\documentclass[12pt]{book}

\usepackage{moresize}

\usepackage[english]{babel}
\usepackage[utf8]{inputenc}
\usepackage{amssymb}
\usepackage{amsmath}
\usepackage{amsthm}
\usepackage{epstopdf}
\usepackage{hyperref}

\usepackage[dvips]{graphicx}
\DeclareGraphicsExtensions{.pdf,.png,.jpg}

\usepackage{color}

\input pictex.tex
\usepackage{graphicx}
\newcommand{\p}{\rho}

\newcommand{\clo}{\mathrm{S}^1}
\newcommand{\sg}{\mathrm{sign}}

\newcommand{\ger}{\mathcal{G}_+(\mathbb{R},0)}

\newcommand{\s}{\sigma}

\newcommand{\precede}{\preceq}
\newcommand{\Z}{\mathbb{Z}}
\newcommand{\esp}{\hspace{0.1cm}}
\newcommand{\mg}{\mathcal{G}}
\newcommand{\vsp}{\vspace{0.1cm}}
\newcommand{\R}{\mathbb R}
\newcommand{\F}{\mathrm{F}}
\newcommand{\id}{\mathrm{id}}
\newcommand{\pa}{[\![}
\newcommand{\pc}{]\!]}

\theoremstyle{definition}

\newtheorem{thm}{Theorem}[section]

\newtheorem{prop}[thm]{Proposition}

\newtheorem{lem}[thm]{Lemma}

\newtheorem{rem}[thm]{Remark}

\newtheorem{question}[thm]{Question}

\newtheorem{cor}[thm]{Corollary}

\newtheorem{ex}[thm]{Example}

\newtheorem{ejem}[thm]{Example}

\newtheorem{ejer}[thm]{Exercise}

\usepackage{makeidx}
\makeindex

\input epsf

\author{\sc { Bertrand Deroin}\\
{\small CNRS and CY Cergy Paris Universit\'e}\\ \\
{\sc Andr\'es Navas}\\
{\small Universidad de Santiago de Chile}\\ \\
{ \sc Crist\'obal Rivas}\\
{\small Universidad de Chile}}

\date{}

\title{\Huge{Groups, Orders,} {\large and} {\Huge Dynamics}
\vspace{-1cm}

\vspace{2cm}}

\begin{document}

\maketitle

\newpage


\pagenumbering{roman}
\tableofcontents


$\mbox{}$

\newpage


\pagenumbering{arabic}

\noindent{\Large {\bf INTRODUCTION}}

\addcontentsline{toc}{chapter}{INTRODUCTION}

\thispagestyle{empty}

\vspace{0.95cm}

The theory of bi-orderable groups is a venerable subject in Algebra that has been extensively developed over the last century, 
stemming from the seminal works of Dedekind, Hölder, and Hilbert. Although initially less developed, the theory of left-orderable 
groups has gained prominence in recent years for several reasons, primarily related to the discovery of new examples, 
key questions concerning group orders, and the application of group orderability in other contexts.

Regarding the discovery of new examples, many researchers have identified several geometrically significant groups that 
are (at least partially) left-orderable. These include braid groups \cite{DDRW}, certain groups of contact diffeomorphisms 
\cite{polte,giroux}, right-angled Artin groups \cite{DuchampThibon}, and, up to finite index, the fundamental groups of 
closed hyperbolic 3-manifolds \cite{Agol,BergeronWise,HaglundWise,KahnMarkovic}. 
Concerning structural questions, key open problems include determining whether certain families of groups are 
left-orderable, notably, groups possessing Kazhdan's property (T) \cite{BHV,kazhdan}. A related question asks 
whether bi-orderable groups are a-(T)-menable \cite{valette}. Furthermore, important problems, such as the 
Boyer-Gordon-Watson conjecture linking left-orderability to Heegaard L-spaces \cite{BGW,Ju}, remain unresolved. 
Finally, the field has attracted diverse mathematical interest due to its growing impact and application. This is 
evidenced by the discovery of new, relevant examples within the framework of orderability, such as finitely-generated, 
left-orderable, simple groups \cite{hyde-lodha-simple, hyde-lodha-finitely presented, MT} and left-orderable, 
finitely presented, non-amenable groups without free subgroups \cite{lodha,monod}. Orderability has also proven 
instrumental in solving long-standing problems (for example, Neumann's conjecture \cite{mineyev} and 
Wiegold's problem on perfect groups \cite{lodha-nuevo}) and providing new 
insight into the study of codimension-1 foliations (particularly $\mathbb{R}$-covered foliations).

There are several classical references on the topic of (left-)orderable groups, such as \cite{botto,glass,kopytov}. 
However, these primarily focus on algebraic aspects, and both the terminology and the approach are somewhat dated 
with respect to the most important lines of research nowadays. This monograph arose from the necessity of placing the 
classical results of the theory in a new and modern perspective to provide a solid background for pursuing research on 
the subject. Quite naturally, in many sections, there is a large intersection of this text with the aforementioned books. 
However, our presentation is new in two aspects. 
On the one hand, we strongly emphasize examples of both groups and phenomena in which orderability plays a 
crucial role. Naturally, many of these examples rely on geometric, combinatorial, topological, or even probabilistic 
insights, but the overarching principle is that orderable groups are dynamical objects, a fact we exploit whenever possible. 
On the other hand, the order we have chosen for the topics is neither logical nor historical, departing from most (if not all) 
known references for the theory. Though this may cause some minor problems for reading (our exposition is not always `linearly ordered"), 
we think that ultimately, this presentation is more appropriate for our purposes. Indeed, the interest of the subject matter lies precisely in its 
several facets —algebraic, geometric, dynamical— not easily reconciled, requiring us to shift from perspective to perspective depending on 
the demands of the particular questions we address. Moreover, our exposition allows the possibility of reading later sections without 
necessarily having mastered the details of earlier ones (though consulting notation will always be necessary).

We stress that there are many other texts that may be considered for complementary reading. In particular, we mention a couple of 
remarkable recent books: {\em Ordering Braids}, by Dehornoy, Dynnikov, Rolfsen, and Wiest \cite{DDRW}, and {\em Ordered Groups 
and Topology}, by Clay and Rolfsen \cite{CR}. These specific works are the reason that topics relating orderability to braid groups 
and low-dimensional topology, including many important examples, are not fully developed in this book. Another subject not treated here is 
circular orders on groups. This topic has attracted significant interest in recent years, and important contributions 
can be found in \cite{baik-s, calegari-circular, calegari-alden, clay-mann-rivas, ghys, ghys-reseaux, malicet-mann-rivas-triestino, mann-rivas, matsu-circ-1, matsu-circ-2}. Despite this, many questions remain open. 
Last but not least, three relatively recent books touch on closely related matters: {\em Foliations and the Geometry of 3-Manifolds}, by Calegari \cite{calegari-libro}; 
{\em Groups of Circle Diffeomorphisms}, by the second-named author of this book \cite{yo}; and the very recent monograph {\em Structure and Regularity of 
Group Actions on One-Manifolds}, by Kim and Koberda, which covers similar (yet updated) topics \cite{KK-book}.

\vsp
  
Let us next briefly describe the content of each chapter of this book.
  
In Chapter 1, we review the basic definitions and treat several relevant examples, such as solvable groups, 
Thompson's groups, and free groups. We also discuss some of the general properties of groups admitting 
orders with different invariance properties, as well as certain closely related combinatorial issues. We close 
the chapter with a result by Gromov concerning the (linear) isoperimetric profile of left-orderable groups.

In Chapter 2, we show that, besides the fact that many groups admit left-orders, they generally admit a multitude of them. 
To better study this phenomenon, we introduce the notion of the space of left-orders associated with a left-orderable group, 
and we discuss some of its properties. As a concrete example, we treat the case of the free group from several points of view. 
Moreover, we present examples of left-orderable groups having uncountably many left-orders but whose associated spaces of 
left-orders contain isolated points, and we give Tararin's description of the groups admitting only finitely many left-orders.

In Chapter 3, we place some classical results of the theory within a dynamical framework and present new developments 
achieved via this approach. We begin with the classical H\"older theorem characterizing group left-orders satisfying an 
Archimedean-type property. We then move to the theory of Conradian left-orders. We first review the classical approach 
by Conrad, and then we provide an alternative dynamical approach, which leads to applications in the study of the topology 
of the space of left-orders. In particular, we give a complete characterization of the groups admitting finitely many Conradian 
left-orders, as well as a description of the space of left-orders of countable solvable groups. We close the chapter with a general 
decomposition of the space of left-orders of finitely-generated, left-orderable groups into three canonical subsets according to 
their dynamical properties. 

Chapter 4 is devoted to several recent results relying on techniques with a probabilistic flavor. We begin with Witte Morris' theorem, 
asserting that left-orderable, amenable groups are locally indicable. We also provide the details of Monod's remarkable construction 
of non-amenable left-orderable groups without rank-two free subgroups ({\em i.e.}, left-orderable counter-examples to the von Neumann 
conjecture), and we also present the explicit (and beautiful~!) finitely-generated version of this due to Lodha and Moore. We then consider 
actions by almost-periodic homeomorphisms and provide a construction of a space involving all of them, which somewhat replaces the 
space of left-orders. Using this, we offer an alternative proof of Witte Morris' theorem mentioned above. Next, we study 
random walks on finitely-generated, left-orderable groups, showing recurrence-type properties and the existence of harmonic 
functions of dynamical origin. More importantly, we explain how probabilistic arguments provide canonical coordinates for 
almost-periodic actions on the line. This is one of the main ingredients for the recent solution of the long-standing problem 
concerning the non-left-orderability of lattices in higher-rank simple Lie groups, obtained by the first-named author together 
with Hurtado \cite{DH}. The proof of this result is not developed here since it also heavily relies on the theory of semi-simple 
Lie groups and symmetric spaces, particularly on the ideas coming from Margulis super-rigidity and Zimmer's program, which 
are too far from the aim of this book. Nevertheless, we hope that this chapter will provide a solid background to immediately 
engage with the reading of \cite{DH}. We conclude with a recent construction by Matte Bon and Triestino, who provide groups of 
piecewise-dyadic homeomorphisms of the line that, on the dynamical side, are almost-periodic, and on the algebraic side, are 
left-orderable, simple, and finitely-generated.

\vsp

Most of the results presented in the book are entirely self-contained. However, some basic knowledge of geometric group theory is desirable 
to fully appreciate the beauty and depth of some of the ideas. In any case, the necessary background not fully developed here can be easily 
grasped by looking at basic books or even on the internet with the right keywords. The text is also complemented with many exercises, which 
sometimes correspond to minor results in the literature. More importantly, several open problems are spread throughout the text. We hope 
that some of these are of genuine mathematical value and will inspire future research in the subject. (A complementary list of open questions, 
mostly concerning classical achievements of the theory, may be found in \cite{glass-lista}; see also \cite[Chapter XVI]{DDRW}.) It is worth 
mentioning that, due to the long delay in the publication of this book since its first appearance online, some of these problems have been (at least partially) solved. However, we decided to still mention them and provide a short discussion and the corresponding references in each case where due.

This text started growing from notes that the second-named author wrote
for mini-courses at the Third Latin American Congress of Mathematicians (2009),  
the Uruguayan Colloquium of Mathematics (2009), and the School Young Geometric
Group Theory II (Haifa, 2013). A first draft was posted on arXiv in 2014; unfortunately, this 
contained many misprints and a couple of small mathematical mistakes that were 
pointed out by several colleagues and are (hopefully) fixed in this version. The three authors 
would like to express their gratitude to the anonymous referees and  
 L.~Bartholdi, J.~Brum, D.~Calegari, M.~Calvez, A.~Clay, Y.~de~Cornulier, P.~Dehornoy, A.~Erschler, 
 \'E.~Ghys, A.~Glass, R.~Grigorchuk, F.~Haglund, T.~Hartnick, S.~Hurtado, J.~ Hyde, T.~Ito, D.~Kielak, 
 S.~ Kim, V.~Kleptsyn, T.~Koberda, J.~Lodha, K.~Mann, I.~Marin, N.~ Matte Bon, G.~Metcalfe, 
 D.~Witte Morris, L.~Paris, F.~Paulin, D.~Rolfsen, F.~le~Roux, Z.~
$\check{\mathrm{S}}$uni\'c, R.~Tessera, M.~Triestino, L.~Vendramin and B.~Wiest,
as well as all the participants of the meetings Orderable Groups, held at Caj\'on del Maipo
(2014), Ordered Groups and Rigidity in Dynamics and Topology, held at Casa 
Matem\'atica Oxaca (2019), and Big Mapping Class Groups and Diffeomorphism Groups, held at 
CIRM, Luminy (2022), for valuable discussions, comments, corrections and suggestions. 
We also thank the Research in Groups Program of CY Advanced Studies for hosting the 
authors during the final stage of the writing. Last but not least, we would like to express our 
deep gratitude to Ina Mette for her extraordinary work and patience during the editing process.

All the three authors acknowledge the funding from the CONICYT PIA 1103 Project DySyRF  
(Center of Dynamical Systems and Related Fields). 
B.~Deroin was also partially supported by ANR-08-JCJC-0130-01, ANR-09-BLAN-0116, ANR-13-BS01-0002 and 
CY Initiative (ANR-16- IDEX-0008).
A. Navas would like to acknowledge the support and 
hospitality of the ERC starting grant 257110 RaWG / Institut Henri Poincar\'e, 
as well as the support of the the CONICYT PIA 1415 Project (Geometry at the Frontier). 
C.~Rivas was also partially supported by a CONICYT's grant Inserci\'on  79130017 and FONDECYT 1241135.


\newpage

$\vsp$

\vspace{8cm}

$\hfill\mbox{\Large To Ana\"{\i}s, Rapha\"el,}$ \\

\vsp\vsp

$\hfill\mbox{\Large Corita, Ch\'{\i}o, Nachito,}$ \\

\vsp\vsp

$\hfill\mbox{\Large Adela, and Benito.}$


\newpage

\noindent{\Large{\bf Notation}}

\vspace{0.4cm}


Because of our dynamical approach, group elements will often be denoted by
the letters $f,g,h$, and compositions are often considered from right to left. 
We sometimes use the letters $u,v,w$, as well as $a,b,c$, particularly when  
dealing with specific groups (as for example free groups, fundamental
groups of surfaces, braid groups, etc). Following the classical notation, 
we also use the letter $\sigma$ to denote certain elements within braid groups. 
Groups are generically denoted by $\Gamma$, though
$G,H$ are sometimes used, as well as $C$ (for convex subgroups), $R$
(for nilpotent radicals) and $T$ (for Tararin groups). In many cases, we implicitly
assume that the groups in consideration are nontrivial; for left-orderable
groups, this is equivalent to being infinite. Similarly, when we work with actions
on the line, we implicitly assume that these are actions by orientation-preserving
homeomorphisms. In general, real-valued function will be denoted by $\phi,\psi$, whereas
group representations by $\Phi,\Psi$. This notation is coherent because certain functions to be
constructed will turn out to be representations, that is, homomorphisms into the additive group of 
real numbers.

Furthermore, when $\Gamma$ is  generated by a finite set $\mathcal{G}\subseteq \Gamma$, we say that 
$\Gamma$ is finitely-generated. In this case, we denote by $||g||$ the {\bf{\em word-length}} of $g \in \Gamma$, 
which by definition is the minimum $m$ for which $g$ which can be written in the form \esp $g = g_{i_1} g_{i_2} \cdots g_{i_m}$,
\esp with $g_{i_j} \in \mathcal{G}$. We also denote by $B_n(id) := \{ f\in \Gamma\mid ||f||\leq n\}$  the {\bf{\em ball
of radius n}} (centered at $id$) with respect to $\mathcal{G}$. Usually, we will consider $\mathcal{G}$ to be 
{\bf{\em symmetric}}, meaning that $g^{-1} \!\in\! \mathcal{G}$ for all $g \in \mathcal{G}$.

\index{word length} \index{word metric}

\vsp

Below, we list some other notation used throughout this text:

\vsp

\noindent
{\boldmath
$id$}:
the trivial (identity) element of a group.
\index{Group!identity}

\vsp

\noindent
{\boldmath
$\langle g_1,g_2,\ldots \rangle$}:
the group generated by $g_1, g_2,\ldots$
\index{Group!group}

\noindent
{\boldmath
$\langle g_1,g_2,\ldots \rangle^+$}:
the semigroup generated by $g_1, g_2,\ldots$
\index{Group!semigroup}

\noindent
{\boldmath
$\Gamma_1 * \Gamma_2$}:
the (non-Abelian) free product of $\Gamma_1$ and $\Gamma_2$.

\vsp

\noindent
{\boldmath
$\mathcal{LO}(\Gamma)$}:
the space of left-orders of $\Gamma$.

\vsp

\noindent
{\boldmath
$\mathcal{BO}(\Gamma)$}:
the space of bi-orders of $\Gamma$.

\vsp

\noindent
{\boldmath
$\mathcal{CO}(\Gamma)$}:
the space of Conradian orders of $\Gamma$.

\vsp

\noindent
{\boldmath
$C_{\preceq}(\Gamma)$}:
the Conradian soul of a left-ordered group $(\Gamma,\preceq)$.
\index{Conradian!soul}

\vsp

\noindent
{\boldmath
$P_{\preceq}^+$}:
the positive cone of a left-order $\preceq$.

\vsp

\noindent
{\boldmath
$\mathbb{N}$} :
$ \!\!\! =  \{1, 2, \ldots\}$.

\vsp

\noindent
{\boldmath
$\mathbb{N}_0$} :
$ \!\!\! =  \{0, 1, 2, \ldots\}$.

\vsp

\noindent
{\boldmath
$(\mathbb{R},+)$}:
the group of real numbers under addition.

\vsp

\noindent
{\boldmath
$\mathbb{R}^*$}:
the group of positive real numbers under multiplication.

\vsp

\noindent
{\boldmath
$\mathrm{Aff}(\mathbb{R})$},
{\boldmath
$\mathrm{Aff}_+ (\mathbb{R})$}:
the group of affine homeomorphisms of the real line and
the subgroup of orientation-preserving ones, respectively.



\vsp

\noindent
{\boldmath
$\mathrm{PSL}(2,\mathbb{R})$},
{\boldmath
$\widetilde{\mathrm{PSL}}(2,\mathbb{R})$}:
the group of orientation-preserving projective homeomorphisms of the circle and
the group of their lifts to the real line, respectively.

\vsp

\noindent
{\boldmath
$\mathrm{Homeo}_+ (\mathbb{R})$, $\mathrm{Homeo}_+ (\mathbb{S}^1)$}:
the group of orientation-preserving homeomorphisms of the line and the circle, respectively.

\vsp

\noindent
{\boldmath
$\widetilde{\mathrm{Homeo}}_+ (\mathbb{R})$: 
the group of lifts to the real line of the orientation-preserving homeomorphisms of the circle.}

\vsp

\noindent
{\boldmath
$\mathrm{PAff}_+(\mathbb{R})$, $\mathrm{PAff}_+(\mathbb{S}^1)$}: the group of orientation-preserving homeomorphisms that are piecewise affine of the real line and the circle respectively.
\index{Group!piecewise affine}

\vsp

\noindent
{\boldmath
$\mathrm{PP}_+(\mathbb{R})$}: the group of orientation-preserving homeomorphisms 
of the real line that are piecewise in $\mathrm{PSL}(2,\mathbb{R})$.
\index{Group!piecewise projective}

\vsp

\noindent
{\boldmath
$\mathrm{F}$}:
Thompson's group of piecewise-affine,
dyadic, orientation-preserving homeomorphisms of the interval.
\index{Thompson's group}

\vsp

\noindent
{\boldmath
$\mathbb{F}_n$}:
the free group on $n$ generators (we will implicitly assume that $n \geq 2$).
\index{Group!free}

\vsp

\noindent
$\mathbb{B}_n$:
the braid group in $n$ strands.
\index{Group!braid}

\vsp

\noindent
{\boldmath
$\mathrm{P} \mathbb{B}_n$}:
the pure braid group in $n$ strands.
\index{Group!pure braid}

\vsp

\noindent
{\boldmath
$BS(1,\ell)$}:
the Baumslag-Solitar group \esp $\langle a,b \!: aba^{-1} = b^{\ell} \rangle$.
\index{Group!Baumslag-Solitar}

\vsp

\noindent
{\boldmath
$G_{m,n}$}:
the torus-knot group $\langle a,b \!: a^m = b^n \rangle$.

\vsp

\noindent
{\boldmath
$\sqcup $}:
the union symbol for disjoint sets.


\chapter{SOME BASIC AND NOT SO BASIC FACTS}

\section{General Definitions}
\label{general}

\index{Order!bi-order}
\index{Order!left-order}
\index{Order!right-order}

\hspace{0.45cm} An order relation $\preceq$ on a group $\Gamma$ is
{\bf{\em left-invariant}} (resp. {\bf{\em right-invariant}}) if for all
$g,h$ in $\Gamma$ such that $g \preceq h$, one has $fg \preceq fh$
(resp. $gf \preceq hf$) for all $f \!\in\! \Gamma$. The relation is
{\bf\textit{bi-invariant}} if it is simultaneously invariant by the left
and by the right. To simplify, we will use the term {\bf{\em left-order}}
(resp. {\bf{\em right-order}}) for a left-invariant total order on a group, 
and {\bf{\em bi-order}} for a bi-invariant total order.
We will say that a group $\Gamma$ is {\bf{\em left-orderable}}
(resp. {\bf{\em right-orderable}})  if it admits a total
order that is invariant by the left (resp. by the
right), and that is is {\bf {\em bi-orderable}} if it admits 
a total order that is simultaneously invariant by the left and right.

\begin{small} \begin{ejem} \label{productos}
Clearly, every subgroup of a left-orderable group is left-orderable. More interestingly,
an arbitrary product $\Gamma$ of left-orderable groups $\Gamma_{\lambda}$ ($\lambda \in \Lambda$) is left-orderable.
(This also holds for bi-orderable groups with the same proof.) Indeed, fixing a 
total well-order on the set of indices $\Lambda$ and a left-order $\preceq_{\lambda}$
on each $\Gamma_{\lambda}$, let $\preceq$ be the associated {\bf{\em lexicographic}}
order. This means that $(g_{\lambda}) \prec (h_{\lambda})$ if the least 
$\lambda \in \Lambda$ such that $g_{\lambda} \neq h_{\lambda}$ satisfies
$g_{\lambda} \prec_{\lambda} h_{\lambda}$. It is easy to check that
$\preceq$ is total and left-invariant.
\end{ejem} \end{small}
\index{Order!lexicographic}


\subsection{Positive and negative cones}
\label{general-1}

\index{Cone!positive}
\index{Cone!negative}

\hspace{0.45cm} If $\preceq$ is an order on $\Gamma$,
then $f \!\in\! \Gamma$ is said to be {\bf\textit{positive}}
(resp. {\bf {\em negative}}) if $f \succ id$ (resp. $f \prec id$). Note that
if $\preceq$ is total, then every nontrivial element is either positive or
negative, and $f \succ id$ if and only if $id \succ f^{-1}$.  
(To get the second inequality, it suffices to multiply on the left 
each term of the first one by $f^{-1}$.) 
Moreover, if $\preceq$ is left-invariant and
$P^{+} = P^{+}_{\preceq}$ (resp. $P^{-}= P^{-}_{\preceq}$)
denotes the set of positive (resp. negative) elements in $\Gamma$
(usually called the {\bf{\em positive}} (resp. {\bf{\em negative}})
{\bf{\em cone}}), then $P^{+}$ and $P^{-}$
are semigroups, and $\Gamma$ is the disjoint union of $P^{+}$, $P^{-}$ and $\{id\}$. 
(Recall that a {\bf {\em semigroup}} is a set endowed with an associative  multiplication.)

Conversely, to every decomposition of $\Gamma$ as a disjoint union of semigroups
$P^{+}$, $P^{-}$, and $\{id\}$ such that $P^{-} = (P^+)^{-1} := \{f\! : f^{-1} \in P^{+} \}$,
there corresponds a left-order $\preceq$ defined by $f \prec g$ whenever $f^{-1}g \in P^{+}$.
Note that $\Gamma$ is bi-orderable exactly when these semigroups may be taken
invariant by conjugation ({\em i.e.}, when they are {\em normal} subsemigroups).

\begin{small}\begin{rem}
The characterization in terms of positive and negative cones shows
immediately the following: If $\preceq$ is a left-order on a group $\Gamma$, then the
{\bf {\em reverse order}} $\overline{\preceq}$ defined by \esp $g \esp \overline{\succ} \esp id$
\esp if and only if \esp $g \prec id$ \esp is also left-invariant and total.
\label{la-primera}
\end{rem}

\begin{rem}
Given a left-order $\preceq$ on a group $\Gamma$, we may define an order $\preceq^*$ by letting
$f \preceq^* g$ whenever $f^{-1} \succ g^{-1}$. Then the order $\preceq^*$ turns out to be
{\em right-invariant}. One can certainly go the other way around, producing left-orders from
right-orders. As a consequence, a group is left-orderable if and only if it is right-orderable.
Since our view is mostly dynamical, we prefer to work with left-orders, yet most of the 
classical literature on the subject deals with right-orders.
\end{rem}

\begin{rem} 
It is worth mentioning that $f \prec g$ for a left-order 
$\preceq$ does not imply that $f^{-1} \succ g^{-1}$. 
(See Example \ref{ejem:no-inverse} on this point.)
Actually, it is easy to check that the implication 
$$f \prec g \implies f^{-1} \succ g^{-1}$$
holds for all elements if and only if $\preceq$ is a bi-order.
\end{rem}\end{small}
\index{Order!bi-order}


\subsection{A characterization involving finite subsets}
\label{general-2}

\hspace{0.45cm} A group $\Gamma$ is left-orderable if and only if for every finite family
$\mathcal{G}$ of nontrivial elements, there exists a choice of ({\em compatible}) exponents 
$\epsilon \!: \mathcal{G} \to \{-1,+1\}$ such that the identity element \esp $id$ \esp does not belong 
to the semigroup generated by the elements $g^{\epsilon(g)}$, $g \in \mathcal{G}$. Indeed, the
necessity of the condition is clear: it suffices to fix a left-order $\preceq$ on $\Gamma$
and choose each exponent $\epsilon (g)$ so that $g^{\epsilon (g)}$ becomes a positive element.
Conversely, assume that for each finite family $\mathcal{G}$ of elements in $\Gamma$
different from the identity, there is such a choice of compatible exponents
$\epsilon \! : \mathcal{G} \to \{-1,+1\}$, and let $\mathcal{X}(\mathcal{G},\epsilon)$
denote the subset of $\{-1,+1\}^{\Gamma \setminus \{id\}}$
formed by the functions \esp $\sg$ \esp satisfying
$$\sg (h) = +1 \quad \mbox{and} \quad \sg (h^{-1}) = -1 \quad \mbox{ for every }
h \in \big\langle g^{\epsilon(g)}, g \in \mathcal{G} \big\rangle^+.$$
(With a slight abuse of notation here and in what follows, 
given a family of group elements $\mathcal{F}$, we let
$\langle \mathcal{F} \rangle^+$ be the semigroup spanned by them.) 
By hypothesis, the set $\mathcal{X}(\mathcal{G},\epsilon)$ is nonempty. Moreover, it is a closed 
subset of $\{-1,+1\}^{\Gamma \setminus \{id\}}$ when this set is endowed with the product topology.

Let $\mathcal{X}(\mathcal{G})$ be the union of all the sets of the form
$\mathcal{X}(\mathcal{G},\epsilon)$ for some choice of compatible exponents
$\epsilon$ on $\mathcal{G}$.
Note that, if \esp $\{\mathcal{X}_i := \mathcal{X}(\mathcal{G}_i),
\esp \esp 1 \leq i \leq n \}$ \esp is a finite family of subsets of this
form, then the intersection \esp $\mathcal{X}_1 \cap \ldots \cap \mathcal{X}_n$
\esp contains the (nonempty) set \esp
$\mathcal{X}(\mathcal{G}_1 \cup \ldots \cup \mathcal{G}_n)$,
\esp and it is therefore nonempty. Since $\{-1,+1\}^{\Gamma \setminus \{id\}}$ is compact,
a direct application of the Finite Intersection Property shows that the intersection
$\mathcal{X}$ of all sets of the form $\mathcal{X}(\mathcal{G})$ is (closed and)
nonempty. Finally, each point in $\mathcal{X}$ corresponds in an obvious way
to a left-order on $\Gamma$.

\medskip

Analogously, one can show that a group is bi-orderable if and only
if for every finite family $\mathcal{G}$ of nontrivial elements,
there exists a choice of exponents $\epsilon \!: \mathcal{G} \to \{-1,+1\}$ such that
\esp $id$ \esp does not belong to the smallest semigroup which simultaneously
satisfies the next two properties:\\

\vspace{0.1cm}

\noindent -- It contains all the elements \esp $g^{\epsilon(g)}$;

\vspace{0.1cm}

\noindent -- For all $f,g$ in the semigroup, both \esp
$fgf^{-1}$ \esp and \esp $f^{-1}gf$ \esp also belong to it.

\vspace{0.1cm}

\noindent We leave the proof to the reader. As a corollary, we obtain that left-orderability
and bi-orderability are {\bf{\em local properties}}; that is, if they are satisfied by
every finitely-generated subgroup of a given group, then they are satisfied by the
whole group. 

Similarly, left-orderability and bi-orderability  are {\bf{\em residual}} properties: if for every nontrivial
element there is a surjective group homomorphism into a group with that property
mapping the prescribed element to a nontrivial one, then the group inherits the
property. Indeed, assuming that a group $\Gamma$ is residually (left- or bi-) orderable, we see from Example \ref{productos} 
that any finite subset $\mathcal G\subseteq \Gamma$ can be mapped injectively into a (left- or bi-) orderable group. Thus, there is a 
compatible choice of exponents for the elements of $\mathcal G$, and so $\Gamma$ is (left- or bi-) orderable by the above criterion.

\begin{small}\begin{ejer} 
\label{ejer-se-puede-extender}
Let $f_1,\ldots, f_k$ be finitely many nontrivial elements of a group $\Gamma$. Suppose that for every finite 
family of elements $g_1,\ldots,g_n$ in $\Gamma$, there exists a left-order (resp. bi-order) on 
$\langle f_1,\ldots,f_k,g_1,\ldots,g_n\rangle$ such that all the $f_i$'s are 
positive. Prove that there exists a left-order (resp. bi-order) on $\Gamma$ for which all the $f_i$'s remain positive. 
\label{ejer-finitos-forced}
\end{ejer}\end{small}


\subsection{Left-orderable groups and actions on ordered spaces}
\label{general-3}

\hspace{0.45cm} If $\Gamma$ is a left-orderable group, then $\Gamma$
acts faithfully on a totally ordered space by order-preserving transformations.
Indeed, fixing a left-order $\preceq$ on $\Gamma$, we may consider the action of
$\Gamma$ by left-translations on the ordered space $(\Gamma,\preceq~\!\!)$.
Conversely, if $\Gamma$ faithfully acts on a totally ordered space $(\Omega,\leq)$ by
order-preserving transformations, then we may fix an arbitrary well-order
$\leqslant_{wo}$ on $\Omega$ and define a left-order $\preceq$ on $\Gamma$
by letting \esp $f \succ id$ \esp if and only if \esp $f (w_f) > w_f$,
\esp where $w_f = \min_{\leqslant_{wo}}\{ w \!: f(w) \neq w\}$.
More generally, if we also have a function \esp $\mathrm{sign} \!:
\Omega \rightarrow \{-,+\}$, \esp we may associate to it the left-order
$\preceq$ for which \esp $f \succ id$ \esp if either $\mathrm{sign}(w_f) = +$
and $f (w_f) > w_f$, or $\mathrm{sign}(w_f) = -$ and $f (w_f) < w_f$. These
left-orders will be referred to as {\bf {\em dynamical-lexicographic}} ones.
\index{Order!dynamically-lexicographic}

\vsp\vsp\vsp

\index{Order!preorder}
\noindent{\bf Left-orders obtained from preorders.} Recall that a
{\bf{\em preorder}} $\preceq$ on a group $\Gamma$ is a reflexive and
transitive relation for which both $f \preceq g$ and $g \preceq f$ may hold
for different $f,g$. The existence of a total, left-invariant preorder
is equivalent to the existence of a semigroup $P$ (containing
the identity) such that $P \cup P^{-1} = \Gamma$. Indeed,
having such a $P$, one may declare $f \preceq g$ if and
only if $f^{-1}g \in P$. Conversely, a preorder $\preceq$
as above yields  the semigroup $P := \{g \!: g \succeq id\}$.
Using the dynamical characterization of left-orders, we next show
that if $\Gamma$ admits sufficiently many total preorders so that
different elements can be ``distinguished'', then it is left-orderable.

\vsp

\begin{prop} \label{preorders}
{\em Let $\Gamma$ be a group and $\{ P_{\lambda}$,
$\lambda \in \Lambda\}$ a family of subsemigroups such that:}

\noindent (i) \esp $P_{\lambda} \cup P_{\lambda}^{-1} = \Gamma$,
{\em for all} $\lambda \in \Lambda$;

\noindent (ii) \esp {\em The intersection} \esp
$P := \bigcap_{\lambda \in \Lambda} P_{\lambda}$
{\em satisfies} $P \cap P^{-1} = \{id\}$.

\noindent {\em Then $\Gamma$ is left-orderable.}
\end{prop}

\noindent{\bf Proof.} For each $\lambda \in \Lambda$, let
$\Gamma_{\lambda} = P_{\lambda} \cap P_{\lambda}^{-1}$. Fix
a total order on the set of indices $\Lambda$, and let $\Omega$ be
the space of all cosets $g \Gamma_\lambda$, where $g \in \Gamma$ and
$\lambda \in \Lambda$. Define an order $\leq$ on $\Omega$ by letting
$g \Gamma_{\lambda} \leq h \Gamma_{\lambda'}$ if either $\lambda$ is smaller
than $\lambda'$, or $\lambda = \lambda'$ and $g^{-1} h \in P_{\lambda}$ (this
does not depend on the chosen representatives $g,h$). By property (i), this order
is total.
The group $\Gamma$ acts on $\Omega$ by  $f (g\Gamma_\lambda) = fg \Gamma_{\lambda}$.
This action preserves $\leq$. Moreover, if $f$ acts trivially, then $f$ lies in
$\Gamma_{\lambda}$ for all $\lambda$. Hence, by property (ii) above, $f \in
\bigcap_{\lambda \in \Lambda} \big( P_{\lambda} \cap P_{\lambda}^{-1} \big)
= P \cap P^{-1} = \{id\}$. This shows that the action is faithful,
hence $\Gamma$ is left-orderable. $\hfill\square$

\vsp\vsp

\index{Order!preorder}
\begin{small}
\begin{ejer} 
Let $\Gamma$ be a group.

\noindent 
(i) If $\Gamma$ is endowed with a left-invariant total preorder $\preceq$, show that  
$H \!:=\! \{ h \!: id \preceq h \preceq id\}$ is a subgroup of $\Gamma$. Show also that $\preceq$ 
induces a total order on the space of classes $\Gamma / H$, which is left-invariant under the 
action of $\Gamma$.

\noindent 
(ii) Conversely, show that every total order on the set of classes $\Gamma / H$ with respect to 
a subgroup $H$ induces a total preorder on $\Gamma$, and that this is left-invariant if and only if 
the $\Gamma$-action on $\Gamma / H$ preserves the order on it.
\end{ejer}

\begin{ejer} \label{preorders-orders}
Let $P := \{ g \!: g \succeq id\}$ be the semigroup of non-negative
elements of a left-invariant total preorder $\preceq$ on a group $\Gamma$. For each
$h \in \Gamma$, let $P_{h} := \{h^{-1}gh \! : g \in P \}$.

\noindent (i) Show that each $P_h$ induces a total preorder on $\Gamma$.

\noindent (ii) Let $H := \bigcap_{h \in \Gamma} ( P_h \cap P_{h}^{-1})$. Show
that $H$ is a normal subgroup of $\Gamma$.

\noindent (iii) Show that $\Gamma / H$ is a left-orderable quotient of $\Gamma$, which 
is nontrivial provided there are at least two nonequivalent elements for $\preceq$.

\noindent{\underline{Hint.}} Although everything can be directly checked, a dynamical view
proceeds as follows: the quotient space $\Gamma / \!\! \sim$ obtained by idendification of
$\preceq$-equivalent points ({\em i.e.}, elements $f,g$ such that $f \preceq g \preceq f$)
is totally ordered, the group $\Gamma$ acts on it by left-translations preserving this order, 
and the kernel of this action corresponds to the subgroup $H$.

\end{ejer}\end{small}

\vsp\vsp

The analogue of the preceding proposition does not hold for partial left-orders. (Recall that a 
{\bf{\em{partial order}}} is a reflexive, transitive, and antisymmetric relation on a set.) Indeed, in 
\S \ref{no-torsion}, we will see many examples of torsion-free groups that are not left-orderable.  
However, these groups admit many partial orders, as shown in the following exercise. 

\begin{small}\begin{ejer} Show that a group is {\em torsion-free} if and only if
it admits a family $\{\preceq_{\lambda} : \lambda \in \Lambda \}$ of partial,
left-invariant orders such that, for each $f \neq g$, there exists
$\lambda \in \Lambda$ satisfying $g \prec_{\lambda} f$.

\noindent{\underline{Hint.}} If $\Gamma$ acts (faithfully) on a set $X$ and $Y \subset X$
has trivial stabilizer, then one may define a partial, left-invariant left-order $\preceq$
on $\Gamma$ by letting $h \succ id$ if and only if $h (Y) \subset Y$. If $\Gamma$ is
torsion-free and $f \neq g$, then this procedure yields a partial left-order for which 
$f \succ g$ by letting $X := \Gamma$ and $Y := \{ h, h^2, \ldots \}$, where $h := g^{-1}f$. 
\end{ejer}\end{small}

\vsp\vsp

\noindent{\bf On group actions on the real line.}
For {\em countable} left-orderable groups, one may take the real line as
the ordered space on which the group acts. (The first reference we found
on this is \cite{holland-dos}; a more modern one is \cite{ghys}.)

\vsp

\begin{prop} \textit{Every left-orderable countable group faithfully
acts on the real line by orientation-preserving homeomorphisms.}
\label{recta}
\end{prop}

\noindent{\bf Proof.} Let $\Gamma$ be a countable group admitting
a left-invariant total order $\precede$. Choose a numbering
$(g_i)_{i \geq 0}$ for the elements of $\Gamma$, set $t(g_0) \! = \! 0$,
and define $t(g_k)$ by induction in the following way: assuming that
$t(g_0),\ldots,t(g_i)$ have been already defined, if $g_{i+1}$ is larger 
(resp. smaller) than all $g_0,\ldots,g_i$ then let $t(g_{i+1})$ be
$\mathrm{max} \{ t(g_0),\ldots,t(g_i) \} + 1$ (resp. $\mathrm{min}
\{ t(g_0),\ldots,t(g_i) \} - 1$), and if $g_m \prec g_{i+1}
\prec g_n$ for some $m,n$ in $\{0,\ldots,i\}$ and no $g_j$
is between $g_m$ and $g_n$ for any $0 \leq j \leq i$
then put $t(g_{i+1}) := (t(g_m) + t(g_n))/2$.

Note that $\Gamma$ acts naturally on $t(\Gamma)$ by $g (t(g_i)) = t(gg_i)$.
We leave to the reader to check that this action extends continuously to the closure
of the set $t(\Gamma)$. Finally, one can extend the action to the entire line by
extending the maps $g$ affinely to each interval in the complement of the
closure of $t(\Gamma)$. $\hfill\square$

\vsp

\begin{small}\begin{rem}  The choice of midpoints in the 
construction above was done to ensure continuity. Many other choices actually 
work, but not arbitrary ones. The important property is the following: for each increasing sequence
of elements $g_1 \prec g_2 \prec \ldots$ smaller than a certain $g$, if every element $h \prec g$ is 
eventually smaller than some $g_n$, then $t(g_n)$ converges to $t(g)$.
\end{rem}\end{small}

\vsp

It is worth analyzing the preceding proof carefully. If $\precede$
is a left-order on a countable group $\Gamma$ and $(g_i)_{i \geq
0}$ is a numbering of the elements of $\Gamma$, then the action of
$\Gamma$ on $\mathbb{R}$ constructed in this proof will be called 
the (associated) {\bf{\em dynamical realization}}. It is easy to
see that this realization has no global fixed point (unless $\Gamma$
is trivial). Moreover, if $f$ is an element of $\Gamma$ whose
dynamical realization has two fixed points $a\!<\!b$ (which may be
equal to $\pm \infty$) and has no fixed point in $]a,b[$, then there
must be some point of the form $t(g)$ inside $]a,b[$. Finally, it is
not difficult to show that the dynamical realizations associated to
different numberings of the elements of $\Gamma$ are all
topologically conjugate.\footnote{A group representation (action) 
$\Phi_1 \!: \Gamma
\rightarrow \mathrm{Homeo}_+(\mathbb{R})$ is {\bf{\em topologically
conjugate}} to $\Phi_2$ if there
exists an orientation-preserving homeomorphism  $\varphi$ of the real
line onto itself such that \esp $\varphi \circ \Phi_1 (g) = \Phi_2 (g)
\circ \varphi$, \esp for all $g \in \Gamma$. Note that 
conjugacy classes yield an equivalence relation; see \S \ref{section-semiconjugacy} 
for more on this.} Therefore, we can
speak of any dynamical property of the dynamical realization
without referring to a particular numbering.

\begin{small}

\begin{rem} \label{empieza-en-cero}
Throughout the text, in most cases we will assume that, in the dynamical 
realization, our numbering of group elements starts at $g_0 := id$, which yields $t(id)=0$.
\end{rem}

\begin{ejer}

Show that for any dynamical realization of a left-order, the set of points in the real line with 
a free orbit is $G_{\delta}$-dense (that is, it contains a countable intersection of dense open subsets). 
\end{ejer}

\index{Order!dynamically-lexicographic}
\begin{rem} Note that, to define a dynamical-lexicographic left-order on the group
$\mathrm{Homeo}_+(\mathbb{R})$ of orientation-preserving homeomorphisms of the
real line, it is not necessary to order all the points in $\mathbb{R}$: it is enough
to consider a well-order on a dense subset (in particular, a dense sequence suffices).
Clearly, $\mathrm{Homeo}_+(\mathbb{R})$ admits uncountably
many left-orders of this type. However, there are left-orders
that do not arise in this manner; see Example \ref{no-son-todos}.
\end{rem}

\begin{rem} \label{rem-germs} The group $\ger$ of germs at the origin of
orientation-preserving homeomorphisms of the real line is left-orderable.
Perhaps the easiest way to show this is by using the characterization in
terms of finite subsets above. Let $\hat{g}_1,\ldots,\hat{g}_k$ be nontrivial
elements in $\ger$, and let $g_1,\ldots,g_k$ be representatives of them. Take
a sequence $(x_{n,1})$ of points converging to the origin in the line so that, 
for each $n$, at least one of the $g_i$'s moves $x_{n,1}$. Passing to a
subsequence if necessary, we may assume that, for each $i \!\in\! \{1,\ldots,k\}$,
either $g_i (x_{n,1}) > x_{n,1}$ for {\em all} \esp $n$, or $g_i (x_{n,1}) < x_{n,1}$
for {\em all} \esp $n$, or $g_i (x_{n,1}) = x_{n,1}$ for {\em all} \esp $n$.
In the first case we let $\epsilon_i := +1$, and in the second case we let
$\epsilon_i := -1$. In the third case, $\epsilon_i$ is still undefined.
However, this may happen only for $k-1$ of the $g_i$'s above. For
these elements, we may repeat the procedure by considering another sequence
$(x_{n,2})$ converging to the origin... In at most $k$ steps, all the
$\epsilon_i$'s will be defined. We claim that this choice is compatible.
Indeed, given an element
$\hat{g} = \hat{g}_{i_1}^{\epsilon_{i_1}} \cdots \hat{g}_{i_{\ell}}^{\epsilon_{i_{\ell}}}$,
the choice above implies that \esp
$g_{i_1}^{\epsilon_{i_1}} \cdots g_{i_{\ell}}^{\epsilon_{i_{\ell}}} (x_{n,1}) \geq x_{n,1}$
\esp for all $n$, where the inequality
is strict if some of the $g_{i_j}$'s ``moves'' some of (equivalently,
all) the points $x_{n,1}$ (meaning that $g_{i_j} (x_{n,1}) \neq x_{n,1}$). If this is the case, then 
$\hat{g}$ cannot be the identity. If not, then we may repeat the argument
with the sequence $(x_{n,2})$ instead of $(x_{n,1})$... Proceeding this
way, we conclude that $\hat{g}$ is nontrivial.

A nice consequence of the claim above is that every countable group of germs
at the origin of homeomorphisms of the real line admits a realization (but not
necessarily an ``extension''!) as a group of homeomorphisms of the interval. Note
that, in the opposite direction, $\mathrm{Homeo}_+([0,1])$ embeds into $\ger$. (This 
embedding is not obtained by looking at the germs of elements of $\mathrm{Homeo}_+([0,1])$
near the origin --the homomorphism thus-obtained is not injective--, but by taking infinitely 
many copies of $\mathrm{Homeo}_+([0,1])$ on intervals accumulating the origin). However, 
despite this embedding, Mann proved in \cite{mann} that the groups $\ger$ and $\mathrm{Homeo}_+([0,1])$ 
are non-isomorphic (see also Exercise \ref{ejer:otra} below). Actually, she proved that there is no nontrivial 
homomorphism from $\ger$ into $\mathrm{Homeo}_+([0,1])$. (See Example \ref{ex-rivas} for
another --much simpler--  example of an uncountable left-orderable group that
has no nontrivial action on the real line.)
\end{rem}

\begin{ejer}\label{ejer:otra} Recall that a {\bf \em polish group} is a topological group that admits a 
complete and separable metric that is compatible with the group topology. The group $\mathrm{Homeo}_+([0,1])$ 
endowed with the compact-open topology is easily seen to be polish. Following the steps below (taken from 
the work of Sch\"oner \cite{unpublished}), show that $\ger$ does not admit a polish topology, 
hence it does not embed into $\mathrm{Homeo}_+([0,1])$.

\noindent (i) Show that for every subgroup $H$ of a polish group $G$, there exists a countable 
subgroup $\hat{H}$ of $H$ such that the centralizers of $H$ and $\hat{H}$ in $G$ coincide.

\noindent (ii) Consider the group $H$ of $\ger$ defined as follows: For each positive integer $n$, 
let $\{T_{n,t}\}$ be a (nontrivial) topological flow supported on $I_n := [1/(n+1), 1/n]$ (for instance, the 
flow associated to a nonzero smooth vector field supported on $I_n$). Then define $H$ as the set of  
elements of $\ger$ having a representative that, for each $n \geq 1$, restricts to a time-$t_n$ map of 
te flow $\{ T_{n,t} \}$ for some $t_n \in \mathbb{Q}$.  
(Note that $H$ is isomorphic to the natural image of the direct product $\mathbb{Q}^{\mathbb{N}}$ 
in $\ger$.) Show that for every countable subgroup $\hat{H}$ of $H$, the centralizer of $\hat{H}$ in 
$G$ is strictly larger than that of $H$.

\noindent{\underline{Hint.}} Fix a numbering $h_1,h_2,\ldots$ for representatives of the elements of $\hat{H}$, where each 
$h_i$ acts by a rational time of the flow $\{ T_{n,t} \}$ on every interval $I_n$. For simplicity, assume also that $h_1$ acts 
nontrivially on every $I_i$. For each positive integer $n$, consider the restrictions of $h_1,\ldots,h_n$ to the interval $I_n$. 
These restrictions generate a cyclic group; denote by $T_{n,t_n}$ a generator of it. Let $g_n$ be a nontrivial homeomorphism 
of $I_n$ that commutes with the map $T_{n,\, t_n}$, but not with $T_{n,\, t_n/2}$. Let $g \in \ger$  be the element represented by 
a homeomorphism whose restriction to each $I_n$ coincides with $g_n$. Prove that $g$ centralizes $\hat{H}$, and show that 
$g$ does not centralize $H$ by explicitly exhibiting an element $h \in H \setminus \hat{H}$ that does not commute with $g$.
\end{ejer}
\end{small}

We do not know whether there is an analogue of Remark \ref{rem-germs}  in higher dimensions.

\vsp

\begin{question} Does there exist a finitely-generated group of germs at the origin
of homeomorphisms of the plane having no realization as a group of homeomorphisms
of the plane~?
\end{question}

\vsp

Note that the results and techniques of \cite{calegari-circular} show that
such a group cannot arise as a group of germs of $C^1$ diffeomorphisms. 
However, imposing such a regularity condition for a group action may lead to very
serious algebraic restrictions (see \cite{yo} for a general panorama on this topic, 
\cite{BMNR,cale,tst} for later developments, and  \cite{Cal-Rolf} for examples of 
a related nature). We will come back to this point in \S \ref{s: Lipschitz conjugation}

To close this discussion, let us mention a related recent result of Hyde \cite{hyde} (see also \cite{triestino-hyde}) 
that solves in the negative a question of Calegari \cite{calegari-blog}: the group of orientation-preserving homeomorphisms 
of the disc that are the identity on the boundary is not left-orderable. It is worth pointing out that this group is torsion-free, as 
follows from a classical result of Kerékjártó \cite{kerej} (see also \cite{kolev}).


\subsection{Semiconjugacy in $\mathrm{Homeo}_+(\mathbb R)$}
\label{section-semiconjugacy}

\index{Action!semiconjugate}
\hspace{0.45cm} Two (non necessarily injective) 
group representations $\Phi_1,\Phi_2: \Gamma \to \mathrm{Homeo}_+(\R)$ will be said to be 
{\bf {\em semiconjugate}} if there is a non-decreasing map $\varphi: \R\to \R$ that is proper ({\em i.e.}, the preimage of 
every compact set is bounded or, equivalently --since $\varphi$ is monotone--, the set $\varphi(\R)$ is unbounded in both 
directions), and such that for all $g\in \Gamma$,  
\begin{equation}\label{eq semiconj} 
\varphi \circ \Phi_1(g) = \Phi_2(g) \circ \varphi.
\end{equation}

In most of the literature, in the definition above, one also requires the continuity of the map $\varphi$. 
However, this extra condition causes more problems than it solves. For instance, if we insist on 
continuity, then actions without global fixed points admitting a discrete minimal invariant set may not be 
semiconjugate to a $\Z$-action by translations. However, with our definition of semiconjugacy, these actions 
are always semiconjugate. More importantly, by dropping the continuity assumption, we have the next 
key fact.

\vsp

\begin{prop} {\em Semiconjugacy is an equivalence relation.} 
\end{prop}

\noindent{\bf Proof.} Reflexivity is obvious and transitivity is easy to check. Below we prove symmetry. 
To do this, suppose that (\ref{eq semiconj}) holds. Since $\varphi$ is proper, we may define 
$$\psi (x) := 
\sup \varphi^{-1}((-\infty,x])=\sup\{y\mid \varphi (y) \leq x\}.$$
From the last equality, the fact that $\psi$ is non-decreasing is obvious. Furthermore, 
the properness of $\psi$ easily follows from the properness of $\varphi$.  
Finally, for all $g \in \Gamma$ and all $x \in \R$, we have 
\begin{eqnarray*}
\Phi_1(g)(\psi (x))
&=& \sup \{\Phi_1(g) (y)\mid \varphi (y) \leq \psi(x) \}\\
&=& \sup\{ z\mid \varphi  (\Phi_1(g)^{-1}(z))\leq \psi(x)\} \\
&=& \sup\{ z\mid \Phi_2(g)^{-1}(\varphi(z))\leq x \} \\
&=& \sup\{ z\mid \varphi (z)\leq \Phi_2(g)(x)\}\\
&=& \psi (\Phi_2(g) (x)).  
\end{eqnarray*}
Therefore, $\psi$ satisfies the semiconjugacy relation. $\hfill\square$

\vspace{0.25cm}

Below we list a couple of exercises and one remark concerning the notion of semiconjugacy that we adopt. 
We refer to \cite{KKM} for further developments on this.

\begin{small}
\begin{ejer} Let $\Gamma$ be a countable group of orientation-preserving homeomorphisms of
the real line. Using its action, produce a dynamical-lexicographic order $\preceq$ on $\Gamma$.
Show that the original action is semiconjugate to the dynamical realization of $\preceq$.
Give examples for which this semiconjugacy is not a conjugacy.
\end{ejer}

\begin{ejer} 
Let $(\Gamma,\preceq)$ be a countable left-ordered group and $\Gamma_0$ a
subgroup. Suppose that, in the dynamical realization of $\preceq$, the subgroup $\Gamma_0$ 
acts with no global fixed point (for instance, this happens if $\Gamma_0$ has finite index in 
$\Gamma$). Show that the restriction to $\Gamma_0$ of the dynamical realization of $\preceq$
is semiconjugate to the dynamical realization of the restriction of $\preceq$ to $\Gamma_0$.
\end{ejer}

\begin{rem} Since our definition of semiconjugacy still involves non injective representations, it applies to 
actions of different groups, provided these actions factor throughout the action of the same group. 
\end{rem}
\end{small}


\section{Some Relevant Examples}

\hspace{0.45cm} At first glance, it may seem surprising that many (classes of) torsion-free
groups turn out to be left-orderable. Here, we give a brief discussion of some of them.

\subsection{Abelian and nilpotent groups}
\label{ejemplificando-1}

\hspace{0.45cm} The simplest bi-orderable groups are the torsion-free, Abelian ones. Obviously, there 
are only two bi-orders on $\mathbb{Z}$. The case of $\mathbb{Z}^2$ is more interesting. According
to \cite{robbin,sikora,teh}, there are two different types of bi-orders on $\mathbb{Z}^2$.
Bi-orders of {\em irrational type} are completely determined by an irrational number
$\lambda$. For such an order $\preceq_{\lambda}$, an element $(m,n)$ is positive if and only
if \esp $\lambda m + n$ \esp is a positive real number. Bi-orders of {\em rational type} are characterized
by two data, namely a pair $(x,y) \in \mathbb{Q}^2$ up to multiplication by a positive real number,
and the choice of one of the two possible bi-orders on the ``kernel'' subgroup \esp
$\{(m,n) \!: mx + ny = 0\} \sim \mathbb{Z}$. \esp Thus, an element
\esp $(m,n) \in \mathbb{Z}^2$ \esp is
positive if and only if either \esp $mx + ny$ \esp
is a positive real number, or \esp $mx + ny = 0$ \esp and
$(m,n)$ is positive with respect to the chosen bi-order on the kernel subgroup. 
The set of left-orders on $\mathbb{Z}^2$ naturally
identifies with the Cantor set (see \S \ref{space-left-orders} for more on this).

The description of all bi-orders on $\mathbb{Z}^n$ for larger $n$ continues inductively.
(A good exercise is to show this using the results of \S \ref{conrad-general}.)
For a general torsion-free, Abelian group, recall that the {\bf{\em rank}} (sometimes 
also called the {\bf{\em torsion-free rank}}) is the minimal
dimension of a vector space over $\mathbb{Q}$ in which the group embeds. The reader
should have no problem showing that, in particular, a torsion-free, Abelian group
of rank$\esp \geq 2$ admits uncountably many left-orders.
\index{rank}

\vspace{0.65cm}


\beginpicture
\label{Cono z2}

\setcoordinatesystem units <0.65cm,0.65cm>

\putrule from -2 7 to 8 7
\putrule from -2 -1 to 8 -1
\putrule from 7 -2 to 7 8
\putrule from -1 -2 to -1 8
\putrule from 7 -2 to 7 8
\putrule from 1 -2 to 1 8
\putrule from 3 -2 to 3 8
\putrule from 2 -2 to 2 8
\putrule from 5 -2 to 5 8
\putrule from -2 1 to 8 1
\putrule from -2 3 to 8 3
\putrule from -2 4 to 8 4
\putrule from -2 5 to 8 5
\putrule from 0 -2 to 0 8
\putrule from 6 -2 to 6 8
\putrule from -2  6 to 8 6
\putrule from 4 -2 to 4 8
\putrule from -2 2 to 8 2
\putrule from 0 -2 to 0 8
\putrule from -2 0 to 8 0

\put{$\bf{id}$} at 3 3
\put{} at -8 0

\plot -2 0.05 2.7 2.8 /

\plot 3.3 3.15 8 5.95 /

\put{Figure 1: The positive  cone \index{Cone!positive} of a left-order on $\mathbb{Z}^2$.} at 3 -2.65


\put{$\bullet$} at -1 -1
\put{$\bullet$} at -1 0

\put{$\bullet$} at 0 -1
\put{$\bullet$} at 0 0
\put{$\bullet$} at 0 1

\put{$\bullet$} at 1 0
\put{$\bullet$} at 1 1
\put{$\bullet$} at 1 -1

\put{$\bullet$} at 2 2
\put{$\bullet$} at 2 1
\put{$\bullet$} at 2 0
\put{$\bullet$} at 2 -1

\put{$\bullet$} at 3 2
\put{$\bullet$} at 3 1
\put{$\bullet$} at 3 0
\put{$\bullet$} at 3 -1

\put{$\bullet$} at 4 0
\put{$\bullet$} at 4 1
\put{$\bullet$} at 4 2
\put{$\bullet$} at 4 3

\put{$\bullet$} at 5 0
\put{$\bullet$} at 5 1
\put{$\bullet$} at 5 2
\put{$\bullet$} at 5 3
\put{$\bullet$} at 5 4

\put{$\bullet$} at 6 0
\put{$\bullet$} at 6 1
\put{$\bullet$} at 6 2
\put{$\bullet$} at 6 3
\put{$\bullet$} at 6 4

\put{$\bullet$} at 7 0
\put{$\bullet$} at 7 1
\put{$\bullet$} at 7 2
\put{$\bullet$} at 7 3
\put{$\bullet$} at 7 4
\put{$\bullet$} at 7 5

\put{$\bullet$} at 7 -1
\put{$\bullet$} at 4 -1
\put{$\bullet$} at 5 -1
\put{$\bullet$} at 6 -1
\put{$\bullet$} at 6 -1

\endpicture


\vspace{0.65cm}

\index{Group!isolator}
Torsion-free, nilpotent groups are also bi-orderable. Indeed, let $\Gamma_i$
denote the $i^{\mathrm{th}}$-term of the {\bf{\em lower central series}} of a group
$\Gamma$ (that is, $\Gamma_1 := \Gamma$ and $\Gamma_{i+1} := [\Gamma,\Gamma_i]$),
and let $H_i (\Gamma)$ be the {\bf{\em isolator}} of $\Gamma_i$ defined by
$$H_i (\Gamma) := \big\{ g \in \Gamma \!: g^n \in \Gamma_i \mbox{ for some }
n \in \mathbb{N} \big\}.$$
If $\Gamma$ is {\bf{\em nilpotent}} ({\em i.e.}, if
$\Gamma_{k+1} = \{id\}$ for a certain $k$), then each $H_{i}(\Gamma)$
is a normal subgroup of $\Gamma$, and $H_i(\Gamma) / H_{i+1} (\Gamma)$
is a torsion-free, central subgroup of $\Gamma / H_{i+1}(\Gamma)$ (see
\cite{kru} for the details). Note that, if $\Gamma$ is also torsion-free,
then $H_{k+1} (\Gamma) =\{ id \}$.
\index{Group!nilpotent}
\index{Lower central series}

Let $P_{i}$ be the positive cone of any left-order on
(the torsion-free Abelian group) $H_i(\Gamma) / H_{i+1}(\Gamma)$,
and let $G_{i}$ be the set of elements in $H_{i}(\Gamma)$ that project to an
element in $P_{i}$ when taking the quotient by $H_{i+1}(\Gamma)$. Using the fact
that each $H_{i}(\Gamma) / H_{i+1}(\Gamma)$ is {\em central} in $\Gamma / H_{i+1}(\Gamma)$,
one may easily check that the semigroup
\esp $P := G_{k-1} \cup G_{k-2} \cup \ldots \cup G_1$
\esp is the positive cone of a bi-order on $\Gamma$.

\index{Group!Heisenberg}
\begin{small}
\begin{ejem} \label{ex-ext} The Heisenberg group
$$H = \big\langle f,g,h \! : [f,g]=h^{-1}, [f,h]= id, [g,h]=id \big\rangle$$
is a non-Abelian nilpotent group of nilpotence degree 2. It may also be seen
as the group of lower-triangular matrices with integer entries so that each
diagonal entry equals 1 via the identifications
$$f =
\left(
\begin{array}
{ccc}
1 & 0 & 0 \\
1 & 1 & 0 \\
0 & 0 & 0 \\
\end{array} \right), \quad
g =
\left(
\begin{array}
{ccc}
1 & 0 & 0 \\
0 & 1 & 0 \\
0 & 1 & 1 \\
\end{array} \right), \quad
h =
\left(
\begin{array}
{ccc}
1 & 0 & 0 \\
0 & 1 & 0 \\
1 & 0 & 1 \\
\end{array} \right).$$
Note that the linear action of $H$ on $\mathbb{Z}^3$ fixes the hyperplane
$\{1\} \times \mathbb{Z}^2$ and preserves the lexicographic order on it. The
left-orders on $H$ induced from this restricted action (see S \ref{general-3})
are (total but) not bi-invariant. This example can be seen as a kind of evidence of the
following nice result due to Darnel, Glass, and Rhemtulla \cite{DGR}: If all left-orders
of a left-orderable group are bi-invariant, then the group is Abelian.

Both the sets of left-orders and bi-orders of countable, torsion-free, nilpotent
groups which are not rank-1 Abelian naturally identify with the Cantor set; see Theorem
\ref{nilpotent-cantor} for left-orders and \cite{witte-nilpotent} for bi-orders. Moreover, a remarkable 
theorem of Malcev \cite{malcev} (resp. Rhemtulla \cite[Chapter 7]{botto}) establishes that every 
bi-invariant (resp. left-invariant) {\em partial} order on a torsion-free, nilpotent group can be
extended to a bi-invariant (resp. left-invariant) {\em total} order.
\end{ejem}

\begin{ejer} \label{heisenberg}
Let $\bar{H}$ be the subgroup of the Heisenberg group formed by the matrices
of the form
$$\left(
\begin{array}
{ccc}
1  & 0  &  0 \\
2x & 1  &  0 \\
z  & 2y &  1 \\
\end{array} \right), \quad x,y,z \mbox{ in } \mathbb{Z}.$$

\noindent (i) Show that the commutator subgroup $[\bar{H},\bar{H}]$ 
is formed by the matrices of the form
$$\left(
\begin{array}
{ccc}
1 & 0 & 0 \\
0 & 1 & 0 \\
4z & 0 & 1 \\
\end{array} \right), \quad z \in \mathbb{Z}.$$

\noindent (ii) Conclude that $\bar{H} / [\bar{H},\bar{H}]$ is isomorphic to \esp
$\mathbb{Z}^2 \times \mathbb{Z} / 4 \mathbb{Z}$, \esp hence has torsion.
\end{ejer}

\begin{ejer} Show that a group $\Gamma$ is residually torsion-free nilpotent if and only
if \esp \esp $\bigcap_{i} H_i(\Gamma) = \{ id \}$. \esp\esp (Since bi-orderability is
a residual property, such a group is necessarily bi-orderable.)

\noindent{\underline{Remark.}} Quite surprisingly, torsion-free, residually nilpotent
groups do not necessarily satisfy this property. Actually, such a group may fail to
be bi-orderable; see \cite{BGG}.
\end{ejer}\end{small}


\subsection{Subgroups of the affine group}
\label{ejemplificando-2}

\index{Group!affine}
\hspace{0.45cm} Let $\mathrm{Aff}_+(\mathbb{R})$ denote the group of
orientation-preserving affine homeomorphisms of the real line (the {\bf{\em affine group}},
for short). For each $\varepsilon \neq 0$, a
partial order $\preceq_{\varepsilon}$ may be defined by declaring that $f$ is
positive if and only if $f(1/\varepsilon) > 1/\varepsilon$. This means that
$$P^+_{\preceq_{\varepsilon}} = \Big\{f =
\left(
\begin{array}
{cc}
u & v  \\
0 & 1  \\
\end{array}
\right) \!: \esp u + v \varepsilon > 1 \Big\}.$$
These orders were introduced (in a more algebraic way) by Smirnov in \cite{smirnov}.

For a finitely-generated subgroup $\Gamma$
of $\mathrm{Aff}_+(\mathbb{R})$, the corresponding action on the line has 
(uncountably many) free orbits.
Thus, one may choose $\varepsilon$ so that $\preceq_{\varepsilon}$ is a total order.
If $\preceq_{\varepsilon}$ is only a partial order, one may ``complete'' it so that it 
becomes total (see \S \ref{general-3}). Consequently, non-Abelian subgroups 
of $\mathrm{Aff}_+(\mathbb{R})$ admit uncountably many left-orders.
\index{Group!Baumslag-Solitar}

As a concrete and relevant example,
for each integer $\ell \geq 2$, the {\bf {\em Baumslag-Solitar
group}} $BS(1,\ell) = \langle g,h \!: \esp hgh^{-1} = g^{\ell} \rangle$ embeds into
the affine group by identifying $g$ and $h$ to $x \mapsto x+1$ and $x \mapsto \ell x$,
respectively. (See Exercise \ref{BS-ex} below.) Note that, for an irrational
$\varepsilon \neq 0$, the associated order $\preceq_{\varepsilon}$ is total. If
one chooses a rational $\varepsilon$, then it may happen that $\preceq_{\varepsilon}$
is only a partial order. However, in this case, the stabilizer of the point
$1 / \varepsilon$ is isomorphic to $\mathbb{Z}$, and thus $\preceq_{\varepsilon}$
can be completed to a total left-order of $BS(1,\ell)$ in exactly two different ways.
Observe also that the reverse orders $\overline{\preceq}_{\varepsilon}$ may be retrieved 
by the same procedure but starting with the embedding $g: x \mapsto x-1$ and
$h: x \mapsto \ell x$, and changing $\varepsilon$ by $-\varepsilon$.

\begin{small}\begin{ejer}\label{BS-ex}
Prove that the map from $BS(1,\ell) = \langle g,h \!: \esp hgh^{-1} = g^{\ell} \rangle$ into
the affine group that makes correspond $g$ and $h$ to $x \mapsto x+1$ and $x \mapsto \ell x$,
respectively, is an embedding.

\noindent{\underline{Hint}.} Prove that the conjugates of $g$ commute, and then write every element
of $BS(1,\ell )$ in normal form as a power of $h$ followed by a product of conjugates of $g$.
\end{ejer}

\begin{ejer} Show that every embedding of $BS(1,\ell)$ into $\mathrm{Aff}_+ (\mathbb{R})$ is obtained by
letting $g,h$ correspond, respectively, to any nontrivial translation and an homothety of ratio $\ell$.
\end{ejer}\end{small}

\begin{small}\begin{ex} \label{ex-bi-order-BS}
Still another way to order $BS(1,\ell)$ is as follows: 
$BS(1,\ell)$ can be thought of as the semidirect product
$\mathbb{Z} [\frac{1}{\ell}] \rtimes \Z$ coming from the exact sequence
$$0 \longrightarrow \mathbb{Z} \Big[ \frac{1}{\ell} \Big] \longrightarrow
BS(1,\ell) \longrightarrow \Z \longrightarrow 0.$$
Using this, one may define the bi-orders $\preceq$, $\preceq'$ by letting
$(\frac{m}{\ell^n},k) \succ id$ (resp. $(\frac{m}{\ell^n},k) \succ' id)$
if and only if either $k > 0$, or $k = 0$ and $\frac{m}{\ell^n} > 0$ (resp.
$k>0$, or $k=0$ and $\frac{m}{\ell^n} < 0$). Together with the reverse
orders $\overline{\preceq}$ and $\overline{\preceq}'$, this completes the
list of all bi-orders on $BS(1,\ell)$ (see Example \ref{solo-cuatro}).
\end{ex}\end{small}

\vsp

Another nice group that embeds into the affine group is  $\Z\wr\Z:=\bigoplus_\Z \Z\rtimes \Z$, the wreath product of $\Z$ 
with itself. Here, the conjugation action of $\Z$ on $\bigoplus_\Z \Z$ is by shifting the indexes. In particular, it is not hard 
to see that 
$$\Z\wr\Z\simeq \langle a,b_0 \mid b_i:=a^i b_0 a^{-i} , [b_i,b_j] = id \, \mbox{ for all } i,j \mbox{ in }\Z\rangle.$$
We leave to the reader to show that, if $\lambda\in \R$ is a {\em non-algebraic} ({\em i.e.}, trascendental) number, 
then the identification of $a$ and $b_0$ to $x\mapsto \lambda x$ 
and $x\to x+1$ induces an homomorphic embedding of $\Z\wr\Z$ into $\mathrm{Aff}_+(\mathbb{R})$.

As for $BS(1,\ell)$ (see Example \ref{ex-bi-order-BS} above), we can use the short exact sequence 
$$\{id\}\to \bigoplus_\Z \Z \to \Z\wr\Z\to \Z\to \{id\}$$
to produce many left-orderings on $\Z\wr\Z$. Indeed, any left-order on $\bigoplus_\Z \Z$ can be extended lexicographically to a left-order on $\Z\wr \Z$. 
However, these orders are not always bi-invariant, since for instance $b_0$ and $b_1$ may not have the same sign. (However, all these orders 
enjoy a slightly weaker property, called the Conradian property, that will be extensively studied in \S \ref{general-Conrad}.) To produce a bi-invariant 
order on $\Z\wr\Z$ using this procedure, we need to consider orders on $\bigoplus_\Z \Z$ that are invariant after shifting the indexes. One such 
order is the lexicographic order on $\bigoplus_\Z \Z$:  an element $g=b_i^{n_i} b_{i+1}^{n_{i+1}} \cdots b_{i+k}^{n_{i+k}}$,  
where $i \in \Z$, $k\geq 0$ and $n_i\neq 0$, is declared  to be positive if $n_i > 0$.

\vsp\vsp

The discussion above can be extended to all non-Abelian subgroups of $\mathrm{Aff}_+(\mathbb{R})$.  In the 
terminology of \S \ref{space-left-orders}, the associated spaces of left-orders identify with the Cantor set. The study
of more general solvable left-orderable groups is more involved, yet it crucially relies on the case of affine groups.
We will come back to this point in \S \ref{section finite-rank-solvable} and \S \ref{section solvable (general)}.

\vsp


\subsection{Free and residually free groups}
\label{ejemplificando-3}

\index{Vinogradov's theorem}
\index{Magnus' expansion}
\hspace{0.45cm} The free group $\mathbb{F}_2$ is bi-orderable. (As a consequence, 
since the commutator subgroup $[\mathbb{F}_2,\mathbb{F}_2]$ is isomorphic to 
$\mathbb{F}_{\infty}$, every non-Abelian free group is bi-orderable as well.) 
Although this result is originally due to Shimbireva \cite{sh}, it is sometimes
attributed to Vinogradov \cite{vinogradov}, and more usually to Magnus.  
Below we first sketch Magnus' construction, which covers Shimbireva approach. 
An alternative argument (which actually applies to free products of arbitrary bi-orderable
groups) will be developed  in \S \ref{sub-vinogradov}. 

Consider the (non-Abelian) ring $\mathbb{A} \! = \! \mathbb{Z}
\langle X_0,X_1 \rangle$ formed by the formal power series with integer coefficients in two
independent variables $X_0$, $X_1$. Denoting by $o(k)$ the subset of $\mathbb{A}$ formed by
the elements all of whose terms have degree at least $k$, one easily checks that
$$\mathbb{F} := 1 + o(1) = \big\{1 + S \! : \esp S \in o(1) \big\}$$
is a subgroup (under multiplication)
of $\mathbb{A}$. Moreover, if $f,g$ are (free) generators of $\mathbb{F}_2$, the map
$\Phi$ sending $f$ (resp. $g$) to the element $1 +X_0$ (resp. $1 + X_1$)
in $\mathbb{A}$ extends in a unique way to a homomorphism
$\Phi \!: \mathbb{F}_2 \rightarrow \mathbb{F}$. Note that $\Phi(f^{-1})=\Phi(f)^{-1}$ is the infinite power series $1-X_0+X_0^2-\ldots$ (which lies in $\mathbb{F}$). 
We claim that $\Phi$ is an injective homomorphism.  To see this, note that for $n \in \mathbb N$,
$$\Phi(f^n)=1+ {n\choose 1}X_0+ \ldots +X_0^n \,\,  \text{ and } \,\, \Phi(f^{-n})= 1- {n \choose 1} X_0+ \ldots. $$
Now, for a reduced word $w\in F_2$, for instance $w=f^{n_1}g^{m_1}\ldots f^{n_k}g^{m_k}$, with $n_i,m_i$ in $\Z \setminus \{0\}$, we have that $\Phi(w)$ 
contains the term $n_1m_1 \cdots n_km_k X_0X_1\ldots X_0X_1$, and therefore $\Phi(w)$ is a nontrivial power series (non-commutativity 
between $X_0$ and $X_1$ is crucial for this argument). 

Next, we observe that $\mathbb F$ can be lexicographically ordered. To do this, we first need to order the monomials of degree $k\geq 1$. 
For this, we consider 
$$ \{0,1\}^k =\{\varphi :\{1,\ldots,k\}\to\{0,1\}\},$$ 
and for each $\varphi \in \{0,1\}^k$ we write $X_\varphi := X_{\varphi(1)} X_{\varphi(2)} \ldots X_{\varphi(k)}$. Now, given $\varphi$ 
and $\psi$ in $\{0,1\}^k$, we declare $X_\varphi \prec_{k} X_\psi$ if $\varphi(i) = 0$ and $\psi(i)=1$, where $i\in \{1,\ldots,k\}$ is the least integer where $\varphi$ 
and $\psi$ differ.  With this notation, for a power series $P$ without constant term, we 
have that $P$ belongs to  $o(k)\setminus o(k+1)$ for some $k\geq 1$ and  hence 
$$P=\sum_{\varphi \in \{0,1\}^k} a_\varphi X_\varphi + T,$$
where $T\in o(k+1)\cup \{0\}$. We declare $ P$ to be {\em positive} if  $ a_{\varphi_{_P}} > 0$, where 
$\varphi_{_P}=\min_{\preceq_k} \{\varphi \mid  \varphi\in \{0,1\}^k \, , \; a_\varphi\neq 0\}$.  

We can finally introduce an order relation $\preceq$ on $\mathbb F$ by letting 
$1+S\prec 1+S'$ \, if $S'-S$ is a {\em positive} power series (with no constant term). It is 
immediate that $\preceq$ is a total order of $\mathbb F$ that is invariant under left and right multiplication. 
Therefore, $(\mathbb F,\preceq)$ is a bi-ordered group containing an isomorphic copy of $\mathbb F_2$.

\begin{small}
\begin{ejer}
Give an explicit description of the positive cone of the order built above and show directly that it is invariant by conjugacy.
\end{ejer}
\end{small}


\begin{small}
\begin{rem}\label{standard bi-orders} The above technique for embedding $\mathbb F_2$ into $\mathbb F$ --called the {\bf {\em Magnus expansion}}-- actually shows
that $\mathbb{F}_2$ is residually torsion-free nilpotent. Indeed, if $\Gamma=\mathbb F_2$ and $\Gamma_i$ denotes $i^{\mathrm{th}}$-term of its lower central series, then  
it is not hard to check that, for every $i \geq 0$, the group $\Phi (\Gamma_i)$ is contained in $1 + o(i+1)$ but not in $1+o(i+2)$.  This implies that the successive quotients 
$\Gamma_i/\Gamma_{i+1}$ are Abelian groups without torsion and that $\bigcap_i \Gamma_i=\{id\}$. Hence, $\mathbb F_2$ is residually torsion-free nilpotent. 
We refer to \cite{MKS} for more details on this. 

Observe that we can use the filtration $\Gamma_i$ to produce many bi-invariant orders on $\mathbb F_2$. Namely, for each $i\geq 0$, take a bi-order 
$\preceq_i$ on the (Abelian) quotient $\Gamma_i/\Gamma_{i+1}$, and then declare an element $g\in \mathbb F_2$ to be {\em positive} if $g$ belongs 
to $\Gamma_i\setminus \Gamma_{i+1}$ and $\Gamma_{i+1} \prec_i g\Gamma_{i+1}$.  Clearly, the set $P$ of {\em positive} elements defines a positive 
cone for a left-order of $\mathbb F_2$. The fact that $P$ is also a normal semigroup follows from the fact that the subgroups in the lower central series 
are normal subgroups of $\mathbb F_2$ and, hence, $g$ and $fg f^{-1}$ define the same element in $\Gamma_i/\Gamma_{i+1}$.
\end{rem}
\end{small}

\vspace{0.25cm}

We next explain a different, more dynamical, approach to produce bi-orders on $\mathbb{F}_2$.
This is done by building actions by homeomorphisms on the real line. To do this, we must consider  {\bf{\em  ping-pong actions}}. 
Given a set $X$, we say that two bijections $f$ and $g$ of $X$ have a ping-pong configuration if there are disjoint sets $A^+,A^-,B^+$ 
and $B^-$ of $X$ such that
$$f (X\setminus A^-)\subset A^+\,, \quad f^{-1}(X\setminus A^+)\subset A^-,$$
$$g (X\setminus B^-)\subset B^+\, , \quad g^{-1}(X\setminus B^+)\subset B^-.$$
We then say that the action of the group generated by $f$ and $g$ is a ping-pong action. The relevance of this 
notion comes from the clever observation, due to Klein, contained in the next exercise.

\begin{small}
\begin{ejer}
\label{ejer-ping-pong}
Show that if $f$ and $g$ are bijections of a set $X$ having a ping-pong configuration, then the group $\langle f, g\rangle$ 
generated by them is isomorphic to $\mathbb F_2$. Show further that any point $x_0$ not contained in $A^+\cup A^-\cup B^+\cup B^-$ 
has a free orbit under $\langle f, g \rangle$. (See \cite{harpe} in case of problems with this.)
\end{ejer}
\end{small}

The archetypical example of a ping-pong action is the action by circle homeomorphisms given 
by (powers of)  two topologically hyperbolic elements having disjoint sets of fixed points; see Figure 2. Recall that 
$ f \in \mathrm{Homeo}_+ (\mathbb{S}^1)$ is said to be {\bf {\em topologically hyperbolic}} if it has exactly two fixed 
points $r_f$ and $a_f$, the first of which is topologically repelling and the other topologically attracting. 

\begin{small}
\begin{ejer}
Let $f$ and $g$ be topologically hyperbolic homeomorphisms of $\mathbb S^1$ with disjoint sets of fixed points. 
Let $A^+$ (resp. $B^+$) be neighborhoods of $a_f$ (resp. $a_g$), and let $A^-$ (resp. $B^-$) be 
neighborhoods of $r_f$ (resp. $r_g$), all of them small enough   
so that $A^+,A^-,B^+$ and $B^-$ are two-by-two disjoint. Show that there is 
$N \in \mathbb N$ such that $A^+,A^-,B^+$ and $B^-$ yield a ping pong configuration for $f^n,g^n$ for each $n \geq N$. 
\end{ejer}
\end{small}

\vsp\vsp


\beginpicture

\setcoordinatesystem units <0.78cm,0.78cm>

\circulararc 360 degrees from 3 0
center at 0 0

\circulararc -101 degrees from 1.85 -2.1
center at 3.34 0.7
\plot 0.52 2.6 0.87 2.7 /
\plot 0.9 2.35 0.87 2.7 /

\put{$a_f$} at 1 3.2
\put{$\bullet$} at 1 2.8

\put{$r_f$} at 2.4 -2.6
\put{$\bullet$} at 2 -2.24

\put{$f$} at 1.1 -0.2

\circulararc -110 degrees from -2.8 0.2
center at -2.8 -1.9
\plot -2.5 -0.1 -2.8 0.2 /
\plot -2.5 0.4 -2.8 0.2 /

\put{$\bullet$} at -3 0.15
\put{$\bullet$} at -0.9 -2.9
\put{$g$} at -1.4 -1.4

\put{$r_g$} at -0.9 -3.3
\put{$a_g$} at -3.45 0.1

\put{} at -9.2 3.5
\put{Figure 2: A ping-pong action on the circle.} at 0 -4

\endpicture


\vspace{0.5cm}

To build a left-invariant order on $\mathbb F_2$, we consider two orientation-preserving, topologically hyperbolic circle 
homeomorphisms $f,g$ such that small neighborhoods of its fixed points yield a ping-pong configuration for 
$\langle f, g\rangle$. We can assume that $f,g$ are piecewise-affine homeomorphisms of $\mathbb S^1$. 
Fix  a point $x_0\in \mathbb S^1$ having free orbit under $\langle f,g\rangle$.

Now let $\tilde f$ and $\tilde g$ be lifts of $f$ and $ g$ to the real line, respectively. This implies that the group $ \langle \tilde f, \tilde g \rangle$ 
commutes with the unit translation  $T_1: x\mapsto x+1$ and is isomorphic to $\mathbb F_2$. Moreover, we can chose the lifts so that they 
have at least one fixed point on the real line (hence infinitely many fixed points). If $\tilde x_0\in \R$ denotes a lift of $x_0$, its orbit under 
$\langle \tilde{f} , \tilde{g} \rangle\simeq \mathbb F_2$ is free. Using this point and the action on the line, we can induce a left-invariant 
order on $\mathbb F_2$ via the dynamical-lexicographic procedure.

\vspace{0.1cm}

We remark that in the previous construction, for each $n\in \Z$, the following property holds: 

\begin{enumerate}

\item[(*)] { Each interval of the form $[n,n+1]$ contains a fixed point of $\tilde f$ and a fixed point of $\tilde g$,  
as well as a point $\tilde x_n$ having free orbit under $\langle \tilde f, \tilde g \rangle$. }

\end{enumerate}

To produce a bi-invariant order on $\mathbb F_2$, we need to modify the previous construction. 
Let $u \in [0,1]$  be a fixed point of $\tilde  f$. We change the homeomorphism $\tilde  f$ so that it fixes every point $x\leq u$,  
while we keep $\tilde  f$ intact on $[u,\infty)$. Analogously, for a fixed point $v\in [0,1]$ of $\tilde  g$, we change $\tilde  g$ so 
that it fixes every point $x\leq v$, while we keep $\tilde  g$ intact on $[v,\infty)$. We denote the resulting homeomorphisms 
by $a$ and $b$. Observe that these are piecewise-linear homeomorphisms. 

\vsp\vsp
\noindent \underbar {Step I.} {\em The homeomorphisms $a$ and $b$ generate a group isomorphic to $\mathbb F_2$.} 

\vsp\vsp

Indeed, let $w$ be a nontrivial reduced word in $\mathbb F_2$, and let $w(\tilde f, \tilde g)$ and $w(a,b)$ be its corresponding evaluations 
in the actions given by $\langle \tilde f, \tilde g \rangle$ and $\langle a, b \rangle$, respectively. It follows from the construction that 
$w(a,b)(x) = w(\tilde f, \tilde g)(x)$ for all large enough $x\in \R$. In particular, referring to (*) above, if we let $x=\tilde{x}_n$ for $n$ 
sufficiently large, we obtain that $w(a,b)(\tilde{x}_n)\neq \tilde{x}_n$. (In fact, as the reader can easily check, it suffices to take $n$ 
equal to 1 plus the word-length of $w$.) Therefore, $w(a,b)$ acts nontrivially, hence 
$w(a,b)$ is not the identity. This shows that $\langle a , b \rangle \simeq \mathbb F_2$.

\vsp\vsp

\noindent \underbar {Step II.} {\em The action of $\langle a,b \rangle $ allows inducing a bi-order on $\mathbb F_2$.} 

\vsp\vsp

Instead of directly defining the bi-order on $\mathbb F_2$, it is easier to define its positive cone. We first observe that, from 
the construction of $ \langle a,b  \rangle$, the set of break points of each $c \in \langle a,b  \rangle$ is bounded from 
below. Thus, every $c \in  \langle a,b \rangle$ has a least break point, which we denote by $x_c$. Remark that for a nontrivial 
$c$, this implies that  $c (x)=x$ for all $x\leq x_c$ and, since $c$ is piecewise affine,  there is a  small  right-neighborhood 
$V_c$ of $x_c$ such that either $c (x)>x$ or $c (x)<x$ for all $x\in V_c$. Equivalently, either $D_+ c (x_c) > 1$ or 
$D_+ c (x_c) < 1$, where $D_+ (\cdot)$ stands for the derivative on the right. Having this, we define 
$$P := \big\{ c \in \langle a,b \rangle \mid c (x) > x \text{ for all } x\in V_c \big\}.$$ 

In other words, since each nontrivial $c \in \langle a,b \rangle$ is piecewise affine, it has a {\em first} break point, on the right of which 
$c$ has a definitive {\em sign}. Clearly, this sign is invariant under conjugation by any homeomorphism of the real line, so $P$ is a normal 
subset of $\mathbb F_2$. Moreover, if $c$ is nontrivial, then exactly one of $c$ or $c^{-1}$ belongs to $P$, so $\mathbb F_2 = P\cup P^{-1}\cup \{id\}$ 
and  $P\cap P^{-1}=\emptyset$. Finally, if $c_1$ and $c_2$ are elements of $P$, then is easy to check that $c_1c_2$ is also an element of $P$, so $P$ is 
a subssemigroup of $\mathbb{F}_2$. This implies that $P$ is the positive cone of a bi-invariant order of $\mathbb F_2$.

\begin{small}\begin{rem}By performing the construction
above appropriately, one may obtain a bi--order on a free group $\mathbb{F}_2
= \langle a,b \rangle$ for which both $a$ and $b$ are positive but the product \esp
$a \esp [a,b]$ \esp is negative. Note that this cannot happen for the bi-orderings coming 
from the Magnus expansion. In fact, this cannot happen for any bi--order on $\mathbb{F}_2$ 
obtained via the lower central series as in Remark \ref{standard bi-orders}.
\end{rem}\end{small}

\index{Group!fully residually free}
\index{Group!surface group}
\noindent{\bf Surface groups.} Surface groups are residually free, hence
bi-orderable (see the end of \S \ref{general-2}). Actually, as we show below,
these groups are fully residually free, which is a
stronger property (see Remark \ref{no-es-lo-mismo} below). Recall that, if P is
some group property, then a group $\Gamma$ is said to be {\bf{\em fully residually}}
P if for every finite subset $\mathcal{G} \subset \Gamma \setminus \{id\}$, there exists
a surjective group homomorphism from $\Gamma$ into a group $\Gamma_{\mathcal{G}}$
satisfying P such that the image of every $g \in \mathcal{G}$ is nontrivial. Equivalently, for
every finite subset $\mathcal{G} \subset \Gamma$, there is an homomorphism $\Phi$ into 
a group satisfying P whose restriction to $\mathcal{G}$ is injective.

\begin{small}\begin{rem}\label{no-es-lo-mismo}
Obviously, the direct product $\mathbb{F}_2 \times \mathbb{F}_2$ is {\em residually free}, as any 
{\em single} nontrivial element is detected by projections. However, $\mathbb{F}_2 \times \mathbb{F}_2$ 
is not fully residually free, because given any distinct $f,g,h$ in $\mathbb{F}_2$, no homomorphism
from $\mathbb{F}_2 \times \mathbb{F}_2$ into a free group maps the elements \esp $(id,id)$, $(f,id)$,
$(g,id)$, $([f,g],id)$, and $(id,h)$, to five different ones. Indeed, as $(id,h)$ commutes
with $(f,id)$ and $(g,id)$, a separating homomorphism must send these three
elements into a cyclic subgroup. However, if this is the case, then
$([f,g],id)$ is mapped to the identity.
\end{rem}\end{small}

Below we deal with the case of surface groups for even genus (the case of odd genus
easily follows from this). The following lemma, due to Baumslag, will be crucial for us. 
The geometric proof that we give appears in \cite{gelander}. For the argument recall that, 
given a group $\Gamma$ generated by finitely many elements $g_1,\ldots,g_k$, its 
{\bf \em Cayley graph} is the graph whose vertices are the group elements, 
two of which are joined by an edge if they differ by left multiplication by some 
generator $g_i$ or its inverse. This graph has a natural metric space structure 
(edges are assumed to have length $1$). Moreover, the natural action of 
$\Gamma$ on itself induces an action by isometries of its Cayley graph.

\vspace{0.05cm}

\begin{lem} \label{clave}
{\em Let $g_1,\ldots,g_k$ be elements in a free group $\mathbb{F}_n$, and let $f$ 
be another element that does not commute with any of them. Then there exists
$N \!\in\! \mathbb{N}$ such that, for every $|n_i| \geq N$, $m \in \mathbb{N}$,
and $j_i \in \{1,\ldots,k\}$,} 
$$g_{j_1} f^{n_1} g_{j_2} f^{n_2} \ldots g_{j_{m}} f^{n_m} \neq id.$$
\end{lem}

\noindent{\bf Proof.} The Cayley graph of $\mathbb{F}_n$ with respect to the canonical system of generators 
naturally identifies with an homogeneous tree $\mathcal{T}_{2n}$ with valence $2n$ at each vertex. This tree 
has a natural {\em boundary at infinity}, that we denote by $\partial \mathcal{T}_{2n}$. One easily shows that 
every nontrivial $f\in \mathbb{F}_n$ acts on this tree as a translation along an axis  (that we denote by 
$\mathrm{axis}(f)$), which has (different) endpoints $a^- = a^-(f)$ and $a^+(f)=a^+$ in $\partial \mathcal{T}_{2n}$. 
If $g_i$ does not commute with $f$, one may show that 
$\{g_i (a^-), g_i(a^+)\} \cap \{a^-,a^+\} = \emptyset$, for every
$i \in \{1,\ldots,k\}$. Let $U^-, U^+$ be neighborhoods in
$\partial \mathbb{F}_n$ of $a^-$ and $a^+$, respectively, satisfying
$$g_i (U^- \cup U^+) \cap (U^- \cup U^{+}) =
\emptyset \esp\esp\esp \mbox{ for each } \esp i.$$
There exists $N \in \mathbb{N}$ such that, for all $r \geq N$,
$$f^r (\partial \mathbb{F}_n \setminus U^{-}) \subset U^{+}, \qquad
f^{-r} (\partial \mathbb{F}_n \setminus U^+) \subset U^-.$$
A ping-pong type argument (see Exercise \ref{ejer-ping-pong}) 
then shows the lemma. $\hfill\square$

\vspace{0.32cm}

Let $\Gamma = \Gamma_{2n}$ be the $\pi_1$ of an orientable surface $S_{2n}$
of genus $2n$ ($n \geq 1$). Let us consider the standard presentation
$$\Gamma = \big\langle g_i, g_i', h_i, h_i', 1 \leq i \leq n :
[g_1, g_1'] \cdots [g_n, g_n' ] \cdot [h_n',h_n] \cdots [h_1',h_1]= id \big\rangle.$$
Following \cite{breuillard et all}, let $\sigma$ be the automorphism of $\Gamma$
that leaves the $g_i$'s and $g_i'$'s fixed but sends $h_i$ to $f h_i f^{-1}$ and 
$h_i'$ to $f h_i' f^{-1}$ for each $i$, where
$f := [g_1, g_1'] \cdots [g_n, g_n']$. (Geometrically, this corresponds to the
Dehn twist along the closed curve obtained from a simple curve that joins
the first and the $2n^{th}$ vertices of the hyperbolic $4n$-gon that yields
$S_{2n}$.) Finally, let $\varphi$ be the surjective homomorphism
from $\Gamma$ to the free group $\mathbb{F}_{2n}$ with free generators
$a_1,\ldots, a_n,a_1',\ldots,a_n'$ defined by $\varphi (g_i ) = \varphi (h_i ) = a_i$
and $\varphi (g_i' ) = \varphi (h_i' ) = a_i'$. We claim that the sequence of
homomorphisms $\varphi \circ \sigma^k$ is {\bf{\em eventually faithful}}, 
in the sense that given any nontrivial elements 
$f_1,\ldots,f_m$ in $\Gamma$, there exists $N \in \mathbb{N}$
so that for all $k \geq N$, the image under $\varphi \circ \sigma^k$ 
of each $f_i$ is nontrivial (hence $\Gamma$ is fully residually free).

To show the claim above, given $g \in \Gamma \setminus \{id\}$, let us write it in the form
$$g = w_1 (g_i, g_i' ) \cdot w_2 (h_i, h_i' )
\cdots w_{2p-1} (g_i, g_i' ) \cdot w_{2p} (h_i, h_i' ),$$
where each $w_j (g_i, g_i')$ and $w_j (h_i, h_i')$ 
are reduced words in $2n$ letters (the first and/or the last $w_j$ may be trivial).
Up to modifying the $w_{2j-1}$'s, we may assume that each $w_{2j}$ (where $1 \leq j \leq p$) 
is such that $w_{2j} (h_i , h_i' )$ is not a power of $f$. Note that the centralizer of $f$ in
$\Gamma$ is the cyclic group generated by $f$. By regrouping several $w_j$'s into a longer word
if necessary, unless $g$ itself is a power of $f$, we may also assume that $w_{2j-1} (g_i , g_i' )$
is not a power of $f$. Let $\bar{f}$ be the image of $f$ under $\varphi$. We have
$$\varphi \circ \sigma^k (g) =
w_1 \bar{f}^k w_2 \bar{f}^{-k} \cdots w_{2p-1} \bar{f}^k w_{2p} \bar{f}^{-k},$$
where $w_j = w_j (a_i , a_i' )$. Since $\bar{f}$ does not commute with any of these $w_j$'s,
Lemma \ref{clave} implies that $\varphi \circ \sigma^k (g)$ is nontrivial.

\begin{rem} 
It is worth pointing out that, in contrast to nilpotent groups (see the discussion on Rhemtulla's theorem discussed at the end of 
\S \ref{ejemplificando-1}), there are partial left-orders on the free group that cannot be extended to total left-orders. Indeed, 
this holds for the partial left-order whose positive elements are those lying in the semmigroup generated by $f^2$, $g^2$ and 
$f^{-1} g^{-1}$ in $\mathbb{F}_2 = \langle f, g \rangle$ (this example is taken from \cite{colacito-metcalfe}). The same holds 
for the semigroup generated by $fg, f^{-1}g^{-1},fg^{-1}$ and $f^{-1}g$ (this last example was kindly communicated to us by 
Metcalfe; its interest comes from that the generators of the semigroup are not powers of other elements).
\end{rem}


\subsection{Thompson's group $\mathrm{F}$}
\label{ejemplificando-4}
\index{Thompson's group}

\hspace{0.45cm} Thompson's group F is perhaps the simplest example of a bi-orderable group
that is not residually nilpotent. For the definition, recall that a {\bf {\em dyadic number}} is a rational number 
of the form \( p / 2^q\), where \(p\) and \(q\) are integers.  We will say that an orientation-preserving 
homeomorphism between intervals of the real line is a {\bf {\em piecewise-dyadic homeomorphism}} 
if it is piecewise-affine with dyadic numbers as break points and derivative equal to some integer power of $2$  
at each regular point. Thompson's group \(F\) is by definition the group of piecewise-dyadic homeomorphisms of the interval \([0,1]\).

 This group is far from being residually nilpotent because its commutator subgroup 
$\mathrm{F}' = [\mathrm{F},\mathrm{F}]$ is simple (see Theorem \ref{Teo F' is simple} further 
on). To see that it is bi-orderable, for 
each nontrivial $f \!\in\! \mathrm{F}$ we denote by $x^{-}_{f}$ (resp. $x^{+}_f$) the leftmost
point $x^-$ (resp. the rightmost point $x^+$) for which $D_+f (x^{-}) \neq 1$ (resp.
$D_{-} f (x^{+}) \neq 1$), where, as before, $D_+f$ and $D_{-}f$ stand for the corresponding lateral
derivatives. One can immediately visualize four different bi-orders on (each subgroup of)
$\mathrm{F}$, namely the bi-order $\preceq_{x^{-}}^{+}$ (resp. $\preceq_{x^{-}}^{-}$,
$\preceq_{x^{+}}^{+}$, $\preceq_{x^{+}}^{-}$) for which \esp $f$ is positive \esp
if and only if \esp $D_+f (x^{-}_f) > 1$ \esp (resp. \esp $D_+ f (x^{-}_f) < 1$,
\esp $D_{-} f (x^{+}_f) < 1$, \esp $D_{-} f (x^{+}_f) > 1$). (Compare the construction in 
Step II of \S \ref{ejemplificando-3}.) Although $\mathrm{F}$
admits many more bi-orders than these (see theorem \ref{thm-F} below), the case
of $\mathrm{F}'$ is quite different. The result below is
essentially due to Dlab \cite{dlab} (see also \cite{F}).

\vspace{0.05cm}

\begin{thm} {\em The only bi-orders on $\mathrm{F}'$
are $\preceq_{x^{-}}^{+}$, $\preceq_{x^{-}}^{-}$, $\preceq_{x^{+}}^{+}$ and
$\preceq_{x^{+}}^{-}$.}
\end{thm}

\vspace{0.05cm}

Remark that there are four other ``exotic'' bi-orders on F, namely:\\

\vspace{0.1cm}

\noindent -- The bi-order $\preceq_{0,x^{-}}^{+,-}$ for which $f$ is positive if and only
if either $x^{-}_f = 0$ and $D_+ f (0) > 1$, or $x^{-}_f \neq 0$ and $D_+ f (x^{-}_f) < 1$;\\

\vspace{0.1cm}

\noindent -- The bi-order $\preceq_{0,x^{-}}^{-,+}$ for which $f$ is positive if and only
if either $x^{-}_f = 0$ and $D_+ f (0) < 1$, or $x^{-}_f \neq 0$ and $D_+ f (x^{-}_f) > 1$;\\

\vspace{0.1cm}

\noindent -- The bi-order $\preceq_{1,x^{+}}^{+,-}$ for which $f$ is positive if and only
if either $x^{+}_f = 1$ and $D_+ f (1) < 1$, or $x^{+}_f \neq 1$ and $D_{-} f (x^{+}_f) > 1$;\\

\vspace{0.1cm}

\noindent -- The bi-order $\preceq_{1,x^{+}}^{-,+}$ for which $f$ is positive if and only
if either $x^{+}_f = 1$ and $D_+ f (1) > 1$, or $x^{+}_f \neq 1$ and $D_{-} f  (x^{+}_f) < 1$.\\

\vspace{0.1cm}

\noindent Note that, when restricted to $\mathrm{F}'$, the bi-order $\preceq_{0,x^{-}}^{+,-}$
(resp. $\preceq_{0,x^{-}}^{-,+}$, $\preceq_{1,x^{+}}^{+,-}$, and $\preceq_{1,x^{+}}^{-,+}$)
coincides with $\preceq_{x^{-}}^{-}$ (resp. $\preceq_{x^{-}}^{+}$, $\preceq_{x^{+}}^{-}$,
and $\preceq_{x^{+}}^{+}$). Let us denote the set of the previous eight bi-orders on F
by $\mathcal{BO}_{Isol} (\mathrm{F})$.

\vspace{0.15cm}

There is another natural procedure to create bi-orders on F. For this, recall the well-known fact that 
$\mathrm{F}'$ coincides with the subgroup of $\mathrm{F}$ formed by the elements $f$ satisfying 
$D_+ f (0) = D_{-} f (1) = 1$ (see Exercise \ref{F:ejer-conmutator} for this).
Now let $\preceq_{\mathbb{Z}^2}$ be any bi-order on $\mathbb{Z}^2$, and let
$\preceq_{\mathrm{F}'}$ be any bi-order on $\mathrm{F}'$. It readily follows from
Dlab's theorem that $\preceq_{\mathrm{F}'}$ is invariant under conjugation by elements
in $\mathrm{F}$. Hence, one may define a bi-order $\preceq$ on F by declaring that
$f \succ id$ if and only if either $f \notin \mathrm{F}'$
and $\big( \log_2 (Df_{+}(0)), \log_2 (Df_{-}(1))\big)
\succ_{\mathbb{Z}^2} \big( 0,0 \big)$, or $f \in \mathrm{F}'$
and $f \succ_{\mathrm{F}'} id$ \esp (see \S \ref{section-convex-extension}
for more details on this type of construction).

All possible ways of left-ordering finite-rank, Abelian groups were described
in \S \ref{ejemplificando-1}. Since there are only four possibilities for
$\preceq_{\mathrm{F}'}$, the preceding procedure gives us four
sets (which we will coherently denote by $\Lambda_{x^{-}}^{+}$, $\Lambda_{x^{-}}^{-}$,
$\Lambda_{x^{+}}^{+}$, and $\Lambda_{x^{+}}^{-}$) naturally homeomorphic to the Cantor
set (in the sense of \S \ref{space-left-orders}) inside the set of bi-orders of F. The
main result of \cite{F} establishes that these bi-orders, together with the eight special
bi-orders previously introduced, are all the possible bi-orders on F. The proof is a 
straightforward application of Conrad's theory to be extensively 
developed in \S \ref{classic-conrad}.

\vspace{0.1cm}

\begin{thm} \label{thm-F}
{\em The set of all bi-orders of \esp $\mathrm{F}$ consists of
the disjoint union of $\mathcal{BO}_{Isol} (\mathrm{F})$ and the sets
$\Lambda_{x^{-}}^{+}$, $\Lambda_{x^{-}}^{-}$, $\Lambda_{x^{+}}^{+}$, and
$\Lambda_{x^{+}}^{-}$.}
\end{thm}

\vspace{0.1cm}

Thompson's group F is remarkable in many aspects. Among its most relevant properties, we can mention that it 
is finitely presented (this is very well explained in \cite{CFP}; see also Exercise \ref{F:fin-pres} below for a sketch 
of proof), it contains no free subgroup (this was first proved in \cite{BS}; see \S \ref{ejemplificando-relatives of F}  
for a proof of a slightly generalized version of this fact) and its commutator subgroup $F'$ is a simple group (see 
Theorem \ref{Teo F' is simple}). At the time of publication of this book, the challenging 
question of the amenability of this group remains open.

\begin{small}\begin{ejer} \label{F:fin-pres}
The goal of this exercise is to provide the main steps to derive the next two presentations of 
Thompson's group F:
$$\mathrm{F}_1 = \big\langle a,b \!: [ a^{-1} b, aba^{-1} ] = [ a^{-1} b, a^2 b a^{-2} ] =  id \big\rangle,$$ 
$$\mathrm{F}_2 = \big\langle c_0, c_1, c_2, \ldots : c_k c_n c_k^{-1} = c_{n+1} \mbox{ for all } k < n \big\rangle.$$ 
Here, $a$ and $b$ are in correspondence with $c_0$ and $c_1$, respectively, and correspond to the 
elements $f_0, f_1$ in F whose graphs are drawn below.

\vspace{0.54cm}


\beginpicture

\setcoordinatesystem units <0.55cm,0.55cm>

\putrule from 0 0 to 6 0

\putrule from 0 0  to 0 6

\putrule from 6 0 to 6 6

\putrule from 0 6 to 6 6

\put{$f_0$} at 1.85 4 
\put{$f_1$} at 11.45 2.1

\put{Figure 3: The graphs of $f_0$ and $f_1$.} at 8 -0.8

\putrule from 10 0 to 16 0
\putrule from 10 6 to 16 6
\putrule from 10 0 to 10 6
\putrule from 16 0 to 16 6

\plot
0  0
1.5  3 /

\plot
1.5 3 
3 4.5  /

\plot
3 4.5 
6   6 /


\plot
10 0
13 3 /

\plot
13   3
13.75  4.5  /

\plot
13.75  4.5 
14.5  5.25  /

\plot
14.5  5.25
16  6 /

\setdots

\putrule from 3 0 to 3 6
\putrule from 1.5 0  to 1.5 6
\putrule from 0 4.5  to 6 4.5
\putrule from 0 3 to 6 3

\putrule from 13 0 to 13 6
\putrule from 14.5 0 to 14.5 6
\putrule from 13.75 0 to 13.75 6
\putrule from 10 3 to 16 3
\putrule from 10 5.25 to 16 5.25
\putrule from 10 4.5 to 16 4.5

\put{} at -5 0 

\endpicture


\vspace{0.54cm}

\noindent (i) Let $\Phi \!: \mathrm{F}_1 \to \mathrm{F}_2$ be the map sending $a$ to $c_0$ and $b$ to $c_1$. Show that $\Phi$ 
extends to a group homomorphism.

\noindent{\underbar{Hint.}} What is to be checked is that the following relations are satisfied in $\mathrm{F}_2$: 
$$[ c_0^{-1} c_1, c_0 c_1 c_0^{-1} ] = [ c_0^{-1} c_1, c_0^2 c_1 c_0^{-2} ] =  id.$$
To do this, just note that from 
$$c_0 c_2 c_0^{-1} = c_3  = c_1 c_2 c_1^{-1} \qquad \mbox{(resp. }  c_0 c_3 c_0^{-1} = c_4 = c_1 c_3 c_1^{-1}),$$
one gets that $c_0^{-1} c_1$ conmutes with $c_2$ (resp. $c_3$).

\vsp

\noindent (ii) Let $\hat\Phi \!: \mathrm{F}_2 \to \mathrm{F}_1$ be the map sending $c_0$ to $a$ and $c_1$ to $b$. Show 
that $\hat\Phi$ extends to a group homomorphism (hence, by (i), to a group {\em isomorphism}, with $\hat\Phi = \Phi^{-1}$).

\noindent{\underbar{Hint.}} Set $a_0 := a$ and $a_n := a^{n-1} b a^{-(n-1)}$ for $n \geq 1$. The task is to show that 
\begin{equation}
a_k a_n a_k^{-1} = a_{n+1} \quad \mbox{for all } k < n.
\label{F:primera}
\end{equation}
One can show this by simultaneously proving that
\begin{equation}
[b a^{-1}, a_j] = id \quad \mbox{for all } j \geq 3.
\label{F:segunda}
\end{equation}
First, note that the second condition above holds for $j = 3$ and $j=4$, since
$$[a^{-1} b, a b a^{-1}] \!=\! id \implies [ba^{-1}, a^2 b a^{-2}] \!=\! id \implies [ba^{-1},a_3] \!=\! id$$
and 
$$[a^{-1} b, a^2 b a^{-2}] \!=\! id \implies [ba^{-1}, a^3 b a^{-3}] \!=\! id \implies [ba^{-1},a_4] \!=\! id.$$
Assume that (\ref{F:primera}) holds for $k \leq n \leq k + i - 3$ and that (\ref{F:segunda}) holds for $3 \leq j \leq i$. 
Then $ba^{-1}$ commutes with both $a_3$ and $a_i$, hence with $a_{i+1} = a_3 a_i a_3^{-1}$, and therefore 
(\ref{F:segunda}) holds for $j = i+1$. Moreover,  
\begin{eqnarray*}
a_{n+1} a_k 
&=& a^n b a^{-n} a^{k-1} b a^{-(k-1)} 
\esp = \esp a^{k-1} a^{n-k+1} b a^{-(n-k+1)} b a^{-1} a^{-k+2} \\
&=& a^{k-1} a_{n-k+2} (ba^{-1}) a^{-k+2}
\esp = \esp a^{k-1} (ba^{-1}) a_{n-k+2} a^{-k+2} \\
&=& (a^{k-1} b a^{-(k-1)} ) (a^{k-2} a_{n-k+2} a^{-(k-2)})
\esp = \esp a_k a_n,
\end{eqnarray*}
and therefore (\ref{F:primera}) holds for $n = k+i-2$. The proof can be completed via 
an induction argument.

\vsp

\noindent (iii) Every element $f \in \mathrm{F}$ can be identified in an obvious way 
with a map that sends in an ordered way the terminal points (called {\em leaves}) of a dyadic rooted tree 
to those of another dyadic rooted tree having the same number of leaves, and conversely. In this view, elements $f_0$ 
and $f_1$ correspond to the following diagrams:

\vspace{0.54cm}


\beginpicture

\setcoordinatesystem units <0.55cm,0.55cm>

\plot 
0  0
-1  -1 /

\plot 
0 0 
1 -1  /

\plot 
-1 -1 
-1.5 -2 /

\plot 
-1 -1 
-0.5 -2 /

\put{$f_0$} at 2.5 -0.5 
\put{$\longrightarrow$} at 2.5 -1 

\plot 
5  0
4  -1 /

\plot 
5 0 
6 -1  /

\plot 
6 -1 
5.5 -2 /

\plot 
6 -1 
6.5 -2 /


\plot 
13  0
12  -1 /

\plot 
13 0 
14 -1  /

\plot 
14 -1 
13.5 -2 /

\plot 
14 -1 
14.5 -2 /

\plot 
13.5 -2 
13.2 -3 /

\plot 
13.5 -2 
13.8 -3 /

\put{$f_1$} at 16.2 -0.5 
\put{$\longrightarrow$} at 16.2 -1 

\plot 
19  0
18  -1 /

\plot 
19 0 
20 -1  /

\plot 
20 -1 
19.5 -2 /

\plot 
20 -1 
20.5 -2 /

\plot 
20.5 -2 
20.2 -3 /

\plot 
20.5 -2 
20.8 -3 /

\put{Figure 4: The diagrams of the elements $f_0$ and $f_1$.} at 10 -3.8
\put{} at -3.3 0 

\endpicture


\vspace{0.54cm}

Given $n \geq 2$, let $f_n := f_0^{n-1} f_{1} f_{0}^{-(n-1)}$. Check that the tree diagram associated to $f_n$ is 
the following (the right-side tree below will be denoted by $\mathcal{T}_n$):

\vspace{0.4cm}


\beginpicture

\setcoordinatesystem units <0.55cm,0.55cm>

\plot 
0  0
-1  -1 /

\plot 
0 0 
1 -1 /

\plot 
1 -1 
0 -2 /

\plot 
1 -1 
2 -2 /

\plot
2 -2 
1 -3 /

\plot 
2 -2 
3 -3 /

\plot
3 -3 
2 -4 /

\plot 
3 -3 
4 -4 /

\plot 
2 -4 
1.5 -5 /

\plot 
2 -4 
2.5 -5 /

\plot 
5 -2
7.5 -2 /

\plot 
7.5 -2 
7.2 -2.2 /

\plot 
7.5 -2 
7.2 -1.8 /

\put{$f_n$} at 6.2 -1.5 

\plot 
10  0
9  -1 /

\plot 
10 0 
11 -1 /

\plot 
11 -1 
10 -2 /

\plot 
11 -1 
12 -2 /

\plot
12 -2 
11 -3 /

\plot 
12 -2 
13 -3 /

\plot
13 -3 
12 -4 /

\plot 
13 -3 
14 -4 /

\plot 
14 -4 
13.5 -5 /

\plot 
14 -4 
14.5 -5 /

\put{$n$} at 2 -1 

\plot 
0.6 0.4 
1.7 -0.7 /

\plot 
2.4 -1.4 
3.5 -2.5 /

\put{$n$} at 12 -1 

\plot 
10.6 0.4 
11.7 -0.7 /

\plot 
12.4 -1.4 
13.5 -2.5 /

\put{Figure 5: The diagram of the element $f_n := f_0^{n-1} f_1 f_0^{-(n-1)}$.} at 7 -5.9
\put{} at -6 0 

\endpicture


\vspace{0.5cm}

\noindent Note that, although 
the diagram representing an element $f \!\in\! \mathrm{F}$ is not unique, there is a unique {\em reduced} one, in the 
sense that any other representative diagram can be obtained from this just by adding {\em carets} ($\bigwedge$) at the 
leaves of the source tree and the image leaves of the target tree. (We keep right-to-left notation for group multiplication, 
which is opposite to most of the literature on the subject, including \cite{CFP}.)

\vsp

\noindent (iv) Show that $f_0, f_1, f_2, f_3, \ldots$ (hence $f_0,f_1$) generate F. 
To do this, show that every element $f \in \mathrm{F}$ may be written in the form 
$$f_0^{-r_0} f_1^{-r_1} \cdots f_n^{-r_n} f_n^{s_n} \cdots f_1^{s_1} f_0^{s_0},$$
where $r_i \geq 0$ and $s_i \geq 0$. Besides, for a nontrivial element, such a writing is unique when respecting the 
next two properties:  exactly one of $r_n,s_n$ is zero, and if $r_k,s_k$ are both positive for a certain 
$k < n$, then at least one of $r_{k+1},s_{k+1}$ is positive.
 
\noindent{\underbar{Hint.}} Let $f \in \mathrm{F}$ be an element represented by a tree diagram in which the target tree 
is $\mathcal{T}_n$. For each leaf $v_i$ of the source tree (which are numbered starting from 0), let $s_i$ be the length 
of the largest path along (unit) left branches that starts at $v_i$ and does not touch the right side of the tree. Show 
that $f = f_{n+2}^{s_{n+2}} f_{n+1}^{s_{n+1}} \cdots f_{1}^{s_{1}} f_0^{s_0}$. (Note that $s_{n+1} = s_{n+2} = 0$.) 
See the figure below for an example.

\vspace{0.4cm}


\beginpicture

\setcoordinatesystem units <0.6cm,0.6cm>

\plot 
0  0
-1  -1 /

\plot 
-1 -1 
-1.4 -2 /

\plot 
-1 -1 
-0.6 -2 /

\plot 
0 0 
1 -1 /

\plot 
1 -1 
0.6 -2 /

\plot 
1 -1 
1.4 -2 /

\plot
0.6 -2  
0.3 -3 /

\plot 
0.6 -2 
0.9 -3 /

\plot 
0.3 -3 
0.1 -4 /

\plot 
0.3 -3 
0.5 -4 /

\plot 
3.7 -2.3
6.8 -2.3 /

\plot 
6.8 -2.3
6.5 -2.1 /

\plot 
6.8 -2.3 
6.5 -2.5 /

\put{$f_2^2 f_0$} at 5.2 -1.5 

\plot 
10  0
9  -1 /

\plot 
10 0 
11 -1 /

\plot 
11 -1 
10 -2 /

\plot 
11 -1 
12 -2 /

\plot
12 -2 
11 -3 /

\plot 
12 -2 
13 -3 /

\plot
13 -3 
12 -4 /

\plot 
13 -3 
14 -4 /

\plot 
14 -4 
13.5 -5 /

\plot 
14 -4 
14.5 -5 /

\put{$v_0$} at -1.4 -2.5
\put{$v_1$} at -0.6 -2.5 
\put{$v_2$} at 0 -4.5
\put{$v_3$} at 0.6 -4.5 
\put{$v_4$} at 1  -3.5
\put{$v_5$} at 1.5 -2.5

\put{$\downarrow$} at -1.4 -5.43
\put{$\downarrow$} at -0.6 -5.43 
\put{$\downarrow$} at 0 -5.43
\put{$\downarrow$} at 0.6 -5.43 
\put{$\downarrow$} at 1 -5.43
\put{$\downarrow$} at 1.5 -5.43

\put{$1$} at -1.4 -6.5
\put{$0$} at -0.6 -6.5 
\put{$2$} at 0 -6.5
\put{$0$} at 0.6 -6.5 
\put{$0$} at 1 -6.5
\put{$0$} at 1.5 -6.5

\put{Figure 6: The tree diagram of the map $f_2^2 f_0$, as predicted by the claim above.} at 6.6 -7.89
\put{} at -6.5 0 

\endpicture


\vspace{0.4cm}

\noindent (v) Show that F is isomorphic to both $\mathrm{F}_1$ and $\mathrm{F}_2$. 

\noindent{\underbar{Hint.}} Show that the map $a \mapsto f_0$ and $b \mapsto f_1$ extends to a 
surjective group homomorphism $\Phi_1 \!: \mathrm{F}_1 \to \mathrm{F}$, hence to a surjective group 
homomorphism $\Phi_2 \!: \mathrm{F}_2 \to \mathrm{F}$. To show that the latter is injective, use the relations 
$$c_n c_k^{-1} = c_k^{-1} c_{n+1}, \quad
c_k c_n^{-1} = c_{n+1}^{-1} c_k, \quad 
c_k c_n = c_{n+1} c_k, \quad \mbox{for }   n > k,$$
to transform an arbitrary expression in the $c_i$'s into one where only negative exponents appear  
on the left, only positive exponents appear on the right, and the subindices are ordered, say  
$$c_0^{-r_0} c_1^{-r_1} \cdots c_n^{-r_n} c_n^{s_n} \cdots c_1^{s_1} c_0^{s_0}.$$
Besides, if both $r_k,s_k$ are positive and both $r_{k+1},s_{k+1}$ are zero for a certain $k$, then using the 
relation $c_k^{-1} c_{n+1} c_k = c_n$ for $n > k$, one may decrease the subindex of each entry between 
$c_k^{-r_k}$ and $c_k^{s_k}$. Proceeding this way as much as possible, we get either the empty word or 
an expression as in (iv) above, which was shown to correspond (under $\Phi_1$) to a nontrivial element of F.
\end{ejer}

\begin{ejer}\label{F:ejer-conmutator} Show that $\mathrm{F}'=[\F,\F]$ coincides with 
the set of elements of $\F$ whose {\bf {\em suppor}t} (that is, the closure of the 
set of points that are moved by some group element) is strictly contained in $]0,1[$. 
To do this, use the fact that $\F$ is generated by two elements, and that the base-2 logarithm 
of the derivatives at $0$ and $1$ provides a surjective homomorphism $\Phi$ from $F$ onto $\Z^2$.

\noindent{\underline{Hint.}} Show that if $w$ is a word in the generators $f_0,f_1$ from Exercise \ref{F:fin-pres} that represents an 
element for which $\Phi$ vanishes, then the total exponents of $f_0$ and $f_1$ vanish. Then use this fact to write $w$ as a product 
of commutators. As a matter of example, for $w = f_0^2 f_1 f_0^{-1} f_1^{-2} f_0^{-1} f_1,$ one has
$$w \!=\! [f_0^2,f_1] (f_1 f_0^2) f_0^{-1} f_1^{-2} f_0^{-1} f_1 = [f_0^2,f_1] f_1 f_0 f_1^{-2} f_0^{-1} f_1 
= [f_0^2,f_1] [f_1, f_0] (f_0 f_1) f_1^{-2} f_0^{-1} f_1,$$ 
hence
$w= [f_0^2,f_1] [f_1, f_0] [f_0, f_1^{-1}] .$
\end{ejer}

\end{small}

We now show that the commutator subgroup $\mathrm{F}'=[\mathrm{F}, \mathrm{F}]$ is a simple group. The standard reference for 
this is \cite{CFP}, where the normal forms in $\mathrm{F}$ obtained in item iv) of \ref{F:fin-pres} are crucial. Here, we use instead a 
more dynamical approach that goes along the lines of the arguments to be exploited in \S \ref{s: finitely generated simple left orderable}.  

For a closed {\bf {\em dyadic interval}} $I\subseteq [0,1]$ (that is, an interval whose endpoints are dyadic rationals), we denote by $\mathrm{F}_I$ 
the subgroup made up of the elements of $\mathrm{F}$ with support contained in $I$. A piecewise-dyadic homeomorphism from \([0,1]\) to \( I\) 
gives a conjugacy between  $\mathrm{F}$ and $\mathrm{F}_I$. The existence of such a homeomorphism follows from the next exercise.

\begin{small}\begin{ejer} \label{ej: transitivity on dyadic interval}
Prove that given two dyadic intervals $J, K$ in the line, there exists a piecewise-dyadic homeomorphism 
sending $J$ onto $K$. 

\noindent{\underline{Hint.}} It suffices to assume that \(J= [0,1]\). By applying a dyadic affine map of the form \( x \mapsto 2^p x + q \) 
for some integers \( p\in  \mathbb N \) and \( q \in \mathbb Z\), one can further assume that \( K \) is of the form \( [0,n ]\). Write $n$ 
in dyadic expansion, say \(n = \varepsilon_0 +2\varepsilon _1+\ldots + \varepsilon _ k 2^k  \), where  \( \varepsilon _i \in \{0,1\}\) 
and \(\varepsilon _k=1\). Now observe that  \([0,n]\) is the image of an interval of the form \(   [0, n']\), with  \(0< n' \leq k < n\),  
by a piecewise-dyadic homeomorphism.  
 \end{ejer}\end{small}

The following lemma is a special case of the famous {\bf {\em Higman's trick}}, 
which applies to general groups of transformations.

\begin{lem}
\label{lem Higman} 
{\em If $N\lhd F$ is a nontrivial normal subgroup, 
then there is a nonempty dyadic interval $I\subset [0,1]$ such that $(\mathrm{F}_I)'$ is contained in $ N$.}
\end{lem}

\noindent{\bf Proof.}
 Let $c \in N\setminus \{\id\}$, and let $I\subseteq [0,1]$ be a closed dyadic interval such that $c(I)\cap I=\emptyset$. We claim that 
 $(\mathrm{F}_I)'$ is contained in $ N$. To show this, first note that, for every $a \in \mathrm{F}_I$, 
$$[a,c](x) := a c a^{-1} c^{-1}(x)= \left\{\begin{array}{ll } a(x) & \text{for } x\in I, \\ c a^{-1} c^{-1}(x) & \text{for } x\in u(I), \\ x & \text{otherwise}. \end{array}\right. $$
This easily implies that, for all $b \in \mathrm{F}_I$ and all $x \in [0,1],$
$$[[a,c],b](x)=[a,b](x).$$
Therefore, $[[a,c],b] =[a,b]$. Since the subgroup $N$ is normal in $\mathrm{F}$ and $c \in N$, we  have that $[a,c]$ belongs to $N$, 
which in its turn implies that $[a,b]=[[a,c],b]$ also belongs to $N$. Since $a,b$ were arbitrary elements of  $\mathrm{F}_I$, 
we have that $(\mathrm{F}_I)'$ is contained in $N$, as claimed. $\hfill\square$

\begin{cor} \label{cor cocientes de F} {\em  If $N$ is a nontrivial normal subgroup of $\mathrm{F}$, then 
$\mathrm{F}'\subset N$. In other words, every proper quotient of $\mathrm{F}$ is Abelian.}

\end{cor}

\noindent{\bf Proof.}  
We know from Lemma \ref{lem Higman} that there is a closed dyadic interval $I$ such that $(\mathrm{F}_I)'\subseteq N$. 
Now, for an arbitrary $f \in \mathrm{F}$, we have that $f \, (\mathrm{F}_I)' f^{-1} = (\mathrm{F}_{f(I)})'$. Since $N$ is normal, this 
implies that $ (\mathrm{F}_{f(I)})'$ is also contained in $N$. 

Let now $(f_n)$ be a sequence of elements in $\mathrm{F}$ such that, for $I_n := f_n (I)$, one has  $\bigcup_n I_n=(0,1)$. 
By the previous discussion, $N$ contains all the groups $(\mathrm{F}_{I_n})'$. 
It follows from Exercises \ref{F:ejer-conmutator} and \ref{ej: transitivity on dyadic interval} 
that $N$ contains $\mathrm{F}'$ as well. $\hfill\square$ 

\vspace{0.2cm}

\begin{thm} \label{Teo F' is simple} 
{\em The commutator subgroup $\mathrm{F}'$ is a simple group.}
\end{thm}

\noindent{\bf Proof.}  Let $N$ be a nontrivial normal subgroup of $F'$. Fix any closed dyadic interval $I \subset (0,1)$. 
We first claim that $N$ contains $(\mathrm{F}_I)'$. Assuming this, one finishes the proof as in the previous corollary 
just taking care of selecting the conjugating elements $f_n$ in $\mathrm{F}'$.

Now, to see that $N$ contains $(\mathrm{F}_I)'$, choose a nontrivial element $f \in N$.  Since $f$ is in $\mathrm{F}'$, 
its support is strictly contained in $(0,1)$. Letting $h \in \mathrm{F}'$ be an element that sends the support of $f$ into 
$I$, we have that $h f h^{-1} \in N$ has support strictly contained in $I$. By Exercise \ref{ej: transitivity on dyadic interval},  
we have that $\mathrm{F}_I$ is isomorphic to $\mathrm{F}$, hence 
from Exercise \ref{F:ejer-conmutator} we conclude that $hfh^{-1}$ belongs to $(\mathrm{F}_I)'$. Thus, $N \cap (\mathrm{F}_I)'$ 
is a nontrivial subgroup of $\mathrm{F}_I \subset \mathrm{F}'$. By Corollary \ref{cor cocientes de F} again,  
$N \cap (\mathrm{F}_I)'$ contains $(\mathrm{F}_I)'$, as claimed. 
$\hfill\square$

\vspace{0.25cm}

Using the finite presentation from Exercise \ref{F:fin-pres} and the fact that every proper quotient of $\mathrm{F}$ is Abelian, 
one can show that many groups of homeomorphisms of the line are isomorphic to $\mathrm{F}$.  This idea can be traced back 
to \cite{brin-ubiquity} and has been largely developed and exploited  in \cite{KKL}.  The key starting point is the so-called 
{\bf {\em chain lemma}}, which is the content of the next exercise.

\begin{small}
\begin{ejer} Let $I=[u,v]$ and $I'=[u',v']$ be intervals such that $u<u'<v<v'$, and let $f$ and $g$ 
be two homeomorphisms of the real line having $I$ and $I'$ as support, respectively. 
Show that if $gf(u') \geq v$, then the group $\langle f, g\rangle$ is isomorphic to $\mathrm{F}$.  

\noindent \underbar{Hint}. Use the finite presentation from Exercise \ref{F:fin-pres} to show that the map  
$a\mapsto gf$ and $b\mapsto g$ extends to a homomorphism from $\mathrm{F}$ to $\langle f,g\rangle$. 
Then use Corollary \ref{cor cocientes de F} to conclude that this is an isomorphism.

\noindent \underbar{Remark.} Under the hypothesis above, one can actually show that the action of $\langle f,g \rangle$ 
on $[u,v']$ is semiconjugate to (the canonical action of) $\mathrm{F}$. 
\end{ejer}
\end{small}

We close this section with an exercise concerning another realization of the group $\mathrm{F}$ that will be useful 
in \S \ref{ss:fin-pres} (see Exercise \ref{hurwitz} therein).

\begin{small}
\begin{ejer} \label{F:ejer-conj}
Let us consider the (binary) Cantor set $\{0,1\}^{\mathbb{N}}$ endowed with the product topology. 

\noindent (i) Let $\hat{\mathrm{F}}$ be the group of homeomorphisms of $\{ 0,1 \}^{\mathbb{N}}$ generated by 
the maps
$$\hat{a} (\xi) := \left \{ \begin{array}{l}
0 \eta \hspace{0.8cm} \mbox {if } \xi = 0 0 \eta,
\\ \\
1 0 \eta \hspace{0.62cm} \mbox{if } \xi = 01 \eta, 
\\ \\
11 \eta \hspace{0.62cm} \mbox {if } \xi = 1 \eta,
\end{array} \right.
\qquad \mbox{ and } \qquad 
\hat{b} (\xi) := \left \{ \begin{array}{l}
\xi \hspace{1.2cm} \mbox {if } \xi = 0 \eta,
\\ \\
1 0 \eta \hspace{0.8cm} \mbox{if } \xi = 100 \eta, 
\\ \\
11 0 \eta \hspace{0.62cm} \mbox {if } \xi = 101 \eta,
\\ \\
1 1 1 \eta \hspace{0.62cm} \mbox{if } \xi = 11 \eta.
\end{array} \right.$$
Show that $\hat{\mathrm{F}}$ is isomorphic to F.

\noindent{\underbar{Hint.}} Note that $\hat{a}$ and $\hat{b}$, respectively,  may be represented by the 
same tree diagrams of the elements $a \sim f_0$ and $b \sim f_1$ of $\mathrm{F}_1 \sim \mathrm{F}$.  

\noindent (ii) Let $\phi_2 : \{0,1\}^{\mathbb{N}} \to [0,1]$ be defined by 
$$\phi_2 ( \xi) = \sum_{j \geq 1} \frac{i_j}{2^{j}}, \qquad \mbox{ where } \xi = (i_1,i_2,\ldots), \, i_j \in \{0,1\}.$$ 
Check that $\phi_2$ is one-to-one except at points that correspond to dyadic rational numbers. Besides, show that 
$\phi_2$ semiconjugates the action of $\hat{F}$ on $\{0,1\}^{\mathbb{N}}$ to that of F on $[0,1]$, in the sense that 
$$ f (\phi_2 (\xi)) = \phi_2 (\hat{f} (\xi))$$
holds for each $f \!\in\! \mathrm{F}$ and all $\xi \!\in\! \{0,1\}^{\mathbb{N}}$, 
where $\hat{f} \!\in\! \hat{\mathrm{F}}$ denotes the element corresponding to $f$.
\end{ejer}

\begin{ejer}\label{F:ejer-extended}
Given a finite binary sequence $s$, we let $\hat a_s$ be the map that consists of the action 
of $\hat a$ localized at the subtree starting at the terminal vertex of the path $s$. In precise terms, 
$$\hat a_s (\xi ) := \left \{ \begin{array}{l}
s \hat a (\eta) \hspace{0.9cm} \mbox {if } \xi = s \eta,
\\ \\
\xi \hspace{1.64cm} \mbox{otherwise}. 
\end{array} \right.$$
Note that $\hat a_1 = \hat b$.

\noindent (i) Prove that all elements $\hat a_s \in \mathrm{F}$ 
are conjugate to one of $\hat a, \hat a_0, \hat a_1, \hat a_{10}$, 
and that none of these elements is conjugate to another one in this list. 
(Note that $\hat a_0 = \hat a^{-1} \hat a_1^{-1} \hat a^2$.)

\noindent (ii) For each pair of finite binary sequences $s,t$, we let $\hat a_t (s)$ be the image of $s$ under 
$\hat a_t$ in case it is defined, which happens either when $s$ and $t$ are incompatible (which 
means that none of them extends the other one, in which case $\hat a_t (s) = s$) or when $s$ starts with 
$t00$, $t01$ or $t1$. Show that F, viewed as a group generated by the elements $\hat a_s$, admits the presentation 
$$\big\langle \hat a_s \!: \hat a_t \hat a_s \hat a_t^{-1} = a _{\hat a_t (s)} 
\mbox { for all } s,t \mbox{ such that } \hat a_t(s) \mbox { is defined} \big\rangle.$$

\noindent{\underline{Hint.}} First check that all these relations are satisfied in F. Moreover, by identifying 
$\mathrm{F} \sim \mathrm{F}_2$, note that the presentation above contains that of $\mathrm{F}_2$, as $c_0$ 
identifies with $\hat a$ and $c_k$ to $\hat a_{1^k}$ for $k \geq 1$ (where $1^k$ stands for a $1$ repeated $k$ 
times), and $\hat a_{1^k} (1^n) = 1^{n+1}$ for all $k < n$. 
\end{ejer}

\end{small}


\subsection{Some relatives of F}
\label{ejemplificando-relatives of F} 
\index{Group!piecewise analytic}

\hspace{0.45cm}

The group of piecewise real-analytic, orientation-preserving diffeomorphisms of the interval is
bi-orderable. (Note that this group contains F.) Indeed, we may let $f$ to be positive if and only
if the point $x_f^- := \inf \{x \!: f (x) \neq x \}$ is such that $f(y) > y$ for every $y > x$ sufficiently
close to $x_f^-$. Restricted to F, this bi-order coincides with $\preceq_{x^-}^+$. Extensions of
$\preceq_{x^-}^-$, $\preceq_{x^+}^+$ and $\preceq_{x^+}^-$ can be defined in an analogous way.
Similarly, groups of piecewise real-analytic, orientation-preserving diffeomorphisms of the real
line that behave nicely close to infinity are also bi-orderable.

\index{Group!piecewise projective}
The group of piecewise analytic diffeomorphisms contains a remarkable subgroup, namely that of {\bf \em piecewise projective 
diffeomorphisms}. Recall that  the projective line $\mathbb P^1 (\mathbb R)$ is  the set of lines of $\mathbb R^2$ passing 
throughout the origin. It identifies with $\mathbb R \cup \{ \infty \}$, where a real $y\in \mathbb R$ corresponds to the line 
$\mathbb R (y,1)$, and the point $\infty$ to the line $\mathbb R (0,1)$. The group $\text{PGL} (2,\mathbb R)$ acts 
on $\mathbb P^1(\mathbb R)$: via the identification $\mathbb P^1 (\mathbb R) \simeq \mathbb R \cup \{ \infty \}$, 
this action is given by 
\[ \left( \begin{array}{cc}
a & b  \\
c & d
 \end{array} \right) y = \frac{ay + b} {cy + d} .\]
Note that this action factors throughout $\mathrm{PSL} (2,\mathbb{R})$.

Let $\mathrm{PP}_+ (\mathbb{R})$ be the subgroup of $\text{Homeo}_+ (\mathbb R)$ consisting of the homeomorphisms 
that coincide with a projective map on each piece of a subdivision of $\mathbb R$ into finitely many intervals. This 
group contains the group $\mathrm{PAff}_+ (\mathbb{R})$ of {\bf \em piecewise affine homeomorphisms}, which itself 
contains $\mathrm{F}$. The next classical theorem was established by Brin and Squier \cite{BS} for the group 
$\mathrm{PAff}_+(\mathbb{R})$. The extension to $\mathrm{PP}_+ (\mathbb{R})$ appears in the work 
of Monod \cite{monod}.

\vsp 

\index{Brin-Squier Theorem}
\begin{thm} \label{th:BS}
{\em The group $\mathrm{PP}_+ (\mathbb{R})$  does not contain any non-Abelian free subgroup.}
\end{thm}

\noindent {\bf Proof.} For each $f \in \mathrm{PP}_+ (\mathbb{R})$, 
denote by $supp_{\mathrm{o}} (f)$ the {\em open support} of $f$, that is, the set of points $x$ such that $f (x)\neq x$. 
This is a finite union of disjoint open intervals. Given $g,h$ in $\mathrm{PP}_+ (\mathbb{R})$, the union 
$supp_{\mathrm{o}} (g) \cup supp_{\mathrm{o}} (h)$ is also a finite number of disjoint 
open intervals $I_1,\ldots, I_n$. We claim that the following property holds: If $[a,b]$ is a compact interval contained in one of these intervals $I_k$, 
with $1 \leq k \leq n$, then there exists a word $w$ in $g$ and $h$ for which $w ([a,b])$ is disjoint from $[a,b]$ (observe that $w ([a,b])$ is still 
contained in $I_k$). Otherwise, the supremum of the orbit of $a$ under $\langle g,h \rangle$ would be smaller than or equal to 
$b$, which is absurd due to the definition of open supports.

Assume for a contradiction that there exist two elements $g, h$  in $\mathrm{PP}_+ (\mathbb{R})$ such that any reduced 
(nontrivial) word in $g$ and $h$ is nontrivial. 
Observe that the map $f_0 : = [ [g,h], [g^2 , h]]$ is the identity close to the endpoints of each $I_k$. Indeed, if we change the projective coordinates on a 
neighborhood of an endpoint of $I_k$ so that this is moved to infinity, then the maps $g$ and $h$ become affine on each half of this neighborhood. 
Thus,  $[g,h]$ and $[g^2, h]$ become translations, hence commute. 

Note that the reduced expression of $f_0$ in $g$ and $h$ is 
$$ f_0 = g h g^{-1} h^{-1} g^2 h g^{-1} h^{-1} g^{-1} h g^2 h^{-1} g^{-2}.$$ 
Thus, $f_0$ is nontrivial. Let hence $f$ be a nontrivial element in $\langle g,h \rangle$ that is the identity on neighborhoods of the endpoints of each $I_j$, 
and such that the number of components $I_j$ intersecting the support of $f$ is minimal among all elements verifying these properties. Choose one 
of these components $I_k$, and let $[a,b] \subset I_k$ be a compact subinterval of $I_k$ such that $f$ is the identity on $I_k \setminus [a,b]$. By 
the claim above, there exists a word $w$ in $g$ and $h$ such that $w ([a,b])$ is disjoint from $[a, b]$. Inside $I_k$, the support of $w f w^{-1}$ is 
hence disjoint from the support of $f$, and thus the restrictions of $f$ and $wfw^{-1}$ to $I_k$ generate a subgroup isomorphic to $\mathbb Z^2$.  
As a consequence, the number of components $I_j$ in restriction to which  $[f, w f w^{-1}]$ is not the identity is strictly smaller that the 
corresponding number for $f$. Since $[f, wfw^{-1} ]$ is the identity on a neighborhood of the endpoints of each $I_j$, we conclude by 
minimality that $[f, w f w^{-1}]$ is the identity everywhere, and therefore $f$ and $w f w^{-1}$ generate a group isomorphic to $\mathbb Z^2$. 
Nevertheless, such a subgroup cannot arise inside a non-Abelian free group on two generators. $\hfill\square$

\vspace{0.3cm}

We will see in \S \ref{sec-amenable-ex} that the group $\mathrm{PP}_+ (\mathbb{R})$ contains many interesting finitely-generated 
subgroups. Among them, the most remarkable are those that are both non-amenable and finitely presented, the existence of which 
has been recently proved in  \cite{lodha}. Actually, these are the first examples of non-amenable, finitely-presented, torsion-free 
groups containing no free subgroup in two generators. (Examples of finitely-presented, non-amenable groups without free 
subgroups but containing many torsion elements were already known; see
\cite{OS}.)


\subsection{Braid groups}
\label{ejemplificando-5}

\index{Group!braid}
\hspace{0.45cm} One of the most relevant examples of left-orderable groups
are the braid groups $\mathbb{B}_n$. Recall that $\mathbb{B}_n$ has a presentation of the form
\begin{small}
$$\mathbb{B}_n \!=\! \big\langle \sigma_1,\ldots,\sigma_{n-1}\!\!: \esp
\sigma_i \sigma_{i+1} \sigma_i= \sigma_{i+1} \sigma_i \sigma_{i+1}
\mbox{ for } 1 \leq i \leq n-2, \esp \sigma_i \sigma_j =
\sigma_j \sigma_i \mbox{ for } |i-j| \geq 2 \big\rangle.$$
\end{small}Following Dehornoy \cite{dehornoy-libro},
for $i \!\in\! \{1,\ldots,n-1\}$, an element of $\mathbb{B}_n$ is
said to be $i$-positive if it may be written as a word of the form
\esp $w_1 \sigma_i w_2 \sigma_i \cdots w_k \sigma_i w_{k+1},$
\esp where the $w_i$'s are (perhaps trivial) words on \esp
$\sigma_{i+1}^{\pm 1}, \ldots, \sigma_{n-1}^{\pm 1}$ (and
$\sigma_i$ appears at least once).
An element in $\mathbb{B}_n$ is said to be $D$-positive if it is
$i$-positive for some $i \!\in\! \{1,\ldots,n-1\}$. A remarkable result
of Dehornoy establishes that the set of $D$-positive elements form the
positive cone of a left-order order $\preceq_{_D}$ on $\mathbb{B}_n$. In other words:

\vspace{0.1cm}

\noindent -- For every nontrivial $\sigma \in \mathbb{B}_n$, either $\sigma$ or $\sigma^{-1}$
is $i$-positive for some $i$. (Actually, Dehornoy provides an algorithm, called
{\em handle reduction}, to recognize positive elements and put 
them in the form above.)

\vspace{0.1cm}

\noindent -- If $\sigma \in \mathbb{B}_n$ is nontrivial, then $\sigma$ and $\sigma^{-1}$
cannot be simultaneously $D$-positive.

\vspace{0.1cm}
\index{Order!Dehornoy left-order}

\noindent We call $\preceq_D$ the {\bf\em Dehornoy left-order} of $\mathbb{B}_n$. 
Note that $\mathbb{B}_n$ is not bi-orderable, as it contains nontrivial elements
that are conjugate to their inverses, as for example:
$$\s_1 \s_2^{-1} = (\s_1\s_2\s_1)^{-1}  (\s_1 \s_2^{-1})^{-1}  (\s_1\s_2\s_1) .$$
In spite of this, $\preceq_{_D}$ satisfies an important
weak property of bi-invariance called {\bf{\em subword property}}:
All conjugates of the generators $\s_i$ are $\preceq_{_D}$-positive.

None of the statements above is easy to prove; see for example \cite{DDRW}.
In \S \ref{fin-gen}, we will give a short proof for the case of $\mathbb{B}_3$.

\vspace{0.25cm}

\index{Magnus' expansion}
\noindent{\bf Pure braid groups.} According to Falk and Randell \cite{FR}, pure braid
groups $\mathrm{P}\mathbb{B}_n$ are residually torsion-free nilpotent, hence bi-orderable.
An alternative approach to the bi-orderability of $\mathrm{P}\mathbb{B}_n$ using the Magnus
expansion was proposed
by Rolfsen and Zhu in \cite{bi-rolfsen} (see also \cite{Marin}). Let us point out,
however, that these bi-orders are quite different from the Dehornoy left-order.
Indeed, we will see in \S \ref{classic-conrad} that, for $n \geq 5$, no bi-order on
$\mathrm{P}\mathbb{B}_n$ can be extended to a left-order of $\mathbb{B}_n$.

For nice bi-orderable groups which are a mixture of pure braid groups and
Thompson's groups, see \cite{burillo}.


\section{Other Forms of Orderability}

\subsection{Lattice-orderable groups}
\label{lattices}

\index{Order!lattice} 
\hspace{0.45cm} A {\bf{\em lattice-ordered group}} (or
{\bf {\em $\ell$-ordered group}}) is a partially ordered group
$(\Gamma,\preceq)$ such that $\preceq$ is left and right invariant, and for each
pair of group elements $f$, $g$, there is a minimal (resp. maximal) element $f \vee g$
(resp. $f \wedge g $) simultaneously larger (resp. smaller) than $f$ and $g$.
(Note that $f \wedge g = (f^{-1} \vee g^{-1})^{-1}$.) For instance, the group
$\mathcal{A} (\Omega,\leq)$ of {\em all} order automorphisms of a totally ordered
space $(\Omega,\leq)$ is $\ell$-orderable, as one may define $f \succeq g$ whenever
$f,g$ in $\mathcal{A} (\Omega,\leq)$ satisfy $f(w) \geq g(w)$ for \esp
{\em all} \esp $\omega \in \Omega$. In this case, we have
$f \vee g \esp (w) = \mathrm{max}_{\leq} \{ f(w),g(w) \}$
and $f \wedge g \esp (w) = \mathrm{min}_{\leq} \{ f(w),g(w) \}$.

\begin{small}\begin{ex} For the group $\mathrm{Homeo}_+(\mathbb{R})$, this
and its reverse order (which is also an $\ell$-order) are the only possible
$\ell$-orders; see \cite{hollito} for a beautiful proof of this nice result 
(compare Example \ref{example-Mul}).
\end{ex}\end{small}

Conversely to the construction above, there is the following important theorem due to Holland
\cite{holland-dos}. (The proof below is taken from \cite[Chapter VII]{botto}; see
also \cite[Chapter 7]{glass}, \cite[Appendix I]{glass-antiguo}, and \cite{holland}.)

\vsp

\begin{thm} \label{holl}
{\em Every $\ell$-ordered group $(\Gamma,\preceq)$ acts by automorphisms of a totally
ordered space $(\Omega,\leq)$ in such a way that $f \preceq g$ implies $f(w) \leq g(w)$
for all $w \!\in\! \Omega$, and $f \vee g \esp (w) = \mathrm{max}_{\leq} \{ f(w),g(w) \}$
and $f \wedge g \esp (w) = \mathrm{min}_{\leq} \{ f(w),g(w) \}$. In particular, every
$\ell$-orderable group is left-orderable.}
\end{thm}

\vsp

To motivate the proof, we start by noting the following: 
Let $(\Gamma,\preceq)$ be an $\ell$-subgroup of $\mathcal{A} (\Omega,\leq)$ 
for a totally ordered space $(\Omega,\leq)$, and let $P$ be the set of {\em non-negative} 
elements for the associate order. If for each $w \in \Omega$ we denote by $P_{w}$
the semigroup $\{f \in \Gamma \!: f(w) \geq w \}$, then:

\noindent (i) \esp \esp $\bigcap_{w \in \Omega} P_{w} = P$,

\noindent (ii) \esp \esp
$P_{w} \bigcup P_{w}^{-1} = \Gamma$, for all $w \in \Omega$.

\noindent This turns natural the following version of Proposition \ref{preorders}.

\vsp

\begin{lem} {\em Let $(\Gamma,\preceq)$ be an $\ell$-ordered group
with set of non-negative elements $P$. Assume that $\Gamma$ contains a family of
subsemigroups $P_{\lambda}$, $\lambda \!\in\! \Lambda$, satisfying} (i) {\em and}
(ii) {\em above. Then the conclusion of Theorem} \ref{holl} {\em is satisfied.}
\end{lem}

\noindent{\bf Proof.} Proceed as in the proof of Proposition \ref{preorders}.
Since $\preceq$ is bi-invariant, for each $f \!\in\! P$, $g \in \Gamma$, and
$\lambda \in \Lambda$, we have $g^{-1} f g \in P \subset P_{\lambda}$. By
definition, this implies that $f (g \Gamma_{\lambda}) \geq g \Gamma_{\lambda}$.
Finally, if $f \notin P$ then, by (ii), we have
$f \in P_{\lambda}^{-1} \setminus P_{\lambda}$ for some $\lambda$,
which yields $f \Gamma_{\lambda} < \Gamma_{\lambda}$. The claims
concerning $f \vee g$ and $f \wedge g$ are left to the reader.
$\hfill\square$

\vspace{0.35cm}

\noindent{\bf Proof of Theorem \ref{holl}.} Denote by $P$ the set of non-negative
elements of $\preceq$, and for each $h \in \Gamma \setminus P$ choose a maximal
$\ell$-subsemigroup $P_h$ of $\Gamma$ containing $P$ but not $h$.
We obviously have
$$\bigcap_{h \in \Gamma  \setminus  P} P_h = P,$$
so that condition (i) above is satisfied. The proof of condition (ii) is by contradiction.
Assume throughout that for certain $h \in \Gamma \setminus P$ and $g \in \Gamma$,
we have $g \notin P_h$ and $g \notin P_h^{-1}$.

\vsp\vsp

\noindent{\underline{Claim (i).}} Neither $id \wedge g$ nor
$id \wedge g^{-1}$ belong to $P_h$.

\vsp\vsp

Indeed, note that $g = [g (id \wedge g)^{-1}](id \wedge g) = (id \vee g)(id \wedge g)$.
Since $id \vee g$ belongs to $P \subset P_h$, if $id \wedge g$ were contained in $P_h$,
then this would imply that $g$ also belongs to $P_h$, contrary to our hypothesis. A
similar argument applies to $id \wedge g^{-1}$.

\vsp\vsp

\noindent{\underline{Claim (ii).}} There exist $n_1,n_2$ in $\mathbb{N}$ and $h_1,h_2$ in $P_h$
such that $[(id \wedge g) h_1]^{n_1} \preceq h$ and $[(id \wedge g^{-1}) h_2]^{n_2} \preceq h$.

\vsp\vsp

By the maximality of $P_h$, the element $h$ belongs to the smallest $\ell$-subsemigroup
$\langle P_h, id \wedge g \rangle_{\ell}$ (resp. $\langle P_h, id \wedge g^{-1} \rangle_{\ell}$)
containing $P_h$ and $id \wedge g$ (resp. $P_h$ and $id \wedge g^{-1}$). Thus, the claim follows
from the following fact: For each $f \prec id$, the semigroup $\langle P_h, f \rangle_{\ell}$
is the set $S$ of elements which are larger than or equal to
$(f \bar{f})^n$ for some $\bar{f} \in P_h$ and some $n \!\in\! \mathbb{N}$.
To show this, first note that this set is an $\ell$-semigroup. Indeed, if
$(f \bar{f}_1 )^{n_1} \preceq g_1$ and $(f \bar{f}_2)^{n_2} \preceq g_2$,
with $\bar{f}_1,\bar{f}_2$ in $P_h$ and $n_1,n_2$ in $\mathbb{N}$, then
both $g_1$ and $g_2$ are larger than or equal to $(f \bar{f})^n$, where
$\bar{f} := id \wedge \bar{f}_1 \wedge \bar{f}_2 \in P_h$ and
$n \!:=\! \max \{n_1,n_2 \}$. Hence, $(f \bar{f})^{2n} \preceq g_1 g_2$ and
$(f \bar{f})^n \preceq g_1 \wedge g_2$, and therefore $g_1 g_2$ and $g_1 \wedge g_2$
(as well as $g_1 \vee g_2$) belong to $S$. Since $f \in S$ and $P_h \subset S$,
this shows that $\langle P_h,f \rangle_{\ell} \subset S$. Finally, we also have
$S \subset \langle P_h,f \rangle_{\ell}$. Indeed, if $\bar{g} \succ (f \bar{f})^n$
for some $\bar{f} \in P_h$ and $n \!\in\! \mathbb{N}$, then since $(f \bar{f})^n
\in \langle P_h, f \rangle_{\ell}$, we have
$\bar{g} := \big( \bar{g} (f \bar{f})^{-n} \big) (f \bar{f})^n \in P \!\cdot\!
\langle P_h, f \rangle_{\ell} = \langle P_h, f \rangle_{\ell}$. This shows the claim.

\vsp\vsp

\noindent{\underline{Claim (iii).}} Let $n := \max \{n_1,n_2\}$ and
$f := id \wedge h_1 \wedge h_2$,
where $h_1,h_2$ and $n_1,n_2$ are as in (ii). Then the element
$\hat{f} := \big( [(id \wedge g)f] \vee [(id \wedge g^{-1})f]\big)^{2n-1}$
is smaller than or equal to $h$.

\vsp\vsp

Indeed, since $f$, $id \wedge g$, and $id \wedge g^{-1}$ lie in $P^{-1}$,
we have $[ (id \wedge g) f]^n \preceq h$ and $[ (id \wedge g^{-1}) f]^n \preceq h$.
Now, as $id \wedge g$ and $id \wedge g^{-1}$ commute (their product equals
$id \wedge g \wedge g^{-1}$), we easily check that $\hat{f}$ may be rewritten as
\begin{footnotesize}
$$[(id \wedge g)f]^{2n-1} \vee \big( [(id \wedge g)f]^{2n-2}[(id \wedge g^{-1})f] \big)
\vee \big( [(id \wedge g)f]^{2n-3}[(id \wedge g^{-1})f]^{2} \big) \vee \ldots \vee
[(id \wedge g^{-1})f]^{2n-1}\!.$$
\end{footnotesize}Each term of this $\vee$-product contains either
$[ (id \wedge g) f]^n$ or $[ (id \wedge g^{-1}) f]^n$ together with
non-positive factors. The claim follows.

\vsp\vsp

To conclude the proof, note that from \esp $(id \wedge g) \vee (id \wedge g^{-1}) = id$,
\esp it follows that $\hat{f} = f^{2n - 1}$. Since $f \!\in\! P_h$, the same holds for
$\hat{f}$. Nevertheless, as $\hat{f}^{-1} h \in P \subset P_h$, this implies that
$h = \hat{f} (\hat{f}^{-1} h) \in P_h$, which is a contradiction. $\hfill\square$

\vspace{0.35cm}

\index{Conjugate roots property (C.R.P.)}
\noindent{\bf Left-orderable groups vs. $\ell$-orderable groups.} Let us
point out that $\ell$-orderability is a stronger property than left-orderability.
For instance, $\ell$-orderability is a non-local property \cite[Theorem 2.D]{glass}.
A more transparent difference concerns roots of elements, as all $\ell$-orderable
groups satisfy the {\bf {\em conjugate roots property (C.R.P.)}}: Any two elements
$f,g$ satisfying $f^n = g^n$ for some $n \!\in\! \mathbb{N}$ are
conjugate. Indeed, if for such $f,g$ we let \esp
$$h := f^{n-1} \vee f^{n-2} g \vee f^{n-3} g^2 \vee \ldots \vee g^{n-1},$$
then we have
$$fh = f^{n} \vee f^{n-1} g \vee f^{n-2} g^2 \vee \ldots \vee fg^{n-1} =
f^{n-1} g \vee f^{n-2} g^2 \vee \ldots \vee f g^{n-1} \vee g^n = hg.$$
This property fails to be true for left-orderable groups, as shown
by the next exercise.

\index{Group!Klein-bottle}
\begin{small} \begin{ejer} The $\pi_1$ of the Klein bottle may be presented in the form
$\langle a,b \!: bab = a\rangle$. (This is nothing but the infinite dihedral group.) This group
is easily seen to be left-orderable (see \S \ref{fin-gen} for a discussion on this). Prove that this
group is not $\ell$-orderable by showing that the elements $x= ba$ and $y= a$ satisfy
$x^2 = y^2$ but are not conjugate. \label{klein-1}

\vspace{0.05cm}

\noindent\underline{Remark.} Despite this example, note that every left-orderable group $\Gamma$
embeds into a lattice-orderable group, namely the group of all order permutations of $\Gamma$
endowed with a left-order.
\end{ejer} \end{small}
\index{Order!lattice}

\begin{small}
\begin{rem} Left-orderable groups satisfying the C.R.P. are not necessarily
$\ell$-orderable. Concrete examples are braid groups: in \S \ref{ejemplificando-5}, we
will see that these groups are left-orderable, the C.R.P. for them is shown in \cite{juan},
and the fact that $\mathbb{B}_n$ is not $\ell$-orderable (for $n \!\geq\! 3$) is proved in
\cite{medve-lattice}.
\end{rem} \end{small}


\subsection{Locally-invariant orders and diffuse groups}
\label{LOG}
\index{Order!locally invariant}
\hspace{0.45cm} Following \cite{chiswell} and the references therein, a partial order
relation $\preceq$ on a group $\Gamma$ is said to be {\bf{\em locally invariant}} if for
every $f,g$ in $\Gamma$, with $g \neq id$, either $f g \succ f$ or $f g^{-1} \succ f$.
Obviously, every left-order is a locally-invariant order. Examples of non left-orderable
groups admitting a locally-invariant order have been recently given: see \S \ref{next}
(see also Theorem \ref{LIO-amenable} for the case of amenable groups).

\begin{small}\begin{ejer} \label{ej-LOG}
Show that a group $\Gamma$ admits a locally-invariant order if
and only if there exist a partially ordered space $(\Omega,\leq)$ and a
map $\varphi: \Gamma \rightarrow \Omega$ such that for every $f,g$ in $\Gamma$,
with $g \neq id$, either $\varphi (fg) > \varphi(f)$ or
$\varphi (fg^{-1}) > \varphi(f)$.
\end{ejer}

\begin{ejem} \label{hyper}
Based on \cite{bow,del,hair}, it is shown in
\cite{chiswell} that many groups with hyperbolic properties
admit locally-invariant orders. More precisely, let $(X,d)$ be
a geodesic $\delta$-hyperbolic metric space \cite{ghys-harpe} and $\Gamma$
a group acting on $X$ by isometries so that \hspace{0.02cm} $d (x,g(x)) > 6\delta$
\hspace{0.02cm} holds for all $x \!\in\! X$. Then the function $g \mapsto d (x_0,g(x_0))$ 
satisfies the property of the preceding exercise for every prescribed $x_0 \in X$. 
In particular, $\Gamma$ admits a locally-invariant order.

This construction applies to many groups. In particular, if $\Gamma$ is a
residually finite Gromov-hyperbolic group (as for instance the $\pi_1$ of a compact
hyperbolic manifold), then $\Gamma$ contains a finite-index subgroup admitting a
locally-invariant order. Similarly, a group acting isometrically and freely 
on a real-tree has a locally-invariant order.
\end{ejem}\end{small}

\vsp

\index{Cone!field}
At first glance, the notion of locally-invariant order may look strange. Perhaps
a more clear view is provided by an equivalent formulation in terms of cones. More
precisely, given a group $\Gamma$, denote by $P (\Gamma)$ the family of subsets
(cones) $P \subset \Gamma$ such that $id \notin P$ and, for all $g \neq id$, at
least one of the elements $g,g^{-1}$ lies in $P$. A {\bf \em field of cones} is a map
$f \to P_f$ from $\Gamma$ into $P (\Gamma)$. This field will be said
to be {\em equivariant} if the following condition holds (see Figure 7): 
$$\mbox{if} \quad g \in P_f \quad \mbox{and} \quad h \in P_{fg},
\quad \mbox{ then } \quad gh \in P_{f}.$$

\vsp\vsp



\begin{center}
\includegraphics[scale=0.124]{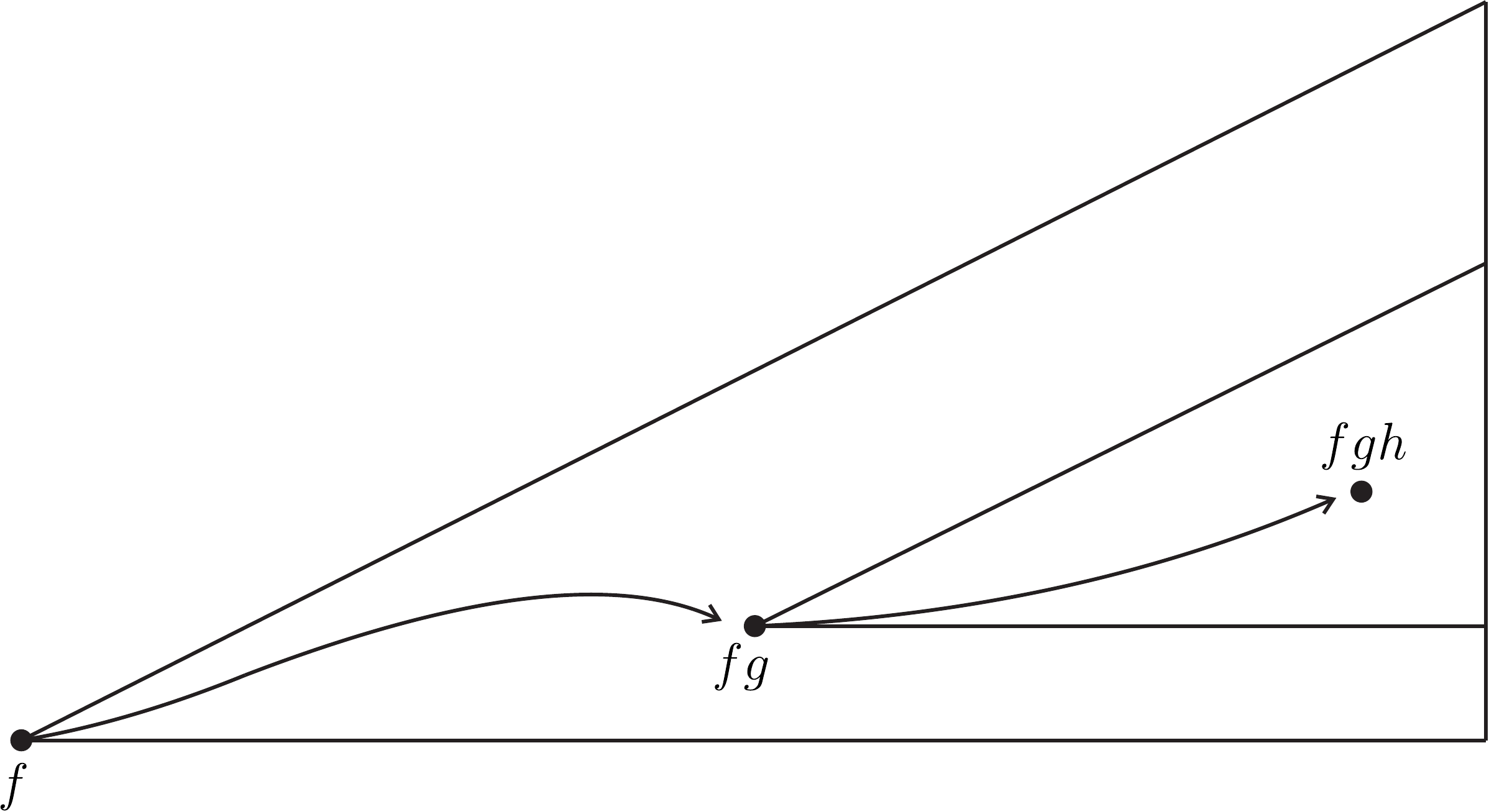}
\hspace{0.3cm}
\includegraphics[scale=0.124]{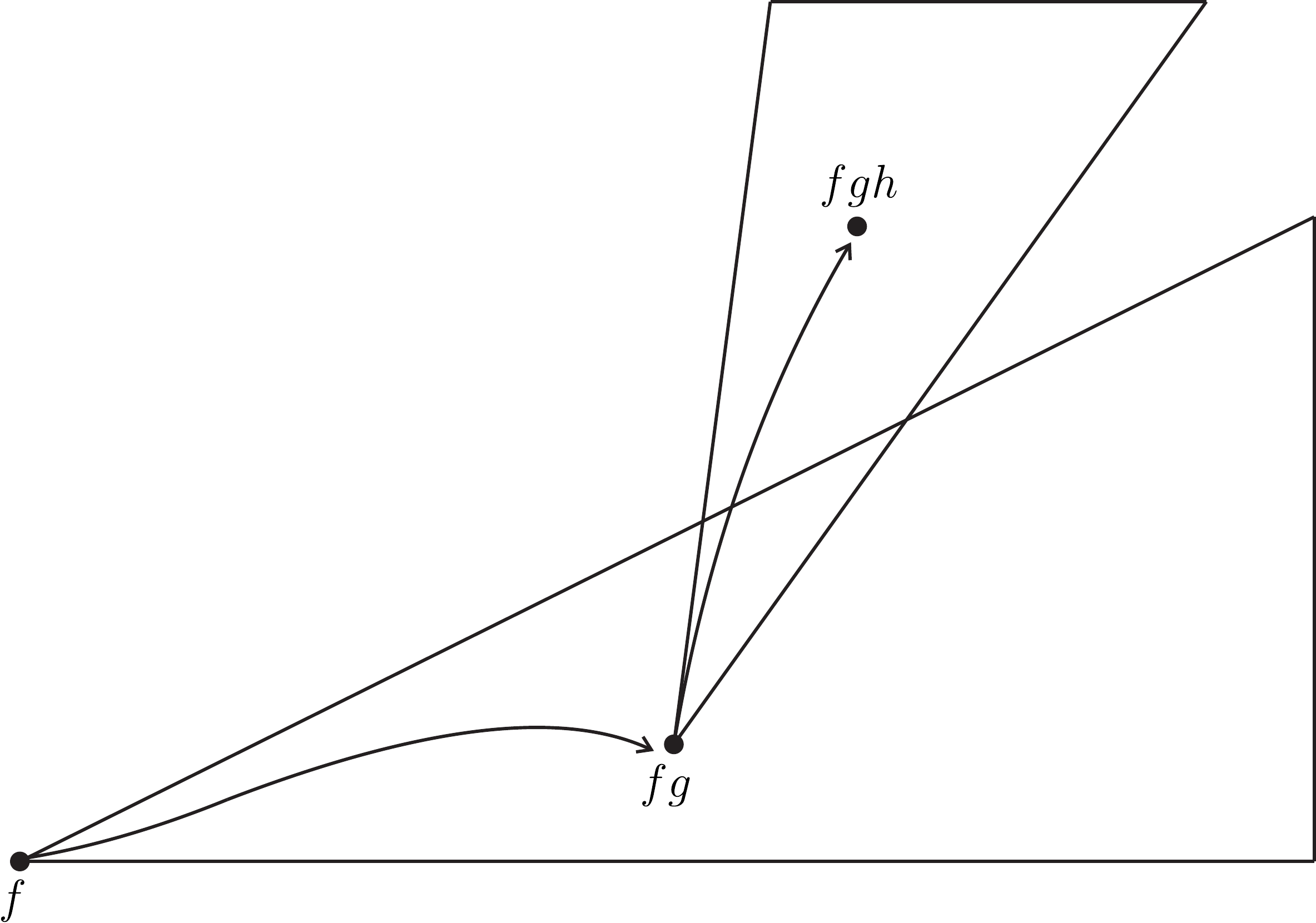}
\end{center}

\begin{center}
Figure 7: The cone condition (image on the left) and its negation (on the right) for locally-invariant orders.
\end{center}


\vspace{0.05cm}


It turns out that locally-invariant orders and equivariant fields of cones are equivalent
notions. Indeed, assume that $\preceq$ is a locally-invariant order on a group $\Gamma$.
For each $f \in \Gamma$, define $P_f$ by letting
$$g \in P_f \quad \mbox{ if and only if } \quad fg \succ f.$$
By definition, each $P_f$ belongs to $P (\Gamma)$.
We claim that the field $f \to P_f$ is equivariant. Indeed, the conditions $g \in P_f$
and $h \in P_{fg}$ mean, respectively, that $fg \succ f$ and $fgh \succ fg$. Hence, by
transitivity of $\preceq$, we have $fgh \succ f$, that is, $gh \in P_f$, as desired.

Conversely, let $f \mapsto P_f$ be an equivariant field of cones. Define a relation
$\preceq$ on $\Gamma$ by letting $f \succ g$ whenever $g^{-1} f \in P_g$. We claim
that this is a locally-invariant order. To see that $\preceq$ is antisymmetric, assume
$f \succ g$ and $g \succ f$. Then $g^{-1}f \in P_g$ and $f^{-1}g \in P_f$. By equivariance,
this implies that $id = (g^{-1}f)(f^{-1}g) \in P_g$, which is a contradiction.
To see that $\preceq$ is transitive, assume $f \succ g$ and $g \succ h$. Then
$g^{-1} f \in P_g$ and $h^{-1} g \in P_h$. By equivariance,
$h^{-1}f = (h^{-1}g) (g^{-1}f) \in P_h$, which means that $f \succ h$.
Finally, given $f \in \Gamma$ and $g \neq id$, we have either $f^{-1}gf \in P_f$,
or $f^{-1}g^{-1}f \in P_f$. In the former case, $gf \succ f$, and in the latter,
$g^{-1}f \succ f$.

\begin{small}\begin{ejer}\label{LIOs-on-Z}
Associated to each $\ell \in \mathbb{Z}$ there is a locally-invariant order $\preceq_{\ell}$
on $\mathbb{Z}$ defined by $m \succ_{\ell} n$ if and only if either $n > m \geq \ell$
or $n < m \leq \ell-1$. (See Figure 8.) Show
that every locally-invariant order on $\mathbb{Z}$ either is the canonical one, its reverse,
or contains one of the orders $\preceq_{\ell}$. (Note that we may enlarge
$\preceq_{\ell}$ by defining non-contradictory inequalities
between integers $m,n$ such that $m > \ell > n$.)

\begin{small}\begin{ejer}\label{facilito}
Show that every group admitting a locally-invariant order is torsion-free.
\end{ejer}\end{small}

\vspace{0.95cm}


\beginpicture

\setcoordinatesystem units <0.86cm,0.86cm>

\putrule from 0 0 to 16 0

\plot 12.15 0
11.85 -0.13 /
\plot 12.15 0
11.85 0.13 /

\plot 10.15 0
9.85 -0.13 /
\plot 10.15 0
9.85 0.13 /

\plot 14.15 0
13.85 -0.13 /
\plot 14.15 0
13.85 0.13 /

\plot 5.85 0
6.15 -0.13 /
\plot 5.85 0
6.15 0.13 /

\plot 3.85 0
4.15 -0.13 /
\plot 3.85 0
4.15 0.13 /

\plot 1.85 0
2.15 -0.13 /
\plot 1.85 0
2.15 0.13 /

\put{$\bullet$} at 1 0
\put{$\bullet$} at 3 0
\put{$\bullet$} at 5 0
\put{$\bullet$} at 7 0
\put{$\bullet$} at 9 0
\put{$\bullet$} at 11 0
\put{$\bullet$} at 13 0
\put{$\bullet$} at 15 0

\put{$\ell\!-\!4$} at 1 -0.5
\put{$\ell\!-\!3$} at 3 -0.5
\put{$\ell\!-\!2$} at 5 -0.5
\put{$\ell\!-\!1$} at 7 -0.5
\put{$\ell$}       at 9 -0.5
\put{$\ell\!+\!1$} at 11 -0.5
\put{$\ell\!+\!2$} at 13 -0.5
\put{$\ell\!+\!3$} at 15 -0.5

\put{} at -0.4 0
\put{Figure 8: A locally-invariant order on $\mathbb{Z}$.} at 8 -1.4
\endpicture


\end{ejer}\end{small}


\vspace{0.55cm}

There is a closely related notion to locally-invariant orders
introduced by Bowditch in \cite{bow}. Namely, given a subset $A$
of a group $\Gamma$, an {\bf{\em extremal point}} of $A$ is a point $f \in A$ such that,
if $f g \in A$ and $f g^{-1} \in A$ for some $g \in \Gamma$, then $g = id$. A
group $\Gamma$ is said to be {\bf{\em weakly diffuse}} if every nonempty
finite subset has an extremal point.
\index{Group!weakly diffuse}

\vsp

\begin{prop} \label{trivia}
{\em A group admits a locally-invariant order if and only if it is weakly diffuse.}
\end{prop}

\noindent{\bf Proof.} Let $\Gamma$ be a group admitting a locally-invariant order
$\preceq$. Given a nonempty, finite subset $A$ of $\Gamma$, let $f$ be a maximal
element (with respect to $\preceq$) of $A$. We claim that $f$ is an extremal
point of $A$. Indeed, let
$g \in \Gamma$ such that $fg \in A$ and $f g^{-1} \in A$. If $g$ were nontrivial
then we would have either $fg \succ f$ or $fg^{-1} \succ f$. However, this
contradicts the maximality of $f \in A$.

For the proof of the converse implication, see Exercise \ref{diffuse-implies-LIO}.
$\hfill\square$

\vsp

\index{Group!diffuse}
\begin{small}\begin{ejer} According to \cite{LIOs-on-amenable}, every weakly diffuse group
is {\bf {\em diffuse}}, that is, every finite subset of cardinality larger than one has at least
two extremal points. Show this by contradiction.

\noindent{\underline{Hint.}} Assume that $A$ is a finite subset having only the identity
as an extremal point. Then the same holds for $A^{-1}$. If $A$ has more than one point,
show that $B := A \cup A^{-1}$ has no extremal point.
\label{weakly diffuse implies diffuse}
\end{ejer}\end{small}

\section{General Properties}
\label{next}

\subsection{Left-orderable groups are torsion-free}
\label{no-torsion}

\index{Group!Promislow}
\hspace{0.45cm} Indeed, if $f \succ id$ (resp. $f \prec id$) for
some left-order $\preceq$, then for all $n \in \mathbb{N}$ we have
$$f^n \succ \ldots f^2 \succ f \succ id \qquad (\mbox{resp. }
f^n \prec \ldots \prec f^2 \prec f \prec f \prec id).$$
As we have seen in \S \ref{ejemplificando-1}, the converse
is true for Abelian and more generally for nilpotent groups,
but does not hold for Abelian-by-finite groups: a classical
relevant example is the {\em Promislow group}, which is the crystallographic group 
\begin{equation}\label{cristal}
\Gamma = \langle a, b \!: a^2 b a^2 = b, b^2 a b^2 = a \rangle.
\end{equation} \index{Group!4-Klein}Here we give some properties of this group (for further details, see
\cite{bow,posmi} as well as \cite[Chapter 13]{passman}). If we let
$c:= (ab)^{-1}$, then the subgroup $\langle a^2, b^2 ,c^2 \rangle$ is
torsion-free, rank-3 Abelian, and normal. The corresponding quotient
is isomorphic to the 4-Klein group. (An order-4, non-cyclic group). The
crystallographic action on $\mathbb{R}^3$ is given by
$$a(x,y,x) = (x+1,1-y,-z),$$
$$b(x,y,z) = (-x,y+1,1-z),$$
$$c(x,y,z) = (1-x,-y,z+1).$$ 
To see that $\Gamma$ is torsion-free, first note that, since every element in the 4-Klein
group has order 2, a nontrivial, finite-order element of $\Gamma$ must have order 2.
Now let $w \in \Gamma$ be nontrivial, say $w = a^{2i} b^{2j} c^{2k} a$ (the cases
where the last factor is either $b$ or $c$ are similar). Then
$$w^2 \esp = \esp a^{2i} b^{2j} c^{2k} a a^{2i} b^{2j} c^{2k} a \esp = \esp
a^{2i} b^{2j} c^{2k} a^{2i} b^{-2j} c^{-2k} a^2 \esp = \esp a^{4i + 2}
\esp \neq \esp id.$$
Finally, note that for any choice of exponents $\varepsilon,\delta$ in
$\{-1,+1\}$, the defining relations of $\Gamma$ yield
$$(a^{\varepsilon} b^{\delta})^2 (b^{\delta} a^{\varepsilon})^2 =
a^{\varepsilon} b^{-\delta} b^{2 \delta} a^{\varepsilon}
b^{2\delta} a^{\varepsilon} b^{\delta} a^{\varepsilon}
= a^{\varepsilon} b^{-\delta} a^{2\varepsilon} b^{\delta} a^{2\varepsilon} a^{-\varepsilon}
= a^{\varepsilon} b^{-\delta} b^{\delta} a^{-\varepsilon} = id.$$
Obviously, this implies that no compatible choice of signs for $a,b$ exists, hence
$\Gamma$ is not left-orderable. (For more conceptual proofs of a different nature,
see either \S \ref{sec-upp} or Example \ref{este-funciona2}.)

\begin{small}\begin{ejer} Consider the set $G$ of triplets of the form $(u, v, w)$, where
each $u, v, w$ is either an integer or of the form $\hat{m}$, with $m \in \mathbb{Z}$.

\noindent (i) Show that the rule
$$(u_1, v_1, w_1) (u_2, v_2, w_2) := (u_1 \oplus u_2, v_1 \oplus v_2, w_1 \oplus w_2),$$
where for $m,n$ in $\mathbb{Z}$,
$$m \oplus n \!:=\! m+n, \quad
m \oplus \hat{n} \!:=\! \widehat{m+n}, \quad
\hat{m} \oplus n \!:=\! \widehat{m-n}, \quad \hat{m} \oplus \hat{n} :=  m - n,$$
endows $G$ with a group structure.

\noindent (ii) Show that Promislow's group $\Gamma$ above identifies with the subgroup of
$G$ generated by $a := (1 \hat{0} \hat{0})$ and $b := (\hat{0} 1 \hat{1})$.
\label{identifying}\end{ejer}\end{small}

\subsection{Unique roots and generalized torsion}
\label{gen-torsion}

\index{Generalized torsion}
\index{Unique root property}
\hspace{0.45cm} Bi-orderable groups have a stronger property
than absence of torsion, namely they have no {\bf{\em generalized
torsion}}: If $f \neq id$, then no product of conjugates of $f$ is the identity.
(In particular, no nontrivial element is conjugate to its inverse.) These groups also have
the {\bf{\em unique root property}}: If $f^n = g^n$ for some integer $n$, then $f = g$.
Once again, none of these properties characterizes bi-orderability (for classical examples, 
see \cite[Example 4.3.1]{botto} and \cite{baumslag,bludov-1}, respectively). Actually, 
they do not even imply left-orderability. For the latter property, a  
concrete example (taken from \cite[Chapter VII]{botto}) is
$$\Gamma_n \esp = \esp \big\langle f, g  \!: f g f g^2 \cdots f g^n
= f^{-1} g f^{-1} g^2 \cdots f^{-1} g^n = id \big\rangle,
\quad \mbox{where } n \mbox { is ``large''}.$$
For the former property, an example has been recently constructed by Cai and Clay in \cite{cai-clai}

\index{Group!Klein-bottle}
\begin{small}\begin{ejem} As we saw in Example \ref{klein-1}, the Klein bottle group is
left-orderable but does not satisfy the C.R.P., hence it is not bi-orderable. Another way to
contradict bi-orderability consists in noting that it has generalized torsion: $(a^{-1}ba) b = id$.
Moreover, the unique root property fails: $(ba)^2 = a^2$, though $ba \neq a$.
\end{ejem}

\begin{ejer} Prove that in any bi-orderable group, the following holds:
If $f$ commutes with a nontrivial power of $g$, then it commutes with
$g$. Show that this is no longer true for left-orderable groups.
\end{ejer}

\begin{ejer} Show that for bi-orderable groups, the normalizer of any finite subset coincides
with its centralizer. Again, show that this is no longer true for left-orderable groups.
\label{ex-normalizer}
\end{ejer}
\end{small}

\subsection{The Unique Product Property (U.P.P.)}
\label{sec-upp}

\index{Unique Product Property (U.P.P.)}
A group $\Gamma$ is said to have the U.P.P. if given any two finite
subsets $\{g_i\}$, $\{h_j\}$, there exists $f \in \Gamma$ that may
be written in a unique way as a product \esp $g_i h_j$.

Every left-orderable group has the U.P.P. Indeed, given two finite subsets
$A~\!\!:=~\!\!\{g_1,\dots,g_n\}$ and $B := \{h_1,\ldots,h_m\}$, let $f := g_i h_j$
be the element of $AB$ that is maximal with respect to a fixed left-order on $\Gamma$.
If $f$ were equal to $g_{i'} h_{j'}$ for some $i',j'$, then $h_j \succeq h_{j'}$, as
otherwise $f = g_i h_j \prec g_i h_{j'}$ would contradict the maximality of $f$.
Similarly, $h_{j'} \succeq h_{j}$, as otherwise $f = g_{i'} h_{j'} \prec g_{i'} h_j$. 
Thus $h_j = h_{j'}$, which yields $g_i = g_{i'}$.

Note that the minimum element in $AB$ has also a unique expression as above.
This is coherent with a result from \cite{st} reproduced below.

\begin{small}\begin{ejer}
Show that U.P.P. implies a ``double'' U.P.P., in the sense that
given any two finite subsets $A,B$ such that \esp $|A| \! + \! |B| > 2$, \esp
there exist at least {\em two} elements in $AB$ which may be written in a
unique way as a product $ab$, with $a \in A$ and $b \in B$. (Compare Exercise
\ref{weakly diffuse implies diffuse}.)

\noindent{\underline{Hint.}} Assume that a group $\Gamma$ has the U.P.P. but only $ab \in AB$
has a unique representation in $AB$, and let $C := a^{-1}A$, $D := Bb^{-1}$, $E := D^{-1}C$,
and $F := DC^{-1}$. Using the fact that, in $CD$, only $id = id \cdot id$ has a unique
representation, show that, in $EF$, no element has unique representation.
\end{ejer}

\begin{ejer}
Show that a group satisfies the U.P.P. if and only if for any finite subset
$A$ there exist at least one element in $A^2$ (actually, two) that may
be written in a unique way as a product $ab$, with $a,b$ in $A$.
\end{ejer}\end{small}

Let us remark that groups with the U.P.P. are torsion-free. Indeed, if
$f^n~=~id$ for some $f \neq id$, then the U.P.P. fails for 
$A = B = \{id,f,\ldots,f^{n-1}\}$. The converse to this remark is false.
Indeed, as Promislow showed in \cite{posmi}, the crystallographic group
of \S \ref{no-torsion} does not satisfy the U.P.P. (see Exercise
\ref{viene-viene} below for the details; see also \cite{rips}
for a different example using small cancellation techniques,
and \cite{K-R,steenbock} for more recent developments.) 
Also, U.P.P groups are non-necessarily left-orderable, as 
it was recently proved in  \cite{K-R}. We discuss this point 
in more detail below.

\vspace{0.21cm}

\index{Group!diffuse}
\index{Order!locally invariant}
\index{Unique Product Property (U.P.P.)}
\noindent{\bf Locally-invariant orders, diffuse groups, and the U.P.P.} The
U.P.P. is satisfied by all weakly diffuse groups (hence, by groups admitting a
locally-invariant order; see Proposition \ref{trivia}). Indeed, given nontrivial finite
subsets $A,B$ of a weakly diffuse group $\Gamma$, let $f \!\in\! AB$ be an extremal point
of $AB$. We claim that $f$ may be written in a unique way as $gh$, with $g \in A$ and
$h \in B$. Indeed, if $f = g_1 h_1 = g_2 h_2$, with $g_1,g_2$ in $A$ and $h_1,h_2$ in
$B$, then letting $h := h_1^{-1} h_2$ we have $fh = g_1h_2 \in AB$ and
$fh^{-1} = g_2 h_1 \in AB$. Since $f$ is an extremal point of $AB$, this
implies that $h = id$, which yields $h_1 = h_2$ and $g_1 = g_2$.

\vspace{0.1cm}

Below we elaborate on an example of a ``large'' group that satisfies the 
U.P.P but is not left-orderable. (For amenable groups the situation is unclear,
due to Theorem \ref{LIO-amenable}.) First note that, by Example \ref{hyper}, isometry
groups of hyperbolic metric spaces with ``large displacement'' admit locally-invariant
orders, and hence satisfy the U.P.P. (Actually, a combination of remarkable recent results
establishes that the $\pi_1$ of every
closed, hyperbolic 3-manifold contains a finite-index group that is bi-orderable; see
\cite{Agol,BergeronWise,DuchampThibon, HaglundWise, KahnMarkovic}.) This
motivates the following question.

\begin{question} Does there exist a sequence of compact,
hyperbolic 3-manifolds whose injectivity radius converges
to infinity and whose $\pi_1$ are non left-orderable~? (Examples
of non left-orderable 3-manifold groups appear in \cite{CD,prit}.)
\end{question}

This question seems to have an affirmative but difficult solution. Indeed,
it is not very hard to prove that, if $\Gamma$ is the $\pi_1$ of a compact,
hyperbolic 3-manifold with nontrivial first Betti number, then $\Gamma$ is
left-orderable (it is actually Conrad-orderable, in the terminology of
\S \ref{general-Conrad}; see \cite{fourier}). A sequence of compact,
hyperbolic 3-manifolds with trivial first Betti number and whose injectivity
radius converge to infinite appears in \cite{frank-cal-dunfield}. However,
it seems hard to adapt the methods therein to show that the $\pi_1$
of infinitely many of these manifolds are non left-orderable. Actually,
an obvious difficulty comes from the fact that they are virtually orderable,
as was mentioned above.

Despite the above, it was cleverly noted by Dunfield and included in the work of Kionke
and Raimbault  (see \cite{K-R}) that there is an hyperbolic 3-manifold whose $\pi_1$ is
known to be non left-orderable and for which a lower estimate of its injectivity radius allows
applying the results described in Example \ref{hyper}. This is enough to conclude that it admits 
a locally-invariant order, hence it satisfies the U.P.P. 
As Kionke and Raimbault point out,
the next question remains open.

\index{Group!weakly diffuse}
\begin{question} Does there exist a U.P.P.-group that is not weakly diffuse~?
\end{question}

\vsp


\begin{small}\begin{ejer} \label{viene-viene} Consider the subset
$$A = B = \big\{ (b a)^2, (a b)^2, a^2 b, a b a^{-1}, b, a b^{-1} a,
b^{-1}, a b a, a b^{-2}, b^2 a^{-1}, a (b a)^2, b a b, a, a^{-1} \big\}$$
of the crystallographic
group $\Gamma = \langle a,b \!: a^2 b a^2 = b, b^2 a b^2 = a \rangle$
introduced in \S \ref{no-torsion}.

\vspace{0.25cm}

\noindent (i) Show that, via the identification of Exercise \ref{identifying},
this set becomes
\begin{small}
$$(0 0 2),
\! (0 0 \underline 2),
\! (\hat 2 1 \hat 1),
\! (\hat 2 \underline 1 \underline{\hat 1}),
\!(\hat 0 1 \hat 1),
\! (\hat 0 1 \underline{\hat 1}),
\! (\hat 0 \underline 1 \hat 1),
\! (\hat 0 \underline 1 \underline{\hat 1}),
\! (1 \hat 2 \hat 0),
\! (\underline 1 \hat 2 \hat 0),
\! (1 \hat 0 \underline {\hat 2}),
\! (\underline 1 \hat 0 \hat 2),
\! (1 \hat 0 \hat 0),
\! (\underline 1 \hat 0 \hat 0),$$
\end{small}where $\underline{m}$ (resp. $\hat{\underline{m}}$)
is written instead of $-m$ (resp. $\widehat{-m}$).

\vspace{0.2cm}

\noindent (ii) In the (partial) multiplication table below, check that each value
corresponding to the product of a pair of elements appears at least twice.
\vspace{0.2cm}
$$\begin{footnotesize}
\begin{array}{l|llllllllllllll}
& \!\!(0 0 2) \!\!
& \!\!(0 0 \underline 2) \!\!
& \!\!(\hat 2 1 \hat 1) \!\!
& \!\!(\hat 2 \underline 1 \underline{\hat 1}) \!\!
& \!\!(\hat 0 1 \hat 1) \!\!
& \!\!(\hat 0 1 \underline{\hat 1}) \!\!
& \!\!(\hat 0 \underline 1 \hat 1) \!\!
& \!\!(\hat 0 \underline 1 \underline{\hat 1}) \!\!
& \!\!(1 \hat 2 \hat 0) \!\!
& \!\!(\underline 1 \hat 2 \hat 0) \!\!
& \!\!(1 \hat 0 \underline {\hat 2}) \!\!
& \!\!(\underline 1 \hat 0 \hat 2) \!\!
& \!\!(1 \hat 0 \hat 0) \!\!
& \!\!(\underline 1 \hat 0 \hat 0) \!\!\\
\hline \\
\vspace{0.1cm}
\!\! (0 0 2) \hspace{0.1cm}
& \!\! (0 0 4) \!\!
& \!\! (0 0 0) \!\!
& \!\! (\hat 2 1 \hat 3) \!\!
& \!\! (\hat 2 \underline 1 \hat 1) \!\!
& \!\! (\hat 0 1 \hat 3) \!\!
& \!\! (\hat 0 1 \hat 1) \!\!
& \!\! (\hat 0 \underline 1 \hat 3) \!\!
& \!\! (\hat 0 \underline 1 \hat 1) \!\!
& \!\! (1 \hat 2 \hat 2) \!\!
& \!\! (\underline 1 \hat 2 \hat 2) \!\!
& \!\! (1 \hat 0 \hat 0) \!\!
& \!\! (\underline 1 \hat 0 \hat 4) \!\!
& \!\! (1 \hat 0 \hat 2) \!\!
& \!\! (\underline 1 \hat 0 \hat 2) \!\!  \\
\vspace{0.1cm}
\!\! (0 0 \underline 2) \!\!
& \!\! (0 0 0) \!\!
& \!\! (0 0 \underline 4) \!\!
& \!\! (\hat 2 1 \underline {\hat 1}) \!\!
& \!\! (\hat 2 \underline 1 \underline {\hat 3}) \!\!
& \!\! (\hat 0 1 \underline {\hat 1}) \!\!
& \!\! (\hat 0 1 \underline {\hat 3}) \!\!
& \!\! (\hat 0 \underline 1 \underline {\hat 1}) \!\!
& \!\! (\hat 0 \underline 1 \underline {\hat 3}) \!\!
& \!\! (1 \hat 2 \underline {\hat 2}) \!\!
& \!\! (\underline 1 \hat 2 \underline {\hat 2}) \!\!
& \!\! (1 \hat 0 \underline {\hat 4}) \!\!
& \!\! (\underline 1 \hat 0 \hat 0) \!\!
& \!\! (1 \hat 0 \underline {\hat 2}) \!\!
& \!\! (\underline 1 \hat 0 \underline {\hat 2}) \!\!  \\
\vspace{0.1cm}
\!\! (\hat 2 1 \hat 1) \!\!
& \!\! (\hat 2 1 \underline {\hat 1}) \!\!
& \!\! (\hat 2 1 \hat 3) \!\!
& \!\! (0 2 0) \!\!
& \!\! (0 0 2) \!\!
& \!\! (2 2 0) \!\!
& \!\! (2 2 2) \!\!
& \!\! (2 0 0) \!\!
& \!\! (2 0 2) \!\!
& \!\! (\hat 1 \hat 3 1) \!\!
& \!\! (\hat 3 \hat 3 1) \!\!
& \!\! (\hat 1 \hat 1 3) \!\!
& \!\! (\hat 3 \hat 1 \underline 1) \!\!
& \!\! (\hat 1 \hat 1 1) \!\!
& \!\! (\hat 3 \hat 1 1) \!\!  \\
\vspace{0.1cm}
\!\! (\hat 2 \underline 1 \underline {\hat 1}) \!\!
& \!\! (\hat2 \underline 1 \underline {\hat 3}) \!\!
& \!\! (\hat 2 \underline 1 \hat 1) \!\!
& \!\! (0 0 \underline 2) \!\!
& \!\! (0 \underline 2 0) \!\!
& \!\! (2 0 \underline 2) \!\!
& \!\! (2 0 0) \!\!
& \!\! (2 \underline 2 \underline 2) \!\!
& \!\! (2 \underline 2 0) \!\!
& \!\! (\hat 1 \hat 1 \underline 1) \!\!
& \!\! (\hat 3 \hat 1 \underline 1) \!\!
& \!\! (\hat 1 \underline {\hat 1} 1) \!\!
& \!\! (\hat 3 \underline {\hat 1} \underline 3) \!\!
& \!\! (\hat 1 \underline {\hat 1} \underline 1) \!\!
& \!\! (\hat 3 \underline {\hat 1} \underline 1) \!\!  \\
\vspace{0.1cm}
\!\! (\hat 0 1 \hat 1) \!\!
& \!\! (\hat 0 1 \underline {\hat 1}) \!\!
& \!\! (\hat 0 1 \hat 3) \!\!
& \!\! (\underline 2 2 0) \!\!
& \!\! (\underline 2 0 2) \!\!
& \!\! (0 2 0) \!\!
& \!\! (0 2 2) \!\!
& \!\! (0 0 0) \!\!
& \!\! (0 0 2) \!\!
& \!\! (\underline {\hat 1} \hat 3 1) \!\!
& \!\! (\hat 1 \hat 3 1) \!\!
& \!\! (\underline {\hat 1} \hat 1 3) \!\!
& \!\! (\hat 1 \hat 1 \underline 1) \!\!
& \!\! (\underline {\hat 1} \hat 1 1) \!\!
& \!\! (\hat 1 \hat 1 1) \!\!  \\
\vspace{0.1cm}
\!\! (\hat 0 1 \underline {\hat 1}) \!\!
& \!\! (\hat 0 1 \underline {\hat 3}) \!\!
& \!\! (\hat 0 1 \hat 1) \!\!
& \!\! (\underline 2 2 \underline 2) \!\!
& \!\! (\underline 2 0 0) \!\!
& \!\! (0 2 \underline 2) \!\!
& \!\! (0 2 0) \!\!
& \!\! (0 0 \underline 2) \!\!
& \!\! (0 0 0) \!\!
& \!\! (\underline {\hat 1} \hat 3 \underline 1) \!\!
& \!\! (\hat 1 \hat 3 \underline 1) \!\!
& \!\! (\underline {\hat 1} \hat 1 1) \!\!
& \!\! (\hat 1 \hat 1 \underline 3) \!\!
& \!\! (\underline {\hat 1} \hat 1 \underline 1) \!\!
& \!\! (\hat 1 \hat 1 \underline 1) \!\!  \\
\vspace{0.1cm}
\!\! (\hat 0 \underline 1 \hat 1) \!\!
& \!\! (\hat 0 \underline 1 \underline {\hat 1}) \!\!
& \!\! (\hat 0 \underline 1 \hat 3) \!\!
& \!\! (\underline 2 0 0) \!\!
& \!\! (\underline 2 \underline 2 2) \!\!
& \!\! (0 0 0) \!\!
& \!\! (0 0 2) \!\!
& \!\! (0 \underline 2 0) \!\!
& \!\! (0 \underline 2 2) \!\!
& \!\! (\underline {\hat 1} \hat 1 1) \!\!
& \!\! (\hat 1 \hat 1 1) \!\!
& \!\! (\underline {\hat 1} \underline {\hat 1} 3) \!\!
& \!\! (\hat 1 \underline {\hat 1} \underline 1) \!\!
& \!\! (\underline {\hat 1} \underline {\hat 1} 1) \!\!
& \!\! (\hat 1 \underline {\hat 1} 1) \!\!  \\
\vspace{0.1cm}
\!\! (\hat 0 \underline 1 \underline {\hat 1}) \!\!
& \!\! (\hat 0 \underline 1 \underline {\hat 3}) \!\!
& \!\! (\hat 0 \underline 1 \hat 1) \!\!
& \!\! (\underline 2 0 \underline 2) \!\!
& \!\! (\underline 2 \underline 2 0) \!\!
& \!\! (0 0 \underline 2) \!\!
& \!\! (0 0 0) \!\!
& \!\! (0 \underline 2 \underline 2) \!\!
& \!\! (0 \underline 2 0) \!\!
& \!\! (\underline {\hat 1} \hat 1 \underline 1) \!\!
& \!\! (\hat 1 \hat 1 \underline 1) \!\!
& \!\! (\underline {\hat 1} \underline {\hat 1} 1) \!\!
& \!\! (\hat 1 \underline {\hat 1} \underline 3) \!\!
& \!\! (\underline {\hat 1} \underline {\hat 1} \underline 1) \!\!
& \!\! (\hat 1 \underline {\hat 1} \underline 1) \!\!  \\
\vspace{0.1cm}
\!\! (1 \hat 2 \hat 0) \!\!
& \!\! (1 \hat 2 \underline {\hat 2}) \!\!
& \!\! (1 \hat 2 \hat 2) \!\!
& \!\! (\hat 3 \hat 1 \underline 1) \!\!
& \!\! (\hat 3 \hat 3 1) \!\!
& \!\! (\hat 1 \hat 1 \underline 1) \!\!
& \!\! (\hat 1 \hat 1 1) \!\!
& \!\! (\hat 1 \hat 3 \underline 1) \!\!
& \!\! (\hat 1 \hat 3 1) \!\!
& \!\! (2 0 0) \!\!
& \!\! (0 0 0) \!\!
& \!\! (2 2 2) \!\!
& \!\! (0 2 \underline 2) \!\!
& \!\! (2 2 0) \!\!
& \!\! (0 2 0) \!\!  \\
\vspace{0.1cm}
\!\! (\underline 1 \hat 2 \hat 0) \!\!
& \!\! (\underline 1 \hat 2 \underline {\hat 2}) \!\!
& \!\! (\underline 1 \hat 2 \hat 2) \!\!
& \!\! (\hat 1 \hat 1 \underline 1) \!\!
& \!\! (\hat 1 \hat 3 1) \!\!
& \!\! (\underline {\hat 1} \hat 1 \underline 1) \!\!
& \!\! (\underline {\hat 1} \hat 1 1) \!\!
& \!\! (\underline {\hat 1} \hat 3 \underline 1) \!\!
& \!\! (\underline {\hat 1} \hat 3 1) \!\!
& \!\! (0 0 0) \!\!
& \!\! (\underline 2 0 0) \!\!
& \!\! (0 2 2) \!\!
& \!\! (\underline 2 2 \underline 2) \!\!
& \!\! (0 2 0) \!\!
& \!\! (\underline 2 2 0) \!\!  \\
\vspace{0.1cm}
\!\! (1 \hat 0 \underline {\hat 2}) \!\!
& \!\! (1 \hat 0 \underline {\hat 4}) \!\!
& \!\! (1 \hat 0 \hat 0) \!\!
& \!\! (\hat 3 \underline {\hat 1} \underline 3) \!\!
& \!\! (\hat 3 \hat 1 \underline 1) \!\!
& \!\! (\hat 1 \underline {\hat 1} \underline 3) \!\!
& \!\! (\hat 1 \underline {\hat 1} \underline 1) \!\!
& \!\! (\hat 1 \hat 1 \underline 3) \!\!
& \!\! (\hat 1 \hat 1 \underline 1) \!\!
& \!\! (2 \underline 2 \underline 2) \!\!
& \!\! (0 \underline 2 \underline 2) \!\!
& \!\! (2 0 0) \!\!
& \!\! (0 0 \underline 4) \!\!
& \!\! (2 0 \underline 2) \!\!
& \!\! (0 0 \underline 2) \!\!  \\
\vspace{0.1cm}
\!\! (\underline 1 \hat 0 \hat 2) \!\!
& \!\! (\underline 1 \hat 0 \hat 0) \!\!
& \!\! (\underline 1 \hat 0 \hat 4) \!\!
& \!\! (\hat 1 \underline {\hat 1} 1) \!\!
& \!\! (\hat 1 \hat 1 3) \!\!
& \!\! (\underline {\hat 1} \underline {\hat 1} 1) \!\!
& \!\! (\underline {\hat 1} \underline {\hat 1} 3) \!\!
& \!\! (\underline {\hat 1} \hat 1 1) \!\!
& \!\! (\underline {\hat 1} \hat 1 3) \!\!
& \!\! (0 \underline 2 2) \!\!
& \!\! (\underline 2 \underline 2 2) \!\!
& \!\! (0 0 4) \!\!
& \!\! (\underline 2 0 0) \!\!
& \!\! (0 0 2) \!\!
& \!\! (\underline 2 0 2) \!\!  \\
\vspace{0.1cm}
\!\! (1 \hat 0 \hat 0) \!\!
& \!\! (1 \hat 0 \underline {\hat 2}) \!\!
& \!\! (1 \hat 0 \hat 2) \!\!
& \!\! (\hat 3 \underline {\hat 1} \underline 1) \!\!
& \!\! (\hat 3 \hat 1 1) \!\!
& \!\! (\hat 1 \underline {\hat 1} \underline 1) \!\!
& \!\! (\hat 1 \underline {\hat 1} 1) \!\!
& \!\! (\hat 1 \hat 1 \underline 1) \!\!
& \!\! (\hat 1 \hat 1 1) \!\!
& \!\! (2 \underline 2 0) \!\!
& \!\! (0 \underline 2 0) \!\!
& \!\! (2 0 2) \!\!
& \!\! (0 0 \underline 2) \!\!
& \!\! (2 0 0) \!\!
& \!\! (0 0 0) \!\!  \\
\vspace{0.1cm}
\!\! (\underline 1 \hat 0 \hat 0) \!\!
& \!\! (\underline 1 \hat 0 \underline {\hat 2}) \!\!
& \!\! (\underline 1 \hat 0 \hat 2) \!\!
& \!\! (\hat 1 \underline {\hat 1} \underline 1) \!\!
& \!\! (\hat 1 \hat 1 1) \!\!
& \!\! (\underline {\hat 1} \underline {\hat 1} \underline 1) \!\!
& \!\! (\underline {\hat 1} \underline {\hat 1} 1) \!\!
& \!\! (\underline {\hat 1} \hat 1 \underline 1) \!\!
& \!\! (\underline {\hat 1} \hat 1 1) \!\!
& \!\! (0 \underline 2 0) \!\!
& \!\! (\underline 2 \underline 2 0) \!\!
& \!\! (0 0 2) \!\!
& \!\! (\underline 2 0 \underline 2) \!\!
& \!\! (0 0 0) \!\!
& \!\! (\underline 2 0 0) \!\! \\
\end{array}
\end{footnotesize}$$
\end{ejer}\end{small}

\vspace{0.1cm}

\begin{small}
\begin{rem}
By the preceding exercise, the crystallographic group $\Gamma$, though being torsion-free, does 
not admit a locally-invariant order. It is worth mentioning that this is actually the case of ``most" 
finitely-presented groups in a very precise random model for groups; see \cite{orlef}.
\end{rem}

\index{Unique Product Property (U.P.P.)}
\begin{rem}
As it was shown by Witte Morris, finite-index subgroups
of $\mathrm{SL}(3,\mathbb{Z})$ are non left-orderable (see Theorem \ref{witte-SL}).
For large index, these
groups are torsion-free, and it seems to be unknown whether they satisfy the
U.P.P. By Exercise \ref{ej-LOG}, the following question makes sense:
Does there exist a ``norm'' on $\mathrm{SL}(3,\mathbb{Z})$ such that
for finite but ``large'' index subgroups $\Gamma$ one has either $\|fg\| > \|f\|$
or $\|fg^{-1}\| > \|f\|$ for every $f,g$ in $\Gamma$, with $g \neq id$~?
\end{rem}\end{small}

\vsp

\index{Kaplansky zero divisor conjecture}
\noindent{\bf On Kaplansky's conjecture.} A famous question due to Kaplansky (commonly
referred to as the {\em Kaplansky zero divisor conjecture}) asks whether the group algebra of a torsion-free
group over a ring $\mathbb{A}$ has no zero-divisors provided $\mathbb{A}$ has no zero-divisors.
(Even the case where $\mathbb{A} = \mathbb{Z}$ is open.) The restriction on the torsion is natural. Indeed,
\begin{equation}\label{div-obvio}
f^n = id \quad \implies \quad
(f - 1)(f^{n-1} + f^{n-2} + \cdots + f + 1) = f^n - 1 = 0.
\end{equation}
It easily follows from the definitions that every group satisfying the U.P.P. also
satisfies the conclusion of the Kaplansky conjecture. For the crystallographic
group considered above, Kaplansky's conjecture is known to be true by different methods
(see for instance \cite{one,two,three}; see also \cite{passman,lueck}).

\begin{small}\begin{ejem}\label{burnside-ex}
Consider the {\bf{\em free Burnside group}}
$$B(m,n) := \big\langle a_1,\ldots,a_m \!: W^n \!=\! id \mbox{ for every word } W \big\rangle.$$
It is known that for $m \geq 2$ and $n$ odd and large enough, $B(m,n)$ is infinite (actually,
it is non-amenable; see \cite{ad1}). Of course, every element in this group has finite order.
However, it is still interesting to look for zero-divisors in its group algebra that are ``nontrivial''
({\em i.e.}, that do not arise from an identity of the form (\ref{div-obvio})). For instance, according to
\cite{ivanov}, this is the case for
$$A:= (1+c+\ldots+c^{n-1}) (1-aba^{-1}), \quad B:= (1-a) (1+b+\ldots+b^{n-1}),$$
where $a:=a_1$, $b:=a_2$, and $c:=aba^{-1}b^{-1}$. (Checking that $AB=0$ is an easy
exercise.)
\end{ejem}\end{small}


\subsection{More Combinatorial Properties}

\hspace{0.45cm}
Recently, orderable groups have been considered as a natural framework to extend certain basic
results of Additive Combinatorics (see \cite{Fr,ad2,ad3} as general references). One of the
most elementary ones is the inequality for product sets
\begin{equation}
|AB| \geq |A| + |B| - 1,
\label{kem}
\end{equation}
which holds for any finite subsets $A,B$ of the integers (this is an easy exercise). In this regard,
it is worth mentioning that this readily extends to finite subsets of left-orderable groups. Indeed,
modulo multiplying $B$ on the right by the largest possible element of type $h^{-1}$, where
$h \in B$, we may assume that $id$ is the smallest element of $B$. Then, if we order the
elements in $A$ (resp. $B$) in the form $g_1 \prec \ldots \prec g_n$ (resp.
$id = h_1 \prec \ldots \prec h_m$), we have
$$g_1 \prec g_2 \prec \ldots \prec g_n \prec g_n h_2 \prec \ldots \prec g_n h_m.$$

Less trivially, (\ref{kem}) still holds for finite subsets of torsion-free groups,
as it was proved by Kemperman in \cite{kemperman}.

\begin{thm} \label{kem-th}
{\em For all finite subsets $A,B$ of a torsion-free group, we have}
$$|AB| \geq |A| + |B| - 1.$$
\end{thm}

\noindent{\bf Proof.} First note that the claim of the theorem trivially holds if
either $|A|$ or $|B|$ equals 1. Moreover, changing $A$ by $g^{-1} A$ and $B$ by $B h^{-1}$
for $g \in A$ and $h \in B$, we reduce the general case to that where $id \in A \cap B$.
Assume for a contradiction that $A,B$ are finite subsets that do not satisfy (\ref{kem})
and for which the value of $m := |AB|$ is minimal, that of $n :=|A| + |B|$ is maximal
while $|AB| = m$, and that of $|A|$ is maximal while $|AB|=m$ and $|A|+|B|=n$ (the
extremal properties being realized among subsets containing $id$).

As $id \!\in\! A \cap B$, we also have
$$|AB| \geq |A| + |B| - |A \cap B|.$$
Hence, $|A \cap B| \geq 2$.
Let $H$ be the subsemigroup generated by $A \cap B$. We consider
two different cases.

\vsp\vsp

\noindent{\underline{Case I.}} We have $Af \subset A$ for all $f \in A \cap B$.

\vsp

Then, as $id \in A$, this implies
$H \subset A$. Therefore, $H$ is a finite subsemigroup of a group, hence a (finite) subgroup.
As $|H| \geq 2$, this produces torsion elements.

\vsp\vsp

\noindent{\underline{Case II.}} There exists $f \in A \cap B$ such that $Af$ is not contained in $A$.

\vsp

Fixing such an $f$, let $A' := \{ g \in A \!: gf \notin A \}$ and $B' := \{ h \in B \! : f h \notin B \}$.
There are two subcases to consider.

If $|A'| \geq |B'|$, then let $A^* \!:=\! A \cup A' f$ and $B^* \!:=\! B \setminus B'$.
(Note that $B \setminus B' \neq \emptyset$ since $id \notin B'$.)
One easily checks that $A^* B^* \subset A B$, hence $|A^*B^*| \leq |AB|$.
Moreover, $|A^*| = |A| + |A' f| = |A| + |A'|$ and $|B^*| = |B| - |B'|$, thus
$|A^*| + |B^*| \geq |A| + |B|.$ Finally, $|A^*| > |A|$, as $A'$ is nonempty.
Therefore, by the choice of $A,B$, we must have
$$|A^*B^*| \geq |A^*| + |B^*| - 1,$$
hence
$$|AB| \geq |A^*B^*| \geq |A^*| + |B^*| - 1 \geq |A| + |B| -1,$$
which is a contradiction.

If $|A'| < |B'|$, then let $A^* \!:=\! A \setminus A'$ (which is nonempty as $id \notin A'$)
and $B^* := B \cup f B$. Again, $A^* B^* \subset AB$, hence $|A^*B^*| \leq |AB|$.
Moreover, $|A^*| = |A| - |A'|$ and $|B^*| = |B| + |B'|$ yield $|A^*| + |B^*| > |A| + |B|$.
By the choice of $A,B$, this implies
$$|A^*B^*| \geq |A^*| + |B^*| - 1,$$
hence
$$|AB| \geq |A^*B^*| \geq |A^*| + |B^*| - 1 > |A| + |B| -1,$$
which is again a contradiction. $\hfill\square$

\vsp\vsp

\begin{small}\begin{ejem}\label{BF-example}
By pursuing on the technique of proof above, Brailovsky and Freiman proved in \cite{B-F} that
equality arises if and only if $A$ and $B$ are {\em geometric progressions} on different sides, 
that is, if there exist group elements $f,g,h$ and non-negative integers $n,m$ such that
$$A = \{ g, gf, \ldots, gf^{n-1} \}, \qquad B = \{ h, fh, \ldots, f^{m-1}h \}.$$
Showing such a claim for left-orderable groups is a straighforward exercise.
\end{ejem}\end{small}

Below we present another proof of Theorem \ref{kem-th} following the ideas of Hamidoune
\cite{hami-solo} that is somewhat closer to the techniques of the next section. We refer to
\cite{hami} for more details and furter developments, including an alternative proof
of the Brailovsky-Freiman theorem above.

\vsp\vsp\vsp

\noindent{\bf Another proof of Theorem \ref{kem-th}.} Given a finite subset $B$ of a group $\Gamma$,
for each finite subset $A \subset \Gamma$ we let $\partial^B  A := AB \setminus A$ (compare
(\ref{a-comparar})). Given a positive integer $k$, we say that a subset $C$ is {\bf {\em $(B,k)$-critical}}
if $|C| \geq k$ and
$$|\partial^B C| = \min \big\{ |\partial^B A | \!: |A| \geq k \big\} .$$
We say that $C$ is a {\bf{\em $(B,k)$-atom}} if it is a $(B,k)$-critical set of smallest cardinality.

\vsp\vsp

\noindent{{\underline{Claim (i).}} If $C$ is a $(B,k)$-atom and $C'$ is $(B,k)$-critical,
then either $C \subset C'$ or $| C \cap C' | \leq k-1$.}

\vsp

Indeed, assume $C$ is not contained in $C'$ and $|C \cap C'| \geq k$.
Then, by definition,
$$|\partial^B C| < |\partial^B (C \cap C')|.$$
Let $C_*$ (resp. $C'_*$) be the complement of $C \cup \partial^B C$
(resp. $C' \cup \partial^B C'$). On the one hand, we have
\begin{small}\begin{eqnarray*}
|\partial^B C \cap C'| \!+\! |\partial^B C \cap \partial^B C'| \!+\! |\partial^B C \cap C'_*|
&=&
|\partial^B C| \\
\!\!&<&\!\! |\partial^B (C \cap C')|\\
\!\!&\leq&\!\!
|C \cap \partial^B C'| \!+\! |\partial^B C \cap C'| \!+\! |\partial^B C \cap \partial^B C'|,
\end{eqnarray*}
\end{small}hence $|\partial^B C \cap C'_*| < |C \cap \partial^B C'|$.
On the other hand, we have
\begin{small}\begin{eqnarray*}
|\partial^B C' \cap C| \!+\! |\partial^B C' \cap \partial^B C| \!+\! |\partial^B C' \cap C_*|
&\leq&
|\partial^B C'| \\
&\leq& |\partial^B (C' \cap C)|\\
&\leq&
|C'_* \cap \partial^B C| \!+\! |\partial^B C' \cap C_*| \!+\! |\partial^B C' \cap \partial^B C|,
\end{eqnarray*}
\end{small}hence
$| \partial^B C' \cap C | \leq | C'_* \cap \partial^B C|$. These two conclusions are certainly
in contradiction.

\vsp\vsp

\noindent{{\underline{Claim (ii).}} If $C$ is a $(B,k)$-atom and $g \neq id$, then
$|C \cap gC| \leq k-1$.}

\vsp

Indeed, the set $gC$ is a $(B,k)$-atom as well. Moreover, we cannot have $gC \subset C$,
otherwise $g$ would be a torsion element. (If $gC \subset C$, then $g C = C$, so that
$g$ acts as a permutation of $C$ and therefore $g^n h \!=\! h$ for all $h \in C$.)

\vsp\vsp

\noindent{{\underline{Claim (iii).}} For all finite sets $A,B$, we have $|AB| \geq |A| + |B| - 1$.}

\vsp

Indeed, we may assume that $B$ contains $id$. Let $C$ be a $(B,1)$-atom. Again, we may
assume $id \in C$. If $C$ contains another element $g$, then $| C \cap g C | \geq 1$, a
contradiction to (ii). Hence, $C = \{ id \}$. Therefore, for every (nonempty) finite subset $A$,
$$|AB| - |A| = |AB \setminus A| \geq |CB \setminus C| = |B \setminus \{id\}| = |B| - 1,$$
which shows the claim. $\hfill\square$

\vsp\vsp

\begin{small}\begin{rem}
It is conjectured that $k$-atoms have cardinality equal to $k$ for torsion-free groups.
This holds for instance  for groups satisfying the U.P.P. (This is an easy exercise; see
\cite[Lemma 4]{hami} in case of problems.)
\end{rem}\end{small}

\vsp\vsp

A direct consequence of the preceding theorem is the inequality $|A^2| \geq 2|A|-1$
for all finite subsets $A$ of torsion-free groups. The next result from \cite{FHLM}
improves this inequality for non-Abelian bi-orderable groups.

\vsp

\begin{thm} {\em Let $A$ be a finite subset of a bi-orderable group. If $|A^2| \leq 3|A|-3$,
then the subgroup generated by the elements of $A$ is Abelian.}
\end{thm}

\noindent{\bf Proof.} The proof is by induction on $|A|$. If $|A| = 2$, say $A = \{f_1,f_2\}$,
then $|A^2| \leq 3 |A| - 3 =3$ implies $A^2 = \{f_1^2, f_1 f_2 = f_2 f_1, f_2^2\}$, since we 
cannot have $f_1^2 = f_2^2$ for $f_1 \neq f_2$ in a bi-orderable group. Therefore, the group
generated by $f_1,f_2$ is Abelian. Assume that the theorem holds for subsets of cardinality $\leq k$,
and let $A := \{ f_1, \ldots, f_{k+1} \}$, where $f_i \prec f_j$ holds whenever $i < j$ for a fixed bi-order
$\preceq$ on the underlying group. We let $i$ be the maximal index for which the subgroup generated
by $B := \{ f_1, \ldots, f_i \}$ is Abelian, and we assume that $i \leq k$. Then $f_{i+1}$ does not belong
to the subgroup generated by $B$. Moreover, there is $f \in B$ not commuting with $f_{i+1}$; we
let $f_j$ be the maximal such element. We also let $C := \{f_{i+1}, \ldots, f_{k+1} \}$.
Assume throught that $|A^2| \leq 3|A| - 3$.

\vsp\vsp

\noindent{\underline{Claim (i).}} We have $|C^2| \leq 3|C| - 3$.

\vsp

Indeed, using $\preceq$ (see also Exercise \ref{ex-normalizer}) and the fact that $f_{i+1}$ 
does not belong to the group generated by $B$, one readily checks that
\begin{equation}\label{relaciones}
B^2 \cap (f_{i+1} B \cup B f_{i+1}) = \emptyset, \quad f_{i+1} B \neq B f_{i+1},
\quad C^2 \cap (B^2 \cup f_{i+1} B \cup B f_{i+1}) = \emptyset.
\end{equation}
Therefore,
\begin{eqnarray*}
|C^2|
&\leq&
|A^2| - |B^2| - |f_{i+1} B \cup B f_{i+1}|\\
&\leq& (3|A| - 3) - (2|B|-1) - (|B| + 1)
\hspace{0.18cm} = \hspace{0.18cm}
3 (|A| - |B|) - 3
\hspace{0.18cm} = \hspace{0.18cm}
3 |C| - 3,
\end{eqnarray*}
as claimed.

By the inductive hypothesis, the group generated by $C$ is Abelian. As $f_j$ and $f_{i+1} \in C$
do not commute, we must have
\begin{equation}\label{otra}
C^2 \cap (f_{j} C \cup C f_{j}) = \emptyset.
\end{equation}

\vsp\vsp

\noindent{\underline{Claim (ii).} We have $B f_{i+1} \cap f_j C = \{f_j f_{i+1} \}$. In particular,
$|B f_{i+1} \cup f_j C| = k$.

\vsp

Indeed, assume $f_m f_{i+1} = f_j f_n$, with $f_m \in B$ and $f_n \in C$. If $f_m \prec f_j$,
then $f_{i+1} \succ f_n$, which is impossible since $f_{i+1}$ is the smallest element in $C$.
If $f_m \succ f_j$, then $f_m$ commutes with $f_{i+1}$, and so does $f_j = f_m f_{i+1} f_n^{-1}$,
which is absurd.

\vsp

\vsp\vsp

\noindent{\underline{Claim (iii).} We have $B^2 \cap (f_j C \cup C f_j) = \emptyset$.

\vsp

Indeed, assume $f_m f_n = f_j f_{\ell}$ holds for $f_m,f_n$ in $B$ and $f_{\ell} \in C$. Since
$f_n \prec f_{\ell}$, we must have $f_m \succ f_j$. Moreover, as $B$ generates an Abelian
group, $f_m f_n = f_n f_m$, hence also $f_n \succ f_j$. Therefore, both $f_m,f_n$ commute
with $f_{i+1}$, and so does $f_j = f_m f_n f_{\ell}^{-1}$, which is absurd. This shows that
$B^2 \cap f_j C = \emptyset$. That $B^2 \cap Cf_j = \emptyset$ is proved similarly.

\vsp\vsp

\noindent{\underline{Claim (iv).} We have $A^2 = B^2 \cup C^2 \cup B f_{i+1} \cup f_j C$.

\vsp

It follows from the above that
\begin{eqnarray*}
|B^2 \cup C^2 \cup B f_{i+1} \cup f_j C|
&=&
|B^2| + |C^2| + |B f_{i+1} \cup f_j C|\\
&\geq&
(2i-1) + (2(k-i+1)-1) + k
\hspace{0.18cm}=\hspace{0.18cm}
3 (k+1) - 3.
\end{eqnarray*}
By the hypothesis $|A^2| \leq 3|A|-3$, this implies the claim.

\vsp

Note that $f_{i+1} f_j \notin B^2$ and $f_{i+1} f_j \notin C^2$, by (\ref{relaciones}).
A contradiction is then provided by the two claims below.

\vsp\vsp

\noindent{\underline{Claim (v).} We have $f_{i+1} f_j \notin B f_{i+1}$.

\vsp

Indeed, assume $f_{i+1} f_j = f_m f_{i+1}$ for $f_m \in B$. If $f_m \succ f_j$, then $f_m$ commutes
with $f_{i+1}$. Thus, $f_{i+1} f_j = f_{i+1} f_m$, hence $f_j = f_m$, which is absurd. Suppose
$f_m \prec f_j$. By Exercise \ref{ex-normalizer}, there exists $f_n \in B$ such that
$f_{i+1} f_n \notin B f_{i+1}$. Thus, necessarily, $f_n \neq f_j$.  We cannot have
$f_n \succ f_j$, otherwise $f_{i+1} f_n = f_n f_{i+1} \in B f_{i+1}$, a contradiction.
Therefore, $f_n \prec f_j$.

By (\ref{relaciones}), $f_{i+1} f_n \notin B^2 \cup C^2$. Hence, by Claim (iv), we have $f_{i+1} f_n \in f_j C$,
so that there is $f_{\ell} \in C$ such that $f_{i+1} f_n = f_j f_{\ell}$. As $B$ generates an Abelian subgroup,
$$f_j f_{\ell} f_j = f_{i+1} f_n f_j = f_{i+1} f_j f_n.$$
Since $f_n \prec f_j$, this implies $f_j f_{\ell} \prec f_{i+1} f_j = f_m f_{i+1}$.
However, this is impossible, because $f_j \succ f_m$ and $f_{\ell} \succeq f_{i+1}$.

\vsp\vsp

\noindent{\underline{Claim (vi).} We have $f_{i+1} f_j \notin f_{j} C$.

\vsp

Assume $f_{i+1} f_j = f_j f_m$ holds for a certain $f_m \in C$. By Exercise \ref{ex-normalizer},
there exists $f_n \!\in\! C$ such that $f_n f_j \notin f_j C$. Since $f_{i+1} f_j \in f_j C$, it holds
$f_{i+1} \prec f_n$. Moreover, by (\ref{otra}) and Claim (iii), $f_n f_j \notin B^2 \cup C^2$.
Thus, by Claim (iv), we have $f_n f_j \in B f_{i+1}$. Let $f_{\ell} \in B$ be such that
$f_n f_j = f_{\ell} f_{i+1}$.

Note that $f_{\ell} \neq f_{j}$, otherwise $f_n f_j$ would belong to $f_j C$. If $f_{\ell} \succ f_j$,
then it commutes with $f_{i+1}$, and so does $f_j = f_n^{-1} f_{i+1} f_{\ell}$, which is absurd.
If $f_{\ell} \prec f_j$, then, as $f_{i+1}$ and $f_n$ commute,
$$f_{i+1} f_{\ell} f_{i+1} = f_{i+1} f_n f_j = f_n f_{i+1} f_j.$$
As $f_{i+1} \prec f_n$, this implies $f_{\ell} f_{i+1} \succ f_{i+1} f_j = f_j f_m$. However,
this is impossible, since $f_{\ell} \prec f_j$ and $f_{i+1} \preceq f_m$. $\hfill\square$

\vsp\vsp

\begin{small}\begin{ejem}
Following \cite{FHLM}, let $A = A_k$ be the subset of the Baumslag-Solitar
group $BS(1,2) := \langle a,b \!: aba^{-1} = b^2 \rangle$ given by
$A := \{a,ab,ab^2,\ldots,ab^{k-1} \}$. Check that $|A^2| = 3 |A| - 2$,
yet $BS(1,2)$ is bi-orderable and non-Abelian.
\end{ejem}

\begin{ejer}
Let $\Gamma$ be the Klein bottle group $\langle a,b \!: aba^{-1} = b^{-1} \rangle$
(see Example \ref{klein-1}).
Check that the set $A = A_k := \{a,ab,ab^{-1}, \ldots ab^{k-1}\}$ satisfies $|A^2| = 2 |A| - 1$,
yet $\Gamma$ is left-orderable and non-Abelian.
\end{ejer}

\begin{ejer} Using Brailovsky-Freiman's theorem (see Example \ref{BF-example}), prove that
if $A$ is a subset of a torsion-free group satisfying $|A^2| = 2|A|-1$, then $A$ generates either
an Abelian subgroup or a group isomorphic to the Klein bottle group.
\end{ejer} \end{small}


\subsection{Isoperimetry and Left-Orderable Groups}

\hspace{0.45cm} The aim of this section is to develop some ideas 
introduced by Gromov in \cite{gromov}. Let us begin with
the notion of isoperimetric profile, due to Vershik.

Let $\Gamma$ be a finitely-generated group acting on a set $X$, and
let $\mathcal{G}$ be a finite generating system containing $id$. For a 
subset $Y \subset X$, its {\bf {\em boundary}} (with respect to $\mg$) is
defined as
\begin{equation}\label{a-comparar}
\partial_{\mg} Y := \mg  Y \setminus Y,
\end{equation}
where $\mg Y := \{g(y) \!: g \in \mg, \esp y \in Y\}$.
The maximal function $I\!: \mathbb{N} \rightarrow \mathbb{N}$
satisfying
$$|\partial_{\mg} Y| \geq I \big( |Y| \big)$$
for all finite $Y \subset X$ is called the {\bf {\em combinatorial isoperimetric profile}}
of the $\Gamma$-action, and will be denoted by $I_{(X;\Gamma,\mathcal{G})}$.
\index{Isoperimetric profile ! combinatorial profile}

An important case arises when $X = \Gamma$ is endowed with the action by left-translations.
In this case, the isoperimetric profile is denoted $I_{(\Gamma,\mg)}$. We list below some
important properties.

\vspace{0.35cm}

\noindent{\bf{Subaditivity.}} \esp \esp If $\Gamma$ is infinite, then for all $r_1,r_2$ we have
$$I_{(\Gamma,\mg)} (r_1 + r_2) \esp \leq \esp I_{(\Gamma,\mg)}(r_1) + I_{(\Gamma,\mg)}(r_2).$$

Indeed, choose $Y_i \subset \Gamma$ such that $|Y_i| = r_i$ and
$\big| \partial_{\mg} Y_i \big| = I_{(\Gamma,\mg)}(r_i)$, with $i \in \{1,2\}$.
Since $\Gamma$ is infinite, after ``moving'' $Y_2$ keeping $Y_1$ fixed, 
we may assume that $Y_1$ and $Y_2$ are disjoint and \esp
$\partial_{\mg} (Y_1 \sqcup Y_2) = \partial_{\mg} Y_1 \sqcup \partial_{\mg} Y_2$. \esp 
(Note that
$h (\partial_{\mg} Y) = \partial_{\mg} (hY)$, for all $h \in \Gamma$ and all $Y \subset \Gamma$.)
This yields
$$I_{(\Gamma,\mg)} (r_1 + r_2) \esp \leq \esp \big| \partial_{\mg} (Y_1 \sqcup Y_2) \big| \esp
= \esp \big| \partial_{\mg} (Y_1) \big| + \big| \partial_{\mg} (Y_2) \big| \esp
= \esp I_{(\Gamma,\mg)} (r_1) + I_{(\Gamma,\mg)} (r_2).$$

\vspace{0.35cm}

\noindent{\bf{$I$ is non-decreasing under extensions.}} \esp \esp
If $\Gamma \subset \Gamma_1$ are infinite groups
and $\mg \subset \mg_1$, then, for all $r$,
$$I_{(\Gamma_1,\mg_1)}(r) \geq I_{(\Gamma,\mg)}(r).$$

Indeed, any finite subset
$Y \!\subset\! \Gamma_1$ may be decomposed as a disjoint union $Y = \bigsqcup_{i=1}^k Y_i$,
where the points in each $Y_i$ are in the same class modulo $\Gamma$. Since $\mg \subset \Gamma$,
$$\partial_{\mg_1} Y \supset \partial_{\mg} Y = \bigsqcup_{i=1}^k \partial_{\mg} Y_i.$$
Thus,
$$\big| \partial_{\mg_1} Y \big| \geq \sum_{i=1}^k \big| \partial_{\mg} Y_i \big|
\geq \sum_{i=1}^k I_{(\Gamma,\mg)} \big( |Y_i| \big)
\geq I_{(\Gamma,\mg)} \Big( \sum_{i=1}^k |Y_i| \Big)
= I_{(\Gamma,\mg)} \big( |Y| \big).$$

\vspace{0.35cm}

\noindent{\bf{$I$ is non-increasing under homomorphisms.}} \esp \esp If
$\Phi \!: \Gamma \rightarrow \underline{\Gamma}$ is a surjective
group homomorphism and $\underline{\mg} = \Phi(\mg)$, then,
for all $r$,
$$I_{(\Gamma,\mg)} (r) \geq I_{(\underline{\Gamma},\underline{\mg})}(r).$$

Indeed, given a finite subset $Y \subset \Gamma$, we let \esp
$\underline{Y}_m :=
\big\{f \in \underline{\Gamma} \!: | \Phi^{-1} (f) \cap Y | \geq m \big\}$.
Clearly, \esp $|Y| = \sum_{m \geq 1} |\underline{Y}_m|$. \esp If we are able
to show that
\begin{equation}
|\partial_{\mg} Y | \geq \sum_{m \geq 1} \big| \partial_{\underline{\mg}} \underline{Y}_m \big|,
\label{astuta}
\end{equation}
then this would yield
$$|\partial_{\mg} Y| \esp \geq \esp
\sum_{m \geq 1} I_{(\underline{\Gamma},\underline{\mg})} \big( |\underline{Y}_m| \big) \esp
\geq \esp I_{(\underline{\Gamma},\underline{\mg})} \Big( \sum_{m \geq 1} |\underline{Y}_m| \Big)
\esp = \esp I_{(\underline{\Gamma},\underline{\mg})} \big( |Y| \big),$$
thus showing our claim. Now, to show (\ref{astuta}), every $\underline{f}$
in $\partial_{\underline{\mg}} \underline{Y}_m$ may be written
as $\Phi(g) \underline{h}$ for some $g \!\in\! \mg$ and
$\underline{h} \!\in\! \underline{Y}_m$. By definition,
$\big| \Phi^{-1}(\underline{h}) \cap Y \big| \geq m$, and
$\big| \Phi^{-1}(\underline{f}) \cap Y \big| < m$. Thus, there
must be some $g_f$ in $\mathcal{G}$ and $h \in \Phi^{-1}(\underline{h})$
such that $f = \Phi(g_f h)$, with $h \in Y$ and $g_f h \notin Y$.
The correspondence $\underline{f} \mapsto g_f h$ from
$\bigcup_m \partial_{\underline{\mg}} \underline{Y}_m$ to $\partial_{\mg} Y$
is injective, because $\Phi(g_f h) = \underline{f}$. This shows (\ref{astuta}).

\vspace{0.3cm}

Suppose that $\Gamma$ acts on a linear space $\mathbb{V}$. Given a subspace
$D \subset \mathbb{V}$ and a finite generating set $\mg \subset \Gamma$
containing $id$, we define its {\bf {\em boundary}} as the quotient space
$$\partial_{\mg} D : = \mathcal{G} \cdot Y \esp \big/ \esp D,$$
where $\mg \cdot D$ is the subspace generated by $\{g(v)\!: g \in \mg, v \in D\}$.
Using now the notation $|\cdot|$ for the dimension of a vector space, we define
the {\bf{\em linear isoperimetric profile}} \index{Isoperimetric profile ! linear profile}
 of the $\Gamma$-action on $\mathbb{V}$ as the
maximal function $I$ satisfying, for all finite dimensional subspaces $D \subset \mathbb{V}$,
$$|\partial_{\mg} D| \geq I \big( |D| \big).$$
We denote this function by $I^*_{(\mathbb{V};\Gamma,\mg)}$. In the special case where
$\mathbb{V}$ is the
group algebra $\mathbb{R}(\Gamma)$ (viewed as the vector space of finitely-supported,
real-valued functions on $\Gamma$), we simply use the notation $I^*_{(\Gamma,\mg)}$.

As is the case of $I_{(\Gamma,\mg)}$, the function $I^*{(\Gamma,\mg)}$ is
subadditive, as well as non-increasing under group extensions. It is unclear
whether it is non-increasing under group homomorphisms. However, we will
see that if the target group is left-orderable, then this property holds
(see Proposition \ref{nigh}).

There is a simple relation between $I$ and $I^*$ for all finitely-generated groups.

\vspace{0.1cm}

\begin{prop} \label{a-dualizar}
{\em For every finitely-generated group $\Gamma$ and all $r \geq 0$,}
\vspace{-0.12cm}
$$I_{(\Gamma,\mg)} (r) \esp \geq \esp I^*_{(\Gamma,\mg)} (r).$$

\end{prop}

\noindent{\bf Proof.} To each finite subset $Y \subset \Gamma$ we may associate the subspace
$D_{Y} := Y^{\mathbb{R}}$ formed by all the functions whose support is contained in $Y$.
We clearly have $|Y| = |D_Y|$ and $|\partial_{\mg} Y| = |\partial_{\mg} D_Y|$,
which easily yields the claim. $\hfill\square$

\vspace{0.4cm}

The opposite inequality is not valid for all groups, as the following example shows.

\vsp

\begin{small}\begin{ex}
If $\Gamma$ is a group containing a nontrivial finite subgroup $\Gamma_0$
and $\Gamma_0 \subset \mg$, then for any finite dimension subspace
$D \subset \mathbb{R} (\Gamma)$
of finitely-supported functions which are constant along the cosets of $\Gamma_0$,
we have $I^*_{(\Gamma,\mg)} \big( |D| \big) < I_{(\Gamma,\mg)} \big( |D| \big)$.
If $\Gamma$ is infinite, this yields \esp
$I^*_{(\Gamma,\mg)} ( r ) < I_{(\Gamma,\mg)} ( r )$, \esp for all $r \geq 0$.
\end{ex}\end{small}

Despite the preceding example, the equivalence between $I$ and $I^*$ holds
for left-orderable groups. (It is an open question whether this remains true
for torsion-free groups; for groups with torsion, see Example \ref{laurent}.) 
The proof of this fact (due to Gromov) is reproduced below.

\vsp

\begin{thm} \label{gromov-thm}
{\em If $\Gamma$ is a finitely-generated left-orderable group, then
for every finite generating system $\mg$ containing $id$, one has
$I_{(\Gamma,\mg)} = I^*_{(\Gamma,\mg)} $.}
\end{thm}

\vsp

To show this theorem, we will use a somewhat ``dual argument'' 
to that of Proposition \ref{a-dualizar}.

\vspace{0.3cm}

\index{Isometric domination (ID)}
\noindent{\bf{Isoperimetric Domination (ID).}} \esp \esp Let $\Gamma$ be a
group acting on a set $X$ and on a vector space $\mathbb{V}$. Suppose there
exists an {\em equivariant} map $D \mapsto Y_D$ from the Grassmanian $\mathrm{Gr}_\mathbb{V}$
of finite dimensional subspaces of $\mathbb{V}$ to the family of subsets of $X$ such that:

\vsp

\noindent (i) \esp \esp $|D| = |Y_D|$, \esp for all \esp $D \in \mathrm{Gr}_{\mathbb{V}}$;

\vsp

\noindent (ii) \esp \esp $|\mathrm{span} (\bigcup_i D_i)| \geq |\bigcup_i Y_{D_i}|$,
\esp for every finite family \esp $\{D_i\} \subset \mathrm{Gr}_{\mathbb{V}}$.

\vsp

\noindent We claim that, in this case, for every finite generating set
$\mg$ containing $id$ and all $ r\geq 0$,
\begin{equation}
I_{(X;\Gamma,\mg)} (r) \geq I^*_{(\mathbb{V};\Gamma,\mg)} (r).
\label{mayor}
\end{equation}
Indeed, taking any $D$ so that $|D| = r$, we have $|Y_D| = r$ and
\begin{multline*}
|\partial_{\mg} D| = \big| \mg \cdot D / D \big|
= \Big| \mathrm{span} \big(\!\bigcup_{g \in \mg} g D \big) \Big| - |D| \geq \\
\geq \Big| \bigcup_{g \in \mg} Y_{gD} \Big| - |D|
\geq \Big| \bigcup_{g \in \mathcal{G}} g (Y_D) \Big| - |Y_D|
= |\partial_{\mg} Y_D|,
\end{multline*}
which easily yields (\ref{mayor}).

\vspace{0.3cm}

\noindent{\bf{ID for left-ordered groups.}} \esp \esp In view of the above discussion,
in order to prove Theorem (\ref{gromov-thm}) it suffices to exhibit an ID from
$\mathrm{Gr}_{\mathbb{R}(\Gamma)}$ to $2^{\Gamma}$.
The construction proceeds as follows. Fix a
left-order $\preceq$ on $\Gamma$. To each finitely-supported, real-valued function
$\varphi$ on $\Gamma$, we may associate the minimum $g \!\in\! \Gamma$ in its
support (where the {\em minimum} is taken with respect to $\preceq$). Denote this
point by $g_{\varphi}$. Now, if $D \subset \mathbb{V}$ is a finitely-dimensional subspace,
then the number of points $g_{\varphi}$ which may appear for some $\varphi \in D$ is finite.
In fact, a simple ``passing to a triangular basis'' argument using the left-order shows that
the cardinality of this subset $Y_D \subset \Gamma$ equals $|D|$, so property (i) above is
satisfied. Property (ii) is also easily verified, thus concluding the proof.

\vsp

\begin{prop}
\label{nigh}
{\em Let $\Phi \!: \Gamma \rightarrow \underline{\Gamma}$ be a
surjective group homomorphism. If $\underline{\Gamma}$ is left-orderable,
then denoting $\underline{\mg} = \Phi(\mg)$ we have, for all $r \geq 0$,}
$$I^*_{(\Gamma,\mg)} (r) \esp \geq \esp
I^*_{(\underline{\Gamma},\underline{\mg})}(r).$$
\end{prop}

\noindent{\bf Proof.} Fix a left-order $\preceq$ on $\underline{\Gamma}$,
and for each $\underline{g} \in \underline{\Gamma}$ denote
$$\mathbb{A}_{\prec \underline{g}}
:= \{\varphi \in \mathbb{A}(\Gamma) \!: g_{\varphi} \prec \underline{g} \},
\qquad \mathbb{A}_{\preceq \underline{g}}
:= \{\varphi \in \mathbb{A}(\Gamma) \!: g_{\varphi} \preceq \underline{g}\}.$$
Given a finitely dimensional subspace $D \subset \mathbb{A}(\Gamma)$, define
$U = U_D \!: \Gamma \rightarrow \mathbb{N}_0$ by
$$U (g) := dim \Big( \mathbb{A}_{\preceq \Phi(g)} \cap D
\big/ \mathbb{A}_{\prec \Phi(g)} \cap D \Big).$$
Let $S_{U}$ be the subgraph of $U$, that is,
$$S_U := \big\{ (g,n) \in \Gamma \times \mathbb{N} \! : \esp U(g) \geq n \big\}.$$
Since $\Gamma$ naturally
acts on $\Gamma \times \mathbb{N}$ and the action is free on each level, we have
\begin{eqnarray*}
|\partial_{\mg} S_U|
= \sum_{m \geq 1} \big| \partial_{\mg} \big( S_U \!\cap\! (\Gamma \times \{ m \}) \big) \big|
\!\!\! &\geq& \!\!\!
\sum_{m \geq 1} I_{(\Gamma,\mg)} \big( \big| S_U \!\cap\! (\Gamma \times \{ m \} ) \big| \big)\\
\!\!\! &\geq& \!\!\! I_{(\Gamma,\mg)} \Big(\! \sum_{m \geq 1}
|S_U \!\cap\! (\Gamma \times \{ m \})| \! \Big)
\! = I_{(\Gamma,\mg)} \big( |S_U| \big).
\end{eqnarray*}
Moreover, one easily convinces that $|S_U| = |D|$. \esp
Putting all of this together, we obtain
$$|\partial_{\mg} D| \esp
= \esp |\partial_{\mg} S_U| \esp
\geq \esp I_{(\Gamma,\mg)} \big( |S_U| \big) \esp
= \esp I_{(\Gamma,\mg)} \big( |D| \big) \esp
= \esp I^*_{(\Gamma,\mg)} \big( |D| \big),$$
where the last equality comes from Theorem \ref{gromov-thm}. $\hfill\square$

\vsp

\index{Group!lamplighter}
\begin{small}\begin{ex}
Following \cite{bartholdi}, we consider the {\bf{\em lamplighter group}}
$\Gamma := \mathbb{Z} \wr \mathbb{Z} / 2\mathbb{Z}
= \mathbb{Z} \ltimes \oplus_{i \in \mathbb{Z}} \mathbb{Z} / 2\mathbb{Z}$,
where the action of $\mathbb{Z}$ consists in shifting coordinates. We view elements
of $\Gamma$ as pairs $(t,f)$, where $t \in\mathbb{Z}$ and $f$ is a finitely-supported
function from $\mathbb{Z}$ into $\mathbb{Z} / 2\mathbb{Z}$. As a generating set we consider
$\mathcal{G} := \{id, (0,\delta_0), (\pm 1,0) \}$, where $\delta_0$ stands for the Dirac delta at $0$.
The subspaces
$$D_n := \left\langle \sum_{supp(f) \subset \{1,\ldots,n\}} (t,f) \!: t \in \{1,\ldots,n\} \right\rangle$$
satisfy $|D_n| = n$ and $|\partial_{\mathcal{G}} D_n| = 2$. However, every finite subset
$Y \!\subset\! \Gamma$ for which $|\partial_{\mathcal{G}} Y| / |Y| \leq 2/n$ must have
at least $2^{\lambda n}$ points for a certain constant $\lambda > 0$. Indeed, this follows
from that the ball of radius $2n+2$ in $\Gamma$ has more than $2^{n}$ points as an application
of the Saloff-Coste's isoperimetric inequality \cite{coste}: If $Y$ satisfies
$|\partial_{\mathcal{G}} Y| / |Y| \leq 1/n$, then its cardinal is greater than or
equal to a half of the cardinal of a ball of radius $n/2$.
(See \cite{gromov} for an elementary proof of this inequality.)
\label{laurent}
\end{ex}

\index{Isoperimetric profile ! linear profile}
\begin{rem} In the example above, the group $\Gamma$ not only contains torsion elements but is
also amenable. In this direction, let us point out that a nice theorem due to Bartholdi \cite{bartholdi}
establishes that for non-amenable groups, the linear isoperimetric profile cannot behave sublinearly
along any subsequence (see \cite{gromov} for an alternative proof using orderings~!).
\end{rem}\end{small}


\chapter{A PLETHORA OF ORDERS}

\section{Producing New Left-Orders}

\subsection{Convex extensions}
\label{section-convex-extension}

\index{Convex!subgroup}

\hspace{0.45cm} A subset $S$ of a left-ordered group $(\Gamma,\preceq)$ is said to be
{\bf{\em convex}} (with respect to $\precede$) if, for all $f \prec g$ in $S$, every element
$h \!\in\! \Gamma$ satisfying $f \prec h \prec g$ belongs to $S$. If $S$ is a subgroup,
this is equivalent to requiring that $g \in S$ for all $g \in \Gamma$ such that $id \prec g \prec f$
for some $f \in S$.

The family of $\preceq$-convex subgroups is linearly ordered (by inclusion). More precisely,
if $\Gamma_0$ and $\Gamma_1$ are convex (with respect to $\preceq$), then either \esp
$\Gamma_0 \!\subset\! \Gamma_1$ \esp or \esp $\Gamma_1 \!\subset\! \Gamma_0$. Moreover,
the union and the intersection of any family of convex subgroups is a convex subgroup.

\index{Convex!jump}
\begin{small}
\begin{ex} \label{salta}
For each $g \in \Gamma$, it is usual to denote $\Gamma_g$ (resp. $\Gamma^g$)
the largest (resp. smallest) convex subgroup which does not (resp. does) contain $g$.
The inclusion $\Gamma_g \subset \Gamma^g$ is called the {\bf{\em convex jump}}
associated to $g$. In general, $\Gamma_g$ fails to be normal in $\Gamma^g$. Normality
holds for bi-orders, and in this case the quotient $\Gamma^g / \Gamma_g$ is Abelian
(see \S \ref{conrad-general} for the study of left-orders for which this holds for every $g$).
\end{ex}

\index{Order!lexicographic}
\begin{ex} \label{actuar-convexo}
It is not difficult to produce examples of group left-orders without maximal
{\em proper} convex subgroups: consider for instance a lexicographic left-order
on $\mathbb{Z}^{\mathbb{N}}$. Nevertheless, if the underlying group $\Gamma$
is finitely-generated, such a maximal subgroup always exists. Indeed, given a
system of generators $id \prec g_1 \prec \ldots \prec g_k$, let $\Gamma_0$
be the maximal convex subgroup that does not contain $g_k$. Then
$\Gamma_0 \subsetneq \Gamma$. Moreover, if $\Gamma_1$ is a convex subgroup
containing $\Gamma_0$, then, by definition, $g_k \in \Gamma_1$. By convexity,
all the $g_i$'s belong to $\Gamma_1$, hence $\Gamma_1 = \Gamma$. Thus,
$\Gamma_0$ is the maximal proper convex subgroup.
\end{ex}\end{small}

In the dynamical terms of \S \ref{general-3}, convex subgroups are characterized by the following 
proposition.

\vsp

\begin{prop} \label{prop convex subgroup dynamics}
{\em Let $(\Gamma,\preceq)$ be a countable left-ordered group, and let $\Gamma_*$ be a convex subgroup.
Then, in the dynamical realization of $\preceq$, there is a bounded, $\Gamma_*$-invariant interval $I$
with the property that $g(I)\cap I=\emptyset$ for every $g\in\Gamma\setminus \Gamma_*$.

Conversely, let $\Gamma$ be a group acting by orientation-preserving homeomorphisms of the real line
without global fixed points. Suppose that there is an interval $I$ with the property that, for all $g\in \Gamma$,
the intersection $g(I)\cap I$ either is empty or coincides with $I$. Then, in any dynamical-lexicographic order
induced from a sequence $(x_n)$ starting with a point $x_1 \!\in\! I$, the stabilizer $Stab_\Gamma(I)$
is a proper convex subgroup.}
\end{prop}

\noindent{\bf Proof.} Suppose $(\Gamma,\preceq)$ is a countable left-ordered group having $\Gamma_*$ as a
proper convex subgroup. Consider its dynamical realization, and let $a := \inf\{t(h) \mid h\in \Gamma_*\}$ and
$b := \sup\{ t(h) \mid h\in \Gamma_* \}$. By Remark \ref{empieza-en-cero}, one has $t(h) = h(0)$, 
which implies that  $I : = (a,b)$ is a bounded interval fixed by $\Gamma_*$.
Moreover, if $g\in \Gamma$ is such that $g(I)\cap I\neq \emptyset$, then there is $h\in \Gamma_*$ such
that $gh(0) \in I$. Therefore, $t(gh) \in I$. By convexity, this implies that $gh\in \Gamma_*$, 
which yields $g\in \Gamma_*$.

Conversely, suppose that for a $\Gamma$-action on the line, there is a bounded interval $I$
satisfying that $g(I)\cap I$ either is empty or coincides with $I$, for each $g \in \Gamma$.
Let $\preceq$ be a  left-order induced from a sequence $(x_n)$ starting with a point 
$x_1\in I$. If $g \in \Gamma$ satisfies $id\prec g\prec h$ for some $h\in Stab_\Gamma(I)$,
then by definition we have $x_1\leq g(x_1)\leq h(x_1)$. By our hypothesis this implies 
$g(I) = I$. Therefore, $Stab_\Gamma(I)$ is $\preceq$-convex.
$\hfill\square$

\vspace{0.4cm}

\index{Convex!extension}
\noindent{\bf The convex extension procedure.}
Let $\Gamma_*$ be a $\precede$-convex subgroup of $\Gamma$, and
let $\precede_*$ be any left-order on $\Gamma_*$. The {\bf{\em extension of $\precede_*$
by $\precede$}} is the order relation $\precede'$ on $\Gamma$ whose positive cone is
$(P_{\precede}^+ \setminus \Gamma_*) \cup P^+_{\precede_*}$.

One easily checks that $\precede'$ is also a left-invariant total order relation,
and that $\Gamma_*$ remains convex in $\Gamma$ with respect to $\precede'$.
Moreover, the family of $\preceq'$-convex subgroups of $\Gamma$ is formed
by the $\preceq_*$-convex subgroups of $\Gamma_*$ and the $\preceq$-convex
subgroups of $\Gamma$ that contain $\Gamma_*$.

\begin{small}\begin{ex} \label{ex-flipping}
Let $(\Gamma,\preceq)$ be a left-ordered group, and $\Gamma_*$
a $\preceq$-convex subgroup. The extension of (the restriction to $\Gamma_*$ of)
$\overline{\preceq}$ by $\preceq$ will be referred as the left-order obtained
by {\bf {\em flipping}} the convex subgroup $\Gamma_*$. An important case
of this seemingly innocuous construction arises for braid groups;
see the end of \S \ref{fin-gen}.
\end{ex}

\begin{rem}\label{rem-ce}
As we have already pointed out, convex subgroups are not necessarily normal. In the case of a 
normal convex subgroup, the left-order passes to the corresponding quotient. Conversely, if
$\Gamma$ contains a normal subgroup $\Gamma_*$ such that both $\Gamma_*$ and
$\Gamma / \Gamma_*$ are left-orderable, then $\Gamma$ admits a left-order for
which $\Gamma_*$ is convex. Indeed, letting $\preceq_*$ and $\preceq_0$ be
left-orders on $\Gamma_*$ and $\Gamma / \Gamma_*$, respectively, we may
define $\preceq$ on $\Gamma$ by letting $f \prec g$ if either
$f \Gamma_* \prec_0 g \Gamma_*$, or $f \Gamma_* = g \Gamma_*$
and $f^{-1} g$ is $\preceq_*$-positive.

Thus, the extension of a left-orderable group by another left-orderable group is
left-orderable. Using Example \ref{productos}, this implies that the {\bf {\em wreath product}}
\hspace{0.08cm} $\Gamma_1 \wr \Gamma_2 := ( \bigoplus_{\Gamma_2} \Gamma_1 ) \rtimes \Gamma_2$
\hspace{0.08cm} of two left-orderable groups is left-orderable.
\end{rem}\end{small}
\index{Wreath product}

\vsp

In dynamical terms, convex subgroups are relevant because of the next remark.

\vsp 

\begin{small}\begin{rem} \label{ssh} Let $(\Gamma,\preceq)$
be a left-ordered group, and let $\Gamma_*$ be a $\preceq$-convex subgroup.
The space of left cosets $\Omega = \Gamma / \Gamma_*$ carries a natural
total order $\leq$, namely $f \Gamma_* < g \Gamma_*$ if $f h_1 \prec g h_2$
for some $h_1,h_2$ in $\Gamma_*$ (this definition
is independent of the choice of $h_1, h_2$ in $\Gamma_*$). Moreover, the action
of $\Gamma$ by left-translations on $\Omega$ preserves this order. An important
case (to be treated in \S \ref{structure-section}) arises when $\Gamma_*$ is the
maximal proper convex subgroup (whenever it exists); see Example \ref{actuar-convexo}.
\end{rem}\end{small}

The preceding construction allows showing the following very useful proposition.

\begin{prop} \label{rel-convex}
{\em Let $\Gamma$ be a left-orderable group, and
let $\{\Gamma_{\lambda} \!: \lambda \in \Lambda\}$
be a family of subgroups each of which is convex with respect
to a left-order $\preceq_{\lambda}$. Then there exists a left-order on $\Gamma$
for which the subgroup \esp $\bigcap_{\lambda} \Gamma_{\lambda}$ \esp is convex.}
\end{prop}

\vsp

For the proof, we need a lemma that is interesting by itself.

\vsp

\begin{lem} {\em Let $\Gamma$ be a group acting faithfully on a totally ordered space
$(\Omega,\leq)$ by order-preserving transformations. Then for every $\overline{\Omega}
\subset \Omega$, there is a left-order on $\Gamma$ for which the stabilizer of \esp
$\overline{\Omega}$ is a convex subgroup.}
\end{lem}

\noindent{\bf Proof.} Proceed as in \S \ref{general-3} using a well-order $\leq_{wo}$
on $\Omega$ for which $\overline{\Omega}$ is an initial segment. $\hfill\square$

\vsp\vsp\vsp\vsp

\noindent{\bf Proof of Proposition \ref{rel-convex}.} As we saw in Example \ref{ssh}, each
space of cosets $\Gamma / \Gamma_{\lambda}$ inherits a total order $\leq_{\lambda}$ that
is preserved by the left action of $\Gamma$. Fix a well-order $\leqslant_{wo}$ on $\Lambda$,
and let $\Omega := \prod_{\lambda \in \Lambda} ( \Gamma / \Gamma_{\lambda} ) \times \Gamma$ be
endowed with the associate dynamical-lexicographic total order $\leq$. This means that
$([g_{\lambda}],g) \leq ([h_{\lambda}],h)$ if either the smallest (according to
$\leqslant_{wo}$) index $\lambda$ such that $[g_{\lambda}] \neq [h_{\lambda}]$
is such that $[g_{\lambda}] >_{\lambda} [h_{\lambda}]$, or the classes of
$g_{\lambda}$ and $h_{\lambda}$ (with respect to $\Gamma_{\lambda}$) are
equal for every $\lambda$ and $g \preceq h$. The left action of $\Gamma$ on
$\Omega$ is faithful and preserves this order. Since the stabilizer of
$ ([id]_{\lambda})_{\lambda \in \Lambda} \times \Gamma$ coincides with
\esp $\bigcap_{\lambda} \Gamma_{\lambda}$, \esp the proposition
follows from the preceding lemma. $\hfill\square$


\subsection{Free products}
\label{sub-vinogradov}
\index{Vinogradov's theorem}

\hspace{0.45cm} As we have seen in \S \ref{ejemplificando-3}, free groups are bi-orderable.
Actually, a much more general statement involving free products holds. The result below was 
first established by Vinogradov \cite{vinogradov}; see \cite{LM} for a kind translation of the original reference.

\vsp

\begin{thm} \label{vino}
{\em The free product of an arbitrary family of bi-orderable groups is bi-orderable.
Moreover, given bi-orders on each of the free factors, there is a bi-order
on the free product that extends these bi-orders.}
\end{thm}

\vsp

Let us point out that a similar statement holds for left-orderability. However, the proof is much simpler. Indeed, let
$\Gamma = * \Gamma_{\lambda}$ be a free product of
left-orderable groups. Then the direct sum $\oplus_{\lambda} \Gamma_{\lambda}$ is
left-orderable. Moreover, the kernel of the natural homomorphism from $\Gamma$ to
$\oplus_{\lambda} \Gamma_{\lambda}$ is very well known to be a free group 
(see for instance \cite{MKS}). Since free groups are left-orderable, $\Gamma$ itself is left-orderable. (An alternative
--dynamical-- argument is contained in the proof of Theorem \ref{alternative-free-products}.)

The statement concerning bi-orderability is more subtle. For instance, the argument above
does not apply, as the bi-orders in the free kernel are not necessarily invariant under
conjugacy by elements of $\Gamma$. Although we are mostly concerned with left-orders
here, we next reproduce the proof of Theorem \ref{vino} given by Bergman in \cite{bergman-2}, which is close to Vinogradov's original approach. 

\begin{small}
\begin{ejer} Show that the free product of groups with the U.P.P.
has the U.P.P. (See \cite{st} in case of problems.) Show an analogous
statement for groups admitting a locally-invariant order.
\end{ejer}\end{small}

\vsp\vsp\vsp

\noindent{\bf Bi-ordering on groups induces from ordered rings.}  The key idea is to embed a free product of groups inside a multiplicative subgroup 
of the ring of 2-by-2 matrices over a suitable {\bf {\em orderable ring}}.  Here, rings are associative with unity.  We say that a ring $(R,+,\cdot)$ is {\em orderable}  
if it admits a total order $\leq$ such that the underlying Abelian group $(R,+)$ is an ordered group such that the set of positive elements (that is, elements $> 0$) is 
a multiplicative subsemigroup. Note that an orderable ring cannot have zero divisors. In particular, the direct product of two orderable rings is not orderable. 

Basic examples of orderable rings are subrings of the set of real numbers. Another important example  is given by the following proposition.

\begin{prop} \label{prop ejemplo anillo ordenable}{\em The group ring $R[G]$ of a left-orderable group $G$ over an orderable ring 
$R$ is an orderable ring.}
\end{prop}

\noindent{\bf Proof.} Recall that the elements of $R[G]$ are formal finite sums of the form $p=\sum r_i g_i$, where $g_i\in G$  and $r_i\in R$. 
Alternatively, one can model $R[G]$ as the set of finitely-supported functions from $G$ to $R$. 

Fix a left-order $\preceq$ on $G$ and an order 
$\leq$ on $R$. In the latter formalism, we declare $p$ to be positive if the image $r_i$ of the least element $g_i$ (according to $\preceq$) 
in the support of $p$ is positive in $R$ (according to $\leq$). Clearly, this makes $(R[G],+)$ an ordered group. Moreover, 
the fact that the product of two positive elements $p,q$ is still positive easily follows, since the least element in the support of  $p\cdot q$ is 
precisely the multiplication (in $G$) of the least element in the support of $p$ times the least element in the support of $q$ 
(compare the argument at the beginning of \S \ref{sec-upp}).  $\hfill\square$ 

\vsp\vsp

We now turn to rings of matrices. If $R$ is a ring, we denote by $R[t]$ the ring of polynomials with coefficients in $R$.

\begin{prop}\label{prop U_R} {\em For an ordered ring $R$, let $M_2(R)$ be the ring of $2$-by-$2$ matrices with coefficients in $R$. 
Let $U_R\subseteq M_2(R)[t]$ be the set of polynomials whose constant terms are diagonal matrices with entries that are positive in $R$. 
Then $U_R$ is a (multiplicative) semigroup admitting a total order which is invariant under left and right multiplication.}

\end{prop}

\noindent{\bf Proof.} Since the product of two diagonal matrices with positive entries is still a diagonal matrix with positive entries, 
we have that $U_R$ is a multiplicative subsemigroup of $M_2(R)[t]$.  We need to show that $U_R$ is {\em bi-orderable}. 
To do this, for each $n\in \mathbb N$, let us order the $R$-submodule $t^nM_2(R)\subseteq M_2(R)[t]$ by choosing an arbitrary order among 
the four positions in the 2-by-2 matrix coefficient, and declare a non-zero element of $t^nM_2(R)$ to be {\em pre-positive} if in the {\em first} non-zero entry 
is a positive element of $R$. Note that if $a\in t^nM_2(R)$ is pre-positive and $d\in M_2(R)$ is a diagonal matrix with positive entries, then the product of 
$d$ and $ a$ is also a pre-positive element of $t^nM_2(R)$.

Let now $a$ and $b$ be two different elements in $U_R$, and let $n\geq 0$ be the least exponent such that $t^n$ appears with a nonzero 
matrix coefficient in $b-a$. Denote this matrix by  $A\in M_2(R)$, and write  $a\prec b$ if  $A\, t^n$ is a pre-positive element in $t^nM_2(R)$. 
Clearly, $\preceq$ is a total order on $U_R$. Moreover, from the above observation that pre-positivity is preserved under multiplication 
by diagonal matrices with positive entries, we readily obtain that $\preceq$ is invariant under left and right multiplication. 
This endows $U_R$ with a bi-order, as desired. $\hfill\square$ 

\vsp\vsp

We are now in position to prove that free products of bi-orderable groups are bi-orderable in full detail.

\vsp\vsp

\noindent {\bf  Proof of Theorem \ref{vino}.} Let $G$ and $H$ be bi-orderable groups. We first show that their free product $G*H$ is bi-orderable as well. 

Let $R=\Z[G*H]$. Since $G*H$ is left-orderable, by Proposition \ref{prop ejemplo anillo ordenable}, we have that $R$ is an orderable ring. With this, 
Proposition \ref{prop U_R} builds  a bi-orderable semigroup $U_R\subseteq M_2(R)[t]$. We claim that $U_R$ contains a subgroup isomorphic to $G*H$, 
and hence it is bi-orderable. 

To show the claim above, we first note that, because of the natural isomorphism $M_2(R)[t]\simeq M_2(R[t])$, we can consider the inclusion of $G$ inside 
$U_R$ given by 
$$\varphi:g\mapsto \left(\begin{array}{cc} 1 & -t \\ 0 & 1\end{array}\right) \left(\begin{array}{cc} g & 0 \\ 0 & 1\end{array}\right) \left(\begin{array}{cc} 1 & t \\ 0 & 1\end{array}\right) =\left(\begin{array}{cc} g & t(g-1) \\ 0 & 1\end{array}\right).$$
Similarly, we can embed $H$ into $U_R$ via 
$$\psi:h\mapsto \left(\begin{array}{cc} 1 & -t \\ 0 & 1\end{array}\right) \left(\begin{array}{cc} 1 & 0 \\ 0 & h\end{array}\right) \left(\begin{array}{cc} 1 & t \\ 0 & 1\end{array}\right) =\left(\begin{array}{cc} 1 & 0 \\ t(h-1) & h\end{array}\right).$$
To conclude, observe that the images of $G$ and $H$ generate a free  product. Indeed, the image of any (reduced) word $g_1h_1\ldots g_nh_n\in G*H$ is nontrivial in 
$U_R$, since the matrix $\varphi(g_1)\psi(h_1)\ldots\varphi(g_n)\psi(h_n)$ moves the vector $(1,1)$ (this is straightforward to check). 

Arguing by induction, one obtains that the free product of any finitely many bi-orderable groups is bi-orderable as well. Further, since bi-orderability is a 
local property (see \S \ref{general-2}), one concludes that 
any arbitrary free product of bi-orderable groups is bi-orderable.

To close the proof, it remains to show that any free product $\Gamma= *\Gamma_{\lambda}$ of bi-ordered groups admits a bi-order 
that extends the orders of the factors. To do this, we use again the exact sequence
$$\{id\} \to G \to  *\Gamma_{\lambda}\to  \oplus \Gamma_{\lambda}\to \{id\},$$
for which $G$ is a free group. By the first part of the proof, there is a bi-order on $G$ that is invariant under conjugation 
by all elements of $\Gamma$.  Using the convex extension procedure (see Remark \ref{rem-ce}), we can 
thus build an order on $\Gamma$  that  extends the bi-orders on the factors and is bi-invariant. $\hfill\square$


\subsection{Left-orders from bi-orders}

\hspace{0.45cm} As we already pointed out, a left-orderable group all
of whose left-orders are bi-invariant is necessarily Abelian \cite{DGR}. This
suggests the existence of natural procedures to create left-orders starting
with bi-orders on groups. Here we briefly discuss two of them.

\vsp\vsp\vsp

\noindent{\bf Left-orders from the sequence of convex subgroups.} Let $\{ \Gamma_i \}$
be the family of convex subgroups for a bi-order $\preceq$ on a group $\Gamma$.
Since $\preceq$ is bi-invariant, for every $g \in \Gamma$, each subset of
the form $g \Gamma_i g^{-1}$ is also convex. Given any
well-order $\leq_{wo}$ on the set of indices $i$, we may define a left-order
$\preceq'$ on $\Gamma$ as follows: Given $g \in \Gamma$, we look for the minimal (with
respect to $\leq_{wo}$) index $i$ such that $g \Gamma_i g^{-1} \neq \Gamma_i$, and we let
$j(i)$ so that $g \Gamma_i g^{-1} = \Gamma_{j(i)}$. If $j(i) >_{wo} i$ (resp.
$j(i) <_{wo} i$), then we let $g \succ' id$ (resp. $g \prec' id$); if $g$ fixes
each $\Gamma_i$, then we let $g \succ' id$ if and only if $g \succ id$.

One easily checks that $\preceq'$ is well-defined and left-invariant.
Note that $\preceq'$ coincides with the original bi-order $\preceq$ if every convex subgroup is normal.

\begin{small} \begin{ejem} Let us consider the bi-order $\preceq_{x^+}^-$ on Thompson's
group $\mathrm{F}$ (see \S \ref{ejemplificando-4}). Let $(x_i)$ be a numbering of all dyadic,
rational numbers of $]0,1[$. Each $x_i$ gives raise to a convex subgroup $\Gamma_i$ formed
by the elements $g$ such that $g(x) \!=\! x$ for all $x \!\in\! [x_i,1]$. Although there are
more convex subgroups than these, this family is invariant under the conjugacy action.
By performing the construction above, we get the left-order $\preceq$ on $\mathrm{F}$
for which $f \succ id$ if and only if $f (x_i) > x_i$ holds for the smallest integer
$i$ such that $f(x_i)  \neq x_i$. (Compare \S \ref{general-3}.)
\end{ejem} \end{small}

\vsp\vsp

\noindent{\bf Combing elements with trivial conjugacy action on a certain left-order.}
Proposition \ref{prop-combing} below appears in \cite{LRR}, yet it was already implicit in \cite{DGR}.

\vsp

\begin{lem}\label{lema combing} {\em Suppose $\preceq$ is a left-order on a group $\Gamma$
admitting a normal, convex subgroup $\Gamma_*$, and let $g \in \Gamma \setminus \Gamma_*$.
If conjugation by $g$ preserves the left-order on $\Gamma_*$ (that is, \esp
$g ( P_{\preceq}^+ \cap \Gamma_*) g^{-1} = P_{\preceq}^+ \cap \Gamma_*$), then
there exists a left-order on the subgroup $\langle g,\Gamma_* \rangle$ that has
$g$ as minimal positive element and coincides with $\preceq$ on $\Gamma_*$.}
\end{lem}

\noindent{\bf Proof.} Since $\Gamma_*$ is normal in $\Gamma$,
every element in $\langle g, \Gamma_* \rangle$ may be written in a
unique way in the form $g^n h$, with $n \!\in \mathbb{Z}$ and $h \in \Gamma_*$.
Define $\preceq_*$ on $\langle g,\Gamma_* \rangle$ by letting $g^n h \succeq_* id$ if
and only if either $h \in P_{\preceq}^+$ or $h = id$ and $n > 0$. Invariance of $\preceq$
under conjugation by $g$ shows that this is a well-defined left-order on
$\langle g,\Gamma^*\rangle$. That $\preceq_*$ coincides
with $\preceq$ on $\Gamma_*$ follows from the definition. Finally, the fact that $g$ is the
minimal positive element of $\preceq_*$ also follows from the definition. $\hfill\square$

\vsp\vsp\vsp

Combined with the convex extension technique, this lemma
allows us to produce many interesting left-orders. Invoking
Example \ref{salta}, this is summarized in the next proposition.

\vsp

\begin{prop} \label{prop-combing} {\em Let $(\Gamma,\preceq)$ be a bi-ordered group,
and let $\Gamma_g \!\subset\! \Gamma^g$ be the convex jump associated to an element
$g \in \Gamma$. Assume that the quotient $\Gamma^g / \langle g, \Gamma_g \rangle$
is torsion-free. Then there exists a left-order $\preceq'$ on $\Gamma$ having
$g$ as minimal positive element and such that $\preceq'$ coincides with $\preceq$ on $\Gamma_g$.}
\end{prop}

\noindent{\bf Proof.} First note that both $\Gamma_g$ and $\Gamma^g$ are invariant
under conjugation by $g$. As $\preceq$ is bi-invariant, conjugacy by $g$ preserves the positive
cone of $\Gamma_g$. Thus, we are under the hypothesis of the preceding lemma, which allows
to produce a left-order on $\langle g, \Gamma_g \rangle$ having $g$ as minimal positive element.
This left-order may be extended to a left-order $\preceq_*$ on $\Gamma^g$, as $\Gamma^g / \langle g,
\Gamma_g \rangle$ is assumed to be torsion-free (recall that $\Gamma^g/\Gamma_g$ is Abelian; see
Example \ref{salta}). Finally, we let $\preceq'$ be the extension of $\preceq_*$ by $\preceq$. $\hfill\square$

\vsp\vsp

\begin{small}
\begin{ex} Given an element $g$ in the free group $\Gamma := \mathbb{F}_n$, let
$k = k(g) \in \mathbb{N}$ be such that $g \in \Gamma_{k} \setminus \Gamma_{k+1}$,
where $\Gamma_i$ denotes the $i^{\mathrm{th}}$-term of the lower central series.
If $g \Gamma_{k}$ has no nontrivial root in $\Gamma_{k} / \Gamma_{k+1}$, then
we are under the hypothesis of Proposition \ref{prop-combing} for any bi-order
on $\mathbb{F}_n$ obtained from the series $\Gamma_i$. Thus, $g$ appears
as the minimal positive element for a left-order on $\mathbb{F}_n$.
\end{ex}

\begin{ex} \label{ex corona laminar action}  Let $\Z\wr\Z:=\bigoplus_\Z \Z\rtimes \Z$ be the wreath product of $\Z$ 
with itself. Recall that the conjugation action of $\Z$  on $\bigoplus_\Z \Z$ is  by shifting the indexes. 
Let $H :=\bigoplus_\Z \Z$. We saw at the end of \S \ref{ejemplificando-2} that we can  use the exact sequence
$$\{id \} \longrightarrow H \longrightarrow \Z\wr\Z\longrightarrow \Z \longrightarrow \{id \},$$
to produce a bi-order $\preceq'$ on $\Z \wr \Z$ as the convex extension of the lexicographic order on $H$ by
one of the two possible orders on the cyclic factor $\Z = \langle a \rangle$. Note that given any element $h\in H$, we can find
two elements $h_\pm$ in $H$ such that $h_-\prec' h^n \prec' h_+$ holds for all $n \!\in\! \Z$. Moreover, $H$ is the maximal proper
$\preceq'$-convex subgroup. Thus, using Lemma \ref{lema combing}, we can produce a left-order $\preceq$ on $\Z \wr \Z$
that has $a$ as its smallest positive element and that coincides with $\preceq'$ on $H$. 

\end{ex}

\begin{ejer} For the items below, we refer to the notations from the preceding Example \ref{ex corona laminar action}. 

\noindent (i) Show that in the dynamical realization of $\preceq$, every element of $H$ acts with fixed points. 

\noindent \underbar{Hint.} Use the following fact already stressed above: for every $h\in H$, there are
 elements $h_\pm$ in $H$ such that $h_-\prec' h^n \prec' h_+$ holds for all $n \!\in\! \Z$. 

\noindent (ii) Show that for any element $g\in \Z\wr\Z$, there is $h\in H$ such that $g\preceq h$. 
Deduce from this that in the dynamical realization of $\preceq$,  the subgroup $H$ has no global fixed points. 

\noindent (iii) Show that the dynamical realization of $\preceq$ is semiconjugate to 
Plante's action from Example \ref{ejemplo de plante}. 
\end{ejer}
\end{small}


\section{The Space of Left-Orders}
\label{space-left-orders}

\index{Order!isolated}
\index{Space of! left-orders}
\hspace{0.45cm} Following Ghys \cite{ghys-ord} and Sikora \cite{sikora}, given
a left-orderable group $\Gamma$, we denote by $\mathcal{LO}(\Gamma)$ the
set of all left-orders on $\Gamma$. This {\bf{\em space of left-orders}} carries
a natural (Hausdorff and totally disconnected) topology whose sub-basis is the
family of sets of the form \esp $U_{f,g} \!=\! \{\preceq : \esp f \!\prec\! g\}$.
Due to left-invariance, another sub-basis is the family of sets \esp
$V_f = \{ \preceq : \esp id \prec f \}$.  In particular, a left-order $\preceq$ is {\bf{\em isolated}} if 
there is a finite set $S\subset \Gamma$ such that $\preceq$ is the only left-order satisfying 
that $id\preceq f$ for all $f\in S$. (Some authors call such an order {\em finitely determined}; 
see for instance \cite{SM85}.) 

To better understand the topology on $\mathcal{LO}(\Gamma)$, one may proceed as in
\S \ref{general-2} by identifying left-orders on $\Gamma$ to certain points in \esp
$\{-1,+1\}^{\Gamma \setminus \{id\}}$. Nevertheless, to later cover also the case of partial
left-orders (see Exercise \ref{space-of-LIOs}), it is better to model $\mathcal{LO}(\Gamma)$ as a subset of \esp
$\{-1,+1\}^{\Gamma \times \Gamma \setminus \Delta}$, \esp namely the one
formed by the functions \esp $\varphi$ \esp satisfying:

\vsp

\noindent -- ({\em Reflexivity}) \esp \esp $\varphi (g,h) = +1$ if and only if $\varphi (h,g) = -1$;

\vsp

\noindent -- ({\em Transitivity}) \esp \esp if $\varphi (f,g) = \varphi (g,h) = +1$, then
$\varphi (f,h) = +1$;

\vsp

\noindent -- ({\em Left-invariance}) \esp \esp $\varphi(fg,fh) = \varphi(g,h)$
for all $f$ and $g \neq h$ in $\Gamma$.

\vsp

\noindent Indeed, every left-order $\preceq$ on
$\Gamma$ leads to such a function $\varphi_{\preceq}$, namely
$\varphi_{\preceq} (g,h) = +1$ if and only if $g \succ h$. Conversely, every $\varphi$
with the above properties induces a left-order $\preceq_{\varphi}$ on $\Gamma$,
namely $g \succ_{\varphi} h$ if and only if $\varphi(g,h) = +1$. Now, if we endow
$\{-1,+1\}^{\Gamma \times \Gamma \setminus \Delta}$ with the product topology
and the subset above with the subspace one, then the induced topology on
$\mathcal{LO}(\Gamma)$ coincides with the one previously defined by
prescribing the sub-basis elements. As a consequence, since
$\{-1,+1\}^{\Gamma \times \Gamma \setminus \Delta}$ is a compact space and
the subspace above is closed, the topological space $\mathcal{LO}(\Gamma)$
is compact.

\vsp

As an example regarding the convex extension procedure (see \S \ref{section-convex-extension}), 
the reader should easily be able to show the next proposition.

\vspace{0.01cm}

\begin{prop} \label{prop convex and LO}{\em Let $\preceq$ be a left-order on $\Gamma$ and
$\Gamma_*$ a $\preceq$-convex subgroup. Then the map from $\mathcal{LO}(\Gamma_*)$
into $\mathcal{LO}(\Gamma)$ that sends $\preceq_*$ to its convex extension by $\preceq$
is a continuous injection. If, in addition, $\Gamma_*$ is normal, then
there is a continuous injection from $\mathcal{LO}(\Gamma_*) \times
\mathcal{LO}(\Gamma / \Gamma_*)$ into $\mathcal{LO}(\Gamma)$
having $\preceq$ in its image.}
\end{prop}

\vspace{0.01cm}

\begin{small}\begin{ex} \label{no-son-todos}
The subspace of dynamical-lexicographic left-orders on
$\mathrm{Homeo}_+(\mathbb{R})$ (see \S \ref{general-3}) is not
closed inside $\mathcal{LO} (\mathrm{Homeo}_+(\mathbb{R}))$. 
To show this, let $(y_k)$ be a dense sequence of real
numbers, and let $(x_n)$ be a monotone sequence converging to a point
$x \!\in\! \mathbb{R}$. For each $n$, define a sequence $(y_{n,k})_k$ by
$y_{n,1} = x_n$ and $y_{n,k} = y_{k-1}$ for $k \!>\! 1$. This gives raise
to a sequence of left-orders $\preceq_n$ (the sign of each point $y_{n,k}$ is
chosen to be $+$). Passing to a subsequence if necessary, we may assume that
$\preceq_n$ converges to a left-order $\preceq$ on $\mathrm{Homeo}_+(\mathbb{R})$.
We claim that $\preceq$ is not an order of  dynamical-lexicographic type. 
Indeed, let $\preceq'$ be an arbitrary  left-order, and let $x'$ be
the first point different from $x$ for the well-order leading to $\preceq'$
(thus, $x'$ may be the first or the second term of this well-order). Let
$f \!\in\! \mathrm{Homeo}_+(\mathbb{R})$ be such that $f(x) = x$ and
$f(y) > y$ for all $y \neq x$. Let $g,h$ be elements in
$\mathrm{Homeo}_+(\mathbb{R})$ that coincide with $f$ in a neighborhood
of $x$ and $g (x') > x' > h (x')$. By definition, the signs of $g,h$ with
respect to $\preceq'$ are different. However, since
$g (y_{n,1}) = g (x_n) > x_n  = y_{n,1}$
and $h (y_{n,1}) > y_{n,1}$ for all $n$ sufficiently large, both $g,h$ are
$\preceq_n$-positive. Passing to limits, both $g,h$ become
$\preceq$-positive, thus showing that $\preceq$ cannot
coincide with $\preceq'$.
\end{ex}
\index{Order!dynamically-lexicographic}

\vspace{0.1cm}

\begin{ex} \label{example-Mul} 
In an earlier version of this book, we asked whether the set of dynamical-lexicographic left-orders on
$\mathrm{Homeo}_+(\mathbb{R})$ is dense in the corresponding space of left-orders. This was recently answered in the negative  
by Muliarchyk in \cite{Muli}, who gave the following brilliant example. Consider the homeomorphisms of the real line 
$f_1,f_2$ defined by $f_1 (x) = x+1$ and 
\begin{equation}
f_2(x) = \left\lbrace
\begin{array}{ll}
x+1 & \mbox{if } x < 1,\\
2x & \mbox{otherwise}.
\end{array}
\right.
\end{equation}
We claim that for all elements $g_1,\ldots,g_n$ in $\mathrm{Homeo}_+(\mathbb{R})$, there exists an order on 
$\langle f_1,f_2,g_1,\ldots,g_n \rangle$ for which $f_1$ is positive but $f_2$ is negative. 
Assuming this, by Exercise \ref{ejer-finitos-forced}, 
there exists an order $\preceq$ on the whole group $\mathrm{Homeo}_+(\mathbb{R})$ for which $f_1$ is positive but $f_2$ 
is negative. However, such an order cannot be approached by dynamical-lexicographic orders, since for any such order,  
obviously the elements $f_1$ and $f_2$ either are both positive or both negative.

To prove the claim, we fix an order on the group of germs at infinity similar to those built on $\ger$ in 
Remark \ref{rem-germs} (alternatively, conjugate the groups of germs via the map $x \mapsto 1/x$). We restrict this 
order to (a perhaps partial order on) $\langle f_1,f_2,g_1,\ldots,g_n \rangle$. If this order is not total, we  extend it 
arbitrarily to a total order (via a convex extension procedure). We thus get an order $\preceq'$ for which both $f_1$ and $f_2$ are positive 
and $f_1^k \prec' f_2$ for every $k \in \mathbb{Z}$. Consider the dynamical realization of this order on the line. The latter condition 
translates into the following: for the point $q := \sup \{ f_1^k (t(0)) : k \in \mathbb{Z} \}$, one has $f_1 (q) = q$ and $f_2 (t(0)) > q$. 
Consider now any dense sequence $(x_n)$ on the line such that $x_1 := q$ and $x_2 := t(0)$, and let $\preceq$ be the 
associated dynamical-lexicographic order on $\langle f_1,f_2,g_1,\ldots,g_n \rangle$ built from this sequence with all 
signs positive except for the first one. One easily checks that $f_1 \succ id$ but $f_2 \prec id$, as claimed.
\label{dense-homeo}
\end{ex}
\end{small}

\vspace{0.1cm}

If $\Gamma$ is a countable left-orderable group, then the natural topology of
$\mathcal{LO}(\Gamma)$ is metrizable. Indeed, if $\mathcal{G}_0 \subset \mathcal{G}_1
\subset \ldots$ is a complete exhaustion of $\Gamma$ by finite sets, then we can define
the distance between two different left-orders \esp $\leq$ \esp and \esp $\preceq$ \esp by
letting \esp $d(\leq,\preceq) = 2^{-n}$, \esp where $n$ is the maximum non-negative
integer such that $\leq$ and $\preceq$ coincide on $\mathcal{G}_n$. An equivalent
metric \esp $d'$ \esp is obtained by letting \esp $d' (\leq,\preceq) = 2^{-n'}$,
\esp where $n'$ is the maximum non-negative integer such that the positive cones
of $\leq$ and $\preceq$ coincide on $\mathcal{G}_{n'}$, that is, \esp
$P_{\leq}^+ \cap \mathcal{G}_{n'} = P_{\preceq}^+ \cap \mathcal{G}_{n'}$. \esp
One easily checks that these metrics are ultrametric. Moreover, the fact that
$\mathcal{LO}(\Gamma)$ is compact becomes more transparent in this case,
as it follows from a Cantor diagonal type argument.

\vspace{0.1cm}

\index{Word length}
When $\Gamma$ is finitely-generated, it is natural to choose $\mathcal{G}_n$ as being the {\em ball
of radius n} (centered at $id$) with respect to some finite, symmetric system of generators $\mathcal{G}$
of $\Gamma$. Such a ball is  usually denoted by $B_n(id)$, or simply as $B_n$.  Recall that by 
{\em symmetric} we mean that $g^{-1} \!\in\! \mathcal{G}$ for all
$g \in \mathcal{G}$, and the ball $B_n$ is the set of elements having word-length
at most $n$, where the {\em word-length} $\|g\|$ of $g \in \Gamma$ is the minimum
$m$ for which $g$ which can be written in the form \esp $g = g_{i_1} g_{i_2} \cdots g_{i_m}$,
\esp with $g_{i_j} \in \mathcal{G}$. 

\begin{small}\begin{rem}\label{rem-HDQ}
The choice $2^{-n}$ for the distances above could be replaced by any other decreasing sequence 
of positive numbers converging to $0$ (for example, $1/n$ would also work). However, for 
finitely-generated groups, an exponential choice seems to be the right one for several reasons. 
One is treated in Exercise \ref{ejer:este-si} below. Another one comes for the easy-to-check fact 
that, for such a choice, the metrics on $\mathcal{LO}(\Gamma)$
resulting from two different finite systems of generators are not only topologically but also H\"older 
equivalent. As a consequence, although the Hausdorff dimension of the space of left-orders 
can change when varying the system of generators, whether its value is zero, positive, or 
infinite makes sense independently of the system. We will come back to this 
interesting issue for the case of the free group in \S \ref{case-free}
\end{rem}
\end{small}

\vsp

\index{Space of bi-orders}
\begin{small}\begin{ejer} \label{bi-inva}
Given a bi-orderable group $\Gamma$, denote by $\mathcal{BO}(\Gamma)$ the
{\bf{\em space of bi-orders}} of $\Gamma$. Show that $\mathcal{BO}(\Gamma)$
is closed inside $\mathcal{LO}(\Gamma)$, hence compact.
\end{ejer}

\index{Space of!locally invariant orders}
\begin{ejer} \label{space-of-LIOs}
Given a group $\Gamma$ admitting a locally-invariant order (see \S \ref{LOG}), denote by
$\mathcal{LIO}(\Gamma)$ the set of all locally-invariant orders on $\Gamma$. Consider
the topology on $\mathcal{LIO}(\Gamma)$ having as a sub-basis the family of sets
\esp $U_{f,g} \!=\! \{\preceq : \esp f \!\prec\! g\}$. Show that, endowed with
this topology, $\mathcal{LIO}(\Gamma)$
is compact. Conclude that a group $\Gamma$ admits a locally-invariant order
if and only if each of its finitely-generated subgroups admits such an order.
(Compare \cite[Theorem 2.4]{chiswell}.)

\noindent{\underline{Hint.}} As a model of $\mathcal{LIO}(\Gamma)$ consider the subset
of $\{-1,0,+1\}^{\Gamma \times \Gamma \setminus \Delta}$ formed by the functions
$\varphi$ such that $\varphi (g,h) = +1$ if and only if $\varphi(h,g) = -1$, and
such that for every $g \neq id$ and $h \in \Gamma$ one has either $\varphi(hg,h) = +1$
or $\varphi(hg^{-1},h) = +1$. (Two elements $g,h$ that are incomparable for a locally
invariant order will then satisfy $\varphi_{\preceq} (g,h) = 0$...)
\end{ejer}

\begin{ejer} \label{diffuse-implies-LIO}
Complete the proof of Proposition \ref{trivia} by showing that every weakly diffuse group 
admits a locally-invariant order. (See \cite{LIOs-on-amenable} in case of problems.)

\noindent{\underline{Hint.}} By a compactness type argument, it is enough to
show the following: For each finite subset $A$ of $\Gamma$, there exists a partial
order $\preceq$ such that for all $f \in A$ and each nontrivial element $g \in \Gamma$
such that both $fg$ and $fg^{-1}$ lie in $A$, either $fg \succ f$ or $fg^{-1} \succ f$.
To construct such a $\preceq$, proceed by induction, the case where $A$ is a
single element being evident. Now, given an arbitrary $A$, by the weakly diffuse
property there is $h \in A$ such that for each nontrivial element $g \in \Gamma$,
either $h g \notin A$ or $hg^{-1} \notin A$. By the induction hypothesis,
$A \setminus \{ h \}$ admits an order as requested. Extend $\preceq$ to
all $A$ by declaring $h$ to be larger than all other elements.
\end{ejer}\end{small}

The group $\Gamma$ (continuously) acts on
$\mathcal{LO}(\Gamma)$ by conjugacy (equivalently, by right multiplication): given an order
$\preceq$ with positive cone $P^+$ and an element $f \!\in\! \Gamma$, the image of $\preceq$
under $f$ is the order $\preceq_{f}$ whose positive cone is \esp $f \esp P^+ f^{-1}$.
\esp In other words, one has \esp $g \preceq_f h$ \esp if and only if \esp $f^{-1} g f
\preceq f^{-1} h f$, \esp which is equivalent to \esp $g f \preceq h f$.
Also note that the map sending $\preceq$ to $\overline{\preceq}$ from Example
\ref{la-primera} is a continuous involution of $\mathcal{LO}(\Gamma)$.

\begin{small}
\begin{ex} If a group left-order is obtained via an action on a totally ordered space $\Omega$, then the
conjugacy action corresponds to changing the order of the comparison points. More precisely, in the notation
of \S \ref{general-3}, if $\preceq$ comes from a well-order $\leq_{wo}$ on $\Omega$, then $\preceq_f$ is
obtained from the same action using the well-order $f_* (\leq_{wo})$ given by $\omega_1 \hspace{0.14cm}
f_*(\leq_{wo}) \hspace{0.14cm} \omega_2$ whenever $f (\omega_1) \leq_{wo} f (\omega_2)$. In particular,
for countable subgroups of $\mathrm{Homeo}_+(\mathbb{R})$, if $\preceq$ is induced from a dense
sequence $(x_n)$ in $\mathbb{R}$, then $\preceq_f$ is induced from the sequence $(f(x_n))$.
\end{ex}

\begin{rem} \label{action-outer}
If $\Gamma$ is a left-orderable group, then the whole group of automorphisms of $\Gamma$
(and not only the group of inner automorphisms) acts on $\mathcal{LO}(\Gamma)$. This
is useful to study bi-orderable groups. Indeed, since the fixed points for the right
action of $\Gamma$ on $\mathcal{LO}(\Gamma)$ correspond to the bi-invariant left-orders,
the group $Out(\Gamma)$ of outer automorphisms of $\Gamma$ acts on the corresponding
space of bi-orders $\mathcal{BO}(\Gamma)$. The reader is referred to
\cite{koberda} for some applications of this idea to the case of free groups.
\end{rem}
\end{small}

In general, the study of the dynamics of the action of $\Gamma$ on $\mathcal{LO}(\Gamma)$
should reveal useful information. Very simple questions on this were already formulated in \cite{order}.

\index{Action!minimal}
\begin{question} For which finitely-generated, left-orderable groups having an infinite space of left-orders
is the action of $\Gamma$ on $\mathcal{LO}(\Gamma)$ uniformly equicontinuous or distal~?
The same question makes sense for {\em minimality}\footnote{Recall that an action is
said to be {\bf {\em minimal}} if every orbit is dense.}, or for having a dense orbit (the latter
is the case of free groups, as we will see in \S \ref{case-free}; for the former,
we do not know any example).
\end{question}

\begin{small}
\begin{ejer}
Give an example of a countable group $\Gamma$ whose action on $\mathcal{LO}(\Gamma)$ 
is minimal. (See \cite{clay-lattice} in case of problems with this.)
\end{ejer}

\begin{ejer} \label{ejer:este-si}
Show that, if the space of left-orders of a finitely-generated group $\Gamma$ is endowed 
with the natural metric (see Remark \ref{rem-HDQ}), then the action of $\Gamma$ on $\mathcal{LO} (\Gamma)$ is by 
bi-Lipschitz homeomorphisms.
\end{ejer}\end{small}


\subsection{Finitely many or uncountably many left-orders}
\label{fin-uncount}

\hspace{0.45cm} We now state the first nontrivial general theorem concerning
the space of left-orders of a left-orderable group. This result was first
obtained by Linnell \cite{linnell-new} by
elaborating on previous ideas of Smirnov, Tararin, and Zenkov. Let us
point that no analogue for spaces of bi-orders holds \cite{butts};
see however \S \ref{section-cristobal}.

\vspace{0.1cm}

\begin{thm} {\em If the space of left-orders of a
left-orderable group is infinite, then it is uncountable.}
\label{linnell-general}
\end{thm}

\vspace{0.1cm}

The starting point to show this result is the following. Let $\Gamma$ be a left-orderable
group and $M$ a {\bf{\em minimal subset}} of $\mathcal{LO}(\Gamma)$; that is, a nonempty,
closed subset that is invariant under the conjugacy action of $\Gamma$ and does not properly
contain any nonempty, closed, invariant set. Since the set $M'$ of accumulation points of
$M$ is both closed and invariant, we must have either $M'=M$ or $M'=\emptyset$. In other words,
either $M$ has no isolated points, or it is finite. In the former case, a well-known result
in General Topology asserts that $M$ must be uncountable (see \cite[Theorem 2-80]{HY}).
In the latter case, the stabilizer of any point $\preceq$ of $M$ is a finite-index
subgroup of $\Gamma$, restricted to which $\preceq$ is bi-invariant. Theorem
\ref{linnell-general} then follows from the following proposition.

\vsp

\begin{prop} \label{linda}
{\em Let $(\Gamma,\preceq)$ be a left-ordered group containing a finite-index
subgroup $\Gamma_0$ restricted to which $\preceq$ is bi-invariant. If $\preceq$ has a
neighborhood in $\mathcal{LO}(\Gamma)$ containing only countably many left-orders, then
$\mathcal{LO}(\Gamma)$ is finite.}
\end{prop}

\vsp

The proof of this proposition uses results and techniques from the theory of
Conradian orders. Hence, we postpone the (end of the) proof of Theorem
\ref{linnell-general} to \S \ref{section-cristobal}.

\vsp

The argument above distinguishes conjugate left-orders, even though one would like to consider them as being 
``equal'' (for instance, they share all dynamical properties). This leads to the couple of natural questions below that 
are contained in \cite{order} and reproduced in an earlier version of this book. Both have been recently answered.

The first question is whether 
for a finitely-generated, left-orderable group $\Gamma$, the space of orbits $\mathcal{LO}(\Gamma) / \Gamma$ can be a 
{\bf {\em non-standard Borelian space}} ({\em i.e.}, a space which is not measurable isomorphic to $[0,1]$). It turns out that the 
answer is affirmative in many cases, as it was proved by Calderoni and Clay in \cite{CC1} (see \cite{CC2,CC3} for further 
developments and examples, as well as \cite{molina} for the particular case of nilpotent groups).

The second question is whether the set of isolated left-orders of a left-orderable group $\Gamma$ is always finite modulo the 
conjugacy action of $\Gamma$. Here, the answer turns out to be negative. This is the case for instance of the free abelian product 
$\mathbb{Z} \times \mathbb{F}_{2n}$ for every integer $n \geq 1$ and the braid group $\mathbb{B}_3$, as it follows from 
the tools and methods from \cite{malicet-mann-rivas-triestino,mann-rivas} and \cite{matsu-circ-1}, respectively. 

Despite the recent results above, the next question taken from \cite{mann-rivas} remains open. 

\vsp

\index{Order!isolated}
\begin{question}
Does there exist a left-orderable group $\Gamma$ with a left-order which neither 
is isolated nor belongs to a Cantor subset of $\mathcal{LO} (\Gamma)$~?
\end{question}

\vsp

\index{Group!Tararin}
\index{Subnormal series}
\index{Normal series}
In the rest of this section, we give a beautiful characterization (due to Tararin
\cite{tararin-paper}) of groups having finitely many left-orders.
(We will refer to them as {\bf {\em Tararin groups.}}) To do this, recall that a
{\bf{\em rational series}}  \index{Rational series} for a group $\Gamma$ is a finite sequence of subgroups
\begin{equation}
\label{rat-ser}
\{id\} = \Gamma^k \lhd \Gamma^{k-1} \lhd \ldots \lhd \Gamma^0 = \Gamma
\end{equation}
that is {\bf{\em subnormal}} (that is, each $\Gamma^{i}$ is normal in $\Gamma^{i-1}$,
but not necessarily in $\Gamma$), and such that each quotient $\Gamma^{i-1} / \Gamma^i$
is torsion-free rank-1 Abelian. Such a series is said to be {\bf{\em normal}} if
each $\Gamma^i$ is normal in $\Gamma$. Note that a repeated application of the 
convex extension procedure shows that every group admitting 
a rational series is left-orderable (see Remark \ref{rem-ce}).

\vspace{0.1cm}

\index{Rational series}
\begin{thm} \label{tararin-theorem}
{\em Every left-orderable group having only finitely many left-orders admits a unique
rational series}
$$\{id\} = \Gamma^k \lhd \Gamma^{k-1} \lhd \ldots \lhd \Gamma^0 = \Gamma.$$
{\em This series is normal and no quotient $\Gamma^{i-2} / \Gamma^i$ is bi-orderable.
Conversely, if a group $\Gamma$ admits a normal rational series such that
no quotient $\Gamma^{i-2} / \Gamma^i$ is bi-orderable, then ($\Gamma$ is left-orderable
and) its space of left-orders $\mathcal{LO}(\Gamma)$ is finite. In this situation, 
for every left-order on $\Gamma$, the convex subgroups are exactly
$\Gamma^0, \Gamma^1, \ldots, \Gamma^k$, the number of left-orders on $\Gamma$
is $2^k$, and each left-order is uniquely determined by the sequence of signs of any
family of elements \esp $g_i \!\in\! \Gamma^{i-1} \setminus \Gamma^{i}$.}
\end{thm}

\vsp

\index{Group!Klein-bottle}
\begin{small}\begin{ejem}
\label{ups-ups} The Klein bottle group $K_2 = \langle a,b \!: a b a^{-1}  = b^{-1} \rangle$
admits exactly four left-orders whose positive cones are $\langle a,b \rangle^+$,
$\langle a,b^{-1} \rangle^+$, $\langle a^{-1},b \rangle^+$, and
$\langle a^{-1},b^{-1} \rangle^+$, respectively (see \S \ref{fin-gen} for details).
The associate rational series is \esp $\{ id \} \lhd \langle b \rangle \lhd K_2$.
\esp More generally, let us consider the group
$$K_k = \langle a_1,\ldots,a_k \!: a_{i+1} a_i a_{i+1}^{-1} = a_i^{-1},
\esp a_i a_j = a_j a_i \mbox{ for } |i-j| \geq 2 \rangle.$$
One can easily check (either using Theorem \ref{tararin-theorem} above or by a direct computation)
that $K_k$ admits $2^k$ left-orders, each of which is determined by the signs
of the $a_i$'s. The corresponding rational series is
$$\{id\} \lhd \langle a_1 \rangle \lhd \langle a_1,a_2 \rangle
\lhd \ldots \lhd \langle a_1,a_2,\ldots,a_k \rangle.$$
\end{ejem}

\begin{ejem} \label{ex-prep}
A dynamical counterpart of having finitely many left-orders for a group is that, up to semiconjugacy,
there may arise only a few actions on the real line. For the case of the group $K_2$ above, this 
translates into the two items below. (See Example \ref{ex-rivas} for an application.)

\vsp

\noindent{\underline{Claim (i).}}  Suppose $K_2 =  \langle a,b \!: a b a^{-1} = b^{-1} \rangle$ acts on the real
line and there is $x\in \mathbb{R}$ such that $x \leq a(x)$. Then $b$ has a fixed point in $I=[x,a(x)]$.

\vsp

Otherwise, changing $b$ by its inverse if necessary, we may assume that $b(z) \!>\! z$ for all $z\in I$. In
particular, $a (x) < b a (x)$, hence $x < a^{-1} b a (x) = b^{-1}(x)$. Therefore, $b(x) < x$, a contradiction.

As a consequence, every open interval $I$ fixed by $a$ on which $a$ acts freely is also fixed by $b$.
Moreover, $b$ has infinitely many fixed points in $I$.

\vsp

\noindent{\underline{Claim (ii).}}  For every open interval $J$ fixed by $b$ and containing no fixed point
of $b$ inside, we have $a(J) \cap J=\emptyset$.

\vsp

Indeed, as $\langle b\rangle$ is normal in $K$, we have that $a(J)\cap J$ is either $J$ or empty. But the first
possibility cannot occur, since in that case $b$ would have fixed points in $J$, due to Claim (i).
\end{ejem}
\end{small}

\vsp\vsp\vsp

The proof of Theorem \ref{tararin-theorem} will be divided into several parts, some
of which involve notions and results contained in the beginning of the next chapter.

\vsp

\begin{lem} \label{evidentemente} {\em If a left-orderable group admits
only finitely many left-orders, then all of them are Conradian.}
\end{lem}
\index{Order!Conradian}

\noindent{\bf Proof.} Let $\Gamma$ be a left-orderable group whose space of left-orders
is finite. For a finite-index subgroup $\Gamma_*$ of $\Gamma$, the conjugacy
action on $\mathcal{LO}(\Gamma)$ is trivial. This means that every left-order
of $\Gamma$ is bi-invariant (hence Conradian) when restricted to $\Gamma_*$.
The lemma then follows from Proposition \ref{listailor} proved later on. $\hfill\square$

\vsp\vsp\vsp\vsp

We may now proceed to show the first claim contained in Theorem \ref{tararin-theorem}.

\vsp

\index{Rational series}
\begin{prop} \label{first-tararin}
{\em Let $\Gamma$ be a left-orderable group admitting only finitely
many left-orders. Then, for every left-order $\preceq$ on $\Gamma$,
the chain of $\preceq$-convex subgroups is a finite rational
series.}
\end{prop}

\noindent{\bf Proof.} To show finiteness of the chain of convex subgroups, 
let us fix $n \!\in\! \mathbb{N}$
such that the number of left-orders on $\Gamma$ is strictly smaller than
$2^{n}$. Following Zenkov \cite{zenkov}, we claim that the family of
$\preceq$-convex subgroups has cardinality $\leq n$. Otherwise, if
$$\{id\} = \Gamma^0 \subsetneq \Gamma^{1} \subsetneq \ldots \subsetneq
\Gamma^n = \Gamma$$
is a chain of distinct $\preceq$-convex subgroups, then
for each $\iota \! = \! (i_1,\ldots,i_n) \!\in\! \{-1,+1\}^n$
we may define the left-order $\preceq_{\iota}$ as being equal to
$\preceq_n$, where $\preceq_1,\preceq_2,\ldots,\preceq_n$ are the left-orders
on $\Gamma^1,\ldots,\Gamma^n$, respectively, which are inductively defined by:

\vsp

\noindent -- If $i_1 = 1$ (resp. $i_1 = -1$), then $\preceq_1$ is the
restriction of $\preceq$ (resp. $\overline{\preceq}$) to $\Gamma^1$;

\vsp

\noindent -- For $n \geq k \geq 2$, if $i_k = 1$ (resp. $i_k = -1$), then
$\preceq_k$ is the extension of $\preceq_{k-1}$ by the restriction of $\preceq$
(resp. $\overline{\preceq}$) to $\Gamma^{k}$.

\vsp

\noindent Clearly, the left-orders $\preceq_{\iota}$ are different
for different choices of $\iota$, which shows the claim.

\vsp

Now let $$\{ id \} = \Gamma^k \subsetneq \Gamma^{k-1} \subsetneq
\ldots \subsetneq \Gamma^0 = \Gamma$$ be the chain of {\em all}
$\preceq$-convex subgroups of $\Gamma$. In the terminology of
Example \ref{salta}, the inclusion $\Gamma^i \subsetneq
\Gamma^{i-1}$ is the convex jump associated to any element in
$\Gamma^{i-1} \setminus \Gamma^i$. By Theorem \ref{super-conrad},
$\Gamma^{i}$ is normal in $\Gamma^{i-1}$, and the induced left-order
on $\Gamma^{i-1} / \Gamma^i$ is Archimedean. By H\"older's theorem
(see \S \ref{section-holder}), this quotient is torsion-free
Abelian. Finally, its rank must be 1, as otherwise it would admit
uncountably many left-orders (see \S \ref{ejemplificando-1}),
which would allow to produce --by convex extension- 
uncountably many left-orders on $\Gamma$.
$\hfill\square$

\vspace{0.25cm}

\index{Rational series}
\begin{prop} \label{second-tararin}
{\em A left-orderable group admitting finitely many left-orders has
a unique (hence normal) rational series.}
\end{prop}

\noindent{\bf Proof.} If
$\{ id \} = \Gamma^k \lhd \Gamma^{k-1} \lhd \ldots \lhd \Gamma^0 = \Gamma$
is a rational series for a group $\Gamma$, then for every $h \in \Gamma$, the conjugate series
$$\{ id \} = h \Gamma^k h^{-1} \lhd h \Gamma^{k-1} h^{-1} \lhd \ldots \lhd  h \Gamma^0 h^{-1}
= \Gamma$$
is also rational. Therefore, the uniqueness of such a series implies its normality.

To show the uniqueness, let us consider two rational series
$$\{ id \} = G^k \lhd G^{k-1} \lhd \ldots \lhd G^1 \lhd G^0 = \Gamma, \qquad
\{id\} = H^{k'} \lhd H^{k'-1} \lhd \ldots \lhd H^1 \lhd H^0 = \Gamma,$$
where $\Gamma$ is supposed to admit only finitely many left-orders. Both
$G^1$ and $H^1$ are normal in $\Gamma$, and the quotients $\Gamma/H^1$ and
$\Gamma/G^1$ are torsion-free Abelian. This easily implies that $G^1 \cap H^1$
is also normal in $\Gamma$ and the quotient $\Gamma / (G^1 \cap H^1)$
is torsion-free Abelian. Since $G^1 \cap H^1$
is convex with respect to some left-order on $\Gamma$ (see Proposition
\ref{rel-convex}), the rank of $\Gamma / (G^1 \cap H^1)$ must be 1; otherwise,
this quotient would admit uncountably many left-orders, thus yielding --by convex
extension-- uncountably many left-orders on $\Gamma$. We conclude that, for
every $g \!\in\! G^1$ (resp. $h \in H^1$), one has $g^n \in G^1 \cap H^1$ (resp.
$f^n \in G^1 \cap H^1$) for some $n \in \mathbb{N}$. However, both $G^1$ and
$H^1$ are stable under roots, hence $g \in H^1$ (resp. $h \in G^1$).
This easily implies that $G^1 = H^1$.

Arguing similarly but with $G^1 = H^1$ instead of $\Gamma$, we
obtain $G^2 = H^2$. Proceeding in this way finitely many times, we
conclude that the rational series above coincide. $\hfill\square$

\vspace{0.35cm}

The structure of the quotients $\Gamma^{i-2} / \Gamma^i$ is given by
the (proof of the) next proposition.

\vsp

\index{Rational series}
\begin{prop} \label{third-tararin}
{\em Let $\Gamma$ be a left-orderable group having finitely many left-orders. If}
$$\{ id \} = \Gamma^k \lhd \Gamma^{k-1} \lhd \ldots \lhd \Gamma^0 = \Gamma$$
{\em is the unique rational series of \esp $\Gamma$, then no quotient
$\Gamma^{i-2} / \Gamma^i$ is bi-orderable.}
\end{prop}

\noindent{\bf Proof.} The group $\Gamma^{i-1} / \Gamma^i$ is normal in $\Gamma^{i-2} / \Gamma^i$.
Hence, $\Gamma^{i-2} / \Gamma^i$ acts by conjugacy on the torsion-free, rank-1, Abelian group
$\Gamma^{i-1} / \Gamma^i$. Now it is easy to see that every automorphism of a torsion-free,
rank-1, Abelian group is induced by the multiplication by a real number. As a consequence,
the non-Abelian group $\Gamma^{i-2} / \Gamma^i$ embeds into the affine group 
$\mathrm{Aff}(\mathbb{R})$. The
non bi-orderability of $\Gamma^{i-2} / \Gamma^i$ is thus equivalent to that the image of
this embedding is not contained in $\mathrm{Aff}_+ (\mathbb{R})$. (This is also equivalent
to that some element is conjugate to a negative power of itself.) But if this were not the
case, then, according to \S \ref{ejemplificando-2}, the quotient $\Gamma^{i-2} / \Gamma^i$
(hence $\Gamma$) would admit uncountably many left-orders. $\hfill\square$

\vspace{0.35cm}

We next proceed to show the converse statements.

\vsp

\index{Rational series}
\begin{prop} \label{fourth-tararin}
{\em Let $\Gamma$ be a group admitting a normal rational series}
$$\{ id \} = \Gamma^k \lhd \Gamma^{k-1} \lhd \ldots \lhd \Gamma^0 = \Gamma$$
{\em such that no quotient $\Gamma^{i-2} / \Gamma^i$
is bi-orderable. For each $i \!\in\! \{1,\ldots,k\}$, let us choose
$g_i \in \Gamma^{i-1} \setminus \Gamma^i$. Then every left-order on $\Gamma$
is determined by the signs of the $g_i$'s. Moreover, for any
such choice of signs, there exists a left-order on $\Gamma$ realizing it.}
\end{prop}

\noindent{\bf Proof.} The realization of signs $\iota \!\in\! \{-1,+1\}^k$ proceeds
as the proof of the first claim of Proposition \ref{first-tararin}, and we leave
the details to the reader. As before, we will denote by $\preceq_{\iota}$ the
left-order that realizes the corresponding signs.

Now let $\preceq$ be a left-order on $\Gamma$, and let $\iota = (i_1,\ldots,i_n)$
be the associate sequence of signs of the $g_i$'s. To prove that the positive
cones of $\preceq$ and $\preceq_{\iota}$ coincide, it suffices to show that
$P_{\preceq_{\iota}}^+ \subset P_{\preceq}^+$ (see Exercise \ref{contains=}).
As we saw in the proof of
the preceding proposition, after changing $g_i$ by a root if necessary, we
may assume that $g_{i-1}^{-1} g_i g_{i-1} = g_i^{r_i}$ for a {\em negative}
rational number $r_i$. Since $g_k^{i_k} \in \Gamma^{k-1}$ belongs to both
$P_{\preceq_{\iota}}^+$ and $P_{\preceq}^+$, and since $\Gamma^{k-1}$
is rank-1 Abelian, we have
$$ P_{\preceq_{\iota}}^+ \cap \esp \Gamma^{k-1} \subset P_{\preceq}^+.$$
Now every element $g \in \Gamma^{k-2} \setminus \Gamma^{k-1}$ may be written as
$g = g_k^{s i_k} g_{k-1}^{t i_{k-1}}$ for some rational numbers $s$ and $t \neq 0$;
such an element is $\preceq_{\iota}$-positive if and only if $t \!>\! 0$. If besides
$s \geq 0$, then $g$ is also $\preceq$-positive. Otherwise, $s<0$, and $g$ may be
rewritten as $g = g_{k-1}^{ t i_{k-1}} g_k^{r_k s i_k}$, and since $r_k s > 0$,
this shows that $g$ is still $\preceq$-positive. Therefore, we have
$$P_{\preceq_{\iota}}^+ \cap \esp \big( \Gamma^{k-2} \setminus \Gamma^{k-1} \big)
\esp \subset \esp P_{\preceq}^+,$$
hence
$$ P_{\preceq_{\iota}}^+ \cap \esp \Gamma^{k-2} \esp \subset \esp P_{\preceq}^+.$$
Proceeding in this way finitely many times, we conclude that
$P_{\preceq_{\iota}}^+ \subset P_{\preceq}^+$. $\hfill\square$

\vspace{0.2cm}

\begin{small}
\begin{ejer}\label{tararin-romain}
Show that every Tararin group admits a unique nontrivial torsion-free Abelian
quotient, namely the quotient with respect to the maximal proper convex subgroup.
\end{ejer}

\begin{ejer} Let $\Gamma$ be a left-orderable group for which the whole
family of subgroups that are convex for some left-order on $\Gamma$ is
finite. Show that $\Gamma$ admits only finitely many left-orders.

\noindent{\underline{Remark.}} This result is also due to Tararin;
see \cite[\S 5.2]{kopytov} in case of problems.
\end{ejer}
\end{small}

We next provide a quite clarifying result on the dynamics of
the action of a Tararin group on its space of left-orders.

\vsp

\begin{prop}\label{prop two orbits of Tararin}
{\em The action of a Tararin group $\Gamma$ on its space of left-orders has two orbits. Moreover, for any two left-orders 
$\preceq$ and $\preceq'$ on $\Gamma$, there is $g \in \Gamma$ such that $\preceq_g$ and $\preceq'$ coincide on the 
subgroup $\Gamma^1$ of its rational series {\em (\ref{rat-ser})} (which corresponds to  the maximal convex subgroup
of any of its left-orders). Furthermore, if we let $h$ be any element in $\Gamma \setminus \Gamma^{1}$  
acting on $\Gamma^{1} / \Gamma^{2}$ as the multiplication by a negative number (see Proposition 
{\em \ref{third-tararin}} and its proof), 
then $g$ can be taken either in $\Gamma^1$ or in $h \Gamma^1$.}
\end{prop}

\noindent{\bf Proof.} Choose elements $g_i \in \Gamma^{i-1} \setminus \Gamma^i$, where $i \in \{2,\ldots,k\}$.
In case the signs of $g_2$ under $\preceq$ and $\preceq'$ are the same, let $h_2 := id$; otherwise, let $h_2$
be the element $h$ above. Then the sign of $g_2$ for $\preceq_{h_2}$ is the same as that for $\preceq'$.

In case the signs of $g_3$ for $\preceq_{h_2}$ and $\preceq'$ coincide, let $h_3 := id$. Otherwise, let $h_3$ be an element
in $\Gamma^1 \setminus \Gamma^2$ acting on $\Gamma^2 / \Gamma^3$ as the multiplication by a negative number.
Then the signs of both $g_2,g_3$ for $\preceq_{h_3 h_2}$ and $\preceq'$ are the same.

Continuing this way, we obtain an element $g := h_k \cdots h_2$ such that the signs of all $g_i$'s for $\preceq_g$
and $\preceq'$ coincide, where $i \in \{2,\ldots,k\}$. This certainly implies that $\preceq_g$ and $\preceq'$ are
the same when restricted to $\Gamma^1$.
$\hfill\square$

\vsp

\begin{small}

\begin{ejem}\label{ex-rivas}
Let us consider the group
$$\Gamma := \big\langle a_s, s \in \mathbb{R} \!: a_s^{-1} a_t a_s = a_t^{-1} \mbox{ whenever } t < s \big\rangle.$$
We claim that $\Gamma$ is left-orderable but has no nontrivial action on the line. Proving that $\Gamma$ is
left-orderable is easy. Indeed, every $g \in \Gamma$ may be written in normal form as
$$g = a_{s_1}^{n_1} \cdots a_{s_k}^{n_k}\, \text{, with } \;s_1>s_2>\ldots >s_k\, , \; n_i\neq0.$$
We may then declare such a $g \!\in\! \Gamma$ to be positive if $n_1 \!>\! 0$, thus getting a left-order on
$\Gamma$ (details are left to the reader). Next, assume for a contradiction that $\Gamma$ acts nontrivially
on the real line. Then there is $t  \in \mathbb{R}$ such that $a_t$ acts nontrivially.
Let $I_t$ be an open interval fixed by $a_t$ containing no fixed point of $a_t$.
By Example \ref{ex-prep}, for each  $s \!>\! t$, we have that $a_s$ has no fixed point in the closure of $I_t$,
and that $a_s (I_t)\cap I_t = \emptyset$.
Let $I_s$ be the minimal open interval fixed by $a_s$ that contains $I_t$. Example \ref{ex-prep} again
implies that for each pair of real numbers $s_1> s_2$ larger than $t$, we have
$a_{s_1}(I_{s_2})\cap I_{s_2}=\emptyset$. We thus obtain that $\{a_s(I_t)\}_{s>t}$
is an uncountable collection of disjoint open intervals, which is absurd.
\end{ejem}
\end{small}


\subsection{The space of left-orders of the free group}
\label{case-free}

\hspace{0.45cm} The space of left-orders of the free group $\mathbb{F}_n$ (with $n \geq 2$) 
is known to be homeomorphic
to the Cantor set. This is a result of McCleary essentially contained in \cite{SM85}, though
an alternative (dynamical) proof appears in \cite{order}. The general strategy of the latter 
reference proceeds as follows:

\vsp

\noindent -- Associated to a given left-order on $\mathbb{F}_n$,
let us consider the corresponding dynamical realization.

\vsp

\noindent -- If we perturb the generators of this realization (as
homeomorphisms of the line), we still have an action of the free group,
which is ``in general'' faithful, thus yielding a new left-order on $\mathbb{F}_n$.

\vsp

\noindent -- If the perturbation above is ``small'', then the new left-order is close
to the original one.

\vsp

\noindent -- Finally, the perturbation can be made so that the resulting left-order differs 
from the original one, as otherwise the original action would be ``structurally stable''
(meaning that actions that are ``close'' to it are semiconjugate), which is easily seen
to be impossible.

\vsp

\index{Order!isolated}
\begin{small} \begin{rem}
The fact that isolated left-orders induce structurally stable actions holds in full generality. The converse, 
however, is false, as it is shown, for example, by the Baumslag-Solitar group (see \S \ref{section finite-rank-solvable}). 
Nevertheless, under certain natural assumptions, the converse is still true. See \cite{mann-rivas} for more on all of this. 
\end{rem}\end{small}

As we will see in Theorem \ref{alternative-free-products} below, 
a similar but more careful argument shows that the space of left-orders
of the free product of two finitely-generated left-orderable groups is a Cantor set.

In another direction, using results from \cite{kopi1,kopi2,SM85}, Clay has shown the
existence of a left-order on $\mathbb{F}_n$ whose orbit under the conjugacy action is
dense \cite{clay-lattice}. Using this, he deduces that $\mathcal{LO}(\mathbb{F}_n)$
is homeomorphic to the Cantor set by means of the argument contained in the following
exercise.

\begin{small}\begin{ejer} Let $\Gamma$ be a countable group having a left-order
whose orbit under the conjugacy action is dense. Show that $\mathcal{LO}(\Gamma)$
is a Cantor set.

\noindent{\underline{Hint.}} If there is an isolated left-order $\preceq$, then its
reverse left-order $\overline{\preceq}$ is also isolated. If there is a left-order of dense orbit,
this forces the existence of $g\in\Gamma$ so that $\preceq_g = \overline{\preceq}$. However,
this is impossible, since the signs of $g$ for both $\preceq$ and $\overline{\preceq}$ coincide.

\noindent{\underline{Remark.}} For the Baumslag-Solitar group $BS(1,\ell)$, no left-order
has dense orbit. However, from the description given in \S \ref{ejemplificando-3}, it
readily follows that, for each irrational number $\varepsilon \neq 0$, the orbit of
$\preceq_{\varepsilon}$ under the action of the whole group of {\em automorphisms}
is dense.
\end{ejer}\end{small}

Actually, the fact that an orbit is dense is not rare but {\bf \em generic}} in the space of left-orders of 
the free group. This means that it holds on a $G_{\delta}$-set of such orders, that is, on a set that is 
a countable intersection of dense open sets (which, according to Baire's theorem, is a dense set). 

\vspace{0.1cm}

\begin{prop} \label{classical}
{\em Let $\Gamma$ be a countable left-orderable group. If $\Gamma$
admits a dense orbit under the conjugacy action, then this is the case for the
orbits of a $G_{\delta}$-set of points in $\mathcal{LO}(\Gamma)$.}
\end{prop}

\noindent{\bf Proof.} Consider an arbitrary finite family of elements
$f_1,\ldots,f_k$ in $\Gamma$ for which the basic open set
$V_{f_1} \cap \ldots \cap V_{f_k}$ is nonempty. (Recall that
$V_f := \{\preceq : id \prec f\}$.) Let
$\mathcal{LO}(\Gamma; f_1, \ldots, f_k)$ be the subset of
$\mathcal{LO}(\Gamma)$ formed by the left-orders $\preceq$ for
which there exists $g \in \Gamma$ such that $\preceq_g$ belongs to
$V_{f_1} \cap \ldots \cap V_{f_k}$. Then the set
$\mathcal{LO}(\Gamma; f_1, \ldots, f_k)$ is a union of open basic sets,
hence open. Moreover, since we are assuming the existence of a dense orbit,
$\mathcal{LO}(\Gamma; f_1, \ldots, f_k)$ is also dense.

Now let $\mathcal{LO}^*(\Gamma)$ be the (countable) intersection of all the sets $\mathcal{LO}(\Gamma; f_1, \ldots, f_k)$ 
obtained above. Then $\mathcal{LO}^*(\Gamma)$ is a $G_{\delta}$-subset of $\mathcal{LO}(\Gamma)$, 
and the the definition easily yields that every left-order in $\mathcal{LO}^*(\Gamma)$ has a dense orbit.
$\hfill\square$

\vspace{0.3cm}

A  dynamical proof of the existence of a left-order with dense orbit in $\mathcal{LO}(\mathbb{F}_n)$
is given in \cite{rivas-free}. This proof is based on the following construction (which is closely related 
to the ideas of \cite{SM85}):

\vsp

\noindent -- Choose a countable dense set of left-orders $\preceq_k$ in $\mathcal{LO}(\mathbb{F}_n)$,
and for each $k$ consider the dynamical realization $\Phi_k$ of $\preceq_k$.

\vsp

\noindent -- Fix a sequence of positive integers $n(k)$ converging to infinity very fast,
and a family of disjoint intervals $[r(k),s(k)]$ that is unbounded in both directions.

\vsp

\noindent -- For each $k$, take a conjugate copy on $[r(k),s(k)]$ of the restriction of
$\Phi_k$ to $[-n(k),n(k)]$.

\vsp

\noindent -- Take extensions of the generators of
$\mathbb{F}_n$ so that they become homeomorphisms of the line.

\vsp

Roughly, the resulting action encodes all the possible ``finite information'' of every 
left-order of $\mathbb{F}_n$. By carefully performing the construction, we can ensure 
there is a single orbit that contains the ``center'' of every $[r(k),s(k)]$. This construction 
therefore yields a new left-order $\preceq$ on $\mathbb{F}_n$. Furthermore, 
through suitable conjugacies inside $\mathbb{F}_n$, this left-order ``captures'' 
all the aforementioned information. In concrete terms, the orbit of
$\preceq$ under the conjugacy action is dense.

\begin{small}
\begin{ejem} For Thompson's group F,
no description of {\em all} left-orders is available. (For bi-orders,
see \S \ref{ejemplificando-4}.) Actually, it is unknown whether its space of left-orders
is a Cantor set. This question is actually open for all non-solvable groups of piecewise-affine
homeomorphisms of the interval.
\end{ejem}

\begin{ejem} It can be checked in many ways that the $\pi_1$ of orientable surfaces are left-orderable.   
We refer to \S \ref{ejemplificando-3} for this. Alternatively, note that  
these groups are torsion-free and 1-relator, and all such groups are locally indicable, a 
property that, as discussed in \S \ref{classic-conrad}, is stronger than left-orderability (see Remark 
\ref{1-relator-LI}). Following the lines of the proof above, it has been recently shown in \cite{ABR} that 
the spaces of left-orders of compact hyperbolic surface groups are homeomorphic to the Cantor set; 
actually, these groups also admit left-orders with dense orbits under the conjugacy action.
\end{ejem}
\end{small}

\vsp

\noindent{\bf The case of free products.} Following a short argument from \cite{rivas-free}, 
we next show that the free product of two arbitrary left-orderable groups can be ordered in many ways.

\vsp

\index{Order!isolated}
\begin{thm} {\em The space of left-orders of the free product
of any two left-orderable groups has no isolated point.}
\label{alternative-free-products}
\end{thm}

\vsp

\noindent{\bf Proof.} We assume that the factors $\Gamma_1$, $\Gamma_2$ of our
free product $\Gamma := \Gamma_1 * \Gamma_2$ are finitely-generated; for
the general case, see Exercise \ref{completar}. Given a left-order $\preceq$
on $\Gamma$, we consider the associated dynamical realization (see Proposition
\ref{recta} and the comments after it). Fix a finite system of generators of
$\Gamma$, and for each $n \geq 1$, let $f_n$ (resp. $g_n$) be the element 
of the ball of radius $n$ that is the smallest (resp. largest) with respect to $\preceq$.
 Let $\varphi_n$ be an orientation-preserving homeomorphism of the
real line that is the identity on $[t(f_n),t(g_n)]$. Consider the following
action of $\Gamma$ on the real line: for $g \in \Gamma_1$, the action is
the conjugate of its $\preceq$-dynamical realization under $\varphi_n^{-1}$; 
for $g \in \Gamma_2$, the action is its $\preceq$-dynamical realization.
(Note that this yields an action since $\Gamma$ is a free product; however,
this action may fail to be faithful.) We claim that if $(x_k)$ is a dense sequence
of points starting at $x_1 := t(id)$, then the positive cone of the induced dynamical
lexicographic left-order $\preceq_{\varphi_n}$ coincides with that of $\preceq$
on the ball of radius $n$. Indeed, by construction, an element $h \in \Gamma$
belongs to $P_{\preceq_{\varphi_n}}$ if and only if $\varphi_n^{-1} h \varphi_n
(t(id)) > t(id)$. Now, since $\varphi_n|_{[t(f_n),t(g_n)]}$ is the identity map, this is
equivalent to $\varphi_n^{-1} h (t(id)) > t(id)$, that is, $\varphi_n^{-1} (t(h)) > t(id)$.
Now, if $h$ belongs to the ball of radius $n$, then $t(h)$ lies in $[t(f_n),t(g_n)]$,
hence $\varphi_n^{-1} (t(h)) = t(h)$, and this is bigger than $t(id)$ if and only if
$h$ is $\preceq$-positive.

Now fix $n \in \mathbb{N}$ and let us perform the preceding construction
with a map $\varphi_n$ such that $\varphi_n (s_2) = s_1$ for $s_1,s_2$
satisfying $t(g_n) < s_1 < t(h_{1,n}) < t(h_{2,n}) < s_2$,
where $h_{i,n}$ is in $\Gamma_i$. Since $t(h_{1,n}) < t(h_{2,n})$,
we have $h_{1,n} \prec h_{2,n}$. Now, from \esp $\varphi_n^{-1} (t(h_{1,n}))
> \varphi_n^{-1} (s_1) = s_2 > t(h_{2,n})$ \esp we obtain \esp
$\varphi_n^{-1} h_{1,n} \varphi_n (t(id)) > h_{2,n} (t(id))$, \esp
which by construction is equivalent to \esp $h_{1,n} \succ_{\varphi_n} h_{2,n}$.

Although the left-order $\preceq_{\varphi_n}$ may be partial (this arises when
the new action of $\Gamma$ is unfaithful), it can be extended (using the convex
extension procedure) to a left-order $\preceq_n$. By construction, the positive
cones of $\preceq$ and $\preceq_n$ coincide on the ball of radius $n$, though
$\preceq_{n}$ and $\preceq$ are different. This concludes the proof. $\hfill\square$

\vsp

\begin{small}
\begin{ejer}\label{completar} Provide the details of the proof of the preceding theorem
for factors which are not finitely-generated.

\noindent{\underline{Hint.}} Use a compactness type argument.
\end{ejer}

\index{Order!isolated}
\begin{rem} \label{r:mann-rivas}
The preceding theorem doesn't hold for direct products. Indeed, using dynamical 
methods, it is shown in \cite{mann-rivas} that the space of left-orders of $\mathbb{F}_2 \times \mathbb{Z}$ 
has isolated points. However, we currently don't know of any example of left-orderable groups $\Gamma_1$ and 
$\Gamma_2$ such that their individual spaces of left-orders are the Cantor set, but the space of left-orders 
of the product $\Gamma_1 \times \Gamma_2$ is not a Cantor set.
\end{rem}
\end{small}
\vsp

\noindent{\bf A geometric/combinatorial proof.} The fact that the space of left-orders of the free group is a 
Cantor set can also be established by an alternative argument which is roughly summarized in the two steps below:

\vsp \vsp

\index{Order!isolated}
\noindent  {\underline {Step I.}}  If a left-order is an isolated point in the space of left-orders
of the free group, then its positive cone must be finitely-generated as a semigroup.

\vsp \vsp

\noindent {\underline{Step II.}} There is no finitely-generated positive cone in the free group.

\vsp \vsp

Concerning Step I, it is not hard to see that a finitely-generated positive cone yields an isolated
point in the space of left-orders (see Proposition \ref{paja} below), yet the converse is not
necessarily true (see Example \ref{pos-riv}). However, the converse can be directly 
established for free groups, as was cleverly shown by Clay and Smith in \cite{clay-smith}. 

\vsp

\begin{question} For which families of groups, isolated left-orders are forced to have finitely-generated 
positive cones~? In particular, is this the case for left-orderable hyperbolic groups~?
\end{question}

\vsp

Below we reproduce the (quite involved) proof of Clay and Smith. Let us stress that it would be 
desirable to get more transparent geometric/combinatorial arguments that apply to other groups, 
for instance small cancelation or hyperbolic groups, as suggested above.

\vsp\vsp

\index{Order!isolated}
\begin{thm} {\em If a left-order on $\mathbb{F}_n$ is isolated in the space of left-orders,
then its positive cone must be finitely-generated as a semigroup.}
\label{thm-fin-gen-free}
\end{thm}

\noindent{\bf Proof.} Let $B_N := B_N(id)$ denote the ball of radius $N$ with respect
to the canonical system of generators. Say that a subset $S \subset \mathbb{F}_n$ is
{\em total at length} $N$ if it is {\em antisymmetric} ({\em i.e.}, 
$g \in S \implies g^{-1} \notin S$) and for all $g \in B_N \setminus \{id\}$, either $g \in S$
or $g^{-1} \in S$. (Note that $id \notin S$.) The crucial point of the proof is the following 
claim. 

\vsp\vsp

\noindent{\underline{Claim (i).}} If $S \subset B_N$ is total at length $N-1$ and satisfies
$S = \langle S \rangle^+ \cap B_N$, then for every element $g$ of length $N$ not
lying in $S \cup S^{-1}$, the semigroup $\langle S, g \rangle^+$ remains antisymmetric.

\vsp

Let us assume this for a while, and let $\preceq$ be an isolated left-order on $\mathbb{F}_n$.
Let $f_1,\ldots,f_k$ be finitely many $\preceq$-positive elements such that $\preceq$ is
the only left-order on $\mathbb{F}_n$ for which all these elements are positive. If $P_{\preceq}$
is not finitely-generated as a semigroup, then there must exist an increasing sequence of
integers $N_m$ such that each set $S := P_{\preceq} \cap B_{N_m}$ is total at length
$N_m$ though there is $g=g_m$ of length $N_m + 1$ that is not contained in $S \cup S^{-1}$.
By Claim (i), the semigroup $\langle S,g \rangle^+$ is antisymmetric. Since it is total
of length $N_m$, using Claim (i) inductively, we may extend it to an antisymmetric
semigroup which, together with its inverse, covers $\mathbb{F}_n \setminus \{ id \}$, thus
inducing a left-order on $\mathbb{F}_n$. Obviously, the same procedure can be carried out
starting with $g^{-1}$ instead of $g$. Now, if $N_m$ is sufficiently large so that $f_1,\ldots,f_k$
are all contained in $B_{N_m}$, then the procedure above would yield at least two different
left-orders with all these elements positive (one with $g$ positive, the other with $g$ negative).
This is a contradiction.

Let us now proceed to the proof of Claim (i). To do this, let's say that a finite subset $S \subset \mathbb{F}_n$
is {\em stable} if for all $f,g$ in $S$, the product $fg$ lies in $S$ whenever $\| fg \| \leq \max\{ \|f\|, \|g\|\}$.
(Here and in what follows, $\| \cdot \|$ stands for the word-length on $\mathbb{F}_n$.)

\vsp\vsp

\noindent{\underline{Claim (ii).}} If $S \subset \mathbb{F}_n$ is stable and $g \in \langle S \rangle^+$
is written in the form $g = h_1 \cdots h_{k}$, with each $h_i \in S$ and $k$ minimal, then for each
$1 \leq i \leq k-1$,
$$\|h_{i+1} \cdots h_k\| < \| h_i h_{i+1} \cdots h_k \|.$$

\vsp


Write $h_i = f_i \bar{f}_i$, $h_{i+1} = \bar{f}_i^{-1} f_{i+1}$, with no cancellation in
$f_if_{i+1} = h_i h_{i+1}$. We claim that $\| \bar{f}_i \| < \| h_i \| / 2$ and
$\| \bar{f}_i \| < \| h_{i+1}\| / 2$. Indeed, if $\| \bar{f}_i \| \geq \| h_i \| / 2$ then
\begin{eqnarray*}
\|h_ih_{i+1}\| = \|f_i\| + \|f_{i+1}\|
&=& \big( \|h_i\| - \|\bar{f}_i\| \big) + \big( \|h_{i+1}\| - \|\bar{f}_i\| \big)\\
&\leq&  \Big( \|h_i\| - \frac{\|h_i\|}{2} \Big) +
             \Big( \|h_{i+1}\| - \frac{\|h_i\|}{2} \Big) = \|h_{i+1}\|,
\end{eqnarray*}
which forces $h_ih_{i+1} \in S$, thus contradicting the minimality of $k$. The
proof of the inequality $\| \bar{f}_i \| < \| h_{i+1} \| / 2$ proceeds similarly.

We may hence write $h_i = \bar{f}_{i-1}^{-1} g_i \bar{f}_i$, where $g_i$ is not the empty word
and $h_i h_{i+1} = \bar{f}_{i-1}^{-1} g_i g_{i+1} \bar{f}_{i+1}$, without cancelation for all $i$.
It follows that $h_{i+1} \cdots h_k = \bar{f}_i^{-1} g_{i+1} \cdots g_k \bar{f}_k$, without cancelation.
Since $\|\bar{f}_i\| < \| h_i \| - \| \bar{f}_i \| = \|\bar{f}_{i-1}\| + \|g_i\|$, finally we have
\begin{eqnarray*}
\|h_{i+1} \cdots h_k\|
&=& \|\bar{f}_i\| + \|g_{i+1}\| + \ldots + \|g_k\| + \|\bar{f}_k\|\\
&<& \|\bar{f}_{i-1}\| + \|g_i\| + \|g_{i+1}\| + \ldots + \|g_k\| + \|\bar{f}_k\|
\hspace{0.1cm} = \hspace{0.1cm} \|h_i \cdots h_k\|.
\end{eqnarray*}

\vsp\vsp

\noindent{\underline{Claim (iii).}} For every subset $S \subset B_N$, the equality
$S = \langle S \rangle^+ \cap B_N$ holds if and only if for each $f,g$ in $S$ such
that $\| fg \| \leq N$, the element $fg$ lies in $S$.

\vsp

The forward implication is obvious. For the converse, given $g \in \langle S \rangle^+ \cap B_N$,
write it in the form $g = h_1 \cdots h_k$, with each $h_i$ in $S$ and $k \leq N$. Since the
hypothesis implies that $S$ is stable, we may apply Claim (ii), thus yielding
$$\|h_k\| \leq \| h_{k-1} h_k \| \leq \ldots \| h_2 \cdots h_k \| \leq \| h_1 \cdots h_k \| = \|g\| \leq N.$$
Again, since $S$ is stable,
this implies $h_{k-1} h_k \in S$; consequently, by induction $h_{k-2} h_{k-1}h_k \in S$, and so on, 
until finally $h_2 \cdots h_{k-1}h_k \in S$ and $g = h_1\cdots h_{k-1} h_k \in S$, as claimed.

\vsp\vsp

\noindent{\underline{Claim (iv).}} If $f,g$ are reduced words in $\mathbb{F}_n$, with
$\|fg\| = N$, $\|f\| \leq N$, $\|g\| = N$, then $\|f\|$ must be even and exactly half of
$f$ must cancel in the product $fg$. Moreover, after cancelation, at least the right half of
$fg$ must be the same as the right half of $g$.

\vsp

Indeed, write $f = h_1 \bar{h}$ and $g = \bar{h}^{-1} h_2$, so that $fg = h_1 h_2$, without cancelation.
Then
$$\|h_1\| + \|\bar{h}\| \leq N, \quad \|h_2\| + \|\bar{h}\| = N, \quad \|h_1\| + \|h_2\| = N.$$
The last two eaqualities yield $\|h_1\| = \|\bar{h}\|$. Therefore, $\|f\| = 2\|h_1\|$ is even, and
$\|h_1\| = \|f\| / 2$, so that exactly half of $f$ disappears in the product $fg$. Moreover, from the
first two relations we obtain $\|h_1\| \leq \|h_2\|$, hence $\|\bar{h}\| \leq \|h_2\|$, which shows
that at least the right half of $g$ survives in the product $fg$.

\vsp\vsp

We can finally finish the proof of Claim (i). Let us begin by letting $S_1 := S \cup \{ g \}$ and, for $i > 0$,
$$S_{i+1} = S_i \bigcup \hspace{0.02cm}
\{ fg \!: f,g \mbox{ in } S_i \mbox{ and } \| fg \| \leq N, fg \notin S_i \}.$$
Obviously, there must exist an index $j$ such that $S_j = S_{j+1}$. By Claim (iii), for such a
$j$, we have $S_j = \langle S, g \rangle^+ \cap B_N$, and $S_j$ is stable.

Assume for a contradiction that $\langle S,g \rangle$ is not antisymmetric. Since Claim (ii) easily implies that the semigroup
generated by an stable set excluding $id$ is antisymmetric, we must have $id \in S_j$, hence there is a smallest index
$k$ such that  both $h,h^{-1}$ belong to $S_k$ for a certain element $h$.

Suppose that $h \in S$ and $h^{-1} \in S_k$, and write $h^{-1} = h_1 h_2$, with $h_1,h_2$ in $S_{k-1}$.
Either $h_1 \notin S$ or $h_2 \notin S$ (otherwise, $h^{-1}$ would be in $S$). Let us consider the first case
(the other is analogous). Then $h_1^{-1} = h_2 h$ belongs to $S_k$, as $h_2 \in S_{k-1}$ and $h \in S \subset S_{k-1}$.
However, $h_1^{-1} \notin S$, otherwise $h_1,h_1^{-1}$ would be both in $S_{k-1}$, thus contradicting the minimality of $k$.
Summarizing, we have that $h_1$ and $h_1^{-1}$ are both in $S_k$, though $h_1 \notin S$ and $h_1^{-1} \notin S$.

The preceding argument allows reducing the general case to that where $h \notin S$ and $h^{-1} \notin S$. Since $S$
is total at length $N-1$, by the minimality of $k$, every element in $S_{k-1} \setminus S$ must have length $N$.

\vsp\vsp

\noindent{\underline{Claim (v).}} Every element in $S_{k-1} \setminus S$, as well as $h$ and $h^{-1}$,
may be written in the form $h_1 g h_2$, where $h_1,h_2$ lie in $S \cup \{id\}$ and have both even length,
and where exactly the left (resp. right)  half  of $h_2$ (resp. $h_1$) cancels in the product $h_1 g h_2$
above. (Note that this implies $\| h_1 g h_2 \| = \|h_1g\| = \|gh_2\| = N$.)

\vsp

The proof is made by induction on $i \leq k$ for elements $f \in S_i \setminus S$ with $\| f \| = N$. In the case $i = 1$,
such an element $f$ corresponds to $g$, which is written in the desired form. For the induction step, we must consider
three different cases:

\noindent -- Assume $f$ is a product $f = h_1gh_2 \bar{h}_1g\bar{h}_2$, with both $h_1gh_2$ and $\bar{h}_1 g \bar{h}_2$
in $S_{i-1} \setminus S$ of length $N$. By Claim (iv), $N$ must be even, and exactly the right half of $h_1gh_2$ must
cancel with the left half of $\bar{h}_1 g \bar{h}_2$ in the product. By the induction hypothesis, the former is nothing
but the right half of $g$ followed by $h_1$, and the latter is $\bar{h}_1$ followed by the left half of $g$. Thus
$h_1 g h_2 \bar{h}_1 g \bar{h}_2 = h_1 g \bar{h}_2$ after cancelation, so that $f$ has the desired form.

\noindent -- Suppose $f = f h_1 g h_2$, where $f \in S$, $h_1gh_2 \in S_i \setminus S$, $\|f\| \leq N$, $\|h_1gh_2 \| = N$.
By Claim (iv), exactly the right half of $f$ cancels in $fh_1gh_2$. If $\|f\| \leq \|h_1\|$, this implies that this cancelation
happens in the product $f h_1$, so that $\|fh_1\| = \|f\| \leq N$, thus yielding $fh_1 \in S$, because $S$ is stable.
If $\|h_1\| \leq \|f\|$, then the entire left half of $h_1$ cancels in $fh_1$, so that $\|fh_1\| \leq \|f\| \leq N$,
yielding again $fh_1 \in S$. That $fh_1$ has even length and half of it cancels in the product $fh_1gh_2$ now
follows from Claim (iv).

\noindent -- Finally, the case where $f = h_1 g h_2 f$, with $f \in S$, $h_1gh_2 \in S_i \setminus S$,
$\|f\| \leq N$, $\|h_1gh_2 \| = N$, can be treated in a similar way to that of the preceding one.

\vsp\vsp

To conclude the proof of Theorem \ref{thm-fin-gen-free}, let us finally
write $h = h_1 g h_2$ and $h^{-1} = \bar{h}_1 g \bar{h}_2$
as in Claim (v). There are two cases to consider:

\noindent -- If $N$ is even, write $g = g_1 g_2$, where $\| g_1 \| = \| g_2 \| = \|g\|/2$. Then
$$id = hh^{-1} = h_1g_1g_2h_2 \bar{h}_1 g_1 g_2 \bar{h}_2.$$
In this product, the right half of $h$ must cancel against the left
half of $h^{-1}$, that is, $g_2 h_2 \bar{h}_1 g_1 = id$.
Therefore, $id = h_1 g_1 g_2 \bar{h}_2 = h_1 g \bar{h}_2$. But this implies
$g^{-1} = \bar{h}_2 h_1 \in \langle S \rangle^{+} \cap B_N = S$, which is a contradiction.

\noindent -- If $N$ is odd, write $g = g_1 f g_2$, where $f$ is the generator of $\mathbb{F}_n$ that appears in the central
position when writing $g$ in reduced form. Proceeding as before, we get $id = hh^{-1} = h_1 g_1 f^2 g_2 \bar{h}_2$
with no further cancelation. However, this is absurd. $\hfill\square$

\vsp\vsp\vsp

\begin{small}
\begin{rem} There are uncountably many left-orders on $\mathbb{F}_n$ for which the
canonical generators $f_i$ are positive. Indeed, this is a direct consequence of the preceding
theorem, though it can be proved in a much more elementary way. What is less
trivial is that there are left-orders on $\mathbb{F}_n$ that extend the lexicographic
order on $\langle f_1, \ldots, f_n \rangle^+$. This is proved in \cite{Sunic} via a
concrete realization of the free group as a group of homeomorphisms of the real
line. The order thus obtained is perhaps the simplest left-order definable 
on $\mathbb{F}_n$, and can be described as follows. To define it, let
$\varphi: \mathbb{F}_n \to \mathbb{R}$ be the function
\begin{multline*}
\varphi(f) \hspace{0.15cm} = \hspace{0.15cm}
\big| \{ \mbox {subwords of } f \mbox{ of the form } f_j f_i^{-1}, j > i\} \big|
\\
- \big| \{ \mbox {subwords of } f \mbox{ of the form } f_j^{-1} f_i, j > i\} \big|
\\
+ \frac{1}{2} \left \{ \begin{array} {l}
1 \hspace{0.67cm} \mbox{ if } f \mbox{ ends with } f_1, f_2, \ldots f_n,\\
-1 \hspace{0.35cm} \mbox{ if } f \mbox{ ends with } f_1^{-1}, f_2^{-1},\ldots,f_n^{-1},\\
0 \hspace{0.67cm} \mbox{ if } f \mbox{ is trivial}, \end{array} \right.
\end{multline*}
where $f$ is a reduced word on $f_i^{\pm 1}$. 
Then declare that $f \succ id$ if and only if $\varphi(f) > 0$.

A slight modification of the method above actually provides a Cantor
set of left-orders, each of which extends the lexicographic order.
\end{rem}
\end{small}

\vsp

Step II. above can be deduced from a the work of S.~Hermiller and Z.~$\check{\mathrm{S}}$uni\'c \cite{HSunic}, 
who showed that no positive cone in a free product of groups is finitely-generated. In fact, more generally, they proved 
that no such positive cone can be described by a regular language\footnote{Informaly, a subset  of a 
finitely-generated 
group is describable by a regular language if it corresponds to a specific set of paths in a finite  graph labelled with the 
group generators. Note that finitely-generated subsemigroups can certainly be described by regular languages. }. 
This was subsequently re-interpreted (and further generalized) in geometric terms in \cite{AABR}, which is the 
approach we adopt here. We begin with a general observation presented in the following exercise.

\begin{small}
\begin{ejer} \label{ejer large balls} 
Let $P$ be the positive cone of a left-order on a finitely-generated group $\Gamma$. 
Show that both $P$ and $P^{-1}$ contain balls of arbitrarily large radius, that is, for each 
$k \in \mathbb{N}$, they contain sets of the form $B_k (g) := g B_k (id)$.
\end{ejer}

\noindent{\underline{Hint.}} 
Look at the dynamical realization: if $g$ is an element that moves the origin far away to the right 
(resp. left), then a big ball $B_k (g) := g B_k (id)$ around $g$ is contained in $P$ (resp. $P^{-1}$). 

\end{small}

\vspace{0.3cm}

The argument is by contradiction. 
Assume that the positive cone $P$ of a left-order on $ \mathbb F_n$ is finitely-generated, say 
$P=\langle g_1,\ldots, g_m \rangle^+$. Let $K := \max \{|| g_i|| \!: i \in \{ 1,\ldots,m \} \}$, where $|| \cdot ||$ 
stands for the word-length with respect to the standard generating system of $\mathbb F_n$. By 
Exercise \ref{ejer large balls}, the negative cone $P^{-1}$ contains a ball $B$ of radius $K+1$ 
centered at some element $g \in P^{-1}$.

Now, since the Cayley graph of the free group with respect to the standard generating system is a tree (see  
\S \ref{ejemplificando-3}), the ball $B$ disconnects $\mathbb F_n$. Denote by $B^{id}$ the vertex set of the 
connected component of $\mathbb F_n\setminus B$ containing $id$, and by $B^\infty$ the vertex set of the 
complement of $B^{id}$ in $\mathbb F_n\setminus B$. 
}

\vsp\vsp\vsp

\noindent {\underline{Claim.} The set $B^\infty$ intersects $P$.

\vsp\vsp

Assuming this, it is easy to obtain the desired contradiction. Indeed, if $w \in B^\infty$ is a positive element, then $w$ can 
be represented as a word $g_{i_1}\ldots g_{i_j}$ in the generators of $P$. By the definition of $K$, the distance between 
two successive prefixes $g_{i_1}$, $g_{i_1}g_{i_2}$, $g_{i_1} g_{i_2} g_{i_3}$, $\ldots$ of $w$ is at most $K$. Since $B$ 
disconnects $\mathbb F_n$, this implies that some of these prefixes must lie inside $B$. However, this contradicts the fact 
that $B$ is contained in the negative cone.

\vsp

To show the claim above, we first choose an arbitrary reduced word $u$ in $B^\infty$. By looking at the initial and the final 
generators appearing in $u$, one easily checks that 
there is a generator $v$ of $\mathbb F_n$ such that $u v u^{-1}$ and $u v^{-1} u^{-1}$ are also reduced 
words. Since both $u v u^{-1}$ and $u v^{-1} u^{-1}$ begin with $u$, they lie inside $B^\infty$. Finally, since 
one of them is positive and the other one is negative, this finishes the proof of the claim and, hence, that of Step II.

\begin{small}
\begin{rem} 
The reader will observe that some arguments from Step II actually prove that, on the free group, no positive cone $P$ is 
{\bf{\em coarsely connected}}, meaning that there is no positive integer $K$ such that the $K$-neighborhood of $P$ 
is connected. (The $K$-neighborhood of $S \subset \mathbb{F}_n$ is the set of elements that differ from an element 
of $S$ by a right factor consisting of at most $K$ factors drawn from the set of generators.) Indeed, the very last step 
of the proof above only used the finite generation of $P$ to infer its coarsely connectedness, thereby deriving a 
contradiction. The phenomenon that positive cones of left-orders cannot be coarsely connected can be also 
shown to hold for free products of  left-orderable groups, fundamental groups of hyperbolic surfaces and 
more generally for {\em limit groups} in the sense of Sela   \cite{AABR}.  
\end{rem}
\end{small}

\begin{question} 
Is it true that in a left-orderable hyperbolic group, no positive cone can be described by a regular language~? 
What about the fundamental group of an hyperbolic 3-manifold~? 
\end{question}

\vsp

A refinement of the argument above, due to Kielak \cite{kielak}, applies in 
more generality to groups of fractions of finitely-generated semigroups inside
groups with infinitely many ends. Recall that given a semigroup $P$ inside a group
$\Gamma$, we say that $\Gamma$ is the {\bf{\em group of fractions}} of $P$ if every element
in $\Gamma$ can be written in the form $gh^{-1}$, where both $g,h$ lie in $P \cup \{1\}$. An
illustrative example is given below. 
For the statement, recall that a semigroup $S$ is said to be a {\bf \em {free semigroup}} if 
for every nontrivial $f,g$ in $S$, different words in positive powers of $f$ and $g$ yield 
different elements in $S$. It is a good exercise to show that nilpotent groups do not 
contain free subsemigroups, so that the next proposition applies to them.
\index{Group!of fractions}

\vsp

\begin{prop} {\em Let $\Gamma$ be a group generated by a finite set of elements
$f_1,\ldots,f_k$, and let $P$ be the semigroup generated by them together with
$id$. If $\Gamma$ has no free sub-semigroup, then each of its elements may
be written in the form $fg^{-1}$ for certain $f,g$ in $P$.}
\label{fractions-free-semigroup}
\end{prop}

\noindent{\bf Proof.} We first claim that, given any $f,g$ in $P$, there
exist $\overline{f}, \overline{g}$ in $P$ such that $g^{-1} f = \bar{f} \bar{g}^{-1}$,
that is $f \bar{g} = g \bar{f}$. Otherwise, we would have $f P \cap g P = \emptyset$,
and this implies that the sub-semigroup generated by $f$ and $g$ is free. Indeed,
if $h_1$ and $h_2$ are different words in positive powers of $f,g$, to see that
$h_1 \neq h_2$ we may assume that $h_1$ begins with $f$ and $h_2$ with $g$ 
(since if they start with the same letter we may cancel it...). Then
the condition $f P \cap g P = \emptyset$ implies that $h_1 \neq h_2$, since $h_1
\in fP$ and $h_2 \in gP$.

Now let $h := f_1 g_1^{-1} f_2 g_2^{-1} \cdots f_k g_k^{-1}$ be an arbitrary element
in $\Gamma$, where all $f_i,g_i$ belong to $P$. By the discussion above, we may
replace $g_{k-1}^{-1} f_{k}$ by $\bar{f}_k \bar{g}_{k-1}^{-1}$, thus obtaining
$$h = f_1 g_1^{-1} f_2 g_2^{-1} \cdots f_{k-1} \bar{f}_k \bar{g}_{k-1}^{-1} g_k^{-1}.$$
Now, we may replace $g_{k-2}^{-1} f_{k-1} \bar{f}_k$ by an expression of the form
$\bar{f}_{k-1} \bar{g}_{k-2}^{-1}$, thus obtaining
$$h = f_1 g_1^{-1} f_2 g_2^{-1} \cdots f_{k-2}
\bar{f}_{k-1} \bar{g}_{k-2}^{-1} \bar{g}_{k-1}^{-1} g_k^{-1}.$$
Repeating this argument no more than $k-1$ times, we finally get an expression
of $f$ of the form $fg^{-1}$, where both $f$ and $g$ belong to $P$. $\hfill\square$

\vsp

\begin{small}\begin{ejer} Given $\Gamma := \mathbb{Z} \wr \mathbb{Z} =
\mathbb{Z} \ltimes \oplus_{\mathbb{Z}} \mathbb{Z}$, let $a$ be a
generator of the (left) factor $\mathbb{Z}$, and let $b$ be a generator
of the $0^{th}$ factor $\mathbb{Z}$ in the right. Show that $a,b$ generate
$\Gamma$, though the semigroup $P$ generated by them and $id$ satisfies
$aP \cap bP = \emptyset$.
\end{ejer}\end{small}

Proposition \ref{fractions-free-semigroup} is false in a very strong way for free groups.
This is the content of the next result from \cite{kielak}. As the reader may note, the 
arguments actually apply to any group having infinitely many ends.

\vsp

\begin{thm} {\em If $P$ is a finitely-generated, proper subsemigroup of $\mathbb{F}_n$,
then $\mathbb{F}_n$ is not the group of fractions of $P$.}
\label{kielak-thm}
\end{thm}

\noindent{\bf Proof.} Let $P$ be a finitely-generated subsemigroup for which $\mathbb{F}_n$
is the group of fractions. Let us consider the finitely many generators of $P$ as a
system of generators of $\mathbb{F}_n$, and let us look at the corresponding
Cayley graph. Since free groups have infinitely many ends, there exists a radius $N$
such that the complement of $B_N (id)$ has at least three connected components. For
simplicity, let us denote just by $B$ the ball centered at the identity and radius
$N$. We will show that $P = \mathbb{F}_n$, hence $P$ is not proper.

\vsp\vsp

\noindent{\underline{Claim (i).}} One has $B (P^{-1} \cup \{id\}) = \mathbb{F}_n$.

\vsp

One easily checks that this claim 
is equivalent to that $P\cup\{id\}$ intersects every ball $B_N (f) = fB$. Let us thus suppose that
\hspace{0.01cm} $(P \cup \{ id \}) \cap fB = \emptyset$ \hspace{0.01cm} for a certain $f \in \mathbb{F}_n$.
Let $E_0$ be a connected component of the complement of $fB$ not containing the identity. Let $h$ be an
arbitrary element in the complement of $B$, and let $E$ be the connected component of
$\mathbb{F}_n \setminus B$ containing $h$. Using the dynamical properties of the action of $\mathbb{F}_n$
on its space of ends (roughly, transitivity and local contraction), one can easily convince oneself that there exists
$g \in \mathbb{F}_n$ such that $gh \in E_0$, $gB \subset E_0$, and $g E$ does not contain $B$.

Write $gh$ in the form $h_1 h_2^{-1}$, with both $h_1,h_2$ in $P \cup \{id\}$.
Since $fB$ does not intersect $P \cup \{id\}$, the element $h_1$ must lie in the
connected component of $\mathbb{F}_n \setminus fB$ containing $id$. Starting from
the point $h_1$, the path obtained by concatenation with $h_2^{-1}$ must cross $fB$ as well
as $gB$. In particular, there is an element in $P^{-1}$ (namely, a terminal subword of
$h_2^{-1}$) joining some point in $gB$ to $h_1 h_2^{-1} = gh$. Thus, there is an element
of $P^{-1}$ joining an element of $B$ to $h$, which shows that $h$ belongs to $BP^{-1}$.

The preceding conclusion was established for all elements $h \in \mathbb{F}_n \!\setminus\! B$.
This obviously implies that $\mathbb{F}_n = B (P^{-1} \cup \{ id \})$, as desired.

\vsp\vsp

\noindent{\underline{Claim (ii).}} One has $P = \mathbb{F}_n$.

\vsp

Let $A$ be a subset of minimal cardinality such that $A (P^{-1} \cup \{id\}) = \mathbb{F}_n$. We
claim that $A$ must be a singleton. Indeed, if $A$ contains two elements $f \neq g$, then
we may write $f^{-1}g = h_1 h_2^{-1}$ for certain $h_1,h_2$ in $P \cup \{ id \}$. Hence
$f,g$ both belong to $f h_1 P^{-1}$. Therefore, letting $A' := A \cup\{ fh_1 \} \setminus \{f,g\}$,
we still have $A' (P^{-1} \cup \{id\}) = \mathbb{F}_n$. However, this contradicts the minimality
of the cardinality of $A$.

We thus conclude that for a certain element $h$, we have $h (P^{-1} \cup \{id\}) = \mathbb{F}_n$.
Certainly, this implies that $P^{-1} \cup \{id\} = \mathbb{F}_n$. In particular, letting $f$ be any
nontrivial element, both $f,f^{-1}$ must belong to $P^{-1}$, hence their product $ff^{-1} = id$
is also in $P^{-1}$. We thus obtain $P^{-1} = \mathbb{F}_n$, and taking inverses, this yields
$P = \mathbb{F}_n$.$\hfill\square$

\vsp\vsp

It is important to note that the preceding proof still leaves open the following question, 
for which we conjecture a negative answer.

\vsp

\begin{question} Do there exist $k > 2$ and a finitely-generated, proper subsemigroup $P$
of $\mathbb{F}_n$ such that  every element of $\mathbb{F}_n$ can be written in the form
$f_1 f_2^{-1} f_3 \cdots f_k^{(-1)^{k+1}}$, with all $f_1, \ldots, f_k$ belonging to
$P \cup \{id\}$~?
\end{question}

\vsp

A direct corollary of Theorem \ref{kielak-thm} is that $\mathbb{F}_n$ does not admit an order
with a finitely-generated positive cone. Together with Theorem \ref{thm-fin-gen-free}, this
yields again that $\mathcal{LO}(\mathbb{F}_n)$ has no isolated point, hence it is a Cantor set.

\vsp

So far we have seen two different ways to show that $\mathcal{LO}(\mathbb{F}_2)$ is a Cantor set.
A natural problem that arises is whether these can be pursued or combined to obtain finer information
of this space in geometric terms. The following question seems to be fundamental in this regard.
In order to solve it, one would need explicit estimates on the speed of
approximation of a given left-order.

\vsp\vsp

\begin{question} Is the Hausdorff dimension of the space of left-orders of the free group $\mathbb{F}_2$ 
zero, finite, or infinite ? 
\end{question}

\vsp\vsp

Another direction of research concerns algorithmic properties for orders. Indeed, using recursive functions
(in the sense of computability theory), it is not hard to construct left-orders on $F_n$ such that the problem
of elements comparison is undecidable, that is, there is no algorithm that, for every input consisting of a pair 
of elements $f,g$ in $F_n$, can decide whether $f \prec g$ or not. Of course, this cannot be the case for 
left-orders on higher-rank Abelian groups, yet in both cases, the spaces of left-orders are Cantor sets.

\vspace{0.32cm}

\noindent{\bf Left-orders v/s bi-orders.} 
Much progress has been recently made in the 
understanding of the space of bi-orders of the free group. In particular, the next important 
result was recently proved by Dovhyi and Muliarchyk \cite{bi-ordenes-free} (see also \cite{koberda}): 
The space of bi-orders of a (finitely-generated, non-Abelian) free group is homeomorphic to the Cantor set. 
This theorem solves in the affirmative a conjecture of McCleary (see \cite[page 127]{wolfe}). 
It should be compared with an analogous result for nilpotent groups, which is described 
at the end of Example \ref{ex-ext}. 


\subsection{Finitely-generated positive cones}
\label{fin-gen}

\hspace{0.45cm} We first recall a short argument due to Linnell
showing that, if a left-order $\preceq$ on a group $\Gamma$ is
non-isolated in $\mathcal{LO}(\Gamma)$, then its positive
cone is not finitely-generated as a semigroup.

\vsp

\index{Order!isolated}
\begin{prop}
{\em If the positive cone of a left-order $\preceq$ on a group
$\Gamma$ is finitely-generated as a semigroup, then $\preceq$ is
isolated in $\mathcal{LO}(\Gamma)$.} \label{paja}
\end{prop}

\noindent{\bf Proof.} If \esp $g_1,\ldots,g_k$ \esp generate $P_{\preceq}^{+}$, then
the only left-order on $\Gamma$ that coincides with $\preceq$ on any set containing
these generators and the identity element is $\preceq$ itself (see the
exercise below). $\hfill\square$

\begin{small}\begin{ejer} \label{contains=} Show that if
two left-orders $\preceq$ and $\preceq'$ on the same group satisfy
$P_{\preceq}^+ \subset P_{\preceq'}^+$, then they coincide.
\end{ejer}
\end{small}

The converse to the preceding proposition is not true. For instance, the dyadic rationals
admit only two left-orders, though none of them has a finitely-generated positive cone.
One may easily modify this example in order to obtain a finitely-generated one, as it is shown below.

\index{Order!isolated}

\begin{small}
\begin{ejem}\label{pos-riv} The group of presentation
$\Gamma := \langle a,b \!: aba^{-1} = b^{-2} \rangle$
is covered by Theorem \ref{tararin-theorem}: it has exactly four left-orders
(compare Example \ref{ups-ups}), which depend only on the signs of the generators $a,b$.
However, the positive cone of none of these orders is finitely-generated, since they all contain
a copy of $\mathbb{Z} [\frac{1}{2}]$ inside $\langle \! \langle b \rangle \! \rangle$
(the smallest normal subgroup containing $b$).
\end{ejem}\end{small}

\vsp

Besides the obvious case of $\mathbb{Z}$, perhaps the simplest example of a 
finitely-generated positive cone for a group left-order occurs for the Klein bottle group
$K_2 = \langle a,b: bab = a \rangle$: one may take  $\langle a,b \rangle^+$
as such a cone (see Figure 8 below). 

\vspace{0.5cm}


\beginpicture

\setcoordinatesystem units <0.9cm,0.9cm>

\putrule from -2 7 to 8 7
\putrule from -2 -1 to 8 -1
\putrule from 7 -2 to 7 8
\putrule from -1 -2 to -1 8
\putrule from 7 -2 to 7 8
\putrule from 1 -2 to 1 8
\putrule from 3 -2 to 3 2.7
\putrule from 3 3.25 to 3 8
\putrule from 2 -2 to 2 8
\putrule from 5 -2 to 5 8
\putrule from -2 1 to 8 1
\putrule from -2 3 to 2.7 3
\putrule from 3.3 3 to 8 3
\putrule from -2 4 to 8 4
\putrule from -2 5 to 8 5
\putrule from 0 -2 to 0 8
\putrule from 6 -2 to 6 8
\putrule from -2  6 to 8 6
\putrule from 4 -2 to 4 8
\putrule from -2 2 to 8 2
\putrule from 0 -2 to 0 8
\putrule from -2 0 to 8 0

\put{Figure 8: The positive cone $P^+ \! = \langle a,b \rangle^+$
on $K_2 = \langle a,b \!: a b a^{-1} = b^{-1} \rangle$.} at 3.1 -2.8
\put{$\bf{id}$} at 3 3
\put{} at -5 0


\plot 3.55 7 3.4 7.05 /
\plot 3.55 7 3.4 6.95 /

\plot 3.55 6 3.4 6.05 /
\plot 3.55 6 3.4 5.95 /

\plot 3.55 5 3.4 5.05 /
\plot 3.55 5 3.4 4.95 /

\plot 3.55 4 3.4 4.05 /
\plot 3.55 4 3.4 3.95 /

\plot 3.55 3 3.4 3.05 /
\plot 3.55 3 3.4 2.95 /


\plot 4.55 7 4.4 7.05 /
\plot 4.55 7 4.4 6.95 /

\plot 4.55 6 4.4 6.05 /
\plot 4.55 6 4.4 5.95 /

\plot 4.55 5 4.4 5.05 /
\plot 4.55 5 4.4 4.95 /

\plot 4.55 4 4.4 4.05 /
\plot 4.55 4 4.4 3.95 /

\plot 4.55 3 4.4 3.05 /
\plot 4.55 3 4.4 2.95 /

\plot 4.55 2 4.4 2.05 /
\plot 4.55 2 4.4 1.95 /

\plot 4.55 1 4.4 1.05 /
\plot 4.55 1 4.4 0.95 /

\plot 4.55 0 4.4 0.05 /
\plot 4.55 0 4.4 -0.05 /

\plot 4.55 -1 4.4 -0.95 /
\plot 4.55 -1 4.4 -1.05 /


\plot 5.55 7 5.4 7.05 /
\plot 5.55 7 5.4 6.95 /

\plot 5.55 6 5.4 6.05 /
\plot 5.55 6 5.4 5.95 /

\plot 5.55 5 5.4 5.05 /
\plot 5.55 5 5.4 4.95 /

\plot 5.55 4 5.4 4.05 /
\plot 5.55 4 5.4 3.95 /

\plot 5.55 3 5.4 3.05 /
\plot 5.55 3 5.4 2.95 /

\plot 5.55 2 5.4 2.05 /
\plot 5.55 2 5.4 1.95 /

\plot 5.55 1 5.4 1.05 /
\plot 5.55 1 5.4 0.95 /

\plot 5.55 0 5.4 0.05 /
\plot 5.55 0 5.4 -0.05 /

\plot 5.55 -1 5.4 -0.95 /
\plot 5.55 -1 5.4 -1.05 /


\plot 6.55 7 6.4 7.05 /
\plot 6.55 7 6.4 6.95 /

\plot 6.55 6 6.4 6.05 /
\plot 6.55 6 6.4 5.95 /

\plot 6.55 5 6.4 5.05 /
\plot 6.55 5 6.4 4.95 /

\plot 6.55 4 6.4 4.05 /
\plot 6.55 4 6.4 3.95 /

\plot 6.55 3 6.4 3.05 /
\plot 6.55 3 6.4 2.95 /

\plot 6.55 2 6.4 2.05 /
\plot 6.55 2 6.4 1.95 /

\plot 6.55 1 6.4 1.05 /
\plot 6.55 1 6.4 0.95 /

\plot 6.55 0 6.4 0.05 /
\plot 6.55 0 6.4 -0.05 /

\plot 6.55 -1 6.4 -0.95 /
\plot 6.55 -1 6.4 -1.05 /


\plot 7.55 7 7.4 7.05 /
\plot 7.55 7 7.4 6.95 /

\plot 7.55 6 7.4 6.05 /
\plot 7.55 6 7.4 5.95 /

\plot 7.55 5 7.4 5.05 /
\plot 7.55 5 7.4 4.95 /

\plot 7.55 4 7.4 4.05 /
\plot 7.55 4 7.4 3.95 /

\plot 7.55 3 7.4 3.05 /
\plot 7.55 3 7.4 2.95 /

\plot 7.55 2 7.4 2.05 /
\plot 7.55 2 7.4 1.95 /

\plot 7.55 1 7.4 1.05 /
\plot 7.55 1 7.4 0.95 /

\plot 7.55 0 7.4 0.05 /
\plot 7.55 0 7.4 -0.05 /

\plot 7.55 -1 7.4 -0.95 /
\plot 7.55 -1 7.4 -1.05 /


\plot 4 7.45 4.05 7.6 /
\plot 4 7.45 3.95 7.6 /

\plot 6 7.45 6.05 7.6 /
\plot 6 7.45 5.95 7.6 /

\plot 4 6.45 4.05 6.6 /
\plot 4 6.45 3.95 6.6 /

\plot 6 6.45 6.05 6.6 /
\plot 6 6.45 5.95 6.6 /

\plot 4 5.45 4.05 5.6 /
\plot 4 5.45 3.95 5.6 /

\plot 6 5.45 6.05 5.6 /
\plot 6 5.45 5.95 5.6 /

\plot 4 4.45 4.05 4.6 /
\plot 4 4.45 3.95 4.6 /

\plot 6 4.45 6.05 4.6 /
\plot 6 4.45 5.95 4.6 /

\plot 4 3.45 4.05 3.6 /
\plot 4 3.45 3.95 3.6 /

\plot 6 3.45 6.05 3.6 /
\plot 6 3.45 5.95 3.6 /

\plot 4 2.45 4.05 2.6 /
\plot 4 2.45 3.95 2.6 /

\plot 6 2.45 6.05 2.6 /
\plot 6 2.45 5.95 2.6 /

\plot 4 1.45 4.05 1.6 /
\plot 4 1.45 3.95 1.6 /

\plot 6 1.45 6.05 1.6 /
\plot 6 1.45 5.95 1.6 /

\plot 4 0.45 4.05 0.6 /
\plot 4 0.45 3.95 0.6 /

\plot 6 0.45 6.05 0.6 /
\plot 6 0.45 5.95 0.6 /

\plot 4 -0.55 4.05 -0.4 /
\plot 4 -0.55 3.95 -0.4 /

\plot 6 -0.55 6.05 -0.4 /
\plot 6 -0.55 5.95 -0.4 /

\plot 4 -1.55 4.05 -1.4 /
\plot 4 -1.55 3.95 -1.4 /

\plot 6 -1.55 6.05 -1.4 /
\plot 6 -1.55 5.95 -1.4 /


\plot 3 3.55
3.05 3.4 /

\plot 3 3.55
2.95 3.4 /

\plot 3 4.55
3.05 4.4 /

\plot 3 4.55
2.95 4.4 /


\plot 3 5.55
3.05 5.4 /

\plot 3 5.55
2.95 5.4 /


\plot 3 6.55
3.05 6.4 /

\plot 3 6.55
2.95 6.4 /


\plot 5 3.55
5.05 3.4 /

\plot 5 3.55
4.95 3.4 /


\plot 5 4.55
5.05 4.4 /

\plot 5 4.55
4.95 4.4 /


\plot 5 5.55
5.05 5.4 /

\plot 5 5.55
4.95 5.4 /


\plot 5 6.55
5.05 6.4 /

\plot 5 6.55
4.95 6.4 /


\plot 5 2.55
5.05 2.4 /

\plot 5 2.55
4.95 2.4 /


\plot 5 1.55
5.05 1.4 /

\plot 5 1.55
4.95 1.4 /


\plot 5 0.55
5.05 0.4 /

\plot 5 0.55
4.95 0.4 /


\plot 5 -0.45
5.05 -0.6 /

\plot 5 -0.45
4.95 -0.6 /


\plot 7 3.55
7.05 3.4 /

\plot 7 3.55
6.95 3.4 /


\plot 7 4.55
7.05 4.4 /

\plot 7 4.55
6.95 4.4 /


\plot 7 5.55
7.05 5.4 /

\plot 7 5.55
6.95 5.4 /


\plot 7 6.55
7.05 6.4 /

\plot 7 6.55
6.95 6.4 /


\plot 7 2.55
7.05 2.4 /

\plot 7 2.55
6.95 2.4 /


\plot 7 1.55
7.05 1.4 /

\plot 7 1.55
6.95 1.4 /


\plot 7 0.55
7.05 0.4 /

\plot 7 0.55
6.95 0.4 /


\plot 7 -0.45
7.05 -0.6 /

\plot 7 -0.45
6.95 -0.6 /


\plot 5 7.55
5.05 7.4 /

\plot 5 7.55
4.95 7.4 /


\plot 7 7.55
7.05 7.4 /

\plot 7 7.55
6.95 7.4 /


\plot 3 7.55
3.05 7.4 /

\plot 3 7.55
2.95 7.4 /


\plot 5 -1.45
5.05 -1.6 /

\plot 5 -1.45
4.95 -1.6 /


\plot 7 -1.45
7.05 -1.6 /

\plot 7 -1.45
6.95 -1.6 /


\put{$\bullet$} at 3 6
\put{$\bullet$} at 3 4
\put{$\bullet$} at 3 5

\put{$\bullet$} at 4 0
\put{$\bullet$} at 4 1
\put{$\bullet$} at 4 2
\put{$\bullet$} at 4 3
\put{$\bullet$} at 4 4
\put{$\bullet$} at 4 5
\put{$\bullet$} at 4 6

\put{$\bullet$} at 5 0
\put{$\bullet$} at 5 1
\put{$\bullet$} at 5 2
\put{$\bullet$} at 5 3
\put{$\bullet$} at 5 4
\put{$\bullet$} at 5 5
\put{$\bullet$} at 5 6

\put{$\bullet$} at 6 0
\put{$\bullet$} at 6 1
\put{$\bullet$} at 6 2
\put{$\bullet$} at 6 3
\put{$\bullet$} at 6 4
\put{$\bullet$} at 6 5
\put{$\bullet$} at 6 6

\put{$\bullet$} at 7 0
\put{$\bullet$} at 7 1
\put{$\bullet$} at 7 2
\put{$\bullet$} at 7 3
\put{$\bullet$} at 7 4
\put{$\bullet$} at 7 5
\put{$\bullet$} at 7 6

\put{$\bullet$} at 3 7
\put{$\bullet$} at 4 7
\put{$\bullet$} at 5 7
\put{$\bullet$} at 6 7
\put{$\bullet$} at 7 -1
\put{$\bullet$} at 7 7
\put{$\bullet$} at 4 -1
\put{$\bullet$} at 5 -1
\put{$\bullet$} at 6 -1
\put{$\bullet$} at 6 -1


\begin{footnotesize}

\put{$a$} at 3.42 7.2
\put{$a$} at 3.42 6.2
\put{$a$} at 3.42 5.2
\put{$a$} at 3.42 4.2
\put{$a$} at 3.42 3.2

\put{$a$} at 4.42 7.2
\put{$a$} at 4.42 6.2
\put{$a$} at 4.42 5.2
\put{$a$} at 4.42 4.2
\put{$a$} at 4.42 3.2
\put{$a$} at 4.42 2.2
\put{$a$} at 4.42 1.2
\put{$a$} at 4.42 0.2
\put{$a$} at 4.42 -0.8

\put{$a$} at 5.42 7.2
\put{$a$} at 5.42 6.2
\put{$a$} at 5.42 5.2
\put{$a$} at 5.42 4.2
\put{$a$} at 5.42 3.2
\put{$a$} at 5.42 2.2
\put{$a$} at 5.42 1.2
\put{$a$} at 5.42 0.2
\put{$a$} at 5.42 -0.8

\put{$a$} at 6.42 7.2
\put{$a$} at 6.42 6.2
\put{$a$} at 6.42 5.2
\put{$a$} at 6.42 4.2
\put{$a$} at 6.42 3.2
\put{$a$} at 6.42 2.2
\put{$a$} at 6.42 1.2
\put{$a$} at 6.42 0.2
\put{$a$} at 6.42 -0.8

\put{$a$} at 7.42 7.2
\put{$a$} at 7.42 6.2
\put{$a$} at 7.42 5.2
\put{$a$} at 7.42 4.2
\put{$a$} at 7.42 3.2
\put{$a$} at 7.42 2.2
\put{$a$} at 7.42 1.2
\put{$a$} at 7.42 0.2
\put{$a$} at 7.42 -0.8


\put{$b$} at 3.65 7.46
\put{$b$} at 3.65 6.46
\put{$b$} at 3.65 5.46
\put{$b$} at 3.65 4.46
\put{$b$} at 3.65 3.46
\put{$b$} at 3.65 2.46
\put{$b$} at 3.65 1.46
\put{$b$} at 3.65 0.46
\put{$b$} at 3.65 -0.54
\put{$b$} at 3.65 -1.54

\put{$b$} at 5.65 7.46
\put{$b$} at 5.65 6.46
\put{$b$} at 5.65 5.46
\put{$b$} at 5.65 4.46
\put{$b$} at 5.65 3.46
\put{$b$} at 5.65 2.46
\put{$b$} at 5.65 1.46
\put{$b$} at 5.65 0.46
\put{$b$} at 5.65 -0.54
\put{$b$} at 5.65 -1.54


\put{$b$} at 2.65 7.6
\put{$b$} at 2.65 6.6
\put{$b$} at 2.65 5.6
\put{$b$} at 2.65 4.6
\put{$b$} at 2.65 3.6

\put{$b$} at 4.65 7.6
\put{$b$} at 4.65 6.6
\put{$b$} at 4.65 5.6
\put{$b$} at 4.65 4.6
\put{$b$} at 4.65 3.6
\put{$b$} at 4.65 2.6
\put{$b$} at 4.65 1.6
\put{$b$} at 4.65 0.6
\put{$b$} at 4.65 -0.4
\put{$b$} at 4.65 -1.4

\put{$b$} at 6.65 7.6
\put{$b$} at 6.65 6.6
\put{$b$} at 6.65 5.6
\put{$b$} at 6.65 4.6
\put{$b$} at 6.65 3.6
\put{$b$} at 6.65 2.6
\put{$b$} at 6.65 1.6
\put{$b$} at 6.65 0.6
\put{$b$} at 6.65 -0.4
\put{$b$} at 6.65 -1.4

\end{footnotesize}

\endpicture


\vspace{0.3cm}

\begin{small}
\begin{ejem} \label{ejem:no-inverse}
Show that for the left-order induced by the positive cone above, if we let $f:= ba$ and $g:= ab$, then 
$f \prec g$ and $f^{-1} \prec  g^{-1}$. 
\end{ejem}\end{small}

Actually, $K_2$ admits exactly four left-orders, and 
each of these has a finitely-generated positive cone. (The other cones
are $\langle a,b^{-1} \rangle^+$, $\langle a^{-1},b \rangle^+$, and
$\langle a^{-1},b^{-1} \rangle^+$.)

Rather surprisingly, finitely-generated positive cones also occur on braid groups,
according to a beautiful result due to Dubrovina and Dubrovin \cite{dub}

\vsp

\begin{thm}\label{thm-DD}
{\em For each $n \geq 3$, the braid group $\mathbb{B}_n$ admits the decomposition}
$$\mathbb{B}_n = \big\langle a_1,\ldots,a_{n-1} \big\rangle^+ \sqcup
\big\langle a_1^{-1},\ldots,a_{n-1}^{-1} \big\rangle^+ \sqcup \{id\},$$
{\em where $a_1 := \s_1 \cdots \s_{n-1}$, $a_2 := (\s_2 \cdots \s_{n-1})^{-1}$,
$a_3 := \s_3 \cdots \s_{n-1}$, $a_4 := (\s_4 \cdots \s_{n-1})^{-1}$, $\ldots$,
and $a_{n-1} := \s_{n-1}^{(-1)^{n-1}}.$}
\end{thm}
\index{Order!DD-order}

\vspace{0.1cm}

Note that this theorem also holds for $n \!=\! 2$, yet it is trivial in 
this case, as $\mathbb{B}_2$ is isomorphic to $\mathbb{Z}$. 
For the case of $\mathbb{B}_3$, the theorem states that the semigroup
$P_{_{DD}} = \langle \sigma_1 \sigma_2, \sigma_2^{-1} \rangle^+$
is the positive cone of a left-order $\preceq_{_{DD}}$. This can be visualized in
Figure 9 below, where we depict the Cayley graph of $\mathbb{B}_3$ (essentially, a
product of a {\em quasi-isometric} copy of $\mathbb{Z}^2$ by a dyadic rooted tree). 
See \cite{Na-ICM} for a more clear picture (in colors). 

\begin{figure}[h!]
\begin{center}

\vspace{0.2cm}
\hspace{-0.5cm} 
\includegraphics[scale=0.3]{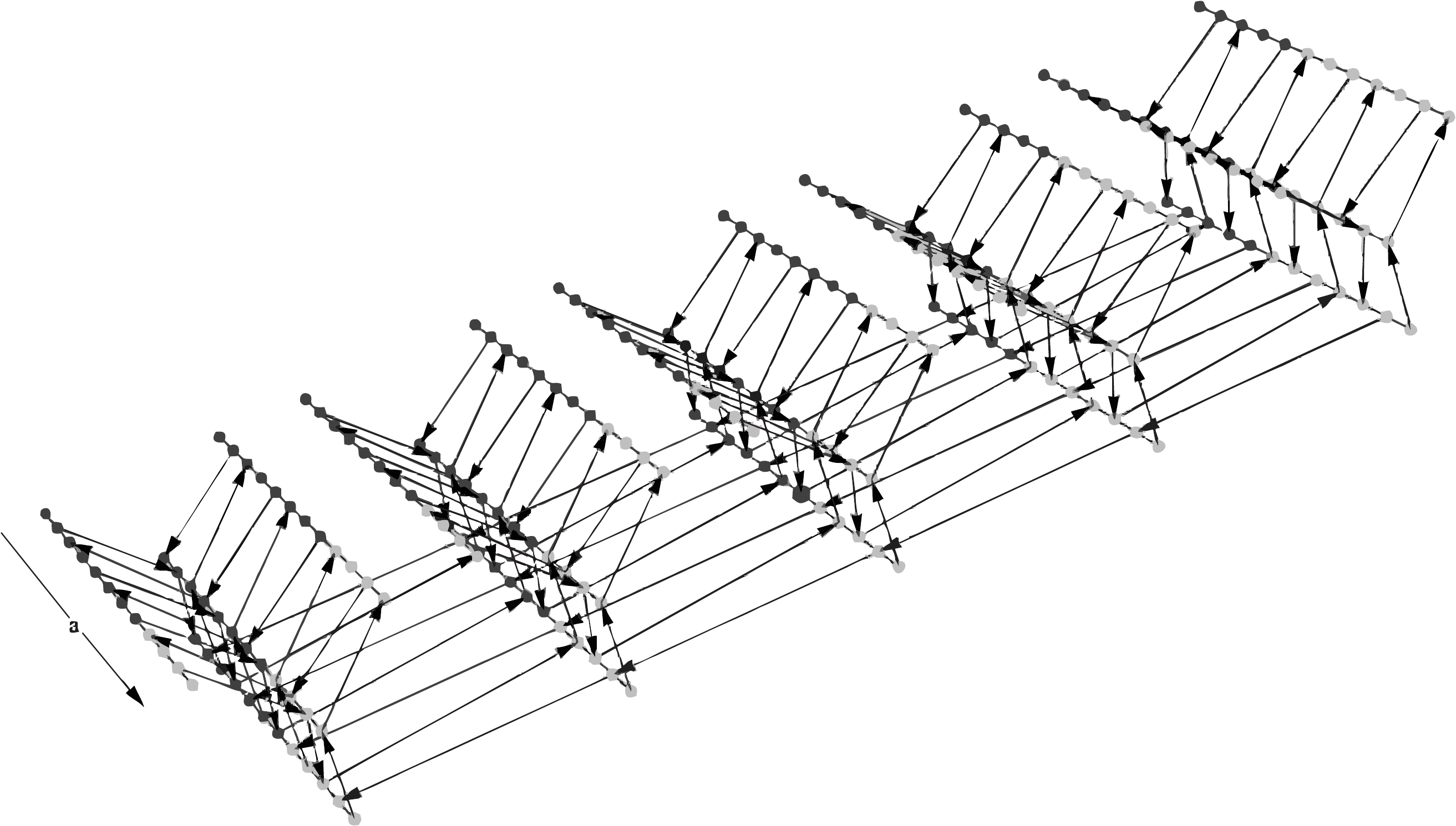}
\end{center}

\vspace{-0.5cm}
\begin{center}
Figure 9: The DD-positive cone on $\mathbb{B}_3$.
\end{center}
\end{figure}

Note that for the generators $a = a_1 := \s_1\s_2$ and $b = a_2 := \s_2^{-1}$ 
of $\mathbb{B}_3$, the presentation becomes \hspace{0.1cm}
$\mathbb{B}_3 = \big\langle a,b \! : b a^2 b = a \big\rangle.$ \hspace{0.1cm}
Thus, in the picture above,
an arrow pointing from left to right should be added to every diagonal edge of the graph.
These arrows represent multiplications by $a$, while all arrows explicitly appearing represent
multiplications by $b$. Starting at the identity, every positive element can be
reached by a path that follows the direction of the arrows. Conversely, every negative 
element can be reached by a path starting at the identity following a direction opposite to
that of the arrows. Finally, no nontrivial element can be reached both ways.

Quite remarkably, this particular example was already known for a long time (and seems to be folklore,
at least for a certain community); see \cite{D-original}.

\vspace{0.35cm}


\index{Order!DD-order}

\noindent{\bf Retrieving the $DD$-order  from the $D$-order.} The proof of Theorem \ref{thm-DD}
strongly uses Dehornoy's theorem dicussed in \S \ref{ejemplificando-5}. To begin with, note
that it readily follows from the definition that for each $j \!\in\! \{1,\ldots,n-1\}$, the subgroup \esp
$\langle \sigma_j,\ldots,\sigma_{n-1}\rangle \sim \mathbb{B}_{n-j+1}$ \esp of
$\mathbb{B}_n$ is $\preceq_{_D}$-convex.

In particular, for the case of $\mathbb{B}_3$, the cyclic subgroup \esp $\langle \sigma_2 \rangle$ \esp is $\preceq_{_D}$-convex.
One can hence define the order $\preceq_3$ on $\mathbb{B}_3$ as being the extension by $\preceq_{_D}$ of the
restriction to \esp $\langle \sigma_2 \rangle$ \esp of the reverse order $\overline{\preceq}_{_D}$. \esp (This
corresponds to flipping $\preceq_{_D}$ on $\langle \sigma_2 \rangle$, as discussed in Example \ref{ex-flipping}.)
We claim that the positive cone of $\preceq_3$ is generated by the elements \esp $a_1 \!:=\! \sigma_1 \sigma_2$
\esp and \esp $a_2 \!:=\! \sigma_2^{-1}$, thus showing the theorem in this particular case. \esp Indeed, by definition,
these elements are positive with respect to $\preceq_3$. Thus, it suffices to show that for every \esp $c \neq id$ in
$\mathbb{B}_3$, either $c$ or $c^{-1}$ belongs to \esp $\langle a_1,a_2 \rangle^+$. \esp Now, if $c$ or $c^{-1}$
is $2$-positive, then there exists an integer $m \neq 0$ such that \esp $c = \sigma_2^m = a_2^{-m}$, \esp
and therefore \esp $c \in \langle a_2 \rangle^+ \subset \langle a_1,a_2 \rangle^+$ \esp if $m < 0$,
and \esp $c^{-1} \in \langle a_2 \rangle^+ \subset \langle a_1,a_2 \rangle^+$ \esp if $m > 0$. If
$c$ is $1$-positive, then for a certain choice of integers $m_1'', \ldots,m_{k''+1}''$, one has
$$c = \sigma_2^{m_1''} \sigma_1 \sigma_2^{m_2''} \sigma_1
\cdots \sigma_2^{m_{k''}''} \sigma_1 \sigma_2^{m_{k''+1}''}.$$
Using the identity \esp $\sigma_1 \!=\! a_1 a_2$, \esp
this allows us to rewrite $c$ in the form
$$c = a_2^{m_1'} a_1 a_2^{m_2'} a_1 \ldots a_2^{m_{k'}'} a_1 a_2^{m_{k'+1}'}$$
for some integers $m_1',\ldots,m_{k'+1}'$. Now, using several times the (easy to check)
identity \esp $a_2 a_1^2 a_2 = a_1$, \esp one may easily  express $c$ as a product
$$c = a_2^{m_1} a_1 a_2^{m_2} a_1 \ldots a_2^{m_{k}} a_1 a_2^{m_{k+1}}$$
in which all the exponents $m_i$ are non-negative. This shows that $c$ belongs to
\esp $\langle a_1, a_2 \rangle^+$. \esp Finally, if $c^{-1}$ is $1$-positive, then
$c^{-1}$ belongs to \esp $\langle a_1,a_2 \rangle^+$.

The extension of the preceding argument to the general case proceeds inductively as follows. Let
us see $\mathbb{B}_{n-1} = \langle \tilde{\sigma}_1, \ldots, \tilde{\sigma}_{n-2} \rangle$ as a subgroup
of \esp $\mathbb{B}_n = \langle \sigma_1, \ldots, \sigma_{n-1} \rangle$ \esp via the homomorphism
\esp $\tilde{\sigma}_i \mapsto \sigma_{i+1}$. \esp Then $\preceq_{n-1}$ induces an order on
$\langle \sigma_2,\ldots,\sigma_{n-1} \rangle \subset \mathbb{B}_n$, which we still denote by
$\preceq_{n-1}$. We then let $\preceq_{n}$ be the extension of $\overline{\preceq}_{n-1}$ by
the $D$-order $\preceq_{_D}$. Then, using the inductive hypothesis as well as the
remarkable identities (that we leave to the reader)
$$(a_2 a_3^{-1} \cdots a_{n-1}^{(-1)^{n-1}}) a_1^{n-1}
(a_2 a_3^{-1} \cdots a_{n-1}^{(-1)^{n-1}}) = a_1$$
and
$$(a_2 a_3^{-1} \cdots a_{n-1}^{(-1)^{n-1}})^2 = a_2^{n-1},$$
one may check as above that the positive cone of the order $\preceq_n$
coincides with the semigroup $\langle a_1, \ldots, a_{n-1} \rangle^+$, thus
showing the theorem.

\index{Order!DD-order}
\begin{small}\begin{ejer}
Prove that the only convex subgroups of \esp $\mathbb{B}_{n}$ \esp
for both the $D$-order and the $DD$-order are
\esp $C^{0}\!:=\! \mathbb{B}_{n}$,
\esp $C^{1}\!:=\!\langle a_2,\ldots,a_{n-1} \rangle\!=\!
\langle \sigma_2,\ldots,\sigma_{n-1} \rangle$,
\ldots,
$C^{n-1}\!:=\!\langle a_{n-1} \rangle \!=\! \langle \sigma_{n-1} \rangle$
\esp and  \esp $C^n\!:=\!\{id\}$.
\end{ejer}\end{small}

Let us emphasize that assuming Theorem \ref{thm-DD}, we can follow the arguments above backwards
and retrieve the Dehornoy's order on $\mathbb{B}_n$. (The details are left to the reader.) A more conceptual 
approach to this phenomenon was proposed by Ito in \cite{ito-1}, and it is developed in the next exercise. 

\index{Order!Dehornoy-like}
\begin{small}\begin{ejer} Let $g_1,\ldots,g_k$ be  finitely many  generators of a group
$\Gamma$. For each $i \in \{1,\ldots,k\}$, let $h_i := (g_i g_{i+1} \cdots g_k)^{(-1)^{k+1}} \! ,$
and denote by $P_i$ the semigroup generated by $g_i,\ldots,g_k$. Assume that the following
condition (called {\bf{\em Property (F)}} in \cite{ito-1}) holds: For each $i \in \{1, \ldots, k-1\}$,
both $g_i P_{i+1} g_i^{-1}$ and $g_i^{-1} P_{i+1} g_i$ are contained in the semigroup
$P_i^{-}$ consisting of the inverses of the elements in $P_i$.

\noindent (i) Prove that $h_1,\ldots,h_k$ generate the positive cone of a left-order on $\Gamma$ if and
only if the $g_i$'s define a {\bf{\em Dehornoy-like order}}, which means that every nontrivial element
may be written as a product of elements $g_i,\ldots,g_k$ so that $g_i$ appears with only positive
exponents, and no $g \in \Gamma$ is such that both $g$ and $g^{-1}$ may be written in such a way.

\noindent (ii) Referring to Theorem \ref{thm-DD}, check that property (F) holds for $g_i := \sigma_i$ and $h_i := a_i$.
\end{ejer}\end{small}

\vspace{0.35cm}


\index{Group!Torus-knot}
\noindent{\bf Torus-knot groups.} We next give an elementary proof of that the torus-knot
groups $G_{m,n} := \langle c, d \!: c^m = d^n \rangle$ do admit left-orders with finitely-generated
positive cones. This is closely related to what was previously shown for braid groups, since for $(m,n) = (3,2)$
we retrieve the braid group $\mathbb{B}_3$ for the generators $c \sim \s_1\s_2$  and $d \sim \s_1\s_2\s_1$.
In this case, the positive cone given by Theorem \ref{thm-DD} is generated by $a :=  c \sim \s_1 \s_2$ and
$b := c^{-2} d  \sim \s_2^{-1}$, with respect to which the presentation becomes $G_{3,2} = \langle a,b \!: b a^2 b = a \rangle$.
Also, note that for $(m,n) = (2,2)$ we retrieve the Klein bottle group. In this case, the generating system of the positive
cone consists of $a := c$ and $ b := c^{-1} d$, for which the presentation becomes $K_4 = \langle a,b \!: bab = a \rangle$.

After some computations, one easily sees that the natural extension of this corresponds to the presentation
$$G_{m,n} = \big\langle a,b \!: (ba^{m-1})^{n-1} b = a \big\rangle,$$
where $a := c$ and $b := c^{-(m-1)} d$.  The following result appears in \cite{pos-1} for $n = 2$ and in \cite{ito-1}
for the general case.

\vsp\vsp

\begin{thm} \label{thm:torus}
{\em For each $m > 1$ and $n > 1$, the group $G_{m,n}$ can be decomposed as}
$$G_{m,n} = \langle a,b \rangle^+ \sqcup \langle a^{-1},b^{-1} \rangle^+ \sqcup \{ id \}.$$
\end{thm}

\vsp\vsp

Quite naturally, the proof of this theorem involves two issues:

\vsp\vsp

\noindent{\underline{Step I.}} Every nontrivial element lies in
$\langle a,b \rangle^+ \cup \langle a^{-1}, b^{-1} \rangle^+$.

\vsp\vsp

\noindent{\underline{Step II.}} No nontrivial element lies in both
$\langle a,b \rangle^+$ and $\langle a^{-1}, b^{-1} \rangle^{+}$.

\vsp\vsp

In what follows, we only consider the case $(m,n) \neq (2,2)$, because 
the choice $(m,n) = (2,2)$ corresponds to the Klein bottle group $K_4$, 
as previously explained. (Some of the arguments below do not apply
in this case.) For Step I, we begin with a simple yet crucial claim. 

\vsp\vsp

\noindent{\underline{Claim (i).}} The element $\Delta := a^m$ belongs to the center of $G_{m,n}$.

\vsp\vsp

Indeed, from $(b a^{m-1})^{n-1} b \!=\! a$, it follows that $(b a^{m-1})^n \!=\! (a^{m-1} b)^n \!=\! a^m$.
Thus,
$$b \Delta = ba^m = b (a^{m-1}b)^n = (b a^{m-1})^n b = a^m b = \Delta b.$$
Moreover, \esp $a \Delta = a^{m+1} = \Delta a$.

\vsp

A word in (positive powers of) $a,b$ (resp. $a^{-1},b^{-1}$) will be said to be {\em positive} 
(resp. {\em negative}). Also, we will say that it is {\em non-positive} (resp. {\em non-negative})
if it is either trivial or negative (resp. either trivial or positive).

\vsp\vsp

\noindent{\underline{Claim (ii).}} Every element $w \!\in\! G_{m,n}$ may be written in the form \esp
$\bar{u} \Delta^{\ell}$ \esp for some non-negative word $\bar{u}$ and $\ell \in \mathbb{Z}$.

\vsp\vsp

Indeed, in any word representing $w$, we may
rewrite the negative powers of $a$ and $b$ using the relations
$$a^{-1} = a^{m-1} \Delta^{-1}, \qquad b^{-1} = a^{-1} (b a^{m-1})^{n-1}
= a^{m-1}\Delta^{-1} (b a^{m-1})^{n-1},$$
and then use the fact that $\Delta$ belongs to the center of $G_{m,n}$.

\vspace{0.1cm}

Note that since $a^m = \Delta$,
every $u \!\in\! \langle a,b \rangle^+$ may be written in the form
$$u = b^{s_0} a^{r_1} b^{s_1} \cdots b^{s_{k-1}} a^{r_k} \Delta^{\ell},$$
where $s_i > 0$ for $i \!\in\! \{1,\ldots,k-1 \}$, $s_0 \geq 0$,
$r_i \!\in\! \{1,\ldots,m-1\}$ for $i \in \{1,\ldots,k-1\}$,
$r_k \geq 0$, and $\ell \!\geq\! 0$.
Therefore, by Claim (ii), every $w \in G_{m,n}$ may be written as
\begin{equation}\label{p}
w \esp
= \esp b^{s_0} a^{r_1} b^{s_1} \cdots b^{s_{k-1}} a^{r_k} \Delta^{\ell} \esp
= \esp b^{s_0} a^{r_1} b^{s_1} \cdots b^{s_{k-1}} a^{r_k + m \ell},
\end{equation}
where the properties of $r_i$ and $s_i$ above are satisfied, and $\ell \in \mathbb{Z}$.
Such an expression will be said to be a {\bf {\em normal form}} for $w$
if $k$ is minimal. (Note that, {\em a priori}, these normal
forms may be non-unique for a given $w$.)

There are two cases to consider. If $r_k + m \ell \geq 0$, then $w$ is obviously non-negative.
Therefore, Step I is concluded by the next claim.

\vsp\vsp

\noindent{\underline{Claim (iii):}} If $r_k + m \ell < 0$, then $w$ is negative.

\vsp\vsp

The proof is by induction on the length $k$ of the normal form. To begin with,
note that $(b a^{m-1})^{n-1} b = a$ yields $b a^{-1} = (a^{-m+1} b^{-1})^{n-1}$.
Thus, for $k = 1$, that is, for $w = b^{s_0} a^{r_1 + m \ell}$, we have
$$w 
= b^{s_0} a^{-1} a^{r_1 + m \ell + 1} 
= a^{-1} \big[ a b a^{-1} \big]^{s_0} a^{r_1 + m \ell + 1}
= a^{-1} \big[ a (a^{-m+1} b^{-1})^{n-1} \big]^{s_0} a^{r_1 + m \ell + 1},$$
and the last expression is easily seen to be negative just by noticing that
$$a (a^{-m+1} b^{-1})^{n-1} = a^{-m+2} b^{-1} (a^{-m+1} b^{-1})^{n-2}.$$

Assume the claim holds up to $k \!-\! 1$. Proceeding as before, expression (\ref{p}) becomes
$$w = b^{s_0} a^{r_1} \cdots b^{s_{k-1} - 1} [b a^{-1}] a^{r_k + m \ell + 1} =
b^{s_0} a^{r_1} \cdots b^{s_{k-1} - 1} [(a^{-m+1} b^{-1})^{n-1}] a^{m_k + m \ell + 1}.$$
If $s_{k-1} > 1$, then writing
$$w = b^{s_0} a^{r_1} \cdots b^{s_{k-1} - 2} [ba^{-1}] a^{-m+2} b^{-1} (a^{-m+1} b^{-1})^{n-2} a^{r_k + m \ell + 1},$$
we see that we can repeat the process changing $ba^{-1}$ by $(a^{-m+1} b^{-1})^{n-1}$. Otherwise,
\begin{eqnarray*}
w
&=& b^{s_0} a^{r_1} \cdots a^{r_{k-1}} [(a^{-m+1} b^{-1})^{n-1}] a^{r_k + m \ell + 1} \\
&=& b^{s_0} a^{r_1} \cdots b^{s_{k-2}}
a^{r_{k-1} - m + 1} b^{-1} (a^{-m+1} b^{-1})^{n-2} a^{r_k + m \ell + 1}\\
&=&
b^{s_0} a^{r_1} \cdots b^{s_{k-2}-1} [ba^{-1} ]
a^{r_{k-1} - m + 2} b^{-1} (a^{-m+1} b^{-1})^{n-2} a^{r_k + m \ell + 1},
\end{eqnarray*}
so that in case $r_{k-1} < m-1$ (and $s_{k-2} > 0$), we may repeat the procedure. Continuing this way,
one easily convinces oneself that, unless
\begin{equation}\label{malos}
s_{k-1} = \ldots = s_{k-(n-1)} = 1, \quad s_{k-n} > 0, \quad r_{k-1}  = \ldots = r_{k-(n-1)} = m-1,
\end{equation}
the expression for $w$ above may reduced to one of the form
$$w = b^{s_0 }\cdots b^{s_i}a^{-1} \overline{\omega}$$
for certain $s_i > 0$, $i < k$, and a non-positive word $\overline{\omega}$. As $i < k$, the induction hypothesis applies
to $b^{s_0 }\cdots b^{s_i}a^{-1}$, which is hence negative, and so does $w$. Assume otherwise that (\ref{malos}) holds.
Then since $b (a^{m-1}b)^{n-1} = a$, replacing $(a^{m-1} b)^{n-1}$ by $b^{-1} a$ and canceling $b^{-1}$, we obtain a new
expression for $w$ of the form
$$w = b^{s_0} \cdots b^{s_{k-n} - 1} a^{r_k + m \ell + 1},$$
which contradicts the minimality of the length of the normal form. This closes the proof.

\vsp\vsp

Step II of the proof of Theorem \ref{thm:torus} can be established via several approaches.
Here, we chose the dynamical one, based on the fact that $G_{m,n}$ embeds into
$\widetilde{\mathrm{PSL}}(2,\mathbb{R})$. To see this, let us first come back
to the presentation
$$G_{m,n} \esp = \esp \big\langle c,d \!: c^{m} = d^n \big\rangle,$$
which exhibits $G_{m,n}$ as a central extension of the group
$$\overline{G}_{m,n} = \big\langle \bar{c},\bar{d} \!: \bar{c}^m = \bar{d}^n =id \big\rangle.$$
A concrete realization of $\overline{G}_{m,n}$ inside $\mathrm{PSL}(2,\mathbb{R})$
arises when identifying $\bar{c}$ to the circle rotation of angle $\frac{2\pi}{m}$,
and $\bar{d}$ to an hyperbolic rotation of angle $\frac{2\pi}{n}$ centered at a point different
from the origin in such a way that, if we let $p_0 \!:=\! p,
p_1 \!:=\! \bar{c}(p), \ldots, p_{m-1} \!:=\! \bar{c}^{m-1}(p)$ and $q_0 \!:=\! p, q_1 \!:=\! \bar{d}(p),
\ldots, q_{n-1} \!:=\! \bar{d}^{n-1}(p)$ for a certain $p \in \clo$, we have that all the points $q_i$'s lie 
between $p_0$ and $p_1$, and $q_{n-1}=p_1$.
This realization allows embedding $G_{m,n}$ into $\widetilde{\mathrm{PSL}}(2,\mathbb{R})$
by identifying $c \!\in\! G_{m,n}$ to the lifting of $\bar{c}$ to the real line given by
$x \mapsto x + \frac{2\pi}{m}$, and $d$ to the unique lifting of
$\bar{d}$ to the real line satisfying \esp $x \leq \bar{d}(x) \leq x + 2\pi$
\esp for all $x \in \mathbb{R}$. (Actually, the arguments
given so far only show that the above identifications induce a group
homomorphism from $G_{m,n}$ into $\widetilde{\mathrm{PSL}}(2,\mathbb{R})$,
and the injectivity follows from the arguments given below.)

The dynamics of the action of $\overline{G}_{m,n}$ on the circle is illustrated
in Figure 10. Passing to the generators $a,b$, we have that
\esp $\bar{b} =\! \bar{c}^{-(m-1)} \bar{d} = \bar{c} \bar{d}$
\esp is a parabolic M\"obius transformation fixing $p_1$, while $\bar{a} = \bar{c}$.
Using this, we next proceed to show that no element $w$ in
$\langle a,b \rangle^+ \subset G_{m,n}$ represents the identity.
By taking inverses, this will imply that no element in
$\langle a^{-1}, b^{-1} \rangle^+$ represents the
identity, thus completing the proof.

\vspace{0.35cm}


\beginpicture

\setcoordinatesystem units <0.78cm,0.78cm>

\circulararc 360 degrees from 3 0
center at 0 0

\circulararc 85 degrees from -1.8 2.1
center at -1.34 3.4

\plot -0.09 2.82 -0.11 2.6 /
 \plot -0.09 2.82 -0.25 2.67 /

\circulararc 175 degrees from 1 3
center at 0 3
\plot -1 3 -0.86 3.3 /
\plot -1 3 -1.12 3.3 /
\put{$\bar{c} = \bar{a}$} at 0 4.3

\circulararc 175 degrees from 2.6 1.8
center at 1.9 2.34
\plot 1.25 2.9 1.5 3 /
\plot 1.25 2.9 1.4 3.2 /
\put{$\bar{c}$} at 2.7 3.3

\circulararc -175 degrees from -2.6 1.8
center at -1.9 2.34
\plot -2.6 1.8 -2.88 1.95 /
\plot -2.6 1.8 -2.58 2.05 /
\put{$\bar{c}$} at -2.6 3.3

\put{$\bar{d}$} at -0.8 1.75

\put{$p = p_0$} at 0.1 3.2
\put{$p_1$} at -2.05 2.55
\put{$p_2$} at -3.1 1.1
\put{$p_{m-1}$} at 2.07 2.8
\put{$p_{m-2}$} at 3.45 1.1
\put{$p_{m-3}$} at 3.44 -1.1
\put{$p_3$} at -3.1 -1.1
\put{$p_{m-4}$} at 2.4 -2.6

\put{$\bullet$} at 0 3
\put{$\bullet$} at -2 2.2
\put{$\bullet$} at -2.9 0.75
\put{$\bullet$} at 2 2.2
\put{$\bullet$} at 2.9 0.75
\put{$\bullet$} at -2.9 -0.75
\put{$\bullet$} at 2.9 -0.75
\put{$\bullet$} at 2 -2.24

\put{} at -9.2 3.5
\put{Figure 10} at 0 -3.7

\endpicture


\vspace{0.46cm}

We begin by writing $w$ in the form
$$w = b^{s_0} a^{r_1} b^{s_1} \cdots a^{r_k} a^{m \ell}, \quad \ell \geq 0,$$
with the corresponding restrictions on the exponents. Here, we may assume that
no expression $(ba^{m-1})^{n-1}b$ appears, since otherwise we may replace it by $a$.

Assume that $w$ is not a power of $a$, and let us consider its reduction
$$\bar{w} = \bar{b}^{s_0} \bar{a}^{r_1} \bar{b}^{s_1} \cdots \bar{a}^{r_k}
\in \widetilde{\mathrm{PSL}}_2(\mathbb{R}).$$
Using $\bar{b} = \bar{a} \bar{d}$ and simplifying $\bar{a}^m = id$, we may rewrite this in the form
$$\bar{w} = \bar{d}^{s'_0} \bar{a}^{r'_1} \bar{d}^{s'_1} \cdots \bar{a}^{r'_{k'}}
\in \widetilde{\mathrm{PSL}}_2(\mathbb{R}),$$
with similar restrictions on the exponents $r_i',s_i'$. What is crucial here is that
the fact that no expression \, $(ba^{m-1})^{n-1}b$ \, appears in the original form
implies that this new expression is nontrivial, as it can be easily checked.
(Indeed, no cancellation $\bar{d}^n = id$ will be performed.)

Unless $\bar{w}$ is a power of $\bar{d} \bar{a}$, we may conjugate it to either
some $\bar{w}' \in \langle \bar{a}, \bar{b} \rangle^+$ beginning and finishing by
$\bar{a}$ and so that all the exponents of $\bar{a}$ lie in $\{1,\ldots,m-1\}$,
or to some $\bar{w}'' \in \langle \bar{a}, \bar{d} \rangle^+$ beginning and
finishing with $\bar{d}$ with the same restriction on the exponents of
$\bar{a}$. An easy ping-pong type argument then shows that
$\bar{w}' \big( ]p_0,p_1[ \big) \!\subset ]p_1,p_0[$ and
$\bar{w}'' \big( ]p_1,p_0[ \big) \!\subset ]p_0,p_1[$,
hence $\bar{w}' \neq id$ and $\bar{w}'' \neq id$.

Thus, to conclude the proof, we need to check that neither $a$ nor $da$ are torsion elements.
That $da$ has infinite order follows from that $\bar{d} \bar{a}$ sends $[p_0,p_1]$
into the strict subinterval $[p_0,\bar{d}(p_2)]$, hence no iterate of it can equal the
identity. Finally, to see that $a$ also has infinite order, just note that it identifies with the
translation by $\frac{2\pi}{m}$ in $\widetilde{\mathrm{PSL}}(2,\mathbb{R})$.


\vspace{0.3cm}

\index{Order!isolated}
\noindent{\bf Some other examples.} The search for more examples of finitely-generated positive cones in
groups with infinitely many left-orders has become a topic of much activity over the last years. The examples given
above as well as the techniques used in proofs have been pursued in three directions. First, there is the close relation
with Dehornoy-like orders in which the previous examples fit, as described in \cite{pos-1} and later in \cite{ito-1}. (See
also Remark \ref{rem:dehornoy-like}.) Second, there is an approach based on partial cyclic amalgamation, which is
fully developed in \cite{ito-2}. This allows iterative implementation, thus establishing for instance that the groups
$$G_{m_1,m_2,\ldots,m_n}
:= \big\langle a_1, \ldots, a_n \!: a_1^{m_1} = a_2^{m_2} = \ldots = a_n^{m_n} \big\rangle$$
do admit finitely-generated positive cones. This approach was somewhat complemented in \cite{ito-3}; however, the orders constructed therein
are only ensured to be isolated, and knowing whether their positive cones are finitely-generated remains an interesting question. Finally,
there is a more combinatorial approach starting from group presentations introduced in \cite{deh}. Roughly, in case these presentations
have a {\em triangular form}, finitely-generated positive cones naturally appear. As a random example, we can mention that the groups
$$H_{m,n} := \big\langle a,b,c \!:  a = ba^2 (b^2a^2)^m c, b = c (ba^2)^n ba \big\rangle$$
fall in this category.

We do not pursue this nice subject here; we just refer the reader to the works mentioned above for the announced results and further
developments (see also Remark \ref{r:mann-rivas}). Nevertheless, let us mention that none of these approaches has provided 
a new proof of Dehornoy's theorem concerning
the $D$-order on $\mathbb{B}_n$ for $n \geq 4$. This issue seems to be beyond the scope of these methods.


\chapter{ORDERABLE GROUPS AS DYNAMICAL OBJECTS}

\section{H\"older's Theorem}
\label{section-holder}

\index{H\"older's theorem}
\index{Order!Archimedean}
\hspace{0.45cm} The results of this section --essentially due to H\"older-- are
classical and perhaps correspond to the most beautiful elementary theorems of the
theory. They characterize group left-orders satisfying an Archimedean  type
property: the underlying ordered group must be ordered
isomorphic to a subgroup of $(\mathbb{R},+)$. For the statement, 
a left-order $\precede$ on a group $\Gamma$ will said to be
{\bf {\textit Archimedean}} if for all $g, h$ in $\Gamma$ such that $g \!\neq\! id$,
there exists $n \!\in\! \mathbb{Z}$ satisfying $g^n \!\succ\! h.$\\

\begin{thm} \textit{Every group endowed with an Archimedean left-order
is order-isomorphic to a subgroup of \esp $(\mathbb{R},+)$.}
\label{completa}
\end{thm}

\vsp

H\"older proved this theorem under the extra assumption that the group is Abelian.
However, his arguments work verbatim without this hypothesis but assuming that
the left-order is bi-invariant. That this hypothesis is also superfluous
was first remarked by Conrad in \cite{conrad}.

\vsp\vsp

\begin{lem} {\em Every Archimedean left-order on a group is bi-invariant.}
\label{lema-conrad}
\end{lem}\index{Order!Archimedean} 

\noindent{\bf Proof.} Let $\preceq$ be an Archimedean left-order on a group
$\Gamma$. We need to show that its positive cone is a normal semigroup.

Suppose that $g \in P^+_{\preceq}$ and $h \in P^{-}_{\preceq}$ are such that
$hgh^{-1} \notin P^+_{\preceq}$. Let $n$ be the smallest positive integer
for which $h^{-1} \prec g^n$. Since \esp $hgh^{-1} \prec id$, \esp we have
\esp $h^{-1} \prec g^{-1} h^{-1} \prec g^{n-1},$ \esp which contradicts
the definition of $n$. We thus conclude that $P_{\preceq}^+$
is stable under conjugation by elements in $P_{\preceq}^{-}$.

Assume now that $g, h$ in $P_{\preceq}^+$ verify \esp
$hgh^{-1} \notin P_{\preceq}^+$. \esp In this case, $hg^{-1}h^{-1} \succ id$,
and since $h^{-1} \in P_{\preceq}^{-}$, the first part of the proof yields
\esp $h^{-1}(hg^{-1}h^{-1})h \in P_{\preceq}^+$, \esp that is,
$g^{-1} \!\in P_{\preceq}^+$, which is absurd. Hence, $P_{\preceq}^+$
is also stable under conjugation by elements in $P_{\preceq}^{+}$,
which concludes the proof. $\hfill\square$

\vspace{0.23cm}

\begin{small}\begin{ejer}
Prove the preceding lemma by using dynamical realizations (see \S \ref{general-3}). 
More precisely, show that the dynamical realization of every Archimedean left-order 
on a countable group is a subgroup of $\mathrm{Homeo}_+(\mathbb{R})$ that acts  
freely on the line (compare Example \ref{ex-free-action}).
\end{ejer}\end{small}\index{Order!Archimedean} 

\vsp\vsp

\noindent{\bf Proof of Theorem \ref{completa}.} Let $\Gamma$ be a group endowed
with an Archimedean left-order $\precede$. By Lemma \ref{lema-conrad}, this order 
is bi-invariant. Fix a positive element $f \in \Gamma$,
and for each $g \in \Gamma$ and each $p \in \mathbb{N}$, consider the
unique integer $q=q(p)$ such that $\hspace{0.15cm} f^q \precede g^p \prec f^{q+1}$.

\vspace{0.25cm}

\noindent{\underline{Claim (i).}} The sequence \esp $\big( q(p) / p \big)$
\esp converges to a real number as \esp $p$ \esp goes to infinity.

\vspace{0.15cm}

Indeed, if $\hspace{0.15cm} f^{q(p_1)} \precede g^{p_1} \prec
f^{q(p_1)+1} \hspace{0.15cm}$ and $\hspace{0.15cm} f^{q(p_2)}
\precede g^{p_2} \prec f^{q(p_2)+1}, \hspace{0.15cm}$ then
the bi-invariance of $\preceq$ yields
$$f^{q(p_1)+q(p_2)} \precede g^{p_1+p_2} \prec f^{q(p_1)+q(p_2)+2}.$$
Therefore,
$\hspace{0.15cm} q(p_1) + q(p_2) \leq q(p_1+p_2) \leq
q(p_1)+q(p_2)+1. \hspace{0.15cm}$ The convergence of the sequence
$(q(p) / p)$ then follows from Exercise \ref{subaditivo}.

\vspace{0.25cm}

\noindent{\underline{Claim (ii).}} The map $\phi: \Gamma
\rightarrow (\mathbb{R},+)$ is a group homomorphism.\\

\vspace{0.15cm}

Indeed, let $g_1,g_2$ be arbitrary elements in $\Gamma$. Suppose that
$g_1g_2 \precede g_2g_1$ (the case where $g_2g_1 \precede g_1g_2$ is analogous).
Since $\preceq$ is bi-invariant, if \esp $f^{q_1} \precede g_1^p \prec f^{q_1 +1}$
\esp and \esp $f^{q_2} \precede g_2^p \prec f^{q_2 +1}$, \esp then
$$\hspace{0.15cm} f^{q_1+q_2} \precede g_1^p g_2^p \precede
(g_1g_2)^p \precede g_2^p g_1^p \prec f^{q_1 + q_2 + 2}{}_.$$
From these relations one concludes that
$$\phi(g_1) + \phi(g_2) = \lim\limits_{p \rightarrow \infty}
\frac{q_1+q_2}{p} \leq \phi(g_1g_2) \leq \lim\limits_{p
\rightarrow \infty} \frac{q_1+q_2 +1}{p} = \phi(g_1)+\phi(g_2),$$
and therefore $\phi(g_1g_2) = \phi(g_1) + \phi(g_2).$

\vspace{0.25cm}

\noindent{\underline{Claim (iii).}} The homomorphism $\phi$ is one-to-one
and order-preserving.\\

\vspace{0.15cm}

That $\phi$ is order-preserving (in the sense that 
$\phi(g_1) \leq \phi(g_2)$ if $g_1 \precede g_2$) follows from the definition.
To show injectivity, first note that $\phi(f)=1.$ Let $h \in \Gamma$ be such 
that $\phi(h)=0$. Assume that $h \neq id$. Since $\preceq$ is Archimedean, 
there exists $n \!\in\! \mathbb{Z}$ such that $h^n \succeq f$. Consequently,
$0 = n \phi(h) = \phi(h^n) \geq \phi(f) = 1,$ which is absurd.
Therefore, if $\phi(h) = 0$, then $h = id$. $\hfill\square$

\vspace{0.1cm}

\begin{small} \begin{ejer}
Let $(a_n)_{n \in \mathbb{Z}}$ be an integer-indexed sequence of real
numbers. Assume that there exists a constant $C \in \mathbb{R}$ such
that, for all $m,n$ in $\mathbb{Z}$,
\begin{equation}
|a_{m+n} - a_m - a_n| \leq C.
\label{subaditividad}
\end{equation}
Show that there exists a unique $\theta \in \mathbb{R}$
such that the sequence $ \big( |a_n - n\theta| \big)$ is bounded.
Check that this number $\theta$ is equal to the limit of the sequence
$(a_n / n)$ as $n$ goes to \esp $\pm \infty$ \esp (in particular,
this limit exists).

\noindent{\underline{Hint.}} For each $n \!\in\! \mathbb{N}$ let
$I_n := \big[ (a_n - C)/n,(a_n + C)/n \big]$. Check that $I_{mn}$
is contained in $I_n$ for every $m,n$ in $\mathbb{N}$. Conclude that
$I := \bigcap_{n \in \mathbb{N}} I_n$ is nonempty (any $\theta$ in
$I$ satisfies the desired property).
\label{subaditivo}
\end{ejer}

\begin{ex} \label{ex-free-action} \index{Order!Archimedean} 
Groups acting freely on the real line are examples of groups admitting 
Archimedean left-orders. Indeed, from such an action one may define
$\precede$ on $\Gamma$ by letting $g \prec h$ if $g(x) < h(x)$ for some
(equivalently, for all) $x \!\in\! \mathbb{R}$. This order relation is
total, and using the fact that the action is free, one readily shows
that it is Archimedean (as well as bi-invariant).

Note that, by the proof of Theorem \ref{completa}, the left-order $\preceq$ above
induces an embedding $\phi$ of $\Gamma$ into $(\mathbb{R},+)$. If $\phi(\Gamma)$ is
isomorphic to $\mathbb{Z}$, then the action of $\Gamma$ is conjugate to the action
by integer translations. Otherwise, unless $\Gamma$ is trivial, $\phi(\Gamma)$ is
dense in $(\mathrm{R},+)$. For each point $x$ in the line, we may then define
$$\varphi(x) = \sup \big\{ \phi(h) \in \mathbb{R} \!: h(0) \leq x \big\}.$$
It is easy to see that $\varphi \!: \mathbb{R} \rightarrow \mathbb{R}$
is a non-decreasing map. Moreover, it satisfies \esp
$\varphi(h(x)) = \varphi(x) + \phi(h)$ \esp for all $x \!\in\! \mathbb{R}$
and all $h \!\in\! \Gamma$. Finally, $\varphi$ is continuous, as otherwise
the set $\mathbb{R} \setminus \varphi(\mathbb{R})$ would be a nonempty
open set invariant under the translations of $\phi(\Gamma)$, which is impossible.
In summary, every free action on the line is (continuously) semiconjugate to an 
action by translations.
\end{ex} \end{small}


\section{The Conrad Property}
\label{general-Conrad}

\subsection{The classical approach revisited}
\label{classic-conrad}

\hspace{0.45cm} A left-order $\preceq$ on a group $\Gamma$ is said to
be {\bf{\em Conradian}} (a {\bf{\em $C$-order}}, for short) if for all positive
elements $f,g$, there exists $n \!\in\! \mathbb{N}$ such that $fg^n \succ g$. 
Groups admitting a $C$-order are called $C$-orderable.
\index{Order!Conradian}

Bi-invariant left-orders are Conradian, as $n=1$ works in the preceding inequality 
for bi-orders. In this direction, it is quite remarkable that one may actually take
\esp $n \!=\! 2$ \esp in the general definition above, as the next proposition
shows. The nice proof we give below, due to Jim\'enez, is taken from \cite{leslie}.

\vspace{0.1cm}

\begin{prop} \label{n=2} {\em If $\preceq$ is a Conradian order on
a group, then $fg^2 \succ g$ holds for all positive elements $f,g$.}
\label{lesla}
\end{prop}
\index{Order!Conradian}

\noindent{\bf Proof.} Suppose that two positive elements $f,g$ for a left-order $\preceq'$
on a group $\Gamma$ are such that $fg^2 \precede' g$. Then $(g^{-1}fg) g \precede' id$, and
since $g$ is a positive element, this implies that $g^{-1}fg$ is negative, and therefore
$fg \prec' g$. Now for the positive element $h := fg$ and every $n \in \mathbb{N}$, one has
\begin{multline*}
f h^{n} = f (fg)^n = f (fg)^{n-2} (fg) (fg) \prec'
f (fg)^{n-2} (fg) g \\
= f (fg)^{n-2} fg^2 \precede' f (fg)^{n-2} g
= f (fg)^{n-3} fg^2 \precede' f (fg)^{n-3}g \precede' \ldots\\
\precede' f (fg) g = f fg^2 \precede' fg = h.
\end{multline*}
This shows that $\preceq'$ does not satisfy the Conrad property. $\hfill\square$

\vsp\vsp\vsp\vsp\vsp

The following is an easy (but important) corollary to the
previous proposition, and we leave its proof to the reader.
(Compare Exercise \ref{bi-inva}.)

\vspace{0.1cm}

\begin{cor} \label{cerrado}
{\em For every left-orderable group, the subspace $\mathcal{CO}(\Gamma)$
of Conradian orders is closed inside the space of left-orders.
Moreover, this subspace is invariant under the conjugacy action.}
\end{cor}

\vsp

Perhaps the most important theorem concerning $C$-orderable groups is the next one.
The direct implication is due to Conrad \cite{conrad}; the converse is due to Brodski
\cite{brodski}, yet it was independently rediscovered by Rhemtulla and Rolfsen
\cite{RR}. We postpone the proof of the first part, and for the second we offer
an elementary one taken from \cite{order}. Recall that a group is said to be
{\bf{\em locally indicable}} if each nontrivial finitely-generated subgroup
admits a nontrivial homomorphism into $(\mathbb{R},+)$.
\vsp\vsp
\index{Group!locally indicable}
\begin{thm} \label{st-th}
{\em A group $\Gamma$ is $C$-orderable if and only if it is locally indicable.}
\end{thm}
\index{Order!Conradian}

\vsp\vsp

To show that local indicability implies $C$-orderability (the converse
will be proved in \S \ref{conrad-general}), we will need the following
lemma, the proof of which is left to the reader. (Compare \S \ref{general-2}.)

\vsp\vsp

\begin{lem} {\em A group $\Gamma$ is $C$-orderable if and only if for every finite family
$\mathcal{G}$ of elements in $\Gamma \setminus \{ id \}$, there exists a choice of exponents
$\epsilon \!: \mathcal{G} \to \{-1,+1\}$ such that $id$ does not belong to the smallest
subsemigroup $\langle \! \langle \mathcal{G} \rangle \! \rangle$ satisfying:

\vsp

\noindent -- It contains all the elements \esp $g^{\epsilon(g)}$, with $g \in \mathcal{G}$;

\vsp

\noindent -- For all $f,g$ in the semigroup,
the element \esp $g^{-1}fg^{2}$ \esp also belongs to it.}
\end{lem}

\vsp\vsp

\noindent{\bf Local indicability implies $C$-orderability.} We need to
check that every locally indicable group $\Gamma$ satisfies the condition of
the preceding lemma. Let $\{g_1,\ldots,g_k\}$ be a finite family of elements
in $\Gamma$ different from the identity. By hypothesis, there is a nontrivial
homomorphism $\phi_1\!: \langle g_1,\ldots,g_k \rangle \rightarrow (\mathbb{R},+)$.
Let \esp $i_1,\ldots,i_{k'}$ \esp be the indices (if any) such that $\phi_1 (g_{i_{j}}) = 0$.
Again by hypothesis, there exists a nontrivial homomorphism
$\phi_2\!: \langle g_{i_1},\ldots,g_{i_{k'}} \rangle \rightarrow (\mathbb{R},+)$. Letting
\esp $i_1',\ldots,i_{k''}'$ \esp be the indices in $\{i_1,\ldots,i_{k'}\}$ for which
$\phi_2 (g_{i_j'}) = 0$, we may choose a nontrivial homomorphism
$\phi_3\!: \langle g_{i_1'},\ldots,g_{i_{k''}'} \rangle \rightarrow (\mathbb{R},+)$... Note that
this process must finish in a finite number of steps (indeed, it stops in at most $k$ steps).
Now, for each $i \!\in\! \{1,\ldots,k\}$, choose the (unique) index $j(i)$ such that $\phi_{j(i)}$
is defined at $g_i$ and $\phi_{j(i)} (g_i) \neq 0$, and let $\epsilon_i := \epsilon(g_i)\!\in\! \{-1,+1\}$ be so that
$\phi_{j(i)} (g_i^{\epsilon_i}) > 0$. We claim that this choice of exponents $\epsilon_i$ is ``compatible".
Indeed, for every index $j$ and every $f,g$ for which $\phi_j$ are defined, one has
\esp $\phi_j(f^{-1} g f^2) = \phi_j (f) + \phi_j (g).$ \esp
\hspace{0.07cm} Therefore, \esp $\phi_1(h) \geq 0$ \esp for every \esp
$h \!\in\! \langle \!\langle g_1^{\epsilon_1},\ldots,g_k^{\epsilon_k} \rangle \!\rangle$.
\esp Moreover, if \esp $\phi_1(h) = 0$, \esp then $h$ actually belongs to \esp
$\langle\!\langle g_{i_1}^{\epsilon_{i_1}},\ldots,g_{i_{k'}}^{\epsilon_{i_{k'}}} \rangle\!\rangle$.
\esp In this case, the preceding argument shows that \esp $\phi_2(h) \!\geq\! 0$, \esp
with equality if and only if
$h \!\in\! \langle\!\langle g_{i_1´}^{\epsilon_{i_1'}},
\ldots,g_{i_{k''}'}^{\epsilon_{i_{k''}'}} \rangle\!\rangle$...
Continuing in this way, one concludes that $\phi_j (h)$ must be strictly
positive for some index $j$. Thus, the element $h$ cannot be equal to
the identity, and this finishes the proof.

\vspace{0.3cm}

If a group $\Gamma$ contains a normal subgroup $\Gamma_*$ so that
both $\Gamma_*$ and $\Gamma / \Gamma_*$ are locally indicable, then $\Gamma$
itself is locally indicable. Equivalently, the extension of a $C$-orderable
group by a $C$-orderable group is $C$-orderable. This is made more 
precise in the next exercise. 

\begin{small}
\begin{ejer} \label{ejercicito}
Let $(\Gamma,\preceq)$ be a $C$-ordered group, and let $\Gamma_*$ be a convex
subgroup. Show that for any $C$-order $\preceq_*$ of $\Gamma_*$, the
extension of $\preceq_*$ by $\preceq$ is still Conradian. In particular,
every left-order obtained from a $C$-left-order by flipping a convex
subgroup is Conradian (see \S \ref{section-convex-extension}).
\end{ejer}

\begin{ejem} \label{1-relator-LI}
A remarkable theorem independently obtained by Brodski \cite{brodski}
and Howie \cite{howie} asserts that torsion-free, 1-relator groups are locally
indicable. Also, all knot groups in $\mathbb{R}^3$ are locally indicable
(see \cite[Lemma 2]{howie-short}).
\end{ejem}\end{small}

\vsp\vsp

\index{Bergman-Thurston example}
\noindent{\bf Examples of left-orderable, non $C$-orderable groups.} Only a few examples
are known. Historically, the first was exhibited (in a slightly different context) by Thurston
\cite{Th}, and rediscovered some years later by Bergman \cite{bergman-1}. It corresponds to
the lifting to $\widetilde{\mathrm{PSL}}(2,\mathbb{R})$ of the $(2,3,7)$-triangle group, and
has the presentation \esp
$$\Gamma = \big\langle f,g,h \!: f^2 = g^3 = h^7 = fgh \big\rangle.$$
Left-orderability follows from that $\widetilde{\mathrm{PSL}}(2,\mathbb{R})$ is a subgroup
of $\mathrm{Homeo}_+(\mathbb{R})$. The fact that $\Gamma$ is not $C$-orderable is a
consequence of the fact that it has no nontrivial homomorphism into $(\mathbb{R},+)$, which 
may be easily deduced from the presentation above. Actually, $\Gamma$ is the $\pi_1$ of an
homological sphere, and this was the motivation of Thurston for dealing with this group
in his generalization of the famous Reeb stability theorem for codimension-1
foliations. We strongly recommend reading \cite{Th} for all of this; 
see also Exercise \ref{ejer-tst} further on.

Below we elaborate on a different and quite important example,
namely braid groups $\mathbb{B}_n$ for $n \geq 5$. Another
example is the lifting $\widetilde{\mathrm{G}}$ of Thompson's
group G to the real line; see \cite{CFP} for more details.

\begin{small}
\begin{ex} \label{sobre-b3-b4}
The braid groups $\mathbb{B}_3$ and $\mathbb{B}_4$ are locally
indicable. For $\mathbb{B}_3$, this may be easily deduced from the exact sequence
$$0 \longrightarrow [\mathbb{B}_3,\mathbb{B}_3] \sim \mathbb{F}_2 \longrightarrow \mathbb{B}_3 \longrightarrow
\mathbb{B}_3 / [\mathbb{B}_3,\mathbb{B}_3] \sim \mathbb{Z} \longrightarrow 0,$$
where the isomorphism $[\mathbb{B}_3,\mathbb{B}_3] \!\sim\! \mathbb{F}_2$ may be shown by looking the action
on the circle of $\mathbb{B}_3 \!\sim\! \widetilde{\mathrm{PSL}}(2,\mathbb{Z})$, and
$\mathbb{B}_3 / [\mathbb{B}_3,\mathbb{B}_3] \sim \mathbb{Z}$ appears by taking ``total exponents''.
For $\mathbb{B}_4$, there is an exact sequence
$$0 \longrightarrow \mathbb{F}_2 \longrightarrow \mathbb{B}_4 \longrightarrow
\mathbb{B}_3 \longrightarrow 0.$$
Here, the homomorphism from $\mathbb{B}_4$ to $\mathbb{B}_3$ is the one that sends $\sigma_1$
and $\sigma_3$ to $\sigma_1$, and $\sigma_2$ to $\sigma_2$. Its kernel
is generated by $\sigma_1 \sigma_3^{-1}$ and
$\sigma_2 \sigma_1 \sigma_3^{-1} \sigma_2^{-1}$. To show that these
elements are free generators, one may consider the homomorphism
$\phi \!: \mathbb{B}_4 \rightarrow \mathrm{Aut} (\mathbb{F}_2)$ defined by
\esp $\phi(\sigma_1)(a) \!:=\! a, \esp
\phi(\sigma_1)(b) \!:=\! ab, \esp
\phi(\sigma_2)(a) \!:=\! b^{-1}a, \esp
\phi(\sigma_2)(b) \!:=\! b, \esp
\phi(\sigma_3)(a) \!:=\! a, \esp
\phi(\sigma_3)(b) \!:=\! ba,$ \esp
and note that \esp $\phi(\sigma_1 \sigma_3^{-1})$ \esp (resp.
$\phi(\sigma_2 \sigma_1 \sigma_3^{-1} \sigma_2^{-1})$) \esp
is the conjugacy by $a$ (resp. $b^{-1}a$).
\end{ex} 

\begin{ejer} The homomorphism from $\mathbb{B}_4$ to $\mathbb{B}_3$ referred to in the preceding example induces a 
homomorphism from the symmetric group $S_4$ to $S_3$. Check that the latter homomorphism arises as follows: If $S_4$ acts 
by permutations of a set $S$ consisting of four objects $s_1,s_2,s_3,s_4$, then for each element in $S_4$, the induced element in $S_3$ 
acts by permuting the elements $t_1,t_2,t_3$ of the set $T$ consisting of the three different manners of splitting $S$ into two pairs. 
(More specifically, $t_1 := \{(s_1,s_4),(s_2,s_3)\}$, $t_2 := \{(s_1,s_3),(s_2,s_4)\}$, and $t_3 := \{(s_1,s_2),(s_3,s_4)\}$.) 
See \cite{benoit} for more on this beautiful observation.
\end{ejer}\end{small}

\vsp

\noindent{\bf Incompatibility between bi-orders on $\mathrm{P}\mathbb{B}_n$
and left-orders on $\mathbb{B}_n$.} In contrast to $\mathbb{B}_3$ and $\mathbb{B}_4$,
the groups $\mathbb{B}_n$ fail to be locally indicable for $n \!\geq\! 5$. Indeed, for
$n \!\geq\! 5$, the commutator subgroup $[\mathbb{B}_n,\mathbb{B}_n]$ is (finitely-generated 
and) {\bf {\em perfect}} ({\em i.e.}, it coincides with its own commutator subgroup),
as shown below.

\begin{small}\begin{ex} \label{n geq 5}
As is well-known (and easy to check), the commutator subgroup
$[\mathbb{B}_n,\mathbb{B}_n]$ is generated by the elements of the form $\s_{i,j} := \s_i \s_j^{-1}$. 
Also, recall that all the generators $\s_i$ of $\mathbb{B}_n$ are conjugate between them.
Indeed, letting $\Delta := \s_1 \s_2 \cdots \s_{n-1}$, one readily checks that
$\s_i \Delta = \Delta \esp \s_{i-1}$. Thus, for all $i \!\in\! \{1,\ldots,n-3\}$,
the equality
$$\s_{i,i+2} =
(\s_{i} \s_{i+1})^{-1} \esp [\s_{i,i+2}, \s_{i+1,i}] \esp (\s_i \s_{i+1})$$
shows that $\s_{i,i+2}$ belongs to $\mathbb{B}_n''$. We will close the proof
by showing that, for $n \!\geq\! 5$, the normal closure $H$ (in $\mathbb{B}_n$) of the
family of elements $\s_{i,i+2}$ (equivalently, of each $\s_{i,i+2}$)
is $\mathbb{B}_n'$. To do this, note that $\s_{i,j}$ and $\s_{i,j'}$ are
conjugate whenever $\{j,j'\} \cap \{ i-1,i+1\} = \emptyset$. Indeed, one
may perform a conjugacy between $\s_j$ and $\s_j'$ as above but inside the subgroup
$\mathbb{B}_{n-2}' \!\subset\! \mathbb{B}_n$ consisting of braids for which the $i$ and $i+1$ strands
remain ``fixed''; such a conjugacy does not change $\s_i$. Therefore, $\s_{i,j}$
belongs to $H$ for all $j \notin \{i-1,i+1\}$. Moreover, since for all
$j \notin \{i-1,i,i+1,i+2\}$ (resp. $j \notin \{i-2,i-1,i,i+1\}$),
$$\s_{i,i+1} = \s_{i,j} \s_{j,i+1} \qquad (\mbox{resp. }
\s_{i,i-1} = \s_{i,j} \s_{j,i-1}),$$
the elements $s_{i,i+1}$ and $s_{i,i-1}$ also belong to $H$.
This shows that $H$ coincides with $\mathbb{B}_n'$.

We recommend \cite{M-Rolfsen} for more details on this example, 
as well as generalizations in the context of Artin groups.
\end{ex}\end{small}

A nice consequence of the example above is that the bi-orders on
$\mathrm{P}\mathbb{B}_n$ do not extend to left-orders on $\mathbb{B}_n$
for any $n \geq 5$. (This fact was established, independently,
in \cite{dub} and \cite{RR}.) Indeed, we have the following 
proposition.

\vsp

\begin{prop} \label{listailor}
{\em Let $\Gamma_0$ be a finite-index subgroup of a left-orderable
group $\Gamma$. If $\preceq$ is a left-order on $\Gamma$ whose restriction
to $\Gamma_0$ is Conradian, then $\preceq$ is Conradian.}
\end{prop}

\noindent{\bf Proof.} Let $f \succ id$ and $g \succ id$ be
elements in $\Gamma$. One has $f^m \in \Gamma_0$ and
$g^n \in \Gamma_0$ for some positive $n,m$ smaller
than or equal to the index of $\Gamma_0$ in $\Gamma$.
Hence, $f^m g^{2n} \succ g^n \succ g$. We claim that this implies that
either $fg \succ g$ or $fg^{2n} \succ g$. Otherwise, $g^{-1} fg \prec id$
\, and \, $g^{-1} f g^{2n} \prec id$, thus yielding
$$id \esp \prec \esp g^{-1} f^m g^{2n} \esp = \esp
(g^{-1} f g)^{m-1} (g^{-1} f g^{2n}) \esp \prec \esp id,$$
which is absurd. $\hfill\square$

\vsp\vsp\vsp

\noindent{\bf A criterion of non left-orderability.} Proposition \ref{listailor}
allows to show that certain ``small'' groups cannot be left-ordered. In concrete
terms, we have the following result due to Rhemtulla \cite[Chapter 7]{botto}.

\vsp

\begin{prop} \label{unico-criterio}
{\em Let $\Gamma$ be a finitely-generated group containing a finite-index subgroup
$\Gamma_0$ all of whose left-orders are Conradian. If $\Gamma$ has no nontrivial
homomorphism into $(\mathbb{R},+)$, then $\Gamma$ is not left-orderable.}
\end{prop}

\vsp

Indeed, if $\Gamma$ were left-orderable then, by Proposition \ref{listailor}, every left-order
on it would be Conradian. Since $\Gamma$ is finitely-generated, Theorem \ref{st-th}
would provide us with a nontrivial homomorphism into $(\mathbb{R},+)$.

\begin{small}\begin{ejem} \label{este-funciona2}
In \S \ref{no-torsion}, we introduced the group
$$\Gamma = \big\langle a,b \!: a^2 b a^2 = b, b^2 a b^2 = a \big\rangle,$$
which contains an index-4 Abelian subgroup, namely $\langle a^2,b^2,(ab)^2 \rangle
\sim \mathbb{Z}^3$. From the presentation, it follows that $\Gamma$ admits
no nontrivial homomorphism into $(\mathbb{R},+)$. Since bi-invariant left-orders are
Conradian, Theorem \ref{unico-criterio} implies that $\Gamma$ is not left-orderable.
\end{ejem}\end{small}


\subsection{An approach via crossings}
\label{section-crossings}

\hspace{0.45cm} An alternative --dynamical-- approach to the theory of
Conradian orders has been recently developed in \cite{order,crossings}.
We begin with the definition of the notion of {\bf {\em crossing}}, which
is the most important tool in this approach.\footnote{It should be noted
that an equivalent notion --namely that of {\bf{\em overlapping elements}}--
was introduced by Glass in his dynamical study of lattice-orderable groups
\cite{glass-antiguo}, though no conexion with the Conrad property
is exhibited therein.}
\index{Crossing}
\index{Crossing!reinforced}

\vspace{0.45cm}


\beginpicture

\setcoordinatesystem units <0.8cm,0.8cm>


\putrule from 1.5 -2.5 to 6.5 -2.5
\putrule from 1.5 2.5 to 6.5 2.5
\putrule from 1.5 -2.5 to 1.5 2.5
\putrule from 6.5 -2.5 to 6.5 2.5

\plot
1.5 0.2
6.5 1.96 /


\plot
6.5 0.2
1.5 -2.02 /


\setdots

\plot 1.5 -2.5
6.5 2.5 /

\putrule from 2.32 -1.68 to 5.68 -1.68
\putrule from 2.32 -1.68 to 2.32 1.68
\putrule from 2.32 1.68 to 5.68 1.68
\putrule from 5.68 -1.68 to 5.68 1.68

\put{Figure 11: A reinforced crossing.} at 4 -3.65
\put{} at -5.2 0

\small


\put{$u$} at 1.5 -2.9
\put{$v$} at 6.5 -2.9
\put{$w$} at 4   -2.9
\put{$f^N v$} at 2.8 -2.9
\put{$g^M u$} at 5.2 -2.9
\put{$\bullet$} at 1.5 -2.5
\put{$\bullet$} at 6.5 -2.5
\put{$\bullet$} at 4   -2.5
\put{$\bullet$} at 2.8 -2.5
\put{$\bullet$} at 5.2 -2.5

\put{$f$} at 4.7 -0.9
\put{$g$} at 3.4 1.1

\endpicture


\vspace{0.45cm}

Let $\preceq$ be a left-order on a group $\Gamma$. Following \cite{crossings}, we say
that a 5-tuple $(f,g;u,v,w)$ of elements in $\Gamma$ is a {\bf{\em crossing}} (resp.
{\bf{\em reinforced crossing}}) for $(\Gamma,\preceq)$ if the following conditions
are satisfied:

\vsp

\noindent -- \esp $u \prec w \prec v$;

\vsp

\noindent -- \esp $g^n u \prec v$ \esp and \esp $f^n v \succ u$
\esp for every $n \in \mathbb{N}$ \hspace{0.1cm}
(resp. also $f u \succ u$ \esp and \esp $g v \prec v$);

\vsp

\noindent --  \esp $f^N v \prec w \prec g^M u$ \esp holds for certain $M,N$ in $\mathbb{N}$.

\vsp

Clearly, every reinforced crossing is a crossing. Conversely, if $(f,g;u,v,w)$ is
a crossing, then one easily checks that $(f^{N} g^{M},g^{M} f^{N};f^N w, g^M w, w)$
is a reinforced crossing.

\vsp

\index{Resilient pair}
An equivalent notion to the above ones is that of a {\bf{\em resilient pair}},
namely a 4-uple of group elements $(f,g;u,v)$ satisfying
$$u \prec fu \prec fv \prec gu \prec gv \prec v.$$
Indeed, if $(f,g;u,v,w)$ is a reinforced crossing, then $(f^N,g^M;u,v)$
is a resilient pair for the corresponding exponents $M,N$. Conversely,
if $(f,g;u,v)$ is a resilient pair, then $(f^2,g; u, v, fv)$ is a
reinforced crossing.

\vspace{0.1cm}


\begin{center}
\hspace{-0.8cm}
\includegraphics[scale=0.35]{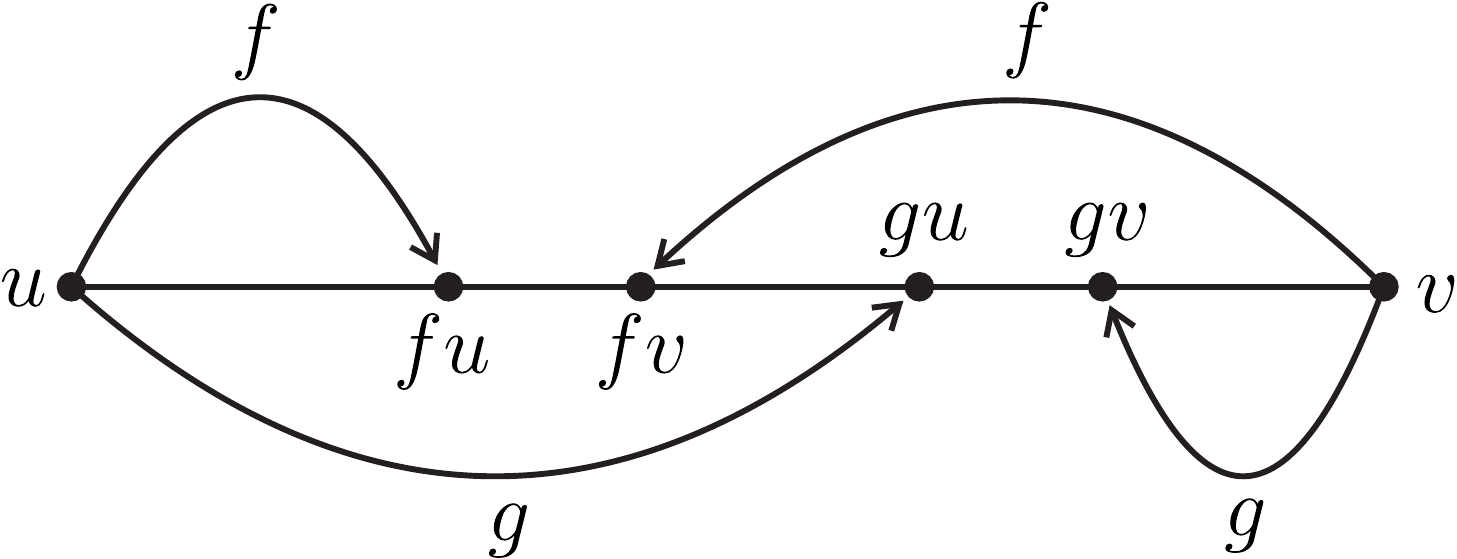}
\end{center}

\vspace{-0.8cm}

\begin{center}
Figure 12: A resilient pair.
\end{center}


\vspace{0.1cm}

\begin{thm} \label{Conrad=noCrossing}
{\em The left-order $\preceq$ is Conradian
if and only if $(\Gamma,\preceq)$ admits
no (reinforced) crossing.}
\end{thm}
\index{Order!Conradian}
\index{Crossing}
\index{Crossing!reinforced}

\noindent{\bf Proof.} Suppose that $\preceq$ is not Conradian, and
let $f,g$ be positive elements such that $fg^n \prec g$ for every
$n \in \mathbb{N}$. We claim that $(f,g;u,v,w)$ is a crossing for
$(\Gamma,\preceq)$ for the choice $u := 1$, $v := f^{-1}g$,
$w := g^2$. Indeed:

\vsp

\noindent -- From $fg^2 \prec g$ one obtains
$g^2 \prec f^{-1} g$, and since $g \succ 1$, this yields \esp
$1 \prec g^2 \prec f^{-1}g,$  \esp that is, \esp $u \prec w \prec v$;

\vsp

\noindent -- From $f g^n \prec g$ it follows that \esp $g^n \prec f^{-1} g$,
\esp that is, \esp $g^n u \prec v$ \esp (for every $n \in \mathbb{N}$);
moreover, since both $f,g$ are positive, we have $f^{n-1} g \succ 1$, and
thus \esp $f^n (f^{-1}g) \succ 1$, \esp that is, \esp $f^n v \succ u$
\esp (for every $n \in \mathbb{N}$);

\vsp

\noindent -- The relation \esp $f (f^{-1}g) = g \prec g^2$ \esp may be read as
\esp $f^N v \prec w$ \esp for $N=1$; finally, the relation \esp
$g^2 \prec g^3$ \esp is \esp \esp $w \prec g^M u$ \esp for $M = 3$.

\vsp\vsp

Conversely, let $(f,g;u,v,w)$ be a crossing for $(\Gamma,\preceq)$
for which for certain $M,N$ in $\mathbb{N}$,
$$f^N v \prec w \prec g^M u.$$
We will prove that $\preceq$ is not Conradian by showing that, for $h := g^Mf^N$ and
$\bar{h} := g^M$, both elements $w^{-1} h w$ and $w^{-1} \bar{h} w$ are
positive, but
$$(w^{-1} h w) (w^{-1} \bar{h} w)^n \prec w^{-1} \bar{h} w,
\quad \mbox{ for all } \esp n \in \mathbb{N}.$$
To do this, first note that \esp $gw \succ w$, \esp as otherwise
$$w \prec g^N u \prec g^N w \prec g^{N-1} w \prec \ldots \prec gw \prec w,$$
which is absurd. Clearly, the inequality \esp $gw \succ w$ \esp implies \esp
$g^M w \succ w,$ \esp hence
\begin{equation}
w^{-1} \bar{h} w = w^{-1} g^M w \succ 1.
\label{zero}
\end{equation}
Moreover, \esp
$h w = g^M f^N w \succ g^M f^N f^N v
= g^M f^{2N} v \succ g^{M} u \succ w,$
\esp thus
\begin{equation}
w^{-1} h w \succ 1.
\label{one}
\end{equation}
Now note that, for every $n \in \mathbb{N}$,
$$h \bar{h}^n w = h g^{Mn} w \prec h g^{Mn} g^M u
= h g^{Mn+M} u \prec h v = g^M f^N v \prec g^M w = \bar{h} w.$$
After multiplying by the left by $w^{-1}$, the last inequality becomes
$$(w^{-1} h w) (w^{-1} \bar{h} w)^n =
w^{-1} h \bar{h}^n w \prec w^{-1} \bar{h} w,$$
as we wanted to check. Together with (\ref{zero}) and (\ref{one}),
this shows that $\preceq$ is not Conradian.
$\hfill\square$

\vsp\vsp

\begin{small}\begin{ejer} Using the characterization of the Conrad property
in terms of resilient pairs, show that the subspace of $C$-left-orders is closed
inside the space of left-orders of a group (see Corollary \ref{cerrado}).
\end{ejer}

\begin{ejer} Using the notion of crossings, give
an alternative proof for Proposition \ref{listailor}.
\end{ejer}

\noindent{\underline{Hint.}} If $(f,g; u,v)$ is a resilient pair, then the same is true for
$(f^n,g^n ; u,v)$, for all $n \geq 1$.

\begin{ejer} Proceed similarly with Proposition \ref{lesla}.

\noindent{\underline{Hint.}} Show that, if $f,g$ are positive elements for which
$fg^2 \prec g$, then $(f,fg;id,fg,g)$ is a crossing for $M=N=2$ (see Figure 13 below).
\end{ejer}

\vspace{0.5cm}


\beginpicture

\setcoordinatesystem units <1.2cm,1.2cm>


\putrule from 1.5 -2.5 to 6.5 -2.5
\putrule from 1.5 2.5 to 6.5 2.5
\putrule from 1.5 -2.5 to 1.5 2.5
\putrule from 6.5 -2.5 to 6.5 2.5

\plot
1.5 1
1.625 1.05296
1.75 1.109
1.875 1.16056
2 1.216
2.125 1.26056
2.25 1.3009
2.375 1.35296
2.5 1.4 /

\plot
2.5 1.4
2.625 1.45296
2.75 1.509
2.875 1.56056
3 1.616
3.125 1.66056
3.25 1.7009
3.375 1.75296
3.5 1.8 /

\plot
3.5 1.8
3.625 1.85296
3.75 1.909
3.875 1.96056
4 2.016
4.125 2.06056
4.25 2.1009
4.375 2.15296
4.5 2.2 /

\setquadratic
\plot
4.5 2.2
4.75 2.329
5 2.5 /

\plot
1.5 -1.1
3.8 -0.18
5 0.2 /

\setlinear



\setquadratic
\plot
1.5 -1.76
5 -1.12
6.5 0.2 /

\setlinear






\setdots

\plot 1.5 -2.5
6.5 2.5 /

\putrule from 1.5 0.2 to 6.5 0.2

\putrule from 5 -2.5 to 5 2.5
\putrule from 1.5 1 to 6.5 1
\putrule from 1.5 -1.12 to 6.5 -1.12

\putrule from 2.35 -0.18 to 3.8 -0.18
\putrule from 2.35 -0.18 to 2.35 -1.68
\putrule from 2.35 -1.68 to 3.8 -1.68
\putrule from 3.8  -0.18 to 3.8 -1.68
\putrule from 2.88 -2.5 to 2.88 -1.12

\put{Figure 13: The $n\!=\!2$ condition.} at 4 -3.4
\put{} at -2.4 0

\small


\put{$id$} at 1.5 -2.8
\put{$g^2$} at 6.5 -2.8
\put{$\bullet$} at 1.5 -2.5
\put{$\bullet$} at 2.88 -2.5
\put{$\bullet$} at 1.5 0.2
\put{$fg^2$} at 1.2 0.2

\put{$f$} at  1.2 -1.78
\put{$fg$} at 1.2 -1.12
\put{$fg$} at 2.88 -2.8
\put{$g$} at  1.2 1

\put{$fg$} at 4.75 -0.07

\put{$g$} at 5 -2.8
\put{$\bullet$} at 1.5 1
\put{$\bullet$} at 1.5 -1.12
\put{$\bullet$} at 1.5 -1.76
\put{$\bullet$} at 5 -2.5
\put{$\bullet$} at 6.5 -2.5

\put{$f$} at 4.5 -1.47
\put{$g$} at 3 1.8

\endpicture


\vspace{0.45cm}

\begin{ex} \label{deh-no}
The Dehornoy left-order $\preceq_{_D}$ on the braid group $\mathbb{B}_n$ (where
$n \!\geq\! 3$) is not Conradian. Indeed, as we next show, $(f,g; u,v,w)
:= (\sigma_2^{-1},\sigma_1,\sigma_2,\sigma_2 \sigma_1, \sigma_2^{-1} \sigma_1)$
is a crossing for $\prec_{_D}$ with $M = N = 1$ (see \cite{NW} for an
alternative argument):

\vsp

\noindent -- It holds that $\sigma_2 \prec_{_D} \sigma_2^{-1} \sigma_1$ is $u \prec_{_D} w$;
moreover, one easily checks that
$\sigma_2 \sigma_1 \succ_{_D} \sigma_1 \succ_{_D} \sigma_2^{-1} \sigma_1$,
hence $w \prec_{_D} v$.

\vsp

\noindent -- For all $k > 0$, we have
$g^k (u) = \sigma_1^k (\sigma_2) \prec_{_D} \sigma_2 \sigma_1 = v$,
where the middle inequality follows from \esp
$\sigma_1^{-1} \sigma_2^{-1} \sigma_1^k \sigma_2 = \sigma_1^{-1} \sigma_1 \sigma_2^k \sigma_1^{-1}
= \sigma_2^k \sigma_1^{-1} \prec_{_D} 1;$ \esp
analogously, for $k \in \mathbb{N}$, we have $f^k(v) = \sigma_2^{-k} (\sigma_2 \sigma_1)
= \sigma_2^{-(k-1)} \sigma_1 \prec_{_D} \sigma_2 \sigma_1$,
where the last inequality follows from
$$\s_1^{-1} \s_2^{k-1} \s_2 \s_1 = \s_2 \s_1^{k} \s_2^{-1} \prec_{_D} id.$$

\vsp

\noindent -- We have
$f(v) = \sigma_2^{-1} (\sigma_2 \sigma_1) = \sigma_1 \succ_{_D} \sigma_2^{-1} \sigma_1 = w$
\esp and \esp $g(u) = \sigma_1 (\sigma_2) \succ_{_D} \sigma_1 \succ_{_D} \sigma_2^{-1} \sigma_1 = w$.
\end{ex}

\begin{ejer} Show that the isolated left-order on the group $G_{m,n}$ constructed in \S \ref{fin-gen} is
not Conradian for $(m,n) \neq (2,2)$.
\end{ejer}

\index{Resilient pair!double}
\index{Crossing!double}
\begin{rem} \label{double-crossing}
The dynamical characterization of the Conrad property should
serve as inspiration for introducing other relevant properties for group left-orders.
(Compare \cite[Question 3.22]{order}.)
For instance, one may say that a $6$-uple $(f,g;u_1,v_1,u_2,v_2)$ of
elements in an ordered group $(\Gamma,\preceq)$ is a
 {\bf{\em double resilient pair}} \esp if both $(f,g;u_1,v_1)$ and
$(g,f^{-1};u_2,v_2)$ are resilient pairs and $u_1 \prec u_2 \prec v_1$
(see Figure 14). Finding a simpler algebraic counterpart of the property of not
having a double crossing for a left-order seems to be an interesting problem.

The notion of $n$-resilient pair can be analogously defined. This corresponds
to a $(2n+2)$-uple $(f,g; u_1, v_1, u_2, v_2, \ldots, u_n, v_n)$ such that:

\vsp

\noindent -- \esp $(f,g;u_1,v_1)$, $(g,f^{-1};u_2,v_2)$, $(f^{-1},g^{-1};u_3,v_3)$,
$(g^{-1},f;u_4,v_4)$, $(f,g;u_5,v_5)$, etc, are all resilient pairs,

\vsp

\noindent -- \esp $u_i \prec u_{i+1} \prec v_i$, for all $i \in \{1,\ldots,n-1\}$.

\vsp

An eventual affirmative answer for the question below would have interesting
consequences; see Proposition \ref{laws}.

\vsp

\begin{question} \label{question-cle}
Let $\Gamma$ be a left-orderable group such that no left-order admits an
$n$-resilient pair for some (large) $n \in \mathbb{N}$. Does $\Gamma$
admit a $C$-order~?
\end{question}
\end{rem}\end{small}

\vspace{0.86cm}


\beginpicture

\setcoordinatesystem units <1cm,1cm>

\putrule from 0 0 to 6 0
\putrule from 0 0 to 0 6
\putrule from 6 0 to 6 6
\putrule from 0 6 to 6 6

\put{$f$} at 4.7 2.7
\put{$g$} at 1.5 3

\put{$u_1$} at 0 -0.6
\put{$v_2$} at 6 -0.6
\put{$u_2$} at 2.3 -0.6
\put{$v_1$} at 3.7 -0.6

\put{$\bullet$} at 0 0
\put{$\bullet$} at 6 0
\put{$\bullet$} at 2.3 0
\put{$\bullet$} at 3.7 0

\put{Figure 14: A double resilient pair.} at 3 -1.2
\put{} at -4.3 -1

\setquadratic

\plot
0 2.5
3 3
6 3.4 /

\plot
0 0.5
4 2.2
6 6.8 /

\setlinear

\setdots

\plot
3.7 0
3.7 3.7 /

\plot
0 3.7
3.7 3.7 /

\plot
2.3 0
2.3 6 /

\plot
0 2.3
6 2.3 /

\plot
0 0
6 6 /

\endpicture


\vspace{0.65cm}

\index{Order!Conradian}
\noindent{\bf Non-Conradian orders yield free subsemigroups.}
Let $\preceq$ be a non-Conradian order on a group $\Gamma$. Let
$(f,g;u,v) \in \Gamma^4$ be a resilient pair for $\preceq$, and denote
$$A := [u,fv]_{\preceq}:= \{w \!: u \preceq w \preceq fv \},
\qquad B := [gu,v]_{\preceq}.$$
Then $A$ and $B$ are disjoint, and for all $n \!\in\! \mathbb{N}$, we have $f^n (A \cup B) \subset A$ and
$g^n (A \cup B) \subset B$. This easily implies that the semigroup generated by $f$ and $g$ is free by 
an application of a ``positive version'' of the ping-pong lemma (see Exercise \ref{ejer-ping-pong}; 
see also \cite{harpe} in case of problems).

This shows in particular that all left-orders on torsion-free, virtually-nilpotent
groups are Conradian, a fact first established in \cite{longo} by different methods.
(This is no longer true for left-orderable polycyclic groups, even for
metabelian ones; see \cite[Corollary 7.5.6]{botto}.) Similarly, the equality
$\mathcal{LO}(\Gamma) = \mathcal{CO}(\Gamma)$ holds for left-orderable groups
$\Gamma$ with subexponential growth, as for example Grigorchuk-Mach\`{\i}'s group
\cite{GM,growth} (see Exercise \ref{close to 1} for a precise definition).
As a consequence of Proposition \ref{simple-alcanza}, we obtain the 
following result.

\vsp

\begin{thm} \label{nilpotent-cantor}
{\em The space of left-orders of a countable, torsion-free, virtually-nilpotent
group with infinitely many left-orders is homeomorphic to the Cantor
set. The same holds for countable, left-orderable groups without
free subsemigroups and having infinitely many left-orders.}
\end{thm}

\vsp

Note also that all left-orders on Tararin groups
({\em i.e.}, groups with finitely many left-orders; see \S \ref{fin-uncount})
are Conradian. Indeed, if $\Gamma$ has finitely many left-orders, then for every $g \!\in\! \Gamma$
and every left-order $\preceq$ on $\Gamma$, the left-order $\preceq_{g^{-n}}$ must coincide with
$\preceq$ for some finite $n$ (actually, for an \esp
$n$ \esp smaller than or equal to the cardinality of
$\mathcal{LO}(\Gamma)$). Thus, \esp $f \succ_{g^{n}} id$ \esp holds
for every $\preceq$-positive element $f$, that
is, \esp $g^{-n} f g^n \succ id$. \esp In particular, if $g \succ id$, then
$f g^n \succ g^n \succ g$, which shows that $\preceq$ is Conradian.

\vsp

\begin{question} Suppose that all left-orders on a finitely-generated,
left-orderable group are Conradian. Must the group be residually almost-nilpotent~?
\end{question}

\vsp\vsp\vsp

\noindent{\bf Every non-Conradian order leads to uncountably many left-orders.}
Using the notion of crossings, we show a refined version of Theorem
\ref{linnell-general} in the presence of non-Conradian orders.

\vsp

\begin{lem} {\em If $\preceq$ is a non-Conradian order on a group
$\Gamma$, then there exists $(f,g,h;u,v)$ in $\Gamma^5$ such that}
$$u \prec fu \prec fv \prec hu \prec hv \prec gu \prec gv \prec v.$$
\end{lem}

\noindent{\bf Proof.} Let $(\bar{f},\bar{g};u,v)$ be a resilient pair for
$\preceq$, so that \esp
$$u \prec \bar{f} u \prec \bar{f} v \prec \bar{g} u \prec \bar{g} v \prec v.$$
\esp Let $f:= \bar{f}$, \esp $h:= \bar{g} \bar{f}$, \esp $g := \bar{g}^2$. Then:

\vspace{0.15cm}

\noindent -- The inequality $fv \prec hu$ \esp is \esp $\bar{f} v \prec \bar{g} \bar{f} u$,
\esp which follows from
$$\bar{f} u \succ u \implies \bar{g} \bar{f} u \succ \bar{g} u \succ \bar{f} v;$$

\vspace{0.1cm}

\noindent -- The inequality $hv \prec gu$ \esp is \esp
$\bar{g} \bar{f} v \prec \bar{g}^2 u$, \esp which follows from
$$\bar{f} v \prec \bar{g} u \implies \bar{g} \bar{f} v \prec \bar{g}^2 u.$$

\vspace{-0.95cm}

$\hfill\square$

\vspace{0.3cm}

\begin{thm} {\em If $\preceq$ is a non-Conradian order on a group $\Gamma$, then
the closure of the orbit of $\preceq$ in $\mathcal{LO}(\Gamma)$ contains a Cantor set.}
\end{thm}

\noindent{\bf Proof.} Fix $(f,g,h;u,v)$ as in the previous lemma. Let $I$ denote
the closure of the subset $\{ \preceq_{w} : u \preceq w \preceq v \}$
of $\mathcal{LO}(\Gamma)$. (Recall that $\preceq_{w}$ is the
left-order with positive cone $w^{-1} P_{\preceq}^+ w$.)
Let $I^+ := \{\preceq' \in I \!: h \succ' id \}$
and $I^+ := \{\preceq' \in I \!: h \prec' id \}$.
We claim that $f(I) \subset I^+$ and $g (I) \subset I^-$.
Indeed, to show that $f(I) \subset I^+$, we need to check that $h (fw) \succ fw$
for all $u \preceq w \preceq v$. But this follows from
$$h(fw) \succeq h(fu) \succ hu \succ fv \succeq fw.$$
The proof of the containment $g(I) \subset I^-$ is analogous.

Denote $\Lambda := \{0,1\}^{\mathbb{N}}$,
and let $h_0 := f$ and $h_1 := g$. Consider the map
$$\Lambda \rightarrow \mathcal{P} \big( \overline{orb(\preceq)} \big),
\qquad \iota = (i_1,i_2,\ldots) \mapsto
\bigcap_{n \geq 1} h_{i_1} h_{i_2} \cdots h_{i_n} (I) = \iota (I).$$
By the claim above, if $\iota \neq \iota'$, then
$\iota (I) \cap \iota' (I) = \emptyset.$
The theorem follows. $\hfill\square$


\subsection{An extension to group actions on ordered spaces}
\label{conrad-general}

\index{Crossing}
\hspace{0.45cm} Let $\Gamma$ be a group acting by order-preserving
bijections on a totally ordered space $(\Omega,\leq)$. A {\bf{\em crossing}}
for the action of $\Gamma$ on $\Omega$ is a 5-tuple $(f,g;u,v,w)$, where
$f,g$ belong to $\Gamma$ and $u,v,w$ are in $\Omega$, such that:

\vsp

\noindent -- It holds \esp $u \prec w \prec v$;

\vsp

\noindent -- For every $n \in \mathbb{N}$, we have
\esp $g^n u \prec v$ \esp and \esp $f^n v \succ u$;

\vsp

\noindent -- There exist $M,N$ in $\mathbb{N}$ so that \esp
$f^N v \prec w \prec g^M u$.

\vsp

\index{Crossing}
\index{Crossing!reinforced}
\index{Resilient pair}
\noindent Analogous definitions of {\em reinforced crossings} and {\em resilient pairs}

Note that for a left-ordered group $(\Gamma,\preceq)$, the
notions of the preceding section correspond to the above ones for the
left-action on the ordered space $(\Gamma,\preceq)$. This is why we 
will sometimes call a $\Gamma$-action on a totally ordered space 
$\Omega$ with no crossings simply a {{\bf \em Conradian action}}.  
\index{Conradian!action}

For another relevant example, recall from Remark \ref{ssh} that, given a
$\preceq$-convex subgroup $\Gamma_0$ of a left-ordered group $(\Gamma,\preceq)$,
the space of left cosets $\Omega = \Gamma / \Gamma_0$ carries a natural total
order $\leq$ that is invariant by the left-translations. (Taking $\Gamma_0$
as the trivial subgroup, this reduces to the preceding example.)
Whenever this action has no crossing, we will say that $\Gamma$
is a {\bf{\em $\preceq$-Conradian extension}} of $\Gamma_0$.
Of course, this is the case of every convex subgroup
$\Gamma_0$ if $\preceq$ is Conradian.
\index{Conradian!extension}

\vsp

\begin{small}\begin{rem} Let $(\Gamma,\preceq)$ be a left-ordered group, and let
$\Gamma_0$ be a $\preceq$-convex subgroup. Given any left-order $\preceq_*$ on
$\Gamma_0$, let $\preceq'$ be the extension of $\preceq_*$ by $\preceq$. One
readily checks that $\Gamma$ is a $\preceq$-Conradian extension of $\Gamma_0$
if and only if it is a $\preceq'$-Conradian extension of it.
\label{convex-extension}
\end{rem}
\index{Conradian!extension}

\begin{ejer} \label{ejer:irreducible}
Let $\Gamma$ be a subgroup of $\mathrm{Homeo}_+(\mathbb{R})$. Say that an open interval $I$
is an {\bf {\em irreducible component}} of a nontrivial element $g \in \Gamma$ if it is fixed by $g$
and contains no fixed point inside. Equivalently, $I$ is a connected component of the complement
of the set of fixed points of $g$.

\noindent (i) Show that if the action of $\Gamma$ is without crossings, then for any pair of
different irreducible components, either one of them contains the other, or they are disjoint.

\noindent (ii) Show that the converse of (i) also holds.
\end{ejer}\end{small}

\vsp

\index{Action!cofinal}
For a general order-preserving action of a group $\Gamma$ on a totally ordered space
$(\Omega,\leq)$, the action of an element $f \!\in\! \Gamma$ is said to be {\bf{\em cofinal}}
if for all $x < y$ in $\Omega$ there exists $n \!\in\! \mathbb{Z}$ such that $f^n(x) > y$.
Equivalently, the action of $f$ is not cofinal if there exist $x<y$ in $\Omega$ such that
$f^n (x) < y$ for every integer $n$. If $(\Gamma,\preceq)$ is a left-ordered group, then
$f \!\in\! \Gamma$ is {\bf\em cofinal} if it is so for the corresponding left action of $\Gamma$
on itself.
\index{Cofinal element}

\vspace{0.12cm}

\index{Crossing}
\index{Action!cofinal}
\begin{prop} \label{prop:es-un subgrupo}
{\em Let $\Gamma$ be a group acting by order-preserving bijections on a
totally ordered space $(\Omega,\leq)$. If the action has no crossings, then the set
of elements whose action is not cofinal is a normal subgroup of \esp $\Gamma$.}
\label{normality}
\end{prop}

\noindent{\bf Proof.} Let us denote the set of elements whose action is not cofinal
by $\Gamma_0$. This set is normal. Indeed, given $g \in \Gamma_0$, let $x < y$ in
$\Omega$ be such that $g^n(x) < y$ for all $n$. For each $h \in \Gamma$ we have
$g^n h^{-1} (h(x)) < y$, hence $(hgh^{-1})^n (h(x)) < h(y)$ (for all
$n \!\in\! \mathbb{Z}$). Since $h(x) < h(y)$, this shows that $hgh^{-1}$
belongs to $\Gamma_0$.

It follows immediately from the definition that $\Gamma_0$ is stable under
inversion, that is, $g^{-1}$ belongs to $\Gamma_0$ for all $g \!\in\! \Gamma_0$.
The fact that $\Gamma_0$ is stable under multiplication is more subtle. For the
proof, given $x \!\in\! \Omega$ and $g \in \Gamma_0$, we will denote by $I_g (x)$
the {\bf{\em convex closure}} of the set $\{g^n(x)\!\!: n \in \mathbb{Z}\}$, that is,
the set formed by the $y \in \Omega$ for which there exist $m,n$ in $\mathbb{Z}$
so that $g^m (x) \leq y \leq g^n (x)$. Note that $I_g(x) = I_g(x')$ for all
$x' \in I_g (x)$. Moreover, $I_{g^{-1}} (x) = I_g (x)$ for all $g \!\in\! \Gamma_0$
and all $x \!\in\! \Omega$. Finally, if $g (x) = x$, then $I_g (x) = \{ x \}$. We
claim that if $I_g (x)$ and $I_f (y)$ are not disjoint for some $x,y$ in $\Omega$
and $f,g$ in $\Gamma_0$, then one of them contains the other. Indeed,
assume that there exist non-disjoint sets $I_f (y)$ and $I_g (x)$, none of
which contains the other. Without loss of generality, we may assume that $I_g (x)$
contains points to the left of $I_f (y)$ (if this is not the case, just interchange
the roles of $f$ and $g$). Changing $f$ and/or $g$ by their inverses if necessary,
we may assume that $g(x) >x$ and $f(y) < y$, thus $g(x') > x'$ for all
$x' \in I_g (x)$, and $f(y') < y'$ for all $y' \in I_y(f)$. Take
$u \in I_g(x) \setminus I_f(y)$, $w \in I_g(x) \cap I_f(y)$, and
$v \in I_f(y) \setminus I_g(x)$. Then one easily checks that
$(f,g;u,v,w)$ is a crossing, which is a contradiction.

Now, let $g,h$ be elements in $\Gamma_0$, and let $x_1 < y_1$ and $x_2 < y_2$ be points in
$\Omega$ such that $g^n (x_1) < y_1$ and $h^n (x_2) < y_2$, for all $n \!\in\! \mathbb{Z}$.
Set $x := \min\{x_1,x_2\}$ and $y := \max\{y_1,y_2\}$. Then $g^n (x) < y$ and $h^n (x) < y$,
for all $n \in \mathbb{Z}$; in particular, $y$ does not belong to neither $I_g (x)$ nor
$I_h (x)$. Since $x$ belongs to both sets, we have either $I_g(x) \subset I_h(x)$ or
$I_h(x) \subset I_g(x)$. Both cases being analogous, let us consider only the first
one. Then for all $x' \!\in\! I_g(x)$ we have $I_h (x') \subset I_g (x') = I_g(x)$.
In particular, $h^{\pm 1}(x')$ belongs to $I_g (x)$ for all $x' \in I_g (x)$. Since
the same holds for $g^{\pm 1}(x')$, this easily implies that $(gh)^n(x) \in I_g(x)$,
for all $n \in \mathbb{Z}$. As a consequence, $(gh)^n(x) < y$ holds for all
$n \in \mathbb{Z}$, thus showing that $gh$ belongs to $\Gamma_0$. $\hfill\square$

\vsp

\begin{small}\begin{ejer}\label{ejer:es-un-subgrupo}
Using the preceding proposition, show that for a nilpotent group action on the real line, 
the set of elements having fixed points forms a normal subgroup. Show that this holds 
more generally for groups with no free subsemigroups.
\end{ejer}\end{small}
\vsp

Slightly extending Example \ref{salta}, a
{\bf{\em convex jump}} of a left-ordered group $(\Gamma,\preceq)$ is
a pair $(G,H)$ of distinct $\preceq$-convex subgroups such that $H$ is
contained in $G$, and there is no $\preceq$-convex subgroup between them.

\vspace{0.15cm}

\begin{thm} {\em Let $(\Gamma,\preceq)$ be a left-ordered group, and let $(G,H)$
be a convex jump in $\Gamma$. Suppose that $G$ is a Conradian extension of
$H$. Then $H$ is normal in $G$, and the left-order induced by $\preceq$ on the
quotient $G/H$ is Archimedean.}
\label{super-conrad}
\end{thm}\index{Order!Archimedean} 
\index{Conradian!extension}

\noindent{\bf Proof.} Let us consider the action of $G$ on the space of cosets
$G /H$. Each element of $H$ fixes the coset $H$, hence its action is not cofinal.
If we show that the action of each element in $G \setminus H$ is cofinal, then
Proposition \ref{normality} will imply the normality of $H$ in $G$.

Now given $f \in G \setminus H$, let $G_f$ be the smallest convex subgroup of
$G$ containing ($H$ and) $f$. We claim that $G_f$ coincides with the set \esp
$$S_f \!:=\! \{g \in G \!: f^m \prec g \prec f^n \mbox{ for some }
m,n \mbox{ in } \mathbb{Z}\}.$$
Indeed, $S_f$ is clearly a convex subset of $G$ containing $H$ and contained in $G_f$.
Thus, to show that $G_f = S_f$, we need to show that $S_f$ is a subgroup. To do
this, first note that, with the notation of the proof of Proposition \ref{normality},
the conditions $g \in S_f$ and $I_g (H) \subset I_{f} (H)$ are equivalent. Therefore,
for each $g \in S_f$, we have $I_{g^{-1}} (H) = I_g (H) \subset I_f (H)$, thus
$g^{-1} \! \in S_f$. Moreover, if $\bar{g}$ is another element in $S_f$, then
$\bar{g} g H \in \bar{g} (I_f (H)) = I_f (H)$, hence $I_{\bar{g}g} (H)
\subset I_f (H)$. This means that $\bar{g} g$ belongs to $S_f$, thus
concluding the proof that $S_f$ and $G_f$ coincide.

Each $f \in G \setminus H$ leads to a convex subgroup $G_f = S_f$ strictly containing $H$.
Since $(G,H)$ is a convex jump, we necessarily have $S_f = G$. Given $g_1 \prec g_2$ in $G$,
choose $m_1,n_2$ in $\mathbb{Z}$ for which $f^{m_1} \prec g_1$ and $g_2 \prec f^{n_2}$.
Then we have $f^{n_2-m_1} g_1 \succ f^{n_2-m_1} f^{m_1} = f^{n_2} \succ g_2$, hence
$f^{n_2 - m_1} (g_1 H) \geq g_2 H$. \esp This easily implies that the action of
$f$ is cofinal.

We have then showed that $H$ is normal in $G$. The left-invariant total order on the
space of cosets $G / H$ is therefore a group left-order. Moreover, given $f,g$ in $G$,
with $f \notin H$, the previous argument shows that there exists $n \!\in\! \mathbb{Z}$
such that $f^n \succ g$, thus $f^n H \succeq gH$. This is nothing but the
Archimedean property for the induced left-order on $G / H$. $\hfill\square$

\vsp

\begin{cor} {\em Under the hypothesis of Theorem} \ref{super-conrad},
{\em up to multiplication by a positive real number, there exists a unique
homomorphism $\tau \!: G \rightarrow (\mathbb{R},+)$ such that $ker(\tau) \!=\! H$
and $\tau(g) \!>\! 0$ for every positive element $g \in G \setminus H$.}
\label{mas}
\end{cor}

\vsp

\index{Conradian!extension}
The homomorphism $\tau$ above will be referred to as the {\bf {\em Conrad
homomorphism}} associated to the corresponding Conradian extension (jump).
\index{Conrad homomorphism}

\begin{small}\begin{ejer}\label{puntos-fijos-nilp}
Let $\Gamma$ be a subgroup of $\mathrm{Homeo}_+(\mathbb{R})$. Show that the action
of $g \in \Gamma$ is not cofinal if and only if $g$ has fixed points on the line. If $\Gamma$ is
finitely-generated and acts without crossings, show that the normal subgroup formed by the
elements having fixed points has global fixed points. If the action corresponds to the
dynamical realization of a left-order $\preceq$, show that this subgroup coincides
with the kernel of the Conrad homomorphism associated to the convex jump
with respect to the maximal proper $\preceq$-convex subgroup (see 
Example \ref{actuar-convexo}).
\end{ejer}\end{small}
\index{Conrad homomorphism}
\index{Action!cofinal}
\index{Crossing}

\vsp\vsp\vsp

\index{Group!locally indicable}
\noindent{\bf $C$-orderability implies local indicability.} Let $\preceq$ be a
Conradian order on a group $\Gamma$. Let $\Gamma_0$ be a nontrivial subgroup
of $\Gamma$ generated by finitely many positive elements $f_1 \prec \ldots \prec f_k$.
Let $\Gamma_f$ (resp. $\Gamma^f$) be the largest (resp. smallest) convex subgroup
which does not contain $f := f_k$ (resp. which contains $f$). By the
corollary above, there exists a nontrivial homomorphism
$\tau \!: \Gamma^f \rightarrow (\mathbb{R},+)$ such that
$ker(\tau) \!=\! \Gamma_f$. This shows that
$\Gamma$ is locally indicable.
\index{Order!Conradian}

\vsp\vsp\vsp

\begin{small}\begin{rem}\label{la-necesito}
The homomorphism $\tau$ produced above respects orders: if $f \preceq g$, then $\tau(f) \leq \tau(g)$.
Moreover, it is trivial when restricted to the maximal convex subgroup. As commutators are mapped to
zero by $\tau$, we conclude that every element in $[\Gamma,\Gamma]$ is strictly smaller than any other
element $f$ satisfying $\tau(f) > 0$.
\end{rem}\end{small}

\vsp\vsp\vsp

We close this section with the following analogue of Proposition \ref{rel-convex}.

\vsp

\begin{prop} \label{rel-convex-conrad} {\em Let $\Gamma$ be a
$C$-orderable group, and let $\{\Gamma_{\lambda} \!: \lambda \in \Lambda\}$
be a family of subgroups each of which is convex with respect to a $C$-left-order
\esp $\preceq_{\lambda}$. Then there exists a $C$-left-order on $\Gamma$ for which
the subgroup \esp $\bigcap_{\lambda} \Gamma_{\lambda}$ \esp is convex.}
\end{prop}

\vsp

The proof is based on a result concerning left-orders obtained from actions
on a totally ordered space.

\vsp

\begin{prop} \label{de-cristobal}
{\em Let $\Gamma$ be a group acting faithfully by order-preserving transformations
on a totally ordered space $(\Omega,\leq)$. If the action has no crossings, then
all induced left-orders on $\Gamma$ are Conradian.}
\end{prop}
\index{Order!Conradian}
\index{Crossing}

\noindent{\bf Proof.} Suppose that the left-order $\preceq$ on $\Gamma$
induced from the action via a well-order $\leqslant_{wo}$ on $\Omega$
(see \S \ref{general-3}) is not Conradian. Then there are
$\preceq$-positive elements $f,g$ in $\Gamma$ such that $fg^{n} \prec g$,
for every $n \in \mathbb{N}$. This easily implies $f \prec g$. Let
$\bar{w} := \min_{\leqslant_{wo}} \{ w_f, w_g \}$. (Recall that
$w_f := \min_{\leqslant_{wo}} \{w \!: f(w) \neq w \}$, and similarly for $w_g$.)
We claim that $(fg,fg^2; \bar{w},g(\bar{w}), fg^2(\bar{w}))$ is a crossing for the action.
Indeed:

\vsp

\noindent -- From $id \prec f \prec g$ we obtain $\bar{w} = w_g \leqslant_{wo}
w_f$ and $g(\bar{w})>\bar{w}$; moreover, $f(\bar{w}) \geq \bar{w},$
which together with $fg^{n} \prec g \,$ yields
$$\bar{w}< fg^2(\bar{w}) < g(\bar{w}).$$

\vsp

\noindent -- The preceding argument actually
shows that $fg^n(\bar{w})<g(\bar{w})$, for
all $n \in \mathbb{N}$. As a consequence,
$fg^2fg^2(\bar{w})< fg^3(\bar{w})< g(\bar{w})$.
A straightforward inductive argument then shows that \esp
$(fg^2)^n(\bar{w})< g(\bar{w}), \mbox{ for all } n \in \mathbb{N}.$
Moreover, from $g(\bar{w})>\bar{w}\,$ and $\,f(\bar{w})\geq\bar{w}$,
\esp we conclude that \esp $\bar{w}<(fg)^n(g(\bar{w})).$

\vsp

\noindent -- Finally, from $\bar{w} < fg^2(\bar{w})$ we obtain $fg^2(\bar{w}) <
fg^2(fg^2(\bar{w}))=(fg^2)^2(\bar{w})$, while $fg^2(\bar{w}) <
g(\bar{w})$ implies $(fg)^2(g(\bar{w}))= fg (fg^2(\bar{w})) < fg
(g(\bar{w}))=fg^2(\bar{w})$. $\hfill\square$

\vsp\vsp\vsp\vsp\vsp

The proof of Proposition \ref{rel-convex-conrad} proceeds as that of
Proposition \ref{rel-convex}. We consider the left action of $\Gamma$ on
$\Omega := \prod_{\lambda \in \Lambda} \Gamma / \Gamma_{\lambda} \times \Gamma$
endowed with the lexicographic order. The stabilizer of
$([id_{\lambda}])_{\lambda \in \Lambda} \times \Gamma$ coincides
with $\bigcap_{\lambda} \Gamma_{\lambda}$, which may be made
convex for an induced left-order $\preceq$ on $\Gamma$.
Now the main point is that, as the action of $\Gamma$ on each
$\Gamma / \Gamma_{\lambda}$ has no crossings, the same holds for
the action of $\Gamma$ on $\Omega$. By Proposition \ref{de-cristobal},
the left-order $\preceq$ is Conradian, thus concluding the proof.

\vsp\vsp

\index{Crossing}
\begin{small}\begin{ejer}\label{ejer:crossings}
Prove the following converse to Proposition \ref{de-cristobal}: If $(\Gamma,\preceq)$
is a countable $C$-ordered group, then its dynamical realization is an action
on the real line without crossings (see \S \ref{general-3}).
\end{ejer}\end{small}


\subsection{The Conradian soul of a left-order}
\label{C-soul}

\hspace{0.45cm} A subgroup of a left-ordered group $(\Gamma,\preceq)$ is said
to be Conradian if the restriction of $\preceq$ to it is a Conradian order.
The {\bf{\em Conradian soul}} \esp $C_{\preceq}(\Gamma)$ of $(\Gamma,\preceq)$
is the (unique) subgroup that is $\preceq$-convex, $\preceq$-Conradian, and that
is maximal among subgroups verifying these two properties simultaneously.
\index{Conradian!soul}

\vsp

\begin{small}\begin{ex} Recall from Example \ref{sobre-b3-b4}
that the commutator subgroup $[\mathbb{B}_3, \mathbb{B}_3]$ is isomorphic to $\mathbb{F}_2$, with  
$\s_1 \s_2^{-1}$ and $\s_1^2 \s_2^{-2}$ as free generators. Denote by $\preceq$ the restriction of the Dehornoy
left-order to $[\mathbb{B}_3,\mathbb{B}_3]$. As we show below, $\preceq$ has no proper convex subgroups.\footnote{This
example is due to Clay \cite{exotic}. However, the existence of a left-order on $\mathbb{F}_2$ with
no proper convex subgroups also follows from the work of McCleary \cite{SM85}. See also \cite{NW}
for left-orders on braid groups without proper convex subgroups.} Since, as it is easily shown,
$\preceq$ is non-Conradian (compare Example \ref{deh-no}), its Conradian soul is trivial.

Let $C \subset \mathbb{F}_2 = [\mathbb{B}_3, \mathbb{B}_3]$ be a nontrivial convex subgroup.
Clearly, we may choose a 1-positive element $\sigma \in \mathbb{F}_2$. If $\sigma$
commutes with $\s_2$, then one may show that $\sigma$ is of the form
$\sigma = \Delta^{2p} \s_2^q$ for some integers $p, q$ satisfying $3p = -q > 0$,
where $\Delta = \s_1 \s_2 \s_1$. We thus have $\Delta^2 \prec \Delta^{4p} \s_2^{-6p} = \sigma^2$.
Since $\Delta^2$ is cofinal for the Dehornoy left-order and central, $\sigma$ is cofinal as
well. Since $C$ is a convex subgroup containing $C$, it must coincide with $\mathbb{F}_2$.

Suppose now that $\sigma$ and $\s_2$ do not commute. By the Subword Property
(see \S \ref{ejemplificando-5}), for every $k > 0$ the braid
$\sigma \s_2^k \sigma^{-1}$ is $1$-positive, as well as
$\sigma \s_2^k \sigma^{-1} \s_2^{-k}$.
Next, $\s_2^k \sigma^{-1} \s_2^{-k}$ is $1$-negative, so that
$\sigma \s_2^k \sigma^{-1} \s_2^{-k} \prec \sigma$. By convexity,
$\sigma \s_2^k \sigma^{-1} \s_2^{-k}$ must lie in $C$. Since $\s \in C$,
both $ \s_2^k \sigma^{-1} \s_2^{-k}$ and $ \s_2^k \sigma \s_2^{-k}$ belong to $C$.
Now $\s$ may be represented as $\s_2^m \s_1 w$, where $m$ is an integer, and $w$ is a
$1$-positive, $1$-neutral, or empty word. Choose $k>0$ so that $m' = k + m > 0$, and set
$\s' := \s_2^k \s \s_2^{-k}$. We know that $\s'$ lies in $C$, and it may be
represented by the $1$-positive braid word $\s_2^{\ell} \s_1 w \s_2^{-k}$. We will now
proceed to show that $C$ must contain both generators of $\mathbb{F}_2$, thus
$C = \mathbb{F}_2$. First note that
$\s_2 (\s_1^{-1} \s_2^{\ell} \s_1) w \s_2^{-k} =
\s_2 (\s_2 \s_1^{\ell} \s_2^{-1}) w \s_2^{-k}$
is 1-positice. Therefore,
\[ id \prec \s_2 \s_1^{-1} \s_2^{\ell} \s_1 w \s_2^{-k} \quad
\implies \quad \s_1 \s_2^{-1} \prec \s_2^{\ell} \s_1 w \s_2^{-k} = \s' \in C, \]
and since $id \prec \s_1 \s_2^{-1}$, this implies that $\s_1 \s_2^{-1} \in C$ by
convexity. Concerning the second generator $\s_1^2 \s_2^{-2}$, observe that
$$\s_2^2 \s_1^{-1} \s_1^{-1} \s_2^{\ell} \s_1 w \s_2^{-k} =
\s_2^2 \s_1^{-1} \s_2 \s_1^{\ell} \s_2^{-1} w \s_2^{-k} =
\s_2^2 \s_2 \s_1 \s_2^{-1} \s_1^{\ell - 1} \s_2^{-1} w \s_2^{-k}$$
is 1-positive. Thus,
$$id \prec \s_2^2 \s_1^{-2} \s_2^{\ell} \s_1 w \s_2^{-k} \quad
\implies \quad \s_1^2 \s_2^{-2} \prec \s_2^{\ell} \s_1 w \s_2^{-k} = \sigma' \in C,$$
and since $1 \prec \s_1^2 \s_2^{-2}$, we conclude from the convexity
of $C$ that $\s_1^2 \s_2^{-2} \in C$.
\end{ex}

\index{Order!DD-order}
\index{Conradian!soul}
\begin{ex} \label{CS-braids} The Conradian soul of $\preceq_{_D}$
on $\mathbb{B}_n$ is the cyclic subgroup generated by $\sigma_{n-1}$. Indeed,
this follows from the facts that the only
$\preceq_{n}$-convex subgroups of \esp $\mathbb{B}_{n}$ \esp are \esp $\{id\}$, \esp
$\langle \sigma_{n-1} \rangle$, \esp $\langle \sigma_{n-2},\sigma_{n-1} \rangle$,
\esp \ldots \esp, $\langle \sigma_2,\ldots,\sigma_{n-1} \rangle$ \esp and \esp
$\mathbb{B}_{n}$ itself, and that the restriction of $\preceq_{_{D}}$ to \esp
$\langle \sigma_{n-2}, \sigma_{n-1} \rangle \sim \mathbb{B}_3$ \esp
is not Conradian (see Example \ref{deh-no}). Let us
examine the case of $\mathbb{B}_3$ by denoting $a := \sigma_1 \sigma_1$ and
$b := \sigma_2^{-1}$. (For a general $\mathbb{B}_n$, one uses a similar argument
together with Theorem \ref{thm-DD}.) Recall that the family of $\preceq_{_D}$-convex
subgroups coincides with that of $\preceq_{_{DD}}$-convex ones.
Clearly, $\langle b \rangle$ does not properly contain any nontrivial convex
subgroup. Suppose that there exists a $\preceq_{_{DD}}$-convex subgroup $B$ of
$\mathbb{B}_3$ such that \esp $\langle b \rangle \subsetneq B \subsetneq \mathbb{B}_3$. \esp Let
$\preceq'$, $\preceq''$, and $\preceq'''$, be the left-orders defined on
$\langle b \rangle$, $B$, and $\mathbb{B}_3$, respectively, by:

\vspace{0.1cm}

\noindent -- $\preceq'$ is the restriction of $\preceq_{_{DD}}$ to $\langle b \rangle$;

\vspace{0.1cm}

\noindent -- $\preceq''$ is the extension of $\preceq'$ by the restriction of
$\overline{\preceq}_{_{DD}}$ to $B$;

\vspace{0.1cm}

\noindent -- $\preceq'''$ is the extension of $\preceq''$ by $\preceq_{_{DD}}$.

\vspace{0.1cm}

\noindent The left-order $\preceq'''$ is different from $\preceq_{_{DD}}$
(the $\preceq_{_{DD}}$-negative elements in $B \setminus \langle b \rangle$
are $\preceq'''$-positive), but its positive cone still contains the elements
$a,b$. Nevertheless, this is impossible, since these elements
generate the positive cone of $\preceq_{_{DD}}$.
\end{ex}

\begin{ejer}\label{no-tan-facil}
Let $\Gamma_* := C_{\preceq}(\Gamma)$ be the Conradian soul of a left-ordered group
$(\Gamma,\preceq)$. Show that, for any Conradian order $\preceq_*$ on $\Gamma_*$,
the extension of $\preceq_*$ by $\preceq$ has Conradian soul $\Gamma_*$.
\end{ejer}\end{small}
\index{Conradian!soul}

To give a dynamical counterpart of the notion of Conradian soul in terms of
crossings, we consider the set $C^+$ formed by the elements $h \!\succ\! id$
such that $h \preceq w$ for every crossing $(f,g;u,v,w)$ satisfying
$id \preceq u$. Analogously, we let $C^{-}$ be the set formed by the
elements $h \prec id$ such that $w \preceq h$ for every crossing
$(f,g;u,v,w)$ satisfying $v \preceq id$. Finally, we let
$$C := \{ id \} \cup C^{+} \cup C^{-}.$$
{\em A priori}, it is not clear that the set $C$ has a nice structure; 
for instance, it is not at all evident that it is a subgroup.
Nevertheless, we have the following result.

\vspace{0.2cm}

\index{Crossing}
\begin{thm} {\em The Conradian soul of $(\Gamma,\preceq)$ coincides
with the set $C$ above.}
\label{S=CS}
\end{thm}
\index{Conradian!soul}


Before passing to the proof, we give four general lemmas on crossings for
left-orders (note that the first three lemmas still apply to crossings
for actions on totally ordered spaces). The first one allows us replacing the
``comparison element'' $w$ by its ``images'' under positive iterates of either
$f$ or $g$.

\vspace{0.1cm}

\begin{lem} {\em If $(f,g;u,v,w)$ is a crossing,
then both $(f,g;u,v,g^n w)$ and $(f,g;u,v,f^n w)$ are
also crossings, for every $n \!\in\! \mathbb{N}$.}
\label{lema1}
\end{lem}

\noindent{\bf Proof.} We only consider the first 5-tuple
(the other is analogous). Since $g w \succ w$,
for every $n \! \in \! \mathbb{N}$ we have $u \prec w \prec g^n
w$; moreover, $v \succ g^{M+n} u = g^n g^M u \succ g^n w$. Hence,
\esp $u \prec g^n w \prec v$. \esp Furthermore, $f^N v \prec w \prec
g^n w$. Finally, from \esp $g^M u \succ w$, \esp we get \esp
$g^{M+n} u \succ g^n w$. \esp $\hfill\square$

\vspace{0.3cm}

Our second lemma allows replacing the ``limiting'' elements $u$ and $v$
by more appropriate ones.

\vspace{0.15cm}
\index{Crossing}

\begin{lem} {\em Let $(f,g;u,v,w)$ be a crossing. If $f u \succ u$ (resp.
$f u \prec u$) then $(f,g;f^n u,v, w)$ (resp. $(f,g;f^{-n} u, v, w)$) is also
a crossing, for every $n \geq 1$. Analogously, if $g v \prec v$ (resp. $g v \succ v$),
then $(f, g; u, g^n v, w)$ (resp. $(f,g;u,g^{-n}v,w)$) is also a crossing, for every
$n \geq 1$.}
\label{lema3}
\end{lem}

\noindent{\bf Proof.} Let us only consider the first 5-tuple (the second case is
analogous). Suppose that $fu \succ u$ (the case $fu \prec u$ may be treated similarly).
Then $f^n u \succ u$, which yields $g^M f^n u \succ g^M u \succ w$. To show
that $f^n u \prec w$, assume by contradiction that $f^n u \succeq w$. Then
$f^n u \succ f^N v$ yields $u \succ f^{N-n} v$, which is absurd.
$\hfill\square$

\vspace{0.3cm}

The third lemma relies on the dynamical nature of the crossing condition.

\vspace{0.15cm}

\begin{lem} {\em If $(f,g;u,v,w)$ is a crossing, then $(hfh^{-1},
hgh^{-1};hu,hv,hw)$ is also a crossing, for every $h \in \Gamma$.}
\label{lema2}
\end{lem}

\noindent{\bf Proof.} The three conditions to be checked are nothing but
the three conditions in the definition of crossing multiplied by $h$
on the left. $\hfill\square$

\vspace{0.3cm}

A direct application of the lemma above shows that, if $(f,g;u,v,w)$ is a crossing, then the
5-tuples $(f,f^ngf^{-n}; f^n u, f^n v , f^n w)$ and $(g^nfg^{-n}, g ;g^n u, g^n v, g^n w)$
are also crossings, for every $n \in \mathbb{N}$.

\vspace{0.1cm}

\begin{lem} {\em If $(f,g;u,v,w)$ is a crossing and $id \preceq h_1 \prec h_2$ are
elements in $\Gamma$ such that $h_1 \in C$ and $h_2 \notin C$, then there exists a
crossing $(\tilde{f},\tilde{g};\tilde{u},\tilde{v},\tilde{w})$ such that
$h_1 \prec \tilde{u} \prec \tilde{v} \prec h_2$.} \label{lemapro}
\end{lem}

\noindent{\bf Proof.} Since $id\prec h_2  \notin C$, there must be a
crossing $(f,g;u,v,w)$ such that $id\preceq u \prec w \prec h_2$.
Fix $N \in \mathbb{N}$ such that $f^N v\prec w $, and consider the
crossing
$$(f,\bar{g};\bar{u},\bar{v},\bar{w}) :=
(f, f^N g f^{-N}; f^N u, f^N v, f^N w).$$
Note that $\bar{v} = f^N v \prec w \prec h_2$. We claim that
$h_1 \preceq \bar{w} = f^N w$. Indeed, if $f^N u \succ u$ then
$f^n u \succ id$, and by the definition of $C$, we must have $h_1 \preceq \bar{w}$.
If $f^N u \prec u$, then we must have $f u \prec u$, thus by Lemma
\ref{lema3} we know that $(f,\bar{g};u,\bar{v},\bar{w})$ is also a
crossing, which still allows concluding that $h_1 \preceq \bar{w}$.

Now, for the crossing $(f,\bar{g};\bar{u}, \bar{v}, \bar{w})$, there exists
$M \in \mathbb{N}$ such that $\bar{w} \prec \bar{g}^M \bar{u}$. Let us
consider the crossing $(\bar{g}^M f \bar{g}^{-M}, \bar{g}; \bar{g}^M \bar{u},
\bar{g}^M \bar{v}, \bar{g}^M \bar{w})$. If $\bar{g}^M \bar{v} \prec \bar{v}$,
then $\bar{g}^M \bar{v} \prec h_2$, and we are done. If not, then
we must have $\bar{g} \bar{v} \succ \bar{v}$. By Lemma \ref{lema3},
$(\bar{g}^M f \bar{g}^{-M}, \bar{g}; \bar{g}^M \bar{u}, \bar{g}^M \bar{v}, \bar{w})$
is still a crossing, and since $\bar{v} \prec h_2$, this concludes the proof.
$\hfill\square$

\vspace{0.4cm}

\noindent{\bf Proof of Theorem \ref{S=CS}}. The proof is divided into several steps.

\vsp\vsp

\noindent {\underline{Claim (i).}} The set $C$ is convex.

\vsp

This follows directly from the definition of $C$.

\vsp\vsp

\noindent {\underline{Claim (ii).}} If $h$ belongs
to $C$, then $h^{-1}$ also belongs to $C$.

\vsp

Assume that $h \in C$ is positive and $h^{-1}$ does not belong to
$C$. Then there exists a crossing $(f,g; u,v,w)$ such that $h^{-1}
\prec w \prec v \preceq id$.

We first note that, if $h^{-1} \preceq u$, then after conjugating by
$h$ as in Lemma \ref{lema2}, we get a contradiction because $(hgh^{-1},
hfh^{-1}; hu, hv, hw)$ is a crossing with \esp $id \preceq hu $ \esp and
\esp $hw \prec hv \preceq h$. To reduce the case $h^{-1} \succ u$ to this one,
we first use Lemma \ref{lema2} and consider the crossing $(g^Mfg^{-M}, g; g^M u,
g^M v, g^M w)$. Since \esp $h^{-1} \prec w \prec g^M u \prec g^M w \prec g^M v$,
\esp if $g^M v \prec v$ then we are done. If not, Lemma \ref{lema3} shows that
$(g^Mfg^{-M}, g; g^M u , g^M v, w)$ is also a crossing, which still
allows concluding.

In the case where $h \in C$ is negative, we proceed similarly but we
conjugate by $f^N$ instead of $g^M$. Alternatively, since $id \in C$
and $id\prec h^{-1}$, if we suppose that $h^{-1}\notin C$ then Lemma
\ref{lemapro} provides us with a crossing $(f,g;u,v,w)$ such
that $id\prec u\prec w \prec v \prec h^{-1}$, which gives a
contradiction after conjugating by $h$.

\vsp\vsp

\noindent {\underline{Claim (iii).}} If $h$ and $\bar{h}$
belong to $C$, then $h\bar{h}$ also belongs to $C$.

\vsp

First, we show that for every pair of positive elements in $C$,
their product still belongs to $C$. (Note that, by Claim (ii), the
same will be true for pairs of negative elements in $C$.) Indeed,
suppose that $h,\bar{h}$ are positive elements, with $h \!\in\! C$ but
$h \bar{h} \notin C$. Then, by Lemma \ref{lemapro}, we may produce
a crossing $(f,g;u,v,w)$ such that $h \prec u \prec v \prec h \bar{h}$.
After conjugating by $h^{-1}$, we obtain the crossing
$(h^{-1}fh,h^{-1}gh;h^{-1}u,h^{-1}v,h^{-1}w)$
satisfying $id \prec h^{-1}u \prec h^{-1} w \prec \bar{h}$,
which shows that $\bar{h} \notin C$.

Now, if $h \prec id \prec \bar{h}$, then
$h \prec h\bar{h}$. Thus, if $h\bar{h}$ is negative, then
the convexity of $C$ yields $h\bar{h} \in C$. If $h\bar{h}$
is positive, then $\bar{h}^{-1}h^{-1}$
is negative, and since $\bar{h}^{-1} \prec \bar{h}^{-1} h^{-1}$,
the convexity gives again that $\bar{h}^{-1}h^{-1}$, hence
$h\bar{h}$, belongs to $C$. The remaining case
$\bar{h} \prec id \prec h$ may
be treated similarly.

\vsp \vsp

\noindent {\underline{Claim (iv).}} The subgroup $C$ is Conradian.

\vsp

In order to apply Theorem \ref{Conrad=noCrossing}, we need to show that there are no
crossings in $C$. Suppose by contradiction that $(f,g;u,v,w)$ is a crossing such that
$f,g,u,v,w$ all belong to $C$. If $id \preceq w$ then, by Lemma \ref{lema2}, we have
that $(g^n f g^{-n}, g; g^n u, g^n v, g^n w)$ is a crossing. Taking $n \!=\! M$ so
that $g^M u \succ w$, this contradicts the definition of $C$, because
$id \preceq w \prec g^M u \prec g^M w \prec g^M v \in C$. The case
$w \preceq id$ may be treated analogously by
conjugating by powers of $f$ instead of $g$.

\vsp\vsp

\noindent {\underline{Claim (v).}} The subgroup $C$ is maximal
among $\preceq$-convex, $\preceq$-Conradian subgroups.

\vsp\vsp

Indeed, if $H$ is a subgroup strictly containing $C$, then there is
a positive element $h \!\in\! H \setminus C$. By Lemma \ref{lemapro},
there exists a crossing $(f,g;u,v,w)$ such that $id\prec u \prec w
\prec v \prec h$. If $H$ is convex, then $u,v,w$ belong to $H$.
To conclude that $H$ is not Conradian, it suffices to show
that $f$ and $g$ belong to $H$.

On the one hand,
since $id \prec u$, we have either \esp $id \prec g \prec g u \prec v$ \esp \esp or
\esp \esp $id \prec g^{-1} \prec g^{-1} u \prec v$. \esp In both cases, the convexity
of $H$ implies that $g$ belongs to $H$. On the other hand, if $f$ is positive,
then from $f^N \prec f^N v \prec w$ we get $f \in H$, whereas in the case
of a negative $f$, the inequality \esp $id \prec u$ \esp gives \esp
$id\prec f^{-1} \prec f^{-1} u \prec v$, \esp which still
shows that $f \in H$. $\hfill \square$


\subsection{Approximation of left-orders and the Conradian soul}
\label{section-soul}

\hspace{0.45cm} The notion of Conradian soul was introduced in \cite{order} as a tool
for leading with the problem of approximating a group left-order by its conjugates.
We begin with the case of a trivial Conradian soul. (Compare Proposition
\ref{prop affine} and its proof.)
\index{Conradian!soul}
\vsp

\index{Conradian!soul}
\begin{thm} {\em If the Conradian soul of an infinite left-ordered group
$(\Gamma,\preceq)$ is trivial, then $\preceq$ may be approximated by its conjugates.}
\label{primero}
\end{thm}

\vsp

We will give two different proofs for this theorem, each of which gives
some complementary information. The first one, due to Clay \cite{clay},
shows that every left-order that is not approximated by its conjugates admits a
nontrivial, convex, {\em bi-ordered} subgroup. This may also be obtained by using
the method of the second  proof below (which is taken from \cite{crossings}) under
the stronger assumption that $\preceq$ is isolated in $\mathcal{LO}(\Gamma)$.
Nevertheless, though more elaborate than the first (it uses the
results of the preceding section), this second proof
is suitable for generalization in the case where the Conradian soul is ``almost
trivial'' ({\em i.e.}, it is nontrivial but admits only finitely many left-orders; see
Theorem \ref{final} below).

\vsp\vsp\vsp

\noindent{\bf First proof of Theorem \ref{primero}.} Suppose that $\preceq$ cannot be
approximated by its conjugates, and let $g_1,\ldots,g_k$ be finitely many positive
elements such that the only conjugate of $\preceq$ lying in
$V_{g_1} \cap \ldots \cap V_{g_k}$ is $\preceq$ itself. (Recall that
$V_{g}$ denotes the set of left-orders making $g$ a positive element.)
For each index $i \in \{1,\ldots,k\}$, let
$$B_i^+ := \big\{ h \in \Gamma \!: id \preceq h \preceq g_i^n
\mbox{ for some } n \in \mathbb{N} \big\},$$
$$B_i := \big\{ h \in \Gamma \!: g_i^{-m} \preceq h \preceq g_i^n
\mbox{ for some } m,n \mbox{ in } \mathbb{N} \big\}.$$

\noindent{\underline{Claim (i).}} For some $j \!\in\! \{1,\ldots,k\}$ we
have $h^{-1} P^+_{\preceq} h = P_{\preceq}^+$ for every $h \in B_j^+$.

\vsp

If not, then for each $i$ there exists $h_i \in \Gamma$ such that
$id \prec h_i \preceq g_i^{n_i}$ for some $n_i \in \mathbb{N}$ and
$h_i^{-1} P_{\preceq}^+ h_i \neq P_{\preceq}^+$. Let $h:= \min \{h_1,\ldots,h_k\}$.
Then $h^{-1} P_{\preceq}^+ h \neq P_{\preceq}^+$. Moreover,
$h \preceq g_i^{n_i}$ for each $i$, thus $h^{-1} g_i^{n_i} \succeq id$.
Since $h$ is necessarily positive, this yields $h^{-1} g_i^{n_i} h \succ id$,
which implies $h^{-1} g_i h \succ id$, that is, 
$g_i \!\in\! h^{-1} P^+_{\preceq} h$. Since this holds for every $i$, by hypothesis, 
the conjugate left-order $\preceq_{h^{-1}}$ must coincide with $\preceq$, which
is a contradiction.

\vsp\vsp

\noindent{\underline{Claim (ii).}} All elements in $B_j^+$ stabilize
$P^+_{\preceq}$ (under conjugation).

\vsp

Indeed, from \esp $g_j^{-m} \preceq h \preceq g_j^n$ \esp we obtain \esp
$id \preceq g_j^m h \preceq g_j^{m+n}$. \esp Thus, $g_j^m h$ belongs to
$B_j^+$. Since $g^m_j$ also belongs to $B_j^+$, by Claim (i) above we have
$$(g_j^m h)^{-1} P_{\preceq}^+ (g_j^m h) = P_{\preceq}^+ \qquad \mbox{ and } \qquad 
g_j^{-m} P_{\preceq}^+ g_j^m = P_{\preceq}^+.$$
This easily yields $h^{-1} P_{\preceq}^+ h = P_{\preceq}^+$, which
in its turn implies $h P_{\preceq}^+ h^{-1} = P_{\preceq}^+$.

\vsp\vsp

\noindent{\underline{Claim (iii).}} The set $B_j$ is a $\preceq$-convex subgroup
of $\Gamma$, and the restriction of $\preceq$ to it is a bi-order (hence a
$C$-order).

\vsp

The convexity of $B_j$ as a set is obvious. Now, for each $h \in B_j$, the
relations $g_j^{-m} \preceq h \preceq g_j^n$ and \esp $h P_{\preceq}^+ h^{-1}
\! = P_{\preceq}^{+}$ \esp easily yield $g_j^m \succeq h \succeq g_j^n$, thus
showing that $h^{-1} \in B_j$. Similar arguments show that $h_1 h_2$ belongs
to $B_j^+$ for all $h_1,h_2$ in $B_j^+$, as well as that the
restriction of $\preceq$ to $B_j^+$ is bi-invariant. $\hfill\square$

\vsp\vsp\vsp\vsp

\noindent{\bf Second proof of Theorem \ref{primero}.}
Let $f_1 \prec f_2 \prec \ldots \prec f_k$ be finitely many
positive elements of $\Gamma$. We need to show that there exists a conjugate of $\preceq$
that is different from $\preceq$ but for which all the $f_i$'s are still positive.

Since $id \!\in\! C_{\preceq}(\Gamma)$ and
$f_1 \notin C_{\preceq}(\Gamma)$, Theorem \ref{S=CS}
and Lemma \ref{lemapro} imply that there is a crossing $(f,g;u,v,w)$
such that $id \prec u \prec v \prec f_1$. Let $M,N$ in $\mathbb{N}$ be such
that $ f^N v \prec w \prec g^M u$. We claim that $id\prec_{v} f_i$
and $id \prec_{w} f_i$ hold for all $1 \leq i \leq k$, but \esp
$g^M f^N \prec_{v} id$ \esp and \esp
$g^M f^N \succ_{w} id$. \esp Indeed,
since $id \prec v \prec f_i$, we have $v\prec f_i\prec f_i v$, thus $id
\prec v^{-1} f_i v$. By definition, this means that $f_i\succ_{v} id$.
The inequality $f_i \succ_{w} id$ is proved similarly. Now note that
$g^M f^N v \prec g^M w \prec v$, hence $g^M f^N \prec_{v} id$. Finally,
from $g^Mf^N w \succ g^M u \succ w $, \esp we deduce that \esp
$g^Mf^N \succ_{w} id$.

Now the preceding relations imply that the $f_i$'s are still
positive for both $\preceq_{v^{-1}}$ and $\preceq_{w^{-1}}$,
but at least one of these left-orders is different from
$\preceq$. This concludes the proof. $\hfill\square$

\vspace{0.35cm}

We next deal with the case where the Conradian soul is nontrivial but admits
finitely many left-orders ({\em i.e.}, it is a Tararin group; see \S \ref{fin-uncount}).
It turns out that, in this case, the left-order may fail to be an
accumulation point of its conjugates. A concrete example is given by
the $DD$-left-order on $\mathbb{B}_n$. Indeed, its Conradian soul is isomorphic to
$\mathbb{Z}$ (see Example \ref{CS-braids}), though it is an
isolated point of the space of braid left-orders because its positive
cone is finitely-generated (see \S \ref{fin-gen}).
Now the $DD$-left-order has the Dehornoy left-order $\preceq_{_D}$ as a natural
``associate'', in the sense that the latter may be obtained from the
former by successive flipings  along convex jumps. For the case of $B_3$,
this reduces to changing the left-order on the Conradian soul in the unique
possible way. As shown below, $\preceq_{_D}$ is an accumulation point
of its conjugates. Moreover, there is a sequence of conjugates of
$\preceq_{_{DD}}$ that converges to $\preceq_{_D}$ as well.
\index{Order!DD-order}

\begin{small}\begin{ex} The sequence of conjugates $\preceq_j$ of $\preceq_{_D}$
by $\sigma_2^{j} \sigma_1^{-1}$ converges to $\preceq_{_D}$ in a nontrivial
way. Indeed, if $w=\sigma_2^k$ for some $k>0$, then
$$\sigma_1^{-1}\sigma_2^{j} w \sigma_2^{-j}\sigma_1
= \sigma_1^{-1} \sigma_2^k \sigma_1 =
\sigma_2 \sigma_1^k \sigma_2^{-1} \succ_{_D} id.$$
If, on the other hand, $w$ is a $\sigma_1$-positive word, say
$w=\sigma_2^{k_1}\sigma_1\sigma_2^{k_2} \ldots \sigma_2^{k_{\ell-1}}\sigma_1
\sigma_2^{k_\ell}$, then
\begin{footnotesize}
$$\sigma_1^{-1}\sigma_2^j w \sigma_2^{-j}\sigma_1 =
\sigma_1^{-1}\sigma_2^j \sigma_2^{k_1}\sigma_1\sigma_2^{k_2}\ldots
\sigma_2^{k_{\ell-1}}\sigma_1\sigma_2^{k_\ell} \sigma_2^{-i}\sigma_1
= \sigma_2\sigma_1^{j+k_1}\sigma_2^{-1}\sigma_2^{k_2}\ldots
\sigma_2^{k_{\ell-1}}\sigma_1\sigma_2^{k_\ell} \sigma_2^{-n}\sigma_1.$$
\end{footnotesize}Thus, $\sigma_1\sigma_2^{-j}w\sigma_2^j\sigma_1$ is
$1$-positive for sufficiently large~$j$ (namely, for
$j>-k_1$). This proves the desired convergence. Finally,
$\preceq_j$ is different from $\preceq_{_D}$ for each positive integer~$j$,
since its smallest positive element is the conjugate of~$\sigma_2$ by
$\sigma_1 \sigma_2^{j}$, and this is different from the smallest
positive element of $\preceq_{_D}$, namely $\sigma_2$. We leave
to the reader the task of checking that the sequence of conjugates
of $\preceq_{_{DD}}$ by $\sigma_1^{-1} \sigma_2^{j}$ converges to
$\preceq_{_D}$ as well.
\end{ex}

\begin{rem} \label{rem:dehornoy-like}
The $\mathbb{B}_3$-case of the preceding example can be
generalized as follows: For all $m,n$ larger than 1, with
$(m,n) \neq (2,2)$, the left-order $\preceq$
on $G_{m,n} = \langle a,b \!: (ba^{m-1})^{n-1} b = a\rangle$
with positive cone $\langle a,b \rangle^+$
given by Theorem \ref{thm:torus} has Conradian soul
$\langle b \rangle \!\sim\! \mathbb{Z}$.
Flipping this order on the Conradian soul
yields a left-order $\preceq'$ that is
accumulated by its conjugates. Moreover, there is a sequence of
conjugates of $\preceq$ that also converges to $\preceq'$.
See \cite{pos-1} as well as \cite{ito-3,ito-2,ito-1} for more on this and
related examples.
\end{rem}\end{small}

It turns out that the phenomenon described above for braid groups occurs for
general left-ordered groups. To be more precise, let $\Gamma$ be a group
having a left-order $\preceq$ whose Conradian soul admits finitely many
left-orders $\preceq_1,\preceq_2,\ldots,\preceq_{2^n}$, where $\preceq_1$
is the restriction of $\preceq$ to its Conradian soul. Each $\preceq_j$
induces a left-order $\preceq^j$ on $\Gamma$, namely the convex extension
of $\preceq_j$ by $\preceq$. (Note that $\preceq^1$ coincides with $\preceq$.)
All the left-orders $\preceq^j$ share the same Conradian soul (see 
Exercise \ref{no-tan-facil}). Assume throughout that $\preceq$ is not
Conradian, which is equivalent to that $\Gamma$ is not a Tararin group.
\index{Conradian!soul}

\vsp\vsp

\begin{thm} {\em  With the notation above, at least one of the left-orders
$\preceq^j$ is an accumulation point of the set of conjugates of $\preceq$.}
\label{final}
\end{thm}

\vspace{0.01cm}

\begin{cor} {\em At least one of the left-orders
$\preceq^j$ is approximated by its conjugates.}
\label{finalito}
\end{cor}

\noindent{\bf Proof.} Assuming Theorem \ref{final}, we have
that $\preceq^k$ belongs to the set of accumulation points
$acc (orb(\preceq^1))$ of the orbit of $\preceq^1$ for some
$k$ in $\{1,\ldots,2^n\}$. Theorem \ref{final} applied to this
$\preceq^k$ instead of $\preceq$ shows the existence of $k' \in \{1,\ldots,2^n\}$
so that $\preceq^{k'} \in acc (orb(\preceq^k))$, and hence
$\preceq^{k'} \in acc(orb(\preceq^1))$. If $k'$ equals either $1$ or $k$ then we are
done; if not, we continue arguing in this way... In at most $2^n$ steps we will find
an index $j$ such that $\preceq^j \in acc(orb(\preceq^j))$. $\hfill\square$

\vspace{0.35cm}

Theorem \ref{final} will follow from the next

\vsp

\begin{prop} {\em Given an arbitrary finite family
$\mathcal{G}$ of $\preceq$-positive elements in
$\Gamma$, there exists $h \in \Gamma$ and a positive
$\bar{h} \notin C_{\preceq} (\Gamma)$ such that
$id \prec h^{-1} f h \notin C_{\preceq} (\Gamma)$
for all $f \in \mathcal{G} \setminus C_{\preceq} (\Gamma)$,
but \esp $id \succ h^{-1}\bar{h}h \notin C_{\preceq} (\Gamma)$.}
\label{a-probar}
\end{prop}

\vsp

\noindent{\bf Proof of Theorem \ref{final} from Proposition \ref{a-probar}.}
Let us consider the directed net formed by the
finite sets $\mathcal{G}$ of $\preceq$-positive elements. For each
such a $\mathcal{G}$, let $h_{\mathcal{G}}$ and $\bar{h}_{\mathcal{G}}$
be the elements in $\Gamma$ provided by Proposition \ref{a-probar}. After
passing to subnets of $(h_{\mathcal{G}})$ and $(\bar{h}_{\mathcal{G}})$
if necessary, we may assume that the restrictions of
$\preceq_{h_{\mathcal{G}}}$ to $C_{\preceq}(\Gamma)$
all coincide with a single $\preceq_j$. Now the properties
of $h_{\mathcal{G}}$ and $\bar{h}_{\mathcal{G}}$ imply:

\vsp

\noindent -- $f \succ^j id$ \esp \esp and \esp \esp
$f \esp (\succ^j)_{h_{\mathcal{G}}} \esp id$, \esp
\esp for all \esp $f \in \mathcal{G} \setminus C_{\preceq} (\Gamma)$;

\vsp

\noindent -- $\bar{h}_{\mathcal{G}} \succ id$, \esp \esp  but \esp \esp \esp
$\bar{h}_{\mathcal{G}} \esp \esp (\prec^j)_{h_{\mathcal{G}}} \prec id$.

\vsp

\noindent This clearly shows the theorem. $\hfill\square$

\vspace{0.5cm}

For the proof of Proposition \ref{a-probar} we will use some lemmas.

\vsp

\begin{lem} {\em For every $id \prec c \notin C_{\preceq} (\Gamma)$,
there is a crossing $(f,g;u,v,w)$ such that $u,v,w$ do not belong
to $C_{\preceq}(\Gamma)$ and $id \prec u\prec w \prec v \prec c$.}
\label{primer-lema}
\end{lem}

\noindent{\bf Proof.} By Theorem \ref{S=CS} and Lemma \ref{lemapro},
for every $id \preceq s\in C_{\preceq} (\Gamma)$ there exists a
crossing $(f,g;u,v,w)$ such that $s\prec u\prec w\prec v \prec c$.
Clearly, $v$ does not belong to $ C_{\preceq} (\Gamma)$. The element
$w$ is also outside $ C_{\preceq} (\Gamma)$, as otherwise the
element $a := w^2$ would satisfy $w \prec a \in C_{\preceq} (\Gamma)$,
which is absurd. Taking $M > 0$ so that $g^M u \succ w$, this gives \esp
$g^M u \notin C_{\preceq} (\Gamma)$, \esp $g^M w \notin C_{\preceq} (\Gamma)$,
\esp and $g^M v \notin C_{\preceq} (\Gamma)$. \esp Consider
the crossing $(g^Mfg^{-M}, g; g^M u,g^M v, g^M w)$. If $g^M v\prec v$, then
we are done. If not, then $gv \succ v$, and Lemma \ref{lema3} ensures that
$(g^Mfg^{-M}\!, g; g^M u,  v, g^Mw)$ is also a crossing, which still allows
concluding. $\hfill\square$

\vspace{0.1cm}

\begin{lem} {\em Given $id \prec c \notin C_{\preceq} (\Gamma)$, there exists
$id \prec a \notin C_{\preceq}(\Gamma)$ (with $a \prec c$) such that, for all
$id \preceq b \preceq a$ and all \esp $\bar{c} \succeq c$, one has
$id \prec b^{-1} \bar{c} b \notin C_{\preceq} (\Gamma)$.}
\label{segundo-lema}
\end{lem}

\noindent{\bf Proof.} Let us consider the crossing $(f,g;u,v,w)$
such that $id \prec u\prec w \prec v \prec c$ and such that $u,v,w$
do not belong to $ C_{\preceq} (\Gamma)$. We affirm that the lemma
holds for $a := u$. Indeed, if $id \preceq b \preceq u$, then from \esp
$b \preceq u \prec v \prec \bar{c}$ \esp we obtain \esp $id \preceq
b^{-1}u \prec b^{-1}v \prec b^{-1} \bar{c}$, thus the crossing
$(b^{-1}fb,b^{-1}gb;b^{-1}u,b^{-1}v,b^{-1}w)$ shows that $b^{-1}
\bar{c} \notin C_{\preceq} (\Gamma)$. Since $id \preceq b$, we conclude
that $id \prec b^{-1} \bar{c} \preceq b^{-1} \bar{c} b$, and the convexity
of $S$ implies that $b^{-1} \bar{c} b \notin C_{\preceq} (\Gamma)$.
$\hfill\square$

\vsp

\begin{lem} {\em For every $g \in \Gamma$, the set $g \esp C_{\preceq} (\Gamma)$
is convex. Moreover, for every crossing $(f,g;u,v,w)$, one has
$u C_{\preceq} (\Gamma) < w C_{\preceq} (\Gamma) < v C_{\preceq} (\Gamma)$,
in the sense that \esp $uh_1 \prec wh_2 \prec vh_3$ \esp holds for all $h_1,h_2,h_3$
in $ C_{\preceq} (\Gamma)$.}
\label{lema4}\end{lem}

\noindent{\bf Proof.} The verification
of the convexity of $g C_{\preceq} (\Gamma)$ is straightforward.
Suppose next that $uh_1 \succ wh_2$ for some $h_1,h_2$ in
$ C_{\preceq} (\Gamma)$. Then, since $u \prec w$, the convexity of both
left classes $u C_{\preceq} (\Gamma)$ and $w C_{\preceq} (\Gamma)$
gives the equality between them. In particular, there exists
$h \in C_{\preceq} (\Gamma)$ such that $uh = w$. Note that such
an $h$ must be positive, hence $id \prec h = u^{-1} w$. But since
$(u^{-1}fu, u^{-1}gu; id ,u^{-1}v, u^{-1}w)$ is a crossing, this
contradicts the definition of $ C_{\preceq} (\Gamma)$. The proof of
the fact that $w C_{\preceq} (\Gamma) \prec v C_{\preceq} (\Gamma)$
is similar. $\hfill\square$

\vspace{0.35cm}

\noindent{\bf Proof of Proposition \ref{a-probar}.} Indexing the elements of
$\mathcal{G} \!=\! \{f_1,\ldots,f_r\}$ so that $f_1 \prec \ldots \prec f_r$,
let $k$ be such that $f_{k-1} \in C_{\preceq} (\Gamma)$ but
$f_{k} \notin C_{\preceq} (\Gamma)$. Recall that, by Lemma
\ref{segundo-lema}, there exists $id \prec a \notin C_{\preceq} (\Gamma)$
such that, for every $id \preceq b \preceq a$, one has $id \prec b^{-1}
f_{k+j} b \notin C_{\preceq} (\Gamma)$ for all $j \geq 0$. We fix
a crossing $(f,g;u,v,w)$ such that $id \prec u\prec v\prec a$
and $u\notin C_{\preceq} (\Gamma)$. Note that the conjugacy
by $w^{-1}$ yields the crossing $(w^{-1}fw, w^{-1}gw;w^{-1}u,w^{-1}v,id)$.

\vsp\vsp

\noindent{\underline{Case I.}} It holds that $w^{-1}v \preceq a$.

\vsp

In this case, we claim that the proposition holds for the choice
$h := w^{-1}v$ and $\bar{h} := w^{-1}g^{M+1}f^N w$. To show this,
first note that neither $w^{-1}gw$ nor $w^{-1}fw$ belong to $C_{\preceq} (\Gamma)$.
Indeed, this follows from the convexity of $C_{\preceq} (\Gamma)$ and the
inequalities \esp $w^{-1}g^{-M}w \prec w^{-1}u \notin C_{\preceq} (\Gamma)$
\esp and $w^{-1} f^{-N} w \succ w^{-1} v \notin C_{\preceq} (\Gamma)$.
\esp We also have $id \!\prec\! w^{-1}g^{M}f^{N}w$, hence $id \!\prec\!
w^{-1}gw \prec w^{-1} g^{M+1}f^{N}w$, which shows that
$\bar{h} \!\notin\! C_{\preceq} (\Gamma)$. Moreover, the inequality
$w^{-1} g^{M+1}f^{N}w (w^{-1} v) \!\prec~w^{-1}v$ can be written as
$h^{-1}\bar{h}h\prec id$. Finally, Lemma \ref{lema1} applied to the
crossing $(w^{-1}fw, w^{-1}gw; w^{-1}u,w^{-1}v,id)$ shows that, for
every $n \!\in\! \mathbb{N}$, the 5-tuple
$(w^{-1} fw, w^{-1}gw; w^{-1} u, w^{-1} v, w^{-1} g^{M+n}f^{N}w)$
is also a crossing. For $n \geq M$ we have $w^{-1}
g^{M+1}f^{N}w (w^{-1}v)\prec w^{-1} g^{M+n}f^{N}w $. Since
$w^{-1} g^{M+n}f^{N}w \prec w^{-1} v$, Lemma \ref{lema4} easily
implies that $w^{-1} g^{M+1}f^{N}w (w^{-1}v) C_{\preceq} (\Gamma) \prec
w^{-1}vC_{\preceq} (\Gamma)$, which yields
$h^{-1}\bar{h}h\notin C_{\preceq} (\Gamma)$.

\vsp\vsp

\noindent{\underline{Case II.}} One has $a\prec w^{-1}v$
and $w^{-1} g^m w \preceq a$, for all $m > 0$.

\vsp

We claim that, in this case, the proposition holds for the choice $h := a$
and $\bar{h} := w^{-1}g^{M+1}f^N w$. This may be checked in the very same
way as in Case I, by noticing that, if $a\prec w^{-1}v$ but $w^{-1}g^m
w \succeq a$ for all $m > 0$, then $(w^{-1}fw, w^{-1}gw; w^{-1}u,a,id)$
is a crossing.

\vsp\vsp

\noindent{\underline{Case III.}} One has $a\prec w^{-1}v$
and $w^{-1} g^m w\succ a$ for some $m > 0$. (Note
that the first condition follows from the second one.)

\vsp

We claim that, in this case, the proposition holds for the
choice $h:=a$ and $\bar{h}:=w \notin C_{\preceq} (\Gamma)$.
Indeed, we have $g^{m}w\succ ha$ (and $w\prec ha$). Since
$g^{m}w\prec v\prec a$, we also have $wa\prec a$, which means that
$h^{-1}\bar{h}h \prec id$. Finally, from Lemmas \ref{lema1} and
\ref{lema4}, we obtain \esp
$$wa C_{\preceq} (\Gamma)\preceq g^m w C_{\preceq} (\Gamma) \prec
v C_{\preceq} (\Gamma) \preceq a C_{\preceq} (\Gamma).$$
This implies that
\esp $a^{-1} w a C_{\preceq} (\Gamma) \prec C_{\preceq} (\Gamma)$, \esp which
means that $h^{-1} \bar{h} h \notin C_{\preceq} (\Gamma)$. $\hfill\square$

\begin{small}\begin{rem}\label{ex-free-clay}
In the context of Theorem \ref{final}, it is possible that one of the orders $\preceq^j$
may be not approximated by its conjugates despite being non-isolated. An illustrative
example of this fact for free groups is the subject of the Appendix of \cite{crossings}.
\end{rem}\end{small}


\subsection{Groups with finitely many Conradian orders}
\label{section-cristobal}

\vsp

\hspace{0.45cm} The starting point of this section is the following

\vsp

\begin{prop} \label{simple-alcanza}
{\em Let $\Gamma$ be a $C$-orderable group. If $\Gamma$ admits a Conradian
left-order having a countable neighborhood in $\mathcal{LO}(\Gamma)$, then
$\Gamma$ admits finitely many left-orders.}
\end{prop}

\vsp

Before showing this proposition, let us show how it leads to a

\vspace{0.3cm}

\noindent{\bf Proof of Theorem \ref{linnell-general}.} We provide
three different arguments (see \S \ref{section-rec} for still another
one that gives supplementary information).
First, as we saw in \S \ref{fin-uncount}, the proof reduces to show 
Proposition \ref{linda}. So, let $(\Gamma,\preceq)$ be a left-ordered group admitting
a finite-index subgroup restricted to which $\preceq$ is bi-invariant. By Proposition
\ref{listailor}, the left-order $\preceq$ is Conradian. By Proposition \ref{simple-alcanza},
if $\Gamma$ admits infinitely many left-orders, then all neighborhoods of $\preceq$ in
$\mathcal{LO}(\Gamma)$ are uncountable.

An alternative argument proceeds as follows. As was shown in \S \ref{section-crossings},
if a group admits a non-Conradian order, then it has uncountably many left-orders.
Assume that $\Gamma$ is left-orderable and all of its left-orders are Conradian. By
Proposition \ref{simple-alcanza}, if some of them has a countable neighborhood inside
$\mathcal{LO}(\Gamma) = \mathcal{CO}(\Gamma)$ (in particular, if $\mathcal{LO}(\Gamma)$
is countable), then $\Gamma$ admits only finitely many left-orders.

As a final argument, note that Proposition \ref{simple-alcanza} together with a convex
extension argument (see Section \ref{section-convex-extension}) show that, if
$\Gamma$ is a left-orderable group such that
$\mathcal{LO} (\Gamma)$ has an isolated point $\preceq$, then the Conradian
soul $C_{\preceq}(\Gamma)$ cannot have infinitely many left-orders. If $C_{\preceq}(\Gamma)$
is trivial (resp. if it is nontrivial and admits finitely many left-orders), then Proposition
\ref{primero} (resp. Proposition \ref{final}) yields the existence of a left-order $\preceq_*$
on $\Gamma$ that is accumulated by its conjugates. As we have already remarked, the closure
of the orbit under the conjugacy action of such a left-order is uncountable. $\hfill\square$

\vspace{0.385cm}

\noindent{\bf Proof of Proposition \ref{simple-alcanza}.} Let $\Gamma$ be a
group admitting a Conradian order $\preceq$ having a {\em countable} neighborhood
in $\mathcal{LO}(\Gamma)$, say
$$V_{f_1} \cap \ldots \cap V_{f_k} = \big\{ \!\preceq' \!: \esp
f_i \succ' id \esp \mbox{ for all } \esp i \!\in\! \{1,\ldots,k\} \big\}.$$

\vsp\vsp

\noindent{\underline{Claim (i).}} The chain of $\preceq$-convex subgroups is finite.

\vsp

Otherwise, there exists an infinite ascending or descending chain of convex jumps
$\Gamma_{g_n} \lhd \Gamma^{g_n}$ so that $f_m \notin \Gamma^{g_n} \setminus \Gamma_{g_n}$
for every \esp $m,n$. \esp As in the proof of Proposition \ref{first-tararin}, for each \esp
$\iota \!=\! (i_1,i_2\ldots) \!\in\! \{-1,+1\}^\mathbb{N}$ \esp let us
define the left-order \esp $\preceq_{\iota}$ on $\Gamma$ by:

\vsp

\noindent -- \esp\esp $P_{\preceq_{\iota}}^+ \bigcap
\big( \Gamma \setminus (\Gamma^{g_n} \setminus \Gamma_{g_n}) \big) =
P_{\preceq}^+ \bigcap \big( \Gamma \setminus (\Gamma^{g_n} \setminus \Gamma_{g_n}) \big)$,
for each $n \in \mathbb{N}$;

\vsp

\noindent -- \esp\esp $P_{\preceq_{\iota}} \cap (\Gamma^{g_n} \setminus \Gamma_{g_n}) 
= P_{\preceq}^+ \cap (\Gamma^{g_n} \setminus \Gamma_{g_n})$ (resp.
$P_{\preceq_{\iota}} \cap (\Gamma^{g_n} \setminus \Gamma_{g_n}) 
= P_{\preceq}^- \cap (\Gamma^{g_n} \setminus \Gamma_{g_n})$) if $i_n = +1$
(resp. $i_n = -1$).

\vsp

\noindent This yields a continuous embedding of the Cantor
set \esp $\{-1,+1\}^{\mathbb{N}}$ \esp into $\mathcal{LO}(\Gamma)$.
Moreover, since \esp $f_m \notin \Gamma^{g_n} \setminus \Gamma_{g_n}$
\esp for every \esp $m,n$, \esp the image of this embedding is contained in
$V_{f_1} \cap \ldots \cap V_{f_k}$. This proves the claim.

\vsp\vsp

\noindent{\underline{Claim (ii).}} Denote by \esp
$\{id\} = \Gamma^k \lhd \Gamma^{k-1} \lhd \ldots \lhd \Gamma^0 = \Gamma$ \esp
the chain of {\em all} $\preceq$-convex subgroups. Then each quotient
$\Gamma^{i-1} / \Gamma^i$ is torsion-free, rank-1 Abelian.

\vsp

If the rank of some $\Gamma^{i-1} / \Gamma^i$ were larger than 1, then the
induced left-order on the quotient would be non-isolated in the space of
left-orders of $\Gamma^{i-1} / \Gamma^i$. This would allow to produce
--by a convex extension type procedure-- uncountably many
left-orders on any given neighborhood of $\preceq$, which
is contrary to our hypothesis.

\vsp\vsp

\noindent{\underline{Claim (iii).}} In the series above,
the group $\Gamma^{k-2}$ is not bi-orderable.

\vsp

First note that $\Gamma^{k-2}$ is not Abelian. Otherwise, it would have rank 2.
This would imply that every neighborhood of the restriction of $\preceq$ to
$\Gamma^{k-2}$ is uncountable, which implies --by convex extension--
the same property for $\preceq$.

Now as in the case of Proposition \ref{third-tararin}, if $\Gamma^{k-2}$ were bi-orderable,
then it would be contained in the affine group $\mathrm{Aff}_+ (\mathbb{R})$. The space
of left-orders of a non-Abelian countable group inside $\mathrm{Aff}_+ (\mathbb{R})$ was
roughly described in \S \ref{ejemplificando-2}: it is homeomorphic to the Cantor set.
(See also \S \ref{section-solvable-spaces}.)
In particular, no neighborhood of the restriction of $\preceq$ to $\Gamma^{k-2}$ is
countable, which implies --by convex extension-- that the same is true for $\Gamma$.
For sake of completeness, we give an explicit sequence of approximating left-orders.
To do this, note that, for some $q > 0$, the group $\Gamma^{k-2}$ can be
identified with the group whose elements are of the form
$$(k,a) \sim
\left(
\begin{array}
{cc}
q^k & a  \\
0   & 1  \\
\end{array}
\right),$$
where $a \!\in\! \Gamma^{k-1}$ and $k \!\in\! \mathbb{Z}$. Let $(k_1,a_1),\ldots,(k_n,a_n)$
be an arbitrary family of $\preceq$-positive elements indexed in such a way that
$k_1 = k_2 = \ldots = k_r = 0$ and $k_{r+1} \neq 0, \ldots, k_n \neq 0$ for some
$r \!\in\! \{1,\ldots,n\}$. Four cases are possible:

\vspace{0.1cm}

\noindent (i) \esp \esp \esp $a_1 > 0, \ldots, a_r > 0$ \esp \esp and
\esp \esp $k_{r+1} > 0, \ldots, k_n > 0$;

\vspace{0.1cm}

\noindent (ii) \esp \esp \esp $a_1 < 0, \ldots, a_r < 0$ \esp \esp and
\esp \esp $k_{r+1} > 0, \ldots, k_n > 0$;

\vspace{0.1cm}

\noindent (iii) \esp \esp \esp $a_1 > 0, \ldots, a_r > 0$ \esp \esp and
\esp \esp $k_{r+1} < 0, \ldots, k_n < 0$;

\vspace{0.1cm}

\noindent (iv) \esp \esp \esp $a_1 < 0, \ldots, a_r < 0$ \esp \esp and
\esp \esp $k_{r+1} < 0, \ldots, k_n < 0$.

\vspace{0.1cm}

\noindent As in \S \ref{ejemplificando-2}, for each irrational number
\esp $\varepsilon$, \esp let $\preceq_{\varepsilon}$ be the left-order
on $\preceq_{\varepsilon}$ whose positive cone is
$$P_{\preceq_{\varepsilon}}^+ =
\big\{(k,a) \! : \esp \esp \esp q^k + \varepsilon a > 1 \big\}.$$
In case (i), for $\varepsilon$ positive and very small, the left-order
$\preceq_{\varepsilon}$ is different from $\preceq$ but still makes positive all the
elements $(k_i,a_i)$. The same is true in case (ii) for $\varepsilon$ negative and near
zero. In case (iii), this still holds for the order $\bar{\preceq_{\varepsilon}}$ when
$\varepsilon$ is negative and near zero. Finally, in case (iv), one needs to consider
again the order $\bar{\preceq_{\varepsilon}}$, but for $\varepsilon$ positive and
small. Now letting $\varepsilon$ vary over a Cantor set formed by irrational
numbers\footnote{Take for example the set of numbers of the form
$\sum_{i \geq 1} \frac{i_k}{4^k}$, where $i_k \!\in\! \{0,1\}$, and
translate it by $\sum_{j \geq 1} \frac{2}{4^{j^2}}$.} very close to $0$
(and which are positive or negative according to the case), this shows that
the neighborhood of (the restriction to $\preceq$ of) $\preceq$ consisting of
the left-orders on $\Gamma^{k-2}$ that make positive all the elements $(k_i,a_i)$
contains a homeomorphic copy of the Cantor set.

\vsp\vsp

\noindent{\underline{Claim (iv).}} The series of Claim (ii) is normal
(hence rational) and no quotient $\Gamma^{i-2} / \Gamma^i$ is bi-orderable.

\vsp

By Theorem \ref{tararin-theorem}, the group $\Gamma^{k-2}$ admits a unique rational
series, namely \esp $\{ id \} \lhd \Gamma^{k-1} \lhd \Gamma^{k-2}$. \esp Since for
every $h \in \Gamma^{k-3}$ the series \esp $\{ id \} \lhd h\Gamma^{k-1}h^{-1} \lhd
h\Gamma^{k-2}h^{-1}$ \esp is also rational for $\Gamma^{k-2}$, they must coincide.
Hence, the rational series
$$\{ id \} \lhd \Gamma^{k-1} \lhd \Gamma^{k-2} \lhd \Gamma^{k-3}$$
is normal. Moreover, proceeding as in Claim (iii) with the induced left-order
on $\Gamma^{k-3} / \Gamma^{k-1}$, one readily checks that this quotient is
not bi-orderable. Once again, Theorem \ref{tararin-theorem} implies that the
rational series of $\Gamma^{k-3}$ is unique... Arguing in this way,
the claim follows.

\vsp\vsp

We may now conclude the proof of the proposition. Indeed, we have
shown that, if $\Gamma$ is a group having a $C$-order with a countable
neighborhood in $\mathcal{LO}(\Gamma)$, then $\Gamma$ admits a rational series
$$\{id\} = \Gamma^k \lhd \Gamma^{k-1} \lhd \ldots \lhd \Gamma^0 = \Gamma$$
such that no quotient $\Gamma^{i-2} / \Gamma^i$ is bi-orderable. By Theorem
\ref{tararin-theorem}, $\Gamma$ has only finitely many left-orders.
$\hfill\square$

\vsp\vsp\vsp\vsp\vsp

We now turn to the study of the space of Conradian orders. The next result from
\cite{rivas} is the analogue of Tararin's theorem describing left-orderable groups
with finitely many left-orders; see \S \ref{fin-uncount}.
\index{Space of!Conradian orders}

\vsp

\index{Rational series}
\begin{thm} \label{lata}
{\em If a $C$-orderable group $\Gamma$ has only finitely many
$C$-orders, then it has a unique (hence normal) rational series
$\{ id \} \!=\! \Gamma^k \lhd \Gamma^{k-1} \! \lhd \ldots \! \lhd \Gamma^0
\!=\! \Gamma$. In this series, no quotient $\Gamma^{i-2}/\Gamma^{i}$ is
Abelian. Conversely, if $\Gamma$ is a group admitting a normal rational series
$\{ id \} = \Gamma^k \lhd \Gamma^{k-1} \lhd \ldots \lhd \Gamma^0 = \Gamma$
such that no quotient $\Gamma^{i-2} / \Gamma^{i}$ is Abelian, then the number
of $C$-orders on $\Gamma$ is $2^k$.}
\end{thm}

\noindent{\bf Proof.} The proof will be divided into four independent claims.

\vsp\vsp

\noindent{\underline{Claim (i).}} If $\Gamma$ is a $C$-orderable group admitting
only finitely many $C$-orders, then for every $C$-order $\preceq$ on
$\Gamma$, the sequence of $\preceq$-convex subgroups is a rational series.

\vsp

Indeed, for each convex jump $\Gamma_g \lhd \Gamma^g$, we may flip the left-order on
$\Gamma_g$ to produce a new left-order (see Example \ref{ex-flipping}) which
is still Conradian (see Exercise \ref{ejercicito}).
If there were infinitely many $\preceq$-convex subgroups,
then this would allow to produce infinitely many
$C$-orders on $\Gamma$, contrary to our hypothesis. Let then
$$\{id\} = \Gamma^k \subsetneq \Gamma^{k-1}
\subsetneq \ldots \subsetneq \Gamma^0 = \Gamma$$
be the sequence of {\em all} $\preceq$-convex subgroups.
As in the proof of Proposition \ref{first-tararin},
$\Gamma^i$ is normal in $\Gamma^{i-1}$, and $\Gamma^{i-1} / \Gamma^i$
is torsion-free Abelian. The rank of this quotient must be 1, as otherwise
it would admit uncountably many orders, which would allow to produce
--by convex extension-- uncountably many $C$-orders on $\Gamma$.

\vsp\vsp

\noindent{\underline{Claim (ii).}} If a left-orderable group admits only
finitely many $C$-orders, then it has a unique (hence normal)
rational series.

\vsp

The proof is almost the same as that of Proposition \ref{second-tararin}. We just
need to change the word ``left-order'' by ``$C$-order'' along that proof, and
replace the (crucial) use of Proposition \ref{rel-convex} by Proposition
\ref{rel-convex-conrad}.

\vsp\vsp

\noindent{\underline{Claim (iii).}} If a group $\Gamma$ with a normal rational series
\esp $\{id\}=\Gamma^k \lhd \Gamma \lhd \ldots \lhd \Gamma^0 = \Gamma$ \esp
admits only finitely many $C$-orders, then no quotient
$\Gamma^{i-2} / \Gamma^i$ is Abelian.

\vsp

First note that every group admitting a rational series is $C$-orderable. Actually,
using the rational series above, one may produce $2^k$ Conradian orders on $\Gamma$.
If one of the quotients $\Gamma^{i-2} / \Gamma^i$ were Abelian, then it would have
rank 2, hence it would admit uncountably many left-orders. This would allow to
produce --by convex extension-- uncountably many $C$-orders on $\Gamma$.

\vsp\vsp

\noindent{\underline{Claim (iv).}} If a group $\Gamma$ has a normal rational series \esp
$\{id\}=\Gamma^k \lhd \Gamma \lhd \ldots \lhd \Gamma^0 = \Gamma$ \esp such that no
quotient $\Gamma^{i-2} / \Gamma^i$ is Abelian, then this series coincides with that
formed by the $\preceq$-convex subgroups, where $\preceq$ is any $C$-order
on $\Gamma$. In particular, such a series is unique.

\vsp

As we have already seen, the rational series above leads to $2^k$ Conradian
left-orders. We have to prove that these are
the only possible $C$-orders on $G$. To show this, let $\preceq$ be a
$C$-order on $\Gamma$. By Claim (iii), there exist non-commuting elements
$g \in \Gamma^{k-1}$ and $h \in \Gamma^{k-2} \setminus \Gamma^{k-1}$.
Denote the Conrad homomorphism
of the group $\, \langle g,h \rangle \,$
(endowed with the restriction of $\, \preceq \,$) by $\tau$. Then we have
$\tau(g) = \tau(hgh^{-1}) \neq 0$. Since $\Gamma^{k-1}$ is rank-1
Abelian, $hgh^{-1}$ must be equal to $g^s$ for some rational number
$s\not=1$. Hence, $\tau(g)=s\tau(g)$, which implies that $\tau(g)=0$.
Therefore, $\, g^n \prec |h| \,$ for every $n \in \mathbb{Z}$, where
$|h| := \max \{h^{-1},h\}$. Since $\Gamma^{k-2} / \Gamma^{k-1}$ has rank 1, this
actually holds for every $h \neq id$ in $\Gamma^{k-2} \setminus \Gamma^{k-1}$.
Therefore, $\Gamma^{k-1}$ is $\preceq$-convex in $\Gamma^{k-2}$.

\vsp

Repeating the argument above, though now with $\Gamma^{k-2}/\Gamma^{k-1}$
and $\Gamma^{k-3}/\Gamma^{k-1}$ instead of $\Gamma^{k-1}$ and $\Gamma^{k-2}$,
respectively, we see that the rational series we began with is no other
than the series given by the $\preceq$-convex subgroups. Since
each $\Gamma^{i-1}/\Gamma^i$ is rank-1 Abelian, if we choose
$g_i \in \Gamma^{i-1} \setminus \Gamma^i$ for each $i$, then
any $C$-order on $\Gamma$ is completely determined by the
signs of these elements. This shows the claim, and concludes
the proof of Theorem \ref{lata}. $\hfill\square$

\vsp

\index{Rational series}
\begin{small}\begin{ejem} \label{solo-cuatro}
The Baumslag-Solitar group \esp
$BS(1,\ell) = \langle a, b \!: aba^{-1} = b^{\ell} \rangle$, \esp $\ell \geq 2$,
admits the rational series
$$\{ id \} \esp \lhd \esp b^{\mathbb{Z}[\frac{1}{\ell}]} := \langle c \!: c^{\ell^i}
= b \mbox{ for some integer } i > 0  \rangle \esp \lhd \esp BS(1,\ell),$$
which satisfies the conditions of Theorem \ref{lata}. Therefore, $BS(1,\ell)$ 
admits four $C$-orders --all of which are bi-invariant--, though
its space of left-orders is uncountable (see \S \ref{ejemplificando-2}).
The reader is referred to \S \ref{section finite-rank-solvable} for more details on this example.
\end{ejem}

\index{Rational series}
\begin{ejem} Examples of groups having exactly $2^k$ left-orders
(hence $2^k$ Conradian orders) were given in Example \ref{ups-ups}.
Namely, one may consider \esp $K_k = \langle a_1,\ldots,a_k \mid R_k\rangle,$
\esp where the set of relations $R_k$ is
$$\;a_{i+1}^{-1} a_i a_{i+1} = a_i^{-1} \quad \mbox{if} \quad i < k,
\qquad a_ia_j = a_ja_i \quad \mbox{if} \quad |i - j| \geq 2.$$
The existence of groups with $2^k$ Conradian orders but infinitely many
(hence uncountably many) left-orders is more subtle. As we have seen in the
preceding example, for $n=k$ this is the case of the Baumslag-Solitar
groups $BS(1,\ell)$ for $\ell \geq 2$. To construct examples for higher $k$
having $BS(1,\ell)$ as a quotient by a normal convex subgroup, we choose an
{\em odd} integer $\ell \geq 3$, and we let $C_n(\ell)$ be the group
$$\big\langle c,b, a_1,\ldots,a_n \mid  cbc^{-1}=b^{\ell}\,,\;
ca_i=a_ic\,,\; ba_nb^{-1}= a_n^{-1}\,,\; ba_i=a_ib \esp \text{ if } \esp i\not=n
\,  ,\;R_n \big\rangle.$$
This corresponds to the set $\Z \times \Z[\frac{1}{3}] \times \Z^{n}$
\esp endowed with the product rule
\begin{multline*}
\Big( c, \,\frac{m}{\ell^{k}}\,, a_1, \ldots, a_n\Big) \cdot
\Big(c^\prime, \,\frac{m^\prime}{\ell^{k^\prime}}\,, a_1^\prime, \ldots, a_n^\prime \Big)
= \\
= \Big( c + c^\prime, \; \ell^{c}\frac{m^\prime}{\ell^{k^\prime}}+\frac{m}{\ell^{k}}\;,
(-1)^{a_2} a_1^\prime + a_1, \ldots,\, (-1)^{a_n} a^\prime_{n-1} + a_{n-1}, \esp
(-1)^{m} a^\prime_n+a_n \Big).
\end{multline*}
Note that this is well-defined, as \esp $(-1)^{m}\!=\!(-1)^{\bar{m}}$
\esp whenever \esp $m/\ell^k = \bar{m}/\ell^{\bar{k}}$ (it is at this step where the fact that
$\ell$ is odd becomes important). The group $C_n (\ell)$ admits the rational series
$$\{ id \} \lhd \langle a_1 \rangle \lhd \langle a_1, a_2 \rangle \lhd
\ldots \lhd \langle a_1,\ldots,a_n \rangle \lhd \big\langle a_1,\ldots,a_n,
b^{\mathbb{Z}[\frac{1}{\ell}]} \big\rangle \lhd C_n(\ell).$$
By Theorem \ref{lata}, it admits exactly $2^{n+2}$ Conradian orders. However, it
has $BS(1,\ell)$ as quotient by the normal convex subgroup $K_n$. Since $BS(1,\ell)$
admits uncountably many (left) left-orders, the same is true for $C_n (\ell)$.
\end{ejem}\end{small}

\vsp

We close this section with a result (also taken from
\cite{rivas}) to be compared with Theorem \ref{linnell-general}.

\vsp

\begin{thm} \label{cristobal-one}
{\em Every $C$-orderable group admits either finitely many or uncountably
many $C$-orders. In the last case, none of these left-orders
is isolated in the space of $C$-orders.}
\end{thm}
\index{Order!isolated}

\vsp

To prove this theorem, we need the lemmas below.

\vsp

\index{Rational series}
\begin{lem} {\em If $\Gamma$ is $C$-orderable group such that $\mathcal{CO}(\Gamma)$
has an isolated point $\preceq$, then the family of $\preceq$-convex
subgroups (is finite and) is a rational series such that no quotient
of the form $\Gamma^{i-2} / \Gamma^i$ is Abelian.}
\end{lem}

\noindent{\bf Proof.} As in Claim (i) of Proposition \ref{simple-alcanza},
the family of $\preceq$-convex subgroups is finite, say
$$\{id\} = \Gamma^k \subsetneq \Gamma^{k-1} \subsetneq \ldots
\subsetneq \Gamma^0 = \Gamma.$$
Since $\preceq$ is Conradian, $\Gamma^i$ is normal in $\Gamma^{i-1}$ for each $i$.
The proofs of that $\Gamma^{i-1} / \Gamma^i$ has rank-1 and no quotient
$\Gamma^{i-2} / \Gamma^i$ is Abelian are similar to those of Theorem
\ref{lata}, and we leave them to the reader. $\hfill\square$

\vsp\vsp\vsp

\index{Order!isolated}
\begin{lem} {\em For any $C$-orderable group whose space of
$C$-orders has an isolated point $\preceq$, the rational
series formed by the $\preceq$-convex subgroups is normal.}
\end{lem}

\noindent{\bf Proof.} The proof is by induction on the
length $k$ of the rational series. For $k=1$, there is nothing
to prove. For $k=2$, the series is automatically normal.
Assume that the claim of the lemma holds for $k$, and let
\begin{equation}\label{tt}
\{id\} = \Gamma^{k+1} \lhd  \Gamma^k \lhd \ldots
\lhd \Gamma^1 \lhd \Gamma^0 = \Gamma
\end{equation}
be the rational series of length $k+1$ associated to a $C$-order
on a group $\Gamma$ that is isolated in $\mathcal{CO}(\Gamma)$.
Note that the truncated chain of length $k$
\begin{equation}\label{ttt}
\{id\} = \Gamma^{k+1} \lhd  \Gamma^k \lhd \ldots \lhd \Gamma^1
\end{equation}
is a rational series for $\Gamma^1$. Moreover, this series is associated to
a $C$-order on $\Gamma^1$ (namely, the restriction of $\preceq$) that is
isolated in $\mathcal{CO}(\Gamma^1)$ (otherwise, $\preceq$ would be
non-isolated in $\mathcal{CO}(\Gamma)$). By the inductive hypothesis,
this series is normal. By the preceding lemma, for each
$i \!\in\! \{3,\ldots,k+1\}$, the quotient $\Gamma^{i-2} / \Gamma^i$
is non-Abelian. We are hence under the hypothesis of Theorem \ref{lata},
which allows us to conclude that this is the unique rational series of $\Gamma^1$.

Now, since $\Gamma^1$ is normal in $\Gamma$, for each $h \in \Gamma$,
the conjugate series
$$\{id\} = h \Gamma^{k+1} h^{-1} \lhd  h \Gamma^k h^{-1} \lhd \ldots \lhd
h \Gamma^1 h^{-1} = \Gamma^1$$
is also a rational series for $\Gamma^{1}$. By the uniqueness above,
this series coincides with (\ref{ttt}). Therefore, (\ref{tt})
is a {\em normal} rational series. $\hfill\square$

\vsp\vsp\vsp\vsp

The proof of Theorem \ref{cristobal-one} is now at hand. Indeed, the two preceding
lemmas imply that, if a $C$-orderable group admits an isolated $C$-order,
then it has a normal rational series satisfying the hypothesis of Theorem
\ref{lata}, thus it has finitely many $C$-orders. If, otherwise, no $C$-order
is isolated in the space of $C$-orders, then this is a Hausdorff, totally
disconnected, topological space without isolated points,
hence uncountable (see \cite[Theorem 2-80]{HY}).

\vsp

\begin{small} \begin{ejer} By slightly extending the arguments above,
show the following analogue of Proposition \ref{simple-alcanza}: If a
$C$-orderable group admits infinitely many $C$-orders, then every
neighborhood of such a left-order in the space of Conradian orders
is uncountable.
\end{ejer} \end{small}


\section{An Application: Ordering Solvable Groups}
\label{section-solvable-spaces}

\hspace{0.45cm} Following \cite{RT}, we will show 
that the space of left-orders of a countable left-orderable virtually-solvable group has no
isolated point, except for the cases where it is finite, which are described in \S \ref{fin-uncount}.
This result requires both algebraic and dynamical developments. As a major particular case,
in \S \ref{section finite-rank-solvable}, 
we focus on finite-rank solvable groups and their finite-index extensions, for which
the result will follow from a classification (up to semiconjugacy) of all actions on the line
(without global fixed points). As concrete relevant examples, at the end of the subsection, 
we treat the cases of the Baumslag-Solitar groups and the groups $Sol$, for which we 
can give a full description of the corresponding spaces of left-orders. 

In \S \ref{section solvable (general)}, 
we deal with the much more difficult case of infinite-rank groups. The approach we develop therein only requires 
a local description of the dynamics of the left multiplication on a left-ordered solvable group around $id$. However, 
a complete classification (again, up to semiconjugacy) of all actions (with no global fixed point) on the real line 
is still available in that case; see \cite{BMRT-solvable} (see also Remark \ref{rem laminar solvable}). 


\subsection{The space of left-orders of finite-rank solvable groups}
\label{section finite-rank-solvable}
\index{Group!finite-rank solvable}

Recall that a group $\Gamma$ is said to be {\bf{\em  virtually finite-rank solvable}}
if it contains a finite-index subgroup $\tilde{\Gamma}$ that admits a normal series
$$\{id\} = \widetilde\Gamma^n \lhd \widetilde\Gamma^1\lhd \ldots \lhd \widetilde\Gamma^0
= \widetilde\Gamma$$
in which every quotient
$\widetilde{\Gamma}^{i-1}/\widetilde{\Gamma}_i$ is finite-rank Abelian.\footnote{In the case
where such a series can be taken so that each $\Gamma^{i-1} / \Gamma^i$ is cyclic, the
group is said to be {\bf {\em virtually polycyclic}}.}
(Note that such a group $\Gamma$ is necessarily countable.) The number $\sum_i
rank(\widetilde{\Gamma}^{i-1}/\widetilde{\Gamma}^i)$ is independent of both the finite-index
subgroup and the normal series chosen. (In particular, we can --and we will-- take
$\widetilde \Gamma$ as being normal in $\Gamma$.) We call this number the {\bf \em Hirsch rank} 
of $\Gamma$ or simply the {\bf \em rank} of $\Gamma$. (Note that when that restricted to Abelian groups, 
the Hirsch rank coincides with the usual rank.) We leave to the reader the task of checking that this number 
strictly decreases when passing to either an infinite-index subgroup or to a quotient by an infinite subgroup. 
(See \cite{robinson} in case of problems.) 
\index{rank}

\vsp

\begin{small}\begin{ejer}\label{maximal-finite-rank}
Show that every left-order on a virtually finite-rank solvable group admits a maximal proper convex
subgroup (despite the fact that such a group may fail to be finitely-generated).

\noindent\underline{Hint.} Proceed by induction on the rank, noting that if $G \subset H$ are distinct convex
subgroups, then the rank of $G$ is strictly smaller than that of $H$.
\end{ejer}\end{small}

\vsp

The main result of this section states as follows.

\vsp

\begin{thm}\label{teo finite rank solvable} {\em The space of left-orders of
a virtually finite-rank solvable group is either finite or a Cantor
set.}
\end{thm}

\vsp

The first step of the proof concerns left-orders induced from non-Abelian affine
actions. (Compare Theorem \ref{primero}.)

\vsp

\begin{prop}\label{prop affine} {\em Let $\Gamma$ be a subgroup of the affine
group endowed with a left-order $\preceq $ induced (in a dynamical-lexicographic way)
from its affine action on the real line. If $\Gamma$ is non-Abelian, then $\preceq$ is an
accumulation point of its set of conjugates.}
\end{prop}

\noindent{\bf Proof.} First note that, as affine homeomorphisms fix at most one point,
the dynamical-lexicographic order $\preceq$ is completely determined by the first two
comparison points, that we denote $x_1,x_2$. (In the case of a single comparison 
point $x_1$, we let $x_2 : = x_1$.) We will assume that the signs that we chose 
for $x_1$ and $x_2$ (in the sense of \S \ref{general-3}) are both positive, since 
this is the only case we use below (the remaining cases are analogous and 
are left to the reader).

By assumption, $\Gamma$ contains both nontrivial homotheties and
nontrivial translations. It follows that the translations in $\Gamma$ form a
subgroup with dense orbits, hence the set  of points that are fixed by some
nontrivial homothethy in $\Gamma$ is dense in $\mathbb{R}$. Therefore,
given any two distinct points in $\mathbb{R}$, there is a nontrivial homothethy
whose unique fixed point lies between them. As a consequence, for any
pair of comparison points $y_1,y_2$ such that $y_1 \neq x_1$, the
induced left-order $\preceq'$ is different from $\preceq$.

We next show that $y_1,y_2$ may be chosen so that $\preceq'$ is close to $\preceq$.
Given a finite set $\mathcal{G} \subset \Gamma$ of $\preceq$-positive elements, we write
it as a disjoint union $\mathcal{G} = \mathcal{G}_1 \sqcup \mathcal{G}_2$, where
$\mathcal{G}_1$ is the subset of elements of $\mathcal{G}$ lying in the stabilizer of
$x_1$ in $\Gamma$. Let $I$ denote the open interval with endpoints $x_1$ and $x_2$.
On the one hand, since $\mathcal{G}_2$ is finite, there is a small neighborhood $U$ of
$x_1$ such that $f(x)>x$ for every $x\in U$ and every $f \in \mathcal{G}_2$. On the
other hand, for every $f \in \mathcal{G}_1$, we have $f(x)>x$ for every $x \in I$.
(Note that each $f \in \mathcal{G}_1$ is an homothethy.) Thus, if we choose any
$y_1$ in the nonempty open set $I \cap U$ (and $y_2$ arbitrary), then the resulting
left-order $\preceq'$ is such that all elements in $\mathcal{G}$ are still $\preceq'$-positive.
Finally, we can choose such a $y_1$ in the $\Gamma$-orbit of $x_1$, say
$y_1 = h(x_1)$. For this choice (and letting $y_2 := h (x_2)$), we have
that $\preceq'$ is the conjugate of $\preceq$ by $h$, as desired.
$\hfill\square$

\vspace{0.28cm}

\begin{cor}\label{cor:affine}
{\em Let $(\Gamma,\preceq)$ be a countable, left-ordered group. Suppose
there is a homomorphism $\Phi \!: \Gamma \to \mathrm{Aff}_+(\mathbb{R})$
with $\preceq$-convex kernel and non-Abelian image. Suppose further that the
dynamical realization of $(\Gamma,\preceq)$ is semiconjugate to the action
given by $\Phi$. Then $\preceq$ is non-isolated in $\mathcal{LO}(\Gamma)$.}
\end{cor}

\noindent{\bf Proof.} By Proposition \ref{prop convex and LO}, it suffices to deal with
the case where $\Phi$ is injective. Let $\varphi$ denote the semiconjugacy assumed by
hypothesis, and let $\Gamma_0$ be the stabilizer of $\varphi (0)$ in $\Phi (\Gamma)$. This
is an Abelian subgroup of $\Gamma$. We claim that it is $\preceq$-convex. Indeed, if
$h_1\prec g \prec h_2$, then $h_1(0) \leq g(0) \leq h_2(0)$, thus
$\varphi(h_1(0)) \leq \varphi (g(0)) \leq \varphi (h_2(0))$,
and hence $\Phi (h_1)(\varphi(0))\leq \Phi(g) (\varphi(0))\leq \Phi (h_2)(\varphi(0))$.
In particular, if $h_1,h_2$ lie in $\Gamma_0$, then $\Phi (g) (\varphi(0)) = \varphi(0)$,
that is, $g$ also lies in $\Gamma_0$.

If $\Gamma_0$ is trivial, then Proposition \ref{prop affine}
directly applies, since in this case $\preceq$ coincides with the left-order
induced from $\varphi(0)$ in the action given by $\Phi$.  If $\Gamma_0$
has rank $1$, then the restriction of $\precede$ to
$\Gamma_0$ is completely determined by the sign of any nontrivial element
therein, say $\Phi (h) \in \Gamma_0$, with $h \succ id$.
As $\Phi (h)$ is a nontrivial homothethy, there
exists $x \in \mathbb{R}$ such that $\Phi (h)(x) > x$. It follows that $\preceq$
coincides with the left-order induced from the action $\Phi$ using
the comparison points $x_1 := \varphi(0)$ and $x_2 := x$.
Therefore, Proposition \ref{prop affine} still allows concluding that $\preceq$
is non-isolated. Finally, the case where $\Gamma_0$ has rank $>1$ is slightly
different, as we cannot argue that $\preceq$ is completely induced from the affine
action. However, by \S \ref{ejemplificando-1}, the restriction of $\preceq$ to
$\Gamma_0$ is non-isolated. Therefore, by convex extension, $\preceq$
itself is non-isolated, as desired. $\hfill\square$

\vspace{0.3cm}

To proceed with the proof of Theorem \ref{teo finite rank solvable}, we need
some general results on the structure of finite-rank solvable groups. If $\Gamma$
is a virtually finite-rank solvable group that is, moreover, torsion-free, then $\Gamma$ 
contains a finite-index subgroup
$\tilde{\Gamma}$ whose commutator subgroup $[\tilde\Gamma,\tilde\Gamma]$
is nilpotent \cite{Ra,robinson}.
Let $R$ be a maximal nilpotent subgroup of $\widetilde{\Gamma}$. By maximality,
$R$ is a {\bf {\em characteristic subgroup}} of $\widetilde{\Gamma}$ (that is, it remains 
invariant under isomorphisms). In particular, it is normal in $\Gamma$. 
Moreover, it is unique (see Exercise \ref{ejer-prod-gr} below). It is sometimes called the
{\bf{\em nilpotent radical}} of $\widetilde\Gamma$.
\index{Group!nilpotent radical}

\vsp

\begin{small}\label{ejer-prod-gr}
\begin{ejer} Let $\Gamma$ be a group and $G,H$ two normal nilpotent subgroups.
Show that the set $GH := \{gh \!: g \in G,\; h \in H \}$ is a nilpotent subgroup of $\Gamma$.
\end{ejer}
\end{small}

\vsp

Theorem \ref{nilpotent-cantor} implies Theorem \ref{teo finite rank solvable}
in the case where $R$ has finite index in $\widetilde{\Gamma}$ (in particular,
when the rank of $\Gamma$ is $1$). Hence, in what follows, we assume that $\widetilde{\Gamma}/R$ is infinite.
We proceed by induction on the rank of $\Gamma$. Thus, we assume that Theorem \ref{teo finite rank solvable} 
holds for every virtually finite-rank solvable group having smaller rank than that of $\Gamma$. Let $\preceq$ be 
a left-order on $\Gamma$. Consider its dynamical realization, and denote by $\Gamma_0 \subset R$ the set of
elements in $R$ having fixed points. By Exercise \ref{ejer:es-un-subgrupo}, $\Gamma_0$ is a normal subgroup 
of $R$. Since $R$ is normal in $\Gamma$, we have that $\Gamma_0$ is normal in $\Gamma$ as well. 
The following lemma implies that $\Gamma_0$ has a global fixed point.
(Compare Exercise \ref{puntos-fijos-nilp}.)
\vsp

\begin{lem} \label{lem global fix pt for nilp} {\em Assume that a
nilpotent group with finite rank acts by orientation-preserving
homeomorphisms of the real line. If every element admits
fixed points, then there is a global fixed point for the action.}
\end{lem}

\noindent{\bf Proof.}
If the nilpotence length of the underlying group $G$ is 0, then the group is trivial,
and there is nothing to prove. We continue by induction on the nilpotent
length, denoting the center of $G$ by $H$. This is a finite-rank
Abelian group, hence it contains a subgroup $H_0$ isomorphic to a certain
$\Z^d$ such that $H/H_0$ is a torsion group. It follows that the (closed) set
$Fix := Fix (H_0)$ of fixed points of $H_0$ is nonempty. Since
$H / H_0$ is torsion, $Fix$ coincides with the set of fixed points of $H$.
The complement of $Fix$ is a disjoint union $\bigsqcup_i I_i$ of open
intervals $I_i$. Moreover, since $H \lhd G$, we have that $Fix$
is $G$-invariant. In particular, the
intervals in the complement of $Fix$ are permuted by $G$. Furthermore,
since every element of $G$ has fixed points, we have that every element in
$G$ must fix some point in $Fix$. Let us now extend in a piecewise-affine
manner the action of $G$ on $Fix$ to the complementary intervals. Doing this,
we obtain a new action of $G$ on $\mathbb{R}$ which factors throughout $G/H$.
Since this action coincides with the original one on $Fix$, every element
of $G$ admits fixed points. We can hence apply the induction hypothesis,
thus concluding that $G/H$ has a global fixed point in $Fix$, hence $G$
has a global fixed point. $\hfill\square$

\vspace{0.28cm}
Now, since $R$ is nilpotent, every left-order on it is Conradian
(see the discussion before Theorem \ref{nilpotent-cantor}). Using
Corollary \ref{mas} (more precisely, by Exercise \ref{puntos-fijos-nilp}),
we obtain that $\Gamma_0$ contains $[R,R]$, and $R/\Gamma_0$
is torsion-free. Moreover, since $\Gamma_0$ is normal in $\Gamma$, 
the set $Fix(\Gamma_0)$ is $\Gamma$-invariant, hence $\Gamma_0$  
admits an integer-indexed sequence $(z_n)_{n \in \mathbb{Z}}$ of 
global fixed points going from $-\infty$ to $\infty$.

We split the induction argument into two cases.

\vsp\vsp\vsp

\noindent{\underline{Case I.}} Either $R/\Gamma_0$ is trivial, or
it has rank $1$ and the conjugacy action of $\widetilde{\Gamma}$
on it is by multiplication by $\pm1$.

\vsp\vsp

In this case, we start by establishing the following claim.

\vspace{0.25cm}

\noindent{\underline{Claim (i).}}
The quotient $\widetilde{\Gamma}/\Gamma_0$ is Abelian.

\vsp

Indeed, if $R / \Gamma_0$ is trivial, then this follows from that
$[\widetilde{\Gamma},\widetilde{\Gamma}]\subseteq R$. Otherwise,
assume for a contradiction that $\widetilde{\Gamma}$ does not centralize
$\tilde \Gamma / \Gamma_0$. As $\tilde \Gamma$ centralizes $\tilde \Gamma / R$,
this means that $\tilde \Gamma$ does not centralize $R / \Gamma_0$. Hence, there
are $f \in R \setminus \Gamma_0$ and $g \in \widetilde{\Gamma}$ such that,
modulo $\Gamma_0$, one has the equality $g f g^{-1} = f^{-1}$. Now, since
$f$ acts without fixed points, changing $f$ by $f^{-1}$ if necessary, we can
assume that $f (y) > y$ for every $y \in \mathbb{R}$. Thus, if we let $x$ be in
the set of fixed points of $ \Gamma_0$, we have that $g f g^{-1}(x) = f^{-1}(x) < x$,
which implies that $f g^{-1} (x) < g^{-1}(x)$, contrary to our assumption on $f$.

\vsp\vsp

Let $I$ be the smallest closed interval containing the origin whose endpoints are fixed by
$\Gamma_0$, and let $H$ be its stabilizer in $\Gamma$. Since $\Gamma_0$ is normal
in $\Gamma$, for every $g \in \Gamma,$ either $g(I)$ equals $I$ or it is disjoint from it. By
Proposition \ref{prop convex subgroup dynamics}, this implies that $H$ is a convex subgroup.

\vspace{0.25cm}

\noindent{\underline{Claim (ii).}} The subgroup $H$ has smaller rank than $\widetilde{\Gamma}$.

\vsp

Indeed, on the one hand, $H \cap \widetilde{\Gamma}$ cannot be equal
to $\widetilde{\Gamma}$, since the latter does not have global fixed
points. On the other hand, since $\Gamma_0$ is contained in $H \cap
\widetilde{\Gamma}$, Claim (i) above implies that $H \cap
\widetilde{\Gamma}$ is a normal subgroup of $\widetilde{\Gamma}$ and
that $\widetilde{\Gamma} / (H \cap \widetilde{\Gamma})$ is Abelian. Therefore,
as the quotient $\widetilde{\Gamma} / (H \cap \widetilde{\Gamma})$ is
left-orderable, it has rank $>0$, thus showing the claim.

\vsp\vsp

It follows by induction that the space of left-orders of $H$ is either finite or a Cantor
set. Hence, by Proposition \ref{prop convex and LO}, if $\preceq$ is isolated in
$\mathcal{LO}(\Gamma)$, then $H$ is a Tararin group. However, if $H$ is a
Tararin group, then every left-order on $H$ is Conradian (see 
Lemma \ref{evidentemente}). By the convexity of $H \cap \widetilde{\Gamma}$
in $\widetilde{\Gamma}$ and the fact that $\widetilde{\Gamma} / (H\cap\widetilde{\Gamma})$
is Abelian, we have that the restriction of $\preceq$ to $\widetilde{\Gamma}$ is Conradian
(see Exercise \ref{ejercicito}). Therefore, by
Theorem \ref{listailor}, we have that $\preceq$ is a Conradian
order of $\Gamma$. As a consequence, using Proposition \ref{simple-alcanza}, we conclude that,
if $\preceq$ is isolated in $\mathcal{LO}(\Gamma)$, then $\Gamma$ must be a Tararin group, as desired.

\vsp\vsp\vsp

\noindent{\underline{Case II.}} Either $rank(R/\Gamma_0)\geq 2$, or $rank(R/\Gamma_0)=1$
and there exists $g \in \widetilde{\Gamma}$ that does not act on $R/\Gamma_0$ by
multiplication by $\pm 1$. (In particular, $R/\Gamma_0$ is not isomorphic to $\Z$.)

\vsp

\index{Measure!Radon measure}
\index{Translation number}
In this case, Proposition \ref{teo type I preserves measure} and
Remark \ref{teo type I preserves measure extended} provide 
an $R$-invariant Radon\footnote{Recall that a {\bf{\em
Radon measure}} is a measure giving finite mass to compact sets.}
measure $\nu$, to which is associated a {\bf{\em translation
number homomorphism}} $\tau_\nu:R \to (\mathbb{R},+)$ defined
by $\tau_\nu(g) := \nu([x,g(x)[)$ (we use the convention 
$\nu([x,y]) := -\nu([y,x])$ for $y<x$ throughout). 
Note that this definition does not depend on
the choice of the point $x$.

\index{Conrad homomorphism}
\begin{small}\begin{ejer} \label{lo-mismo}
Show that $\tau_{\nu}$ coincides (up to a positive multiple) with the
Conrad homomorphism on $R$ associated to the convex jump with respect
to the maximal proper convex subgroup (see \S \ref{conrad-general})
\end{ejer}

\index{Translation number}
\begin{ejer} \label{ejer traslaciones}
Let $G$ be a subgroup of $\mathrm{Homeo}_+ (\mathbb{R}) $ with no global fixed point. Suppose that 
its action preserves a Radon measure for which the translation number homomorphism has an image 
that is not isomorphic to $\mathbb{Z}$. Prove that $G$-action is semiconjugate to an action by translations.
\end{ejer}

\begin{ejer} \label{ejer-se-usa}
Let $G$ be a subgroup of $\mathrm{Homeo}_+ (\mathbb{R})$ preserving a Radon measure $\nu$.

\noindent (i) Show that the kernel of $\tau_\nu$ coincides with the subset $G_0$ consisting of the elements
having fixed points. Moreover, show that for all $x$ in the support $supp(\nu)$, its stabilizer in $G$
coincides with $G_0$.

\noindent (ii) Conclude that if every element in $G$ has fixed points, then there is a global fixed point for the action.
\end{ejer}
\end{small}

The next proposition (which is interesting in its own right) tell
us that, up to multiplication by a positive constant, $\nu$ is
the unique $R$-invariant Radon measure. This is somewhat a
dynamical counterpart of Exercise \ref{lo-mismo} above.

\vsp
\begin{prop} \label{prop unique measure} {\em  Let $G$ be a subgroup of
$\mathrm{Homeo}_+ (\mathbb{R})$ preserving a Radon measure $\nu$.
Then, for any other (nontrivial) $G$-invariant Radon measure $\nu'$, there
is a positive real number $\kappa$ such that $\kappa \tau_{\nu} = \tau_{\nu'}$.
Moreover, if $\tau_{\nu} (G)$ is dense in $(\mathbb{R},+)$, then $\kappa\nu = \nu'$.}
\end{prop}

\noindent{\bf Proof.} It easily follows from Exercise \ref{ejer-se-usa} that $\tau_{\nu} (G)$ and
$\tau_{\nu'} (G)$ are simultaneously either discrete or dense in $\mathbb R$. In the former case,
the claim of the proposition is obvious. Below we deal with the latter case.

Fix $g\notin G_0$ and a point $x$ that is fixed by $G_0$. Then, as a combination
of Exercises \ref{subaditivo} and \ ref{ejer traslaciones}, we have that for all $f\in G$,
$$ \tau_\nu (f)=\tau_\nu(g)  \lim_{p\to \infty}\Big\{\frac{q}{p} \! : g^q(x)\leq f^p(x)<g^{q+1}(x)\Big\},$$
and the same holds changing $\nu$ by $\nu'$. Therefore, we have
$\tau_{\nu'} (g) \tau_\nu(f)= \tau_{\nu'} (f)\tau_\nu(g)$ for every $f \in G$, hence
$\tau_{\nu'}$ equals $\kappa \tau_{\nu}$ for a certain positive $\kappa$.

Next, we claim that the supports of $\nu$ and $\nu'$ coincide. Indeed, the density of
$\tau_{\nu'}(G)$ implies that $\nu'$ has no atoms and the action of $G$ on $supp(\nu')$
is minimal ({\em i.e.}, every orbit is dense). It follows that if there is a point $x\in supp(\nu')\setminus supp(\nu)$,
then there exists $g\in G$ such that $g(x)>x$ and $\nu([x,g(x)))=0$, contradicting the fact that
$\ker (\tau_{\nu'}) = \ker (\tau_\nu)$. Therefore, $supp(\nu') \subset supp(\nu)$,
and the reverse inclusion is proven analogously.

Finally, let $x<y$ be two points in the (common) supports of $\nu$ and $\nu'$,
and let $g_n\in G$ be such that $g_n(x)$ converges to $y$. Then
$$\nu' \big( [x,y] \big)
=  \lim_{n\to\infty}\nu' \big(  [x,g_n(x) ] \big)
= \lim_{n\to\infty}\tau_{\nu'} (g_n)
= \lim_{n\to\infty} \kappa \tau_\nu(g_n)= \kappa \nu \big( [ x,y] \big),$$
which finishes the proof. $\hfill\square$

\vspace{0.3cm}

Using the $R$-invariant Radon measure, we can describe the action of $\Gamma$
up to semiconjugacy. 

\vspace{0.25cm}

\noindent{\underline{Claim (i).}} There is a homomorphism $\Phi \!: \Gamma \to
\mathrm{Aff}_+(\mathbb{R})$ such that $\Phi(R)$ contains nontrivial translations,
and $\Gamma_0$ coincides with $\ker (\Phi) \cap R$. Moreover, the dynamical
realization of $(\Gamma,\preceq)$ is (continuously) semiconjugate to this affine action.

\vsp

\index{Translation number}
Indeed, let us continue denoting by $\nu$ an $R$-invariant Radon measure.
Since $R / \Gamma_0$ is not isomorphic to $\Z$, Proposition \ref{prop unique measure}
(and its proof) implies that $\nu$ is unique up to a scalar multiple, and that $\Gamma_0$
is the kernel of the translation number homomorphism $\tau_\nu$. As $R$ is normal in
$\Gamma$, this implies that for each $g\in \Gamma$, the measure $g_*(\nu)$
is also $R$-invariant. Thus, for every $g\in
\Gamma$, there is $\lambda_g>0$ such that $g_*(\nu)=\lambda_g\nu$.
This yields a group homomorphism $\lambda \!: \Gamma\to
\mathbb{R}^*$ into the group of positive reals (with multiplication). We then define
$\Phi \!: \Gamma\to \mathrm{Aff}_+(\mathbb R)$ by
$$\Phi (g) (x) := \frac{1}{\lambda_g} \;
x+\nu([0,g(0)]).$$
One can easily check that this is a homomorphism that extends $\tau_\nu$.

To show that the dynamical realization of $\preceq$ is semiconjugate to this
affine action, for each $x\in\mathbb{R}$ we let  $\varphi(x) := \nu([0,x])$.
Then $\varphi$ is a continuous, non-decreasing surjective map, and a direct
computation shows that for all $g \in \Gamma$ and every $x \in \mathbb{R}$,
$$
\varphi(g(x))=\Phi(g) (\varphi(x)),
$$
which shows the announced semiconjugacy.

\vspace{0.25cm}

Next, we let $I_\nu := (a,b)$, where $a := \sup\{x<0 \! : x\in
supp(\nu)\}$ and $b := \inf \{x>0 \! : x\in supp(\nu)\}$. We also
let $\Gamma_\nu$ be the stabilizer in $\Gamma$ of $I_\nu$.
The subgroup $\Gamma_{\nu}$ is easily seen to be convex.
Moreover, $\Gamma_\nu \cap R = \Gamma_0$.

Note that the rank of $\Gamma_{\nu}$ is smaller than that of $\Gamma$.
Thus, by the induction hypothesis, if $\Gamma_\nu$ admits infinitely many left-orders,
then no left-order on it is isolated. By convex extension, we conclude that $\preceq$
is non-isolated in $\mathcal{LO}(\Gamma)$. For the other case, the next claim applies.

\vspace{0.285cm}

\noindent{\underline{Claim (ii).}} If $\Gamma_\nu$ is a
Tararin group, then $\ker (\Phi)$ is convex.

\vspace{0.15cm}

Indeed, since $\Phi(\Gamma_\nu)$ does not contain any nontrivial translation,
it can only contain homotheties centered at $0$; in particular, it is Abelian. If
it is trivial, then $\ker (\Phi) = \Gamma_{\nu}$, so it is convex, as desired.
Assume that $\Phi(\Gamma_{\nu})$ is nontrivial, and let
$\{id\} = \Gamma^n \lhd \Gamma^{n-1} \lhd \ldots \lhd \Gamma^0 = \Gamma_\nu$
be the series of all convex subgroups of the Tararin group $\Gamma_\nu$. (Recall that
$\Gamma^{i-1} / \Gamma^{i}$ has rank $1$ and that the action of $\Gamma^{i-1}$
on $\Gamma^{i-1} / \Gamma^{i}$ is by multiplication by some negative number.)
By Exercise \ref{tararin-romain}, $\Gamma_\nu$ has a unique torsion-free Abelian
quotient, namely $\Gamma_\nu / \Gamma^{1}$. As this must  coincide with
$\Phi (\Gamma)$, we conclude that $\ker (\Phi)$ equals
$\Gamma^{1}$, hence it is convex.

\vspace{0.3cm}

Knowing that $\ker (\Phi)$ is convex, we can proceed to show that $\preceq$ is non-isolated. Indeed,
either $\Gamma/\ker (\Phi)$ is Abelian of rank at least $2$, or it is a non-Abelian subgroup of the
affine group. In the former case, it has no isolated left-orders (see \S \ref{ejemplificando-1}), hence
--by convex extension-- the left-order $\preceq$ is non-isolated in $\mathcal{LO}(\Gamma)$. In the
latter case, we are under the hypothesis of Proposition \ref{prop affine}, which yields the same
conclusion. This finishes the proof of Theorem \ref{teo finite rank solvable}.

\vspace{0.3cm}

\index{Group!Baumslag-Solitar}
\noindent{\bf Left-orders on the Baumslag-Solitar groups.}
Perhaps the simplest examples of finite-rank
solvable groups that are non virtually-nilpotent are the Baumslag-Solitar groups
$BS(1,\ell) := \langle h,g \! : hgh^{-1} = g^{\ell} \rangle$, where $\ell > 1$.  We have seen
that $BS(1,\ell)$ admits only four Conradian orders (see Example \ref{solo-cuatro}),
yet it also admits the left-orders induced from its affine faithful actions on the line
(see \S \ref{ejemplificando-2}). Below we follow the lines of the previous
proof to show that, actually, these are the only possible left-orders on $B(1,\ell)$.

As in \S \ref{ejemplificando-2}, we can see $B(1,\ell)$ as a semidirect product $\Z[\frac{1}{\ell}]\rtimes\Z$,
where the $\Z$-factor acts on $\Z[\frac{1}{\ell}] := \big\{ \frac{k}{\ell^m} \!: k,\,m \mbox { in }\Z \big\}$ by
multiplication by $\ell$. Viewed way, it easily follows that the nilpotent radical of $B(1,\ell)$ is $R:=\Z[\frac{1}{\ell}]$.

Now, given a left-order on $BS(1,\ell)$, we consider its dynamical realization.
Since $R = \Z[\frac{1}{\ell}]$ has rank 1, two cases may arise.

\vspace{0.3cm}

\noindent{\underline{Case I.}} There is a global fixed point for $R$.

\vsp

As $R$ is normal in $BS(1,\ell)$, the set of $R$-fixed points is $BS(1,\ell)$-invariant; thus, it is unbounded in both 
directions. In particular, $R$ must coincide with the convex subgroup $H$ that arises as the stabilizer of the interval that 
contains the origin and is enclosed by two consecutive $R$-fixed points. As a consequence, $\preceq$ is Conradian. 

\vspace{0.3cm}

\noindent{\underline{Case II.}} Every nontrivial element of $R$ is fixed-point free.

In this case, the action of $R$ is continuously semiconjugate to that of a dense group of translations, thus it preserves
a Radon measure $\nu$ without atoms that is unique up to a scalar factor. Moreover, $h$ does not preserve
$\nu$, otherwise we would have $\tau_\nu(g^{\ell-1})=\tau_\nu(g^{-1}hgh^{-1})=0$, contradicting the fact
that $g^{\ell-1}$ acts freely. Thus, the dynamical realization of $\preceq$ is semiconjugate to (the action given
by) a faithful embedding of $BS(1,\ell)$ into $\mathrm{Aff}_+(\mathbb R)$, and the left-order $\preceq$
coincides with a left-order induced from this affine action. 
\vsp

Note that, in Case I above, the fact that $\preceq$ is non-isolated follows from Proposition \ref{simple-alcanza},  since 
$BS(1,\ell)$ is not a Tararin group. More concretely, one readily sees that there is a sequence of affine-like orders coming 
from Case II that approximate $\preceq$; namely, it suffices to choose the first comparison point tending to either $-\infty$ 
or $\infty$. Note also that the fact that  the orders arising in Case II are non-isolated follows from Corollary \ref{cor:affine}.

It is worth pointing out that the description above --as well as its proof-- applies not only to dynamical
realizations of left-orders, but also to general (faithful) actions on the line with no global fixed point.
Such an action is hence either without crossings (with $R$ being the subgroup of elements having
fixed points) or semiconjugate to an affine action. (See \cite{rivas} for more details on this.)

\vsp\vsp\vsp

\index{Group!$Sol$}
\noindent{\bf Left-orders on $Sol$ groups.} Relevant examples of finite-rank solvable
groups that are non virtually-polycyclic are those of the form $Sol := \Z^2 \rtimes_A \Z$, where
$A$ is an hyperbolic automorphism of $\Z^2$ ({\em i.e.}, it is given by a matrix in $\mathrm{SL} (2,\Z)$
with trace larger than $2$, so that it has two irrational eigenvalues). Below we follow the lines of
the previous proof in this particular case to get an accurate description of the space of left-orders
and its subspaces of bi-invariant and Conradian orders. Actually, the methods employed yield
a complete description of all faithful actions on the line with no global fixed point.

We denote by $R$ the commutator subgroup of $Sol$ --which coincides with the $\Z^2$-factor--,
and we denote by $f$ the element of $\Z$ acting on $R$ as $A$. The subgroup $R$ is easily 
seen to coincide with the nilpotent radical of $Sol$.

Given a left-order $\preceq$ on $Sol$, let us consider its dynamical realization.
Since $A^T$ is $\mathbb Q$-irreducible and $R$ is Abelian and finitely-generated,
the next three properties are equivalent:

\vsp

\noindent -- There is an element in $R$ having a fixed point;

\vsp

\noindent -- Every element of $R$ has a fixed point;

\vsp

\noindent -- There is a global fixed point for $R$.

\vsp

\noindent Indeed, on the one hand, having a fixed point for $g \!\in\! R$ is equivalent to $\tau_{\nu} (g) \!=\! 0$
for an $R$-invariant Radon measure $\nu$ (see Exercise \ref{ejer-se-usa}). On the other hand, $A$ also acts at
the level of translation numbers, as is next shown.
\index{Translation number}

\begin{small}\begin{ejer}
Let $g_1,g_2$ be the canonical basis of $R = \mathbb{Z}^2$. Show that
$$\left( \begin{array}
{c}
\tau_{\nu} (f g_1 f^{-1})  \\
\tau_{\nu} (f g_2 f^{-1})  \\
\end{array} \right)
= A^T
\left( \begin{array}
{c}
\tau_{\nu} ( g_1 )  \\
\tau_{\nu} ( g_2 )  \\
\end{array} \right).$$
\end{ejer}\end{small}

Thus, the two cases considered in the proof of Theorem \ref{teo finite rank solvable}
fit with those considered below.

\vspace{0.25cm}

\noindent{\underline{Case I.}} The subgroup $R$ has a global fixed point.

\vsp

Since $Sol$ acts without global fixed points and $R$ is normal in $Sol$, in this case
the set of $R$-fixed points is unbounded in both directions (and $\Gamma$-invariant).
As for $BS(1,\ell)$, this implies  that $\preceq$ is Conradian. To see that $\preceq$ is
non-isolated, one may argue by convex extension by noticing that $R$ is convex and
rank-two Abelian. Alternatively, $Sol$ is not a Tararin group.

\vspace{0.25cm}

\noindent\underline{Case II.} There is no global fixed point for $R$.

\vsp

In this case, $R$ is semiconjugate to a dense group of translations, thus it preserves
a Radon measure without atoms $\nu$ that is unique up to a scalar factor. As $f$
is hyperbolic, it cannot preserve $\nu$: it acts as an homothethy with ratio one of the
eigenvalues of $A^T$. Thus, the dynamical realization of $\preceq$ is semiconjugate 
to (the action given by) a faithful embedding of $\Gamma$ into $\mathrm{Aff}_+(\mathbb R)$.
The fact that $\preceq$ is non-isolated in this case follows from Corollary \ref{cor:affine}.

\vspace{0.3cm}

\begin{figure}[h!]
\begin{center}
\includegraphics[scale=0.7]{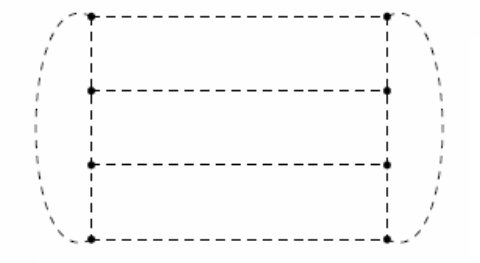}
\end{center}
\vspace{-0.865cm}
\begin{center}
Figure 15: Depicting the space of left-orders of $Sol$.
\end{center}
\end{figure}

\vspace{-0.2cm}

It follows from the previous analysis that, as it was the case for the
Baumslag-Solitar groups, there are two types of left-orders on $Sol$:

\vspace{0.2cm}

\noindent{\underline{Case I.}} Conradian orders.

\vsp

These correspond to those left-orders for which the normal subgroup
$R = \Z^2$ is convex. Thus, $\mathcal{CO}(Sol)$ is made up 
of two copies of the Cantor set $\mathcal{LO}(\Z^2)$, each of which
corresponds to the choice of a sign for $f$. (In Figure 15, these are
represented by the two ``vertical dashed circles''.) Observe that
among all bi-orders on $R$, those that are invariant under
conjugacy by $f$ are those that that correspond (under Conrad's
homomorphisms) to eigendirections of the matrix $A^T$. Since the
corresponding left-orders on $Sol$ are the bi-invariant ones, we
conclude that $Sol$ supports exactly eight bi-orders.

\begin{small}\begin{rem}
A classification of finitely-generated solvable groups admitting only
finitely many bi-orders can be found in \cite{botto,kopytov}.
\end{rem}\end{small}

\vspace{0.2cm}

\noindent{\underline{Case II.}} Left-orders coming from affine actions.

\vsp

These form an open set (which is locally a Cantor set) that complements the
subspace of Conradian orders. They can be described as in \S \ref{ejemplificando-2}.
Note, however, that $Sol$ admits four embeddings into $\mathrm{Aff}_+ (\mathbb{R})$. 
According to this, the affine-like orders are depicted as four ``horizontal dotted lines'' in Figure 15. 
These ``lines'' accumulate at the eight bi-invariant orders in $\mathcal{CO}(Sol)$. This 
is similar to the previously described approximation of the four bi-orders on 
Baumslag-Solitar's groups by affine-like orders.

\vsp

Based on all of this, it is not difficult to describe the dynamics of the conjugacy action
of $Sol$ on its space of left-orders. We leave this as an exercise to the reader.


\subsection{The space of left-orders of (general) solvable groups}
\label{section solvable (general)}

\hspace{0.45cm} In this section, we prove the general result announced in the previous one. 

\vsp

\begin{thm} \label{teo solvable}
{\em The space of left-orders of a countable virtually-solvable group is either
finite or homeomorphic to the Cantor set.}
\end{thm}
\index{Measure!Radon measure}

In the preceding section, we  treated the case of virtually finite-rank solvable groups. To do that, we 
established that these groups admit quasi-invariant measures when acting on the line. For concreteness, recall that 
a Radon measure $\nu$ on the line is {\bf{\em quasi-invariant}} under the action of a group $\Gamma$ if for every 
$g\in \Gamma$, there exists  a positive real number $\lambda_g$ such that $g_*(\nu)=\lambda_g\,\nu$, where 
by definition $g_*(\nu)(X) := \nu(g^{-1}(X))$ for every measurable set $X$. 

Let us briefly recall the argument. Virtually finite-rank solvable groups are virtually nilpotent-by-Abelian, with 
finite-rank nilpotent part. The key point is that finite-rank nilpotent groups preserve a Radon measure when acting by 
orientation-preserving homeomorphisms of the line. Moreover, this measure is unique up to a scalar factor. Since the 
nilpotent part is normal, this yields the announced quasi-invariance. As a consequence, the group action is necessarily 
semiconjugate to an affine action.

It turns out that this nice picture does not longer hold for actions of general solvable groups. Next, we reproduce a classical 
example due to Plante \cite{plante} of an action of $\Z\wr\Z$ for which there is no quasi-invariant measure. Alternatively, 
the reader may check  Example \ref{ex corona laminar action} where we describe a left-order on $\Z\wr\Z$  whose 
dynamical realization is semiconjugate to Plante's action of $\Z\wr\Z$.

\index{Measure!quasi-invariant}

\begin{small}
\begin{ex} \label{ejemplo de plante}
The wreath product $\Z\wr\Z := \bigoplus_{\Z} \Z\rtimes \Z$ is a metabelian
group having $H := \bigoplus_{\Z}\Z$ as its maximal nilpotent subgroup. We next describe an action
of $\Z\wr \Z$ on the real line with the property that for every shift-invariant subgroup of $H$, no global
fixed point arises, although every element therein admits fixed points. This implies in particular
that there is no quasi-invariant measure for $\Z\wr\Z$. Indeed, such a measure would be invariant by the
commutator subgroup $[\Z\wr\Z\, ,\,\Z\wr\Z]$. However, since this subgroup is shift-invariant, 
this is in contradiction with Exercise \ref{ejer-se-usa}.

For the construction, let $f$ denote the homothethy $x\mapsto 2x$. Let $I_0:=[-1,1]$, and
for $i \!\in\! \Z$, denote $I_i:=f^i(I_0)$. Let $h \!: I_0 \to I_0$ be a homeomorphism such that 
$h(-1/2) = 1/2$ and $h(x)>x$ for all $x \!\in\! (-1,1)$.
We define $h_i \!: I_i \to I_i$ by $h_i := f^i h f^{-i}$. Note that this is equivalent to saying that 
$f^{-1} h_i (x) = h_{i-1} f^{-1} (x)$ holds for all $x \in I_i$, that is, $h_i f (y) = f h_{i-1} (y)$
for all $y \in I_{i-1}$.
Below, we extend the definition of each $h_i$ to
the whole line in such a way that $f$ and $h_0$ generate a group isomorphic to $\Z\wr\Z$.

One easily convinces oneself that there is a unique way to extend the maps $h_i$ to commuting homeomorphisms of the
real line. For instance, to ensure commutativity, we must necessarily have $h_{i-1}(x) := h_i^m h_{i-1}h_i^{-m}(x)$
for $x\in h_i^m(I_{i-1})$. The (proof of the uniqueness of the) extension can then be easily achieved
by induction. We continue denoting by
$h_i$ the resulting homeomorphisms. We claim that $fh_if^{-1} = h_{i+1}$ holds. Indeed, this follows from
the definition for $x\in I_{i+1}$. Assume inductively that $fh_if^{-1}(x)=h_{i+1}(x)$ holds for all $x\in I_k$
for a certain $k \geq i+1$, and let $x \in I_{k+1}$. Letting $m \in \mathbb{Z}$ be such that $x = h_{k+1}^m (y)$
for a certain $y \in I_k$, we have
\begin{eqnarray*}
fh_if^{-1}(x)
&=& f h_i f^{-1} (h_{k+1}^m (y))
\hspace{0.18cm} = \hspace{0.18cm} f h_i h^m_{k} f^{-1} (y) \\
&=& f h_k^m h_i f^{-1} (y)
\hspace{0.18cm} = \hspace{0.18cm} h_{k+1}^m f h_i f^{-1} (y)
\hspace{0.18cm} = \hspace{0.18cm} h^m_{k+1} h_{i+1} (y)
\hspace{0.18cm} = \hspace{0.18cm} h_{i+1}(x),
\end{eqnarray*}
where the second and fourth equalities follow from the definition of $h_k$, the third from
the commutativity between $h_i$ and $h_k$, and the fifth from the induction hypothesis.
\end{ex}
\end{small}

\vspace{0.2cm}
\begin{figure}[h!]
\begin{center}
\includegraphics[width=0.45\textwidth]{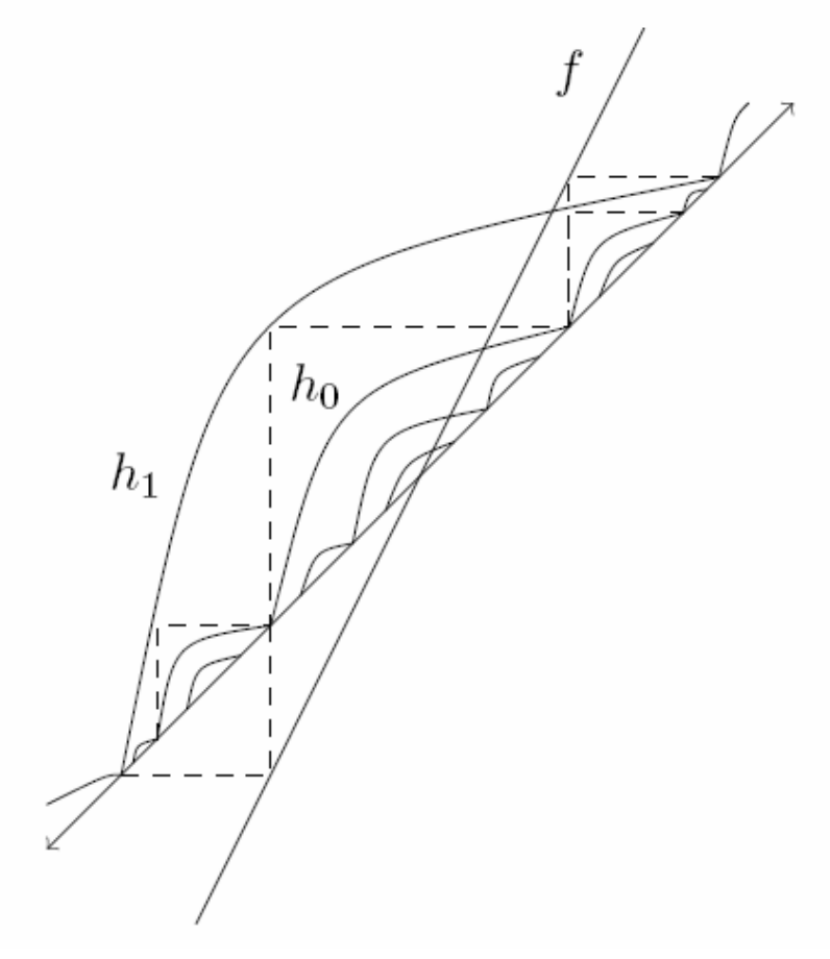}
\end{center}
\vspace{-0.8cm}
\begin{center}
Figure 16: Plante's action of $\Z\wr\Z$.
\end{center}
\end{figure}

\vspace{-0,5cm}
\index{Plante's action}

To deal with the phenomenon illustrated by the preceding example, we will use the machinery developed in 
\S\ref{general-Conrad}. The price to pay is that, unlike \S \ref{section finite-rank-solvable},  we will not give 
a classification --up to semiconjugacy-- of all actions on the real line for general countable solvable groups. Rather, 
we will give a rough local description of the dynamics on the real line that still allows us to conclude Theorem \ref{teo solvable}. 

\begin{small}
\begin{rem} \label{rem laminar solvable} A full classification of solvable group actions on the real line up to 
semiconjugacy was recently obtained in \cite{BMRT-solvable}. Therein, it is shown that such an action is either semiconjugate to an 
affine action (in which case there is a quasi-invariant Radon measure), or it is a {\bf {\em laminar action}}. Roughly, the latter means that 
there is an action on an oriented (real) tree $\mathcal T$ such that the induced action on the (visual) boundary $\partial\mathcal T$ is 
semiconjugate to the original action on the line. The simplest example of a laminar action is Plante's action of $\Z\wr\Z$ above, 
where the  associated tree $\mathcal T$ is a simplicial one (of infinite valence). More precisely, there is a vertex in the tree for 
each interval in the open support of a canonical generator $h_i$ of $\bigoplus_\Z\Z$. Moreover, if $u$ (resp. $v$) is the vertex 
associated to an interval $I$ (resp. $J$) in the open support of $h_i$ (resp. $h_j$), then $u$ is connected to $v$ if and only if 
$|i-j|\leq 1$ and $I$ contains $J$ or vice versa.

The notion of laminarity for actions on the real line was introduced and extensively developed in \cite{BMRT-LM}, and is also 
useful to understand  the space of actions of other groups such as Thompson group $\mathrm{F}$.
\end{rem}\end{small}

We start with an exercise that follows from an easy reformulation of part of the proof of Propositions
\ref{first-tararin} and \ref{simple-alcanza}.

\vsp

\begin{small}\begin{ejer}
Show that every left-order on a group admitting infinitely many convex subgroups is non-isolated in
the corresponding space of left-orders.
\end{ejer}\end{small}

\vsp

Due to the preceding exercise, in order to prove Theorem \ref{teo solvable}, it suffices to consider left-orders
with finitely many convex subgroups. Let  $\preceq$ be such an order on a group $\Gamma$, with
$$\{ id\} = C_n \subsetneq C_{n-1}\subset \ldots \subsetneq C_{0}=\Gamma$$
being the family of convex subgroups. One of these subgroups must coincide with the Conradian soul
$G := C_{\preceq} (\Gamma)$, that is, with the maximal $\preceq$-convex subgroup restricted to which $\preceq$
is Conradian (see \S \ref{C-soul}). If $G$ is not a Tararin group, then by Proposition \ref{simple-alcanza},
the restriction of $\preceq$ to $G$ is non-isolated in $\mathcal{LO}(G)$; by convex extension, $\preceq$ is
not isolated in $\mathcal{LO}(\Gamma)$. Hence, in all what follows, we assume that $G$ is a Tararin group.

If $G = \Gamma$, then we are done: $\Gamma$ admits only finitely many left-orders. If $G$ is trivial,
then $\preceq$ is non-isolated in $\mathcal{LO}(\Gamma)$, due to Theorem \ref{primero}. We hence
suppose that $G$ is a nontrivial, proper subgroup of $\Gamma$, say $G = C_{\ell}$, with $n > \ell > 0$.
We will show that the restriction of $\preceq$ to $C_{\ell-1}$ is non-isolated; by convex extension, this
in turns will imply that $\preceq$ is non-isolated in $\mathcal{LO}(\Gamma)$, as desired. As the claim to
be shown only involves $C_{\ell - 1}$, to simplify will denote this group as $\Gamma$; equivalently, we
will assume that $\ell = 1$, that is, there is no convex subgroup strictly between $G = C_1$ and $\Gamma$.

We consider the dynamical realization of $\preceq$. Since
$G$ is a proper convex subgroup, it admits at least one fixed point on each side of the
origin. We let $I_G$ be the smallest open interval fixed by $G$ that contains the origin. By
Proposition \ref{prop convex subgroup dynamics}, the convexity of $G$ immediately implies 
the following lemma.

\begin{lem}\label{lem:T}
{\em Every element of $\Gamma$ either fixes $I_G$ or moves it to a disjoint interval.
In particular, the stabilizer of $I_G$ coincides with $G$.}
\end{lem}

\vsp

From now on, we assume $\Gamma $ to be virtually solvable. Let $\widetilde{\Gamma}$ be a finite-index,
normal, solvable subgroup of $\Gamma$. We let
$\widetilde{\Gamma}^0 := \widetilde{\Gamma}$ and
$\widetilde{\Gamma}^j := [\widetilde{\Gamma}^{j-1},\widetilde{\Gamma}^{j-1}]$
be the associated derived series:
$$\{ id\}=\widetilde{\Gamma}^k\lhd \widetilde{\Gamma}^{k-1}\lhd \ldots
\lhd \widetilde{\Gamma}^1\lhd \widetilde{\Gamma}^0 = \widetilde \Gamma \lhd \Gamma.$$
Note that each $\widetilde{\Gamma}^j$ is normal in $\Gamma$. We let $i$ be the minimal
index such that $\widetilde{\Gamma}^{i}$ is contained in $G$. Since $G$ is a nontrivial,
proper, convex subgroup, we have that $k > i \geq 1$; see Figure 17 below.

\vspace{0.25cm}


\beginpicture

\setcoordinatesystem units <1.2cm,1.2cm>


\put{$\Gamma$} at 3 2
\put{$\widetilde{\Gamma}^{i-1}$} at 4,5 1,2
\put{$G$} at 1,75 0,7
\put{$G \cap \widetilde{\Gamma}^{i-1}$} at 3,85 -0,5
\put{$\widetilde{\Gamma}^{i}$} at 3,4 -1,2
\put{$\{id\} = \widetilde{\Gamma}^k$} at 3 -2,05

\put{$\bullet$} at 3 1,7
\put{$\bullet$} at 2 0.7
\put{$\bullet$} at 4 1,2
\put{$\bullet$} at 3 -0,5
\put{$\bullet$} at 3 -1,2
\put{$\bullet$} at 3 -1,7

\plot
3 1,7
3 1,5 /

\plot
3 1,5
2 0,7
3 -0.5 /

\plot
3 1,5
4 1,2
3 -0.5
3 -1,7 /

\put{Figure 17: The groups $G$, $\widetilde{\Gamma}^i$ and $\widetilde{\Gamma}^{i-1}$.} at 3 -2.7

\put{} at -3 0


\endpicture


\vspace{0.5cm}

The subgroup $\widetilde{\Gamma}^{i-1}$ will be crucial for our analysis: although it is not always nilpotent, the restriction
of $\preceq$ to $\widetilde{\Gamma}^{i-1}$ will be shown to be Conradian. Thus, dynamically, it will play the role played by the nilpotent
radical in the finite-rank case.  We like to think of it as a kind of {\em Conradian skeleton} of $(\Gamma, \preceq)$.

\vsp

\begin{lem} \label{lem:the restriction is Conrad}
{\em The order $\preceq$ restricted to $\widetilde{\Gamma}^{i-1}$ is Conradian.}
\end{lem}

\noindent{\bf Proof.} By definition, the subgroup $\widetilde{\Gamma}^{i}=[\widetilde{\Gamma}^{i-1},\widetilde{\Gamma}^{i-1}]$
is contained in $G$. Therefore, $\widetilde{\Gamma}^{i-1}\cap G$ is normal in $\widetilde{\Gamma}^{i-1}$, as well as convex therein.
Moreover, as the quotient \hspace{0.05cm}
$\widetilde{\Gamma}^{i-1}/\widetilde{\Gamma}^{i-1}\cap G$ \hspace{0.05cm}
is Abelian, it only admits Conradian orders. Since $\preceq$
restricted to $\widetilde{\Gamma}^{i-1}\cap G$ is Conradian, this implies that $\preceq$ restricted to  $\widetilde{\Gamma}^{i-1}$ is a
convex extension of a Conradian order by a Conradian one, hence Conradian (see Exercice \ref{ejercicito}).$\hfill\square$

\vsp

\begin{lem}\label{lem:the action is cofinal}{\em The action of $\widetilde{\Gamma}^{i-1}$ has no global fixed point.}
\end{lem}

\noindent{\bf Proof.} Let $I$ be the smallest open interval containing the origin that is fixed by $\widetilde{\Gamma}^{i-1}$.
Since $\widetilde{\Gamma}^{i-1}$ is normal in $\Gamma$, the interval $I$ is either fixed or moved to a disjoint interval 
by each $g \in \Gamma$. In particular, the stabilizer $Stab_\Gamma(I)$ of $I$ is a convex subgroup of $\Gamma$
(see Proposition \ref{prop convex subgroup dynamics}). Now, if $I$ were not the whole line, then the maximality of 
$G$ (as a proper, convex subgroup) would imply that  $Stab_\Gamma(I) \subseteq G$, thus yielding
$\widetilde{\Gamma}^{i-1}\subset G$, which is a contradiction. $\hfill\square$

\vspace{0.3cm}

Since $\preceq$ restricted to $\widetilde{\Gamma}^{i-1}$ is Conradian, its action on the real line has no crossings 
(see Exercise \ref{ejer:crossings}). It follows that the set of elements in $\widetilde{\Gamma}^{i-1}$ having fixed
points is a normal subgroup of $\widetilde{\Gamma}^{i-1}$ (actually, of $\Gamma$); see Proposition~\ref{normality}. 
In particular, if $g\in\widetilde{\Gamma}^{i-1}$ does not act freely, then the set of fixed points of $g$ accumulates at 
both $-\infty$ and $+\infty$. Thus, in order to prove Theorem \ref{teo solvable}, we need to analyze two cases.

\vsp\vsp

\noindent {\bf \underbar{Case I.}} The subgroup $\widetilde{\Gamma}^{i-1}$ contains elements without fixed points.
(This case arises for instance if $\widetilde{\Gamma}^{i-1}$ has finite rank.)

\vsp\vsp

We first observe that, in this case, $\widetilde{\Gamma}^{i-1}$ preserves a nontrivial Radon measure $\nu$. Indeed, since the order on
$\widetilde{\Gamma}^{i-1}$ is Conradian, its action on the real line has no crossings. Further, since there is $g_0 \! \in \! \widetilde{\Gamma}^{i-1}$
acting freely, Proposition \ref{normality} easily implies that $\widetilde{\Gamma}^{i-1}$ has a maximal proper convex subgroup, namely
$\{g\in \widetilde{\Gamma}^{i-1} \!: g \text{ has a fixed point}\}$. It then follows from Proposition \ref{teo type I preserves measure} and
Remark \ref{teo type I preserves measure extended} that $\widetilde{\Gamma}^{i-1}$ preserves a nontrivial Radon measure on the line.

Now, from the normality of $\widetilde{\Gamma}^{i-1}$ in $\Gamma$ and Proposition \ref{prop unique measure},
we have that there is a homomorphism $\lambda \!: g\mapsto \lambda_g$ from $\Gamma$ into $\mathbb{R}^*$ 
satisfying $\tau_{g_*\nu}=\lambda_g\tau_{\nu}$. The next lemma comes from the work of Plante \cite{plante}.

\vsp

\index{Measure!Radon measure}
\begin{lem} {\em If the homomorphism $\lambda$ is trivial, then $\Gamma$ preserves a Radon measure on the real line.
Otherwise, $\Gamma$ admits a quasi-invariant Radon measure which is $\widetilde{\Gamma}^{i-1}$-invariant.}
\end{lem}

\noindent{\bf Proof.} Recall that by Proposition \ref{prop unique measure}, if $\tau_\nu(\widetilde{\Gamma}^{i-1})$ is dense, then
$\nu$ is quasi-invariant. This occurs for instance if $\lambda$ is nontrivial. Indeed, choosing $g\in \Gamma$ such that
$\lambda_g < 1$, we have for every $f\in \widetilde{\Gamma}^{i-1}$,
\begin{equation}\tau_\nu \big( g^{-1}fg \big) = \nu \big( [g^{-1}(x), g^{-1}f(x)) \big) =
g_* \nu ([x,f(x))) = \tau_{g_*\nu}(f)=\lambda_g\tau_\nu(f). \label{eq lambda y mu}\end{equation}

We thus assume that $\tau_\nu(\widetilde{\Gamma}^{i-1})\simeq\Z$ is not dense. In particular,
we suppose that $\lambda$ is trivial and that $\tau_\nu(\widetilde{\Gamma}^{i-1})\simeq\Z$. 
We let $H := ker(\tau_\nu)=\{g\in \widetilde{\Gamma}^{i-1}\mid \tau_\nu(g)=0\}$. We have seen that
$H$ consists of the elements in $\widetilde{\Gamma}^{i-1}$ having fixed points (see Exercise
\ref{ejer-se-usa}). Therefore, $H$ is normal not only in $\widetilde{\Gamma}^{i-1}$, but also in 
$\Gamma$. Moreover, the condition $\tau_\nu(\widetilde{\Gamma}^{i-1})\simeq\Z$ translates to 
$\widetilde{\Gamma}^{i-1} / H\simeq \Z$. We claim that $\widetilde{\Gamma}^{i-1} / H$ is
in the center of $\Gamma / H$. Indeed, letting $f \!\in\! \widetilde{\Gamma}^{i-1}$ be a generator
of $\widetilde{\Gamma}^{i-1} / H$, for each $g\in \Gamma$ we have that $g^{-1}fg = f^n h$
holds for certain $h \in H$ and $n \in \mathbb{Z}$. We need to show that $n=1$. Now, by
(\ref{eq lambda y mu}), we have
$$n\, \tau_\nu(f)=\tau_\nu(g^{-1}fg)=\lambda_g\tau_\nu(f) = \tau_\nu (f) \neq 0,$$
which implies that $n=1$, as desired.

Finally, the quotient group $\Gamma/\widetilde{\Gamma}^{i-1}\simeq (\Gamma/H)/(\widetilde{\Gamma}^{i-1}/H)$
acts on the compact quotient $Fix(H)/\!\sim$, where $x \sim f(x)$ for each $f \in \widetilde{\Gamma}^{i-1}$ and all $x \in Fix(H)$. 
This space is easily seen to be homeomorphic to the circle. Therefore, since $\Gamma$ is solvable (hence amenable; see 
\S \ref{super-witte}) it preserves a probability measure on it. Pulling back this measure to the real line, we obtain a 
$\Gamma$-invariant Radon measure on $\mathbb{R}$. $\hfill\square$

\vspace{0.3cm}

We now claim that the dynamical realization of $\preceq$ is semiconjugate to a non-Abelian affine action. As in the preceding 
section, this will follow once we show that the homomorphism $\lambda$ is nontrivial. Assume for a contradiction that $\lambda$ is
trivial. Then by the preceding lemma, there is a $\Gamma$-invariant Radon measure $\nu$. Moreover, as the origin is moved by
every nontrivial element and $\Gamma$ contains elements having fixed points (for instance, those of $G$), the origin does not
belong to the support of $\nu$. Let $I_\nu$ be the connected component of the complement of the support of $\nu$ containing
the origin. The interval $I_\nu$ is either fixed or moved to a disjoint interval by each element of $\Gamma$, hence its stabilizer
$Stab_\Gamma(I_\nu)$ is a convex subgroup of $\Gamma$. Since this subgroup contains $G$ and since $G$ is the
maximal proper convex subgroup of $\Gamma$, we must have $Stab_\Gamma(I_\nu) = G$. Further, $Stab_\Gamma(I_\nu)$
coincides with the kernel of the translation number 
\index{Translation number} homomorphism $\tau_\nu \! : \Gamma \to (\mathbb{R},+)$, thus it is normal
in $\Gamma$. We thus conclude that $G$ is normal and co-Abelian in $\Gamma$. Therefore, $\preceq$ is a convex extension
of a Conradian order by a Conradian one, hence it is Conradian (see Exercise \ref{ejercicito}). However, this contradicts
the fact that $G$ is the Conradian soul of $\Gamma$.

We can finally show that $\preceq$ is non-isolated by invoking Corollary \ref{cor:affine}. Indeed, it easily follows from
the construction of the dynamical realization that the kernel of the induced homomorphism from $\Gamma$ into
$\mathrm{Aff}_+(\mathbb{R})$ is a $\preceq$-convex subgroup, hence the hypotheses of the corollary are fulfilled.

\vspace{0.35cm}

\noindent {\bf \underbar{Case II.}} Every element of $\widetilde{\Gamma}^{i-1}$ admits fixed points.

\vsp\vsp

In this case, we will prove that the approximation scheme by conjugates developed in
\S \ref{section-soul} applies. More precisely, starting from the dynamical realization of
$\preceq$, we will induce a new left-order using a comparison point that is outside but
very close to $I_G$. The main issue here is to ensure that this procedure can be performed in
such a way that the order restricted to $G$ remains untouched (compare Theorem \ref{final}).
In the proof, it will become clear that the action is somewhat similar to the one described
in Example \ref{ejemplo de plante}.

For each nontrivial element $g \in \widetilde{\Gamma}^{i-1}$, let us denote by $I_g$
the connected component of the complement of its set of fixed points that contains the origin. 
Similarly, denote by $I_G$ the connected component of the complement of the set of fixed points 
of $G$ that contains the origin. It follows from Lemma \ref{lem:the action is cofinal} that the union 
of all the $I_g$'s is the whole real line.

\vsp

\begin{lem} \label{lem no strong crossings} 
{\em For each $f\in \Gamma$ and $g\in \widetilde{\Gamma}^{i-1}$, one of the following possibilities occurs:

\vsp

\noindent -- \esp $f(I_g)=I_g$;

\vsp

\noindent -- \esp $f(I_g)$ is disjoint from $I_g$;

\vsp

\noindent --  \esp $\overline{f(I_g)}\subset I_g$ holds up to changing $f$ by its inverse if necessary.}
\end{lem}

\noindent{\bf Proof.} By Lemma \ref{lem:the restriction is Conrad}, the order $\preceq$ restricted to
$\widetilde{\Gamma}^{i-1}$ is Conradian, hence $\widetilde{\Gamma}^{i-1}$ acts without crossings (see 
Exercise \ref{ejer:crossings}). As $\widetilde{\Gamma}^{i-1}$ is normal in $\Gamma$, the lemma easily follows.
$\hfill\square$

\vspace{0.3cm}

The next two lemmas are similar to the preceding one. The first follows from the convexity of $G$,
and the second from the nonexistence  of crossings for the action of $\widetilde{\Gamma}^{i-1}$
and the fact that $\widetilde{\Gamma}^{i-1}$ is normal in $\Gamma$. Details are left to the reader.

\vsp

\begin{lem}\label{lem:T either dilates or fixes} {\em In the preceding lemma, if $g$ does not belong to $G$, then
the second possibility cannot occur for $f\in G$. In other words, for all $g \in \widetilde{\Gamma}^{i-1}\setminus G$  and each
$f \!\in\! G$, either $f$ fixes $I_g$, or (up to changing $f$ by $f^{-1}$ if necessary), it holds that $\overline{I_g}\subset f(I_g)$.}
\end{lem}

\begin{lem}\label{cor: union/intersection} {\em Let $I$ be the intersection of all the
intervals $I_g$, where $g$ ranges over $\widetilde{\Gamma}^{i-1} \setminus G$. Then each
element $f \in \Gamma$ either moves $I$ to a disjoint interval, or up to
replacing it by its inverse, it holds that $I \subset f(I)$.}
\end{lem}

\vsp

The next lemma is a kind of refined version of Theorem \ref{S=CS} knowing
that $G$ has finite rank and/or admits only finitely many left-orders, and that
$\widetilde{\Gamma}^{i-1}$ is normal and $\preceq$-Conradian.

\vsp

\begin{lem}\label{lem:intersection}
{\em The intersection of  all the intervals $I_g$ for $g\in \widetilde{\Gamma}^{i-1}
\setminus G$ coincides with $I_G$.}
\end{lem}

\noindent{\bf Proof.} Since $\widetilde{\Gamma}^{i-1}$ is a $\preceq$-Conradian subgroup
(see Lemma \ref{lem:the restriction is Conrad}), its action has no crossings, which implies
that the family of intervals $I_g$, with $g\in \widetilde{\Gamma}^{i-1} \setminus \{ id \}$, is totally
ordered by inclusion (see Exercise \ref{ejer:irreducible}). Moreover, as $G$ is convex, for
each $g \in \widetilde{\Gamma}^{i-1} \setminus G$ we have that $I_g$ strictly contains $I_G$.
Therefore, letting $I$ be the intersection of all the $I_g$'s
for $g\in \widetilde{\Gamma}^{i-1} \setminus G$, we have
that $I$ is a bounded interval containing $I_G$.

Assume that no element $f \in \Gamma$ is such that $I$ strictly contains $f(I)$. Then, by Lemma 
\ref{cor: union/intersection}, every $f \in \Gamma$ either fixes $I$ or moves it to a disjoint interval. 
It readily follows that the stabilizer of $I$ in $\Gamma$ is a proper convex subgroup of $\Gamma$; 
moreover, this subgroup contains $G$. As $G$ is the maximal proper convex subgroup, this necessarily
implies that $I_G$ equals $I$.

Therefore, it suffices to prove that no $f \! \in \! \Gamma$ satisfies $f(I) \subsetneq I$. Assume otherwise
for a certain $f$. As $I = \bigcap_{g \in \widetilde{\Gamma}^{i-1} \setminus G} I_g$, there must exist
$g \in \widetilde{\Gamma}^{i-1}\setminus G$ such that $I \cap f (I_{g})$ is strictly contained in $I$.
Since $fgf^{-1}$ belongs to $\widetilde{\Gamma}^{i-1}$, by the definition of $I$, we must have that
$f g f^{-1}$ belongs to $G$. Actually, the same holds for $f^n g f^{-n}$, for all $n \geq 1$.

Note that $f^n(I_g)$ is an open interval (not necessarily containing the origin) that is fixed 
by $f^n g f^{-n}$, with no fixed point inside. Moreover, by Lemma \ref{lem no strong crossings},
one has $\overline{f (I_g)} \subset I_g $, which implies that 
$\overline{f^{n+1}(I_g)} \subset f^n (I_g)$ holds for all $n \geq 1$. Together with the fact that the action of $G$
has no crossings, this easily implies that $f^k g f^{-k}$ belongs to $Stab_G (f^n(I_g))$, for all $k \geq n$.
As a consequence, $(Stab_G (f^n(I_g)))$ is a strictly decreasing sequence of convex subgroups
of $G$ for any left-order induced from a sequence starting with a point in the (nonempty)
intersection of the compact intervals $\overline{f^n(I_g)}$. However, this contradicts the
fact that $G$ is a Tararin group. 
$\hfill\square$

\vspace{0.35cm}

The next lemma follows almost directly from Proposition \ref{a-probar} and its proof. Actually,
it is a kind of restatement of it for dynamical realizations. We leave the details to the reader.

\vsp

\begin{lem}\label{lem:approx}
{\em Let $(x_n)$ be a sequence of points outside $I_G$ that converges to one of the endpoints of $I_G$. 
For each $n \geq 1$, let $\preceq_n$ be any left-order on $\Gamma$ obtained in a dynamical-lexicographic 
way from a sequence starting with $x_n$. Then $\preceq_n$ converges to $\preceq$, and differs from 
$\preceq$ for $n$ sufficiently large.}
\end{lem}

\vsp

\begin{small}\begin{ejer} In the context of the preceding lemma, show that $\preceq_n$ differs from $\preceq$ 
for {\em every} $n$. To do this, show that for every $g,h$ in $\tilde{\Gamma}^{i-1} \setminus G$ such that 
$I_g \subset I_h$, there exists $f \in \Gamma$ such that $I_h \subset f (I_g)$. (See \cite[Lemma 5.10]{RT} 
in case of problems with this.)
\end{ejer}\end{small}

\vsp

For the sake of concreteness, the orders $\preceq_n$ in Lemma \ref{lem:approx} may --and will-- 
be taken as those for which the second comparison point is the origin (so that no other
comparison point is necessary). The convergence in the statement means that
for any sequence $(x_n)$ converging to an endpoint of $I_G$ from outside, given
$g \in \Gamma\setminus G$, we have that $g \succ id$ holds if and only if $g \succ_n id$
holds for all sufficiently large $n$.  
Therefore, to prove that $\preceq$ is non-isolated, we
are left to showing that for a well-chosen sequence $(x_n)$ as above, the left-orders
$\preceq_n$ coincide with $\preceq$ when restricted to $G$ for sufficiently large $n$.
This is achieved by the next lemma, which closes the proof of Theorem \ref{teo solvable}.

\vsp

\begin{lem} \label{lem approx over T}{\em There exists a sequence of points
$x_n$ converging to an endpoint of $I_G$ from outside such that the induced
left-orders $\preceq_{n}$ coincide with $\preceq$ on $G$ for all $n$.}
\end{lem}

\vsp

To show this, we need one more general lemma for Tararin groups.

\begin{lem}\label{lem:TararinOrder}
{\em Let $T$ be a Tararin group, with chain of convex subgroups
$$\{id\} = T^k \lhd T^{k-1} \lhd \ldots \lhd T^{1} \lhd T^0 = T,$$
and let $f_T$ be an element in $T \setminus T^1$ acting on $T^{1} /
T^{2}$ as the multiplication by a negative (rational) number. Suppose 
that $T$ acts by orientation-preserving homeomorphisms of the line
in such a way that the sets of fixed points of nontrivial elements
have empty interior (as is the case for dynamical realizations). Let 
$y\in\mathbb{R}$ be a point that is not fixed by $f_T$.  Then, for
every left-order $\preceq$ on $T$,  there exists a point $x$ between 
$f_T^{-2} (y)$ and $f_T^2 (y)$ such that $\preceq$ coincides 
on $T^{1}$ with the restriction of any left-order $\preceq'$ on 
$T$ induced from a sequence starting with $x$.}
\end{lem}

\noindent{\bf Proof.} As the sets of fixed points of nontrivial
elements have empty interior, and since $T$ is countable, there is
a point $z$ between $f_T^{-1}(y)$ and $y$ whose orbit under $T$
is free. Any such point induces --in a  dynamical-lexicographic way--
a left-order $\preceq^*$ on $T$. Since $T^{1}$ is necessarily convex for this
order, there is an open interval $I$ containing $z$ that is fixed by $T^{1}$
and contains no fixed point of $T^1$ inside. Moreover, $I$ is
mapped to a disjoint interval by any nontrivial power of $f_T$; in
particular, $I$ contains at most one point of the orbit of $y$ under
$\langle f_T \rangle$. As a consequence, $I$ strictly lies between
$f_T^{-2}(y)$ and $f_T (y)$. Now, by Proposition \ref{prop two
orbits of Tararin}, there exists an element $g$ in either $T^{1}$ or
$f_T T^1$ such that $\preceq$ and $\preceq^*_g$ coincide on $T^{1}$.
As $\preceq^*_g$ is the dynamical-lexicographic order with comparison
point $x := g (z)$, this point satisfies the conclusion of the lemma. $\hfill\square$

\vspace{0.35cm}

\noindent{\bf Proof of Lemma \ref{lem approx over T}.}
Due to Lemma \ref{lem:intersection} (and since $\widetilde{\Gamma}^{i-1}$ acts without crossings),
there exists a sequence of elements $g_n \in \widetilde{\Gamma}^{i-1}\setminus G$ such that
$I_{g_n}$ converges to $I_G$. As in the preceding lemma, let $f_G$ be an element in
$G \setminus G^1$ acting on $G^1 / G^2$ as the multiplication by a negative number.
According to Lemma \ref{lem:T either dilates or fixes}, we may pass to a subsequence
for which one of the two possibilities below occur.

\vsp\vsp

\noindent{\underline{Subcase (i).}} Each interval $I_{g_n}$ is fixed by $f_G$.

\vsp\vsp

We first claim that $G$ must fix each interval $I_{g_n}$. Indeed, letting $y$ be an endpoint
of any of the $I_{g_n}$'s, we may induce a left-order on $G$ from  a sequence having $y$
as its initial point. For such an order, the stabilizer of $y$ in $G$ is a convex subgroup of $G$
containing $f_G$. It must hence coincide with $G$, and therefore $G$ fixes $y$, as desired.

We may now let $x_n$ be any of the endpoints of $I_{g_n}$, say the right one. As $G$ fixes $x_n$
and the second comparison point of $\preceq_n$ is the origin, the restriction of $\preceq_n$ to
$G$ coincides with that of $\preceq$. Moreover, since $I_{g_n}$ converges to $I_G$, the
points $x_n$ converge (from outside) to the right endpoint of $I_G$.

\vsp\vsp

\noindent{\underline{Subcase (ii).}} For each $n$, we have either 
$ \overline{f_G (I_{g_n})} \subset I_{g_n}$ or $\overline{I_{g_n}} \subset f_G (I_{g_n})$. 

\vsp\vsp

Up to taking a subsequence and changing $f_G$ by its inverse if necessary, we may assume that 
$\overline{I_{g_n}} \subset I_{g_{n-1}}$ and $\overline{I_{g_n}} \subset f_G(I_{g_n})$, for all $n \geq 1$. 
In the case where $f_G \prec id$
(resp. $f_G \succ id$), let $y_n$ be the left (resp. right) endpoint of $I_{g_n}$. Note that $y_n$ converges to a
fixed point of $f_G$ (hence of $G$). Moreover, for all $n \geq 1$, if $f_G \prec id$ (resp. $f_G \succ id$), then
\begin{equation} \label{eq aproximando f_G}
f_G (y_n) < y_n \;\;  \text{(resp. } f_G(y_n) > y_n).
\end{equation}

Now, for each $y_n$, let us consider the point $x_n$ provided by Lemma \ref{lem:TararinOrder}. Then
$\preceq_{n}$ coincides with $\preceq$ in restriction to the maximal proper convex subgroup $G^1$ of $G$.
Moreover, as $x_n$ lies between $f_G^{-2}(y_n)$ and $f_G^{2}(y_n)$, it converges to the same point as $y_n$; 
in particular, it converges to an endpoint of $I_G$ from outside. Furthermore, (\ref{eq aproximando f_G})
obviously holds for $x_n$ instead of $y_n$. This implies that $\preceq_n$ coincides with
$\preceq$ over the whole of $G$, as desired. $\hfill\square$

\vspace{0.35cm}

To close this section, let us mention that it is unclear what is the most general framework for which
Theorem \ref{teo solvable} still holds. This naturally yields to the next

\vsp

\index{Group!amenable}
\begin{question} Is the space of left-orders of a countable amenable group either finite or a Cantor set~?
What about groups without free subgroups in two generators~?
\end{question}

\vsp

It is very likely that the previous methods can be extended to a wide family of amenable groups, namely
that of {\bf{\em elementary amenable}} ones. Roughly, this is the smallest family of groups that contains
all Abelian groups and that is stable under taking extensions, direct limits, quotients and subgroups
(see \cite{chou,mexico} for more on this). A relevant example is considered in described below. 

\vsp

\begin{small}\begin{ex}
The group $\Gamma := \mathbb{Z} \ltimes (... \mathbb{Z} \wr (\mathbb{Z} \wr ... )...)))$ in which the  conjugacy
action of the left factor consists in shifting (the level of) the factors in the right wreath product is obviously
elementary amenable. It is somewhat a variation of Plante's example (yet it is not solvable), and was
simultaneously introduced in \cite{brin} and \cite{mexico}. A natural (faithful) action of this group on the
line comes from identifying the generator of the left factor with the map $x \mapsto 2x$ and a generator
of the $0^{th}$-factor of the wreath product with a homeomorphism with support contained in $[-2,2]$ 
and sending $-1$ to $1$. It is very likely that its space of left-orders is a Cantor set, and it actually 
seems reasonable to try to describe all of its actions on the line.
\end{ex}\end{small}


\section{Verbal Properties of Left-Orders}
\label{verbal-section}

\hspace{0.45cm}
Let $\mathcal{W}$ be the set of reduced words in two letters $a,b$.
(This naturally identifies with the free group in two generators.) We distinguish three subsets of
$\mathcal{W}$, namely $\mathcal{W}^+$, $\mathcal{W}^-$, and $\mathcal{W}^\pm$, the
set of words involving only positive, negative, and mixed exponents in $a$ and $b$, respectively.
Given elements $f,g$ in a group $\Gamma$ and $W \in \mathcal{W}$, we let $W(f,g)$ be the
element in $\Gamma$ obtained from the expression of $W$ by replacing $a$ and $b$ with 
$f$ by $g$, respectively.

For $W \in \mathcal{W}$, a left-order $\preceq$ on a group $\Gamma$ will be said to satisfy {\bf {\em verbal property $W$}},
or that it is a {\bf {\em $W$-order}}, if whenever $f$ and $g$ are $\preceq$-positive, the element $W(f,g)$ is also $\preceq$-positive. 
Note that this defines a nontrivial property only in the case where $W\in \mathcal{W}^\pm$,
hence in the sequel we will only consider these words.

\begin{small}\begin{ex}
For $W(a,b) := b^{-1}ab$, one easily checks that the set of $W$-left-orders coincides with
that of bi-orders.
\end{ex}

\begin{ex}
For $W(a,b) := b^{-1}ab^2$, Proposition \ref{lesla} tells us that the set of $W$-orders
corresponds to that of Conradian ones.
\end{ex}\end{small}
\index{Order!W-order}

The next two questions become natural in this context.

\begin{question} Does there exist a word $W$ such that the $W$-orders are those that satisfy
an specific and relevant algebraic property different from bi-orderability or the Conradian one~?
\end{question}

\index{Crossing!double}
\begin{question} Is the property of not having a double crossing (see Example \ref{double-crossing})
for a left-order equivalent to a verbal property (or to an interesection of finitely many ones)~?
\end{question}

As it is easy to check, the subset of $W$-orders is closed inside $\mathcal{LO}(\Gamma)$, and the
conjugacy action preserves this subset. The next result on free groups is only stated for
two generators, though it can be easily extended to more generators.

\vsp

\begin{thm} {\em The free group on two generators admits left-orders satisfying no
verbal property $W \in \mathcal{W}^{\pm}$. Actually, this is the case of a $G_{\delta}$-dense
subset of $\mathcal{LO} (\mathbb{F}_2)$.}
\label{no-verbal}
\end{thm}

\vsp

Let us first show that the existence of a single left-order satisfying no verbal
property implies that this is the case for most left-orders. Thies relies on Lemma
\ref{classical}, as shown by the next lemma.

\vsp

\begin{lem} {\em Every left-order on $\mathbb{F}_2$ having a dense orbit under the
conjugacy action satisfies no verbal property $W \in \mathcal{W}^{\pm}$.}
\end{lem}

\noindent{\bf Proof.} Otherwise, as the closure of such an orbit only contains
$W$-orders, we would be in contradiction with Theorem \ref{no-verbal}.
$\hfill\square$

\vspace{0.2cm}

\begin{question}
It is a nontrivial fact that the real-analytic homeomorphisms of the line
given by $x \mapsto x+1$ and $x \mapsto x^3$ generate a free group \cite{glass-free}. By analyticity,
a $G_{\delta}$-dense subset $S$ of points in the line have a free orbit under this action. Given a point
$x \in S$, we may associate to it the left-order on $F_2$ defined by $f \succ g$ whenever $f(x) > g(x)$.
Is the set of $x \!\in\! S$ for which the associated order satisfies no verbal property still a
$G_{\delta}$-dense subset of $\mathbb{R}$~?
\end{question}

\vsp\vsp

We next proceed to the proof of the first claim of Theorem \ref{no-verbal}, which is done
via a very simple dynamical argument. Namely, given $W \!\in\! \mathcal{W}^{\pm}$,
we will construct two increasing homeomorphisms of the real line $f, g$, both moving
the origin to the right, such that in the action  of $\mathbb{F}_2$ given by $a\to f$,
$b \to g$, the homeomorphism $W(f,g)$ moves the origin to the left.
Then, any dynamical-lexicographic left-order $\preceq$ associated to a sequence starting at
the origin will be such that $f \succ id$, $g \succ id$, and $W(f,g) \prec id$. This is enough for
our purposes except for that the action we will produce will be not necessary faithful. However,
this is just a minor detail that may be solved in many ways. For instance, one can make the action
faithful by perturbing it close to infinity, as in \S \ref{case-free}; alternatively, one may consider
a convex extension of the order $\preceq$, as in \S \ref{section-convex-extension}.

The construction of the desired action is done as follows. By interchanging $a$ and $b$ if
necessary, we may assume that the word $W=W(a,b)$ writes in the form $W = W_1 a^{-n} W_2$,
where $W_2$ is either empty or a product of positive powers of $a$ and $b$, the integer $n$ is positive,
and $W_1$ is arbitrary.
Let us consider two local homeomorphisms defined on a right neighborhood of the real line
such that $f(0) > 0$, $g(0) > 0$ and $W_2(f,g)(0) < f^n (0)$. This can be easily done by taking
$f(0) \gg g(0)$ and letting $g$ be almost flat on a very large right-neighborhood of the origin.
If $W_1$ is empty, just extend $f$ and $g$ to homeomorphisms of the real line. Otherwise,
write $W_1 = a^{n_k} b^{m_k} \ldots a^{n_2} b^{m_2} a^{n_1} b^{m_1}$, where all $m_i, n_i$ are
nonzero excepting perhaps $n_k$. The extension of $f$ and $g$ to a left-neighborhood of the
origin depends on the signs of the exponents $m_i, n_i$, and is done in a constructive manner.
Namely, first extend $f$ slightly so that $f^{-n} W_2 (f,g)(0)$ is defined and $f$ has a fixed point
$x_1$ to the left of the origin.  Then extend $g$ to a left-neighborhood of the origin so that
$g^{m_1} f^{-n} W_2(f,g)(0) < x_1$ and $g$ has a fixed point $y_1$ to the left of $x_1$. Note
that $m_1 > 0$ forces $g$ to be topologically attracting on the right towards $y_1$ on an interval containing
$f^{-n} W_2(f,g)(0)$, whereas $m_1 < 0$ forces right topological repulsion. Next, extend $f$ to a left
neighborhood of $x_1$ so that $f^{n_1} g^{m_1} f^{-n} W_2(f,g)(0) < y_1$ and $f$ has a fixed point
$x_2$ to the left of $y_1$. Again, if $n_1 > 0$, this forces topological attraction on the right towards
$x_2$, whereas $n_1 < 0$ implies topological repulsion on the right.


\vspace{0.35cm}
\beginpicture
\label{picture verbal orderings}

\setcoordinatesystem units <0.9cm,0.9cm>

\plot 0 -2 0 2 / \plot  -4.3 0 2 0 /

\put{$f$} at 2.3 2

\plot
0 1.55
0.1 1.64
0.2 1.69
0.5 1.8
1 1.89
1.5 1.96
2 2 /

\plot
0 1.55
-0.1 1.23
-0.2 0.35
-0.3 0.03
-0.4 -0.28
-0.6 -0.64
-0.65 -0.65 /

\plot
-0.65 -0.65
-0.7 -0.69
-0.8 -0.78
-0.9 -0.8
-2.2 -0.82
-2.3 -0.825
-2.4 -0.83
-2.5 -0.85 /

\plot
-2.5 -0.85
-2.6 -1
-2.7 -2.7 /

\plot
-2.7 -2.7
-2.71 -3
-2.72 -3.2
-2.75 -3.6
-2.77 -3.78
-2.78 -3.793
-2.79 -3.796
-2.8 -3.9 /

\plot
-2.8 -3.9
-3 -3.96
-3.4 -3.97
-3.8 -3.98
-4.3 -4 /

\put{$g$} at 2.4 0.4

\plot -0.05 0.2
0 0.3
2 0.4 /

\plot -0.05 0.2
-0.2 -1.3
-0.25 -1.45
-0.3 -1.5 /

\plot -0.3 -1.5
-0.35 -1.55
-0.4 -1.6
-1.7 -1.7
-3 -1.8
-3.2 -1.86
-3.3 -1.9
-3.4 -2 /

\plot
-3.4 -2
-3.45 -3
-3.47 -3.4
-3.48 -3.48
-3.49 -3.49
-3.5 -3.5 /

\plot
-3.5 -3.5
-3.8 -3.58
-4 -3.592
-4.1 -3.595
-4.2 -3.598
-4.3 -3.6 /

\put{{\tiny $\bullet$}} at  -0.65 -0.65
\put{{\tiny $\bullet$}} at  -2.7 -2.7
\put{{\tiny $\bullet$}} at  -3.99 -3.99
\put{{\tiny $\bullet$}} at   -1.7 -1.7
\put{{\tiny $\bullet$}} at  -3.5 -3.5


\setdots  \plot -4.3 -4.3 2 2 /

\small

\put{Figure 18 \!: The case $W_1=a^{n_2}b^{m_2}a^{n_1}b^{m_1}$, where $m_1 >0, n_1 < 0, m_2 <0$, and $n_2 >0$.} at -1 -5

\put{ } at -8.7 2.2
\endpicture


\vspace{0.65cm}

Continuing the procedure in this manner (see Figure 18 for an
illustration), we get partially defined  homeomorphisms $f,g$ for which
$$0 > f^{n_k} g^{m_k} \ldots f^{n_2} g^{m_2} f^{n_1} g^{m_1} f^{-n} W_2(f,g)(0) = W(f,g)(0).$$
Extending $f,g$ arbitrarily to homeomorphisms of the real line, we finally obtain the desired action.


\section{A Non Left-Orderable Group, and More}
\label{structure-section}

\subsection{No left-order on finite-index subgroups of $\mathrm{SL}(n,\mathbb{Z})$}

\hspace{0.45cm} Proposition \ref{listailor} gave us a simple
criterium for non left-orderability of certain groups. In the same
spirit, an important result due to Witte Morris \cite{Wi} establishes
that finite-index subgroups of $\mathrm{SL}(n,\mathbb{Z})$ are non left-orderable
for $n \geq 3$. (Note that most of these groups are torsion-free, because of the
classical Selberg lemma \cite{sel}.)

\vspace{0.15cm}

\index{Witte Morris!non-orderability of $SL(n,\Z)$}
\begin{thm} {\em If $\Gamma$ is a finite-index subgroup of $\mathrm{SL}(n,\mathbb{Z})$,
with $n \geq 3$, then $\Gamma$ is non left-orderable.}
\label{witte-SL}
\end{thm}

\noindent{\bf Proof.} Since $\mathrm{SL}(3,\mathbb{Z})$ injects
into $\mathrm{SL}(n,\mathbb{Z})$ for every $n \geq 3$, it suffices
to consider the case $n = 3$. Assume for a contradiction that
$\preceq$ is a left-order on a finite-index subgroup $\Gamma$ of
$\mathrm{SL}(n,\mathbb{Z})$. Note that for large enough
$k \!\in \!\mathbb{N}$, the following elements must
belong to $\Gamma$:
$$\begin{array}
{ccccccccccccccccccccccccccccccccc}
\hspace{1.2cm}
g_1  = \left(
\begin{array}
{ccc}
1 & k & 0 \\
0 & 1 & 0 \\
0 & 0 & 1 \\
\end{array}
\right) ,& &
g_2 = \left(
\begin{array}
{ccc}
1 & 0 & k \\
0 & 1 & 0 \\
0 & 0 & 1 \\
\end{array}
\right) ,& &
g_3  = \left(
\begin{array}
{ccc}
1 & 0 & 0 \\
0 & 1 & k \\
0 & 0 & 1 \\
\end{array}
\right), \\

& & & & & & & & & & & & & & & & & & & & & &   \\

\hspace{1.2cm}
g_4 =
\left(
\begin{array}
{ccc}
1 & 0 & 0 \\
k & 1 & 0 \\
0 & 0 & 1 \\
\end{array}
\right) ,& &
g_5 = \left(
\begin{array}
{ccc}
1 & 0 & 0 \\
0 & 1 & 0 \\
k & 0 & 1 \\
\end{array}
\right) ,& &
g_6 = \left(
\begin{array}
{ccc}
1 & 0 & 0 \\
0 & 1 & 0 \\
0 & k & 1 \\
\end{array}
\right) . \\
\end{array}$$
\noindent It is easy to check that for each $i \in \mathbb{Z}/ 6 \mathbb{Z}$,
the following relations hold:\\
$$g_i g_{i+1} = g_{i+1} g_i,
     \qquad [g_{i-1},g_{i+1}] = g_i^k.$$
In particular, the group generated by $g_{i-1},g_i$ and $g_{i+1}$ is nilpotent.

For $g \!\in\! \Gamma$, we define $|g|\!:=\!g$ if $g \succeq id$,
and $|g|\!:=\!g^{-1}$ in the other case. We also write
$g \gg h$ if $g  \succ  h^n$ for every $n \! \geq \! 1$.
We claim that either \esp $|g_{i-1}| \gg |g_i|$ \esp or \esp $|g_{i+1}| \gg |g_i|$.
\esp Indeed, as $\preceq$ restricted to the subgroup $\langle g_{i-1},g_i,g_{i+1}\rangle$
is Conradian (see Theorem \ref{nilpotent-cantor}) and a power of $g_i$ is a
commutator, this follows from Remark \ref{la-necesito}.

Assume for instance that $|g_1| \ll |g_2|$, the case where $|g_2| \ll |g_1|$
being analogous. Then we obtain
$|g_1| \ll |g_2| \ll |g_3| \ll |g_4| \ll |g_5| \ll |g_6| \ll |g_1|,$
which is absurd.  $\hfill\square$

\vspace{0.35cm}

It follows from an important theorem due to Margulis that for $n \geq 3$,
\index{Margulis!normal subgroup theorem}
every normal subgroup of a finite-index subgroup of $\mathrm{SL}(n,\mathbb{Z})$
either is finite or has finite index (see \cite{Ma2}). As a corollary, we obtain
the following strong version of Theorem \ref{witte-SL}.

\vspace{0.1cm}

\begin{thm} {\em For $n \geq 3$, no torsion-free, finite-index subgroup of
$\mathrm{SL}(n,\mathbb{Z})$ admits a total, nontrivial, left-invariant preorder.}
\label{wittefuerte}
\end{thm}\index{Order!preorder}

\noindent{\bf Proof.} If $\Gamma$ is such a group and admits a nontrivial, total preorder,
then by Exercise \ref{preorders-orders}, there is a nontrivial quotient $\Gamma/H$
that is left-orderable. Since $\Gamma$ is torsion-free, it has no nontrivial finite
subgroup. Therefore, there are only two possible cases: either $H$ is trivial, in
which case we contradict Theorem \ref{witte-SL}, or $\Gamma/H$ is finite and nontrivial,
which is impossible, as no nontrivial finite group admits a nontrivial, left-invariant
preorder. (Indeed, if $f \succ id$ for such a preorder, then $f^n \succ id$ for all
$n \in \mathbb{N}$.)
$\hfill\square$

\vsp\vsp\vsp

In terms of semigroups, this translates into the next result.

\vsp

\begin{cor} {\em If $n \!\geq\! 3$ and $\Gamma$ is a torsion-free, finite-index subgroup
of $\mathrm{SL}(n,\mathbb{Z})$, then there is only one subsemigroup $P$ of $\Gamma$
satisfying $P \cup P^{-1} = \Gamma$, namely $P = \Gamma$.}
\end{cor}

\vsp\vsp\vsp

The results above remained conjecturally true for all lattices in simple Lie groups of rank$\geq 2$ for 
many years, with an important contribution \cite{wit1} by Lifschitz and Witte Morris concerning the case of
higher $\mathbb{Q}$-rank, as well as other non-cocompact lattices. Nevertheless, in the recent breakthrough \cite{DH}, 
Hurtado and the first-named author of this book  managed to produce a complete proof of the  non-left-orderability of 
lattices in higher-rank simple Lie groups. Although this uses deep machinery coming from Lie group theory, the 
most important ingredients of proof from the viewpoint of orderable groups will be discussed in the next chapter. 


\subsection{A canonical decomposition of the space of left-orders}
\label{descomp}

\index{Order!type I}
\hspace{0.45cm} Let $(\Gamma,\preceq)$ be a finitely-generated, left-ordered
group, and let $\Gamma_0$ be its maximal $\preceq$-convex subgroup (see 
Example \ref{actuar-convexo}). The action of $\Gamma$ on $\Gamma / \Gamma_0$ 
may or may not be Conradian. In the first case, we will say that $\preceq$ is of {\bf{\em type I}}.
The next proposition generalizes Corollary \ref{cerrado}.

\vsp

\begin{prop} {\em The set of left-orders of type I is closed inside $\mathcal{LO}(\Gamma)$.}
\end{prop}

\noindent{\bf Proof.} Since $\Gamma$ is finitely-generated, $\Gamma / \Gamma'$ may be written as
$\mathbb{Z}^k \times G$, where $k \!\geq\! 1$ and $G$ is a finite Abelian group. Let $(\preceq_n)$
be a sequence of type-I left-orders on $\Gamma$ converging to a left-order $\preceq$. We
must show that $\preceq$ is also of type I. To do this, note that associated to each
$\preceq_n$, there is a Conrad's homomorphism $\tau_n$, which may be though of as defined
on $\mathbb{Z}^k$. This homomorphism may be chosen {\em normalized}. More precisely, 
if we let $\{g_1,\ldots,g_k\}$ be a family of elements whose representatives generate 
$\Gamma / \Gamma'$ and denote $a_{n,i} := \tau_n (g_i)$, 
then the vector $(a_{n,1},\ldots, a_{n,k})$ belongs to the
$(k-1)$-sphere $\mathrm{S}^{k-1}$, for each $n$.

\vsp\vsp\vsp

\noindent{\underline{Claim (i).}} The points $(a_{n,1},\ldots, a_{n,k})$
converge to some limit $(a_1, \ldots, a_k) \in \mathrm{S}^{k-1}$.

\vsp

Otherwise, there are subsequences $(\tau_{n_i})$ and $(\tau_{m_i})$ such that
the associated vectors converge to two different points of $\mathrm{S}^{k-1}$. Let 
$\mathbb{H}_1, \mathbb{H}_2$ be the ortogonal hyperplanes to these points. 
These hyperplanes divide $\mathbb{R}^k$ into four regions.
Let us pick a point of integer coordinates on each of these
regions, and let $h_1,h_2,h_3,h_4$ be elements of $\Gamma$ which project to
these points under the quotient map $\Gamma \rightarrow \mathbb{Z}^k \times G$. 
For a large enough index $i$, the values of both $\tau_{n_i}(h_j)$ and $\tau_{m_i} (h_j)$
are nonzero for each $j$, but the signs of these numbers must be different for some
$j$. As Conrad's homomorphisms are non-decreasing, after passing to a subsequence
of $(\preceq_{n_i})$ and $(\preceq_{m_i})$ this implies that, for some $j$, the
element $h_j$ will have different signs for $\preceq_{n_i}$ and $\preceq_{m_i}$.
However, this is in contradiction with the convergence of $\preceq_n$. 
This shows the announced convergence.

Note that the vector $(a_{1},\ldots, a_{k})$ gives raise to a 
group homomorphism $\tau \!: \mathbb{Z}^k \rightarrow \mathbb{R}$ 
(which may be though of as defined on $\Gamma$), namely, for 
$g \sim g_1^{n_1} \cdots g_k^{n_k}$ in $\Gamma/\Gamma'$, 
$$\tau (g) := \sum_{i=1}^k a_i n_i.$$

\vsp

\noindent{\underline{Claim (ii).}} The kernel of
$\tau$ is a $\preceq$-convex subgroup of $\Gamma$.

\vsp

Indeed, let $g \in \Gamma$ and $f \!\in\! \ker (\tau)$ 
be such that $id \preceq g \preceq f$. As Conrad's 
homomorphisms are order preserving, for each $n$, we have
$$0 = \tau_n (id) \leq \tau_n (g) \leq \tau_n (f).$$
As $\tau_n$ pointwise converges to $\tau$ and $\tau (f) = 0$, the inequalities
above yield, after passing to the limit, $\tau (g) = 0$. Thus, $g$ belongs to $\ker(\tau)$.

\vsp\vsp

As a consequence of Claim (ii), the maximal $\preceq$-convex subgroup $\Gamma_0$ contains
$\ker (\tau)$. Also, the action of $\Gamma$ on $\Gamma / \Gamma_0$ is order-isomorphic to
that on $\Gamma / \ker(\tau) \big/ \Gamma_0 / \ker (\tau)$. Since the latter is an action by
translations, the former is, in particular, Conradian. Therefore, $\preceq$ is of
type I. $\hfill\square$

\vsp\vsp\vsp\vsp\vsp

The case where the action of $\Gamma$ on $\Gamma / \Gamma_0$ is not Conradian 
is dynamically more interesting. We know by definition that there must exist a
crossing for the action. The question is ``how large'' can be the ``domain of crossing''.
To formalize this idea, for each \, $h \!\in\! \Gamma$, \, let us consider the ``interval''
$$I(h) := \big\{ \bar{h} \in \Gamma \!\!: \mbox{there exists a crossing } (f,g;u,w,v)
\mbox{ such that } fv \prec h \prec \bar{h} \prec gu \big\}\!.$$
By definition, $I(h)$ is a convex subset of $\Gamma$.

\vsp

\begin{lem} {\em If the set $I(h)$ is bounded from above for some
$h \in \Gamma$, then it is bounded from above for all $h \in \Gamma$.}
\end{lem}

\noindent{\bf Proof.} As the notion of crossing is invariant under conjugation, 
it holds that \, $h_1 \big( I(h_2) \big) = I(h_1h_2)$ \, for all $h_1,h_2$ in 
$\Gamma$. The lemma easily follows. $\hfill\square$

\vspace{0.45cm}

If (the action of $\Gamma$ on $\Gamma / \Gamma_0$ is has crossings
and) $I(h)$ is bounded from above for all $h \in \Gamma$, we will
say that $\preceq$ is of {\bf{\em type II}}. Otherwise, $\preceq$ will
be said of {\bf{\em type III}}. We then have a canonical
decomposition of the space of left-orders of $\Gamma$ into three
disjoint subsets (compare \cite[Theorem 7.E]{glass}):
$$\mathcal{LO}(\Gamma) =
\mathcal{LO}_I (\Gamma) \sqcup
\mathcal{LO}_{II} (\Gamma) \sqcup \mathcal{LO}_{III}(\Gamma).$$

\vspace{0.01cm}

\index{Order!type II}
\index{Order!type III}
\begin{small}
\begin{ex} Every Conradian order is of type I. Therefore, by Theorem \ref{st-th},
finitely-generated, locally-indicable groups admit left-orders of type I.
\end{ex}

\begin{ex} Smirnov's left-orders $\preceq_{\varepsilon}$ (with $\varepsilon$ irrational;
see \S \ref{ejemplificando-2})
on subgroups of the affine group are prototypes of type-III left-orders.
However, these groups being bi-orderable, they also admit left-orders of type I,
which actually arise as limits of Smirnov type orders. Moreover, the description given
in \S \ref{ejemplificando-2} shows that these groups do not admit type-II left-orders.
As a consequence, $\mathcal{LO}_{III}(\Gamma)$ is not necessarily closed
inside $\mathcal{LO}(\Gamma)$.

The last remark above can be made more precise.  
Namely, if we choose a sequence $(g_n)$ in (any non-Abelian subgroup of) 
the affine-group so that $g_n^{-1} (\varepsilon)$
tends to $+\infty$, then all conjugate left-orders $(\preceq_{\varepsilon})_{g_n}$ are of
type III, but the limit left-order $\preceq_{\infty}$ is bi-invariant, hence of type I.
Thus, a limit of type-III left-orders in the same orbit of the conjugacy action may
fail to be of type III.
\end{ex}

\begin{ejer} Let $(x_i)$ and $(y_n)$ be two sequences of points in $]0,1[$ so that
$x_i$ converges to the origin and $\{y_n\}$ is dense. For each $i \geq 1$, let
$(z_{n,i})_n$ be the sequence having $x_i$ as its first term and the $y_n$'s as the
next ones. Associated to this sequence there is a dynamical-lexicographic left-order
$\preceq_i$ on Thompson's group F, namely, \esp $f \succ_i id$ \esp if and only if the
smallest $n$ for which $f(z_{n,i}) \neq z_{n,i}$ is such that $f(z_{n,i}) > z_{n,i}$
(see \S \ref{general-3}).
Show that $\preceq_i$ is of type III for all $i$, but any adherence point of the
sequence $(\preceq_i)$ in $\mathcal{LO}(\mathrm{F})$ is of type I.
\end{ejer}

\begin{rem} The subset of left-orders of type II on the free group $F_2$ is dense in
$\mathcal{LO}(\mathbb{F}_2)$. Roughly, the proof proceeds as follows. (Compare \S
\ref{case-free}.) Start with an arbitrary left-order $\preceq$ on $\mathbb{F}_2$ together
with an integer $n \in \mathbb{N}$. Consider the dynamical realization of $\preceq$
as well as a very large compact subinterval $I$ in the real line on which the dynamics
captures all inequalities between elements in the ball $B_n (id)$ or radius $n$ in
$\mathbb{F}_2$. Then consider a new action of $\mathbb{F}_2$ which coincides
with this dynamical realization on $I$ and commutes with a translation of the line
(of very large amplitude). This new action induces a (perhaps partial) left-order
$\preceq_n$, which can be easily completed to a total one (by convex extension)
which is not of type I (by adding crossings). Clearly, the left-orders
$\preceq_n$ are all of type II and converge to $\preceq$.
\end{rem}

\begin{rem} A similar construction allows us to produce a sequence of type-III left-orders 
on $\mathbb{F}_2$ that converges to an order of type-II. Roughly, starting with a type II
left-order $\preceq$, we consider its dynamical realization. We keep it untouched on a 
very large compact interval $I$, and outside $I$ we perturb it by inserting infinitely many
crossings for the generators along larger and larger domains. The new action will then
induce a type-III  left-order on $\mathbb{F}_2$ that is very close to $\preceq$.
We leave the details to the reader.
\end{rem} \end{small}

\vsp







\index{Action!cofinal}
\noindent{\bf Cofinal elements and the type of left-orders.} Recall from \S \ref{conrad-general}
that an element $f$ of a left-ordered group $(\Gamma,\preceq)$ is {\em $\preceq$-cofinal} 
if for any $g \in \Gamma$ there exist integers $m,n$ such that $f^m \prec g \prec f^n$. In terms
of dynamical realizations (see \S \ref{general-3}), for countable groups, this corresponds
to that $f$ has no fixed point on the real line.

\index{Cofinal element}
Following \cite{cofinal}, we say that $f$ is a {\bf{\em cofinal element}} of $\Gamma$
if it is $\preceq$-cofinal for every left-order $\preceq$ on $\Gamma$.
The following should be clear from the discussion above.

\vsp

\begin{prop} {\em If a finitely-generated, left-orderable group $\Gamma$ has
a cofinal, central element, then no left-order on $\Gamma$ is of type III.}
\end{prop}

\vsp

\begin{small}
\begin{ex} In \S \ref{general-Conrad}, we introduced the group
$$\Gamma = \big\langle f,g,h \! : f^2 = g^3 = h^7 = fgh \big\rangle,$$
which is left-orderable but admits no nontrivial homomorphism into $(\mathbb{R},+)$ (hence no
left-order of type I). We claim that the central element $\Delta := fgh$ is cofinal. Indeed, if
$\Delta = f^2 = g^3 = h^7$ has a fixed point for a dynamical realization, then this is fixed
by $f,g,h$, hence by the whole group, which is impossible. As a consequence, every
left-order on $\Gamma$ is of type II.
\end{ex}

\begin{ex}
Recall that the center of the braid group $\mathbb{B}_n$ is
generated by the square of the so-called {\bf {\em Garside element}} $\Delta_n$. Moreover, one has
$$\Delta_n^2 = (\s_1 \s_2 \cdots \s_{n-1} )^n = (\s_1^2 \s_2 \cdots \s_{n-1})^{n-1}.$$
We next reproduce the proof given in \cite{cofinal} of that $\Delta_n^2$ is cofinal in $\mathbb{B}_n$.

\vsp\vsp\vsp

\noindent{\underline{Claim (i).}} If $\preceq$ is a left-order on $\mathbb{B}_n$ for which
$\Delta_n \succ id$, then for any braid $\sigma$ that is conjugate to either \esp
$\alpha_n := \s_1 \s_2 \cdots \s_n$ \esp or \esp $\beta_n := \s_1^2 \s_2 \cdots \s_n$, 
\esp we have \esp $id \prec \sigma \prec \Delta_n^2$.

\vsp\vsp

Indeed, as \esp $\sigma^k = \Delta_n^2$ \esp for $k$ equal to either $n$ or $n-1$,
we must have
$$id \prec \s \prec \s^2 \prec \s^3 \prec \ldots \prec \s^k = \Delta_n^2.$$

\vsp

\noindent{\underline{Claim (ii).}} If $\preceq$ is as above, then \esp
$\Delta_n^{-2} \prec \s_i \prec \Delta_n^2$, \esp for all $i \in \{1,\ldots,n-1\}$.

\vsp

Indeed, since $\Delta_n^2$ is central, by Claim (i) we have,
for all $\delta \in \mathbb{B}_n$,
$$\Delta_n^{-2} \prec \delta \alpha_n^{-1} \delta^{-1} \prec id, \qquad
id \prec \delta \beta_n \delta^{-1} \prec \Delta_n^2.$$
Since \esp $\beta_n \alpha_n^{-1} \!=\! \s_1$, \esp this yields \esp
$\Delta_n^{-2} \!\prec \delta \s_1 \delta^{-1} \!\prec\! \Delta_n^2,$ \esp
and since all the $\s_i$'s are conjugate between them, this shows the claim.

\vsp\vsp

\noindent{\underline{Claim (iii).}} The element $\Delta_n^2$
is cofinal in $\mathbb{B}_n$.

\vsp

Let $\preceq$ be a left-order on $\mathbb{B}_n$. Using again the fact that $\Delta_n^2$
is central, the set $\{\s \in \mathbb{B}_n \!: \Delta_n^{2r} \prec \s \prec \Delta_n^{2s}
\esp \mbox{ for some } r,s \mbox{ in } \mathbb{Z} \}$ is easily seen to be a subgroup of
$\mathbb{B}_n$. By Claim (ii), it contains all the $\sigma_i$'s. Thus, it coincides with
$\mathbb{B}_n$, which concludes the proof.

\end{ex}

\begin{rem} We do not know whether there exist type-III left-orders on the derived groups
$\mathbb{B}_n'$ for $n \geq 5$. Note that these groups do not admit left-orders of type I, since they 
admit no nontrivial homomorphism into the reals (see Example \ref{n geq 5}). By the preceding 
example, the restriction to $\mathbb{B}_n'$ of any left-order on $\mathbb{B}_n$ is of type II.
\end{rem}\end{small}

\vsp\vsp\vsp

\noindent{\bf A dynamical view.} As we showed in \S \ref{general-3} (and used a number of times), finitely-generated, 
left-ordered groups may be realized as groups of orientation-preserving homeomorphisms of the real line.
In what follows, we use this approach to visualize the dynamical differences between orders of
different types. For example, type-I left-orders are characterized as follows.

\vsp

\index{Order!type I}
\index{Measure!Radon measure}
\begin{prop} \label{teo type I preserves measure}
{\em If $\preceq$ is a left-order of type I on a
finitely-generated group $\Gamma$, then its dynamical realization
preserves a Radon measure on $\mathbb{R}$. Conversely,
any left-order induced (in a dynamical-lexicographic manner) 
from a faithful action of $\Gamma$ on the
real line that preserves a (nontrivial) Radon measure is of type I.}
\end{prop}

\noindent{\bf Proof.} Let $\preceq$ be a left-order of type I on $\Gamma$, and 
let $\Gamma_0$ be its maximal proper convex subgroup. Consider the dynamical
realization of $\preceq$. By convexity, $\Gamma_0$ fixes the interval $[a,b],$ where
$a,b$ are, respectively, the infimum and the supremum of the orbit of the origin under
$\Gamma_0$. Moreover, by Theorem \ref{super-conrad} and Corollary \ref{mas}, we have
that $\Gamma_0$ is normal in $\Gamma$, that $\Gamma_0 = \ker (\tau_\preceq)$,
where $\tau_\preceq:\Gamma \to (\mathbb R, +)$ is the Conrad homomorphism,
and that the induced order on $\Gamma/\Gamma_0$ is Archimedean.
In particular, the set $Fix(\Gamma_0)$ of global fixed points of
$\Gamma_0$, is $\Gamma$-invariant (hence infinite), and the
action of $\Gamma/\Gamma_0$ on $Fix(\Gamma_0)$ is free.

Now, if $Fix(\Gamma_0)$ is discrete (equivalently, if
$\tau_\preceq(\Gamma) \sim \Z$), then for each $x\in Fix(\Gamma_0)$,
the measure $\sum_{h\in \Gamma/\Gamma_0} \delta_{h(x)}$ is a
$\Gamma$-invariant Radon measure. If
$Fix(\Gamma_0)$ is non-discrete (equivalently, if
$\tau_\preceq(\Gamma)$ is a dense subgroup of $\mathbb R$), we may
proceed as in Example \ref{ex-free-action} to show that the action
of $\Gamma$ is continuously semiconjugate to an action by translations that factors
throughout $\Gamma/\Gamma_0$. Pulling back the Lebesgue measure by this
semiconjugacy, we obtain a $\Gamma$-invariant Radon measure.

Conversely, assume that an action of $\Gamma$ by orientation-preserving
homeomorphism of the real line preserves a Radon measure $\nu$. Then there is a
translation number 
\index{Translation number} homomorphism $\tau_\nu \!: \Gamma \to (\mathbb{R},+)$
defined by $\tau_\nu(g) := \nu ([y,g(y)) )$. (Recall that the value is independent of $y$,
due to invariance.) We claim that $\ker (\tau_\nu)$ is a convex subgroup for
any left-order induced from the action. Indeed, let $x \in \mathbb{R}$ be
the first reference point for inducing such a left-order on $\Gamma$ (see 
\S \ref{general-3}). On the one hand, if $x$ lies in the support of $\nu$, then
$\ker (\tau_\nu)$ coincides with the stabilizer of $x$, hence it is a convex
subgroup. On the other hand, if $x$ does not belong to the support of $\nu$,
let $I$ be the connected component of the complement of the support of $\nu$
containing $x$. At least one endpoint of $I$ is finite, which easily allows us to show
that for each $g \in \Gamma$, either $g (I) \cap I$ is empty or coincides with
$I$. It follows that the stabilizer of $I$ is a convex subgroup of
$\Gamma$ that coincides with $\ker (\tau_\nu)$.

Note that for all $g,h$ in $\Gamma$, the inequality $\tau_\nu(g) \!>\! \tau_\nu(h)$
implies $g (x) \!>\! h (x)$ for every $x \in \mathbb R$. It easily follows from this and the
discussion above that $\ker (\tau_{\nu})$ is the maximal convex subgroup. Finally, the
action of $\Gamma$ on $\Gamma / \ker (\tau_\nu)$ is Conradian, because it is
order-isomorphic to an action by translations. $\hfill\square$

\vsp\vsp

\begin{small}\begin{rem}\label{teo type I preserves measure extended}
In the proof above, the finite-generation hypothesis was only used in the direct implication
to ensure the existence of a maximal proper convex subgroup. Since this is known to exist
in some other situations (see, for instance, Exercise \ref{maximal-finite-rank}), the proposition
still holds in these cases.
\end{rem}\end{small}
\vsp\vsp

\index{Action!locally contracting}
\index{Action!globally contracting}

To deal with type-II and type-III left-orders, we closely follow \cite{DKN2} (compare \cite{malyu}).
We say that the action of a subgroup $\Gamma$ of $\mathrm{Homeo}_+ (\mathbb{R})$ is
{\bf{\em locally contracting}} if for every $x \in \mathbb{R}$ there is $y > x$ such that the
interval $[x,y]$ can be contracted to a point by a sequence of elements in $\Gamma$. We say
that the action is {\bf{\em globally contracting}} if such a sequence of contractions exists for any
compact subinterval of $\mathbb{R}$. We denote by $\widetilde{\mathrm{Homeo}}_+ (\mathrm{S}^1)$
the group of homeomorphisms of the line that are liftings of orientation-preserving circle
homeomorphisms. The next lemma is to be compared with Example \ref{actuar-convexo}.
\vsp

\begin{lem} \label{minimal-set}
{\em Every finitely-generated subgroup of $\mathrm{Homeo}_+ (\mathbb{R})$ preserves a
nonempty minimal closed subset of the line. This set is unique if no discrete orbit exists.}
\end{lem}

\noindent{\bf Proof.} Assuming that 
there are no global fix points, fix a point $x_0$ and a compact interval $I$ containing $x_0$ as well as all its
images under the (finitely many) generators of the group $\Gamma$. By obvious reasons, every orbit 
must intersect $I$; hence, this must be also the case of every closed
set of the line that is invariant under the action. Therefore, the standard argument
using Zorn's lemma to detect (nonempty) minimal sets may be applied by looking at the (compact) ``traces''
in $I$ of nonempty invariant closed subsets of the line. We leave the details to the reader (see
\cite[Proposition 2.1.12]{yo} in case of problems). 

To prove uniqueness, note that for a
closed invariant subset $K$, the set of acummulation points $K'$ is also closed and
invariant, hence $K' \!=\! K$ if $K$ is not a discrete orbit. Assume that $K$ is not the
whole line (otherwise, the uniqueness is obvious). It is then easy to see that every connected
component of the complement of $K'$ has a sequence of images converging to any point
in $K=K'$.  In other words, every orbit acummulates at $K$, which obviously implies the
uniqueness of the nonempty minimal invariant closed set.
$\hfill\square$

\vsp

\begin{thm}\label{T: topological structure}
{\em Let $\Gamma$ be a finitely-generated subgroup of $\mathrm{Homeo}_+ ({\mathbb{R}})$ whose
action admits no global fixed point. Then one of the following mutually-exclusive possibilities
occur:}

\vspace{0.08cm}

\noindent (i) {\em $\Gamma$ is semiconjugate to a group of translations;}

\vsp

\noindent (ii) {\em $\Gamma$ is semiconjugate to a minimal, locally contracting
subgroup of $\widetilde{\mathrm{Homeo}}_+(\mathrm{S}^1)$;}

\vsp

\noindent (iii) {\em $\Gamma$ is globally contracting.}
\end{thm}

\noindent{\bf Proof.} Assume there is no discrete orbit for the action. By Lemma \ref{minimal-set},
there is a unique minimal nonempty closed $\Gamma$-invariant subset $K$. In case $K$ is not the
whole line, collapse each connected component of the complement of $K$ to a point in order
to continuously semiconjugate $\Gamma$ to a group $\overline{\Gamma}$ whose action is minimal. If
$\Gamma$ preserves a Radon measure then, after semiconjugacy, this measure becomes a
$\overline{\Gamma}$-invariant Radon measure of total support and no atoms. Therefore,
$\overline{\Gamma}$ (resp. $\Gamma$) is conjugate (resp. semiconjugate) to a group of translations.

Suppose next that $\Gamma$ has no invariant Radon measure. Then the action of
$\overline{\Gamma}$ cannot be free. Otherwise, $\overline{\Gamma}$ would be conjugate
to a group of translations (see Example \ref{ex-free-action}), and the pull-back of the
Lebesgue measure by the semiconjugacy would be a $\Gamma$-invariant Radon measure.

Let $\bar{g} \in \overline{\Gamma}$ be a nontrivial element having fixed points,
and let $\bar{x}_0$ be a point in the boundary of the complement of $\mathrm{Fix}(\bar{g})$.
Then there is a left or right neighborhood $I$ of $\bar{x}_0$ that is
contracted to $\bar{x}_0$ under iterates of either $\bar{g}$ or its inverse.
By minimality, every $\bar{x}$ has a neighborhood that can be contracted to
a point by elements in $\overline{\Gamma}$. Coming back to the original action,
we conclude that every $x \!\in\! \mathbb{R}$ has a neighborhood that can be
contracted to a point by elements in $\Gamma$. Note that such a limit point
can be chosen arbitrarily in $K$; in particular, it may be chosen to belong to a
compact interval $I$ that intersects every orbit (as in the proof of Lemma \ref{minimal-set}).

For each $x\in \mathbb{R}$, let $M(x) \in \mathbb{R} \cup \{+\infty\}$ be
the supremum of the $y>x$ such that the interval $(x,y)$ can be contracted to a
point in $I$ by elements of $\Gamma$. Then either $M \equiv +\infty$, in which
case the group $\Gamma$ is globally contracting, or $M (x)$ is finite for every
$x\in \mathbb{R}$. In the latter case, $M$ induces a non-decreasing map
$\overline M \!: \mathbb{R} \rightarrow \mathbb{R}$ that commutes 
with all the elements in $\overline{\Gamma}$. Since the union of the intervals on which
$\overline M$ is constant is invariant under $\overline{\Gamma}$, the minimality of the action implies
that there is no such interval, that is, $\overline M$ is strictly increasing. Moreover,
the interior of $\mathbb{R} \setminus \overline M (\mathbb{R})$ is also invariant, hence
empty because the action is minimal. In other words, $\overline M$ is continuous.
All of this shows that $\overline M$ induces a homeomorphism of $\mathbb{R}$ into its
image. Since the image of $\overline M$ is $\overline{\Gamma}$-invariant, it must be the whole
line. Therefore, $\overline M$ is a homeomorphism from the real line to itself. Observe that
$\overline M (x)>x$ for any point $x$, which implies that $\overline M$ is conjugate to
the translation $x \mapsto x+1$. After this conjugacy, $\overline{\Gamma}$ becomes a
subgroup of $\widetilde{\mathrm{Homeo}}^+(\mathrm{S}^1)$.
$\hfill\square$

\vsp\vsp\vsp\vsp

The next proposition should now be clear to the reader.

\vsp

\index{Order!type I}
\index{Order!type II}
\index{Order!type III}

\begin{prop} \label{bluf}
{\em Let $\Gamma$ be a finitely-generated left-orderable group, and let $\preceq$ be a left-order on
it. Then $\preceq$ is of type I, II, or III if and only if its dynamical realization
satisfies property (i), (ii), or (iii)) above, respectively.}
\end{prop}

\vsp

So far, we haven't given any example of a left-orderable group all of whose left-orders are of type III. 
Actually, in an earlier version of this book, we explicitly asked for the existence of such a group, and we 
provided several (very strong) consequences of the eventual non-existence of them. It turns out, however, 
that these groups exist, though their construction is not at all easy. For explicit examples, we refer to 
\cite{HLNR} and \cite{MT}. It is worth mentioning that the example provided in the latter reference 
corresponds to the group that will be analyzed in detail in the final section of this book.

\chapter{PROBABILITY AND LEFT-ORDERABLE GROUPS}

\section{Amenable Left-Orderable Groups}
\label{super-witte}

\hspace{0.45cm} Starting from the work of von Neumann and Day, amenability became 
one of the deepest notions in the theory of infinite groups. Although there are many equivalent definitions 
(see, for instance, Exercise \ref{ejer:dos-equiv} below), we introduce this concept via von Newmann's 
original approach using means. 

\index{Group!amenable}
A {\bf {\em mean}} on a countable set \( \Gamma \) is a linear functional \(M\) on $\mathcal{L}^{\infty} (\Gamma)$ 
that satisfies: 

\vsp

\noindent -- ({\em Positivity}) If $\phi$ is non-negative, then $M (\phi) \geq 0$;

\vsp

\noindent -- ({\em Normalization}) If $1_{\Gamma}$ denotes the constant function equal to 1 along $\Gamma$, 
then $M (1_{\Gamma}) = 1$.

A countable group $\Gamma$ is said to be {\bf {\em amenable}} if there exists a mean \(M\) on 
$\Gamma$ that is invariant under right multiplication, that is:

\noindent -- ({\em Invariance})  For all $\phi \in \mathcal{L}^{\infty}(\Gamma)$ and all $g \in \Gamma$,  
one has $M (\phi) = M (\mathrm{\phi} \circ R_g)$, where $R$ is the right action of $\Gamma$ on 
$\mathcal {L}^\infty (\Gamma)$, that is, $R_g (\phi)(h) := \phi (hg)$.

\vsp

Among the many equivalent definitions, in the next exercise we highlight the one that will be useful in this section. 

\vsp

\begin{small}\begin{ejer} \label{ejer:dos-equiv}
Prove that a group is amenable if every action by homeomorphisms of a compact metric space admits an 
invariant probability measure. (See Appendix A of \cite{juschenko} or \cite{zimmer} in case of problems.)
\end{ejer}

\begin{ejer} Prove that every subgroup of an amenable group is amenable.
\end{ejer}\end{small}

\vsp

Our aim now is to discuss another nice result due to Witte Morris 
\cite{witte}. The theorem below was conjectured by Linnell in \cite{linnell-ind},
but it was already suggested by Thurston (see \cite[page 348]{Th}). 
\index{Witte Morris!theorem on amenable groups}

\vsp

\begin{thm} \label{witte-bien}
{\em Every amenable, left-orderable group is locally indicable.}
\end{thm}

\vsp

\index{Order!right-recurrent}
For the proof, we will say that a left-order $\preceq$ is {\bf\textit{right-recurrent}}
if for every pair of elements $f,h$ in $\Gamma$ such that $f \!\succ\! id$, there
exists $n \!\in\! \mathbb{N}$
satisfying $fh^n \succ h^n$. Note that every
right-recurrent order is Conradian. (The converse does not hold; see Example
\ref{no-tiene-pelo}.) As subgroups of amenable groups are amenable, 
this implies that Theorem \ref{witte-bien} follows from the next proposition.

\vspace{0.1cm}

\begin{prop} {\em If \esp $\Gamma$ is a finitely-generated, amenable,
left-orderable group, then $\Gamma$ admits a right-recurrent order.}
\label{witte-fund}
\end{prop}

\vspace{0.1cm}

To prove this proposition, we will need the following weak form of the Poincar\'e 
recurrence theorem. We recall the proof for the reader's convenience.
\index{Poincar\'e recurrence theorem}

\vspace{0.1cm}

\begin{thm} {\em If $S$ is a measurable map that preserves a probability
measure $\mu$ on a space $X$, then for every measurable subset $A$
of \esp $X$ and $\mu$-a.e. point $x \! \in \! A$, there
exists $n \! \in \! \mathbb{N}$ such that $S^n(x)$ belongs to $A$.}
\label{rec-poincare}
\end{thm}

\noindent{\bf Proof.} The set $B$ of points in $A$ that do not come back to $A$ under
iterates of $S$ is $A \setminus \bigcup_{n \in \mathbb{N}} S^{-n}(A)$. One easily
checks that the sets $S^{-i}(B)$, with $i \geq 1$, are pairwise disjoint. Since
$S$ preserves $\mu$, these sets have the same measure, and since the total mass of
$\mu$ equals $1$, the only possibility is that this measure equals zero. Therefore,
$\mu(B)=0$, that is, $\mu$-a.e. point in $A$ comes back to $A$ under some
iterate of $S$. $\hfill\square$

\begin{small}\begin{ejer}\label{ejer:unbounded}
In the framework above, show that for $\mu$-a.e. point $x \! \in \! X$, 
the set of positive integers $n$ such that $S^n(x) \!\in\! A$ is unbounded.
\end{ejer}\end{small}

\index{Order!right-recurrent}
\noindent{\bf Proof of Proposition \ref{witte-fund}.} 
By Exercise \ref{ejer:dos-equiv}, if $\Gamma$ is a (countable) left-orderable 
amenable group, then its action on (the compact metric space) $\mathcal{LO}(\Gamma)$
preserves a probability measure $\mu$. We claim that $\mu$-a.e. point in
$\mathcal{LO}(\Gamma)$ is right-recurrent. To show this, for each $g \!\in\! \Gamma$, let us
consider the subset $V_g$ of $\mathcal{LO}(\Gamma)$ formed by the left-orders $\preceq$ on $\Gamma$
such that $g \!\succ\! id$. By the Poincar\'e recurrence theorem, for each $f \in \Gamma$, the set
\esp $B_g(f) := V_g \setminus \bigcup_{n \in \mathbb{N}} f^{-n}(V_g)$ \esp has null $\mu$-measure.
Therefore, the measure of $B_g := \bigcup_{f \in \Gamma} B_g (f)$ is also zero, as is the 
measure of $B := \bigcup_{g \in \Gamma} B_g$. Let us consider an arbitrary element $\preceq$
in the ($\mu$-full measure) set $\mathcal{LO}(\Gamma) \setminus B$. Given $g \succ id$ and
$f \in \Gamma$, from the inclusion $B_g(f) \subset B$ we deduce that $\preceq$ does not belong
to $B_g(f)$. Thus, there exists $n \!\in\! \mathbb{N}$ such that $\preceq$ belongs to $f^{-n} (V_g)$,
hence $\preceq_{f^n}$ is in $V_g$. In other words, one has $g \succ_{f^{n}} id$, that is, $gf^n \succ f^n$.
Since $g \succ id$ and $f \in \Gamma$ were arbitrary, this shows the right-recurrence of $\preceq$.

\index{Witte Morris!a group without recurrent orders}
\begin{small}\begin{ex} \label{no-tiene-pelo}
Following \cite[Example 4.5]{witte}, we next show that there exist
$C$-orderable groups that do not admit right-recurrent orders.
This is the case of the semidirect product
$\Gamma \!=\! \mathbb{F}_2 \ltimes \mathbb{Z}^2$,
where $\mathbb{F}_2$ is any free subgroup of $\mathrm{SL}(2,\mathbb{Z})$ acting
linearly on $\mathbb{Z}^2$. (Such a subgroup may be taken of finite index.)
Indeed, that $\Gamma$ is $C$-orderable follows from the local indicability of
both $\mathbb{F}_2$ and $\mathbb{Z}^2$. Assume throughout that $\preceq$ is a
right-recurrent left-order on $\Gamma$. For a matrix $f \!\in  \mathbb{F}_2$
and a vector $v \!=\! (m,n) \in \mathbb{Z}^2$, let us denote by
$\bar{f}$ and $\bar{v}$ the corresponding elements in $\Gamma$, so that
$f(v) = \bar{f} \bar{v} \bar{f}^{-1}$. Let $\tau$ be the Conrad's homomorphism
associated to the restriction of $\preceq$ to $\mathbb{Z}^2$, so that we have
$v \succ id$ whenever $\tau(v) > 0$, and $\tau(v) \geq 0$ for all $v \succ id$
(see Corollary \ref{mas}). Let $f$ be a hyperbolic matrix in $\mathbb{F}_2$,
with positive eigenvalues $\alpha_1, \alpha_2$ and corresponding eigenvectors $v_1, v_2$
in $\mathbb{R}^2$. Since $v_1$ and $v_2$ are linearly independent, we may assume that
$\tau (v_1) \neq 0$. Furthermore, we may assume that $\tau (v_1) > 0$ and $\alpha_1 > 1$
after replacing $v_1$ with $-v_1$ and/or $f$ with $f^{-1}$,
if necessary. Let $L : \mathbb{R}^2 \rightarrow \mathbb{R}$ be the (unique) linear
functional that satisfies $L (v_1) = 1$ and $L (v_2) = 0$. Given any $v \in \mathbb{Z}^2$
such that $\tau (v) > 0$, right-recurrence provides us with an increasing sequence $(n_i)$
such that $\bar{v} \bar{f}^{-n_i} \succ \bar{f}^{-n_i}$ for every $i$. This implies that
$\bar{f}^{n_i} \bar{v} \bar{f}^{-n_i} \succ id$, hence $\tau (f^{n_i}(v)) \geq 0$. Since
$$\lim_{i \rightarrow \infty} \frac{\tau (f^{n_i}(v))}{\alpha_1^{n_i}} =
\tau \left( \lim_{i \rightarrow \infty} \frac{f^{n_i} (v)}{\alpha_1^{n_i}} \right)
= \tau (L(v) v_1) = L (v) \tau (v_1),$$
we conclude that $L (v) \geq 0$.
Since $v$ is an arbitrary element of $\mathbb{Z}^2$ satisfying $\tau (v) > 0$, this
necessarily implies that \esp $\ker (\tau) \!=\! \ker (L)$ \esp is an eigenspace of $f$.
But $f$ is an arbitrary hyperbolic matrix in $\mathbb{F}_2$, and it is easy to show
that there are hyperbolic matrices in $\mathbb{F}_2$ with no common eigenspace. This
is a contradiction.
\end{ex}
\end{small}

\vsp

\noindent{\bf An extension for left-orderable groups without free subgroups~?} A prototype of 
a non-amenable group is the free group in two generators; see Exercise \ref{ejer:no invariant probability} below. 
Since every subgroup of an amenable group is amenable, 
every group containing a (non-Abelian) free group is also non-amenable. The converse is, however, 
false, even in the framework of bi-orderable, finitely-presented groups (see \S \ref{sec-amenable-ex}).

\begin{small}\begin{ejer}\label{ejer:no invariant probability}
Show that the group generated by two homeomorphisms of the circle having nonempty disjoint 
sets of fixed points does not preserve any probability measure on the circle. Deduce that the 
free group on two generators is not amenable.
\end{ejer}\end{small}

It is unknown whether Theorem \ref{witte-bien} extends to groups without free subgroups,
that is, whether left-orderable groups not containing $\mathbb{F}_2$ are locally indicable.
(See \cite{linnell-free} for an interesting result pointing in the affirmative direction.) 
\index{Engel's condition}
A relevant class of groups that do not contain free subgroups consists of those satisfying
a nontrivial {\bf{\em law}} (or {\bf{\em identity}}). This is a
reduced word $W \!=\! W (x_1,\ldots,x_k)$ in positive and negative powers such
that $W (g_1,\ldots,g_n)$ is trivial for {\em every} $g_1,\ldots,g_n$ in the group. For
instance, Abelian groups satisfy a law, namely $W_1 (x_1,x_2) := x_1 x_2 x_1^{-1} x_2^{-1}$.
Nilpotent and solvable groups also satisfy group laws. Another important but less understood
family is the one given by groups satisfying an {\bf{\em Engel condition}} $W^E_k$, where
$$W^E_k (x_1,x_2) := W^E_1 \big( W^E_{k-1} (x_1,x_2), x_2 \big),
\qquad W^E_1 (x_1,x_2) := x_1 x_2 x_1^{-1} x_2^{-1}.$$
It is an open question whether left-orderable groups satisfying an Engel condition
must be nilpotent. This is known to be true if the group is Conrad-orderable (see
\cite[Theorem 6.G]{glass}). In other words, if $\Gamma$ is an Engel group having a
left-order without resilient pairs, then $\Gamma$ is nilpotent. In this direction, the next
proposition becomes interesting, and shows the pertinence of Question \ref{question-cle}.
The (easy) proof is left to the reader. (See \cite{gol} for more on this.)

\vsp

\begin{prop} \label{laws}
{\em If \esp $\Gamma$ is a left-orderable group satisfying
a law, then there exists $n \!\in\! \mathbb{N}$ such that no left-order
on $\Gamma$ admits an $n$-resilient pair.}
\end{prop}

\vsp\vsp

\index{Order!locally invariant}
\noindent{\bf Locally-invariant orders on amenable groups.} The ideas involved in the proof 
of Theorem \ref{witte-bien} yield interesting results for other type of orders on amenable
groups. The next result is due to Linnell and Witte Morris \cite{LIOs-on-amenable}.

\vsp

\begin{thm} \label{LIO-amenable}
{\em Every amenable group admitting a locally-invariant order is left-orderable (hence locally indicable).}
\end{thm}

\noindent{\bf Proof.} As for Theorem \ref{witte-bien},
we may assume that $\Gamma$ is finitely-generated.

First, it is not hard to extend the claim of Exercise \ref{LIOs-on-Z} to
describe the restriction of a locally-invariant order to any left coset of a cyclic subgroup:
For every $f \in \Gamma$ and $g \neq id$, either
$$f g^n \prec f g^{n+1} \quad \mbox{for all } n \in \mathbb{Z},$$
or
$$f g^{n} \prec f g^{n-1} \quad \mbox{for all } n \in \mathbb{Z},$$
or there exists $\ell \in \mathbb{Z}$ such that
$$fg^n \prec fg^{n+1} \quad \mbox{for all } n \geq \ell \quad \mbox{and}
\quad fg^{n} \prec fg^{n-1} \quad \mbox{for all } n < \ell.$$
We next argue that for amenable groups, there is a locally-invariant order
for which the third possibility never arises.

Recall from Exercise \ref{space-of-LIOs} that the space of locally-invariant orders is a
compact topological space, which is metrizable whenever $\Gamma$ is countable. The
group $\Gamma$ acts on $\mathcal{LIO}(\Gamma)$ by left and right translations.
Since $\Gamma$ is amenable, both actions preserve probability measures on
$\mathcal{LIO}(\Gamma)$. Let $\mu$ be a probability measure that is invariant
under the right action. We leave to the reader the task of showing that a
generic locally-invariant order $\preceq$ is {\em strongly right-recurrent}. More precisely, there is a
subset $A$ of full $\mu$-measure such that for every $\preceq$ in $A$, the following
happens: If $f \prec g$, then given $h \in \Gamma$, the set of integers $n$ such
that $fh^n \prec gh^n$ is unbounded in both directions. (Compare Exercise \ref{ejer:unbounded}.) 
Since in the third case
above this property fails, we conclude that a generic locally-invariant order is either the canonical
one or its reverse whenever restricted to a left-coset of a cyclic subgroup.
\index{Order!right-recurrent}

We next show that every $\preceq$ in $A$ is a left-order. Indeed, by the definition
of locally-invariant order, for every $g \neq id$, either $g \succ id$ or $g^{-1} \succ id$ holds.
Both inequalities cannot hold simultaneously, otherwise the restriction of $\preceq$ to
the cyclic subgroup $\langle g \rangle$ wouldn't be neither the canonical order nor its reverse.
Therefore, the positive cone $P := \{g \!: g \succ id\}$ is disjoint from its inverse, and their
union covers $\Gamma \setminus \{ id \}$.

It remains to show that $P$ is a semigroup. Assume for a contradiction that $g,h$
in $P$ are such that $gh \notin P$. Then $gh \prec id$.  Thus $g \succ id \succ gh$. 
Using the property of a locally-invariant order, one easily checks that, necessarily,
$gh^2 \prec gh$. More generally,
$$g \succ gh \succ gh^2 \succ gh^3 \succ \ldots.$$
Therefore, $gh^n \prec id$, for all $n \geq 1$. However, due to the right-recurrence of $\preceq$, there
is some $n \in \mathbb{N}$ such that $gh^n \succ h^n \succeq h \succ id$. This is a contradiction.
$\hfill\square$

\vsp

The theorem above makes the following question natural. 

\begin{question}
Does there exist an amenable U.P.P. group that is not left-orderable~? (See \S \ref{sec-upp}.)
\end{question}


\section{Non-Amenable, Left-Orderable Groups with no Free Subgroups}
\label{sec-amenable-ex}

\hspace{0.45cm} The natural problem of finding non-amenable groups without (non-Abelian) free subgroups goes back to 
von Neumann. This was first solved by Ol'shanskii in \cite{olshanski} and soon after by Adian \cite{ad1}, at the beginning of 
the eighties. However, their examples do not admit a finite presentation, and also contain many torsion elements. 
Indeed, torsion is fundamental for their constructions; for instance, Adian proves that the free Burnside group 
$B(2,n)$ introduced in Example \ref{burnside-ex} is non-amenable for large enough $n$. (Note that, in this group, 
{\em every} element has finite order.) The first example of a finitely-presented, non-amenable group without free 
subgroups was constructed much later in \cite{OS}; however, this group still has many torsion elements.

Here we present a construction of a bi-orderable (hence torsion-free) non-amenable group having no free subgroups 
which, moreover, admits a finite presentation. This example comes from the beautiful recent work of Lodha and 
Moore \cite{lodha}, who were able to isolate a particular finitely-presented group inside a much larger family of 
groups previously studied by Monod \cite{monod}. 

\index{Group!Lodha-Moore}

\subsection{(Non-)amenable relations}
\label{sec amenable actions}
\index{Group!piecewise projective}
\hspace{0.45cm} In this section, we will show that $\mathrm{PP}_+ (\mathbb{R})$, the group of orientation-preserving, piecewise-projective 
homeomorphisms of the real line (see \S \ref{ejemplificando-relatives of F}), contains many countable subgroups that are non-amenable 
\cite{monod}. The key ingredient is the notion of amenable equivalence relation and  a result of Carri\`ere and Ghys \cite{CG} that we present as Theorem \ref{th:CG} below. 

Consider an action of a countable group $\Gamma$ on a measure space $X$ by measurable maps. 
Assume that the images of zero-measure sets under group elements have zero measure. The {\bf {\em 
orbital equivalence relation}} associated to $\Gamma$ on $X$ is the one whose equivalence classes are the orbits of $\Gamma$. (We denote the orbit of the point $x \in X$ by $\Gamma x$.) Such a relation (and the underlying action)  
is said to be {\bf{\em amenable}} if for almost every $x \!\in\! X$ there is a mean $M_x \!: \mathcal{L}^{\infty} (\Gamma x) \to \mathbb{R}$ that satisfies:
\index{Action!amenable}

\vsp

\noindent -- ({\em Invariance}) For almost every $x \in X$ and every $y\in \Gamma x$, one has $M_x  = M_y$;

\vsp

\noindent -- ({\em Measurability}) If $\psi$ is a bounded measurable function defined on the {\bf {\em graph of the equivalence relation}} \esp $\mathcal{G} := \big\{ (x,y)\in X\times X \!:   y \in \Gamma x \big\},$ \esp 
then $x \mapsto M_x (\psi (x, \cdot))$ is a measurable function from $X$ into $\mathbb{R}$.

\vsp

It is worth to stress that this notion of amenable action only depends on the associated orbital equivalence 
relation.  However, it is a kind of extension of the notion of group amenability, as the next exercise shows. 

\begin{small}\begin{ejer} \label{ej: amenable action}
Show that if \(\Gamma\) is an amenable group, then the orbital equivalence relation induced by any 
measurable action of $\Gamma$ is amenable. 

\noindent{\underline{Hint.}}  
If $M$  is a mean on $\Gamma$, define the family of means on the $\Gamma$-orbits by letting $M_x (\phi) := M (\tilde{\phi}_x)$, 
where $\tilde{\phi}_x (g) := \phi (gx)$. Show that this family satisfies the measurability axiom, and that if \(M\) is right invariant, it also satisfies the invariance axiom. 
\end{ejer}

\begin{rem} 
The class of amenable actions is considerably larger than that of actions of amenable groups. A nice 
example witnessing this is the action of a lattice of \( \text{PSL}(2,\mathbb R) \) on the projective real 
line \(\mathbb P^1(\mathbb R)\) (endowed with the Lebesgue measure). Indeed, such a lattice is 
never amenable (as it contains free subgroups), though its action on \(\mathbb P^1(\mathbb R)\) is 
amenable (see \cite[Corollary 4.3.7]{zimmer}). 
\end{rem}\end{small}

\vsp

The next result gives great insight on the projective action of certain subgroups of $\text{PSL} (2,\mathbb R)$ 
on $\mathbb{P}^1(\mathbb{R}) \!\sim\! \mathbb{R} \cup \{ \infty \}$, where $\infty \!\sim\! [1:0]$. 

\index{Carri\`ere-Ghys Theorem}
\begin{thm} \label{th:CG}
{\em Let $\Gamma\subset \mathrm{PSL}(2,\mathbb R)$ be a countable subgroup. If $\Gamma$ contains a 
non-discrete copy of $\mathbb F_2$, then the equivalence relation given by the $\Gamma$-orbits of its 
projective action on $\mathbb P^1(\mathbb R)$  is non-amenable.}
\end{thm}

\noindent{\bf Proof.}  Assume for a contradiction that the relation induced by the $\Gamma$-action on 
$\mathbb P^1(\mathbb R)$ is amenable. Then the induced action of the free subgroup $\mathbb F_2$ is also amenable.  
Indeed, we may define the family $M'$ of means along the $\mathbb F_2$-orbits by letting $M'_{x} (f) := M_x (\bar{f})$, 
where $\bar{f}$ coincides with $f$ along the $\mathbb F_2$-orbit of $x$ and equals $1$ on 
$\Gamma x \setminus \mathbb F_2 x$. Therefore, we may assume that $\Gamma = \mathbb F_2$.

Let $\{a, b\}$ be a free generating set of $\mathbb{F}_2 = \Gamma$, and let $A$ (resp. $B$)
be the set of elements in $\mathbb{F}_2$ whose reduced forms finish with a nontrivial power of $a$ (resp. $b$).
Assuming the existence of the linear functionals $M_x$ as above, 
let $\psi = \psi_a \!: \mathbb P^1 (\mathbb R) \to [0,1]$ be defined by $\psi (x) := M_x (1_{{A x}})$, where $1_{{A x}}$ stands for the characteristic function
of the corresponding set $Ax := \{h (x) \!: h \in A\}$. 
By applying the ({\em Measurability}) axiom to the function $(x,y) \mapsto 1_{{Ax}} (y)$, we conclude that 
$\psi$ is a measurable function. Moreover, the action of $\mathbb{F}_2$ on $\mathbb P^1(\mathbb R)$ is almost 
everywhere free, because every nontrivial element in $\mathrm{PSL} (2, \mathbb{R})$ fixes at most two points. 
Since the sets $Ab^i$ are pairwise disjoint for $i \in \mathbb{Z}$, by the ({\em Positivity}) and 
({\em Normalization}) axioms of the definition of a mean, we obtain a.e. 
$$0\leq \sum_i M_x(1_{Ab^ix})\leq 1.$$
Since by ({\em Invariance}) we also have
$$\psi (b^i x) = M_{b^i x}(1_{Ab^i x}) = M_x(1_{Ab^i x}),$$
we obtain a.e.
\begin{equation}\label{eq C-G}
0 \leq \sum_i \psi (b^i x) \leq 1.
\end{equation}
In particular, a.e it holds that, if $\psi (x)>1/2$, then $\psi (b^i x)<1/2$ for all $i\neq0$. 

Using ({\em Invariance}) and the fact that the action is almost free, one concludes that $M_x (1_{\{x\}}) = 0$ 
holds for almost every $x$. Therefore, by ({\em Normalisation}), we have  $\psi_b := 1-\psi_a$. This allows us 
to conclude in the very same way as above that, a.e., if $\psi (x)<1/2$, then $\psi (a^i x)>1/2$ for all $i\neq 0$.

These properties fit into the framework of the classical Klein's ping-pong argument (see Exercise 
\ref{ejer-ping-pong}). Namely, for the measurable subsets $P$ and $Q$ of $\mathbb P^1 (\mathbb R)$ 
defined by
$$P:=\{x\in \mathbb P^1 (\mathbb R) \mid \psi(x)<1/2\} \;\; \text{ and } \;\; Q:=\{x\in G\mid \psi(x)>1/2\},$$
we have $a^i P\subset Q$ and $b^jQ\subset P$ for every nonzero $i,j$. 
In particular, for every element $g\in \mathbb F_2$ that (in reduced form) begins and finishes with a power of $a$, 
we have $g (P) \subset Q$.\footnote{Strictly speaking, these containments hold up to sets of null measure, but 
this is enough to make the argument work.} 

By hypothesis, there is a sequence of nontrivial elements $g_n \in \mathbb{F}_2$ converging to the identity in 
$\mathrm{PSL} (2,\mathbb{R})$. We claim that we may take these elements to begin and end with nontrivial 
powers of $a$ (and hence $g_n (P) \subset Q$ for every $n$). Indeed, if $g_n$ begins or finishes by a power 
of $b$, then at least one of the following elements $ag_n a^{-1}$ or $a^{-1}g_n a$ begins and finishes by 
powers of $a$, and such a conjugate stays close to the identity. 

From the Lebesgue Density Theorem, it follows that 
$$\lim_{n \to \infty} \mu \big( g_n (P) \cap P \big) = \mu(P),$$ 
where $\mu$ stands for the Lebesgue measure on $\mathbb P^1(\mathbb R)$. Since $g_n (P) \subset Q$ and $P$ 
and $Q$ are disjoint, we obtain that $\mu(P)=0$. The same argument shows that $\mu (Q)= 0$. In particular, 
almost surely the function $\psi$ takes the value $1/2$, but this contradicts 
the finiteness of the series \eqref{eq C-G}.  $\hfill\square$



\vspace{0.25cm}

Before stating the next result, we introduce some notation. 
Given a subgroup $\Gamma$ of $\text{PSL} (2,\mathbb R)$ acting on $\mathbb{P}^1(\mathbb{R})$, we denote 
$\mathrm{P} (\Gamma)$ the subgroup of $\text{Homeo}_+ (\mathbb R)$ consisting of the homeomorphisms that coincide 
with the restriction of an element of $\Gamma$ on each piece of a division of the real line into finitely many intervals. Observe that $\Gamma$ is not assumed to fix $\infty$, whereas $P(\Gamma)$ fixes it.

\vsp

\begin{prop} \label{p:relation}
{\em Let $\Gamma$ be a subgroup of $\mathrm{PSL}(2,\mathbb R)$. Assume that $\Gamma$ contains a nontrivial translation $x \mapsto x + t$. 
Then the orbit relations induced on $\mathbb R\setminus \Gamma \infty$ by both the actions of $\Gamma$ and $\mathrm{P} (\Gamma)$ coincide a.e.}
\end{prop}

\noindent {\bf Proof.} After conjugating by a suitable affine map, we may assume that the translation $x\mapsto x+1$ is contained in 
$\Gamma$. Let $x,y$ be any pair of points not lying in the $\Gamma$-orbit of $\infty \!\in\! \mathbb P^1(\mathbb R)$ but lying in the 
same $\Gamma$-orbit. 
We need 
to show that there exists an element of $\mathrm{P} (\Gamma)$ sending $x$ to $y$. To do this,  
let $g \in \Gamma$ be such that $y = g(x)$, say $g(z)= \frac{az+b}{cz+d}$ (with $ad-bc\neq 0$). If $c= 0$, then $g$ already fixes $\infty$, and 
the restriction of $g$ to $\mathbb R$ belongs to $\mathrm{P} (\Gamma)$, hence we are done in this case. Assume next that 
$c\neq 0$. For each $n\in \mathbb Z$ of sufficiently large modulus, the equation $g(z) = z-n$ has two solutions, namely 
$$z_\pm = \frac{a-d+cn}{2c} \left( 1 \pm   \sqrt{1+\frac{4c(dn+b)}{(a-d+cn)^2}} \right).$$ 
One easily checks that these satisfy
\[  z_+ \sim_{|n|\rightarrow \infty} n \ \qquad \text{and} \ \qquad \ z_- \rightarrow _{|n|\rightarrow \infty} -\frac{d}{c} .\]
Since $x \neq -\frac{d}{c} = g^{-1}(\infty)$, by choosing $n$ large enough or small enough according to whether $x >  -\frac{d}{c}$ 
or $x < -\frac{d}{c}$, we have that $x$ lies inside the interval $I$ with endpoints $z_-$ and $z_+$. Define $\hat{g}(z) := g(z)$ if $z$ 
belongs to $I$, and $\hat{g}(z) := z-n$ otherwise. Then $\hat{g}$ is an element of $\mathrm{P} (\Gamma)$ that sends $x$ 
to $y$, as desired.  $\hfill\square$

\vspace{0.3cm}

\begin{cor}
{\em If $\Gamma\subset \mathrm{PSL} (2,\mathbb R)$ is a countable subgroup containing a non-discrete copy of $\mathbb F_2$ and a 
nontrivial translation $x\mapsto x+t$, then the group $\mathrm{P} (\Gamma)$ is non-amenable and does not contain any copy of $\mathbb F_ 2$.}
\label{cor:cg}
\end{cor}

\noindent {\bf Proof.} Theorem \ref{th:CG} shows that the relation induced by \(\Gamma\) on \(\mathbb P^1 (\mathbb R)\) is non-amenable, and the same is true on the set  \(\mathbb P^1 (\mathbb R) \setminus \Gamma \infty= \mathbb R \setminus \Gamma \infty\), since the orbit \(\Gamma \infty\)  has null measure. Proposition \ref{p:relation} shows that, on \(\mathbb R \setminus \Gamma \infty\), the orbits of \( \Gamma\) and \( \mathrm{P} (\Gamma)\) are the same, hence the relation induced by \(  \mathrm{P} (\Gamma)\) on \(\mathbb R \setminus \Gamma\infty\) is non-amenable as well. Since the relation induced by any action of an amenable group is amenable (see Exercise \ref{ej: amenable action}), this shows that the group \(\mathrm{P} (\Gamma)\) is not amenable. Finally, the fact that \(  \mathrm{P} (\Gamma)\) does not contain any copy of $\mathbb F_2$ comes from Theorem \ref{th:BS}. $\hfill\square$

\vsp
\begin{rem} The preceding corollary as well as Theorem \ref{th:CG} 
can be extended under the weaker hypothesis of density of $\Gamma$. 
This is due to another result of \cite{CG} claiming that every countable dense 
subgroup of $\mathrm{PSL} (2, \mathbb R)$ contains a non-discrete copy of 
$\mathbb F_2$. It should be noted that, in 
\cite{breuillard gelander}, it is shown that if $\Gamma$ is a countable dense subgroup of a 
connected, real, semi-simple Lie group $G$ (for instance, $\mathrm{PSL}(2,\mathbb{R})$), 
then $\Gamma$ contains a copy of $\mathbb F_2$ which is also dense in $G$.
\end{rem}


\subsection{A finitely-presented version}
\label{ss:fin-pres}

In this section, we let $\Gamma \subset \text{PSL}(2, \mathbb R)$ be the group generated by
$$\tilde a =\left( \begin{array}{cc} 1 & 1\\  0 & 1  \end{array}\right)\, ,\; \tilde b =\left( \begin{array}{cc} 0 & 1\\  -1 & 0  \end{array}\right)\, 
\text{ and } \;\tilde c =\left( \begin{array}{cc} \sqrt{2} & 0\\  0 & 1/\sqrt{2}  \end{array}\right).$$
Note that the first two elements generate $\text{PSL}(2,\Z )$. 
By Theorem \ref{th:BS}, $\mathrm{P}(\Gamma)$ has no free subgroup. However, using the results of the preceding section, we will prove that the group $\mathrm{P} (\Gamma)$ is non-amenable. Moreover, following \cite{lodha}, we 
will show that $\mathrm{P} (\Gamma)$ contains a finitely-presented subgroup that is also non-amenable.

\vsp\vsp\vsp

\noindent{\bf The action of $\mathrm{P} (\Gamma)$ on $\mathbb{P}^1 (\mathbb{R})$ is non-amenable.} Since  
$\tilde{a}$ is the translation by 1, the non-amenability of the action will follow from Theorem \ref{th:CG} provided we check that $\Gamma$ contains a 
copy of $\mathbb F_2$ that is non-discrete in $\text{PSL}(2,\mathbb R)$. To do this, fix $n\geq 1$, and consider
$$g := \left( \begin{array}{cc} 0 & 1\\  -1 & 1/2^n \end{array}\right) = \tilde c^{n} \tilde b \tilde a^{-1} \tilde c^{n} \in \Gamma.$$
Letting $\delta := 1/2^n$, we see that the eigenvalues of $g$ are
$$\lambda_\pm= \frac{\delta}{2} \pm \frac{i \sqrt{4-\delta^2}}{2}.$$
Thus, the element $g \in \mathrm{PSL} (2,\mathbb R)$ is elliptic, because $\lambda_+\lambda_-=1$ and both 
$\lambda_{-}, \lambda_+$ have nontrivial imaginary part. Therefore, the projective action of $g$ is 
conjugate to that of a rotation.

We claim that $g$ acts hyperbolically on $\mathbb Q_2\times \mathbb Q_2$, where $\mathbb Q_2$ denotes the 2-adic rationals 
(see \cite{koblitz} for background). This means that none of the \(2\)-adic norms $|\lambda_\pm|_2$ of $\lambda_\pm$ is equal to 1. To check this, we write
$$\lambda_\pm=\frac{1}{2^{n+1}} \,w_\pm,$$
where $w_\pm = 1 \pm i \sqrt{2^{2n+2}-1}$. Thus, we need to show that $|w_\pm|_2\neq 1/2^{n+1}$. Looking for a contradiction, 
we assume that $|w_-|_2=|w_+|_2=1/2^{n+1}$.  Then
$$1 
= \left| \frac{w_+ w_+}{w_-w_+} \right|_2 
= \left|  \frac{1 + i \sqrt{2^{2n+2}-1}}{2^{2n+1}} + 1 \right|_2
= \left| \frac{w_+}{2^{2n+1}} + 1\right|_2.$$
But for an ultrametric norm (such as the 2-adic norm), we have that $|x|<|y|$ implies that $|x-y|=|y|$. 
Since $\left| \frac{w_+}{2^{2n+1}} \right|_2 = \frac{2^{2n+1}}{2^{n+1}} = 2^n$, letting 
$x:=  \frac{w_+}{2^{2n+1}}  + 1$ and $y :=  \frac{w_+}{2^{2n+1}} $, this yields 
$1 = \left|\frac{w_+}{2^{2n+1}}\right|_2=2^n$, which is the desired contradiction. 

Thus, $g$ acts hyperbolically on $\mathbb Q_2\times \mathbb Q_2$, hence it has infinite order.\footnote{A nice 
consequence of this argument is that $g$ is conjugate to a rotation on $\mathbb P^1(\mathbb R)$ by an angle 
which is an irrational multiple of $\pi$; otherwise, it would have a finite order. Such a rotation will be called an 
{\bf \em irrational rotation}, for short.}  
Now, since $\Gamma$ is a non-solvable group, there is a conjugate $f$ of $g$ such that $f$ and $g$ do not 
share any eigenvector (this is an easy exercise; see \cite{katok-s} in case of problems). By a ping-pong argument 
applied to the action of $\langle f, g \rangle$ on $\mathbb Q_2 \times \mathbb Q_2$, it follows that $g$ and $f$ 
generate a free group. Finally, this free group is non-discrete in $\mathrm{PSL} (2,\mathbb R)$, since $g$ is 
conjugate to an irrational rotation.

\begin{small}\begin{rem}
The fact that the element $g$ above is of infinite order has a nice arithmetic consequence. Namely, 
since the angle of rotation $\theta$ of an element $h \in \mathrm{PSL}(2,\mathbb{R})$ satisfies \, 
$2 \cos (\theta) = \pm tr (h)$, \, one concludes that $\arccos (1/2^{n+1})$ is an irrational multiple 
of $\pi$, for each $n \geq 1$.
\end{rem}\end{small}

\vsp\vsp\vsp

\noindent{\bf A finitely presented, non-amenable subgroup.} This constitutes the main contribution of \cite{lodha}. 
We first provide a crucial construction that will allow us to analyze $P(\Gamma)$ in combinatorial terms.

\begin{small}\begin{ejer} \label{hurwitz}
The Hurwitz application is the map $\phi \!: \{0,1\}^{\mathbb{N}} \to [0,\infty]$ recursively defined by 
$$\phi (0 \xi) := \frac{1}{1+\frac{1}{\phi(\xi)}} \qquad \mbox{and} \qquad \phi (1 \xi) := 1 + \phi (\xi), 
\qquad \mbox{ where } \xi \in \{0,1\}^{\mathbb{N}}.$$

\noindent (i) Check that, for $n_1 \geq 0$ and all positive integers $n_2,n_3,\ldots$,
$$\phi (1^{n_1} 0^{n_2} 1^{n_3} 0^{n_4} \ldots) = n_1 + \frac{1}{n_2+\frac{1}{n_3+\frac{1}{n_4+\ldots}}}.$$

\noindent (ii) Show that $\phi$ is one-to-one, except at points in $\{0,1\}^{\mathbb{N}}$ that become eventually constant 
yet are not constant. Check that these points map under $\phi$ to the rational points in $]0,\infty[$, and that the 
non-injectivity comes from the fact that, for every finite sequence $s$ of 0's and 1's, 
$$\phi (s 0 \underbar{1}) = \phi (s 1 \underbar{0}),$$
where $\underbar 0$ (resp. $\underbar 1$) stands for the (infinite) constant sequence with all entries $0$ (resp. $1$).

\noindent (iii) Show that $\phi$ is increasing (resp. continuous) with respect to the lexicographic order 
(resp. product topology) on $\{0,1\}^{\mathbb{N}}$.

\noindent (iv) Let $\phi_0 \!: \{0,1\}^{\mathbb{N}} \to [0,1]$ be defined by $\phi_0 (\xi) := \phi (0\xi)$. 
Also, recall the map $\phi_2 : \{0,1\}^{\mathbb{N}} \to [0,1]$ from Exercise \ref{F:ejer-conj} defined by 
$$\phi_2 ( \xi) = \sum_{j \geq 1} \frac{i_j}{2^{j}}, \qquad \mbox{ where } \xi = (i_1,i_2,\ldots), \esp i_j \in \{0,1\}.$$ 
Show that $\psi := \phi _0 \circ \phi_2^{-1} : [0,1] \to [0,1]$ is a well-defined homeomorphism that 
sends bijectively the dyadic numbers onto the rationals (in $[0,1]$). 
More accurately, show that $\psi$ sends each point $i/2^n$ to $p_i^n / q_i^n$, where 
\esp $0 = p_0^n / q_0^n < p_1^n / q_1^n  < \ldots p_i^n / q_i^n  < \ldots < p_{2^n}^n / q_{2^n}^n = 1$ 
\esp is the $\mathrm{n}^{th}$ step Farey sequence of rationals recursively defined by  
$p^k_0 := 0, q^k_0 = p^k_1 = q^k_1 := 1$ for all $k \geq 0$, and 
$$p^{k+1}_{2i+1} := p^k_i, \quad 
p_{2i}^{k+1} := p^k_{i-1} + p^k_{i}, \quad q^{k+1}_{2i+1} := q^k_i, \quad
q_{2i}^{k+1} := q^k_{i-1} + q^k_{i}.$$
\index{Farey sequence}

\noindent{\underbar{Hint.}} Although the claim above can be proven by induction, a dynamical argument proceeds as follows. 
Let $H_2 : [0,1] \to [0,1]$ be the map defined by $H_2 (t) := \{ 2t \}$ for $t < 1$ and $H(1)=1$. 
Also, let $H_0: [0,1] \to [0,1]$ be defined by 
$$H_0 (t) := \left \{ \begin{array}{l}
\frac{t}{1-t} \hspace{0.76cm} \mbox {if } 0 \leq t < \frac{1}{2},
\\ \\
\frac{2t-1}{t} \hspace{0.62cm} \mbox{if } \frac{1}{2} \leq t \leq 1. 
\end{array} \right.$$
Show that both $H_0$ and $H_1$ are conjugate to the shift 
$\sigma : \{0,1\}^{\mathbb{N}} \to \{0,1\}^{\mathbb{N}}$ defined by 
$\sigma (i_1,i_2,i_3,\ldots) := (i_2,i_3,\ldots)$. More precisely, show that 
$$H_0 \circ \phi_0 = \phi_0 \circ \sigma, \qquad H_2 \circ \phi_2 = \phi_2 \circ \sigma.$$

\noindent (v) For each $\xi \in \{0,1\}^{\mathbb{N}}$, let $\bar{\xi}$ be the {\em conjugate} of $\xi$, which results from 
$\xi$ by changing all 0's into 1's and vice versa. Check that $\phi (\xi) \phi(\bar{\xi}) = 1$ holds for all sequences $\xi$ 
that are not constant.

\noindent (vi) Let $\Phi: \{ 0,1 \}^{\mathbb N} \to \mathbb{R} \cap \{\infty\} =
\mathbb{P}^1 (\mathbb{R})$ be defined by 
$$\Phi (0 \xi) := -\phi (\bar{\xi}), \qquad \Phi (1 \xi) := \phi (\xi).$$
Translate the properties of $\phi$ above into analogous properties of $\Phi$.
\end{ejer}\end{small}

Recall from Exercise \ref{F:ejer-conj} the homeomorphisms $\hat{a}$ and $\hat{b}$ 
of $\{0,1\}^{\mathbb{N}}$ given by 
$$\hat{a} (\xi) := \left \{ \begin{array}{l}
0 \eta \hspace{0.8cm} \mbox {if } \xi = 0 0 \eta,
\\ \\
1 0 \eta \hspace{0.62cm} \mbox{if } \xi = 01 \eta, 
\\ \\
11 \eta \hspace{0.62cm} \mbox {if } \xi = 1 \eta,
\end{array} \right.
\qquad \mbox{ and } \qquad 
\hat{b} (\xi) := \left \{ \begin{array}{l}
\xi \hspace{1.2cm} \mbox {if } \xi = 0 \eta,
\\ \\
1 0 \eta \hspace{0.8cm} \mbox{if } \xi = 100 \eta, 
\\ \\
11 0 \eta \hspace{0.62cm} \mbox {if } \xi = 101 \eta,
\\ \\
1 1 1 \eta \hspace{0.62cm} \mbox{if } \xi = 11 \eta.
\end{array} \right.$$

\begin{small}\begin{ejer} \label{ejer:a}
Define $a,b$ in $P(\Gamma)$ by letting
$$a (t) := \tilde a (t) = t+1, 
\qquad \mbox{ and } \qquad 
b (t) := \left \{ \begin{array}{l}
t  \hspace{1.35cm} \mbox {if } t \leq 0,
\\ \\
\frac{t}{1-t} \hspace{0.95cm} \mbox{if } 0 \leq t \leq \frac12, 
\\ \\
3 - \frac{1}{t} \hspace{0.64cm} \mbox{if } \frac{1}{2} \leq t \leq 1,
\\ \\
t+1 \hspace{0.78cm} \mbox{if } t \geq 1. 
\end{array} \right.$$
Check that $\Phi (\hat{a} (\xi)) = a ( \Phi (\xi))$ and  
$\Phi (\hat{b} (\xi )) = b (\Phi (\xi))$ hold for all $\xi \!\in\! \{ 0, 1 \}^{\mathbb{N}}$.
\end{ejer}\end{small}

Observe that the element $\tilde c = \left( \begin{array}{cc} \sqrt{2} & 0\\  0 & 1/\sqrt{2}  \end{array}\right)$ above corresponds 
to multiplication by 2. It is hence crucial to encode the action of the map $x \to 2x$ in the coordinates given by $\phi$.

\begin{small}\begin{ejer}
Let $\hat c: \{0,1\}^{\mathbb{N}} \to \{0,1\}^{\mathbb{N}}$ be recursively defined by 
$$\hat c ( \xi ) := \left \{ \begin{array}{l}
0 \hat c (\eta) \hspace{1.6cm} \mbox {if } \xi = 0 0 \eta,
\\ \\
10 \hat{c}^{-1}(\eta) \hspace{1cm} \mbox{if } \xi = 0 1 \eta,
\\ \\
11 \hat c (\eta) \hspace{1.4cm} \mbox{if } \xi = 1 \eta. 
\end{array} \right.
$$
Show that, for all $\xi \in \{ 0 ,1 \}^{\mathbb{N}}$,
$$\phi (\hat c ( \xi )) = 2 \phi (\xi).$$

\noindent{\underline{Hint.}} Check the equality above for sequences $\xi$ that become eventually constant (to 
do this, use an induction argument on the length of the largest finite subword of $\xi$ before it becomes constant). 
Show also that $\hat c$ is a homeomorphism, and conclude by continuity that the equality above holds for all $\xi$. 
\end{ejer}
\end{small}

Similarly to Exercise \ref{F:ejer-extended}, given a finite binary sequence $s$, we let $\hat a_s$ (resp. $\hat c_s$) 
be the map that consists of the action of $\hat a$ (resp. $\hat c$) localized at the subtree starting at the terminal vertex of 
the path $s$. In precise terms, 
$$\hat a_s (\xi ) := \left \{ \begin{array}{l}
s \hat a (\eta) \hspace{0.9cm} \mbox {if } \xi = s \eta,
\\ \\
\xi \hspace{1.64cm} \mbox{otherwise}, 
\end{array} \right.
\qquad \quad
\hat c_s (\xi ) := \left \{ \begin{array}{l}
s \hat c (\eta ) \hspace{0.9cm} \mbox {if } \xi = s \eta,
\\ \\
\xi \hspace{1.64cm} \mbox{otherwise}. 
\end{array} \right.$$
Note that $\hat b = \hat c_1$. Note also that $\hat{c}_s$ is conjugate 
to $\hat c_{10}$ for each non-constant, nonempty, binary sequence $s$.

\begin{small}\begin{ejer} \label{ejer:c}
Let $c$ be the piecewise projective homeomorphism of the real line defined by 
$$ c ( t ) := \left \{ \begin{array}{l}
\frac{2t}{1+t} \hspace{1.2cm} \mbox {if } t \in [0,1],
\\ \\
t \hspace{1.64cm} \mbox{otherwise}. 
\end{array} \right.$$
Check that \esp $\Phi (\hat{c}_{10} (\xi)) = c (\Phi (\xi))$ \esp 
holds for all $\xi \!\in\! \{ 0,1 \}^{\mathbb{N}}$.
\end{ejer}\end{small}

We denote by $G_0$ the subgroup of $\mathrm{P} (\Gamma)$ generated by the elements $a,b,c$ from 
Exercises \ref{ejer:a} and \ref{ejer:c}. Then $G_0$ has no free subgroup, and we would like to understand the 
orbit equivalence relation of $G_0$. 

According to the Proposition \ref{p:relation}, the orbit equivalence relation of $\mathrm{P}(\Gamma)$ coincides a.e. with 
that of $\Gamma$, which is non-amenable by Theorem \ref{th:CG}. A priori, the orbit relation of $G_0$ is finer, meaning 
that elements could be related by $P (\Gamma)$ without being related by $G_0$. We claim that, however, the same orbit 
equivalence relation arises for this smaller group $G_0$, which is therefore non-amenable. The 
argument presented below is taken verbatim from \cite{lodha}.

\begin{small}\begin{ejer} Show that the equivalence relation of $G_0$ coincides with that of $\Gamma$ by 
following the steps below.

\noindent (i) Check that the following relations hold for certain restrictions of elements in $G_0$:
$$a c^{-1} a^{-1} c b \, (x) = 2x = \tilde{c} \, (x) \quad \mbox{ for all} \quad x \in [0,1],$$
$$a^{-3} b \, (x) = -1/x = \tilde{b} \, (x) \quad \mbox{ for all}  \quad x \in [1,2,1],$$
$$aba \, (x) = -1/x = \tilde{b} \, (x) \quad \mbox{ for all} \quad x \in [-1,-1/2].$$

\noindent (ii) Using the first equality above, show that if $x,y$ in $\mathbb{R} \setminus \mathbb{Q}$ are related by 
$\tilde{c}$ (say $y = 2 x$), then they are related by the group $G_0$.

\noindent{\underline{Hint.}} Denote the integer part of $x$ by $k$. 
Then $a^{-k} (x) = x-k$ belongs to $[0,1]$, and 
$$\tilde{c} (x) = 2 x = 2 (x-k) + 2 k = a^{2k} ( ac^{-1} a^{-1} c b \, (a^{-k} (x))).$$

\noindent (ii) Using the second  and third of the equalities above, show that if $x,y$ in $\mathbb{R} \setminus \mathbb{Q}$ 
are related by $\tilde{b}$, then they are related by the group $G_0$.

\noindent{\underline{Hint.}} Assuming that $x > 0$, there exists $n \in \mathbb{Z}$ such that $x' := 2^n x$ belongs to $[1/2,1]$. 
By the previous claim, $x$ and $x'$ are related in $G_0$, hence we can work with $x'$ instead of $x$. Now note that 
$$\tilde{b} (x') = -\frac{1}{x'} = a^{-3} b (x').$$
For negative $x$, proceed similarly by using the third equality above.
\end{ejer}\end{small}

The remaining task now is to give a finite presentation for $G_0$. To do this, we will deal with 
its isomorphic version acting on $\{0,1\}^{\mathbb{N}}$ generated by $\hat a$, $\hat b$ and $\hat c$.

We start by noticing that the subgroup of $G_0$ generated by $\hat a$ and $\hat b = \hat a_1$ is isomorphic to Thompson's 
group F (see Exercise \ref{F:ejer-conj}).

\index{Thompson's group!Thurston realization}
\index{Thompson's group!Ghys-Sergiescu realization}
\begin{small}
\begin{rem} 
The map \esp $\psi = \phi _0\circ \phi_2^{-1}$ \esp  from Exercise \ref{hurwitz} conjugates the standard dyadic action of F to an 
action by piecewise projective homeomorphisms of $[0,1]$. (The latter was first constructed by Thurston, and it is well described 
in \cite{CFP} and \cite{yo}.) It hence corresponds to the Ghys-Sergiescu conjugacy between these two actions; see \cite{GS}. It 
may be proved that $\psi$ is totally singular with respect to the Lebesgue measure; see for instance \cite{DKN-moscow,KS}. This 
function is known as {\em Conway's box function}, and denoted $\framebox{x} := \psi(x)$. Its inverse function was known to 
Minkowski, who denoted it by $? (x) := \psi^{-1}(x)$; it is called the {\em Minkowski question mark function} (and also the 
{\em slippery devil staircase}).
\end{rem}
\end{small}
\index{Conway's box function}
\index{Minkowski's question mark function}

Recall from Exercise \ref{F:ejer-extended} that for finite binary sequences $s,t$, we let $\hat{a}_t (s)$ be the action 
of $\hat a_t$ on $s$ whenever it is defined. The main result of Lodha and Moore \cite{lodha} may be stated as follows. 
\index{Group!Lodha-Moore}
\begin{thm} {\em Consider $G_0$ as a group generated by $\hat a_s$ and $\hat{c}_t$, where $s$ is any (perhaps empty) sequence and $t$ 
any nonempty, non-constant sequence. Then the family of relations below provide a presentation of $G_0$ with respect to these generators:}

\vsp

\noindent (i) {\em If $\hat{a}_t (s)$ is defined, then $\hat a_t \hat a_s = \hat a_{\hat a_t (s)} \hat a_t$;}

\vsp

\noindent (ii) {\em If $\hat{a}_t (s)$ is defined, then $\hat a_t \hat c_s = \hat c_{\hat a_t (s)} \hat a_t$;}

\vsp

\noindent (iii) {\em If $s,t$ are (nonempty, non-constant and) incompatible, then $\hat c_s \hat c_t = \hat c_t \hat c_s$;}

\vsp

\noindent (iv) {\em For each $s$ a nonempty, non-constant, binary sequence,  
$\hat c_s = \hat c_{s11} \hat c_{s10}^{-1} \hat c_{s0} \hat a_s$.}
\label{thm:pres-LM}
\end{thm}

By Exercise \ref{F:ejer-extended}, the first group of relations correspond to those of a presentation of F, 
and by Exercise \ref{F:fin-pres}, these may be summarized into 
$$[\hat a^{-1} \hat b, \hat a \hat b \hat a^{-1}] = [\hat a^{-1} \hat b, \hat a^2 \hat b \hat a^{-2}] = id,$$
where $\hat b := \hat a_1$. Using conjugacy by elements in F, it is not hard to see that the second 
group of relations may be reduced to 
$$[\hat a_{0} , \hat c_{10}] = [\hat a_{01}, \hat c_{10}] = [\hat a_{11}, \hat c_{10}] = [\hat a_{111}, \hat c_{10}] = id.$$
By a similar procedure, the third group reduces to 
$$[\hat c_{10}, \hat c_{01}] = [\hat c_{10}, \hat c_{001}] = id,$$
and the last group to the single non-commuting relation
$$\hat c_{10} = \hat c_{1011} \hat c_{1010}^{-1} \hat c_{100} \hat a_{10},$$
thus providing a finite presentation.

\begin{small}\begin{ejer}
By rewriting the relations above in terms of $a \sim \hat a$, $b \sim \hat a_1$ and $c \sim \hat c_{10}$, 
check that the list of relations above reduces to
$$[a^{-1}  b,  a  b  a^{-1}] = [ a^{-1}  b,  a^2  b  a^{-2}] = id,$$
$$[c, a^{-1} b^{-1} a^2] = [c, a^{-1} b^{-1} a b^{-1} a^{-1} b^2 a] = [c, aba^{-1}] = [c, a^2 b a^{-2}] = id,$$
$$[a^{-1} c a, c] = [a^{-2} c a^2] = id,$$
$$c = b^{-1} a b^{-1} a c a^{-1} c^{-1} b a^{-1} c a b^{-1} a^{-1} b^2,$$
the last one being a simplification of 
$$c = (b^{-1} a b^{-1} a c a^{-1} b a^{-1} b) ( b^{-1} a b^{-1} c^{-1} b a^{-1} b) (b^{-1} c b) (b^{-1} a b^{-1} a^{-1} b^2).$$
\end{ejer}\end{small}

The relations in Theorem \ref{thm:pres-LM} are easy to check. Those of type (i) and (ii) arise by conjugacy. 
Commutativity in case (iii) holds because $\hat c_s$ and $\hat c_t$ correspond to maps with disjoint supports 
whenever $s$ and $t$ are incompatible. Finally, 
relations of type (iv) reduce by conjugacy to \esp $\hat c_{10} = \hat c_{1011} \hat c_{1010}^{-1} \hat c_{100} \hat a_{10}$, 
\esp which is straightforward to check. (Note that, by the definition of $\hat c$, we have 
$\hat c = \hat c_{0} \hat c_{10}^{-1} \hat c_{11} \hat a$.) The goal now is to provide the main ideas of the proof of Theorem 
\ref{thm:pres-LM}. Although the proof itself will be left as an exercise, we will illustrate most of the crucial steps with examples. We refer the reader to \cite{lodha} in case of problems with the formal arguments.

We first consider words representing elements in $G_0$ made up of letters of the form $\hat a_s$ and $\hat c_t$, where $s$ 
is arbitrary (perhaps empty) and $t$ is 
nonempty and non-constant. We say that such a word is {\em standard} if, from right to left, it is the concatenation of a word 
on the $\hat a_s$ and a word on the $\hat c_t$, and whenever both $\hat c_{t_1}$ and $\hat c_{t_2}$ occur with nontrivial 
exponents for certain $t_2 \subsetneq t_1$, the first of these appears before the other. The {\em depth} of such a word 
is the smallest $\ell$ for which there is some $\hat c_t$ appearing in it satisfying $\mathrm{length} (t) = \ell$. (In case no 
$\hat c_t$ appears, the depth is defined to be infinite.) Two words are said to be equivalent if it is possible to 
derive one of them from the other one by applying the relations listed in Theorem \ref{thm:pres-LM}.

\begin{small}\begin{ejer} For most of the claims below, use an inductive procedure for the proof.

\vsp

\noindent (i) Prove that for every nonempty, non-constant, binary sequence $s$ and each $\ell \geq 1$, there are 
standard words $W_1, W_2$ equivalent to $\hat c_s$ and $\hat c_s^{-1}$, respectively, such that:

\noindent -- If $\hat a_{t}$ occurs in $W_i$, then $t$ entends $s$;

\noindent -- If $\hat c_{t}$ occurs in $W_i$, then $t$ extends $s$, has length $\geq \ell$, 
and the exponent of $\hat c_t$ is $\pm 1$;

\noindent -- If $\hat c_{t_1}$ and $\hat c_{t_2}$ occur in $W_i$ for $t_1 \neq t_2$, then $t_1$ and $t_2$ 
are incompatible.

\noindent{\underline{Hint.}} Use the relations of type (iv) above and their inverse versions, namely:
$$\hat c_s^{-1} 
= \hat a_s^{-1} \hat c_{s0}^{-1} \hat c_{s10} \hat c_{s11}^{-1} 
= \hat c_{s00}^{-1} \hat c_{s01} \hat c_{s1}^{-1} \hat a_s^{-1} .$$
For example, for $s = 10$ and $\ell = 5$, we have:
\begin{footnotesize}
\begin{eqnarray*}
\hat c_{10} 
\!\!\! &=& \!\!\!
\hat c_{1011} \hat c_{1010}^{-1} \hat c_{100} \hat a_{10} \\
&=& \!\!\!
( \hat c_{101111} \hat c_{101110}^{-1} \hat c_{10110} \hat a_{1011} ) 
( \hat c_{101000}^{-1} \hat c_{101001} \hat c_{10101}^{-1} \hat a_{1010}^{-1} )
( \hat c_{10011} \hat c_{10010}^{-1} \hat c_{1000} \hat a_{100} ) \hat a_{10} \\
&=& \!\!\!
\hat c_{101111} \hat c_{101110}^{-1} \hat c_{10110}  
\hat c_{101000}^{-1} \hat c_{101001} \hat c_{10101}^{-1}  
\hat c_{10011} \hat c_{10010}^{-1} \hat c_{1000} 
\hat a_{1011} \hat a_{1010}^{-1} \hat a_{100}  \hat a_{10} \\
&=& \!\!\!
\hat c_{101111} \hat c_{101110}^{-1} \hat c_{10110}  
\hat c_{101000}^{-1} \hat c_{101001} \hat c_{10101}^{-1}  
\hat c_{10011} \hat c_{10010}^{-1} \hat c_{100011} \hat c_{100010}^{-1} \hat c_{10000} \hat a_{1000}
\hat a_{1011} \hat a_{1010}^{-1} \hat a_{100}  \hat a_{10}.
\end{eqnarray*}
\end{footnotesize}

\vspace{-0.35cm}

\noindent (ii) Given a word $W$ in the $\hat a_s$, show that there exists $\ell \geq 1$ such 
that for any standard word $W'$ of depth $\geq \ell$, the word $W W'$ is equivalent 
to a standard word of depth $\geq \ell - k$, where $k$ is the word-length of $W'$.

\noindent{\underline{Hint.}} Use relations of type (ii).

\vsp

\noindent (iii) Using (i) and (ii), prove that for every $\ell \geq 1$, each word in the 
$\hat a_s, \hat c_t$  (with $t$ nonempty and non-constant) is equivalent to a standard word.
\end{ejer}
\end{small}

If $W$ is a standard word and $\hat c_s$ occurs in $W$, we say that $s$ is {\em exposed} if there is an infinite 
binary path starting at $s$ along which no $y_t$ occurs in $W$. The word is {\em sufficiently expanded} if,  
whenever $\hat c_s$ occurs in $W$ and $s$ is not exposed, we have that:

\noindent -- $\hat c_{s0}$ occurs in $W$ if $\hat c_s$ occurs with a positive exponent;

\noindent -- $\hat c_{s1}$ occurs in $W$ if $\hat c_s$ occurs with a negative exponent.

\begin{small}\begin{ejer}
Show that every standard word is equivalent to one which is sufficiently expanded.

\noindent{\underline{Hint.}} Use the relations of type (iv) and their inverse versions. For instance, 
this yields
$$\hat c_{0}^3 \hat c_{001} \hat a_1 
= \hat c_{0}^2 \hat c_{00} \hat c_{010}^{-1} \hat c_{011} \hat a_0 \hat c_{001} \hat a_1
= \hat c_{0}^2 \hat c_{00} \hat c_{010}^{-1} \hat c_{011} \hat c_{010} \hat a_0 \hat a_1
= \hat c_{0}^2 \hat c_{00} \hat c_{011} \hat a_0 \hat a_,$$
where the left-side expression is not sufficiently expanded yet the right-side one is.
\end{ejer}\end{small}

For a word $\Lambda$ in the alphabet $\{0,1,\hat c, \hat c^{-1}\}$, we define the process of 
{\em advancing} the occurrence of a certain $\hat c^{\pm 1}$ the result of replacing 
$$\hat c 00 \to 0 \hat c, \qquad \hat c 01 \to 10 \hat c^{-1}, \qquad \hat c 1 \to 11 \hat c,$$
$$\hat c^{-1} 0 \to 00 \hat c^{-1}, \quad \hat c^{-1} 10 \to 01 \hat c, \quad \hat c^{-1} 11 \to 1 \hat c^{-1}.$$
The resulting word $\Lambda'$ after advancing this occurrence of $\hat c^{\pm 1}$ several times is said 
to be an advanced version of $\Lambda$. A potencial cancellation in $\Lambda$ is a concatenation 
of the form $\hat c \hat c^{-1}$ or $\hat c^{-1} \hat c$ obtained in an advanced version of $\Lambda$.

\begin{small} 
\begin{ejer}\label{ejer:advancing}
Let $\Lambda$ be a word in the alphabet $\{0,1,\hat c, \hat c^{-1}\}$ having no potential cancellation. 

\noindent (i) Prove that no advanced version of $\Lambda$ contains potential cancellations.

\noindent (ii) Prove that there exist binary, finite sequences $s,t$ such that $\Lambda s$ can be transformed 
into $t\hat c^n$ by advancing all occurrences of $\hat c^{\pm 1}$ ($n$ is the number of these occurrences).
\end{ejer}

\begin{ejer}\label{ejer:ultimito}
Let $W$ be a nontrivial, sufficiently expanded word on $\hat c_s$, where $s$ ranges over nonempty, non-constant, 
binary sequences. The goal is to define finite binary sequences $u,v$ such that for all $\xi \!\in\!  \{ 0 , 1\}^{\mathbb{N}}$, 
the image of $u \xi$ under $W$ is $v \hat c^n (\xi)$ for some $n > 0$.

\noindent (i) Assume the existence of $u, v$ as above. Show that $W$ sends \, $u 0^{2^n} 1 0^{2^n} 1 \ldots$ \, 
to \, $v 0 1^{2^n} 0 1^{2^n} \ldots$ \, As these two points have nonequivalent tails, conclude that $W$ does 
not represent the identity.

\noindent (ii) Let $y_{t_0}$ be the last entry (from right to left) in $W$. In case $t_0$ is exposed, let $u$ be any 
path extending $t_0$ that points to infinite without passing through the $t_i$'s for which $y_{t_i}$ appears 
in $W$. If $t_0$ is not exposed, extend it as follows: if $y_{t_0}$ appears with a positive exponent, add a 0 
to the right, and a 1 otherwise; continue the procedure until an exposed index is obtained, and conclude the 
construction of a path, denoted by $u_0$. Let $\Lambda$ be the word on $\{0,1,\hat c, \hat c^{-1}\}$ obtained by 
inserting $\hat c^n$ just after $s$ whenever $\hat c^n_s$ appears in $W$. By Exercise \ref{ejer:advancing}, there 
exist $s,t$ such that $\Lambda s$ can be transformed into $t \hat c^n$ by the procedure of advancing. Show that 
the claim holds for $u := u_0 s$ and $v := t$.

\noindent\underline{Hint.} First note that $\Lambda$ contains no potential cancellation, and then check that 
applying $W$ to $u \xi$ corresponds to advancing the word $\Lambda u$. 
\end{ejer}
\end{small}

The conclusion of the proof of Theorem \ref{thm:pres-LM} is now at hand. Namely, every word in $\hat a_s, \hat c_t$ 
can be transformed to a sufficiently expanded standard one using the relations of type (i), (ii), (iii), (iv) above. If we 
assume that such a word represents the identity, then its part in the generators $\hat c_t$ must be the trivial word, 
otherwise by Exercise \ref{ejer:ultimito} there would be a point being sent into another one with a nonequivalent 
tail, which is impossible. As a consequence, $W$ is a word in the generators $\hat a_s$. Now, every such word 
representing the identity can be transformed into the trivial word using the relations of type (i), as these give 
a presentation of $\mathrm{F} \sim \langle \hat a_s \rangle$. 


\section{Almost-Periodicity}
\label{section: almost periodicity}

\hspace{0.45cm} In this section, we develop the notion of almost-periodicity for group actions on the real
line. A left-orderable group being given, the set of such actions equipped with the compact-open topology
can be used as a substitute for the space of left-orders. As an example, we will show how this yields an
alternative proof of Theorem \ref{witte-bien} which does not rely on the theory of Conradian orders. 
The concept of almost-periodic actions has also been used recently to provide new constructions of 
finitely-generated, left-orderable, simple groups (see \S \ref{s: finitely generated simple left orderable}) 
and as a key tool for the proof of the non left-orderability of irreducible lattices in semi-simple real algebraic group \cite{DH}. 

\subsection{Almost-periodic actions}
\label{s: almost periodic}

\hspace{0.45cm} The group of orientation-preserving homeomorphisms of the real line
is equipped with the compact-open topology, which makes it a topological group. The {\bf\textit{translation flow}} on $\mathrm{Homeo}_+ (\mathbb R)$ is defined by 
\[ (s, h ) \in \mathbb R \times \mathrm{Homeo}_+ (\mathbb R) \mapsto \tau_s ^{-1} \circ h \circ \tau_s \in \mathrm{Homeo}_+ (\mathbb R),\]
where $\tau_s (t) := s+t$.  An element $h \in  \mathrm{Homeo}_+ (\mathbb R)$
is said to be an {\bf\textit{almost-periodic homeomorphism}} if its orbit by the translation flow 
is relatively compact in $\mathrm{Homeo}_+({\mathbb R})$. The set of
almost-periodic, orientation-preserving homeomorphisms of the line is denoted $APH_+ ({\mathbb R})$.
\index{Action!almost-periodic}

\begin{small}\begin{ex}\label{flujo-toral}
Certainly, every homeomorphism that is {\em periodic} ({\em i.e.}, it commutes with a nontrivial translation)
is almost-periodic. An example of an almost-periodic homeomorphism that is non-periodic is given by
\begin{equation} \label{eq: example}
\varphi ( t ) := t + \frac{1}{3} \big( \sin (t) + \sin (t \sqrt{2} ) \big).
\end{equation}
To see this, consider the continuous map \( (x,y) \mapsto \varphi_{x,y} \) from the torus \( (\mathbb R/ \mathbb Z)^2 \) to \( \mathrm{Homeo}_+({\mathbb R})\) defined by 
\[  \varphi _{x,y} (t) = t +\frac{1}{3} \big( \sin ( x+t)+\sin ( y + t \sqrt{2} ) \big) \]
Let  $T= \{T_s\}_{s \in \mathbb{R}}$ be the irrational flow on the torus defined by 
$$ T_s (x,y):= (x+s, y+s\sqrt{2} ). $$ 
The map \(\varphi \) is equivariant with respect to the actions of the flow \(T\) on the torus and  of the translation flow on \( \mathrm{Homeo}_+({\mathbb R})\). More precisely, for every \( s\in \mathbb R\), we have  
\[ \varphi _{ T_s (x,y) } = \tau_ s ^{-1} \circ \varphi_{x,y} \circ \tau_s   . \]
In particular, the image set $\{ \varphi_{x,y} (\cdot) \!: (x,y) \in ( \mathbb{R} / \mathbb{Z} )^2\}$ is a compact subset of \( \mathrm{Homeo}_+({\mathbb R})\) which is invariant under the translation flow, hence \(\varphi = \varphi _{0,0} \) is almost-periodic.

\end{ex}\end{small}

\vsp

\begin{lem}
{\em The subset $APH_+ (\mathbb R)$  is a subgroup of $\mathrm{Homeo}_+ (\mathbb R)$.}
\end{lem}

\noindent{\bf Proof.}
This is a consequence of the continuity of the composition and inverse operations on
$\mathrm{Homeo}_+ ({\mathbb R})$ with respect to the compact-open topology.
$\hfill\square$

\vspace{0.35cm}

\index{Action!almost-periodic}
A group action on the real line whose image is contained in $APH_+(\mathbb R)$ will be said to be an  
{\bf \textit{almost-periodic action}}. 
There are several ways to construct faithful almost-periodic actions of a given left-orderable, countable group $\Gamma$
on the line. The simplest one consists in considering a faithful action on the interval and then to extend it to the
whole line so that it commutes with the translation $t\mapsto t+ 1$. This somewhat trivial construction shows
that $APH_+(\mathbb R)$ contains a copy of every left-orderable, countable group. Nevertheless, 
in order to carry out a study that also involves actions that appear as limits of conjugates
of a given one, we are forced to consider actions that may be unfaithful.
This is the reason why we use the notation $\Phi \!: \Gamma \to \mathrm{Homeo}_+(\mathbb{R})$ in what follows.

Starting with an almost-periodic action of a group $\Gamma$ on the real line, we next provide
a compact one-dimensional foliated space together with a $\Gamma$-action on it that preserves
the leaves. It is this construction that lends interest to considering almost-periodic actions.

\vsp

\begin{prop} \label{l:closure under translation flow}
\textit{Let $\Gamma$ be a finitely-generated group and $\Phi_0 \!: \Gamma \rightarrow \mathrm{Homeo}_+({\mathbb R})$ 
an action by orientation-preserving homeomorphisms of the line. Then $\Phi_0$ is almost-periodic if and only if there
exists a topological
flow $ T =
\{ T_s \}_{s\in \mathbb R}$ acting freely on a compact space $X$, an action of $\Gamma$ on $X$
by homeomorphisms preserving every $T$-orbit together with its orientation, and a point $x_0 \!\in\! X$,
such that for every $g\in \Gamma$ and every $t\in {\mathbb R}$,
\begin{equation}\label{eq: G-action}
g \big( T_t (x_0) \big) = T_ {\Phi_0(g)(t)} (x_0).
\end{equation}
Moreover, the flow can be taken so that the $T$-orbit of $x_0$ is dense in $X$.}
\end{prop}

\noindent{\bf Proof.} Let us first show that if there is a compact space
$X$ together with a flow $T$ and a $\Gamma$-action verifying~\eqref{eq: G-action},
then the action $\Phi_0$ is almost-periodic. Indeed, for each $x \! \in \! X$,
we can lift the $\Gamma$-action on $X$ to an action
$\Phi^x \!: \Gamma \rightarrow \mathrm{Homeo}_+ (\mathbb R)$ verifying
\[  g \big( T_t (x) \big) = T_{\Phi^x (g)(t)} (x) ,\]
which is well-defined since the flow $T$ is free.
Moreover, as the $\Gamma$-action on $X$ is by homeomorphisms, for every $g\in \Gamma$,
the map $x \mapsto \Phi^x (g)$ from $X$ into $\mathrm{Homeo}_+(\mathbb R)$
is continuous. Hence, the set of elements $\Phi^x (g)$, where $x \!\in\! X$, is
compact. Now, for every $s, t$ in $\mathbb R$ and every $x\in X$, we have
\[ g\big( T_t( T_s (x) ) \big)
= g\big( T_{t+s} (x) \big)
= T_{\Phi^x(g)(t+s)}(x)
=  T_{\Phi^x(g)(t+s)-s}(T_s(x)), \]
which yields
\[ \Phi^{T_s(x)} (g) = \tau_{-s} \circ \Phi^x(g)  \circ \tau_s .  \]
Therefore, for every $g\in \Gamma$, the conjugates of $\Phi^x (g)$ by the translations
$\tau_s$ stay in a compact set, which proves that $\Phi^x$ is almost-periodic for every
$x\in X$. In particular, for $x=x_0$, we deduce that $\Phi_0 = \Phi^{x_0}$ is almost-periodic.

Conversely, let us start with an almost-periodic action $\Phi_0$, and let us provide the compact space
$X$ together with the flow $T= \{T_s\}_{s\in \mathbb{R}}$ and the $\Gamma$-action verifying~\eqref{eq: G-action}.
Denote by $APA_+(\Gamma)$ the set of almost-periodic actions of $\Gamma$ on the real line. This can
be seen as a closed subset of $APH_+ (\mathbb R)^{\mathcal{G}}$, where $\mathcal{G}$ is a finite
generating set of $\Gamma$. Consider the {\bf \textit{translation flow}}
$T=\{T_s\}_{s\in \mathbb{R}}$ of conjugacies by translations on $APA_+ (\Gamma)$, namely
\begin{equation}\label{eq: translation flow}
T_s (\Phi) (g) : = \tau_{-s} \circ \Phi (g) \circ \tau_s,
\end{equation}
where $\Phi \in APA_+(\Gamma) $ and $g \in \Gamma$. This is a topological flow acting on
$APA_+(\Gamma)$. Denote by $X$ the closure of the $T$-orbit of $\Phi_0$. This is a
compact $T$-invariant subset of $APA_+(\Gamma)$, since $\Phi_0$ is almost-periodic.

We claim that the formula
\begin{equation} \label{eq:def_action}
g (\Phi) : = \tau_{-\Phi(g) (0)} \circ \Phi \circ \tau_{\Phi(g)(0)} \end{equation}
defines an action of $\Gamma$ on $APA_+(\Gamma)$.
One can verify this by a tedious computation, but a more conceptual argument proceeds as follows.
Consider the actions of $\mathbb R$ and $\Gamma$ on the product space $APA_+ (\Gamma) \times \mathbb R$
given by
\[ s (\Phi, t)  := \big( T_s (\Phi) , t-s \big) \ \ \
\ \text{and} \ \ \ \ g  (\Phi , t  ) = \big( \Phi, \Phi(g)(t) \big) . \]
A point in $APA_+(\Gamma) \times \mathbb R$ can be thought of as an almost-periodic action of $\Gamma$
together with a marker. The action of \(\mathbb R\) on $APA_+(\Gamma) \times \mathbb{R}$ corresponds to
translating the marker while conjugating the almost-periodic action by the same translation. The
$\Gamma$-action on $APA_+(\Gamma) \times \mathbb R$ corresponds to acting on the marker
using the action of the first coordinate while leaving the almost-periodic representation unchanged.
An easy computation shows that these two actions commute. Hence, there is a natural action of $\Gamma$
on the quotient of $ APA_+ (\Gamma) \times \mathbb R$ by $\mathbb R$, which naturally identifies with
$APA_+ (\Gamma)$ via the embedding
$$\Phi \in APA_+ (\Gamma) \mapsto (\Phi , 0) \in APA_+ (\Gamma) \times \mathbb R.$$
The action of $\Gamma$ on $APA_+ (\Gamma)$ induced by
this identification is given by the formula~\eqref{eq:def_action}.

{\em A priori}, there is no reason to expect that the flow $T$ on $X$ acts freely. However, it is possible
to change it into a free one by the following procedure: Let $Y$ be any compact space endowed with a
topological flow $S$ acting freely. (For instance, the toral flow of Example \ref{flujo-toral};
see also Exercise \ref{delones} below.) Consider the space $\tilde{X} := X \times Y$ together with
the $\Gamma$-action on it defined as
$$g : ( \Phi, y ) \mapsto \big( g(\Phi), S_{\Phi(g)(0)} (y) \big)$$
and the (diagonal) flow
$$s : ( \Phi, y ) \mapsto \big( T_s (\Phi), S_s (y) \big).$$
Then all the properties above are still satisfied, and moreover the new flow on $\tilde{X}$ is 
free. Hence, we can (and we will) assume that the flow $T$ on $X$ is free.

Equation (\ref{eq: G-action}) is obvious from the construction, as well as the fact that
$X$ can be taken to be the closure of a single point. This closes the proof.
$\hfill\square$

\vspace{0.35cm}

In the context above, we will call an {\bf \em almost-periodic space} for $\Gamma$ a compact 
space $X$ together with a flow $ T = \{ T_s \}_{s\in \mathbb R}$ and a $\Gamma$-action 
by homeomorphisms preserving every $T$-orbit together with its orientation, such that 
$$g \big( T_t (x) \big) = T_ {\Phi_0(g)(t)} (x)$$
holds for every $g\in \Gamma$, $t\in {\mathbb R}$ and $x \in X$. 
Note that we do not impose the condition of freeness for the flow $T$ here, 
and we also relax the hypothesis of denseness of some orbit under $T$. 
We will come back to this point in \S \ref{s:free-aps}.

\vsp

\begin{small}\begin{ejer}
Show that in the case of the almost-periodic homeomorphism \eqref{eq: example} 
of Example \ref{flujo-toral}, the flow \((X, S) \) of Proposition \ref{l:closure under 
translation flow} can be taken to be an irrational flow on the torus. 
\end{ejer}\end{small}

\vsp 

\begin{small}\begin{ejer} \label{delones}
A {\bf {\em Delone set}} $D$ in $\mathbb{R}$ is a subset that is discrete and almost dense in a uniform way. More 
concretely, there exist positive constants $\varepsilon, \delta$ such that $|x-y| \geq \varepsilon$ for all $x \neq y$ in $D$,
and for all $z \in \mathbb{R}$ there is $x \in D$ such that $|x-y| \leq \delta$. Two Delone sets are {\em close} if they
coincide over a very large interval centered at the origin (this induces a topology that is metrizable). A Delone
set is repetitive if for all $r > 0$ there is $R = R(r) > 0$ such that for every pair of intervals $I,J$ of length
$r,R$, respectively, there is a translated copy of $D \cap I$ contained in $D \cap J$.

Assume that $D_0$ is a repetitive yet non-periodic Delone set in $\mathbb{R}$
(it is easy to build such sets). Show that the natural translation flow $T_t \!: D \mapsto D + t$,
restricted to the closure of the orbit of $D_0$, is a minimal flow.
\end{ejer}\end{small}

\begin{small} \begin{ejer} \label{ex: flow box} 
Given a flow \( T =\{T_t\}_{t\in \mathbb {R}} \)  acting continuously on a topological space \(X\), a {\bf \em flow box} is a local homeomorphism 
\( h : U \rightarrow I \times S \), where \(U\subset X\) is an open subset, \( I\subset \mathbb{R}\) is an open interval, and \(S\) is a topological space, such that \( h ( T_t ( x) ) =  \Psi _t ( h(x))\) holds whenever $h$ is defined at both $x$ and $T_t(x)$, where \(\Psi\) is the local flow acting on \(I\times S\) by the formula \( \Psi_t  ( u , s)  = (t+u, s) \).   Prove that, if continuous functions on \(X\) separate points, then every point \( x\in X\) belongs to the domain of a flow box. A {\bf \em plaque} of a flow box is a set of the form \( h ^{-1} (I\times \{s\}) \), where \( s\in S\). (Note that this is a connected subset of an orbit of the flow.) Prove that any connected compact interval in an orbit of the flow is contained in a plaque of a flow box.

\noindent{\underline{Hint.}} Let \( C^1 _T (X)\) be the space of continuous fonctions \(f: X\rightarrow \mathbb{R}\) such that \( \frac{d f\circ T_t ( \cdot) }{dt} \) exists and is continuous. Observe that, by using a kernel on the real line, the space \( C^1 _T (X) \) is dense in \( C^0 (X) \), and that for every point \(x\in X\), there exists a function of \( C^1 _T (X) \) with derivative \( \frac{d f\circ T _t (\cdot) }{dt}(x)  \neq 0\). Use such a function to prove the existence of a flow box whose domain contains \(x\). See \cite{Whitney} or \cite[Lemma 5.1]{DH} for more details. \end{ejer} \end{small}



\subsection{A bi-Lipschitz conjugacy theorem}
\label{s: Lipschitz conjugation}

\hspace{0.45cm} We denote by $\mathrm{BiLip}_+ (\mathbb R)$ the group of orientation-preserving,
bi-Lipschitz homeomorphisms of the real line. For every $h\in \mathrm{BiLip}_+ (\mathbb R)$, we let
$L(h)$ be its {\bf {\em bi-Lipschitz constant}}, that is, the minimum of the numbers $L \geq 1$ such that
\begin{equation} \label{eq: lipschitz}
L^{-1} |y-x| \leq |h(y) -h(x) | \leq L  |y-x| \ \ \ \mbox{ for all } x, y \mbox{ in } {\mathbb R}.
\end{equation}
We equip $\mathrm{BiLip}_+ (\mathbb R)$ with the topology of uniform convergence on compact sets.

\vsp

\begin{thm} \label{t: lipschitz conjugation}
\textit{Every finitely-generated group of homeomorphisms of the real line is topologically
conjugate to a group of bi-Lipschitz homeomorphisms.}
\end{thm}

\noindent{\bf Proof.}
Let $\nu = \lambda(x) dx$ be a probability measure on ${\mathbb R}$ with a smooth, positive density
$\lambda$ such that for $|x|$ large enough, we have $\lambda (x) = 1/x^2$. The following observation
will be central in what follows: If for some constant $L \geq 1$, a homeomorphism $h$ of the real line
satisfies
\begin{equation}
\label{eq: cauchy condition} h_*(\nu) \leq L \nu \quad \text{ and } \quad (h^{-1})_* (\nu) \leq L \nu,
\end{equation}
then $h$ is Lipschitz. To prove this fact, first note that $\nu \big( [x, +\infty) \big) = 1/x$ for
all sufficiently large positive numbers $x$ (and similarly, $\nu \big( (-\infty, x] \big) = 1 / |x|$ if $|x|$
is large enough and $x$ is negative). Thus, the left-side inequality in \eqref{eq: cauchy condition}
shows that $|h(x)| \leq L |x|$ holds for $|x|$ large enough. Since the density of $(h^{-1})_* (\nu)$ is 
given by $Dh (x) \lambda ( h(x) ) $, the right-side inequality in \eqref{eq: cauchy condition} yields
$Dh (x) \leq  L \lambda (x) \big/ \lambda (h(x)) $ for almost every $x$. Thus, up to sets of zero
Lebesgue measure, the derivative $Dh$ is bounded on every compact interval, and for $|x|$
large enough, we have
$$Dh (x) \leq \frac{ L \lambda (x) }{ \lambda (h(x)) }
= \frac{L |h(x)|^2}{|x|^2} \leq L^3.$$
This proves that $Dh$ is a.e. bounded, hence $h$ is Lipschitz,
with Lipschitz constant at most $L^3$.

Next, let $\Gamma$ be a finitely-generated subgroup of $\mathrm{Homeo}_+ ({\mathbb R})$, and let
$\mathcal{G}$ be a finite, symmetric system of generators of $\Gamma$. Let $\phi \in \mathcal{L}^1 (\Gamma)$
be a function taking positive values such that, for every $h \!\in\! \mathcal{G}$, there is a constant $L_h$ satisfying 
$\phi (hg) \leq L_h  \phi(g)$ for all $g \in \Gamma$. For instance, one can take $\phi (g) = \kappa^{||g||}$, where
$\| g \|$ is the word-length of $g$ with respect to $\mathcal{G}$
and $\kappa$ a sufficiently small positive number. (For
$\kappa < 1/|\mathcal{G}|$, one can ensure that $\phi$ belongs to $\mathcal{L}^1(\Gamma)$;
see also Exercise \ref{close to 1} below.) Let us normalize the function $\phi$ so
that $\sum_{g\in \Gamma} \phi(g) = 1$,
and let us define the probability measure $\nu_0$ on ${\mathbb R}$ by letting
\[  \nu_0 : = \sum_{g\in \Gamma} \phi (g) \hspace{0.05cm} g_* (\nu).\]
Note that for each $h\in \Gamma$, we have
\[  h_* (\nu_0) = \sum _{g\in \Gamma} \phi (g) \hspace{0.05cm} (hg)_* (\nu) \leq L_{h^{-1}} \nu_0.\]

 The measure $\nu_0$ has no atoms; indeed since \(\nu\) has a positive density with respect to the Lebesgue measure, we have \( \nu_0 (\{x\}) = \sum \phi(g) \nu (\{g^{-1} (x)\}) = 0\) for every \(x\in \mathbb R\). Also, the measure \(\nu_0\) has full support, since \(\nu_0\geq \phi(e)\nu\) and the support of \(\nu\) is full.  Lastly, the total mass \(\nu_0 (\mathbb R)\) of  \(\nu_0\) is the same as that of \(\nu\):
\[ \nu _ 0 (\mathbb R) = \sum _{g\in G} \phi(g) \nu (g^{-1}(\mathbb R)) =\nu (\mathbb R).\] 
Thus, there exists a homeomorphism $\varphi$ of the real line
sending $\nu_0$ into $\nu$. This homeomorphism is explicitly given by \( \varphi = \varphi_{\nu}^{-1}\circ \varphi_{\nu_0}\) where \( \varphi_{\nu_0}, \varphi_\nu : \mathbb R \rightarrow (0, \nu(\mathbb R) )\) are the homeomorphisms defined for every \(x\in \mathbb R\) by
\[ \varphi_{\nu_0} (x) := \nu_0 ((-\infty, x) ) \quad \text{ and } \quad \varphi_\nu (x) := \nu ((-\infty, x) ) . \]
These homeomorphisms satisfy \( (\varphi_{\nu_0})_* dx = \nu_0 \) and \( (\varphi_{\nu})_* dx= \nu \), where \(dx\) is understood as the Lebesgue measure on \(\mathbb R\).
For each $h \in \Gamma$, we have
\[  (\varphi \circ h \circ \varphi^{-1})_* (\nu)
= \varphi_* h_* (\nu_0) \leq L_{h^{-1}} \hspace{0.05cm} \varphi_* (\nu_0) = L _{h^{-1}}\hspace{0.05cm} \nu . \]
From the discussion at the beginning of the proof, we deduce that the conjugate of $\Gamma$
by $\varphi$ is contained in $\mathrm{BiLip}_+(\mathbb R)$.
$\hfill\square$

\vspace{0.3cm}

We should stress that there is no analogue of the previous theorem in higher dimension,
even for actions of (infinite) cyclic groups; see \cite{Yano}. 

\vsp

\begin{small}\begin{ejer}
Let $D$ be a Delone subset of $\mathbb{R}$ (see Exercise \ref{delones}). Show that there exists
a bi-Lipschitz homeomorphism of the real line that sends $D$ onto $\mathbb{Z}$. (Again, there is no
higher-dimensional analogue of this fact; see \cite{burago-kleiner, mcmullen} as well as
\cite{cortez-navas}.)
\end{ejer}\end{small}

The proof of Theorem \ref{t: lipschitz conjugation} above was taken from \cite{deroin}. In \S \ref{s:random walks}, 
we will give a more conceptual
(yet quite elaborate) proof based on probabilistic arguments (see also Exercise \ref{ejer-asi-tambien} below). 
For analogous results for transverse pseudo-groups of codimension-one foliations or groups acting 
on the circle, see~\cite[Proposition 2.5]{Deroin} and~\cite[Th\'eor\`eme D]{DKN}. 

\vsp

It is worth pointing out that Theorem \ref{t: lipschitz conjugation} holds more generally for {\em countable} groups of homeomorphisms; see 
Exercise \ref{ejer contable es lipschitz}. However, in general, it is not possible to conjugate an arbitrary group of homeomorphisms of the 
real line to a group of bi-Lipschitz transformations. For instance, this is obviously impossible for the whole group $\text{Homeo}_+(\R)$. 
The next exercise (built from a clever remark of Calegari taken from  \cite{cale}) shows that  the Abelian group $\Z^\Z$ consisting of 
maps $\{\alpha :\Z\to\Z\} $ endowed with the pointwise addition, admits an action by homeomorphisms of the real line that cannot be 
conjugated to a group of bi-Lipschitz transformations.

{\small
\begin{ejer} For each $n\in\Z$, let $g_n:[n,n+1)\to [n,n+1)$  be a homeomorphism that acts freely on $(n,n+1)$. Consider the 
(faithful) representation $\Phi:\Z^\Z\to \mathrm{Homeo}_+ ({\mathbb R})$ defined by $\Phi(\alpha)(x)=g_n^{\alpha(n)}(x)$ whenever $x\in[n,n+1)$. 
Show that $\Phi(\Z^\Z)$ is a group that cannot be conjugated so that it fits inside the group of bi-Lipschitz homeomorphisms of the line.
\end{ejer}
}

\begin{small}
\begin{ejer} \label{ejer contable es lipschitz}
In~\cite[Th\'eor\`eme D]{DKN}, it is proved that every countable group of circle homeomorphisms is topologically conjugate to a 
group of bi-Lipschitz maps. Use this fact to extend Theorem \ref{t: lipschitz conjugation} to countable groups of homeomorphisms 
of the line.

\noindent{\underline{Hint.}} Compactify the real line as the projective line \( \mathbb{P} ^1(\mathbb R) = \mathbb R \cup \{\infty\}\), and use 
the fact that a Lebesgue measure on it is of the form \(\lambda (x) dx\), with \(\lambda (x) \sim \frac{1}{x^2} \) as \( x\) tends to infinity. 
Then use the ideas of the proof of Theorem \ref{t: lipschitz conjugation}.
\end{ejer}

\begin{ejer}\label{ejer-asi-tambien} 
Using \cite[Th\'eor\`eme D]{DKN} stated above, show that every countable group of homeomorphisms of the interval is conjugate to a 
group of bi-Lipschitz homeomorphisms. Then conclude that Theorem \ref{t: lipschitz conjugation} holds for countable groups using 
Exercise \ref{ejer conjug} below.
\end{ejer}

\begin{ejer}\label{ejer conjug}
Let $\Gamma$ be a subgroup of $\mathrm{BiLip}_+([0,1])$, and let $\varphi \!: [0,1] \to
\mathbb{R}$ be an orientation-preserving homeomorphism such that $\varphi(x) = -1/x$
for $x$ close to zero, and $\varphi(x) = 1/(1-x)$ for $x$ close to $1$. Check that the
conjugate of $\Gamma$ by $\varphi$ is a subgroup of
$\mathrm{BiLip}_+(\mathbb{R})$.
\end{ejer}
\end{small}

\vspace{0.2cm}

\noindent {\bf On certain actions of the Baumslag-Solitar group $BS(1,2)$.}  
There are many actions for which Theorem \ref{t: lipschitz conjugation} is counterintuitive, as the Lipschitz constant seems to naturally explode 
close to $\pm \infty$ for certain group elements. This is the case, for instance, of the action of $BS(1,2)=\langle a,b \mid aba^{-1}=b^2\rangle $ on 
the line, where $a$ acts as a translation by $1$ and $b$ is a homeomorphism that fixes each $n\in \Z$;  
see the left picture in Figure 19 below. (Note that this action already arose in Example \ref{ex-bi-order-BS}.) Indeed, the action of $b$ on each 
interval $[n-1,n]$ is (a conjugate of) the square of the action of $b$ on $[n,n+1]$, which seems to produce an explosion of its Lipschitz constant. 
However, below we give an explicit realization of this action by bi-Lipschitz diffeomorphisms.

\begin{small}\begin{ex}\label{ex-action-bil} 
Let $f,g$ be the maps $x \mapsto 2 x$ and $x \mapsto x+1$ viewed as the diffeomorphisms of the projective line. Denote by $p$ 
the common fixed point of $f$ and $g$ (that is, the infinity in the standard projective chart). Let $\tilde{f}$ and $\tilde{g}$ be the lifts in 
$\widetilde{\mathrm{Homeo}}_+(\mathbb{R})$ of $f$ and $g$, respectively, both fixing some lift (hence all the lifts) of $p$. 
Finally, let $a,b$ be the real-analytic diffeomorphisms of the real line defined by $a := \tilde{f} T_1$ and $b := \tilde{g}$, 
where $T_1$ is the unit translation on the line. 
It is not hard to see that these define an action of $BS(1,2)$. Indeed, from the relation $fgf^{-1} = g^2$ and the fact that 
both $\tilde{f}$ and $\tilde{g}$ fix the lifts of $p$, one easily concludes that $\tilde{f} \tilde{g} \tilde{f}^{-1} = \tilde{g}^2$. Since 
both $\tilde{f}$ and $\tilde{g}$ commute with $T_1$, this still gives $aba^{-1} = b^2$. Finally, note that $b$ fixes all the 
lifts of $p$, while $a$ acts on this set of lifts as a translation. Using this, it is not hard to see that the action is actually 
faithful, and topologically conjugate to the one previously constructed. Moreover, since $a$ and $b$ are smooth and 
commute with the unit translation, they are actually bi-Lipschitz on the whole real line.    
The dynamics of the group action is depicted on the right in the Figure 19 below.
\end{ex}
\end{small}

\begin{figure}[h!]
\begin{center}
\includegraphics[scale=0.34]{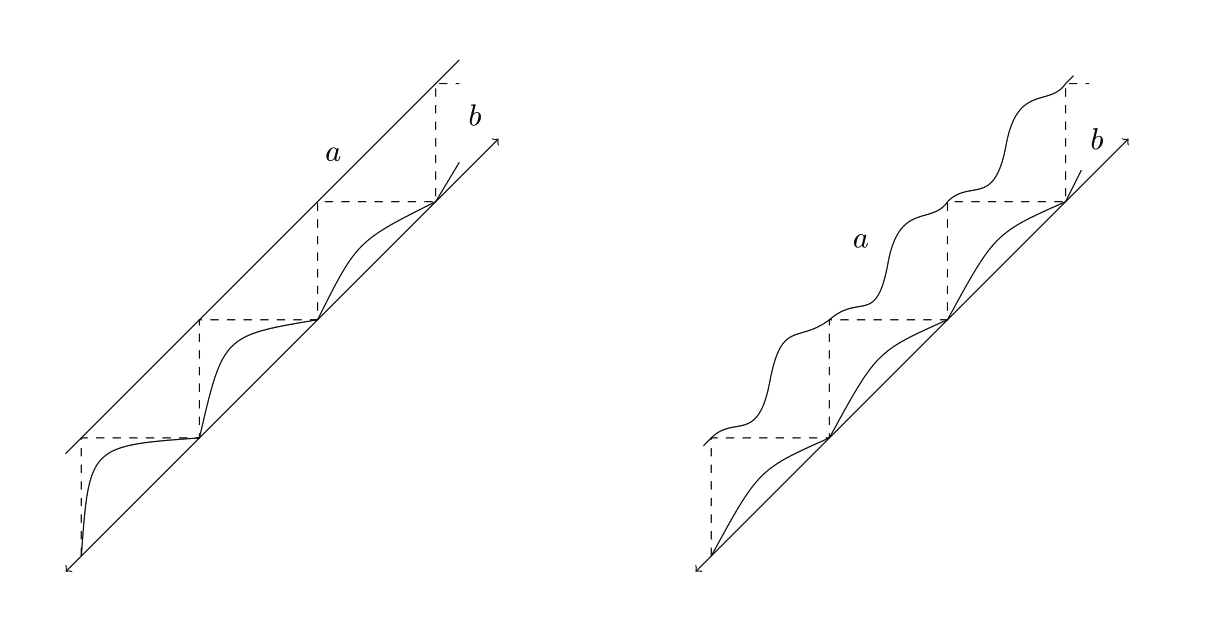}
\end{center}
\vspace{-1.32cm}
\begin{center}
Figure 19: Two conjugate actions of $BS(1,2)$.
\end{center}
\end{figure}

Although the action above is by real-analytic diffeomorphisms, things become very subtle when passing to actions on the closed interval. 
The following seminal result is essentially due to Cantwell and Conlon \cite{conlon} (see also \cite{GL}). 
For the statement, when we speak of semi-conjugacy  for an action on the 
interval, we refer to the action restricted to the interior, which is homeomorphic to the line. However, when dealing with regularity issues, 
usually (but not always) the key phenomena arise near the endpoints, as the reader will notice along the proof below (see also 
Exercise \ref{ejer:dif-borde} further on). 

\vspace{0.2cm}

\begin{thm} \label{thm-CC}
{\em The action of $BS(1,2)$ on the line built in Example \ref{ex-action-bil} is not 
topologically semiconjugate to an action by $C^1$ diffeomorphisms of the closed interval $[0,1]$.}
\end{thm} 

\noindent{\bf Proof.} Assume that we have an action by $C^1$ diffeomorphisms of $[0,1]$ that is semiconjugate to the action above, 
and keep denoting $a,b$ the $C^1$ diffeomorphisms associated to the generators. 
Let $x_0 \in (0,1)$ be a point that is fixed by $b$, 
and denote $x_n := a^{-n} (x_0)$. Let $y_0$ be a point in $(x_1,x_0)$ that is not fixed by $b$, and let $I$ be the interval 
with endpoints $y_0$ and $b(y_0)$. Note that the intervals $\{b^k (I): k \in \mathbb{N} \}$ have two-by-two disjoint interiors. 
Fix $k \in \mathbb{N}$, and write it in dyadic notation: 
$$k = \epsilon_0 2^0 + \epsilon_1 2^1 + \epsilon_2 2^2 + \ldots + \epsilon_{\ell} 2^{\ell}, \quad \mbox{where } \epsilon_i \in \{0,1\}, \, \epsilon_{\ell} = 1.$$
We thus see that
$$b^k = (b^{\epsilon_{\ell}})^{2^{\ell}} \cdots (b^{\epsilon_2})^4 (b^{\epsilon_1})^2 b^{\epsilon_0} = a^{\ell} b^{\epsilon_{\ell}} \cdots a^{-1} b^{\epsilon_2} a^{-1} b^{\epsilon_1} a^{-1} b^{\epsilon_0}.$$ 
Let $z \in I$ be a point such that $|b^k(I)| = Db^k (z) \! \cdot \! |I|$. By the chain rule, the previous relation gives 
\begin{equation}\label{eq:prueba-CC}
\frac{|b^k(I)|}{|I|} 
= \prod_{i=0}^{\ell} D b^{\epsilon_i} (z_i) \cdot \frac{Da^{\ell} (z')}{\prod_{i=0}^{\ell} Da (z_i)} 
= \prod_{i=0}^{\ell} D b^{\epsilon_i} (z_i) \cdot \frac{\prod_{i=1}^{\ell} Da (a^{i-1}(z))}{\prod_{i=1}^{\ell} Da (z_i)},
\end{equation}
where $z_i := a^{-1} b^{\epsilon_{i-1}} \cdots a^{-1} b^{\epsilon_0}(z)$ and $z' = b^{\epsilon_{\ell}} \cdots a^{-1} b^{\epsilon_2} a^{-1} b^{\epsilon_1} a^{-1} b^{\epsilon_0}(z)$.  
Note that, for each $i \geq 1$,  both $z_i$ and $a^{i-1}(z)$ belong to $a^{-i} ([x_1,x_0])$. 

Since $b$ fixes all the intervals $[x_i,x_{i-1}]$, on each of them there is a point where the derivative of $b$ equals 1. By the continuity of $Db$, 
this forces $Db(0) = 1$. Thus, if we fix a very small $\varepsilon > 0$ (actually, any $\varepsilon < 1 - \sqrt{2}/2$ will work for our argument),  
we can take $L \in \mathbb{N}$ large enough so that
$$Db (u) \geq 1-\varepsilon  \quad \mbox{for all } i \geq L \mbox{ and all } u \in a^{-i} ([x_1,x_0]).$$
Also, by the continuity of $Da$, up to slightly enlarging $L$, we may also assume that 
$$\frac{Da (v)}{Da (u)} \geq 1-\varepsilon \quad \mbox{for all } i \geq L \mbox{ and all } u,v \mbox{ in } a^{-i} ([x_1,x_0]).$$ 
Letting $u := z_i$ and $v := a^{i-1}(z)$ in these inequalities and using (\ref{eq:prueba-CC}),  
we obtain that there are  constants $C'$ and $C:=C' (1-2\varepsilon)^{-2L}$ such that 
\begin{equation}\label{eq:CC-explodes}
\frac{|b^k(I)|}{|I|} \geq C' (1-\varepsilon)^{2(\ell-L)} = C (1-\varepsilon)^{2\ell}.
\end{equation}

We now let $k$ vary from $2^{\ell}$ to $2^{\ell+1}-1$. Using (\ref{eq:CC-explodes}), we immediately obtain 
$$\sum_{2^{\ell} \leq k < 2^{\ell+1}} | b^k(I) | \geq 2^{\ell} \, |I| \, C \, (1-2\varepsilon)^{2\ell}.$$
Note that the right-side expression explodes as $\ell$ goes to infinite, due to our choice of $\varepsilon$. However, this is absurd, 
since as the intervals $b^k (I)$ in consideration have two-by-two disjoint interior, the sum of their lengths 
is smaller than 1. $\hfill\square$

\vspace{0.35cm}

Quite surprisingly, despite the theorem above, the action from Example \ref{ex-action-bil} is conjugate to an action on $[0,1]$ for which every 
element is $C^{\infty}$ on the interior and differentiable (but not continuously differentiable\,!) at the endpoints. This is a direct consequence of 
the claims in the next exercise, which is closely related to recent work of Virot (observe that every homeomorphism of the real line that commutes 
with the unit translation satisfies the key property (\ref{eq:bornes}) below).

\begin{small} \begin{ejer}\label{ejer:dif-borde}
Let $\varphi : (0,1) \to \mathbb{R}$ be a $C^{\infty}$ diffeomorphism that coincides with the map $x \mapsto -e^{1/x}$ (resp. $x \mapsto e^{1/(1-x)}$) 
close to $- \infty$ (resp. $+ \infty$).\\

\noindent (i) Show that if $f \in \mathrm{Homeo}_+(\mathbb{R})$ is a homeomorphism satisfying 
\begin{equation}\label{eq:bornes}
x - c \leq f(x) \leq x + c
\end{equation}
for a certain $c>0$ and every $x \in \mathbb{R}$, then the homeomorphism $\hat{f} := \varphi^{-1} \circ f \circ \varphi$, 
viewed as a map from $[0,1]$ to itself, is differentiable (with derivative equal to 1) at the endpoints.

\noindent{\underbar{Hint.}} For $x$ close to $0$, show that (\ref{eq:bornes}) implies 
$$- \frac{c}{e^{1/x}+c} \leq \frac{1}{x} - \frac{1}{\hat{f}(x)} \leq \frac{c}{e^{1/x}-c},$$ 
hence
\begin{equation}\label{eq-raro-raro}
\left| \frac{\hat{f}(x)}{x} - 1 \right| = \left| \hat{f}(x) \left[ \frac{1}{x} - \frac{1}{\hat{f}(x)} \right] \right| \leq \frac{c \,\hat{f}(x)}{e^{1/x}-c}.
\end{equation}
Check that the last expression converges to 1 as $x$ goes to the origin. 

\noindent (ii) Show that if $f$ is a bi-Lipschitz homeomorphism that satisfies (\ref{eq:bornes}), 
then $\hat{f}$ is also a bi-Lipschitz homeomorphism.

\noindent{\underbar{Hint.}} By the chain rule, for almost every point $x$, 
$$D \hat{f}(x) = \frac{D \varphi (x)}{D \varphi (\hat{f}(x))} \cdot D f (\varphi(x)).$$ 
Moreover, since $f$ is bi-Lipschitz, the factor $Df (\varphi (x))$ remains uniformly bounded (almost everywhere). 
For the quotient of derivatives of $\varphi$, check the equality 
$$ \frac{D \varphi (x)}{D \varphi (\hat{f}(x))} 
= \frac{e^{1/x} / x^2}{e^{1/\hat{f}(x)}/\hat{f}(x)^2} = \left(\frac{\hat{f}(x)}{x} \right)^2 \exp \left( \frac{1}{x} - \frac{1}{\hat{f}(x)} \right),$$
and use (\ref{eq-raro-raro}) to show that it remains uniformly bounded as well.
\end{ejer} \end{small}

\begin{question} Does there exist a (finitely-generated) left-orderable group having no differentiable action on the closed 
interval or the real line\,? Here, by {\em differentiable} we mean that every group element is a differentiable function
(with finite derivative at every point), yet the derivative may vary discontinuously.
\end{question}

\noindent{\bf On $C^1$ actions of other left-orderable groups.} 
What is nice about the nonsmoothable action of $BS(1,2)$ previously discussed is that it is Conradian. Indeed, as we will see in the example below 
(built on a seminal remark by Bonatti, Croivisier and Wilkinson), it is not hard to produce non-Conradian actions that are not smoothable on $[0,1]$.

\begin{small}
\begin{ex} \label{ex-BCW}
Let $\Gamma$ be a group of homeomorphisms of $[0,1]$ containing two elements $f,g$ that are linked forming a {\em resilient pair} 
$(f,g;u,v)$. Recall that this means that $u < f(u) <  f(v) < g(u) < g(v) < v$, and that such a pair arises for every non-Conradian action; 
see \S \ref{section-crossings}. Assume that $\Gamma$ also 
contains an element $h$ with no fixed point in $(0,1)$ that commutes with both $f$ and $g$. We claim that the action of $\Gamma$ on 
$[0,1]$ cannot be by $C^1$ diffeomorphisms. Indeed, by the resilience property and the commutativity assumption, for each $n \in \mathbb{Z}$ 
we have that $(f,g; h^n (u),h^n(v))$ is also a resilient pair for $f$ and $g$. On the one hand, this implies that $f$ and $g$ have fixed points 
on each interval $[h^n(u),h^n(v)]$, and since these intervals converge to the endpoints of $[0,1]$ as $n$ goes to $\pm \infty$, this forces 
$Df(0) = Df(1) = Dg(0) = Dg(1) =1$. On the other hand, the resilience condition obviously implies that, on each interval $[h^n(u),h^n(v)]$, 
there must be a point $x_n$ at which either $Df (x_n) < 1/2$ or $Dg(x_n) < 1/2$. This contradicts the continuity of at least 
one of the derivatives $Df$ or $Dg$. 
\end{ex}
\end{small}

Besides the action built in Example \ref{ex-action-bil}, the group $BS(1,2)$ also acts by affine transformations on the projective line, 
which can be thought of as an action on the interval if we cut the circle at the global fixed point. Thus, the smoothness obstruction of 
Theorem \ref{thm-CC} does not apply to all actions of $BS(1,2)$. However, building on the proof technique of the same 
theorem plus some algebraic considerations, it is proved in \cite{tst} that the group $\Gamma \!=\! \mathbb{F}_2 \ltimes \mathbb{Z}^2$ 
discussed in Example \ref{no-tiene-pelo} has no faithful action by $C^1$ diffeomorphisms of the interval. Again, what is nice about this 
group is that it is still locally indicable. Indeed, for non locally indicable groups, no $C^1$ action on the interval can arise, because 
of the celebrated {\em stability theorem} of Thurston; see \cite[Chapter 5]{yo} for a full discussion on this topic. 

\vsp

\begin{small}
\begin{ejer}\label{ejer-tst}
In \S \ref{classic-conrad} we discussed the following example of a non locally-indicable, left-orderable group:
$$\Gamma = \big\langle f,g,h \!: f^2 = g^3 = h^7 = fgh \big\rangle.$$
Although this is a subgroup of $\widetilde{\mathrm{PSL}} (2,\mathbb{R})$, it has no action by $C^1$ diffeomorphisms on the interval, 
because of Thurston's stability theorem referred to above. 
The reader is invited to give an alternative argument by developing the two items below.  

\noindent (i) Prove that for every action of $\Gamma$ on $(0,1)$ with no global fixed point, the central element $fgh$ has no fixed point.

\noindent (ii) Using Example \ref{ex-BCW}, conclude that $\Gamma$ has no nontrivial action by $C^1$ diffeomorphisms of $[0,1]$.
\end{ejer}
\end{small}

Actually, the group $\mathbb{F}_2 \ltimes \mathbb{Z}^2$ is even nicer: it has no faithful action by $C^1$ diffeomorphisms of the line \cite{tst} 
(note that the group $\Gamma$ from Exercise \ref{ejer-tst} is naturally a group of real-analytic diffeomorphisms of the line).  Many 
other groups share this property, but in most cases this is difficult to establish. A particularly interesting example is Higman's group 
$H := \langle a_i \, (i\in \Z/4\Z) \mid a_i a_{i+1}a_i^{-1}=a_{i+1}^2\rangle$.  
Indeed, in \cite{RT-Higman}, 
it is proved that $H$ is left-orderable, though it admits no nontrivial action by $C^1$ diffeomorphisms of the real line.

\vspace{0.35cm}

\noindent{\bf On small groups/actions.} The discussion above suggests that obstructions to $C^1$ smoothability are more common for 
``large'' groups/actions, and tend to disappear for ``small'' ones. The example contained in the next exercise should however give us 
some warning on this claim.

\begin{small}
\begin{ejer}
Another example of a locally indicable group having no $C^1$ action on the interval is the Baumslag-Solitar group 
$BS(1,-2) = \langle aba^{-1} = b^{-2} \rangle$ (yet it admits smooth actions on the real line). Show this by developing the items below:

\noindent (i) Show that for every faithful action of $BS(1,-2)$ on the real line with no global fixed point, the element $a$ acts with no fixed 
point, but $b$ fixes points in each interval with endpoints $x,a(x)$, for all $x \in \mathbb{R}$ (compare Example \ref{ex-prep}). 

\noindent (ii) Using the argument of proof of Theorem \label{thm-CC}, conclude that $BS(1,-2)$ has no faithful action by $C^1$ 
diffeomorphisms of the closed interval.

\noindent{\underline Remark.} This example is less satisfatory than the preceding ones because $BS(1,-2)$ has only a few 
actions on the line. Indeed, it is a Tararin group admitting only four left-orders (see \S \ref{fin-uncount}). 
\end{ejer}
\end{small}

Despite the example described in the exercise above, the following question arises naturally.

\begin{question} Does there exist a finitely-generated, left-orderable group that is not a Tararin group but has only Conradian left-orders 
and admits no $C^1$-smoothable action on the interval or the real line ?
\end{question}

The following remarkable recent result by Kim, Matte Bon, de la Salle and Triestino \cite{KMST} points in the negative direction. 
For the statement, recall that a finitely-generated group $\Gamma$ has {\bf{\em subexponential growth}} if the number of elements in the ball of radius
$n$ with respect to a finite generating system grows subexponentially in $n$ (this property does not depend on the choice of the finite generating 
system). Such a group cannot contain free subsemigroups, hence all of its left-orders (if any) must be Conradian (see \S \ref{section-crossings}). 

\begin{thm} {\em Every action of a group of subexponential growth on the closed interval is semiconjugate to a group of $C^1$ diffeomorphisms. 
Moreover, for every $L > 1$, the semiconjugacy can be chosen so that for all generators $f$ and all $x \in [0,1]$ one has 
$1/L < D f(x) < L$.}
\end{thm}

It is worth stressing that the statement doesn't claim that the original action is conjugate to an action by $C^1$ diffeomorphisms. 
It only deals with semiconjugacies, and passing to a genuine conjugacy seems to be a subtle issue. Despite this, the theorem 
somehow extends and simplifies several works in the literature, for instance \cite{CJN,FF,jorquera,navas-approx,growth,parkhe}. 
In order to get a taste of it, the reader is invited to work through the exercise below.

\begin{small}
\index{Group!of subexponential growth}
\begin{ejer} \label{close to 1}
Assume that a subgroup $\Gamma$ of $\mathrm{Homeo}_+(\mathrm{R})$ has subexponential growth. 
Show that, for every $L \!>\! 1$, it is possible to simultaneously conjugate the generators of $\Gamma$ to $L$-bi-Lipschitz homeomorphisms.

\noindent{\underline{Hint.}} In the proof of Theorem  \ref{t: lipschitz conjugation} above, the (positive) function $\phi \in \mathcal{L}^1 (\Gamma)$
can be taken so that $\phi (hg) \leq L^{1/3} \phi (g)$ for every $h \in \mathcal{G}$ and every $g \in \Gamma$. See \cite{navas-approx} for more on this.
\end{ejer}\end{small}


\subsection{Actions almost having fixed points} \label{s: almost-fixed}

\index{Action!almost having fixed points}
\hspace{0.45cm} Let $\Gamma$ be a finitely-generated group with finite generating set
$\mathcal{G}$, and let $\Phi$ be an almost-periodic action of $\Gamma$ on $\mathbb R$.  (Recall
that we do not assume $\Phi$ to be faithful.) We say that $\Phi$ {\bf \textit{almost has fixed points}}  if
$$\inf _{t\in \mathbb R} \sup _{g\in \mathcal{G}} \big| \Phi(g)(t) - t \big| = 0.$$
An equivalent way to think about this property is that the action of $\Gamma$ on the compact space constructed
in Proposition~\ref{l:closure under translation flow} has a global fixed point in the closure of the orbit
of $\Phi_0$ by the translation flow \(T\).

It is not obvious how to construct almost-periodic actions that do not almost have fixed points.
(Consider, for instance, the case of affine groups.) This is actually the goal of the next result. 

\vsp

\begin{thm}
\label{t: almost-periodic action}
\textit{Every action of a finitely-generated group by orientation-preserving homeomorphisms of the
real line is topologically conjugate to an almost-periodic action. Moreover, if the original action has
no global fixed point, then there is such a conjugate that does not almost have fixed points.}
\end{thm}

\vsp

To prove this result, let $\Gamma$ be a finitely-generated group provided with a finite, symmetric
system of generators $\mathcal{G}$. Given constants $L>1$ and $D > D' > 0$, we denote
$R = R (\Gamma, \mathcal{G}, L, D, D')$ the set of representations
$\Phi \!: \Gamma \rightarrow \mathrm{BiLip}_+(\mathbb R)$
such that every $g\in \mathcal{G}$ satisfies $L (\Phi(g)) \leq L$ and
\begin{equation}\label{eq: displacement condition}
t-D
\leq \min_{g\in \mathcal{G}} \Phi(g)(t)
\leq t - D' \leq t+ D'
\leq \max _{g\in \mathcal{G}} \Phi(g)(t)
\leq t + D
\end{equation}
for all $t \in \mathbb{R}$.
This set can be seen as a closed subset of $\mathrm{BiLip}_+ (\mathbb R)^{\mathcal{G}}$,
and as such is equipped with the product topology. Relations (\ref{eq: lipschitz}) and
(\ref{eq: displacement condition}) imply that $R$ is compact, by the Arzela-Ascoli 
theorem. Moreover, the same relations show that the translation flow $T$
defined by~\eqref{eq: translation flow} preserves $R$. Hence, every element of $R$ is an
almost-periodic action of $\Gamma$, and (\ref{eq: displacement condition})
shows that, moreover, such an element does not almost have fixed points.

\vsp\vsp

\begin{lem} \label{l: R non empty}
\textit{Let $\Phi_0 \!: \Gamma \rightarrow \mathrm{Homeo}_+ (\mathbb R)$ be a
faithful action without global fixed points of a finitely-generated group $\Gamma$.
Then there are constants $L>1$ and $D  \!>\! D' \!>\! 0$, as well as a finite, symmetric
generating system of $\Gamma$, such that the corresponding set $R$  contains
a representation that is conjugate to $\Phi_0$.}
\end{lem}

\noindent{\bf Proof.} We see $\Gamma$
as being contained in $\mathrm{Homeo}_+(\mathbb{R})$ via $\Phi_0$. 
By Theorem \ref{t: lipschitz conjugation}, it is enough to prove the statement 
in the case where $\Gamma$ is a subgroup of $\mathrm{BiLip}_+(\mathbb R)$. 
Let $\mathcal{G}$ be a finite, symmetric generating set of $\Gamma$, and let $L$ be a
constant such that every $g \!\in\! \mathcal{G}$ is $L$-bi-Lipschitz.
Let $(t_n)_{n\in \mathbb Z}$ be the sequence of points in $\mathbb R$ defined by
$t_0 := 0$ and $t_{n+1} := \max _{g\in \mathcal{G}} g(t_n)$. Equivalently,
$t_{n-1} = \min_{g\in \mathcal{G}} g (t_n)$, as $\mathcal{G}$ is symmetric.
Since $\Gamma$ has no fixed point on the real line, we have
\[ \lim _{ n\rightarrow \pm \infty} t_n = \pm \infty. \]
Let $\varphi$ be the homeomorphism of the real line that sends $t_n$ to $n$ and that
is affine on each interval $[t_n, t_{n+1}]$. We claim that the action of $\Gamma$
defined by $\Phi (g) := \varphi \circ g \circ \varphi^{-1}$ belongs to
$R(\Gamma,\overline{\mathcal{G}}, L^6,1,4)$ for the generating set
$\overline{\mathcal{G}} := \mathcal{G} \cup \mathcal{G}^2$.

To prove this, we first note that the distortion of the sequence $(t_n)$ is uniformly bounded. In
concrete terms, if for each $n \!\in\! {\mathbb Z}$ we denote $\delta_n := t_{n+1} - t_n$, then
\begin{equation} \label{eq: distortion}
L^{-1} \delta_{n+1}  \leq  \delta_n \leq L \hspace{0.05cm} \delta_{n+1}.
\end{equation}
Indeed, let $g_n \in \mathcal{G}$ be such that $t_{n+1} = g_{n} (t_n)$. By definition,
$g_n (t_{n+1} ) \leq t_{n+2}$, and since $g_n$ is an $L$-bi-Lipschitz map, we have
\[  t_{n+2} - t_{n+1} \geq g_n (t_{n+1} ) - g_n(t_n) \geq L^{-1} (t_{n+1} - t_n),\]
which yields the right-side inequality in (\ref{eq: distortion}). (The left-side
inequality is obtained analogously.)
Note that, by construction, for all $w,z$ in $[t_{n}, t_{n+1}]$,
\begin{equation}\label{eq: nec}
\big| \varphi (z) -\varphi (w) \big|
= \frac{\hspace{0.05cm} |z-w|}{\delta_n}.
\end{equation}

We next claim that for every $g\in \mathcal{G}$, we have $L (\Phi(g) ) \leq L^3$. To show this,
it suffices to prove that each such $\Phi(g)$ is Lipschitz on every interval 
$[n,n+1]$,
with Lipschitz constant at most $L^3$. To check this, consider two arbitrary points
$x,y$ in $ [n,n+1]$, and define $w := \varphi^{-1} (x)$ and $z:= \varphi^{-1} (y)$. Then $w,z$
both belong to $[t_{n},t_{n+1}]$, which in virtue of (\ref{eq: distortion}) and (\ref{eq: nec}) yields
\begin{small}
$$\big| \Phi(g) (y) - \Phi(g) (x) \big|
= \big| \varphi ( g (z) ) - \varphi (g(w)) \big|
\leq \frac{L\hspace{0.04cm} \big| g(z)-g(w) \big|}{\delta_n}
\leq \frac{L^2\hspace{0.04cm} \big| z-w \big|}{\delta_n}
\leq L^3 |y-x|,$$
\end{small}as desired.

By construction, for every generator $g\in \mathcal{G}$ and all $x \in \mathbb{R}$,
\[ x-2 \leq  \Phi (g) (x) - x \leq x +2.\]
Indeed, the integer points just after and before $x$ are moved a distance less than or equal to $1$ by
$\Phi(g)$. Moreover, as for every $n\in {\mathbb Z}$ we have $\Phi(g_{n+1} g_n) (n) = n+2$, letting 
$n$ be the integer part of $x$, this yields $\Phi(g_{n+1} g_n) (x) \geq x+1$. We have hence
proved that $\Phi $ belongs to $R (\Gamma,\overline{\mathcal{G}}, L^6 , 1, 4 )$.
$\hfill\square$

\vsp\vsp\vsp

Theorem~\ref{t: almost-periodic action} immediately follows from the preceding lemma in the case 
where $\Gamma$ has no global fixed point. If such a point exists, we replace $\Gamma$ by the
free product $\Gamma * \mathbb Z$ (or a quotient of it) and we extend $\Phi_0$ so that the generator
of the $\mathbb{Z}$-factor is mapped to a nontrivial translation. This new group has no fixed point, so
that the preceding lemma applies to it. We thus conclude that the $\Gamma$-action is topologically 
conjugate to an almost-periodic action.


\subsection{Free almost-periodic spaces}  
\label{s:free-aps}

The previous construction allows replacing the space of left-orders 
of a given left-orderable group by an object that is provided with a flow and is more natural when dealing 
with dynamical realizations.

\vsp\vsp\vsp

\begin{cor} \label{c:substitute}
\textit{Let $\Gamma$ be a finitely-generated, left-orderable group. Then there exists a compact space
$X$, a free flow $T=\{T_s\}_{s\in \mathbb R}$ on $X$, and an action of $\Gamma$ on $X$ without global fixed
points which preserves the $T$-orbits together with their orientations.}
\end{cor}

\noindent{\bf Proof.}
Since $\Gamma$ is finitely-generated and left-orderable, it admits a faithful action by orientation-preserving
homeomorphisms of the real line without global fixed point (see \S \ref{general-3}).
By Theorem~\ref{t: almost-periodic action}, this action
is conjugate to an almost-periodic action $\Phi_0$ that does not almost have fixed points. Consider the space
$X$ constructed in the second part of the proof of Proposition~\ref{l:closure under translation flow} together with 
the free flow $T$ and the $\Gamma$-action on it. Because $\Phi_0$ does not almost have fixed points
and $\{ T_s (\Phi_0) \!: s \in \mathbb{R} \}$ is dense in $X$, there is no fixed point for the $\Gamma$-action on $X$.
Moreover, $\Gamma$ stabilizes every $T$-orbit, and preserves the orientation on each of them.
$\hfill\square$

\vsp\vsp\vsp

In the sequel, a space \(X\) together with a \(\Gamma\)-action and a flow \(T\) as in the conclusion of 
the preceding corollary will be called a {\bf \em free almost-periodic space} for the group \(\Gamma\). 
Referring to the terminology introduced in \S \ref{s: almost periodic}, this is an almost-periodic space 
for which the associated flow $T$ is free and the underlying $\Gamma$-action has no fixed point.

\vsp

\begin{small}\begin{ejer}
Prove that under the assumptions of Corollary \ref{c:substitute}, there exists a free almost-periodic space \(X\) that can be endowed with a metric so that the \(\Gamma\)-action and the translation flow \(T\) act by Lipschitz transformations (we point out that, in general, this metric need not be of finite Hausdorff dimension).
 
\noindent{\underline{Hint.}} 
Equip the subset $R = R (\Gamma, \mathcal{G}, L, D, D')\subset APA_+(\Gamma)$ with the distance \(d(\Phi, \Phi')\) defined as the infimum of the real numbers so that 
\[ |\Phi (g) (x) - \Phi' (g) (x) | \leq d(\Phi, \Phi') \ L^{||g||}\  \exp \left(\frac{\log L}{D} |x| \right)  \text{ for every } x\in \mathbb R, \]
where \(||g||\) is the minimum number of elements of \(\mathcal{G}\) needed to write \(g\).
Show that the metric space \( (R, d) \) is compact, and that the \(\Gamma\)-action together with the translation flow 
are actions by Lipschitz homeomorphisms with respect to \(d\).   
\end{ejer}\end{small}

\vsp

There is a kind of converse philosophy to that of Corollary \ref{c:substitute}: a flow on a compact metric space with some 
of the features above can be used to construct interesting left-orderable groups. This is the case, for instance, with the group 
introduced by Matte Bon and Triestino in \cite{MT}, which is fully described in \S \ref{s: finitely generated simple left orderable}. 
In that setting, the translation flow corresponds to the suspension flow over a subshift of finite type. The construction of 
the group then proceeds via a careful study of the group of homeomorphisms of the total space that preserve each orbit 
of the flow individually, acting on them by orientation-preserving, piecewise-dyadic homeomorphisms.


\subsection{Indicability of amenable, left-orderable groups revisited}\label{s: Witte}

\hspace{0.45cm} Based on the previous construction, and following \cite{deroin},
we next give an alternative proof of Theorem \ref{witte-bien}. Let $\Gamma$ be
a finitely-generated, left-orderable group, and let $X$ be an almost-periodic space equipped
with a free flow $T$ 
and a $\Gamma$-action, as described by
Corollary~\ref{c:substitute}. If $\Gamma$ is amenable, then there exists a probability measure
$\mu$ on $X$ that is invariant by $\Gamma$. Consider the conditional measures of $\mu$
along the orbits of the translation flow $T$. These are defined in the following way: in a flow box \( h : U \rightarrow I \times S\) (see Exercise \ref{ex: flow box}), use Fubini's theorem to disintegrate the measure \( h_* \mu \) as an integral \(\int _ S \left(  \mu_{I\times s} \otimes \delta _{s} \right) d\nu (s)\), where \(s\mapsto  \mu_ {I\times s}\) is a measurable family of probability measures on the interval \(I\). The conditional measure on the plaque \( h^{-1} (I\times\{s\}) \) is the measure \( (h^{-1})_* \mu_{I\times s}\). We leave it to the reader to verify that on the intersection of two plaques of two different flow boxes, the two conditional measures differ by multiplication by a constant. In particular, by considering long flow boxes, this enables us to construct Radon measures on $\mu$-a.e.  
$T$-orbit that are well-defined up to multiplication by a positive constant.  We denote
by $\mu_l$ this Radon measure on the $T$-orbit $l$; it is well-defined on the orbit of \(\mu\)-a.e. point. Note that this family of projective Radon measures is canonically associated with the lamination space induced by the flow \(T\), namely it is invariant under reparametrization of the flow. 

The countable group $\Gamma$ preserves $\mu$ and the laminated structure induced by the flow \(T\). Therefore, 
for $\mu$-a.e. $T$-orbit $l$ in $X$,  the measure $\mu_l$ is nonzero, and every $g \in \Gamma$ multiplies it by a certain factor:
\[  g_* (\mu_{l}) = c_{l} (g) \mu_{l}, \ \ \ \mathrm{where}\ \ \ c_{l} (g) >0 .\]
If $\mu_{l}$ is not preserved by $\Gamma$, then the map $g \mapsto \log c_{l}(g)$ is a nontrivial
homomorphism from $\Gamma$ into $({\mathbb R}, +)$. (See Remark \ref{cannot-arise} below concerning this case.)
Otherwise, $\mu_{l}$ is preserved by $\Gamma$. If $\mu_{l}$ has an atom, then its orbit must be discrete,
and $\Gamma$ acts by translations along this orbit, thus giving rise to a nontrivial homomorphism into the
integers. If $\mu_l$ has no atom, then the $\Gamma$-action on $l$ is semiconjugate to an action by
translations, which induces a nontrivial homomorphism into the reals.

\begin{small}\begin{rem} \label{cannot-arise}
It seems that the condition $\mu_l = c_l (g) \mu_l$ for a function $c_l$
that is not identically equal to 1 cannot arise in the context above (with positive measure). 
At least, it doesn't arise for the space $X$ built upon $APA_+(\Gamma)$ in the proof of 
Lemma \ref{l:closure under translation flow}. Indeed, assuming otherwise, 
the action of $\Gamma$ on $l$ is semiconjugate to that of a non-Abelian affine group 
(see \S \ref{section-solvable-spaces}). In particular, there must exist a resilient pair \,
$u \!<\! f(u) \!<\! f(v) \!<\! g(u) \!<\! g(v) \!<\! v$ \, for the action on $l$;  
moreover, there is an element $h \in \Gamma$ whose inverse (and all of its iterates) sends $[u,v]$
into a disjoint interval (hence to a region where no crossing for $f,g$ arises). 
Denoting the associated representation by $\Phi$, we thus have that, for all $n \!\in\! \mathbb{N}$, 
the conjugate representations $h^n (\Phi)$ remain outside a certain neighborhood of
$\Phi$. However, this is in contradiction with the Poincar\'e recurrence theorem.
(We refer to Examples \ref{ex:new-proof} and \ref{ex:BS} for another application of this idea.)
\end{rem}\end{small}


\section{Random Walks on Left-Orderable Groups} \label{s:random walks}

\hspace{0.45cm} Our goal now is to provide a more conceptual proof of the existence of almost-periodic
actions for left-orderable groups based on probabilistic arguments. Throughout this section, $\Gamma$ will denote
a finitely-generated group and $\rho$ a probability measure on $\Gamma$ that is {\bf{\em symmetric}}, in the 
sense that $\rho (g) = \rho (g^{-1})$ for all $g \in \Gamma$, and whose support $\mathcal{G}$ generates $\Gamma$. 
Although otherwise stated, $\mathcal{G}$ will also be assumed to be finite. 
\index{Measure!symmetric}

We start with an emphasis on a particular type of actions, namely, those for which the Lebesgue measure 
is {\em stationary}, {\em i.e.}, invariant on average. These actions, called {\em $\p$-harmonic}, will appear to have many 
nice properties. In particular, we will see that they are 
always almost-periodic. Quite remarkably, we will show that all actions on 
$\mathbb R$ become harmonic under suitable conjugacies, and these conjugacies are unique up to 
post-composition with an affine map. 


\subsection{Harmonic actions and Derriennic's property}
\label{S: derriennic}

\index{Action!$\rho$-harmonic}
\index{Measure!stationary}
\hspace{0.45cm}  Let $\Gamma$ be a subgroup of $\mathrm{Homeo}_+(\mathbb{R})$. 
The action of $\Gamma$ is said to be {\bf \textit{$\p$-harmonic}}
(or just harmonic, if the probability $\p$ is clear from the context) if the Lebesgue measure is
{{\em stationary}}, that is, if for every $x,y$ in $\mathbb R$,
\begin{equation} \label{eq:harmonic}
y - x = \int _{\Gamma} \big( g(y) - g (x) \big) \, d\p (g) =
\sum_{g \in \mathcal{G}} \big( g(y) - g(x) \big) \hspace{0.05cm} \rho(g).
\end{equation}
We will assume throughout that the $\Gamma$-action has no global fixed point. 
However, we will see in Exercise \ref{e:irreducible} below that this assumption is redundant, 
since no group action on the line satisfying property (\ref{eq:harmonic}) above can globally fix a point.

Obviously, $\p$-harmonic actions include those that satisfy, for every $x\in \mathbb R$,
\begin{equation*}
x = \int_{\Gamma} g (x) \, d\p (g).
\end{equation*}
This will be called the {\bf \textit{Derriennic property}}, as it corresponds to a weak form
of a property studied by Derriennic in~\cite{Derriennic} in the more general context of
Markov processes on the line (not necessarily coming from a group action). Quite
surprisingly, as was cleverly noticed by Kleptsyn, all $\p$-harmonic actions 
satisfy this property.
\index{Derriennic's property}

\vsp

\begin{prop}\label{prop: Derr}
\textit{Every $\p$-harmonic action has the Derriennic property. }
\end{prop}

\vsp

For the proof, we need the following elementary lemma.

\vspace{0.1cm}

\begin{lem}\label{l:Phi-int}
\textit{For all $h \in \mathrm{Homeo}_+({\mathbb R})$
and each compact interval $[a,b]$, we have
\begin{equation}\label{eq:drift}
\int_a^b \big[ (h(x)-x)+ (h^{-1}(x)-x) \big] \, dx = \Delta^{h}(b) - \Delta^{h}(a),
\end{equation}
where $\Delta^{h} (x)$ is the non-signed area of the region depicted in Figure 20 below:
$$\Delta^{h}(c) := \left \{ \begin{array} {l}
\int_{ h^{-1}(c)}^{c} \big[ h (s) - c \big] ds,
\hspace{0.67cm} \mbox{ if } \esp h (c) \geq c,\\
\\
\int_{h (c)}^c \big[ h^{-1}(s) - c \big] ds,
\hspace{0.62cm} \mbox{ if } \esp h (c) \leq c.
\end{array} \right.$$}
\end{lem}

\noindent{\bf Proof.} Denoting $|A|$ the Lebesgue measure of a subset
$A \subset \mathbb{R}^2$, we have that  $\int_a^b ( h (x)-x) \, dx \esp$ equals
$$\big| \big\{(x,y) \!: a<x<b, \, x<y<h(x) \big\} \big|
\esp - \esp \big| \big\{ (x,y)  \!: a<x<b, \,  h (x)<y<x \big\} \big|,$$
which may be rewritten as
\begin{small}
$$\big| \big\{(x,y) \!: a\!<\!x\!<\!b, \, b\!<\!y\!<\! h(x) \big\} \big|
+ \big| \big\{ (x,y)  \!: a\!<\!x\!<\!b, \, a\!<\!y\!<\!b, \, x\!<\!y\!<\! h(x) \big\} \big|$$
$$- \big| \big\{ (x,y)  \!: a\!<\!x\!<\!b, \, h(x)\!<\!y\!<\!a \big\}  \big|
-  \big| \big\{ (x,y)  \!: a\!<\!x\!<\!b, \, a\!<\!y\!<\!b, \, h(x)\!<\!y\!<\!x \big\} \big|.$$
\end{small}A similar
equality holds for $h^{-1}$. Now, in the sum
$$\int_a^b (h (x) - x) \, dx + \int_a^b (h^{-1}(x)-x) \, dx,$$
the corresponding
second and fourth terms above cancel each other. Indeed, these terms involve all
couples $(x,y) \in [a,b]^2$, and we have $x<y<h(x)$ if and only if $h^{-1}(y)<x<y$.
Therefore, the second term for $h$ is exactly the negative of
the fourth term for~$h^{-1}$, and vice versa.

As a consequence, the value of
$$\int_a^b \big[ ( h (x)-x)+ ( h^{-1}(x)-x) \big] \, dx$$
equals
$$\big| \big\{ (x,y)  \!: a\!<\!x\!<\!b, \,\, b\!<\! y\!<\! h(x) \big\} \big| \, +
\, \big| \big\{ (x,y)  \!: a\!<\!x\!<\!b, \,\, b\!< \!y\!<\! h^{-1}(x) \big\}  \big|$$
$$- \, \big| \big\{ (x,y)  \!: a\!<\!x\!<\!b, \, h(x)\!<\!y\!<\!a \big\} \big| \, -
\, \big| \big\{ (x,y)  \!: a\!<\!x\!<\!b, \,\, h^{-1}(x) \!<\!y\!<\!a \big\}  \big|,$$
and one can easily check that the expressions above and below are
equal to $\Delta^{h}(b)$ and $\Delta^{h}(a)$, respectively.
This proves the desired equality.
$\hfill\square$

\vspace{0.75cm}


\beginpicture
\setcoordinatesystem units <1cm,1cm>

\putrule from 0 0 to 0 3
\putrule from 3 0 to 3 3
\putrule from 0 3 to 3 3
\putrule from 0 0 to 3 0
\plot 0 0 3 3 /

\putrule from 7 0 to 7 3
\putrule from 10 0 to 10 3
\putrule from 7 3 to 10 3
\putrule from 7 0 to 10 0
\plot 7 0 10 3 /

\setquadratic

\plot
0 1.4 2 3.5 3 4.3 /

\plot
7 0.8 9 1.7 11.2 3 /
\plot
7.8 0 8.7 2 10 4.2 /

\setlinear
\setdots
\plot 3 3 3 4.3 /
\plot 0 4.3 3 4.3 /
\plot 0 3 0 4.3 /
\plot 1.5 0 1.5 3 /

\plot 1.6 3 1.6 3.1 /
\plot 1.7 3 1.7 3.3 /
\plot 1.8 3 1.8 3.3 /
\plot 1.9 3 1.9 3.4 /
\plot 2   3 2   3.5 /
\plot 2.1 3 2.1 3.6 /
\plot 2.2 3 2.2 3.7 /
\plot 2.3 3 2.3 3.8 /
\plot 2.4 3 2.4 3.9 /
\plot 2.5 3 2.5 4 /
\plot 2.6 3 2.6 4 /
\plot 2.7 3 2.7 4.2 /
\plot 2.8 3 2.8 4.2 /
\plot 2.9 3 2.9 4.3 /
\plot 3   3 3   4.3 /

\plot 9.3 3 9.3 3 /
\plot 9.4 3 9.4 3.2 /
\plot 9.5 3 9.5 3.4 /
\plot 9.6 3 9.6 3.5 /
\plot 9.7 3 9.7 3.6 /
\plot 9.8 3 9.8 3.8 /
\plot 9.9 3 9.9 4 /
\plot 10  3 10  4.2 /

\plot 10 0 11.2 0 /
\plot 11.2 0 11.2 3 /
\plot 10 3 10 4.2 /
\plot 10 3 11.2 3 /
\plot 7 2.3 10 2.3 /
\plot 9.3 0 9.3 3 /

\put{$h$} at 0.6 2.4
\put{$c$} at 3 -0.3
\put{$h^{-1}(c)$} at 1.5 -0.3
\put{$h(c)$} at -0.5 4.3
\put{$c$} at -0.3 3
\put{$\longrightarrow \Delta^{h}(c)$} at 4 3.6

\put{$c$} at 10 -0.3
\put{$h^{-1}(c)$} at 11.6 -0.3
\put{$h(c)$} at 9.3 -0.3
\put{$h(c)$} at 6.56 2.3
\put{$h$} at 9.6 1.7
\put{$h^{-1}$} at 8.34 2
\put{$\longrightarrow \Delta^{h}(c)$} at 11 3.6

\put{} at -2.1 0

\put{Figure 20: The definition of $\Delta^{h}(c)$ in the two possible cases.}
at 5.5 -1.1

\endpicture


\vspace{0.875cm}

\noindent{\bf Proof of Proposition \ref{prop: Derr}.}
First, note that, by $\p$-harmonicity, the value of
\[ \int_{\Gamma} \big( g(x) -x \big) d\p (g) \]
is independent of $x$. We call it the {\bf {\em drift}} of the action and denote it
by $Dr (\Gamma ,\p)$. The statement to be proved is hence equivalent to the vanishing
of the drift. To show this, we integrate (\ref{eq:drift}) over $\Gamma$ and use the
symmetry of $\rho$ to obtain, for all $a<b$,
$$
2(b-a) Dr(\Gamma, \p) = \int_{\Gamma} (\Delta^{g} (b) - \Delta^{g}(a) \big) \, d\p (g).
$$
Denoting now $\Delta(c) := \int_{\Gamma} \Delta^{g}(c)  \, d\p (g)$, this yields
$$
2(b-a) Dr(\Gamma, \p) = \Delta(b) - \Delta(a).
$$
The last equality shows that $\Delta$ is an affine function. Moreover,
$\Delta$ is an average of non-negative functions, thus it is non-negative.
Therefore, $\Delta$ must be constant, which implies that
$Dr(\Gamma,\p)=0$, as desired.
$\hfill\square$

\vsp\vsp\vsp

The constant  function \(\Delta\) is then an invariant of the harmonic action. 
In the sequel, by abuse of notation, we will denote its (positive) value by \(\Delta\). Note that \(\Delta\) depends continuously 
on the \(\rho\)-harmonic action if we equip the set of actions with the compact-open topology. The next proposition 
shows the relevance of the Derriennic property in the study of almost-periodic actions. Recall that \(\mathcal{G}\) 
stands for the support of the underlying probability distribution  \(\p\).

\vspace{0.1cm} 

\begin{prop} \label{P: harmonic are almost-periodic}
\textit{Every $\p$-harmonic action is almost-periodic and does not almost have fixed points. More precisely, every 
element is a Lipschitz homeomorphism, and there exist positive constants \( D,D'\) (depending only on \(\p\)) such that 
for all $x \in \mathbb{R}$,
\[x-D\sqrt{\Delta}
\leq \min_{g\in \mathcal{G}} g(x)
\leq x - D'\sqrt{\Delta} \leq x+ D'\sqrt{\Delta}
\leq \max _{g\in \mathcal{G}} g(x)
\leq x + D\sqrt{\Delta}.\]}
\end{prop}

\vspace{0.1cm}

\noindent{\bf Proof.} It suffices to prove the claims for the elements in $\mathcal{G}$. 
Indeed, this is obvious for the Lipschitz property, whereas for the boundedness
of the displacements, this is a consequence of the (cocycle) relation
$$gh (x) - x = \big(gh (x) - h (x) \big) + \big( h (x) - x \big).$$

For the Lipschitz property, note that for every element $g \in \mathcal{G}$ 
and every $x < y$, we have
\[ \p (g) \big[ g (y) - g (x) \big]
\leq \int_{\Gamma} \big[ g (y) - g (x) \big] d \p (g) = y-x.\]
Hence,
\begin{equation} \label{eq: Lipschitz} g (y) - g (x) \leq \frac{y-x}{\p (g)},\end{equation}
proving that $g$ has Lipschitz constant at most $1 / \p (g)$. 

\vsp

We next show that the displacements of the elements of \(\Gamma\) are bounded. 
We will in fact prove that for every Lipschitz homeomorphism \(g\) of the real line, the  
value of $\Delta^{g}(x)$ is comparable to $[ g(x)-x ]^2$ up to a multiplicative constant 
depending on  \(L(g)\), the maximum between the Lipschitz constant of \(g\) and that 
of its inverse.

\vspace{0.1cm}

\begin{lem}\label{L: bounded displacements}
\textit{For every Lipschitz orientation-preserving homeomorphism \(g\) 
of the real line and every \(x\in \mathbb{R}\), it holds that 
\[ \frac{1}{2L} \left( g(x) - x \right) ^2 \leq \Delta^g (x)\leq  L \left( g(x) - x \right) ^2, \]
where $L = L(g)$.} 
\end{lem}

\noindent{\bf Proof.} Assume for simplicity that \(g(x) >x\), so that \(g^{-1}(x) <x\). 
Recall that  \(\Delta^g(x) \) is the area of the region \(R\) defined by 
\(\{ (y, z) \! : \ x\leq z\leq g(y) \}\). Observe that this region is contained 
in the rectangle \( [g^{-1}(x), x] \times [x, g(x) ]\), 
hence
\[ \Delta^g (x) \leq (x- g^{-1} (x) ) (g(x) - x) \leq L(g(x) - x ) ^2 .\]
This is the left-side inequality of the statement. To prove the other inequality, 
let \(x'\in \mathbb{R}\) be the point such that \( x- x' = \frac{g(x) - x}{L}\). 
Since \(g\) is \(L\)-Lipschitz, for every \(y\in [x',x]\),  
it holds that \( g(y) \geq g(x) - L(x-y) \). In geometric terms, this 
means that the region \(R\) contains the triangle whose vertices 
are the three points \( (x', x), \ (x,x)\), and \( (x, g(x) )\). This gives the estimate 
\[ \Delta^g (x) \geq \frac{1}{2} (x-x') (g(x) - x) = \frac{1}{2L} \left( g(x)-x\right) ^2. \] 
Analogous arguments apply in the case where $g (x) \leq x$. $\hfill\square$

\vsp\vsp\vsp

We are now in a position to finish the proof of Proposition \ref{P: harmonic are almost-periodic}. 
To do this, set \esp $L:= \max \{1/\p (h) \!: \p (h) > 0 \}$. By \eqref{eq: Lipschitz}, 
this quantity is a Lipschitz constant for each element in $\mathcal G$. 
We claim that, for every $x\in \mathbb R$, 
\[\frac{\Delta}{2L} = \sum _{g\in \mathcal G, \ g(x) > x} \p(g) \frac{\Delta^g(x)}{L} 
\leq \sum _{g\in \mathcal G, \ g(x) > x} \p(g) \left( g(x)-x\right) ^2 \leq L \Delta . \]  
Indeed, this follows as a direct application of Lemma \ref{L: bounded displacements} taking 
into account that \(\Delta^g(x) = \Delta^{g^{-1}} (x)\) and that \( \Delta^g(x) =0\) if \(g(x)=x\). 

If $h \in \mathcal{G}$ is such that $h(x) = \max_{g \in \mathcal{G}} g(x),$ 
then the right-side inequality above shows that 
$$L \Delta 
\geq \sum _{g\in \mathcal G, \ g(x) > x} \p(g) \left( g(x)-x\right)^2 
\geq \p(h)  \left( h(x)-x\right)^2,$$
Thus, 
$$(h(x) - x)^2 \leq \frac{L \Delta}{\p (h)} \leq L^2 \Delta,$$
which yields 
$$h(x) \leq x + L \sqrt{\Delta}.$$
Similarly, the left-side inequality (together with the symmetry of $\p$) yields 
$$\frac{\Delta}{2L} \leq \sum _{g \in \mathcal G, \ g(x) > x} \p(g) \left( g(x)-x \right)^2 
\leq (h(x)-x)^2 \sum_{g \in \mathcal{G}, \ g(x) > x} \p (g) \leq \frac{1}{2} (h(x)-x)^2,$$
hence
$$h(x) \geq x + \sqrt{\frac{\Delta}{L}}.$$
We have thus established that 
$$ x + \sqrt{\frac{\Delta}{L}} \leq 
\max_{g \in \mathcal{G}} g(x)  \leq x + L \sqrt{\Delta}.$$
Analogous considerations show that 
\[ x- L \sqrt{\Delta}\leq \min _{g\in \mathcal G} g(x) \leq x-   \sqrt{\frac{\Delta}{L}} .\]
This concludes the proof of Proposition \ref{P: harmonic are almost-periodic} with 
$D' = \frac{1}{\sqrt{L}}$ and \( D=L\). $\hfill\square$

\vsp\vsp\vsp\vsp

Let $\Gamma$ be a finitely-generated, left-orderable group equipped with a symmetric probability measure \(\p\) 
supported on a finite generating set. The set of all \(\p\)-harmonic actions with \(\Delta = 1\) is a compact subset 
of the space of all actions (equipped with the compact-open topology), and it is invariant under the translation 
flow defined by \eqref{eq: translation flow}. Moreover, the group \(\Gamma\) acts on this space by the formula 
\eqref{eq:def_action}, with each orbit of the flow being preserved. This almost-periodic space is called the 
{\bf \em harmonic space}, and the translation flow acting on it is called  the 
{\bf \em harmonic flow}. These objects play an important role in the proof of the non-orderability of irreducible lattices 
in semi-simple Lie groups in \cite{DH}. Note that this flow may fail to be free, as the next exercise shows.

\vsp

\begin{small}\begin{ejer} \label{ej: harmonic space abelian group}
The goal of the exercise is to compute the harmonic space of $\mathbb{Z}^d$, the free Abelian group of rank \(d\). 

\noindent (i) Show that every bounded $\p$-harmonic function on \(\Z^d\) is constant. (Here, denoting the group law 
on \(\Z^d\) additively, a function $\phi \!: \Z^d \to \mathbb R$ is said to be a {\bf \em $\p$-harmonic function} 
if $\phi (g) = \sum_{h \in \Z^d} \phi (g+ h) \rho(h)$ holds for all $g \in \Z^d$.)

\noindent{\underline{Hint.}} Assume by contradiction that \(\phi : \Z^d \rightarrow \mathbb R\) is a nonconstant bounded 
$\p$-harmonic function. Then there exists \( g_0\in \Z^d \) such that  \(g \mapsto  \phi ( g+g_0) - \phi (g) \) takes a positive 
value. Denote \(M:=\sup_{g\in \Z^d } (\phi (g + g_0)  - \phi (g) ) > 0\), and let \( g_n \in \Z^d \) be such that  
$\phi (g_n + g_0)  - \phi (g_n)$ converges to $M$ as $n$ goes to infinity. 
Up to extracting a subsequence, we can assume that the bounded sequence of functions \( g\mapsto \phi (g+g_n+g_0) - \phi (g+g_n) \)  
pointwise converges to a function \(\psi\). This is a bounded $\p$-harmonic function having a maximum at \(0\) equal to \(M\). 
By a maximum principle argument, \(\psi\) is constant equal to \(M\). Note that \( \phi (g+g_n+g_0) - \phi (g+g_n) \) converges 
to \(M\) for every \(g\in \Z^d \). In particular, \( \phi (g_n+kg_0) - \phi (g_n) \) converges to 
\( kM\)  for every positive integer \(k\). For $k$ large enough, this contradicts the boundedness of $\phi$. 
(Note that commutativity of the underlying group is crucial for this argument.)

\noindent (ii) Show that a Lipschitz $\p$-harmonic function on \(\Z^d\) is affine, that is, the sum of a 
group homomorphism from $\Z^d$ into $(\mathbb R,+)$ and a constant.

\noindent{\underline{Hint.}} Given a Lipschitz $\p$-harmonic function \( \phi: \Z^d \rightarrow \mathbb R\), 
the function $\psi : \Z^d \to \mathbb{R}$ defined by $\psi(g) := \phi (g+ g_0 ) - \phi (g)\) is bounded and harmonic. 

\noindent (iii) Show that every $\p$-harmonic action of $\Z^d$ on the real line is an action by translations. 

\noindent{\underline{Hint.}} Observe that, for each $x \in \mathbb{R}$, the function 
$g \mapsto g (x)$ is $\p$-harmonic.

\noindent (iv) Deduce that the harmonic space of \(\Z^d \) is homeomorphic to a sphere of dimension  
\(d-1\), and that the harmonic flow acts trivially on it.  

\noindent{\underline{Hint.}} For the first statement, denote by \( g_1, \ldots , g_d\) a 
system of generators of $\Z^d $ (that is, a basis of \(\Z^d \) as a \(\mathbb Z\)-modulus). 
Show that the map \(\Phi \mapsto ( \Phi(g_1) (0), \ldots, \Phi (g_d) (0) ) \) 
induces a homeomorphism from the harmonic space of \(\Z^d \) onto an ellipsoid in \(\mathbb R^d\). 
The second statement is immediate.

 \end{ejer}\end{small}


\subsection{Infiniteness of stationary measures}
\label{S:properties of invariant measures}

\hspace{0.45cm} As we have seen, desirable properties hold for actions for which the Lebesgue measure 
is stationary. Here, we start a broad study of general {\em stationary measures}, assuming their existence. 
The universal properties thus obtained will be crucial in \S \ref{section-rec} to establish a fundamental result, 
namely, that there is {\em always} a stationary measure (provided that the probability distribution $\p$ 
is symmetric and supported on a finite generating set). Actually, we will see that this measure 
is unique up to multiplication by a positive constant and can be transformed to the Lebesgue measure 
by a semi-conjugacy, which will thereby allow us to exploit all the properties discussed so far. 

Throughout this section, $\Gamma$ will continue to denote a finitely-generated subgroup of 
$\mathrm{Homeo}_+(\mathbb{R})$ having no global fixed point and endowed with a symmetric 
probability measure $\rho$ supported on a generating set $\mathcal{G}$. Note, however, that 
we will not assume that $\mathcal{G}$ is finite.

Following \cite{DKN2}, let us introduce the Markov process on the line defined by
$$X_x ^n := g_n  \cdots  g_1  (x),$$
where ${\bf g} = (g_n) \in \Gamma^{\mathbb{N}}$ is a family of independent random variables with law $\p$.
(For a general introduction to the theory of Markov processes, we refer the reader
to the very nice book \cite{DY}.)  The transition probabilities of this process are
$$\p_X (x,y) : = \sum_{y = g(x)} \p (g).$$
The associated Markov operator $P = P_X$ acting on the space of bounded
continuous functions $C_b (\mathbb R)$ is given by
\begin{equation} \label{eq:markov}
P \phi (x) : = \mathbb E ( \phi ( X^x_1) ) = \int _{\Gamma} \phi ( g (x) ) \, d\p (g) .
\end{equation}
The iterates of this operator correspond to the operators associated to the
convolutions of $\p$. More precisely, we have $P^n_{\p} = P_{\p^{*n}}$, where 
$\p^{*n} := \p * \p * \cdots * \p$ ($n$ times) and $*$ stands for the 
{\bf{\em convolution}} of probabilities, which is defined by
$$\p_1 * \p_2 (h) := \sum_{fg = h} \rho_1 (f) \hspace{0.05cm} \rho_2 (g).$$
\index{Measure!convolution}

We will still denote by $P$ the dual action on the space of Radon measures on the line. 
Similarly to (\ref{eq:harmonic}), such a measure will be said to be {\bf {\em stationary}} 
if it is $P$-invariant, that is, $P \nu = \nu$. Equivalently,
$$\nu = \sum_{g \in \Gamma} g_* (\nu) \, \p (g).$$
Note that, by definition, an action is harmonic if and only if the Lebesgue measure is stationary.
\index{Measure!stationary}

\vspace{0.1cm}

\begin{lem} \label{L:bi-infiniteness}
\textit{Every nonzero stationary measure $\nu$ on the real line is bi-infinite
 (i.e., $\nu(x,\infty) = \infty$ and $\nu(-\infty,x) =\infty$, for all $x \in \mathbb R$).}
\end{lem}

\noindent{\bf Proof.}
Suppose that there exists $x \!\in\! \mathbb R$ such that $\nu(x,\infty) < \infty$.
Since we are assuming that the $\Gamma$-action has no global fixed point on $\mathbb{R}$,
for every $y\in \mathbb R$, there is an element $g\in \Gamma$ such that $g(x) < y$. As
the support of $\p$ generates $\Gamma$, we can choose $n > 0$ such that
$p^{* n} (g^{-1}) >0$. Then
\[
\nu (y,\infty) \leq \nu (g(x),\infty) \leq \frac{\nu (x,\infty)}{\p^{* n}(g^{-1})} <\infty.
\]
This shows that $\nu(y,\infty ) <\infty$ holds for all $y \in \mathbb R$.

Next, let $\phi \!: \mathbb R \rightarrow (0,\infty)$ be the function defined by
$\phi (x) := \nu(x,\infty)$.
Since $\p$ is symmetric, this function is harmonic; in other words, we have $P \phi = \phi$,
where $P \phi$ is defined by (\ref{eq:markov}). (This definition still makes sense, though $\phi$
is not necessarily continuous.) Fix a real number $C$ satisfying $0< C < \nu(-\infty,\infty)$,
and let $\psi := \max \{ 0, C - \phi \}$. The function $\psi$ is {\bf {\em subharmonic}},
which means that $\psi \leq P \psi$. Moreover, it vanishes on a neighborhood of $-\infty$
and is bounded on a neighborhood of $\infty$.
This implies that $\psi$ is $\nu$-integrable, and since $\int P \psi \, d\nu = \int \psi \, d\nu$,
the function $\psi$ must be $\nu$-a.e $P$-invariant.
Now, a classical lemma from \cite{Garnett} asserts that a measurable function which is in
$\mathcal{L}^1(\nu)$ and $P$-invariant must be a.e. $\Gamma$-invariant
(see Exercise \ref{ex: Garnett} for a schema of proof).
Thus, $\psi$ is constant on almost every orbit. However, this is impossible, since every
orbit intersects every neighborhood of $-\infty$ (where $\psi$ vanishes) and of
$\infty$ (where $\psi$ is positive). This contradiction establishes the lemma.
$\hfill\square$

\vsp

\begin{small}
\begin{ejer} \label{ex: Garnett}
Let $\Gamma$ be a countable group and $\rho$ a measure on $\Gamma$ whose support generates
$\Gamma$. Assume that $\Gamma$ acts on a probability space $(X, \nu)$ by measurable maps and
that $\nu$ is $\rho$-stationary, meaning that
$$\nu = \int_{\Gamma} g_* (\nu) \hspace{0.03cm} d \rho(g).$$
Prove that any function $\phi \in \mathcal{L}^1 (X,\nu)$ that is a.e.
$P$-invariant is a.e. $\Gamma$-invariant.

\noindent{\underline{Hint.}} Pick a constant $C \!\in\! \mathbb R$ and consider the function
$\psi := \max \{\phi, C\}$.  Observe that $\psi$ belongs to $L^1(\nu)$ and satisfies $\psi \leq P \psi$,
and deduce that $\psi$ is $\rho$-harmonic on almost every $\Gamma$-orbit. Conclude by ranging
$C$ over all rational numbers.
\end{ejer}

\begin{ejer} \label{ejer:stat}
Let $\Gamma$ be a finitely-generated subgroup of $\mathrm{Homeo}_+(\mathbb{R})$, and let $\rho$
be a symmetric probability measure on it with generating support. Prove that for all $x \in \mathbb{R}$,
every compact interval $I$, and almost every sequence $(g_n) \in \Gamma^{\mathbb{N}}$, the set of
integers $n$ for which $X^n_x$ belongs to $I$ has density zero, that is,
$$\lim_{k \to \infty} \frac{1}{k} \big| \big\{ n \in \{ 1, \ldots, k \} \!: X^n_x \in I \big\} \big| = 0.$$

\noindent{\underline{Hint.}} Let $\nu_k$ be the measure on the line defined by
$$\nu_k (I) := \frac{1}{k} \sum_{n=1}^k \rho^{* n} \big( \{ h \!: h (x) \in I \} \big).$$
Assuming that the zero-density above doesn't hold, show that, up to a subsequence,
$\nu_k$ converges to a nonzero, finite, stationary measure, thus contradicting
Lemma \ref{L:bi-infiniteness}.
\end{ejer}

\begin{ejer} \label{e:irreducible}
Prove that no nontrivial group action on the line satisfying (\ref{eq:harmonic})
for all $x,y$ in $\mathbb{R}$ can have a global fixed point.

\noindent{\underline{Hint.}} By definition, the Lebesgue measure is $P$-invariant for a $\p$-harmonic
action. Apply Lemma \ref{L:bi-infiniteness} to the restriction of the action to a connected component
of the complement of the set of global fixed points.
\end{ejer}\end{small}

\vsp

\subsection{Recurrence}
\label{section-rec}

\hspace{0.45cm} As in previous sections, we continue to consider a symmetric
probability measure $\rho$ on a group $\Gamma$ acting on the real line without global fixed points.
We also assume that the support $\mathcal{G}$ of $\rho$ generates $\Gamma$. We start with a 
result concerning {\em oscillation} of random orbits; more precisely, it asserts that almost every 
random orbit escapes to infinity in both directions. Note that this result does not assume 
that $\mathcal{G}$ is finite. 

\vsp\vsp

\begin{prop} \label{P:oscillation}
\textit{For every $x \in \mathbb R$, almost surely we have
\[ \limsup_{n\rightarrow \infty} X_x^n = +\infty\ \ \ \mathrm{and}\ \ \
\liminf_{n\rightarrow \infty} X^n_x = -\infty.\]}
\end{prop}

\noindent{\bf Proof.} Denote $\mathbb{P} := \rho^{\mathbb{N}}$. 
Given points $C$ and $x$ on the real line, let
\[  p_C (x) := \mathbb P \Big[ \limsup_{n\rightarrow \infty} X_x ^n > C \Big] .\]
Since $\Gamma$ acts by orientation-preserving homeomorphisms, for all $x \leq y$, we have
\[
\Big\{ (g_n) \in \Gamma^{\mathbb N} \!:  \limsup_{n\rightarrow \infty} X_x^n > C \Big\}
\subset
\Big\{(g_n)\in \Gamma^{\mathbb N} \!: \limsup_{n\rightarrow \infty} X_y^n >  C \Big\}.\]
In particular, $p_C (x) \leq p_C (y)$, that is, $p_C$ is non-decreasing.
Moreover, since $p_C$ is the probability of the tail event
\[ \Big[ \limsup_{n\rightarrow \infty} X_x^n > C \Big] \]
and $X$ is a Markov chain, $p_C$ is a \textit{harmonic} function,
that is, for every $x\in \mathbb R$ and every integer $n\geq 0$,
\[ p_C (x) = \sum_{g \in \Gamma} p_C \big( g(x) \big) \hspace{0.01cm}
\rho^{* n} (g) = \mathbb E \big( p_C (X_x^n) \big).\]

Now, we would like to see $p_C$ as the distribution function of a \textit{finite} measure on
the line. However, this is only possible when $p_C$ is continuous on the right, which is {\em a
priori} not necessarily the case. We are hence led to consider the right-continuous function
\[ \overline{p_C}(x): = \lim _{y\rightarrow x, \ y>x} p_C (y).\]
This function is still non-decreasing. Therefore, there exists a finite
measure $\nu$ on $\mathbb R$ such that for all $x < y$,
\[  \nu (x,y]  \hspace{0.1cm} = \hspace{0.18cm}
\overline{p_C} (y) - \overline{p_C} (x).\]
Since $p_C$ is harmonic and $\Gamma$ acts by homeomorphisms, $\overline{p_C}$
is also harmonic. Since $\rho$ is symmetric, this yields that $\nu$ is $P$-invariant.
Now recall that Lemma~\ref{L:bi-infiniteness} implies that any $P$-invariant
finite measure identically vanishes (see also \cite[Proposition 5.7]{DKN}). Therefore,
$\overline{p_C}$ is constant; in particular its value does not depend on
the starting point $x$. The 0-1 law then allows to conclude that (for any
fixed~$C$) either $p_C \equiv 0$ or $p_C \equiv 1$.

Let us now show that $p_C$ is identically equal to~$1$ for each $C$. To do this, fix any
$x_0 > C$. Since for every $g\in \mathrm{Homeo}_+(\mathbb R)$, we have either $g(x_0)\ge x_0$
or $g^{-1}(x_0) \ge x_0$, the symmetry of~$\p$ yields that $X_{x_0}^n \ge x_0$ holds with
probability at least~$1/2$, for all $n \in \mathbb{N}$. It is then easy to see that
$$
p_L = p_C (x_0) \ge \mathbb{P} \Big[ \limsup_{n \to \infty} X^n_{x_0} \ge x_0 \Big] \ge 1/2.
$$
As we have already shown that $p_C$ equals ~$0$ or~$1$, this implies that
$p_L$ is identically equal to~$1$.

The latter means that for every $x\in \mathbb R$, the equality
\[  \limsup_{n\rightarrow \infty} X_x^n =+ \infty \]
holds almost surely.
Analogously, for every $x \in \mathbb R$, almost surely we have
\[  \liminf _{n\rightarrow \infty} X_x^n = -\infty. \]
This completes the proof of the proposition.
$\hfill\square$

\vspace{0.4cm}

We are now ready to prove the main result of this section, namely the {\em recurrence}
of the Markov process under the extra hypothesis that $\mathcal{G}$ is finite. 

\vspace{0.1cm}

\begin{cor}
\label{C:recurrence}
\textit{There exists a compact interval $K$ such that, for every $x \in \mathbb{R}$,
almost surely the sequence $(X_x^n)$ intersects $K$ infinitely many times.}
\end{cor}

\noindent{\bf Proof.}
Consider a closed interval $K$ as in the proof of Lemma \ref{L: atomic part 2}, that is, $K \!=\! [A, B]$,
where $A<B$ are such that for every $g$ of the support of $\p$, we have $g(A) < B$. (Recall that
$\rho$ is finitely-supported.) By Proposition \ref{P:oscillation}, for every $x\in \mathbb R$, almost
surely the sequence $(X_x^n)$ will pass from $(-\infty, A]$ to $[B,+\infty)$ infinitely many times. 
The desired conclusion follows from the observation that the choice of $A$ and $B$
imply that every time this happens, $(X_x^n)$ must cross the interval $K$.
$\hfill\square$

\vspace{0.4cm}

\noindent{\bf On left-orders that are generic with respect to a stationary measure.} Given a probability on 
a left-orderable group $\Gamma$ upported on a generating set, we can also consider stationary
probability measures for the action of $\Gamma$ on its space of left-orders (see Example \ref{ex: Garnett}).
By this, we mean a probability measure $\mu$ on $\mathcal{LO}(\Gamma)$ such that
\begin{equation}\label{stat-LO}
\mu = \sum_{g \in \Gamma} g_* (\mu) \rho (g).
\end{equation}
Since $\mathcal{LO}(\Gamma)$ is compact, such a probability measure $\mu$ always exists. This follows
from a direct application of either Kakutani's fixed point theorem \cite{conway} or the Bogoliubov-Krylov 
averaging procedure \cite{sinai}. (Note that we do not require $\rho$ to be finitely supported for this.) 
It seems quite interesting to study the relation of $\mu$ with the algebraic properties of $\Gamma$ as
well as its dependence on $\rho$. Below we give two examples on this.

\begin{small}\begin{ex}\label{ex:new-proof}
We next give still another proof of Theorem \ref{linnell-general} for finitely-generated groups.
To do this, fix $\rho$ and $\mu$ as above. By a standard argument of desintegration into ergodic 
components \cite{OV}, we can assume that $\mu$ is {\em ergodic}, in the sense that it cannot  be written 
as a nontrivial convex combination of two different stationary probability measures. We have two possibilities:

\vsp

\noindent \underline{Case (i).} The measure $\mu$ has an atom.

\vsp

If $\preceq$ is an atom of maximal $\mu$-measure, then (\ref{stat-LO}) easily implies that
its orbit must be finite. (Actually, by ergodicity, this orbit necessarily coincides with the support of $\mu$.)
In particular, $\preceq$ is right-recurrent, hence Conradian (see \S \ref{super-witte}). Thus, if $\Gamma$ 
has infinitely many left-orders, then Proposition \ref{simple-alcanza} implies that $\mathcal{LO}(\Gamma)$
is uncountable, as desired.

\vsp

\noindent \underline{Case (ii).} The measure $\mu$ is non-atomic.

\vsp

By ergodicity, for almost every pair 
$\big(\! \preceq,(g_n) \big)$ in $\mathcal{LO}(\Gamma) \times \Gamma^{\mathbb{N}}$
(endowed with the measure $\mu \times \rho^{\mathbb{N}}$), the sequence
$( \preceq_{\sigma^n({\bf g})},\sigma^n({\bf g}))$ is dense in
$supp(\mu) \times \Gamma^{\mathbb{N}}$,
where $\sigma$ stands for the shift $\sigma ((g_n)) := (g_{n+1})$.
Let us fix such a pair $\big( \! \preceq, (g_n) \big)$, and let $(U_k)$ be a sequence of open subsets
of positive $\rho$-measure in  $\mathcal{LO}(\Gamma)$, none of which contains $\preceq$, but
which do converge to $\preceq$. For each $k$, there exists $n(k) \!\in\! \mathbb{N}$ such that
$\preceq_{\sigma^{n(k)} ({\bf g})}$ belongs to $U_k$. Hence, $\preceq_{\sigma^{n(k)} ({\bf g})}$
converges to $\preceq$, with $\preceq_{\sigma^{n(k)}({\bf g})}$ being distinct from $\preceq$
for all $k$. As a consequence, the closure of the orbit of $\preceq$ under the action of $\Gamma$
is a totally disconnected compact metric space with no isolated point, that is, a Cantor
set. In particular, $\mathcal{LO}(\Gamma)$ is uncountable.
\end{ex}

\index{Conradian!soul}
\begin{rem} Note that the approximation by conjugates in the example above is essentially 
different from that of the original proof of Theorem \ref{linnell-general}. Namely, in \S \ref{section-soul}, 
the conjugating elements are positive but ``small' (as small as possible 
outside the Conradian soul). In the proof above, the conjugating elements 
are ``random''. By Exercise \ref{ejer:stat}, if $\rho$ symmetric, then these elements  are mostly ``near 
the infinite'' (either ``very positive'' or ``very negative''), despite the recurrence of the associated random 
walk on the line. 
\end{rem}

\begin{ex} \label{ex:BS}
According to Example \ref{solo-cuatro}, for each integer $\ell \geq 2$,
the Baumslag-Solitar group $B(1,\ell) = \langle a, b \!:  a b a^{-1} = b^{\ell} \rangle$
admits four Conradian orders, which are actually bi-invariant and come from
the exact sequence
$$0 \longrightarrow \mathbb{Z} \Big[ \frac{1}{\ell} \Big] \longrightarrow B(1,\ell)
\longrightarrow \mathbb{Z} \longrightarrow 0.$$
We claim that, although $\mathcal{LO}(B(1,\ell))$ is a Cantor set
(see \S \ref{section finite-rank-solvable}), for every symmetric probability distribution $\rho$
on $B(1,\ell)$ with finite generating support, every stationary probability measure $\mu$ 
on $\mathcal{LO}(B(1,\ell))$ is supported on these four points. Indeed, we proved in
\S \ref{section finite-rank-solvable} that for every left-order $\preceq$ on $BS(1,\ell)$
that is not bi-invariant, the associated dynamical realization is semiconjugate to a non-Abelian subgroup
of the affine group. In particular, there exist elements whose sets of fixed points are bounded and for
which $-\infty$ and $+\infty$ are topologically repelling fixed points. Let $g$ be such an element
(actually, such a $g$ can be taken as the image of $a^{-1}$), and
denote by $Fix (g)$ its set of fixed points. Let $f_1,f_2$ be in the realization of $B(1,\ell)$ so that $f_1$
(resp. $f_2$) sends the leftmost (resp. the rightmost) fixed point of $g$ to the right (resp.
left) of $0 \!=\! t(id)$. Denote $g_1 := f_1 g f_1^{-1}$ and $g_2 := f_2 g f_2^{-1}$.
If we identify elements in $BS (1, \ell)$ with their realizations, we have
$g_1 \succ id$ and $g_2 \prec id$. Moreover, $h g_i h^{-1} \succ id$ holds
for both $i \!=\! 1$ and $i\!=\!2$ provided $h$ is sufficiently large (say, larger than a
certain element $h_+$). Similarly, $h g_i h^{-1} \prec id$ holds for $i=1$ and $i=2$
provided $h$ is smaller than a certain element $h_{-}$; see Figure 21 below.

\vspace{0.05cm}


\beginpicture

\setcoordinatesystem units <1.05cm,1.05cm>

\putrule from -6.6 0 to 6.6 0
\putrule from 0 0 to 0 3
\putrule from 0 -3 to 0 -0.8


\plot
-5  0.5
5.5  2.2 /

\plot
-5 0.51
5.5 2.21 /

\plot
-5 0.52
5.5 2.22 /

\plot
-5 0.49
5.5 2.19 /

\plot
-5 0.48
5.5 2.18 /


\plot
-5  -2.7
5.5  -1 /

\plot
-5 -2.71
5.5 -1.01 /

\plot
-5 -2.72
5.5 -1.02 /

\plot
-5 -2.69
5.5 -0.99 /

\plot
-5 -2.68
5.5 -0.98 /


\plot
-3 -3
-0.75 -0.75 /

\plot
0 0
3 3 /

\put{$\bullet$} at 0 0
\put{$\bullet$} at 2.5 0
\put{$\bullet$} at -3 0

\put{$t(id) = 0$} at  0 -0.35
\put{$t(h_{-})$} at -3 -0.35
\put{$t(h_{+})$} at  2.5 -0.35
\put{$g_1 = f_1 g f_1^{-1}$} at 4.8 1.75
\put{$g_2 = f_2 g f_2^{-1}$} at 4.8 -1.5

\setdots
\putrule from 1.55 0 to 1.55 1.55
\putrule from -2.25 0 to -2.25 -2.25

\put{} at -7.2 0
\put{Figure 21: The elements $g_1,g_2,h_-$ and $h_+$.} at 0 -3.7
\endpicture


\vspace{0.75cm}

Assume $\mu$ is a stationary probability measure on $\mathcal{LO} (BS(1,\ell))$
that is not fully supported on the four bi-orders. Then any ergodic component of this
measure outside these bi-orders is still stationary, and supported on the complement
of the bi-orders. For simplicity, we still denote this ergodic component by $\mu$. Let $\preceq$
a point in the support of $\mu$. If we perform the construction of the elements
$g_1,g_2,h_-,h_+$ above, then the measure of the open neighborhood \esp
$V_{g_1} \cap V_{g_2^{-1}} = \big\{ \! \preceq' : \esp \esp g_1 \succ' id, \esp g_2 \prec' id \big\}$
\esp of $\preceq$ must be positive, say equal to $\kappa > 0$. A direct application of Birkhoff's 
ergodic theorem then shows that, for a generic random path \esp $(h_n) \in B(1,\ell)^{\mathbb{N}}$,
the set of integers $n$ for which $\preceq_{_{X^n}}$ lies in $V_{g_1} \cap V_{g_2^{-1}}$ 
has density $\kappa$, where $X^n := h_n \cdots h_1$. Nevertheless, among these
integers $n$, \esp with density 1 we have \esp either \esp $X^n \prec h_{-}$ \esp
or \esp $X^n \succ h_{+}$ (see Exercise \ref{ejer:stat}), which is a contradiction.
\end{ex}
\end{small}


\subsection{Further properties of stationary measures}


In this section, we start by studying in detail the case where a group action on the line admits discrete orbits. Obviously, such 
an orbit supports an invariant Radon measure, namely the counting measure. In particular, this measure is $P$-invariant. The 
next two lemmas show that if there exists a discrete orbit, then all $P$-invariant Radon measures lie in the convex closure 
of the set of counting measures along discrete orbits. To do this, we will use the following useful combinatorial lemma.

\begin{lem} \label{L: useful lemma} \textit{Let $\Gamma$ be a group endowed with a symmetric probability measure $\rho$ 
supported on a finite generating set $\mathcal{G}$. Let 
\(\Gamma_*\subsetneq \Gamma\) be a strict subgroup, \(\preceq\) a \(\Gamma\)-invariant 
total order on \(Y=\Gamma / \Gamma_*\), and \(\varphi : Y\rightarrow (0, +\infty) \) a \(\rho\)-harmonic 
function, in the sense that \(\varphi (y) = \int \varphi(g (y)) 
d \rho( g )\) holds for every \(y\in Y\). Assume that for each couple of elements \( y_1\prec y_2 \) in \(Y\), we have 
\[ \sum _{y_1 \preceq y \prec y_2} \varphi (y) < +\infty. \] 
Then \(\varphi \) is constant, the subgroup \( \Gamma_*\) is normal in $\Gamma$, and \( \Gamma/\Gamma_*\) is infinite cyclic. } 
\end{lem}

\noindent{\bf Proof.} We proceed to construct a \(\rho\)-harmonic action from the data. We define a 
collection of nonempty open intervals $  I_y := \{ (a_y, b_y)\}_{y\in Y} \} $ in the real line by the formulae:
$$a_{y_2} - a_{y_1}  := \sum _{y_1\preceq y \prec y_2} \varphi (y) \quad \text{for} \quad y_1\prec y_2, 
\quad \mbox{and} \quad b_y - a_y := \varphi (y).$$
This family is uniquely defined up to translations. Moreover, the intervals \(I_y\) are two-by-two disjoint, 
and their union \(I=\bigcup_{y\in Y}  I_y\) is a subset of the real line of full Lebesgue measure. In particular, 
the complement of \(I\) is a closed set of empty interior. Note that, by construction, the order of the intervals 
\(I_y \) is the one induced by \(\preceq\). 

Define an action of the group \(\Gamma\) on \(I\), by imposing that the element \( g \) maps \(I_y\) onto \(I_{g (y)}\) 
by the orientation-preserving affine map 
$$
x \mapsto \frac{\varphi (g (y))}{\varphi (y)} (x-a_y) + a_{\varphi(y)}.
$$ 
Since \(\Gamma\) preserves the ordering of the intervals \(I_y\)'s, it acts by increasing homeomorphisms of \(I\). 
Moreover, as the complement of this latter has zero Lebesgue measure, the action extends to an action by 
homeomorphisms of the line. Furthermore, since \(\varphi\) is a \(\rho\)-harmonic function and \( \rho\) is 
symmetric, this action is \(\rho\)-harmonic. 

Proposition \ref{P: harmonic are almost-periodic} then shows that 
\begin{equation}\label{eq:boundedd}
\sup_{ x\in \mathbb R, \ g\in \mathcal{G}} (g(x) - x) <+\infty .
\end{equation}
Since the subgroup \(\Gamma_*\) is strict, for every \( y\in Y\) there exists \(g\in \mathcal G\) such that \( y \prec g(y)\). 
As the interval \(I_{g(y)}\) is on the right of \( I_y\), this yields \( \varphi(y) = b_y-a_y\leq g(a_y)-a_y\). From (\ref{eq:boundedd}),  
we thus infer that \(\varphi \) is a bounded harmonic function.  

Let \(K\subset \mathbb R\) be a recurrence interval 
for the \(\Gamma\)-action given by Corollary \ref{C:recurrence}, 
and let \(Y_*\subset Y\) the set of \(y \in Y \) such that \(I_y \cap K\neq \emptyset\). The set 
\(Y_*\) is a recurrent subset, in the sense that for every \(y\in Y\) and \(\rho^{\mathbb N}\)-a.e. 
\( {\bf g} = (g_n) \in \Gamma ^{\mathbb N}\), there exists \(n\in \mathbb N\) such that 
\(g_n \ldots g_0 (y) \in Y_*\). We denote by \( n_{Y_*} (y,{\bf g}) \) the infimum of these integers. The function 
\(n_{Y_*}\) is a {\bf \em stopping time}, namely, \(n_{Y_*} (y, {\bf g'})\) is equal to \( n_{Y_*} (y,{\bf g}) \) if 
\(g_n'=g_n\) for every \( n\leq n_{Y_*} (y,{\bf g})\)). 

Since the union of the intervals \(I_y\)'s for \(y\in Y_*\) is bounded in \(\mathbb R\), we have 
\[ \sum _{y\in Y_*} \varphi (y) <+\infty.\]
In particular, there exists \(y_{\max}\in Y_*\) such that \(\varphi (y_{max})\) is maximal. 
Applying the martingale convergence theorem to the stopping time \(n_{Y_*}\), we conclude that, denoting \( l_n ({\bf g}) := g_n \ldots g_0\), 
we have 
\[ \varphi (y_{max}) = \int _{\Gamma^{\mathbb N} } \varphi (l_{n_{Y_*}(y,{\bf g})}(y_{max})) \ \rho ^{\mathbb N} (d{\bf g}) .\]
Since \(l_{n_{Y_*}(y,{\bf g})} (y_{max})\) belongs to \(Y_*\), we have  \(\varphi (l_{n_{Y_*}(y,{\bf g})} (y_{max}))\leq \varphi(y_{\max}) \). 
Therefore, almost surely, it holds that 
\[ \varphi (l_{n_{Y_*}(y,{\bf g})} (y_{max}))= \varphi(y_{max}) .\]
We claim that this implies that \(\varphi\) is constant on \(Y_*\). 
Indeed, given any \(y_*\in Y_*\), there exist elements \( g_0, \ldots, g_n\in \mathcal G\) such that \(y_*=g_n\ldots g_0 (y_{max})\). 
Let \( 0=n_0\leq n_1\leq \ldots \) be the integers between \(0\) and \(n\) for which \( g_{n_{k}}\ldots g_0 (y_{max})\) belong to 
\( Y_*\). Since \(\mathcal G\) is contained in the support of \(\rho\). Then the arguments above show recursively that 
\( \varphi (  g_{n_{k}}\ldots g_0 (y_{max}) ) = \varphi (y_{max})\), and thus \(\varphi (y_*)= \varphi (g_n \ldots g_0 (y_{max}) ) = \varphi (y_{max})\). 

The conclusion above was obtained for any recurrence interval $K$. Thus, by considering an exhaustion of the real line by 
recurrence intervals, we deduce that the function \(\varphi \) is constant on \(Y\). In particular, this implies that between two 
points \( y_1\prec y_2\) in \(Y\), there are only a finite number of points. As a consequence, the ordered space \((Y, \preceq)\) 
is isomorphic to \( (\mathbb Z, \leq)\). Since the group of automorphisms of \((\mathbb Z, \leq)\) is a cyclic group, the lemma follows.
$\hfill\square$

\vspace{0.15cm}

\begin{lem} \label{L: atomic part 1}
\textit{Let $\Gamma$ be a finitely-generated subgroup of $\mathrm{Homeo}_+(\mathbb{R})$ having 
no global fixed point, which is endowed with a symmetric probability measure $\rho$ supported on a finite generating set 
$\mathcal{G}$. Let $\nu$ be a Radon measure on the line that is stationary for the $\Gamma$-action. If there is a discrete 
orbit, then $\nu$ is supported on the union of discrete orbits, and is totally invariant.}
\end{lem}

\noindent{\bf Proof.} If there is a discrete orbit $\mathcal O$, then $\Gamma$ acts on it
by translating its points. Thus, the normal subgroup $\Gamma_*$ formed by
the elements acting trivially on $\mathcal O$ is recurrent, by Polya's classical theorem
\cite{Polya} (see also Corollary \ref{C:recurrence} for an alternative proof of this fact).
Let $\p_*$ be the (symmetric) measure on $\Gamma_*$ obtained by {\bf{\em balayage}}
of $\p$ to $\Gamma_*$. More precisely, we consider the random walk on $\Gamma$ with 
transition probabilities given by \( p (g,h) =\rho (gh^{-1})  \) starting at the neutral element of $\Gamma$, 
and we stop it at the first moment where it visits $\Gamma _*$. This yields a random variable with 
values in $\Gamma_*$ whose distribution is $\p _*$ (using the notation of the proof of Lemma 
\ref{L: useful lemma}, \(\rho_*\) is the distribution of the random variable 
\({\bf g}\in \Gamma^{\mathbb N}\mapsto l_{n_{\Gamma_*} (e,{\bf g})} \in \Gamma_* \)). 

We claim that the restriction \(\nu_*\) of $\nu$ to a connected component $C$ of
$\mathbb R \setminus \mathcal O$ is a (finite) measure that is \(\rho_*\)-stationary. 
To show this, note that for every Borel subset \(A\subset C\), the function $\phi : \Gamma \to \mathbb{R}$ 
given by $\phi (g) := \nu (g^{-1} A) \in [0,+\infty) $ is {\bf \em right \(\rho\)-harmonic}, which means that 
for every \(g \in \Gamma\), 
\[  \phi (g) = \sum_{f \in \Gamma} \rho (f) \phi (g f) . \]
If we replace \(A\) by \(C\) itself, the function becomes left \(\Gamma_*\)-invariant, namely, 
\( \phi ( g f) = \phi (f) \) for every \( g\in \Gamma_*\) and $f \in \Gamma$.
It hence induces a nonnegative right \(\rho\)-harmonic function on the infinite cyclic group $\Gamma / \Gamma^*$, 
which is necessarily constant; see Exercise~\ref{ej: non negative harmonic functions on Z} below. 
In particular, the \(\nu\)-measures of the images of \(C\) by the elements of \(\Gamma\) 
are bounded, and in particular those of \(A\). Therefore, the claim is a consequence of the martingale convergence 
theorem applied to the stopping time \( n_{\Gamma^*}(e, \cdot)\) and to the bounded harmonic function $\phi$. 

The claim above having been proved for {\em every} connected component $C$ of $\mathbb R \setminus \mathcal O$, 
it follows from Lemma~\ref{L:bi-infiniteness} that \(\nu_*\)  is supported on $\mathrm{Fix}(\Gamma_*) \cap \overline{C}$, 
the set of global fixed points for the group $\Gamma_*$ contained in the closure of $C$. 
(Note that the proof of Lemma \ref{L:bi-infiniteness} does not use finiteness of the support of the underlying probability distribution.) 
As a consequence,  the support of $\nu$ consists of a union of discrete orbits, each one being isomorphic as an ordered space 
to \((\mathbb Z, \leq)\). To see that $\nu$ is invariant, note that for every bounded Borel subset \(B\subset \mathbb R\), the 
function $g\in \Gamma \mapsto \nu (g^{-1} (B))$,  is \(\rho\)-harmonic invariant by \(\Gamma_*\). It hence induces an 
harmonic function on the quotient $\Gamma / \Gamma_*$, and since this is isomorphic to $\mathbb Z$, it must be constant. 
$\hfill\square$

\vsp\vsp

\begin{small}\begin{ejer}\label{ej: non negative harmonic functions on Z}
By developing the items below, show that for any symmetric probability $\p$ with finite generating support on $\Z$, 
every non-negative \(\rho\)-harmonic function is constant. 

\noindent (i) Given a $\rho$-harmonic function $\phi$, write it in the form of a sequence: $a_n : = \phi(n)$. Check that 
$\rho$-harmonicity translates into the following recurrence relation for $(a_n)$: For each $n \in \mathbb{Z}$, 
$$\p (k) [ a_{n+k}+a_{n-k} ] + \p (k-1) [ a_{n+k-1}+a_{n-k+1} ] + \ldots + \p (1) [ a_{n+1} + a_{n-1}  ] + (\p (0) - 1) a_n = 0.$$
where $k \in \mathbb{Z}$ is the maximum element in the support of $\p$. (Note that, by symmetry, $k$ is strictly positive.)

\noindent (ii) Check that the characteristic polynomial $P$ of this recurrence relation 
can be written as $P = x^k Q$, where 
$$Q (x) = \p (k) [x^k + x^{-k} ] + \ldots + \p (1) [x+x^{-1}] + (\p (0) - 1).$$
Using the classical theory of recurring sequences, conclude that $a_n$ can be written in the form
$$a_n = \sum_{\lambda} P_{\lambda} (n) \, \lambda^n,$$
where the sum ranges over the roots $\lambda$ of $P$ and each $P_{\lambda}$ is a polynomial whose degree is 
equal to the multiplicity of $\lambda$ as a root of $P$ minus 1.

\noindent (iii) Check that $0$ is not a root of $P$, and prove that $P$ has no positive real root other than 1.

\noindent{\underbar{Hint.}}  Use the fact that for each positive real number $x$ different from $1$ and all integers $m \geq 1$, 
it holds that $x^m + x^{-m} > 2$, hence $Q(x) > 0$.

\noindent (iv) Check that $1$ is a root of $P$ with multiplicity 2.

\noindent (v) Show that no root of $P$ other than 1 can have norm 1.

\noindent{\underbar{Hint.}} Use that for all $z = e^{i \theta} \neq 1$ of norm 1, the expression 
$z^m + z^{-m} = 2 \cos (m \theta)$ lies in $[-2,2[$, hence $Q (z) < 0$.

\noindent (vi) Show that $P$ has no root with a modulus greater than $1$. 

\noindent{\underbar{Hint.}} Assume otherwise and let $r > 1$ be the maximum modulus of a root of $P$. Let 
$d$ be the maximal degree of the polynomials $P_{\lambda}$ associated with the roots $\lambda$ of $P$ 
of norm $r$. Finally, let $\lambda_1, \ldots, \lambda_{\ell}$ be the roots $\lambda$ of norm $r$ whose 
associated polynomials have degree $d$. Writing $\lambda_j = r e^{i \theta_j}$ and 
$P_{\lambda_j} (x) = c_j x^d + \ldots$ (where the dots stand for lower-order 
terms), check that, as $n$ goes to infinity, 
$$a_n \sim n^d r^n R(n), \quad \mbox{where } \, R(n) = \sum_{j=1}^{\ell} c_j e^{i \theta_j n}.$$
Conclude by contradiction by observing that the trigonometric polynomial $R$ takes negative values 
for  arbitrarily large integers $n$. 

\noindent (vii) Use a similar argument to the one above (with $n$ going to $-\infty$ this time) 
to show that $P$ has no root with a modulus less than $1$. 

\noindent (viii) Conclude that there exist constants $c_1,c_2$ such that $a_n = c_1n+c_2$ holds for all $n$. Using that 
$a_n$ is non-negative for all $n$, conclude that $c_1 = 0$, and therefore $a_n$ is constant.  
\end{ejer}

\begin{ejer} Show that the claim of the preceding exercise does no longer hold if $\p$ fails to be symmetric.
\end{ejer}
\end{small}

\vsp

\begin{small}\begin{ejer}
Assume that a $\Gamma$-action on the line is harmonic and admits a discrete orbit.
Prove that $\Gamma$ coincides with the (cyclic) group generated by
a translation of the line.
\end{ejer}\end{small}

\vsp\vsp

\begin{lem} \label{L: atomic part 2}
\textit{Let $\Gamma$ be a finitely-generated subgroup of $\mathrm{Homeo}_+(\mathbb{R})$ having 
no global fixed point, which is endowed with a symmetric probability measure $\rho$ supported on a finite generating 
set $\mathcal{G}$. Let $\nu$ be a Radon measure on the line that is stationary for the $\Gamma$-action.
If the atomic part of $\nu$ is nontrivial, then $\nu$ is invariant and supported on a union of discrete orbits.}
\end{lem}

\noindent{\bf Proof.} Let \(x\in \mathbb R\) be a point such that \( \nu (\{x\})>0\). Let $\Gamma_*$ be 
the stabilizer of $y$ in $\Gamma$, and let $Y \sim \Gamma / \Gamma_*$ denote the orbit of $x$. 
The function \(\varphi : Y\rightarrow [0, +\infty) \)  defined by \( \varphi (y):= \nu (\{y\})\)  satisfies the assumptions 
of Lemma \ref{L: useful lemma}. As a consequence, the restriction of the measure \(\nu\) to $Y$ is invariant, 
which implies that the orbit of \(x\) is discrete. The conclusion then follows from Lemma \ref{L: atomic part 1}.  
$\hfill\square$

\vspace{0.44cm}

We next turn to the case where $\Gamma$ has no discrete orbits. 
Recall from Lemma \ref{minimal-set} that, if $\Gamma$ is finitely-generated, then 
there is a unique nonempty minimal invariant closed set for the action. In the sequel, 
we denote this set by $ M$.

\vspace{0.1cm}

\begin{lem} \label{L: support}
\textit{Let $\Gamma$ be a finitely-generated subgroup of $\mathrm{Homeo}_+(\mathbb{R})$ having 
no global fixed point, which is endowed with a symmetric probability measure $\rho$ supported on a finite 
generating set $\mathcal{G}$. Assume that $\Gamma$ does not have any discrete orbit on the real line.
Then any stationary measure is supported on the minimal set $M$.}
\end{lem}

\noindent{\bf Proof.} Let $\nu$ be a stationary measure on the real line. 
For all $h \in \mathcal{G}$, we have
\[h_* (\nu) \leq \frac{1} {\p (h)} \sum_{g \in \Gamma} g_* (\nu) \, \p(g)
= \frac{\nu}{\p (h)}.\]
This obviously implies that  $\nu$ quasi-invariant by $\Gamma$. As a consequence, 
the support of $\nu$ is a closed $\Gamma$-invariant subset of the line, hence it
contains $M$. Thus, it suffices to show that $\nu$ does not charge any connected
component of the complement $M^c$.

Assume $M^c$ is nonempty, and collapse each of its connected components
to a point. We thus obtain a topological line carrying a $\Gamma$-action for which all
orbits are dense. The stationary measure $\nu$ can be pushed to a stationary measure
$\overline{\nu}$ for this new action. If a component of $M^c$ had a positive
$\nu$-measure, then $\overline{\nu}$ would have atoms. However, by Lemma~\ref{L: atomic part 2},  
this contradicts the minimality of the $\Gamma$-action obtained after collapsing.
$\hfill\square$

\vsp

\begin{small}
\begin{ejer}
Prove that if the action of $\Gamma$ is harmonic, then $\Gamma$ either is a cyclic group of translations
or acts minimally on the real line.
\end{ejer}
\end{small}


\subsection{Existence of stationary measures}
\label{sec:stationary}

\hspace{0.45cm} Using the recurrence result of the preceding section,
we can now establish the existence of a $P$-invariant Radon
measure via a quite long but standard argument.

\vspace{0.18cm}

\begin{thm}
\label{T:existence of invariant measure}
\textit{Let $\Gamma$ be a finitely-generated subgroup of $\mathrm{Homeo}_+(\mathbb{R})$
endowed with a symmetric probability measure \(\rho\). If the support of $\rho$ is finite and generates 
$\Gamma$, then there exists a (nonzero) $\rho$-stationary measure on the real line.}
\end{thm}

\noindent {\bf Proof.} Fix a continuous compactly-supported function 
$\xi \!: \mathbb R \to [0,1]$ such that $\xi \equiv 1$ on the compact recurrence interval $K$ 
(see Corollary \ref{C:recurrence}). For any initial point $x$, let us stop the process $X^n_x$ at a
\emph{random} stopping time $T$ chosen in a Markovian way so that, for all $n \in \mathbb{N}$,
$$\mathbb P \big[ T=n+1 \mid T\ge n \big] = \xi (X^{n+1}_x).$$
(Here, $T = T({\bf g})$, where ${\bf g} = (g_i)_{i \in \mathbb{N}}$.)
In other words, after each step of the initial random walk arriving to a point $y=X^{n+1}_x$,
we stop with probability $\xi(y)$, and we continue the compositions with probability $1 - \xi(y)$.

Denote by $Y_x$ the random stopping point $X_x^T$, and consider its distribution $\p_x$
(note that $T$ is almost-surely finite since the process $X_x^n$ almost surely visits $K$ and
$\xi \equiv 1$ on $K$). Due to the continuity of $\xi$, the measure $\p_x$ on $\mathbb R$
depends continuously (in the weak topology) on~$x$. Therefore, the corresponding diffusion
operator $P_{\xi}$ defined by
\[ P_{\xi}(\phi)(x) = \mathbb E \big( \phi (Y_x) \big) = \int_{\mathbb R} \phi (y) \, d\p_x (y) \]
acts on the space of bounded continuous functions on~$\mathbb R$, and hence it acts by duality
on the space of probability measures on~$\mathbb R$. Note that for any such probability measure,
its image under $P_{\xi}$ is supported on $\hat{K} := \text{supp} (\xi)$. Thus, by applying the
Bogolyubov-Krylov procedure of time averaging \cite{sinai}, 
we see that there exists a $P_{\xi}$-invariant probability measure~$\nu_0$.

In order to construct a Radon measure that is stationary for the initial process, we proceed as follows: 
For each $x \!\in\! \mathbb R$, let us consider the sum of the Dirac measures supported on its random
trajectory before the stopping time~$T$. In other words, we consider the  random measure
$$m_x ({\bf g}) := \sum_{j=0}^{T(w)-1} \delta_{X^j_x}.$$
Let $m_x$ denote its expectation
$$
m_x := \mathbb{E} \big( m_x({\bf g}) \big) = \mathbb E \left(\sum_{j=0}^{T(w)-1} \delta_{X^j_x}\right),
$$
which is considered as a measure on $\mathbb R$. 
Finally, we integrate $m_x$ with respect to the measure $\nu_0$ on~$x$,
thus yielding a Radon measure $\nu := \int m_x \, d\nu_0(x)$ on $\mathbb R$. 
Formally, this means that for any compactly supported function $\phi$, 
\begin{equation}\label{eq:nudef}
\int_{\bf R} \phi \, d \nu = \int_{\mathbb R} \mathbb E \left(\sum_{j=0}^{T(w)-1} \phi (X_x^j)\right) \, d\nu_0(x).
\end{equation}
Note that the right-side expression in~\eqref{eq:nudef} is well-defined and finite. Indeed, there exist
$N \in \mathbb{N}$ and $p_0 > 0$ such that with probability at least $p_0$ a trajectory starting at
any point of $\text{supp} (\phi)$ hits $K$ in at most $N$ steps. Therefore, the distribution of the measure
$m_x(w)$ on $\text{supp} (\phi)$ ({\em i.e.}, the number of steps that are spent in $\text{supp} (\phi)$
until the stopping time) has an exponentially decreasing tail. Thus, its expectation is finite and
bounded uniformly on $x \in \text{supp}(\phi)$, which implies the finiteness of the integral.

Next, let us check that the measure $\nu$ is $P$-invariant. To do this, let us rewrite
the measure ~$\nu$ as follows. First, note that, by definition, we have
$$
m_x
= \sum_{n\ge 0} \sum_{g_1,\dots,g_n\in G} \left[ \prod_{j=1}^n \p (g_j)
\prod_{j=1}^{n} \big[ 1-\xi(g_j \cdots g_1(x)) \big] \right] \delta_{g_n \cdots  g_1(x)}.
$$
Thus,
\begin{eqnarray*}
P (m_x) \!\!
&=&
\! \sum_{g \in G} \p(g) \hspace{0.06cm} g_* (m_x)
\\
&=&
\! \sum_{g \in G} \p (g) \hspace{0.1cm} g_* \!\! \left( \sum_{n\ge 0} \sum_{g_1,\dots,g_n\in G}
\left[ \prod_{j=1}^n \p (g_j)  \prod_{j=1}^{n} \big[ 1-\xi(g_j \cdots  g_1(x)) \big] \right]
\delta_{g_n\cdots g_1(x)}\right)
\\
&=&
\! \sum_{n\ge 0} \sum_{g_1,\dots,g_n, g\in G} \left( \p (g) \prod_{j=1}^n p(g_j)\right)
\prod_{j=1}^{n} \big[ 1-\xi(g_j \cdots g_1(x)) \big] \hspace{0.1cm} g_* \big( \delta_{g_n \cdots g_1(x)} \big)
\\
&=&
\! \sum_{n\ge 0}\sum_{g_1,\dots,g_n,g_{n+1} \in G} \left( \prod_{j=1}^{n+1} \p (g_j)\right)
\prod_{j=1}^{(n+1)-1} \big[ 1-\xi(g_j \cdots g_1(x)) \big] \hspace{0.1cm}
\delta_{g_{n+1} g_n\cdots g_1(x)}.
\end{eqnarray*}
As before, the last expression equals the expectation of the random measure
$\sum_{j=1}^{T({\bf g})} \delta_{X^j_x}$. In this sum, we are counting
the stopping time, but not the initial one. Therefore,
$$
P m_x = m_x -\delta_x + \mathbb E (\delta_{Y_x}).
$$
By integrating with respect to $\nu_0$, this yields
\begin{multline*}
P \nu = P \Big( \int_{\mathbb  R} m_x \, d\nu_0(x) \Big) = \int_{\mathbb R} P(m_x) \, d\nu_0(x) =  \\
\int_{\mathbb  R} m_x \, d\nu_0(x) - \int_{\mathbb R} \delta_x \, d\nu_0 (x) +
\int_{\mathbb  R} \mathbb E (\delta_{Y^x}) \, d\nu_0(x) =
\nu - \nu_0 + P_{\xi} (\nu_0).
\end{multline*}
Since $\nu_0$ is $P_{\xi}$-invariant, we finally obtain $P\nu=\nu$,
as we wanted to show.
$\hfill\square$

\vsp\vsp\vsp

\begin{thm}\label{c:semiconjugation}
\textit{Let $\Gamma$ be a finitely-generated group endowed with a symmetric 
probability $\rho$. If the support of $\rho$ is finite and generates $\Gamma$, then  
every minimal action $\Phi : \Gamma \rightarrow \text{Homeo}_+ (\mathbb R)$
is topologically conjugate to a $\rho$-harmonic action.  }
\end{thm}

\noindent{\bf Proof.}
 By Theorem \ref{T:existence of invariant measure}, there exists a $P$-invariant Radon measure $\nu$. 
 Due to Lemmas  \ref{L:bi-infiniteness}, \ref{L: atomic part 2}, and \ref{L: support}, respectively, 
 the measure is bi-infinite, has no atoms, and its support support is 
 total. As a consequence, there exists a homeomorphism 
$\varphi \!: \mathbb R \rightarrow \mathbb R$ such that $\varphi _ * (\nu)$ is the Lebesgue 
measure. The conjugate action $\varphi \circ \Phi \circ \varphi ^{-1}$ is then $\p$-harmonic.
$\hfill\square$

\vspace{0.45cm}

When the action of $\Gamma$ admits discrete orbits, we know from 
Lemma \ref{L: atomic part 1} that every stationary 
measure must be $\Gamma$-invariant. However, two such measures
may be supported on different orbits. We next establish the uniqueness (up to
a scalar factor) of the stationary measure in the case where there is no discrete orbit.
Recall that, in this case, there exists a unique nonempty, closed, minimal
$\Gamma$-invariant set $M$ (see Lemma \ref{minimal-set}).

\vspace{0.31cm}

\begin{prop} \label{T:unicity}
\textit{Assume that there is no discrete orbit for the $\Gamma$-action on the line.
Then the $P$-invariant Radon measure $\nu$ is unique up to a scalar factor.}
\end{prop}

\vsp

We begin the proof with some reductions. First, we can assume that the action is minimal, since stationary measures are supported on
$M$ (see Lemma \ref{L: support}), and the action is semiconjugate to a minimal one. Moreover, Theorem \ref{c:semiconjugation} 
allows us (via a topological conjugacy) 
to assume that the Lebesgue measure is stationary, that is, that the action is $\rho$-harmonic.

Recall that a $P$-invariant measure is said to be {\bf{\em ergodic}} if every $\Gamma$-invariant measurable subset of the line either 
has measure $0$ or its complement has measure $0$.  Every $P$-invariant measure decomposes as an integral of ergodic measures 
\cite{sinai}. Thus, to prove Theorem \ref{T:unicity}, it suffices to show that, up to multiplication by a constant, there exists a unique 
ergodic $\p$-stationary measure.  

\vsp

\begin{lem}\label{L:new-lem}
\textit{Assume that the action of $\Gamma$ is minimal and $\rho$-stationary.
Let $\nu$ be an ergodic $P$-invariant measure. Then for all continuous functions $\phi, \psi$
with compact support, with $\phi \geq 0$ and $\phi \equiv 1$ on the recurrence interval $K$
given by Corollary \ref{C:recurrence}, and for every $x \in \mathbb R$, it almost surely holds
\begin{equation}\label{E: convergence}
\frac{S_k \psi ( x, {\bf g} ) }{S_k \phi (x,{\bf g})  }
\longrightarrow \frac{\int \psi \esp d\nu} {\int \phi \esp d\nu}
\end{equation}
as $k$ tends to infinity, where
$S_k \psi  (x,{\bf g}) := \psi (X_x^0) + \psi (X_x^1) + \ldots + \psi (X_x^{k-1})$ (and similarly for $S_k \phi $).}
\end{lem}

\vsp

For the proof, we will apply Hopf's ratio ergodic theorem \cite{Hopf} (see also \cite{KK}) to the system
$(\mathbb R^{\mathbb{N}_0}, \sigma, \hat{\nu})$,
where $\sigma$ is the shift operator $\sigma (X^n ) _n = (X^{n+1})_n$,
and $\hat{\nu}$ is the image of the measure $\nu \times \rho^{\mathbb{N}}$ under the map
$$\big( x, {\bf{g}} = (g_n) \big) \mapsto (X_x^0 = x, X_x^1, \ldots, X_x^n, \ldots). $$
We leave as an exercise to the reader to verify that $\hat{\nu}$ is invariant under $\sigma$.
(Actually, this is nothing but a reformulation of the fact that $\nu$ is $P$-invariant.)

We claim that the system $(\mathbb R^{\mathbb{N}_0}, \sigma, \hat{\nu})$ is {\bf {\em ergodic}},
that is, every measurable $\sigma$-invariant subset $A$ of $\mathbb R^{\mathbb{N}_0}$ has either
zero or full $\hat{\nu}$-measure. Indeed, for such an $A$, and for a fixed $x\in \mathbb R$, let $p_A (x) $ be the
probability that the sequence $(X_x^n)_{n\geq0}$ belongs to $A$. The function $p_A \!: \mathbb R \rightarrow [0,1]$ 
thus defined is measurable. Since $A$ is $\sigma$-invariant, the property of belonging to $A$ depends only on the tail of 
the sequence. It is then straightforward to check that the function $p_A$ is $P$-invariant. We claim that this function is indeed
constant. To prove this, note that we cannot directly apply Exercise \ref{ex: Garnett}, because the function $p_A$ has no reason 
to belong to $\mathcal{L}^1(\mathbb{R},\nu)$. To overcome this difficulty, let us consider a compact interval $I$ containing the 
recurrence interval $K$. Given a point $x \in I$, we denote by $Y_x^m, \ldots, Y_x^m \ldots$ the points of the sequence $(X_x^n)$ 
that belong to $I$. As we are assuming that the Lebesgue measure is
$P$-stationary, the Markov process  $Y$ on $I$ leaves the restriction of the Lebesgue measure
on $I$ invariant. Moreover, the restriction of the function $p_A$ to $I$ is still harmonic for the Markov process $Y$, 
namely, $p_A (x) = \mathbb E ( p_A (Y_x^1) ) $ for every $x \in I$. The Lebesgue measure of $I$ being finite,
an easy extension of Exercise \ref{ex: Garnett} for Markov processes
shows that $p_A$ is almost-surely constant on the intersection of a.e. orbit and $I$. As this is true for every compact
interval $I$ containing $K$, we conclude that $p_A$ is constant on almost-every orbit, and since the measure \(\nu \) is ergodic, \(p_A\) is almost everywhere equal to a constant, as was claimed. Now, the $0-1$ law
shows that this constant is either $0$ or $1$, thus showing that $A$ has measure $0$ or its complement has
measure $0$. This concludes the proof that the system  $(\mathbb R^{\mathbb{N}_0} , \sigma, \hat{\nu})$ is ergodic.

\vsp

Next, let $\phi \!: \mathbb R \rightarrow \mathbb R$ be a non-negative function with compact support such that
$\phi \equiv 1$ on the recurrence interval $K$. Then, letting  $\widehat{\phi} ( x,(X^n)_{n \geq 1} ) := \phi(x)$, the
function $\widehat{\phi}$ belongs to $L^1 (\mathbb{R}^{\mathbb{N}_0}, \widehat{\nu})$, and the recurrence property
implies that for $\hat{\nu}$-a.e. $(x,(X^n))$, we have
$$\sum_{k \geq 0}  \widehat{\phi} \big( \sigma^k (x, (X^n )_n) \big) = \infty.$$
A direct application of Hopf's ratio ergodic theorem then implies that for every function
$\widehat\psi \in \mathcal{L}^1 (\mathbb{R}^{\mathbb{N}_0} ,\widehat{\nu})$,
almost surely we have the convergence
$$\frac{\widehat\psi + \widehat\psi \circ \sigma + \ldots + \widehat\psi \circ \sigma^{k-1} }
{\widehat{\phi} + \widehat{\phi} \circ \sigma + \ldots + \widehat{\phi} \circ \sigma^{k-1}}
\longrightarrow \frac{\int \widehat\psi d\hat{\nu}}{\int \widehat{\phi} d\hat{\nu}}.$$
Applying this to a function of the form  $\widehat{\psi}(x,(X^n)_n ) := \psi (x)$, where
$\psi \!: \mathbb{R} \to \mathbb{R}$ is continuous with compact support, and noting that
$$\int \widehat{\phi}  \esp d \hat{\nu}= \int \phi \esp d\nu \qquad \mbox{and} \qquad
\widehat{\phi} + \widehat{\phi} \circ \sigma + \ldots + \widehat{\phi} \circ \sigma^{k-1}
\big( x, (X^n)_n \big)  = S_k \phi (x,{\bf g})$$
(and similarly for $\psi$),
we conclude that \eqref{E: convergence} holds for $\nu$-a.e. $x\in \mathbb R$.

\vsp

The difficulty now is to extend \eqref{E: convergence} to {\em every} $x \in \mathbb R$.
This will follow from the contraction property for $\rho$-harmonic actions below.

\vsp

\begin{lem} \label{L:weak contraction}
\textit{For any fixed number $0<p<1$ and all $x,y$ in the line,
with probability at least $p$ we have}
\begin{equation*}
\lim _{n\rightarrow \infty}  | X^n_x - X^n_y | \leq \frac {|x-y| } {1-p}.
\end{equation*}
\end{lem}

\noindent{\bf Proof.} For simplicity, assume that $y<x$. Since $\nu$ is $P$-invariant, the sequence of random variables
${\bf{g}} \mapsto X^n_x - X^n_y $ is a {\em positive martingale}. In particular, for every integer
$n \geq 1$, we have
\[ \mathbb E ( X^n_x - X^n_y ) = x - y.\]
By the martingale convergence theorem, the sequence $ (X^n_x - X^n_y) $
almost surely converges to a non-negative random variable $v(x,y)$.
By Fatou's inequality, we have
\[ \mathbb E (v (x,y)) \leq \lim_{n\rightarrow \infty}
\mathbb E ( X^n_x - X^n_y ) = x- y.\]
The lemma then follows from Chebyshev's inequality.
$\hfill\square$

\vspace{0.4cm}

Let now $y \!\in\! {\mathbb R}$ and the functions $\phi$, $\psi$ as in the statement of Lemma
\ref{L:new-lem} be fixed. We claim that, for each $m \geq 1$, with probability at least
$1 - 1/m$ we have
\begin{equation}
\label{eq:almost equidistribution}
\limsup_{k \rightarrow \infty} \left| \frac{S_k \psi (y,{\bf g})}{S_k \phi (y,{\bf g}) }
- \frac{\int \psi \esp d\nu}{\int \phi \esp d\nu} \right| \leq \frac{1}{m}.
\end{equation}
Once this is established, it will obviously imply that (\ref{E: convergence})
holds almost surely at all points, as desired.

To show (\ref{eq:almost equidistribution}), note that we know from Lemma \ref{L:new-lem} that   
for a {\em $\nu$-generic} point $x$, the equality 
\begin{equation}\label{eq:esta-si-si}
\lim_{k \rightarrow \infty} \frac{S_k \bar{\psi} (x,{\bf g})}{S_k \phi (x,{\bf g})}
= \frac{\int \bar{\psi} \esp d\nu}{ \int \phi \esp d\nu},
\end{equation}
holds with probability 1 for any compactly supported function $\bar{\psi}$. Since
$\nu$ has total support, 
such an $x$ may be chosen sufficiently close to $y$, say 
$|x- y| \leq \varepsilon$ for a prescribed $\varepsilon > 0$. 
By Lemma~\ref{L:weak contraction}, with probability at least~$1 - 1/m$ we have that for all 
sufficiently large $k$, say $k \ge k_0({\bf g})$, 
\begin{equation}\label{eq:dist-3e}
|X^k_y - X^k_x | \leq m \varepsilon .
\end{equation}

Now, instead of estimating the difference in~(\ref{eq:almost equidistribution}),
it suffices to obtain estimates of the ``relative errors''
\begin{equation}\label{eq:rel1}
\limsup_{k \rightarrow \infty}
\left| \frac{S_k \psi (y,{\bf g}) - S_k \psi (x,{\bf g})}{S_k \phi (x,{\bf g}) } \right|
\leq \delta_1(\varepsilon)
\end{equation}
and
\begin{equation}\label{eq:rel2}
\limsup_{k \rightarrow \infty}
\left| \frac{S_k \phi (y,{\bf g}) - S_k \phi (x,{\bf g})}{S_k \phi (x,{\bf g}) } \right|
\leq \delta_2(\varepsilon),
\end{equation}
in such a way that $\delta_1 (\varepsilon) \to 0$ and $\delta_2(\varepsilon) \to 0$
as $\varepsilon \to 0$.

Since the estimate~\eqref{eq:rel2} for~$\phi$ is a particular case of the
estimate~\eqref{eq:rel1}, we will only check~\eqref{eq:rel1}. Now,~\eqref{eq:dist-3e}
implies that $\left| S_k \psi (y,{\bf g}) - S_k \psi (x,{\bf g}) \right|$ is at most
\begin{small}
\[ {\rm{mod}}(m \varepsilon, \psi) \hspace{0.05cm} \text{card} \big\{ k_0({\bf g})\le j\le k \! : \text{either} \,
X^j_x \, \text{or} \, X^j_y \,\, \text{is in } \, \text{supp}\psi \big\} + 2 k_0({\bf g}) \max |\psi| \]
\[ \le  {\rm{mod}}(m \varepsilon, \psi) \esp
\text{card}\{j\le k \mid X^j_x \in U_{m \varepsilon}(\text{supp}\psi)\} + \text{const}({\bf g}).
\]
\end{small}Here, $\rm{mod}(\cdot, \psi)$ stands for the modulus of continuity of~$\psi$ with respect
to the distance $d$ on the variable, and $U_{m \varepsilon}(\text{supp} \psi)$ denotes the
$m \varepsilon$-neighborhood of the support of $\psi$, again with respect to $d$.

Let $\xi$ be a continuous function satisfying $0 \le \xi \le 1$ and that is equal to~$1$ on
$U_{m\varepsilon}(\text{supp} \psi)$ and to~$0$ outside~$U_{(m+1) \varepsilon}(\text{supp} \psi)$.
We have
$$
\text{card} \big\{ j\le k \!: X^j_x \in U_{m \varepsilon}(\text{supp} \psi) \big\} \le S_k \xi (x,{\bf g}).
$$
Thus,
\begin{eqnarray*}
\left| \frac{S_k \psi (y,{\bf g}) - S_k \psi (x,{\bf g})}{S_k \phi (x,{\bf g}) } \right|
&\leq&
\frac{\text{const}({\bf g}) + {\text{mod}}(m \varepsilon, \psi) \hspace{0.05cm}
S_k \xi (x,{\bf g})}{S_k \phi (x,{\bf g})} \\
&\longrightarrow&
{\rm{mod}}(m\varepsilon, \psi)\hspace{0.04cm} \frac{\int \xi \, d\nu}{\int \phi \, d\nu}
\hspace{0.38cm} =: \hspace{0.281cm} \delta_1(\varepsilon).
\end{eqnarray*}
Here, we have applied the fact that, by our choice of $x$,
equality~\eqref{eq:esta-si-si} holds with~$\xi$ in the numerator and
$\phi$ in the denominator. Since ${\rm{mod}}(m\varepsilon, \psi)$ tends
to~$0$ as $\varepsilon\to 0$ and the quotient
$$
\frac{\int \xi \, d\nu}{\int \phi \, d\nu} \le
\frac{\nu(U_{(m+1)\varepsilon}(\text{supp} \psi))}{\int \phi \, d\nu}
$$
remains bounded, this yields $\delta_1(\varepsilon)\to 0$ as $\varepsilon \to 0$.
$\hfill\square$

\vspace{0.421cm}

Having Lemma \ref{L:new-lem} at hand, it is now easy to finish the proof of Proposition \ref{T:unicity}.
Indeed, given any two ergodic $P$-invariant Radon measures $\nu_1,\nu_2$, for each
$x \in M$ and every compactly supported,
real-valued function $\psi$, almost surely we have
$$\frac{S_k \psi (x,{\bf g})}{S_k \phi (x,{\bf g})}
\longrightarrow \frac{\int \psi \esp d\nu_i} {\int \phi \esp d\nu_i},$$
where $i \in \{1,2\}$. Thus, $\int \psi \esp d\nu_1 = \lambda \int \psi \esp d \nu_2$,
with $\lambda := \int \phi \esp d \nu_1 / \int \phi \esp d\nu_2$.
This proves that $\nu_1 = \lambda \nu_2$, and concludes the proof of Proposition \ref{T:unicity}.

\vsp

\begin{small}\begin{ejer} Show that the condition on $\phi$ in Lemma \ref{L:new-lem} can
be relaxed to $\phi \geq 0$ and $\phi \not\equiv 0$.
\end{ejer}\end{small}

\vsp

We close with the next result of uniqueness of the conjugation to an harmonic action.

\vspace{0.15cm}

\begin{thm} \label{t: unicity conjugation to harmonic action}
\textit{The conjugacy of a minimal action to a $\p$-harmonic one is unique
up to post-composition with an affine map.}
\end{thm}

\noindent{\bf Proof.}
Given a minimal action $\Phi \!: \Gamma \rightarrow \text{Homeo}_+ (\mathbb R) $
and two homeomorphisms $\varphi_i \!: \mathbb R \rightarrow \mathbb R $ such that	 
each $\varphi_i \circ \Phi \circ \varphi_i^{-1}$ is $\rho$-harmonic (with $i \in \{1,2\}$), the images of the
Lebesgue measure by $\varphi_1^{-1}$ and $\varphi_2^{-1}$ are $\p$-stationary for $\Phi$,
hence they differ by multiplication by a constant. Therefore, $\varphi_2 \circ \varphi_1^{-1}$
sends the Lebesgue measure to a multiple of itself, which means that
$\varphi_2 \circ \varphi_1^{-1}$ is an affine map.
$\hfill\square$

\vspace{0.25cm}

Having established existence and uniqueness (up to post-composition with an affine map) of the $\p$-stationary 
measure, the next exercise gives some insight into what happens when changing the probability distribution $\p$. 

\begin{small}\begin{ejer}
Show that the harmonic flows (see the end of \S \ref{S: derriennic}) corresponding to two different finitely supported symmetric probability measures on \(\Gamma\) whose supports generate \(\Gamma\) are {\bf {\em orbitally conjugate}}, which means that there exists a homeomorphism between the corresponding almost-periodic spaces that exchange trajectories of the harmonic flows (without necessarily preserving their time parametrizations).  
\end{ejer}\end{small}

We close this section with a clever remark by Brum in the form of an exercise.

\begin{small}
\begin{ejer} 
Given a finitely-generated, left-orderable group $\Gamma$, let \(d\in \mathbb N\) denote its first Betti number, and 
let \(\pi : \Gamma \rightarrow \mathbb{Z}^d \) be the quotient of the abelianization of \(\Gamma\) by its torsion subgroup.  
Recall from Exercise \ref{ej: harmonic space abelian group} that $\text{Har} (\mathbb{Z}^d)$, the harmonic space of 
$\mathbb{Z}^d$, is homeomorphic to $\mathbb{S}^{d-1}$. 
Show that the natural map \[ \Phi \in \text{Har} (\Lambda ) \rightarrow \Phi \circ \pi \in \text{Har} (\Gamma )\] 
induces an injection \( \pi _k (\mathbb S^{d-1}) \rightarrow \pi _k (\text{Har} (\Gamma))\) at the level of homotopy groups. 

\noindent{\underline{Hint.}} Let \(g_1,\ldots, g_d$ be elements in $ \Gamma $ whose images generate $\mathbb{Z}^d$. Study 
the map that sends an action \(\Phi\in \text{Har} (\Gamma)\) to  \( (\Phi (g_1)(0), \ldots, \Phi (g_d) (0)) \in \mathbb{R}^d \). 
Observe that this map does not take the value \( (0,\ldots, 0)\),  and use that \(\mathbb R^d\setminus \{(0,\ldots, 0)\}\) 
deformation retracts onto \(\mathbb S^{d-1}\).
\end{ejer}
\end{small}

\thispagestyle{empty}

\section{A finitely-generated, left-orderable, simple group} \label{s: finitely generated simple left orderable}

Finitely-generated, infinite, simple groups are not easy to construct. For example, any finitely-generated matrix group is residually finite, 
that is, the intersection of its finite index normal subgroups is the trivial group \cite{malcev res finite}. In particular, such a group 
is not simple if infinite. Infinite hyperbolic groups form another family where no simple group arises \cite{Delzant}.  
In 1951, Higman built the first example of a finitely-generated, infinite, simple group \cite{Higman}. This was later 
refined by Thompson in unpublished notes dating from 1965, where he introduced what is nowadays called 
{\bf \em Thompson's group $\mathrm{T}.$}}  This is an extension of Thompson's group $\F$ introduced 
in \S \ref{ejemplificando-4},  and consists of the homeomorphisms of the circle that are piecewise-affine 
with powers of 2 as derivatives and dyadic points as break points. Thompson proved that 
this group is both simple and not only finitely-generated but also finitely-presented. 

In the famous Kourovka Notebook \cite[Question 16.50]{kourovka}, Rhemtulla asked whether there exist finitely-generated simple groups 
which are left-orderable. The question was brilliantly answered in the affirmative by Lodha and Hyde \cite{hyde-lodha-simple}. 
Soon after, an easier and illuminating construction was produced in 
\cite{MT}  by Matte Bon and Triestino.\footnote{More recently, more examples of simple finitely-generated, left-orderable groups 
where built in \cite{hyde-lodha-rivas}. However, all the aforementioned examples fail to be finitely-presented, as it is very well 
explained in \cite{fedo-lodha}. However, in \cite{hyde-lodha-finitely presented}, Hyde and Lodha found an example of a 
finitely-presented, simple group that is left orderable.} The key observation 
for their construction is that one can use Thompson's ideas but replacing the action of \(\mathrm{T}\) 
on the circle by a kind of almost-periodic action on another compact laminated space. 
The goal of this closing section is to present the beautiful examples of Matte Bon and Triestino. 

\subsection{The Thompson group of a suspension}

Let \(X\) be a  Cantor set and \( \varphi : X\rightarrow X\) a homeomorphism of \(X\). We will 
refer to the pair $(X,\varphi)$ as a {\bf \em Cantor system}. In the sequel, we will mostly 
assume that \(\varphi\) is minimal, {\em i.e.}, every \(\varphi\)-orbit is dense in \(X\). In such 
a case, we will refer to  $(X,\varphi)$ as a {\bf \em Cantor minimal system}. 

The {\bf \em suspension} of  $(X,\varphi)$ is the space $\hat{X} := (X \times \R)/\Z$, where the quotient is taken with respect 
to the diagonal action of $\Z$ on $X \times  \R$ given by $n \cdot (x, t) = (\varphi^n(x), t+n)$. The natural projection of 
$(x,t)\in X\times \R$ to $Y$ will be denoted $\pa x,t \pc$. 

The suspension $\hat{X}$ is a compact space naturally equipped with the {\bf \em translation flow} 
$S = \{S_s \}_{s \in \mathbb{R}}$ given by 
\begin{equation}\label{eq flow} S_s (\pa x,t \pc) = \pa x,t+s \pc. 
\end{equation}
Note that, by definition, $S_1(\pa x,0 \pc) = \pa \varphi^{-1}(x),0 \pc$. 
Thus, the time-1 map of the suspension flow mimics the inverse of the Cantor set homeomorphism. 

Recall from \S  \ref{ejemplificando-4} that a dyadic number is a rational number whose denominator is a power of \(2\), 
and that a piecewise-dyadic homeomorphism of the real line (or between open subsets of the real line) is an orientation-preserving 
homeomorphism that is piecewise-affine, having powers of \(2\) as derivatives and dyadic numbers as break points. 
 {\bf \em The Thompson group of the suspension} of $\varphi$, denoted \( \mathrm{T}(\varphi) \), is a subgroup of 
 the group of homeomorphisms of \(\hat{X} \) that preserve the \(S\)-orbits, acting on each of them as piecewise-dyadic 
 homeomorphisms with respect to their time parametrization by \(S\).  Specifically, $\mathrm{T}(\varphi)$ 
 is the group of homeomorphism of $\hat{X}$ that are locally of the form 
$$\pa x,t \pc \mapsto \pa x,h(t) \pc,$$
where $h$ is a dyadic homeomorphism of the real line. Note that $S_1$, the time-1 map of the flow, 
as well as all its integer powers, are elements of  $\mathrm{T}(\varphi)$.

Since \(\varphi\) is minimal, the action of the translation flow \( S \) is free. 
As a consequence, every element \( h \in  \mathrm{T} (\varphi) \) lifts to a unique homeomorphism 
\(\tilde{h} \) of \( X \times \R\) of the form \( \tilde{h} (x, t) = ( x, h_x (t)) \), where \( \{h_x\}_ {x\in X}\) 
is a family of dyadic homeomorphisms of the real line that satisfy the following two properties:

\vsp

\noindent -- ({\em Equivariance}) \, For all  $x\in X$ and $t \in \R$, 
\[     h_{\varphi (x) } (t+1) = h_x (t) + 1;\]  

\vsp

\noindent -- ({\em Continuity}) \, The map from $X$ into $\mathrm{Homeo}_+(\mathbb{R})$ 
sending $x$ to $h_x$ 
is locally constant.

\vsp

\begin{prop}\label{prop LO MB-T}
\textit{For every Cantor minimal system $(X,\varphi)$, the associated group \(\mathrm{T}(\varphi)\) is left-orderable.}
\end{prop}

\noindent{\bf Proof.} By looking at the action of $\mathrm{T} (\varphi)$ on the parametrization of any $S$-orbit, we obtain an action 
of $\mathrm{T}(\varphi)$ on the real line by orientation-preserving homeomorphisms. Since all the \(S\)-orbits are dense, this action 
is faithful. Therefore, \(\mathrm{T}(\varphi)\) is left-orderable. $\hfill\square$

\vspace{0.28cm}

One of the fundamental features of the group \( \mathrm{T}(\varphi)\) is that it contains many copies of two of the Thompson 
groups.  We first describe the construction of copies of Thompson's group \( \F\) inside \( \mathrm{T}(\varphi)\) . 

Recall that a  {\em dyadic interval} is a compact interval of the real line whose endpoints are dyadic numbers. 
Recall also that, in \S \ref{ejemplificando-4}, for each dyadic interval $I$ contained in $[0,1]$, we let $\F_I$ be the 
subgroup of $\F$ formed by the elements whose support is contained in $I$. Equivalently, such an $\F_I$ may 
be viewed as the group of piecewise-dyadic homeomorphisms of \(I\). Here, we extend the latter definition and the notation 
to every dyadic interval in the line. Since any dyadic interval is the image of the interval \([0,1]\) by a piecewise-dyadic 
homeomorphism (see Exercise \ref{ej: transitivity on dyadic interval}), 
each group \( \F_I\) is a copy of the classical group \( \F= \F_{[0,1]}\). 

\vsp

\begin{small}\begin{ejer} \label{ej: extension of dyadic homeomorphism}
Prove that if \(I,J,K\) are dyadic intervals of the real line such that both \(J\) and \(K\) are contained in the interior of  \(I\), 
then every piecewise-dyadic homeomorphism  from $J$ onto $K$ can be extended to an element of \( \F_I\). Prove that 
the same conclusion holds if $J$ and $K$ are both contained in $I$, the three intervals share one endpoint, but the 
other endpoints of $I$ and $J$ are both different from that of $K$.

\noindent{\underline{Hint.}}  Apply Exercise \ref{ej: transitivity on dyadic interval} to the right (resp. left) 
components of \( I\setminus  \text{int} (J)\) and \( I\setminus \text{int} (K)\), where $\mathrm{int}(\cdot)$ 
denotes the interior of the corresponding interval.
 \end{ejer}\end{small}

\begin{small}\begin{ejer} \label{ej: fragmentation in F}
Let \(I,J,K\) be dyadic intervals such that \(J\cap K\) has nonempty interior and \(I = J\cup K\). 
Show that \( \F_I\) is generated by the natural copies of \( \F_J\) and \( \F_K\) inside \(\F_I\).

\noindent{\underline{Hint.}}  One can assume that \(I=[0,1]\), \( J= [0,b]\) and \(K=[a,1]\) for some dyadic numbers \(0<a<b<1\).  
Let $f \in \F_I$. If \( f(a) <b\), use Exercise \ref{ej: extension of dyadic homeomorphism} to prove that there exists \(g \in \F_J\) 
whose restriction to \([0,a] \) agrees with \(f\), and then note that \( g^{-1} f \in \F_K\). If \( f(a) \geq b\), choose 
\(h \in \F_K\) such that \( h f(a) <b\), and apply the first case. 
 \end{ejer}\end{small}
 
\vsp

A subset \(C \subset X\) that is simultaneously closed and open is usually called a {\bf \em clopen} set.  
The family of these subsets is a countable base for the topology of \(X\). If \(I\) is a dyadic interval and 
\(C\subset X\) is a clopen subset, it may be the case that the inclusion of \( C \times I \) in \(X \times \R \) 
induces an injective map from $C \times I $ into $ Y$. If this happens, its image is called a {\bf \em dyadic flow box}, 
and it is denoted by \( \pa C\times I \pc \). Note that such a box naturally identifies with \(C\times I \). We also 
say that the pair \((C,I)\) defines a dyadic flow box. 

Two observations are in order. First, if \( (C, I)\) defines a dyadic flow box, then the same holds for 
\( (\varphi^n (C), I+n)\) for all \(n\in \mathbb Z\), and we have 
\begin{equation} \label{eq: equivariance dyadic box} 
\pa \varphi^n (C) \times ( I + n) \pc = \pa C\times I \pc. 
\end{equation} 
Second, denoting by \(l\) the integer part of the length of \(I\), the pair \((C,I)\) defines a dyadic flow box if 
and only if the clopen sets \( C, \varphi (C) , \ldots , \varphi ^{l}(C) \) are two-by-two disjoint. In particular, if 
the length of \(I\) is \(<1\), then \((C,I)\) always defines a dyadic flow box.

\vsp

\begin{small}\begin{ejer} \label{ej: complement dyadic box}
Prove that for any dyadic flow box \( \pa C \times I \pc \subset \hat{X} \), the set \( \hat{X} \setminus \text{int} ( \pa C \times I \pc ) \) 
is a finite disjoint union of dyadic flow boxes whose boundaries are contained in the boundary of \( \pa C \times I \pc \). 

\noindent{\underline{Hint.}} By minimality of the translation flow \( S\) on \(\hat{X}\), every point in the complement of \( \pa C\times I  \pc \) 
belongs to a unique piece of \( S\)-trajectory of the form \( \pa \{x_0\} \times J \pc \) for some dyadic interval \(J= [a,b]\subset \mathbb R\) 
and some \(x_0\in X\). Both boundary points \( \pa x_0, a \pc \) and \( \pa x_0, b \pc \) belong to the boundary of \( \pa C\times I \pc \). 
Note that for every \(x\) sufficiently close to \(x_0\), the set \(\pa x \times \text{int} (J) \pc \) is disjoint from \( \pa C \times I \pc \), 
while both points \( \pa x, a \pc \) and \( \pa x, b \pc \) belong to \( \pa C\times I \pc \). \end{ejer}\end{small}

\vsp

Given an element \( g\in \F_I\), we define \( \F_{C,I} (g)  \) as the homeomorphism of \(\hat{X}\) that acts as 
\( \text{id} \times g\) on \( \pa C\times I  \pc \) and as the identity outside.  The set of all these homeomorphisms 
is denoted by \( \F_{C,I}\). Obviously, this set is a copy of \( \F_I \sim \F\). A fundamental fact  is that these copies of 
\( \F\) inside \(\mathrm{T}(\varphi)\) all together form a generating set. 

\vsp

\begin{prop}  \label{p: generating Thompson copies}
\textit{The groups \( \F_{C,I}\), where \(C\) ranges over all clopen subsets of \(X\) and \(I\) over all dyadic intervals 
of \(\mathbb R\) of length \(<1\), generate the whole group \( \mathrm{T} (\varphi)\). }
\end{prop}

\noindent{\bf Proof.} 
Let \( C\subset X\) be a clopen set and \( I \subset \mathbb R\) a dyadic interval such that \( (C,I)\) defines a dyadic flow box. Let 
\( I_1, \ldots , I_r\) be a family of dyadic intervals contained in $I$, all of length \(<1\), such that each intersection \( I_i \cap I_{i+1}\) 
has nonempty interior and the union $\bigcup I_i$ equals \(I\). An inductive application of Exercise \ref{ej: fragmentation in F} shows 
that the group \( \F_I \) is generated by the natural copies of \( \F_{I_i}\) therein, with \(i\) ranging over \( \{ 1,\ldots , r \} \). 
In particular, \( \F_{C,I}\) is generated by \( \F_{C, I_1}, \ldots, \F_{C,I_r}\). 

As a consequence, it suffices to establish that the groups \( \F_{C,I}\) generate \( \mathrm{T} (\varphi)\) when \(C\) ranges 
over all the clopen subsets of \(X\) and \(I\) over all the dyadic intervals of \(\mathbb R\) such that \( (C,I)\) defines a dyadic flow box. 
To do this, let \(h \in \mathrm{T} (\varphi)\) and \(  \pa x,s \pc \in \hat{X} \). Let \( J \) be a dyadic interval whose interior contains the segment 
\( [s, h_x (s) ]\). By the ({\em Continuity}) property above, there exists a clopen neighborhood \(C\) of \(x\) such that the restrictions 
of the homeomorphisms \( h_y \) to \( J\) do not depend on \(y \in C\). We denote this common homeomorphism by $\bar{h}$. 
Note that one can choose \(C\) small enough so that \( (C, J )\) defines a dyadic flow box.    

Let \(I\subset J\) be a dyadic interval that contains \(s\) and such that \( \bar{h} (I) \) is contained in the interior of \(J\).  
Using Exercise \ref{ej: transitivity on dyadic interval}, one can extend \( \bar{h} \) to a dyadic homeomorphism 
of \( J\), that is, to an element in \( \F_ J \) (still denoted $\bar{h}$). We let \( f := \F_{C,J} ( \bar {h}) \in \F_{C,J}\). 
By construction, the element \( g := f^{-1}  h \in \mathrm{T}(\varphi) \) agrees with the identity on \( \pa C\times I \pc \).   

By Exercise \ref{ej: complement dyadic box}, the complement of the interior of \( \pa C\times I \pc \) is a finite union of dyadic boxes 
\( \pa C_j\times I_j \pc \) whose boundaries are contained in the boundary of \( \pa C \times I \pc \). In particular, since \( g \) preserves 
the trajectories of the flow, each of the dyadic boxes \( \pa C_j\times I_j \pc \) is invariant under \( g \). Again, by the ({\em Continuity}) 
property above, up to taking a subdivision of each of the clopen sets \(C_j\) into a finite number of smaller clopen sets, one can 
assume that the restriction of $g$ to each \( \pa C_j\times I_j \pc \) is an element of $ \F_{C_j,I_j}$.
$\hfill\square$

\vspace{0.35cm}

There is also a copy of another Thompson's group inside \( \mathrm{T}(\varphi)\), namely, the group \(\widetilde{\mathrm{T}}\) 
of piecewise-dyadic homeomorphisms of the real line that commute with the translation \( t\mapsto t+1\). This is a central 
extension of Thompson's group \( \mathrm{T} \), the group of piecewise-dyadic homeomorphisms of the circle.

\vspace{0.15cm}

\begin{prop} 
\textit{There is a copy of \(\widetilde{\mathrm{T}} \) inside \(\mathrm{T}(\varphi)\).} 
\end{prop}

\noindent{\bf Proof.} Fix a point $x_0 \in X$. Given an element \( g\in \widetilde{\mathrm{T}}\), define a (partial) map $h$ on $\hat{X}$ 
by setting \(h (\pa x_0,t \pc) := \pa x_0, g(t) \pc \). Since the orbits under $\varphi$ are dense, this extends to a homeomorphism 
$h$ of $\hat{X}$. Moreover, by the ({\em Equivariance}) property above, this homeomorphism $h$ belongs to $\mathrm{T}(\varphi)$. 
Furthermore, the map that sends \(g\in \widetilde{\mathrm{T}} \) to the element \( h \in \mathrm{T} (\varphi) \) just 
constructed is a group homomorphism $\Phi$. Finally, $\Phi$ is injective, since  the restriction of each $\Phi (g)$
to the \( S\)-orbit of $x_0$ is the natural action of \(g \in \widetilde{\mathrm{T}}\) on the real line.  
$\hfill\square$

\vsp\vsp

\begin{small}\begin{ejer} \label{ej: finite generation tildeT}
Let \(J\) and \(K\) be two dyadic intervals of length \(<1\) whose intersection has nonempty interior and whose union 
is an interval of length \(>1\). Show that the two natural copies of \( \F_J\) and \( \F_K\) inside \(\widetilde{\mathrm{T}}\) together 
with the unit translation \(t\mapsto t+1\) generate \(\widetilde{\mathrm{T}}\). Conclude that \(\widetilde{\mathrm{T}}\) is finitely-generated. 

\noindent{\underline{Hint.}}  One can assume that \(J=[0,b]\) and \(K= [a,c]\), with \( 0<a<b<1<c\). 
Given any element \(f\in \widetilde{\mathrm{T}}\), show that one can multiply \(f\) by an element of the group generated by 
\( \F_J\), \( \F_K\), and the unit translation, so that the resulting element of \(\widetilde{\mathrm{T}}\) fixes the origin~\(0\). 
Then apply Exercise \ref{ej: fragmentation in F} to show that this element belongs to the group generated by 
\( \F_J\) and \( \F_{[a,1]} \subset \F_K\).
\end{ejer}\end{small}

\subsection{Simplicity of \( \mathrm{T}(\varphi)\)} 

This section is devoted to the proof that the group \( \mathrm{T} (\varphi)\) is simple. The action of this group on \( \hat{X} \)  
shares some properties with actions of groups of homeomorphisms of higher-dimensional manifolds, many of which are 
known to be simple. We thus adapt ideas arising in the classical proofs of simplicity of manifold homeomorphism groups. 
These crucially use the concepts of group perfection and {\bf \em fragmentation}}, the latter meaning that general 
elements of the group can be expressed as a product of elements with smaller supports.   

\begin{prop}\label{teo simple MB-T}
\textit{For every Cantor minimal system $(X,\varphi)$, the associated group \( \mathrm{T} (\varphi)\) is simple.}
\end{prop}

\noindent{\bf Proof.}  Let  \(h \neq id \)  be an element of \( \mathrm{T} (\varphi)\). Our goal is to show that every 
element \(f\) of \(\mathrm{T} (\varphi)\) belongs to the normal closure $N(h)$ of \(h\), that is, the smallest normal subgroup of 
\(\mathrm{T}(\varphi)\) containing \( h \). To do this, by Proposition \ref{p: generating Thompson copies}, we can 
assume that \(f\) belongs to some subgroup \( \F_{C,I}\). 

Let \( \hat{x} \in \hat{X} \) be a point such that \(h (\hat{x}) \neq \hat{x}\). Choose an element \(g\in \mathrm{T} (\varphi)\) that 
is the identity on some small neighborhood of \( \hat{x} \) but moves \( h (\hat{x})\), {\em i.e.}, \( g h (\hat{x}) \neq h(\hat{x})\).  
We can then choose a dyadic box  \( \pa D\times J \pc \) containing \(\hat{x}\) sufficiently small so that \(g\) is the identity 
on \( \pa D\times J \pc\), and the three sets  \( \pa D\times J \pc \), \(h (\pa D\times J \pc\)) and 
\( gh (\pa D\times J\pc ) \) are two-by-two disjoint.

First assume that \( \pa C\times I \pc \subset \pa D\times J \pc \). By Exercise \ref{F:ejer-conmutator}.
the commutator subgroup \( \F_J '\) of Thompson's group \( \F_J\) is the group of dyadic homeomorphisms of the 
interval \(J\) acting trivially on a neighborhood of the endpoints of \(J\). One thus deduces that \(f\) belongs to \( \F_{D, J} ' \). 
In particular, we can write $f$ as a product  \( f = [a_1,b_1] \cdots  [a_r, b_r] \), with each \( a_j, b_j \) in \( \F_{D,J}\). 
We can then use a variant of the Highman trick (compare the proof of Lemma \ref{lem Higman}). Namely, 
one can readily check the relation
\[ [a_j, b_j] = \left[ [a_j,h], g [g^{-1} b_j g, h] g^{-1}  \right] .\]
This shows that \(f\) is a product of commutators of \(h \) and $h^{-1}$ with elements of \( \mathrm{T} (\varphi)\). 
Since these commutators belong to \( N(h) \), we are done in this case. 

To close the proof, the idea is to reduce to the previously treated case by fragmenting and conjugating 
 \(f\). Because \(\varphi \) has dense orbits 
and  \(C\) is compact, one can find a finite partition of \( C\) into clopen sets \( C_i\) so that, for each \(i\), there 
exists a certain integer \(n_i\) such that \( \varphi ^{n_i} (C_i ) \subset D \). The dyadic flow boxes \( \pa C_i \times I \pc \) 
form a partition of \( \pa C\times I \pc \). Moreover, one has \( f = \prod f_i \), where \( f_i\in \mathrm{T} (\varphi)\) 
is the map that coincides with \(f\) on \( \pa C_i\times I \pc \) and equals the identity elsewhere. 
To show $f$ is in  $N(h)$, it suffices to show each $f_i$ is in $N(h)$. 
Now note that the conjugate element $g_i := S_{n_i} f_i S_{-n_i}$ has support equal to 
$S_{n_i}(\text{supp}(f_i))$, which is 
$$\text{supp}(g_i) = \pa \varphi^{n_i} (C_i) \times (I+ n_i) \pc$$
Since  $\varphi^{n_i} (C_i) \subset D$, this support is contained in $\pa D \times (I + n_i) \pc$. 
In other words, we have reduced the case to the one where the involved clopen set is contained in $D$. 

Assume hence that $C \subset D$ and 
consider a dyadic  interval \( K\subset \mathbb R\) whose interior contains both \( I\) and \(J\).  Assume 
first that  \((C,K)\) defines a dyadic flow box. There is an element \(a\) of Thompson's group \( \F _K\)  that 
sends \(I \) into \(J\). The element \( \F_{C,K}(a)\in \F_{C,K} \subset \mathrm{T}(\varphi)\) then sends the dyadic 
flow box \( \pa C\times I \pc \) inside \( \pa C\times J \pc \). By construction, the conjugate of \(f\) by \( \F_{C,K} (a) \) 
belongs to the group \( \F_{C, J} \). Since $C \subset D$, we are  in the first situation above. This 
allows us to conclude that this conjugate belongs to the normal closure of \(h\), and therefore the same holds for \(f\).

Finally, it may happen that \( (C, K)\) doesn't define a dyadic flow box. In this case, we can use again a fragmentation trick. 
Namely, we look for a finite partition of \( C\) into sufficiently small clopen subsets \( C_j\) such that each \( (C_j, K)\) defines 
a dyadic flow box, and we write \( f\) as the product of the maps equal to \(f\) on \( \pa C_j, K \pc \) and the identity elsewhere.  
By the case treated immediately above, each such factor must belong to $N(h)$. Therefore this is 
the case of $f$ as well, which closes the proof. $\hfill\square$

\subsection{Finite generation of $T(\varphi)$.} 

Although the group $\mathrm{T}(\varphi)$ is simple for every Cantor minimal system 
$(X,\varphi)$, it is not always finitely generated, as it is shown be the next two exercises.

\vsp

\begin{small}
\begin{ejer} 
Denote by \( X := \mathbb{Z}_2 \) the Cantor set of {\bf $2$-{\em adic integers}}, and let \( \varphi (x) := x + 1\) be the {\bf \em adding machine}. 
Show that $(X,\varphi)$ is a Cantor minimal system.

\noindent{\underline{Hint.}} Write each $x \in \mathbb{Z}_2$ in the form $x= \sum _{k\geq 0} \varepsilon_k 2^k$, with $\varepsilon_k \in \{0,1\}$, 
and compute.
\end{ejer}

\begin{ejer}\label{T-no-fin-gen}
Referring to the Cantor minimal system of the preceding exercise, 
show that there is a strictly increasing sequence of copies of Thompson's group \(\tilde{\mathrm{T}}\) 
whose union is \( \mathrm{T}(\varphi)\). Conclude in particular that \(\mathrm{T} (\varphi)\) is not finitely-generated. 
 
\noindent{\underline{Hint.}}  For each positive integer $n$, let \( \tilde{\mathrm{T}}_n\) be the group of piecewise-dyadic homeomorphisms 
of the real line that commute with the translation \( t \mapsto t+2^n\). Show that every element of $\mathrm{T}(\varphi)$ belongs to 
 \( \tilde{\mathrm{T}}_n\) for some $n$.
\end{ejer}
\end{small}

We next restrict our attention to a special kind of Cantor minimal systems for which the group $\mathrm{T}(\varphi)$ is 
finitely-generated. These are the so-called subshifts that naturally arise in symbolic dynamics. In concrete terms, 
we say that a Cantor system \((X,\varphi)\) is a {\bf \em subshift} if there is a partition of \(X\) into a finite 
number of clopen subsets \(C_1,\ldots, C_d\) so that one can characterize any point \(x\) of \(X\) by its {\em itinerary}, that is, 
the sequence of clopen sets of the partition visited by its \(\varphi\)-iterates. More concretely, this is the sequence 
\( ( i_n) \in \{1,\ldots, d\} ^{\mathbb Z}\) defined by the property \( \varphi^n (x) \in C_{i_n}\). In such a framework, 
the partition $\{C_1,\ldots,C_d\}$ is a {\bf {\em generating partition}}, in the sense that the sets $\varphi^n(C_i)$, 
for $n\in \Z$ and $i \in \{1,\ldots,d\}$, generate the Boolean algebra of clopen subsets of $X$ (and hence its topology as well).  


\begin{small}\begin{ejer} \label{ej: subshift 0}
Show that any subshift can be conjugated to the restriction to a closed invariant 
set of the {\bf \em full shift over a finite alphabet}, where by the latter we mean a system of the form 
\( (A^{\mathbb Z}, \sigma) \), with \(A\) a finite set and \( \sigma ( (a_n) _{n\in \mathbb Z} ) = (a_{n+1})_{n\in\mathbb{Z}}\). 

\noindent{\underline{Hint.}} Take \( A= \{1,\ldots, d\}\) and look at the map \( i: X\rightarrow A^{\mathbb Z}\) defined by \( i (x) = (i_n(x))_n\), 
where \( i_n (x)\in \{1,\ldots, d\}\) satisfies \( \varphi^n (x) \in C_{i_n}\). Then check the conjugacy relation $i \circ \varphi = \sigma \circ i$. 
\end{ejer}\end{small}

In what follows, we will be interested in subshifts that are minimal. However, examples of minimal subshifts 
are not so easy to provide. The first one was found by Morse and Hedlund in their study of symbolic dynamics 
of irrational rotations of the circle \cite{Morse Hedlund}. We reproduce this example in the exercise below.

\begin{small}\begin{ejer} \label{ej: subshift}
Let \( R = R_{\alpha }\) denote the rotation by an angle $\alpha$ on the circle $\mathbb{S}^1$, which we now identify 
with $\mathbb R /\mathbb Z$ for simplicity. Assume throughout that $\alpha$ is irrational. 
Let $I^-$ and $I^+$ be the subsets of $\mathbb{S}^1$  given by 
\(I^- := [0, \alpha)\) and \( I^+ := (0, \alpha]\). 
Let \(i^- = (i_n^-): \mathbb S^1 \rightarrow \{0,1\}^{\mathbb Z}\) be the map defined by 
$$
i_n^-(x) = \left\{
  \begin{array}{l}
    1  \quad \mathrm{ if } \,\, R^n(x) \in I^-, \\
    0  \quad \mathrm{otherwise.}
  \end{array}
\right.
$$
Let  \(i^+ = (i_n^+): \mathbb S^1 \rightarrow \{0,1\}^{\mathbb Z}\)  be the analogously defined map obtained 
by replacing $I^-$ by $I^+$. 

\noindent (i) Show that the subset of \(\{0,1\}^{\mathbb Z}\) defined by 
\[ X : = i ^-(\mathbb S^1) \cup i ^+ (\mathbb S^1) \]
is invariant under the shift map \(\sigma ( (a_n)_n)) = (a_{n+1}) _n  \). 
(Sequences belonging to \(X\) are called {\em Sturmian sequences}.) 

\noindent (ii) Show that \(X\) is closed.

\noindent (iii) Show that \( (X, \sigma ) \) 
is a Cantor minimal system which is a subshift. 

\noindent{\underline{Hint.}}  Note the following facts: 

\noindent - The equality \( i ^- (x) = i^+ (x)   \) holds for every \(x\in \mathbb S^1\) 
that does not belong to the \(R\)-orbit of \( 0\), and at such a point,  \(i^ -\) and \(i^+\) are continuous.

\noindent - At every \(x\) in the \(R\)-orbit of \(0\), the function \(i^- \) is right-continuous (resp. \(i^+\) is left-continuous) and 
\[ \lim _{y \rightarrow x, y < x} i ^- (y) = i ^+ (x) \quad \text{   (resp.     }   \lim _{y \rightarrow x, y > x} i ^+ (y) = i ^- (x)). \]

\noindent - Given any two points $x,y$ in \(\mathbb S^1\), there is a sequence of integers \((k_n)_n\) such that \( R^{k_n}(x) \) 
tends to \(y\) from the right. The same property holds replacing right with left.     
\end{ejer}
\end{small}

\vsp

We next turn to the proof that, for a minimal subshift $\varphi$, the associated group $\mathrm{T}(\varphi)$ is finitely-generated.\footnote{In fact, 
this holds more generally for any subshift, independently of whether it is minimal or not. (The definitions of the suspension 
and the associated Thompson group work {\em verbatim} for that case.) See \cite{MT} for the details.}

\vsp

\begin{prop} \label{prop finite generation} {\em Let $(X,\varphi)$ be a minimal subshift with generating partition $\{C_1,\ldots, C_d\}$.  
If \(I\subset \mathbb R\) is a dyadic interval of length \(<1\), then the group $\mathrm{T} (\varphi)$ is 
generated by the groups $\widetilde{\mathrm{T}}$ and $\F_{C_i,I}$ for $i \in \{1,\ldots,d\}$. 
In particular, $\mathrm{T}(\varphi)$ is finitely-generated.}
\end{prop}

\vsp

\begin{small}\begin{rem} Together with Exercise \ref{T-no-fin-gen}, the preceding proposition implies that the adding machine 
is not a subshift. However, this can be directly checked as follows: Starting with any partition of \(\mathbb Z_2\) 
by clopen subsets, one can choose a sufficiently large integer \(m\) for which this partition can be refined by the one 
whose elements are the clopen sets \( C_i\) indexed by \, $i \in \mathbb{Z} / 2^m \mathbb{Z}$ \, and defined as 
$$C_i := \left\{ x= \sum _{k\geq 0} \varepsilon_k 2^k \in \mathbb Z_2 : \sum _{0\leq k\leq m-1} \varepsilon_k 2^k = i \right\}.$$
One readily sees that \( \varphi (C_i) = C_{i+1}\) for every \( i\in \mathbb Z/ 2^m\mathbb Z\). As a consequence, 
one cannot distinguish two points with the same itinerary with respect to the partition $\{C_1,\ldots,C_{2^m} \}$, 
and hence one cannot do it with the original partition either.
\end{rem}\end{small}

The next general lemma (that applies to all Cantor minimal systems $(X,\varphi)$) will be crucial for the proof of Proposition \ref{prop finite generation}.

\vsp 

\begin{lem} \label{lem intersection} {\em For any dyadic interval \(J\subset \mathbb R\) and any pair of clopen subsets 
$C,D$ of $X$ such that \((C,J)\) and \( (D,J)\) define dyadic flow boxes, the group $\langle \F_{C,J}', \F_{D,J}'\rangle$ 
contains the group $\F'_{C\cap D, J}$. If, moreover, \(C\) and \(D\) are disjoint, and \( (C\cup D , J)\) defines a dyadic 
flow box, then $\langle \F_{C,J}',\F_{D,J}'\rangle$ also contains $\F'_{C\cup D,J}$. }
\end{lem}

\noindent{\bf Proof.} We start with the following general remark: 
Given any two elements $f_0$ and $g_0$ in $\F'_J$, let $f\in \F'_{C,J}$ and $g\in \F'_{D,J}$ be defined by 
\( f := \F_{C,J} (f_0) \) and \( g := \F_{D,J} (g_0)\). One readily verifies that 
\[ [f,g]= \F_{C\cap D, J} ([f_0, g_0]) .\] 

Now let $h_0 \in \F_J'$ be arbitrary. Since \( \F_J'\) is simple (see Theorem \ref{Teo F' is simple}), 
we have $\F_J' = [ \F_J', \F_J']$.  In particular, we may write $h_0$ as a product of commutators 
$ h_0 = [f_1,g_1] \cdots [f_m,g_m]$, with each $ f_i,g_i$ in  $\F_J'$. By the remark just above, we have 
$$\F_{C \cap D, J} (h_0) =  [\F_{C,J} (f_1), \F_{D,J} (g_1) ] \cdots [ \F_{C,J} (f_m) , \F_{D,J} (g_m)),$$
which shows that $\F_{C \cap D, J} (h_0) $ belongs to $\langle \F_{C,J}', \F_{D,J}'\rangle$. Since this 
holds for all $h_0~\in~\F_J'$, we conclude that $\langle \F_{C,J}', \F_{D,J}'\rangle$ contains $\F'_{C\cap D, J}$,  
proving the first assertion.

For the second claim, given \( f\in \F'_{C\cup D, J} \), write it as \( f = \F_{C\cup D, J} (f_0)\), with \( f_0\in \F'_J\). 
Since $C$ and $D$ are disjoint, one obviously has 
\[ f = \F_{C, J} (f_0) \circ \F_{D,J} (f_0), \]
hence $f\in \langle \F_{C,J}',\F_{D,J}'\rangle$. This shows that  $\F'_{C\cup D,J}$ 
is contained in $\langle \F_{C,J}', \F_{D,J}'\rangle$
$\hfill\square$

\vspace{0.5cm}

\noindent{\bf Proof of Proposition \ref{prop finite generation}.} 
Let $H \subset \mathrm{T}(\varphi)$ be the subgroup generated by  $\widetilde{\mathrm{T}} $ and the subgroups $\F_{C_i,I}$ 
for $i \in \{ 1,\ldots,d\}$. We will show that $H$ contains all the subgroups \( \F_{C, J}\) with \(J\) of length \(<1\), 
which by Proposition \ref{p: generating Thompson copies} implies that $H$ coincides with $\mathrm{T}(\varphi)$. 
Note that the property to be proved is obviously equivalent to the property that $H$ contains all the subgroups 
\( \F'_{C, J}\) with \(J\) of length \(<1\), because of Exercise \ref{F:ejer-conmutator}. This is actually what we will show below.

The group \( \widetilde{\mathrm{T}}\) acts transitively on the set of dyadic intervals of \(\mathbb R\) of length \(<1\). Thus, 
conjugating by elements of $\widetilde{\mathrm{T}}$, we deduce that $H$ contains all copies $\F_{C_i,J}$ for every 
$i \in \{ 1,\ldots, d \}$ and every dyadic interval $J$ of length $<1$. The ({\em Equivariance}) property and \eqref{eq: equivariance dyadic box} 
show that $H$ also contains each subgroup of the form $\F_{\varphi^n(C_i),J}$ for every $n\in \Z$, every $i \in \{1,\ldots,d\}$, 
and every dyadic interval $J$ of length $<1$. Using the first part of Lemma \ref{lem intersection}, we deduce that 
$H$ contains each subgroup $\F'_{D,J}$ with $J$ any dyadic interval of length $<1$ and $D$ any clopen set that 
can be written as a finite intersection of the form $D=\bigcap_{j=1}^k \varphi^{n_j}(C_{i_j})$, with $n_j\in \Z$ and 
$i_j \in \{ 1,\ldots,d \}$. Now, the fact that $\{C_1,\ldots,C_d\}$ is a generating partition exactly means that each clopen 
set $C$ can be written as a finite union $C=D_1\cup \ldots \cup D_\ell$ of sets $D_i$ of this form. Hence, the second 
part of Lemma \ref{lem intersection} implies that $\F'_{C,J} \subset H$ for each clopen set $C$ and each dyadic 
interval $J$ of length $<1$, 
as we wanted to prove. 

The fact that $\mathrm{T}(\varphi)$ is finitely-generated follows since $\F_{C_i,I}\simeq \F$ and $\tilde{\mathrm{T}}$ 
are finitely-generated (see Exercise \ref{ej: finite generation tildeT} for the latter group). $\hfill\square$

\vspace{0.35cm} 

Putting together Propositions \ref{prop LO MB-T}, \ref{teo simple MB-T} and \ref{prop finite generation}, we finally conclude 
the following remarkable result of Matte Bon and Triestino.

\vsp 

\begin{thm} {\em If $(X,\varphi)$ is a minimal subshift, then the group 
$\mathrm{T} (\varphi)$ is left-orderable, simple, and finitely-generated.}
\end{thm}


\begin{footnotesize}

\end{footnotesize}

\printindex


\end{document}